\documentclass{book}
\usepackage[english]{babel}

\usepackage{amsmath, amssymb, amscd, amsthm}
\usepackage{xfrac}   
\usepackage{bbm}  

\usepackage{paralist}

\usepackage[colorlinks=true,linkcolor=darkgray,citecolor=darkgray]{hyperref}

\usepackage{xcolor}
\usepackage{tcolorbox}
\tcbuselibrary{theorems}
\tcbuselibrary{breakable}
\tcbuselibrary{skins}
\usepackage{float}

\usepackage{graphicx}
\usepackage{wrapfig}
\usepackage{caption}
\usepackage{subcaption}

\usepackage[abbreviations]{foreign}

\usepackage{tocbibind}

\usepackage{ifthen}
\usepackage{imakeidx}
\usepackage{multicol}

\makeindex

\def\D{\mathbb{D}}
\def\E{\mathbb{E}}
\def\N{\mathbb{N}}
\def\Z{\mathbb{Z}}
\def\Q{\mathbb{Q}}
\def\R{\mathbb{R}}
\def\C{\mathbb{C}}
\def\H{\mathbb{H}}
\def\SS{\mathbb{S}}
\def\T{\mathbb{T}}

\def\cA{\mathcal{A}}
\def\cB{\mathcal{B}}

\def\cD{\mathcal{D}}
\def\cE{\mathcal{E}}
\def\cF{\mathcal{F}}
\def\cI{\mathcal{I}}
\def\cG{\mathcal{G}}
\def\cH{\mathcal{H}}
\def\cL{\mathcal{L}}

\def\cP{\mathcal{P}}
\def\cQ{\mathcal{Q}}
\def\cR{\mathcal{R}}
\def\cS{\mathcal{S}}
\def\cT{\mathcal{T}}

\def\cT{\mathcal{T}}

\def\tr{\operatorname{tr}}
\def\det{\operatorname{det}}

\def\lcm{\operatorname{lcm}}

\def\Re{\operatorname{Re}}
\def\Im{\operatorname{Im}}
\def\Leb{\operatorname{Leb}}
\def\SL{\operatorname{SL}}
\def\GL{\operatorname{GL}}
\def\SO{\operatorname{SO}}
\def\PSL{\operatorname{PSL}}

\def\Sp{\operatorname{Sp}}
\def\P{\operatorname{P}}
\def\Aff{\operatorname{Aff}}
\def\Aut{\operatorname{Aut}}
\def\Inn{\operatorname{Inn}}
\def\Isom{\operatorname{Isom}}

\def\Dil{\operatorname{Dil}}
\def\Out{\operatorname{Out}}
\def\Tr{\operatorname{Tr}}
\def\Sp{\operatorname{Sp}}

\def\dist{\operatorname{d}}
\def\area{\operatorname{Area}}
\def\Area{\operatorname{Area}}

\def\deck{\operatorname{Deck}}
\def\Deck{\operatorname{Deck}}

\def\Sing{\operatorname{Sing}}
\def\Ends{\operatorname{Ends}}
\def\Trans{\operatorname{Trans}}
\def\sign{\operatorname{sign}}
\def\MCG{\operatorname{MCG}}
\def\Fix{\operatorname{Fix}}

\def\height{\operatorname{height}}
\def\circumference{\operatorname{circumference}}

\def\endray{\operatorname{end}}
\def\Emb{\operatorname{Emb}}

\def\MCG{\operatorname{MCG}}

\def\UTan{\operatorname{T}^1}

\def\hol{\operatorname{hol}}
\def\kerhol{\operatorname{kerhol}}

\def\epsilon{\varepsilon}
\def\codim{\operatorname{codim}}
\def\rk{\operatorname{rk}}

\newcommand{\leX}{\mathcal{L}^{\varepsilon}(M)}
\newcommand{\leXp}{\mathcal{L}^{\varepsilon'}(M)}
\newcommand{\eps}{\varepsilon}

\newcommand{\emphdef}[2][]{%
\ifthenelse{ \equal{#1}{} }%
 {\emph{#2}\index{#2|textbf}}%
 {\emph{#2}\index{#1|textbf}}%
}
\newcommand{\emphdefrecall}[2][]{%
\ifthenelse{ \equal{#1}{} }%
 {\emph{#2}\index{#2}}%
 {\emph{#2}\index{#1}}%
}

\newtcbtheorem[number within=section]{tcbtheorem}{Theorem}%
{
breakable,
enhanced,
theorem style=plain apart,
frame hidden,
frame style={fill=gray!20},
title style={color=gray!30},
fonttitle=\bfseries,coltitle=black,
fontupper=\itshape,colupper=black,
borderline west = {5pt}{0pt}{red!70}
}{thm}

\newtcbtheorem[use counter from=tcbtheorem]{tcbproposition}{Proposition}%
{
breakable,
enhanced,
theorem style=plain apart,
frame hidden,
frame style={fill=gray!20},
title style={color=gray!30},
fonttitle=\bfseries,coltitle=black,
fontupper=\itshape,colupper=black,
borderline west = {5pt}{0pt}{orange!70}
}{prop}

\newtcbtheorem[use counter from=tcbtheorem]{tcblemma}{Lemma}%
{
breakable,
enhanced,
theorem style=plain apart,
frame hidden,
frame style={fill=gray!20},
title style={color=gray!30},
fonttitle=\bfseries,coltitle=black,
fontupper=\itshape,colupper=black,
borderline west = {5pt}{0pt}{orange!70}
}{lem}

\newtcbtheorem[use counter from=tcbtheorem]{tcbcorollary}{Corollary}%
{
breakable,
enhanced,
theorem style=plain apart,
frame hidden,
frame style={fill=gray!20},
title style={color=gray!30},
fonttitle=\bfseries,coltitle=black,
fontupper=\itshape,colupper=black,
borderline west = {5pt}{0pt}{orange!70}
}{cor}

\newtcbtheorem[use counter from=tcbtheorem]{tcbdefinition}{Definition}%
{
breakable,
enhanced,
theorem style=plain apart,
frame hidden,
frame style={fill=gray!20},
title style={color=gray!30},
fonttitle=\bfseries,coltitle=black,
fontupper=\itshape,colupper=black,
borderline west = {5pt}{0pt}{blue!50}
}{def}

\newtcbtheorem[use counter from=tcbtheorem]{tcbremark}{Remark}%
{
breakable,
enhanced,
theorem style=plain apart,
frame hidden,
frame style={fill=gray!20},
title style={color=gray!30},
fonttitle=\bfseries,coltitle=black,
fontupper=\itshape,colupper=black,
borderline west = {5pt}{0pt}{green!40}
}{rk}

\newtcbtheorem[use counter from=tcbtheorem]{tcbquestion}{Question}%
{
breakable,
enhanced,
theorem style=plain apart,
frame hidden,
frame style={fill=gray!20},
title style={color=gray!30},
fonttitle=\bfseries,coltitle=black,
fontupper=\itshape,colupper=black,
borderline west = {5pt}{0pt}{green!40}
}{qu}

\newtcbtheorem[use counter from=tcbtheorem]{tcbexercise}{Exercise}%
{
breakable,
enhanced,
frame hidden,
breakable,
frame style={fill=gray!20},
fonttitle=\bfseries,coltitle=black,
colupper=black,
}{exo}

\newtcbtheorem[use counter from=tcbtheorem]{tcbexample}{Example}%
{
breakable,
enhanced,
breakable,
frame hidden,
frame style={fill=gray!10},
fonttitle=\bfseries,coltitle=black,
colupper=black,
}{exa}

\newenvironment{psmallmatrix}
  {\left(\begin{smallmatrix}}
  {\end{smallmatrix}\right)}

\definecolor{darkgreen}{rgb}{0,0.7,0}

\newcounter{EnumerateCounterContinuation}

\title{Infinite Translation Surfaces in the Wild\\ {\large Volume 1: Geometry and Symmetries}}
\author{Vincent Delecroix \\ Pascal Hubert \\ Ferr\'an Valdez}

\begin{document}
\frontmatter
\maketitle
\tableofcontents


\chapter*{Preface}
\label{ch:Preface}

 Most probably, hearts brave enough to open this book have already encountered the charm of polygonal billiards. However, paraphrasing the Captain of the White Tower at the time of the War of the Ring, one does not simply walk into the realm of translation surfaces without speaking about them. Polygonal billiards are dynamical systems describing the frictionless movement of a point inside a polygon $P$ which enjoys elastic collitions with the sides, except at the corners, where movement stops as if there was an infinitesimal pocket. In simpler terms, this means that upon colliding with a side of $P$, the point continues its trajectory following the optical law: \emph{angle of incidence equals angle of reflection}. A point which after a finite number of collisions comes back to its initial position with the same direction on which its movement started describes a \emph{periodic trajectory}. It is fair to say that one of the simplest questions one could ask about billiards is whether for any polygon $P$ periodic billiard trajectories exist. For instance, in every acute triangle, the triangle formed by projecting the orthocenter to the sides yields a periodic billiard trajectory. For right triangles, almost any point hitting the hypotenuse at a 90 degree angle describes a periodic trajectory. Surprisingly, finding periodic trajectories for arbitrary obtuse triangles has remained up to date a very hard task. No one seems to know why. On the matter Richard Schwartz, one of the leading experts in this area, opined:\footnote{For more details, see \href{https://mathoverflow.net/posts/344448/revisions}{https://mathoverflow.net/posts/344448/revisions}}\\
\\
\emph{I am not sure why it is so hard. All I can really say is that, after a lot of experimentation, I can't really see any pattern to it. It might be hard in the same way that building the fountain of youth is hard: nobody has any idea what to do.
Sometimes I have described the problem as being akin to riding your bicycle to the North Pole. You know in advance that something is going to go bad for you, but it is hard to know what that will be exactly.\\}

Billiards and translation surfaces are related by an unfolding procedure, also know as the unfolding trick, which is attributed to Fox and Kershner~\cite{FoxKershner36}. It was rediscovered 40 years later by Katok and Zemljakov~\cite{KatokZemliakov75}. This procedure allows to reinterpret the dynamics of the ball as a flow on a surface $S_P$ tiled by isometric copies of $P$. In particular, the existence of a closed trajectory for this flow implies the existence of a periodic trajectory for the billiard. An eager reader can check the details of the unfolding trick in Section~\ref{ssec:PolygonalBilliards}. There are three important facts to retain about the surface $S_P$. First, if one turns a blind eye to what happens near points corresponding to vertices of $P$, \emph{every transition function of $S_P$ is a translation of the plane}. This is the property that defines translation surfaces. Second, if each interior angle of $P$ is a rational multiple of $\pi$, case in which $P$ is called \emph{rational}, the surface $S_P$ is compact, without boundary and its genus depends on the aforementioned angles. As explained by Masur, one can use this fact and the language of quadratic diffentials to show that every billiard on a rational polygon has a periodic trajectory~\cite{Masur86}. Finally, and most importantly for this book, when $P$ is not rational $S_P$
is an infinite genus translation surface whose topology does not depend on the interior angles of $P$~\cite{Valdez09_2}. Infinite genus translation surfaces are a particular case of infinite-type translation surfaces, see Definition~\ref{def:FiniteTypeInfiniteType}.\\
\\
\indent In short: \emph{behind a typical polygon lies an infinite-type translation surface}.\\
\\
This book is a modest attempt to lay the foundations for the study of infinite-type translation surfaces. It discusses old and recent results on the subject, but also presents new Theorems and their applications. It is divided in four Chapters, whose content we brielfy describe in the following paragraphs.\\

Chapter 1 is an extended introduction and is divided into two parts. It starts with a detailed discussion regarding three different definitions of what a general translation surface is. Roughly speaking, translation surfaces can be defined in three ways: \emph{constructively}, by gluing polygons together; \emph{analytically}, through pairs $(X,\omega)$ where $X$ is a Riemann surface and $\omega$ is a non-identically zero holomorphic 1-form; or \emph{geometrically}, as structures on a topological surface where transition functions are translations. One of the main contributions of this Chapter is to show that, in the infinite-type setting, these definitions are indeed equivalent. We then briefly discuss several objects and structures one can attach to a translation surface $M$, \emph{e.g.} the metric completion  w.r.t. the intrinsic path metric (coming from its natural flat Riemannian metric), singularities or group of affine homeomorphisms. The second part of Chapter 1 is a compendium of examples on which translation surfaces naturally appear. These include polygonal billiards, problems arising from theoretical physics such as wind-tree models or Eaton lenses, and homogeneous holomorphic foliations.

Chapter 2 starts discussing in detail the topological classification of topological infinite-type surfaces, \emph{i.e.} those surfaces whose fundamental groups are not finitely generated.
Then, in Section~\ref{sec:CoveringSpaces} it addresses  infinite-type translation surfaces arising from covering spaces. Special focus is given to those $p:\widetilde{M}\to M$ where $M$ is of finite type, and $\Deck(p)$ is an abelian torsion-free group of finite rank $d$. In particular, we discuss how these are in correspondence with subspaces of dimension $d$ of the first homology group of $M$ with rational coefficients. This kind of coverings are revisited on Chapter 3. Section~\ref{Sec:Existence-translation-structure} raises several natural questions: Does every topological surface $S$ admit a translation surface structure? Can this structure have finite area? Is it possible to freely prescribe a set of conical singularities?
Finite-type translation surfaces are characterized by the presence of only conical singularities. On the other hand, infinite-type translation surfaces can have singularities which are not conical. These are called \emph{wild} singularities and appear in the metric completion w.r.t. the intrinsic path metric coming from its natural flat (Riemannian) metric. Chapter 2 concludes with Section~\ref{sec:Singularities}, a detailed discussion about wild singularities, their affine invariants and some interesting examples.

Chapter 3 deals with symmetries of translation surfaces. This is the largest and probably the densest Chapter of the book. In this chapter, we study elements of $\Aff(M)$, the group of affine homeomorphisms of a translation surface $M$. We also examine $\Trans(M)$, the group of 'translations' of $M$ (which is the subgroup of $\Aff(M)$ formed by those affine homeomorphisms whose derivative is the identity) or $\Gamma_M<\GL(2,\R)$, the group formed by the derivatives of elements of $\Aff(M)$. This latter is also known as the Veech group of $M$ and, in the case when $M$ has finite topological type, gives necessary conditions for the translation flow of $M$ to have 'optimal dynamics', see~\cite{Veech89}.
In Section~\ref{sec:AnalyticAutomorphisms} we give a classification of translation surfaces for which $\Trans(M)$ does not act properly discontinuously. Then, in Section~\ref{sec:EveryGroupIsVeech} we show that if there is no restriction on the area of $M$, then basically any countable subgroup of $\GL(2,\R)$ can be the Veech group of an translation surface of infinite genus. This builds on previous results by Przytycki, Weitze-Schmithuesen and the third author. Section~\ref{sec:AffineGroupsCoverings} discusses Affine and Veech groups in the context of an infinite covering $p:\widetilde{M}\to M$ where the base is a finite-type translation surface. We focus the discussion in the case where $\Deck(p)$ is an infinite torsion free Abelian group. As mentioned before, these are in correspondence with subspaces
$V<H_1(M,\Sigma;\mathbb{Q})$, where $\Sigma$ is the set of conical singularities. In particular we show that every element in $\Aff(M)$ which leaves $V$ invariant has a lift, and thus
the study of Affine and Veech groups in this context amounts to the study stabilizers of vector subspaces of homology. We pay special attention to covering defined by the subspace $\kerhol(M)$, which is the kernel of the map that associates to each homology class $c$ the complex number $\int_c\omega$, where $M=(X,\omega)$. In particular, we give criteria for parabolic elements to lift, \emph{e.g.} the so-called one-cylinder trick, see Theorem~\ref{thm:OneCylinderTrick} for a precise stament. We also discuss an important result of Hooper and Weiss that gives conditions for a covering defined by $\kerhol(M)$ over a genus 2 surface to have a Veech group of the First kind. Section~\ref{sec:AffineGroupsCoverings} finishes discussing Veech groups of wind-tree models, which are also examples infinite-type translation surfaces defined by infinite abelian coverings. To be specific, one of the main aspects we discuss are parameters defining the wind-tree model for which the Veech group is infinitely generated.

As detailed in the first Chapter, there is a classical construction by Thurston and Veech
which allows to construct a translation (or half-translation) surface $M$ on any finite-type surface $S$ with non-abelian fundamental group which admits a horizontal and a vertical cylinder
decomposition where all cylinders involved have the same modulus. As a consequence the surface $M$ has a lot of affine pseudo-Anosov homeomorphisms. Using ideas of Hooper and the theory of harmonic functions on infinite graphs, in Section~\ref{sec:HooperThurstonVeechConstruction} we detail an extension of this construction. It details how to construct, for any infinite-type surface $S$, a translation (or half-translation) surface structure whose Veech group contains at least a Fuchsian group of the second kind generated by two parabolic matrices. We also provide multiple examples to illustrate this construction. Chapter 3 concludes with Section~\ref{ssec:VeechGroupBakersSurface}, in which we
explicitly calculate the Veech group of a family of examples: the Chamanara (or baker) surfaces.

One of the main features of translation surfaces is that for every direction $\theta$ there is a well defined vector field formed by unit vectors parallel to $\theta$. In particular, for every such direction there exists a flow parallel to $\theta$ defined on a subset of full dimension (again, singularities create some small problems). These flows are mentioned through the book, but it is until Chapter~\ref{chap:TranslationFlowsAndIET} that we begin their study. This Chapter is intended to be an intermezzo between this book and~\cite{DHV2}, which will be dedicated to dynamical aspects of infinite type translation surfaces.
In Section~\ref{sec:IETAndTranslationFlows}, we begin by discussing the preliminaries, focusing on infinite interval exchange transformations (IETs). The main result in this direction is that the first return map to a finite length transversal of the translation flow on a general translation surface is conjugate to a (possible infinite) IET. We explain how to encode the dynamics of an infinite IET using two total orders on an at most countable alphabet and a (summable) length sequence. Given that infinite IETs are transformations on finite length intervals, first return maps to infinite length transversal cannot be described with IETs. We also discuss several of these more general first return maps. They appear, for example, when studying the  Hooper-Thurston-Veech surfaces.
Section~\ref{sec:Jungle} is intended to serve as a warning of the following fact, due to Arnoux, Ornstein and Weiss (see Theorem~\ref{thm:VershikArnouxOrnsteinWeiss}) : every aperiodic measure-preserving transformation is measurably isomorphic to a, possibly infinite, interval exchange transformation on $(0,1)$. In plain words: anything can happen in the world of infinite IETs.
The content of this Section provides the necessary tools to understand the proof of Arnoux-Orstein-Weiss's result: Rokhlin's towers, cutting and stacking constructions for IETs and the Kakutani-Rokhlin Lemma (see Lemma~\ref{lem:Rokhlin}). We also discuss briefly Bratteli-Vershik diagrams. As an illustration of how things might be different in the world of infinite IETs, we provide in Section~\ref{sec:CounterexampleKeane} a counterexample to Keane's criterion. More precisely, we explain how to construct an explicit infinite IET without connections which is not minimal. Chapter 4 finishes with Section~\ref{sec:entropy}, where we discuss (metric) entropy in general and then we focus on IETs. The main result of this section is Theorem~\ref{thm:EntropyInfiniteIET}, which provides an upper bound for the entropy of an infinite IET in terms of the lengths of the subintervals it permutes. We give several applications of this formula. The chapter concludes with a short discussion on Abramov's formula, which relates the entropy of the translation flow with the entropy of the first return map to a global transversal.\\


\textbf{Disclaimers}
\begin{enumerate}
\item This volume does not include most of the known results regarding the dynamical aspects of infinite-type translation surfaces. These topics will be extensively covered in the forthcoming second volume, ~\cite{DHV2}.
\item While most theorems in this book are presented in their original form, some statements have been enhanced for clarity or comprehensiveness, such as Theorem~\ref{thm:AnyGroupIsVeech}.
\end{enumerate}

\textbf{Around this book}. More material about infinite translation surfaces, including errata, solution to some exercises and an updated bibliography, is available at \href{https://www.labri.fr/perso/vdelecro/infinite-translation-surfaces-in-the-wild.html}{https://www.labri.fr/perso/vdelecro/infinite-translation-surfaces-in-the-wild.html} .  If you have any kind words (or complaints), solutions to exercises or maybe a poem, you can send them to \href{mailto:tsw@matmor.unam.mx}{tsw@matmor.unam.mx}.

\textbf{Acknowledgements}. We would like to thank the following people for reading early versions of this manuscript: Christoph Karg, Jan Kohlm\"{u}ller, Anja Randecker, Garrett Proffitt, Jenya Sapir, Carlo Schmid, Barak Weiss, Christopher Zhang. We would like to thank the following institutions for their support during all these years this proyect: Centro de Modelacion Matem\'atica, Max-Planck-Institut fuer Mathematik, Centro de Ciencias\linebreak Matem\'aticas UNAM, Centre International de Rencontres Math\'ematiques, Laboratoire Bordelais de Recherche en Informatique, Mathematisches\linebreak Forschungsinstitut Oberwolfach, Simons Laufer Mathematical Sciences Institute, Casa Matem\'atica Oaxaca and Karslruhe Institute of Technology.
The third author would like to thank PAPIIT IN102018 \& IN101422, PASPA de la DGAPA, and
Fondo CONACYT-Ciencia Básica’s project 283960 for financial support, and Marion and Anika Yunuen for emotional support.

\mainmatter

\chapter{Introduction}
\label{CH:INTRODUCTION}
\label{ch:Introduction}
\label{chap:Introduction}

Very much like Tomoyuki Tanaka's \emph{Kingu Gidora}\footnote{Tomoyuki Tanaka, born on April 26, 1910, in Osaka Prefecture, was a film producer. He is known for being one of the creators of \emph{Gojira}, the most famous of all monsters in the \emph{kaiju} genre. Kingu Gidora, one of Tanaka's creation is, according to the Wikipedia, \emph{an armless, bipedal, golden and yellowish-scaled dragon with three heads, two wings in the shape of fans to fly and two tails}.} , translation surfaces have three different, but related, faces. This polyvalence gives the world of translation surfaces a peculiar and rich essence. In the first section of this introductory chapter, we discuss these three different definitions, prove their equivalences and introduce the main associated geometrical and dynamical objects. Unlike Kingu Gidora, translation surfaces are not extraterrestrial or genetically engineered monsters. In the second section, we present infinite translation surfaces in the wild (mathematical) world. They appear in dynamical systems such as polygonal billiards, wind-tree models or baker's map as well as in geometric constructions related to surface homeomorphisms such as the Thurston-Veech's construction. Most results stated of this second section are proven in subsequent chapters of this book or in the second volume~\cite{DHV2} and, hopefully, will serve as a motivation for the reader.

\section{What is a translation surface?}
\label{sec:WhatIsATranslationSurface}
\label{SECC:WhatIsATranslationSurface}

\subsection{Three examples}
\label{ssec:ThreeBabyExamples}

A translation surface is a flat object constructed by gluing polygons along parallel sides of the same length using translations. Before introducing formal definitions let us discuss three illustrative examples.

\smallskip

\noindent \textbf{A flat torus}. Consider the unit square in $\C$ given by $0\leq \Re(z), \Im(z) \leq 1$ and identify parallel sides using translations. That is, if $z=x+\sqrt{-1}y$, we identify the points $x$ in the lower side of the square with $x+\sqrt{-1}$ and the points $\sqrt{-1}y$ in the left side of the square with $1+\sqrt{-1}y$. The result is a surface $M$ homeomorphic to a torus and which has the special property that \emph{every} point $z\in M$ has a small  neighbourhood isometric to a neighbourhood of the origin in $\C$.

The flat torus can be alternatively obtained by considering the quotient of Abelian groups $\C / (\Z \oplus \sqrt{-1} \Z)$.

\bigskip

\begin{figure}[!ht]
\begin{center}
\includegraphics{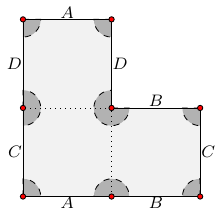} \\
\includegraphics{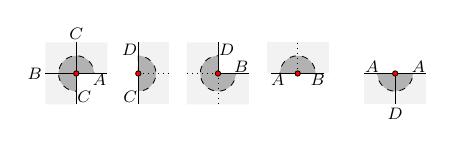}
\caption{A genus 2 translation surface and a neighbourhood of its conical singularity.}
\label{fig:LShapedOrigami}
\end{center}
\end{figure}
\noindent \textbf{A genus 2 surface}. Consider three copies of the unit square, glue them and label their sides as depicted in Figure~\ref{fig:LShapedOrigami}. If we identify edges with the same labels using translations, the result is a genus 2 surface\footnote{This can be shown easily by calculating the Euler characteristic of the surface: $2-2g = V - E + F$ where $g$ is the genus,  $V = 1$ is the number of vertices, $E=6$ is the number of edges and $F=3$ the number of faces.}. In this surface all points \emph{except one} have a small neighbourhood isometric to a neighbourhood of the origin in $\C$. The ``problematic'' point, which we denote by $p_0$, appears because all vertices are merged by the identifications into a single point. In the lower part of Figure~\ref{fig:LShapedOrigami} we illustrate a small neighbourhood $U_{p_0}$ of $p_0$. Remark that $U_{p_0}$ is not isometric to a neighbourhood of the origin in $\C$. It is however isometric to the space obtained by gluing cyclically 3 copies of a neighbourhood of the origin in $\C$, which is an example of a \emph{ramified covering} of degree 3. For this reason, $p_0$ is called a \emph{conical singularity} of total angle $6\pi$.

Both the flat torus and the genus 2 surface are examples of \emph{finite-type} translation surfaces.

\begin{figure}[!ht]
\begin{center}
\includegraphics{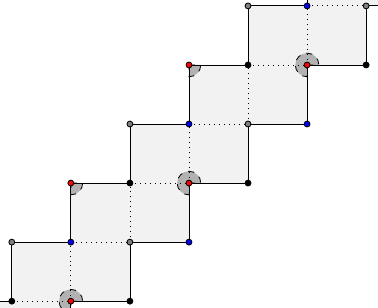}
\end{center}
\caption{The infinite staircase. Opposite sides are identified. There are four infinite degree vertices in the surface.}
\label{fig:StaircaseFirst}
\end{figure}

\noindent \textbf{The \emphdef{infinite staircase}}. Consider a countable family of squares and identify their parallel sides as depicted in Figure~\ref{fig:StaircaseFirst}, where pairs of opposite (parallel) sides are identified using translations. Every point which is not a vertex  has a  neighbourhood that is isometric to a neighbourhood of the origin in $\C$. After identifications all vertices involved merge into four points.
In figure \ref{fig:StaircaseFirst} below we depict with a dashed line the boundary of a small neighbourhood of one on these four points, which we denote by $U_{z_0}$ and $z_0$ respectively. Remark that $U_{z_0}$ can be constructed by gluing cyclically \emph{infinitely} many copies of a neighbourhood of the origin in $\C$. However, $z_0$ is worse than the problematic point from the preceding example because
\emph{$z_0$ does not have a compact neighbourhood}.
In other words, because of ``very problematic'' points like $z_0$ the topological space that
we have constructed \emph{is not locally compact}, and hence not a surface. For this reason we remove all the vertices of the squares involved in the construction. The result is an \emph{infinite-type} translation surface called the infinite staircase. The nomenclature in this case is justified because,
as we will see later, this surface has infinite genus.

These three examples will appear all along this text. Flat tori have been studied since the 19th and early 20th centuries, by L.~Kronecker and H.~Weyl, among others. The surface of genus 2 in the second example is a particular case of a compact translation surface. Such translation surfaces have been intensively studied since the 1970's and their theory is well-developed. This book assumes some basic knowledge of the theory of compact translation surfaces and we refer the reader to the following four references when needed: H.~Masur, S.~Tabachnikov~\cite{MasurTabachnikov02}, A.~Zorich~\cite{Zorich06}, J.-C.~Yoccoz~\cite{Yoccoz10}, G.~Forni, C.~Matheus~\cite{ForniMatheus14} or A. Wright~\cite{Wright15}.

\subsection{Three definitions}
\label{ssec:TranslationSurfaceDefinitions}
In this section we define what a translation surface is in three different ways: constructive, geometric and analytic. Each definition has its own pros and cons depending on the context in which translation surfaces are studied.
\smallskip

We start with the constructive definition, which generalizes the three examples we presented before. This is the definition that will be used to present most examples in Section~\ref{sec:ExamplesInfiniteTranslationSurfaces}. First let us define an \emphdef{Euclidean polygon} as a simply connected and bounded closed set in the Euclidean plane whose boundary is a curve formed by finitely many segments. Let $\cP$ be an at most countable family of Euclidean polygons and  $E(\cP)$ be the set of all the edges in $\cP$. For each edge $e \in E(\cP)$ we consider $n_e$ the (unit) vector normal to $e$ which points toward the interior of the polygon having $e$ as a side. Suppose that there exists $f: E(\cP) \to E(\cP)$ a pairing (that is an involution without fixed point) such that for every $e \in E(\cP)$ the edges $e$ and $e' = f(e)$ differ by a translation $\tau_{e}$ and $n_{e'} = - n_{e}$.

Let $\bigsqcup_{P \in\cP} P$ be the disjoint union of the polygons in $\cP$. For every $e \in E(\cP)$ seen as a subset of $\bigsqcup_{P \in\cP} P$, we identify the points in $e$ with the points in $f(e)$ using $\tau_e$. Note that each point in the interior of an edge $e$ is identified with exactly one point in $f(e)$. After these identifications, we obtain a topological space where the quotient map $\pi: \bigsqcup_{P \in\cP} P \to (\bigsqcup_{P \in\cP} P) / \sim$ is injective in the interior of each polygon and 2-to-1 on the edges. The situation for vertices can be more complicated. For example in the L shaped surface of Figure~\ref{fig:LShapedOrigami} the map $\pi$ sends all $8$ vertices involved in the construction to the same point. On the other hand, for the infinite staircase in Figure~\ref{fig:StaircaseFirst} $\pi$ sends infinitely many vertices to the same point. For this reason we say that a vertex $v \in P \in \mathcal{P}$ is \emph{of finite degree} if $\pi^{-1}(\pi(v))$ is finite and \emph{of infinite degree} otherwise. This notion takes care of vertices which merge into ``very problematic'' points like the ones we created when constructing the infinite staircase (see Lemma~\ref{lem:ConicalPointConstructive}).
\begin{tcbdefinition}{Constructive}{TranslationSurfaceConstructive}
Let $\cP$ be an at most countable set of Euclidean polygons and $f: E(\cP) \to E(\cP)$ a pairing as above. Let $M$ be $ \bigsqcup_{P\in\cP} P/\sim$ with all vertices of infinite degree removed. If $M$ is connected, we call it \emph{the translation surface obtained from the family of polygons $\cP$ and pairing $f$}.
\end{tcbdefinition}

\begin{tcbexercise}{}{}
In this simple exercise we discuss the connectedness of $M$ (in Definition~\ref{def:TranslationSurfaceConstructive}) in terms of the combinatorics of the edge pairing. Let $(\cP, f: E(\cP) \to E(\cP))$ be a family of polygons and a pairing, and let $M$ be the translation surface obtained from the family of polygons $\mathcal{P}$. For each edge $e$ we denote $P_e$ the polygon in $\cP$ to which $e$ belongs.
\begin{compactenum}
\item Prove that $M$ is connected if and only if for each pair $P,Q \in \cP$ there exists a finite sequence of polygons and edges $(P_{e_0},e_0,P_{e_1},e_1, \ldots, e_{n-1},P_{e_n})$ where $P_{e_0} = P$, $P_{e_{n}} = Q$ and for $i \in \{0, 1, \ldots, n-1\}$ we have $P_{f(e_i)} = P_{e_{i+1}}$.
\item Show that the above property is equivalent to the connectivity of a graph built from $(\cP,f)$ and whose vertex set is the set of polygons $\cP$.
\end{compactenum}
\end{tcbexercise}

The following lemma and exercise describe the local geometry of a translation surface around a vertex of finite degree.
\begin{tcblemma}{}{ConicalPointConstructive}
Let $M$ be the translation surface generated by a family of polygons $\cP$. For each vertex $v\in P\in\mathcal{P}$ we denote by $\alpha_v\in(0,2\pi)$ the interior angle of $P$ at $v$. Then for each vertex $v\in P$ of finite degree there exists a positive integer $k_v \in \{1, 2, \ldots\}$ such that
\begin{equation}
\label{eq:AdmissibleVertices}
\sum_{w\in\pi^{-1}(\pi(v))} \alpha_w = 2 k_v\pi.
\end{equation}
If $k_v > 1$, the point $\pi(v) \in M$ is called a \emphdef{conical singularity} of angle $2 k_v\pi $ while if $k_v=1$ it is called a \emphdef{regular point}.
\end{tcblemma}

Since the Euclidean metric $dx^2 + dy^2$ in the plane is invariant under translation, any translation surface built from polygons inherits a Euclidean metric that is well-defined in the complement of conical singularities. In particular,
there is a well-defined notion of distance and area. The following elementary exercise describes the behavior of a metric at a vertex of finite degree.
\begin{tcbexercise}{}{ConicalSingularities}
This exercise discusses conical singularities of a Euclidean metric. It is very much inspired by the first section of~\cite{Troyanov86}. We use $x,y$ for the standard coordinates in $\R^2$, $z = x + i y$ the corresponding number in $\C$ and $(r,\theta)$ for polar coordinates $x = r \cos(\theta)$ and $y = r \sin(\theta)$ (or $z = r \exp(i \theta)$).
\begin{enumerate}
\item Show that the Euclidean metric $dx^2 + dy^2$ can also be written as $dz d\bar{z}$ or $(dr)^2 + (r d\theta)^2$.
\setcounter{EnumerateCounterContinuation}{\theenumi}
\end{enumerate}
Let us consider the punctured plane $\C^* = \C \backslash \{0\}$ with the metric $g_\alpha = (dr)^2 + (\alpha r d\theta)^2$ where $\alpha > 0$ is a parameter. The metric $g_1$ is nothing else than the Euclidean metric.
\begin{enumerate}
\setcounter{enumi}{\theEnumerateCounterContinuation}
\item Show that a circle of radius $r > 0$ around the origin has length $2\pi r \alpha$ in the metric $g_\alpha$.
\item Using the change of coordinates $\theta' = \alpha \theta$ show that $g_\alpha$ is a flat metric, \ie a metric with zero curvature.
\item \label{ques:ZKAsRamifiedCover} When $\alpha = k$ is a positive integer show that the map $z \mapsto z^k$ is a local isometry from $(\C^*, g_k)$ to $(\C^*, g_1)$. Deduce that $g_k = |z|^{2(k-1)} dz d\bar{z}$.
\item Show that for any $\alpha>0$ the metric completion of $\C^*$ with respect to $g_\alpha$ is the whole plane $\C$.
\item \label{exo:item:ProlongationConicalMetric} Show that the metric $g_\alpha$ extends to a metric on $\C$ if and only if $\alpha = 1$.
\setcounter{EnumerateCounterContinuation}{\theenumi}
\end{enumerate}
The origin in $\C$ is called a \emphdefrecall{conical singularity} of angle $2 \pi \alpha$ for the metric $g_\alpha$. More generally, if $S$ is a surface with a flat metric $g$ defined on the complement of a discrete set $\Sigma \subset S$, we say that a point $p \in \Sigma$ is a conical singularity of \emphdefrecall[angle of conical singularity]{angle} $2 \pi \alpha$ for the metric $g$ if there exists an open neighborhood $U$ of $p$ such that $(U,g)$ is isometric to $(D'_r, g_\alpha)$ where $D'_r \subset \C$ is the punctured disk of radius $r$ centered at the origin.
\bigskip

We now provide another construction of the cone $(\C^*, g_\alpha)$.
\begin{enumerate}
\setcounter{enumi}{\theEnumerateCounterContinuation}
\item Let $m$ be the metric on $\C$ defined by $m = |\exp(z)|^2 (dz d\overline{z})$. Show that the exponential map $\exp: \C \to \C^*$ is a local isometry between $(\C, m)$ and $(\C^*, g_1)$.
\item Show that the translations by purely imaginary numbers are isometries in $(\C, m)$.
\item Show that the quotient $\C / (2 \pi i \alpha \Z )$ with the quotient metric $\overline{m}$ is isometric to $(\C^*, g_\alpha)$.
\item Show that the metric completion of $(\C, m)$ is obtained by adding one point $p_\infty$. \textit{Hint: $p_\infty$ is the limit of any horizontal ray $z_t = -t + i y$ with $t \to +\infty$ and y fixed.}
\item Show that the map $\exp$ extends to a map $\overline{\exp}:\C \cup \{p_\infty\}\to\C^* \cup \{0\} $.
\setcounter{EnumerateCounterContinuation}{\theenumi}
\end{enumerate}
The point $p_\infty$ in $(\C, m)$ is called an \emphdefrecall{infinite angle singularity}.
\begin{enumerate}
\setcounter{enumi}{\theEnumerateCounterContinuation}
\item \label{exo:item:FiniteDegreeAndConicalSingularity} Let $M$ be the translation surface generated by a family of polygons $\cP$. Show that for each vertex of finite degree $k_v$, the point $\pi(v)$ is a conical singularity of angle $2 k_v \pi$ of the Euclidean metric on $M$. \textit{Hint: take a look at Lemma~\ref{lem:ConicalPointConstructive}}.
\item Let $M$ be the infinite staircase constructed in Section~\ref{ssec:ThreeBabyExamples}. Show that each vertex gives rise a point which is locally isometric to the infinite angle singularity $p_\infty$ defined above.
\end{enumerate}
\end{tcbexercise}

\begin{tcbremark}{}{}
The map $z \mapsto z^k$ considered in question~\ref{ques:ZKAsRamifiedCover} of Exercise~\ref{exo:ConicalSingularities} is actually the local model for ramified covers of surfaces.
\end{tcbremark}

We now consider the geometric definition of a translation surface. As we saw, a (constructive) translation surface carries a natural flat metric with conical singularities. However, as we will see later, not all surfaces with a flat metric are translation surfaces. The restriction is related to the so-called holonomy that will be discussed in Section~\ref{ssec:Structures}. A \emphdef{translation atlas} on a topological surface $S$ is a set of maps $\cT = \{\phi_i:U_i \to \C\}$ where $(U_i)_{i\in\N}$ forms an open covering of $S$, each $\phi_i$ is a homeomorphism from $U_i$ to $\phi(U_i)$ and for each $i,j$ the transition map $\phi_j \circ \phi_i^{-1}: \phi_i(U_i \cap U_j) \to \phi_j(U_i \cap U_j)$ is a translation in $\C$. Any topological surface with a translation atlas is naturally endowed with a flat metric. Indeed, this metric is obtained by pulling back the (translation invariant) Euclidean metric in $\C$.

\begin{tcbdefinition}{Geometric}{TranslationSurfaceGeometric}
A \emphdef[translation surface]{(geometric) translation surface} is a pair $(S, \cT)$ made of a connected topological surface $S$ and a maximal translation atlas $\cT$ on $S\setminus\Sigma$, where:
\begin{enumerate}
\item $\Sigma$ is a discrete subset of $S$ and
\item every $z\in\Sigma$ is a conical singularity.
\end{enumerate}
The maximal translation atlas $\cT$ is called a \emphdef{translation structure} on $S$ and its charts are called the \emphdef{flat charts} or \emphdef{flat coordinates}.
\end{tcbdefinition}

\begin{tcbexercise}{}{}
Prove that the angle of a conical singularity $z\in \Sigma$ in the definition above is of the form $2\pi k$, for some positive integer $k>1$.
\end{tcbexercise}

\begin{tcbremark}{}{GXStructure}
The pair $(S \setminus \Sigma, \cT)$ is a particular case of a geometric structure in the sense of Thurston locally modeled on the space $\C$ endowed with its group of translations $\C$. For more details see Appendix~\ref{Appendix:GXStructures}.
\end{tcbremark}

Finally, we introduce the analytic definition. It provides a tight link between translation surfaces and complex algebraic geometry that permits the use of very powerful tools. This third definition is central in many important results on compact translation surfaces such as the Eskin-Kontsevich-Zorich formula. The conformal structure already appears in the proof of Theorem~\ref{thm:EquivalencesInDefinitions} and will be important when we study isometries in Section~\ref{sec:AnalyticAutomorphisms}.
The reader not familiar with Riemann surfaces might choose to skip the sections using this point of view or consult one of the many textbooks available on the subject such as~\cite{Forster77}, \cite{FarkasKra80}, \cite{Donaldson11} or~\cite{Hubbard-book1}.

\begin{tcbdefinition}{Analytic}{TranslationSurfaceAnalytic}
An \emphdef[translation surface]{(analytic) translation surface} is a pair $(X,\omega)$ formed by a (connected) Riemann surface $X$ and a holomorphic 1-form (\emph{a.k.a.} Abelian differential) $\omega$ on $X$ which is not identically zero.
\end{tcbdefinition}

\begin{tcbexercise}{}{ConicalSingularityAndZeroOfDifferential}
This exercise is a natural counterpart of Exercise~\ref{exo:ConicalSingularities}.
\begin{enumerate}
\item Let $(X,\omega)$ be an analytic translation surface. Show that the tensor $dz \overline{dz}$ is a flat metric on $X \setminus Z(\omega)$ where $Z(\omega)$ is the set of zeros of $\omega$ and $z$ is a local coordinate for $X$.
\item Show that each $p \in Z(\omega)$ is a conical singularity for the flat metric and determine the angle in terms of
the order of $p$ as a zero of $\omega$. \textit{Hint: there exists a non-negative integer $k$ and a holomorphic coordinate $z$ at $p$ such that $\omega = z^k dz$.}
\item Show that the area form of $(X,\omega)$ is $\frac{i}{2} \omega \wedge \overline{\omega}$.
\textit{Hint: write $\omega \wedge \overline{\omega}$ in a translation chart $z = x + i y$ in terms of $dx$ and $dy$.}
\end{enumerate}
\end{tcbexercise}

\textbf{Notation}.
Unless stated otherwise topological surfaces will be denoted by $S$, Riemann surfaces by $X$ and translation surfaces by $M$.

\subsection{Relations between the three definitions} \label{sec:RelationsBetweenDefinitions}

Before diving into the geometry and dynamics of translation surfaces we prove in this section that the three definitions given above are equivalent.

Let $M$ be a constructive translation surface built from a family $(\cP, f: E(\cP) \to E(\cP))$.  Given an edge $e \in E(\cP)$ belonging to the  polygon $P_e$ consider all possible triangles contained in $P_e$ that have $e$ as one of their edges and such that the third vertex is a vertex of $P_e$. This set of triangles is not empty and we denote it $\cT_e$.  Now, let $(e,e')$ be a pair of identified edges and $\tau$ the translation such that $\tau(e') = e$. Let $t$ and $t'$ be respectively triangles in $\cT_e$ and $\cT_{e'}$. Then $t \cup \tau(t') \subset \C$ is a planar quadrilateral and we denote by $V_{t,t'}$ its interior. Also let $U_{t,t'}$ be the interior of the image of $t \sqcup t'$ in $M$. We have a well defined bijection $\phi_{t,t'}: U_{t,t'} \to V_{t,t'}$.

The following result gathers the natural links between our three definitions.
\begin{tcbtheorem}{}{EquivalencesInDefinitions}
\begin{enumerate}
\item \label{item:ConstructiveGeometric} (\emph{constructive are geometric})
Let $M$ be a constructive translation surface. Denote by $V$ and
$\Sigma$ the image of the vertices in $M$ and the conical singularities respectively (we have
$\Sigma \subset V$). Then the maps $\{\phi_{t,t'}\}_{(t,t')}$ define a
translation atlas on $M \setminus V$ that can be uniquely extended to a translation
atlas on $M \setminus \Sigma$.

\item \label{item:GeometricAnalytic} (\emph{geometric are analytic})
Let $(M,\cT)$ be a geometric translation surface with conical singularities $\Sigma$. Then there is a unique holomorphic structure on $M$ such that translation charts are holomorphic.
Moreover, there is a unique Abelian differential $\omega$ such that for each translation chart $\phi: U \to V$ we have $\phi^* dz = \omega$.

\item \label{item:AnalyticGeometric} (\emph{analytic are geometric})
Let $(X,\omega)$ be an analytic translation surface and $Z(\omega)$ the set of zeros of $\omega$. Then on $X \setminus Z(\omega)$ the local coordinates $z$ with $\omega=dz$
define a translation atlas on $X \backslash Z(\omega)$.

\item \label{item:GeometricConstructive} (\emph{geometric are constructive})
Every geometric translation surface $(M,\mathcal{T})$ can be obtained from a family of polygons and a pairing $(\cP, f: E(\cP) \to E(\cP))$.
\end{enumerate}
\end{tcbtheorem}

\begin{proof}
\noindent \textbf{\ref{item:ConstructiveGeometric}. Constructive are geometric.} The domains $U_{t,t'}$ of the maps $\phi_{t,t'}$ form by definition an open cover of $M \setminus V$. It is also clear that the transition maps between the $\phi_{t,t'}$ are translations. Now, it follows from the questions~\ref{exo:item:ProlongationConicalMetric} and~\ref{exo:item:FiniteDegreeAndConicalSingularity} of Exercise~\ref{exo:ConicalSingularities} that each point in $V \setminus \Sigma$ (\ie vertex of conical angle $2\pi$) has a neighborhood isometric to a flat disk. This isometry is the unique way of extending the atlas.

\smallskip
\noindent \textbf{\ref{item:GeometricAnalytic}. Geometric are analytic.} Let $(M,\mathcal{T})$ be a geometric translation surface. Since translations are holomorphic maps, the translation structure in $M\setminus\Sigma$ is also a Riemann surface structure. By holomorphic continuation, the Riemann structure is unique at each point of $\Sigma$ (recall that $\Sigma$ is discrete by assumption). Now, in each chart we can pull back the Abelian differential $dz$ from the plane. Because $dz$ is invariant under translation, these pull-backs agree on the intersection of the domains of two flat charts and merge into a globally well-defined holomorphic 1-form $\omega$ on $M\setminus\Sigma$. Note that this 1-form is nowhere zero. Exercise~\ref{exo:ConicalSingularityAndZeroOfDifferential} explains why a conical singularity of angle $2 k_v\pi $ corresponds to a zero of degree $v_k-1$ of the differential $\omega$.

\smallskip
\noindent \textbf{\ref{item:AnalyticGeometric}. Analytic are geometric.} Given an analytic translation surface $(X,\omega)$, the integration of $\omega$ on $X\setminus Z(\omega)$ endows $X\setminus Z(\omega)$ with an atlas where transition functions are translations. Indeed, define the chart $\varphi:U\subset X\to\C$ by $\varphi(p)=\int_{p_0}^p\omega$, where $p_0\in U$ and $U$ is a simply connected neighbourhood. Now consider another chart $\psi:V\subset X\to\C$ defined by $\int_{q_0}^p\omega$ ($V$ also simply connected), where $p\in U\cap V$ and the constant $C=\int_{p_0}^{q_0}\omega$. The equation $\int_{p_0}^{q_0}\omega=\int_{p_0}^{p}\omega+\int_{p}^{q_0}\omega$ implies that the defined charts are related by a translation: $\varphi(p)=\psi(p)+C$. The relation between the zeros of $\omega$ and the conical singularities can again be deduced from Exercise~\ref{exo:ConicalSingularityAndZeroOfDifferential}.

\smallskip
\noindent \textbf{\ref{item:GeometricConstructive}. Geometric are constructive.} We now turn to the delicate point of the proof. In what follows we show that every geometric translation surface admits a triangulation $\mathfrak{T}$ whose 1--skeleton in flat coordinates is made of Euclidean segments and whose 0--skeleton contains the set of conical singularities $\Sigma$. The existence of such a triangulation implies that the geometric translation surface can be constructed using Euclidean polygons and hence it is a constructive surface.

Let $(M,\mathcal{T})$ be a geometric translation surface with conical singularities $\Sigma$. We decompose the proof in two steps.
\begin{enumerate}
\item First we show that $M \setminus \Sigma$ admits a triangulation $\mathfrak{T}'$ whose 1-skeleton is formed (in flat coordinates) by Euclidean segments.
\item Next, we explain how the triangulation $\mathfrak{T}'$ can be refined in order to obtain a (flat) triangulation $\mathfrak{T}$ of $M$ such that $\Sigma$ is contained in the vertices of $\mathfrak{T}$.
\end{enumerate}

\emph{Step 1}.
We already know that $M$ inherits a Riemann structure from item \ref{item:AnalyticGeometric}. Hence there exists a triangulation $\mathfrak{T}''$ of $M$ whose edges are smooth and vertices have finite degree (see~\cite{AhlforsSario-book} for a proof that this triangulation always exists). We can assume that the vertices of $\mathfrak{T}''$ are disjoint from $\Sigma$ using an homotopy. Now we modify $\mathfrak{T}''$ to obtain the triangulation $\mathfrak{T}'$ as follows:
\begin{itemize}
\item\textbf{Part A.} Let $Sk_n(\mathfrak{T}'')$ denote the $n=0,1,2$ skeletons of $\mathfrak{T}''$ (\ie, vertices, edges and faces). Every vertex  $p\in Sk_0(\mathfrak{T}'')$ is in the domain of a translation chart $\phi_p:U_p\to\C$ whose image is a disk and such that $U_p\cap Sk_0(\mathfrak{T}'')=p$. Since $Sk_0(\mathfrak{T}'')$ is a discrete subset of $M$, we can further assume that $U_p \cap U_q=\emptyset$ for every $p\neq q$ in $Sk_0(\mathfrak{T}'')$. We can choose $U_p$ so that $\phi_p(U_p\cap Sk_1(\mathfrak{T}''))$  consists of finitely many smooth arcs $\gamma_j$ from $\phi_p(p)$ to the boundary. We can thus replace each $\gamma_j$ by a radius from $\phi_p(p)$ to the boundary of $\phi_p(U_p)$. See Figure~\ref{fig:FlateningTriangulation}.
\begin{figure}[!ht]
\begin{center}
\includegraphics[scale=0.9]{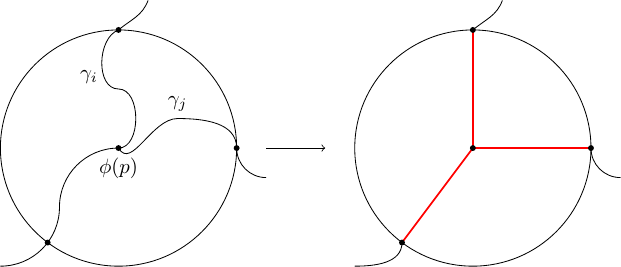}
\end{center}
\caption{Modifying the triagulation $\mathfrak{T}''$ near a vertex.}
\label{fig:FlateningTriangulation}
\end{figure}

\item\textbf{Part B.} Let us denote by $U=\cup_{p\in Sk_0(\mathfrak{T}')} U_p$, where $U_p$ is the disk chosen in the previous step. For every edge $e\in Sk_1(\mathfrak{T}')$ the subset $e\setminus U$ can be covered by finitely many open disks $D_1^e,\ldots,D_{m_e}^e$, each of which is the domain of a flat chart and which together form a chain, that is:
\begin{enumerate}
\item $D_j^e\cap U\neq \emptyset$ if and only if $j=1$ or $j=m_e$.
\item For every $j=2,\ldots,m_e-1$, $D_j^e\cap D_l^e\neq \emptyset$ if and only if $l=j \pm 1$,
\item There are no three disks in $\{D_1^e,\ldots,D_{m_e}^e\}$ sharing a common point,
\item For every $j=1,\ldots,m_e$, the intersection $e \cap D_j^e$ is non-empty smooth-connected arc.
\end{enumerate}
The intersection $e\cap (\cup_{j=1}^{m_e}\partial D_j^e)$ defines a finite sequence of points $p_1,\ldots, p_{2m_e}$. Since the collection of disks is finite and disks are convex, we can replace as in the previous step the arc between $p_i$ and $p_{i+1}$ by a straight line segment. See Figure~\ref{fig:FlateningTriangulationEdges}.
\begin{figure}[!ht]
\begin{center}
\includegraphics[scale=0.8]{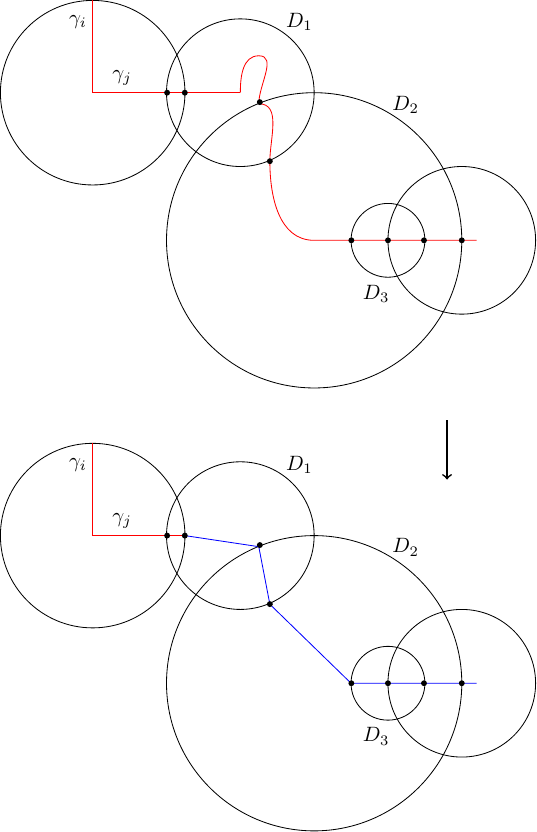}
\end{center}
\caption{Modifying the triagulation $\mathfrak{T}''$ near an edge.}
\label{fig:FlateningTriangulationEdges}
\end{figure}

\end{itemize}
\emph{Step 2}. The second step is a consequence of the following lemma.
\begin{tcblemma}{}{triangulation}
Let $M$ be a simply connected compact translation surface with a boundary that consists of finitely many Euclidean segments. Then $M$ admits a triangulation whose edges are Euclidean segments and whose vertices are exactly the vertices on the boundary and the singularities.
\end{tcblemma}
The proof we propose, which is elementary, uses concepts that we have not introduced yet. Therefore we adress it in Section~\ref{SSEC:Debt}, p.~\pageref{SSEC:Debt}.
\end{proof}

\subsection{Half-translation surfaces}
\label{ssec:HalfTranslationSurfaces}
\begin{wrapfigure}{R}{0.5\textwidth}
\centering
\includegraphics{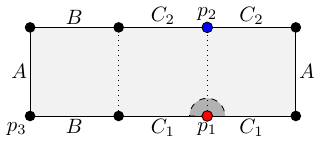}\\ \vspace{1cm}
\includegraphics{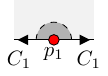}
\caption{A genus 1 half-translation surface and a neighbourhood of a conical singularity of angle $\pi$ (corresponding to a simple pole of the quadratic differential).}
\label{fig:HalfTransSurf}
\end{wrapfigure}
Translation surfaces are part of a bigger class of surfaces called half-translation surfaces or quadratic differentials\footnote{Some authors use the misleading term \emph{flat surface} for a quadratic differential. However a flat Riemannian metric on a surface does not necessarily
come from a quadratic differential.}. Before giving any formal definition let us consider the example of the surface $M$ depicted in Figure~\ref{fig:HalfTransSurf}. This surface is obtained by gluing three squares along their edges using isometries \emph{that respect the orientations shown in the figure}. In particular, the isometry used to glue the two edges labeled $C_1$ is a \emphdef{point reflection} (\ie a map of the form $z \mapsto -z+c$). Following tradition, we call these kind of isometries \emphdef{half-translations}. Note that the side $A$ is glued by translation. A direct calculation of the Euler characteristic shows that $M$ is a torus. The vertices of the three squares merge into three points $p_1$, $p_2$ and $p_3$ in the surface $M$. The points $p_1$ and $p_2$ are conical singularities of angle $\pi$, in the sense of Exercise~\ref{exo:ConicalSingularities}. Indeed, these points have neighbourhoods that can be constructed from a half-disk by gluing opposite radii using a half-translation. On the contrary, $p$ is a conical singularity of angle $4\pi$.

\bigskip

Now let us introduce a more formal definition. Consider a family of polygons
$\cP$ and a pairing
$f:E(\cP)\to E(\cP)$ such that for every $e\in E(\cP)$
the edges $e$ and $e'=f(e)$ differ by a translation. The first difference with
the definition of a translation surface is that now there are two possible
cases to consider: (i) $n_{e'}=-n_e$ and (ii) $n_{e'}=n_e$. We identify the
points in $e$ with the points in $f(e)$ in case (i) using a translation and in case
(ii) using a point reflection. Let $\pi: \bigsqcup_{P \in\cP} P \to
(\bigsqcup_{P \in\cP} P) / \sim$ be the corresponding quotient map.

\begin{tcbdefinition}{Constructive}{HalfTranslationSurfaceConstructive}
Let $\cP$ be an at most countable set of Euclidean polygons and $f: E(\cP) \to E(\cP)$ a pairing as above.  Let $M$ be $ \bigsqcup_{P\in\cP} P/\sim$ with all vertices of infinite degree and all vertices that define a conical singularity of angle $\pi$ removed. If $M$ is connected we call it the \emphdef[half-translation surface]{(constructive) half-translation surface} generated by the family of polygons $\cP$.
\end{tcbdefinition}


\begin{tcbdefinition}{Geometric}{FlatSurfaceGeometric}
A \emphdef[translation surface]{(geometric) half-translation surface} is a pair $(S,\mathcal{T})$ made of a connected topological surface $S$ and a maximal half-translation atlas $\cT$ on $S\setminus\Sigma$, where:
\begin{enumerate}
\item $\Sigma$ is a discrete subset of $S$ and
\item every $z\in\Sigma$ is a conical singularity whose angle is larger than or equal to $2\pi$.
\end{enumerate}
\end{tcbdefinition}
A half-translation atlas is just as a translation atlas, except that we allow change of coordinates to be half-translations. Given that the Euclidean metric in $\C$ is also invariant by half-translations, every geometric half-translation surface inherits a natural Riemannian metric from the plane by pull-back via its flat coordinates.


Finally, we define flat surfaces analytically.
\begin{tcbdefinition}{Analytic}{FlatSurfacesAnalytic}
An \emphdef[half-translation surface]{(analytic) half-translation surface} is a pair $(X,q)$ formed by a (connected) Riemann surface $X$ and a holomorphic non-identically zero quadratic differential $q$.
\end{tcbdefinition}

\begin{tcbremark}{}{NoPoles}
In this text we sometimes consider  pairs $(X,q)$ formed by a (connected) Riemann surface $X$ and a \emph{meromorphic} not identically zero quadratic differential $q$ whose poles are all simple, if any. Every such pair defines the half-translation surface $(X\setminus Poles(q),q)$.
The convention made in this text is \emph{not to consider simple poles (or equivalently conical singularities of angle $\pi$) as part of a half-translation surface}. The reasons behind this convention are rather technical, but will be clarified to the reader when the necessary material has been introduced\footnote{Two examples of these technical reasons are: core curves of cylinders in a sphere with a meromorphic quadratic differential where poles have not been removed are always null-homotopic; for any finite-type hyperbolic Riemann surface, the cotangent space of Teichm\"{u}ller space (at a point) is identified with the space of integrable \emph{holomorphic} quadratic differentials, see Proposition 6.6.2 in \cite{Hubbard-book1}.}.
\end{tcbremark}


Let us now explain how quadratic differentials define geometric half-translation surfaces. Let $\{z_i:U_i \to \C\}_{i\in I}$ be an atlas for the Riemann surface structure of $X$. A quadratic differential $q$ on $X$ is a collection of holomorphic functions $f_i$ defined on $z_i(U_i)$ such that
\begin{equation}
\label{eq:CocycleQuadraticDiff}
\forall p \in U_i\cap U_j, \qquad
f_i(z_i(p))\left(\frac{dz_i}{dz_j} (p)\right)^2=f_j(z_j(p)).
\end{equation}
We denote by $Z(q)$ the set of zeroes of $q$ in $X$. The key point is to look for the local coordinate $\xi$ on $X\setminus Z(q)$ such that $q=d\xi^2$. Indeed, for every $p\in X\setminus Z(q)$ there exists a neighbourhood $U_i=U_i(p)$ on which a square root of $f_i$ exists. If we define for every $z\in U_i$
\begin{equation}
\label{eq:CoordinateFlatStructure}
\xi_i(z)=\int_p^z\sqrt{f_i(w)}dw.
\end{equation}
this coordinate satisfies $q=d\xi_i^2$ on $U_i$. Moreover, since we had to make a choice for the square root of $f$, the coordinate $\xi_i(z)$ is unique up to translation and change of sign. In other words, the coordinates $\xi_i(z)$ define an atlas on $X\setminus Z(q)$ where transition functions are of the form $z\to\pm z + a$. That is,  $\{(U_i,\xi_i)\}_{i\in I}$ is a half-translation  atlas.

On the analytic side, an Abelian differential $\omega$ is just turned into a quadratic differential by considering $q = \omega^2$. However, not every quadratic differential is globally the square of an Abelian differential (\eg $dz^2 / (z(z+1)(z-1))$ on $\C\setminus \{0,1,-1\}$).
\begin{tcbexercise}{}{SingFlatSurface}
Show that the three definitions of half-translation surface that we have given above are equivalent. \emph{Hint:} the proof is \emph{mutatis mutandis}, the same as the proof of Theorem \ref{thm:EquivalencesInDefinitions}.
\end{tcbexercise}

\begin{tcbexercise}{}{ConicalSingularitiesHalfTranslationSurfaces}
This exercise discusses conical singularities of half-translation surfaces.
\begin{enumerate}
\item Let $M$ be a geometrical half-translation surface. Prove that the angle of a conical singularity $z\in\Sigma$ in the definition above is of the form $\pi k$ for some positive integer $k\geq 2$.
\item Let $(X,q)$ be an analytic half-translation surface and $z_0\in Z(q)$ a zero of multiplicity $k\geq 1$. Show that $z_0$ in flat coordinares is a conical singularity of total angle $(k+2)\pi$.
\item Let $M$ be a constructive half-translation surface and $v$ a vertex of a polygon in $\mathcal{P}$. Prove that there exist a positive integer $k_v\geq 2$ for which the following analog of equation (\ref{eq:AdmissibleVertices}) holds:
\begin{equation}
\label{eq:AdmissibleVerticesDF}
\sum_{w\in\pi^{-1}(\pi(v))} \alpha_w =k_v \pi.
\end{equation}
\end{enumerate}
\end{tcbexercise}

\textbf{Translation double covering}. Let $M=(X,q)$ be a half-translation surface and suppose that $q$ is not the square of an Abelian differential. In the next paragraphs we explain how to construct a (ramified) double covering $\pi:M^2\to M$ such that $\pi^*q=\omega^2$, where $\omega$ is an Abelian differential on $M^2$.

As before, let $\{z_i: U_i \to \C\}_{i \in I}$ be holomorphic charts for $X$ and $q$ given locally by $\{f_i(z_i) dz_i^2\}$ (which satifies~\eqref{eq:CocycleQuadraticDiff}). Let us consider only charts that avoid the zeros $Z(q)$ of $q$, that is the charts $z_i: U_i \to \C$ such that $f_i(z_i(p)) \not=0$ for all $p \in U_i$. In each such chart, we can consider $g_i^+(z_i)$ and $g_i^-(z_i)$ the two square root of $f_i(z_i)$ on $U_i$, that is $(g_i^+(z_i))^2 = (g_i^-(z_i))^2 = f_i(z_i)$. We build the translation double covering locally. We first consider two copies $U_i^+$ and $U_i^-$ of $U_i$. Then for each pair $i,j$ such that $U_i\cap U_j\neq\emptyset$ and let $\pi_{i,\pm}:U_i^\pm\to U_i$ be the natural maps. We identify $\pi_{i,+}^{-1}(U_i\cap U_j)$ with $\pi_{j,+}^{-1}(U_i\cap U_j)$  and $\pi_{i,-}^{-1}(U_i\cap U_j)$ with $\pi_{j,-}^{-1}(U_i\cap U_j)$ using the maps $\pi_{j,+}^{-1}\circ\pi_{i,+}$ and $\pi_{j,-}^{-1}\circ\pi_{i,-}$ respectively if
\begin{equation}
\label{eq:CocycleAbelianDiferentialLocal}
g_i^+(z_i(p))\cdot\frac{dz_i}{dz_j}(p)=g_j^+(z_j(p)).
\end{equation}
On the other hand, if $g_i^+(z_i(p))\cdot\frac{dz_i}{dz_j}(p)=g_j^-(z_j(p))$ we identify $\pi_{i,+}^{-1}(U_i\cap U_j)$ with $\pi_{j,-}^{-1}(U_i\cap U_j)$  and $\pi_{i,-}^{-1}(U_i\cap U_j)$ with $\pi_{j,+}^{-1}(U_i\cap U_j)$  using the maps $\pi_{j,-}^{-1}\circ\pi_{i,+}$ and $\pi_{j,+}^{-1}\circ\pi_{i,-}$ respectively. These gluings define an open Riemann surface which is a (non-ramified) double covering of $X\setminus Z(q)$. The collection of expressions $\{g_i^\pm\circ z_i \circ \pi_{i,\pm} d(z_i\circ\pi_{i,\pm})\}_{i\in I}$ defines a holomorphic 1--form
$\omega$ on this open Riemann surface that can be holomorphically extended to a (possibly ramified) covering $\pi:M^2\to M$ called the \emphdef{translation double covering of $M$}. Moreover, by construction $\pi^*q=\omega^2$ and the branching locus\footnote{For a precise definition of branching locus and ramification points see Section~\ref{sec:CoveringSpaces}.} of this covering (which might be empty) is contained in $Z(q)$.

\begin{tcbexercise}{}{RamificationPoints}
Let $M$ be half-translation surface, $\pi:M^2\to M$ its orientation double covering and suppose that $z_0\in M$ is a conical singularity of angle $n\pi$, for some $n\in\N$. Prove that if $n$ is even, then $\pi^{-1}(z_0)$ is formed by two conical singularities of total angle $n\pi$; whereas if $n$ is odd, $\pi^{-1}(z_0)$ is just a conical singularity of total angle $2n\pi$. Deduce that the branching points of the orientation double cover are the conical singularities in $M$ of angle $n\pi$, $n$ odd.
\end{tcbexercise}

We finish this section by presenting a way to define the double covering $M^2$ that only uses the half-translation atlas and elementary algebraic topology. Most of the details are left to the reader on exercise \ref{exo:double cover}. The idea is to use the \emphdefrecall{holonomy} of the half-translation atlas (for more details on this concept see Appendix~\ref{Appendix:GXStructures}) and then take linear parts of the corresponding maps: let $\gamma: [0,1] \to M$ be a closed loop. Then there exists a finite set of charts $\xi_0: U_0 \to \C$, $\xi_1: U_1 \to \C$, \ldots, $\xi_n: U_n \to \C$ and real numbers
 $t_0 = 0 < t_1 < \ldots < t_n = 1$ so that $\{\gamma(t)\}_{t \in [t_{i}, t_{i+1}]} \subset U_i$. On each intersection $U_i \cap U_{i+1}$ the linear
 part of the transition map $\xi_{i+1} \circ \xi_i^{-1} |_{U_i \cap
 U_{i+1}}$ is either $+1$ or $-1$. Let $s_i$ be this sign. We set $f(\gamma) = s_0 s_1 \ldots s_n$.



\begin{figure}[!ht]
\begin{center}
\includegraphics{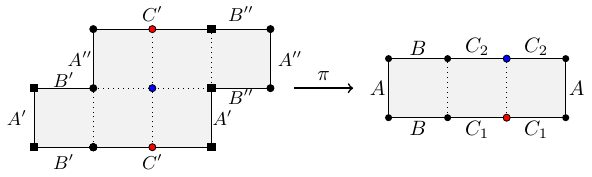}
\end{center}
\caption{The orientation double covering of the half-translation surface of Figure~\ref{fig:HalfTransSurf}.}
\label{fig:DoubleCovering}
\end{figure}

\begin{tcbexercise}{}{double cover}
Show that the map $f$ defined above defines a map\footnote{Formally speaking, this map is the holonomy of the half-translation atlas on $M$ composed with the derivative map, see Appendix \ref{Appendix:GXStructures} for more details.} $\hat{f}: \pi_1(M\setminus \Sigma,x_0) \to \{+1,-1\}$.
\begin{enumerate}
\item Show that the half-translation atlas of $M$ contains a translation atlas if and only if the image of $f$ is $\{+1\}$.
\item If the half-translation atlas of $M$ contains a translation atlas show that it contains exactly two.
\item Show that the double (regular) covering $\mathring{M}^2\to M\setminus \Sigma$ given by the subgroup $\ker{\hat{f}}$ can be endowed with a translation structure. Show that this double regular covering can be extended to a possibly ramified covering $\pi:M^2\to M$ and that the branching locus of $\pi$ is contained in $\Sigma$. In other words, show that $M^2$ is precisely the orientation double covering of $M$.
\item If $M$ comes from a holomorphic quadratic differential $(X,q)$ show that the surface defined by the equation
$\{\omega \in T^*(X):\ \omega^2 = q\}$, where $T^*(X)$ is the (complex) cotangent space of $X$, carries a canonical translation structure.

\end{enumerate}
\end{tcbexercise}
In Figure~\ref{fig:DoubleCovering} we depict an example of an orientation double covering.

%


\subsection{Structures}
\label{ssec:Structures}
In this section we discuss the most important structures and invariants that one can associate with a translation (or half-translation) surface. The material we present will be illustrated with several examples in the following sections of this chapter.
We stress that this material is standard in the context of \emph{compact} translation surfaces, as discussed in several different surveys on the subject such as H. Masur, S. Tabachnikov~\cite{MasurTabachnikov02}, A. Zorich~\cite{Zorich06}, J.-C. Yoccoz~\cite{Yoccoz10}, G. Forni, C. Math\'eus~\cite{ForniMatheus14} or A. Wright~\cite{Wright15}.
\\~\\

\noindent\textbf{Metric completion and singularities}. Let $M$ be a translation or half-translation surface and $\Sigma\subset M$ its conical singularities. As discussed in the preceding section, $M\setminus \Sigma$ can be naturally endowed with a flat Riemannian metric $\mu$ which provides a canonical notion of distance and area. We denote by $\widehat{M}$ the corresponding metric completion of $M\setminus\Sigma$.

\begin{tcbexercise}{}{}
Let $M$ be a translation or half-translation surface. Show that:
\begin{enumerate}
\item $M$ is canonically embedded in $\widehat{M}$ and the image of this embedding is dense in  $\widehat{M}$.
\item Suppose that $M=(X,q)$ is a half-translation surface and that $q$ can be meromorphically extended to a point $z_0\in\widehat{M}$. Show that if $q$ is not holomorphic at $z_0$, then this point is a simple pole and that the Euclidean metric of $M$ can be extended to $z_0$. Show that in the case where $z_0$ is a simple pole this point is a conical singularity of angle $\pi$.
 \item If $M$ is constructed from a family of polygons $\mathcal{P}$, show that there is an inclusion $\left(\bigsqcup_{P\in\cP} P/\sim\right)\hookrightarrow \widehat{M}$. In other words: show that all vertices need to be added to obtain the metric completion of a constructive translation surface. Show also that the total area of $M$ w.r.t the Riemannian metric $\mu$ is the sum of the areas of the polygons in $\mathcal{P}$.
\end{enumerate}
\end{tcbexercise}

Recall that in Exercise \ref{exo:ConicalSingularities} we introduced the notion of conical singularity and infinite angle singularity.


\begin{tcbdefinition}{}{Singularities}
Let $M$ be a translation or half-translation surface and $\widehat{M}$ its metric completion with respect to the natural flat metric. A point $p\in \widehat{M}\setminus M$ is called:
\begin{enumerate}
\item \emphdef[regular point]{regular} if it has a neighbourhood isometric to an open disk in the plane. A point in $\widehat{M}$ that is not regular is called a \emphdef{singularity}.
\item a \emphdef{wild singularity} if $p$ is not a regular point, nor a conical singularity, nor an infinite angle singularity.
\end{enumerate}
The set of singularities is denoted by $\Sing(M)$. A translation surface is called \emphdef[tame translation surface]{tame} if $\Sing(M)$ is formed exclusively by conical and infinite angle singularities, and \emphdef[wild translation surface]{wild} in all other cases.
\end{tcbdefinition}


Wild and infinite angle singularities only appear in the metric completion of the surface. The infinite staircase presents 4 infinite angle singularities and no wild singularity, see Figure~\ref{fig:StaircaseFirst}. Baker's surfaces, which we define in Section~\ref{ssec:BakerBowmanArnouxYoccoz}), are finite area translation surfaces presenting only wild singularities. We study singularities of translation surfaces in detail in chapter \ref{sec:Singularities}.\\


\begin{tcbremark}{}{}
The set of wild singularities can be disjoint from the set of vertices on a constructive translation surface. Indeed, it is not difficult to construct an infinite triangulation $\mathfrak{T}$ of the (open) unit disk $\mathbb{D}=\{z\in\C\hspace{1mm}|\hspace{1mm}|z|<1\}$. Every vertex of $\mathfrak{T}$ is properly contained in $\mathbb{D}$, however $\Sing(\mathbb{D})=\{z\in\C\hspace{1mm}|\hspace{1mm}|z|=1\}$.
\end{tcbremark}

\noindent \textbf{Infinite type.}
A translation surface is compact if and only if it can be built from finitely many polygons. We introduce a slightly more general notion (finite type) which allows us to define the main objects of study of this book.

\begin{tcbdefinition}{}{FiniteTypeInfiniteType}
Let $M=(X,\omega)$ be a translation surface. We say that $M$ is of \emphdef[finite-type (translation surface)]{finite-type} if $M$ has finite area and $X$ is a finite-type\footnote{Recall that a Riemann surface $X$ has \emphdef[finite-type (Riemann surface)]{finite-type} if it is isomorphic (as Riemann surface) to a compact Riemann surface from which we have removed a finite set of points. For example, the plane $\C$ is of finite-type (as it is the sphere $\widehat{\C}$ minus one point) whereas the unit disk $\D$ is not.} Riemann surface. If $M$ is not of finite-type we say it is of \emphdef[infinite-type (translation surface)]{infinite-type}.
\end{tcbdefinition}


We stress that the fundamental group of $X$ does not play a role in Definition~\ref{def:FiniteTypeInfiniteType}. However, a translation surface $M$ whose fundamental group is not finitely generated is always of infinite type.

To define half-translation surfaces of finite type just change $\omega$ to $q$ (a quadratic differential) in the preceding definition. In Section~\ref{ssec:PanovPlanes} we show that there are infinite-type half-translation surfaces of the form $(\C,q)$ with interesting dynamical properties.

\smallskip

\noindent \textbf{Affine maps.} Let $M_1$ and $M_2$ be translation surfaces with conical singularities $\Sigma_1$ and $\Sigma_2$ respectively. A map
$f:M_1\to M_2$ with $f(\Sigma_1)\subset\Sigma_2$ is called an \emphdef{affine map} if the
restriction $f:M_1\setminus\Sigma_1\to M_2\setminus\Sigma_2$ in flat charts is an
$\R$-affine map.
That is, if $z_1$ (respectively $z_2$) denotes a flat local coordinate of $M_1\setminus\Sigma_1$ (resp. $M_2\setminus\Sigma_2$) and we write $z_1 = x_1 + iy_1$ (resp. $z_2 = x_2 + iy_2$) then $f$ in these coordinates is given by an expression of the form:
\begin{equation}
\label{eq:AffineHomeo}
\begin{pmatrix}x_1\\y_1\end{pmatrix} \mapsto
\begin{pmatrix}x_2\\y_2\end{pmatrix} =
\begin{pmatrix}a & b \\ c & d \end{pmatrix}
\begin{pmatrix}x_1 \\ y_1 \end{pmatrix} + \begin{pmatrix}\xi_0 \\ \eta_0 \end{pmatrix}.
\end{equation}
Since two coordinates differ only by a translation, the matrix $\displaystyle \begin{pmatrix}a & b \\ c & d \end{pmatrix}$ does not depend on the coordinates $z_1$ and $z_2$; however the translation vector $\displaystyle \begin{pmatrix}\xi_0\\\eta_0\end{pmatrix}$ does depend.

Let $M$ be a translation surface and let us denote by $M^0:=M\setminus\Sing(M)$. Because $M^0$ has an atlas where transition functions are translations, the tangent bundle of $M^0$ is trivial, that is, $TM^0=M^0\times\C$. For an affine map $f:M_1\to M_2$ we denote by $Tf:TM_1^0\to TM_2^0$ the tangent map of $f$. If $f$ is given in local coordinates by~\eqref{eq:AffineHomeo}, then the tangent map can be written as $Tf(z_1,v)=(f(z_1),Df\cdot v)$, where $Df$ is the (constant) derivative of $f$ given by $\begin{psmallmatrix}a & b \\ c & d \end{psmallmatrix}$.
\begin{tcbdefinition}{}{TypologyAffineMaps}
Let $M_1$ and $M_2$ be translation surfaces. An affine map $f: M_1 \to M_2$ is called a \emphdef{translation} if $Df = 1$ and an \emphdef{isometry} if $Df \in \SO(2,\R)$.

The translation surfaces $M_1$ and $M_2$ are \emphdef{isomorphic} if there exists a 1-to-1 translation $f: M_1 \to M_2$.
\end{tcbdefinition}

Note that translations and more generally isometries are special cases of volume preserving affine maps (\ie affine maps $f$ with $\det(Df) = \pm 1$).

An \emphdef{affine automorphism} is a homeomorphism $f:M\to M$ which is also an affine map. We denote by $\Aff(M)$ the group of affine homeomorphisms of $M$ and by $\Aff^+(M)$ the subgroup of $\Aff(M)$ made of orientation preserving affine automorphisms (\ie their linear part has positive determinant). If $\Sigma\subset M$ we denote by $\Aff(M,\Sigma)$ the subgroup formed by all elements $f\in\Aff(M)$ such that $f(\Sigma)=\Sigma$, \ie that leave $\Sigma$ invariant (not necessarily fixing every point). This notation extends to all subgroups of $\Aff(M)$.  We denote by $\Trans(M)$ (resp. $\Isom(M)$) the groups formed by all affine automorphisms of $M$ which are translations (resp. isometries).

\begin{tcbdefinition}{}{VeechGroup}
The derivative $Df$ of an affine map lives in $\GL(2,\R)$ and it provides a homomorphism $D: \Aff(M) \to \GL(2,\R)$. We define $\Gamma(M):=D(\Aff(M))$, its (at most index two) subgroup $\Gamma^+(M):=D(\Aff^+(M))$ and $\Gamma(M,\Sigma):=D(\Aff(M,\Sigma))$. We call $\Gamma^+(M)$ and $\Gamma(M)$ respectively the \emphdef{Veech group} and the \emphdef{extended Veech group} of $M$.
\end{tcbdefinition}

\begin{tcbdefinition}{}{TypologyAffineAutomorphisms}
An orientation and volume preserving affine automorphism $f: M \to M$ (\ie so that $\det(Df) = 1$) is called \emphdef[elliptic (affine map)]{elliptic} if $|\tr(Df)| < 2$, \emphdef[parabolic (affine map)]{parabolic} if $|\tr(Df)| = 2$ and \emphdef[hyperbolic (affine map)]{hyperbolic} if $|\tr(Df)| > 2$. In the case $M$ is a compact translation or half-translation surface, an affine hyperbolic automorphism is called an \emphdef{affine pseudo-Anosov homeomorphism}.
\end{tcbdefinition}
This nomenclature mimics the one used for matrices in $\SL(2,\R)$. Notice that an automorphism is elliptic if and only if it is an isometry. Note that there are several equivalent definitions of what a pseudo-Anosov homeomorphism is (see \eg \cite{Thurston88}, \cite{FathiLaudenbachPoenaru} or~\cite{FarbMargalit}).

We study affine and Veech groups of infinite-type translation surfaces in detail in Chapter~\ref{ch:Symmetries}.

\smallskip

\noindent \textbf{$\GL(2,\R)$-action}. There is a natural action of the group $\GL(2,\R)$ on the set of all translation surfaces, which is easily defined as follows: if $(S,\mathcal{T})$ is a geometric translation surface, $\mathcal{T}=\{\phi_i:U_i\to\C\}_{i\in I}$ and $A\in\GL(2,\R)$, then $A(S,\mathcal{T})=(S,A\mathcal{T})$, where $A\mathcal{T}=\{A\circ\phi_i:U_i\to\C\}_{i\in I}$. Remark that $A\mathcal{T}$ is indeed a translation atlas for:
$$
(A\circ\phi_j)\circ(A\circ\phi_i)^{-1} = A\circ(\phi_j\circ\phi_i^{-1})(A^{-1}(z)) = A(A^{-1}(z)+c) = z + A(c)
$$

In the context of finite-type translation surfaces, this action is a fundamental ingredient for the study of the dynamics of the translation flow, which we define in the next paragraphs. For a detailed discussion of this action in the context of finite-type surfaces see Section 3 in \cite{Wright15}.

\smallskip

\noindent \textbf{Translation flow.}
For each direction $\theta\in\R/2\pi\Z$ we have a well-defined translation flow $F_{\C,\theta}^t: \C \to \C$ given by $F_{\C,\theta}^t(z) = z + t e^{i \theta}$. This flow defines a constant vector field $X_{\C,\theta}:=\frac{\partial F_{\C,\theta}^t}{\partial t}|_{t=0}(z)$.
Now let $M$ be a translation surface and recall that $\Sing(M) \subset \widehat{M}$ is the set of singularities.
In each flat chart of $M$ we can do a pull-back of $X_{\C,\theta}$ and this define a vector field $X_{M,\theta}$ on $M\setminus\Sing(M)$. For every $z\in M\setminus\Sing(M)$ let us denote by $\gamma_z:I\to M$, where $I\subset\R$ is an interval containing zero, the maximal integral curve of $X_{M,\theta}$ with initial condition $\gamma_z(0)=z$. We define $F_{M,\theta}^t(z):=\gamma_z(t)$ and call it the \emphdef{translation flow} on $M$ in direction $\theta$. Let us remark that formally speaking $F^t_{M,\theta}$ is not a flow because the trajectory of the curve $\gamma_z(t)$ may reach a singulatiry of $M$ in finite time. However for most of the surfaces that we study in this text for every direction $\theta$ we have that $F^t_{M,\theta}$ is a flow on a subset of $M$ of full measure (w.r.t. to the natural area form defined by the translation structure). When there is no need to distinguish translation flows in different translation surfaces we abbreviate $F_{M,\theta}^t$ and $X_{M,\theta}$ by $F_{\theta}^t$ and $X_{\theta}$ respectively.

 The term \emph{translation flow} (\ie without specifying the direction) is reserved for $F_{M,\pi/2}^t$, that is, the translation flow in the vertical direction and is abbreviated as $F_{\pi/2}^t$ or simply by $F^t$ when all other variables in the context are clear. The (vertical) translation flow in $\begin{psmallmatrix} \cos(\theta) & \sin(\theta)\\ -\sin(\theta) & \cos(\theta) \end{psmallmatrix}M$ and the translation flow in direction $\theta+\frac{\pi}{2}$ in $M$ are conjugate.

When $M$ has finite type, for all directions $\theta$, the set of points for which $F_{M,\theta}^t$ is defined for all times forms a subset of full Lebesgue measure. We will say that such a translation flow is \emphdef[complete (flow)]{complete}. In Chapter~\ref{chap:TranslationFlowsAndIET} we will discuss under which conditions the translation flow on a non-compact translation surface is complete. The flat disk $(\D, dz)$ is an example where the flow is not complete in any direction.

The translation flow in any direction $\theta$ also preserves the (flat) Riemannian metric $\mu$ on $M$. In particular, where it is defined, it preserves the area: if $A \subset M$ and $F_{M,\theta}^t$ is well defined on $A$ then the area of $A$ and $F_{M,\theta}^t(A)$ are identical. For this reason many tools from ergodic theory are available to study translation flows. Because of the presence of conical singularities, in general the translation flow does not act by isometries.
The properties of translation flows are studied in a (second) forthcoming volume~\cite{DHV2}
dedicated mostly to dynamical aspects of infinite-type translation surfaces.

\begin{tcbremark}{}{GeodesicFlowVSTranslationFlow}
Let $M$ be a translation surface and $\UTan(M)$ the unit tangent bundle of $M\setminus\Sing(M)$ (w.r.t. the flat metric of $M$). Given that all coordinate changes are translations this bundle is trivial. More precisely, we have a canonical isomorphism $\UTan(M) \simeq M \times \R/2\pi\Z$: the vector field $X_\theta$ corresponding to $e^{i \theta}$ considered above gives the slice $M \times \{\theta\}$ in this product.
Now, the geodesic flow $G^t$ on $\UTan(M)$ in these product coordinates is just $G^t: (z, \theta) \mapsto (F_{M,\theta}^t(z), \theta)$. In other words, the 3-dimensional unit tangent bundle $\UTan(M)$ is foliated by copies of $M$ which are invariant by the geodesic flow.  On each invariant surface $G^t$ is canonically identified to a translation flow on $M$.
\end{tcbremark}

Recall that, given a half-translation surface $M=(X,q)$, the local coordinates $z$ for which the quadratic differential satifies $q = (dz)^2$ are well defined up to a sign. For this reason the pull back of the vector field $e^{i\theta}$ to $M$ does not define a global vector field on $M$. However, this pull back does define a \emph{direction field} (see \eg~\cite{ArnoldODE} for a definition) and the corresponding set of integral curves defines a foliation $\mathcal{F}_\theta$ on $M\setminus Z(q)$ by curves which locally look like straight line segments.
Note that on the translation double covering $M^2$ of $M$, the directional vector fields and translation flows are well defined.

\smallskip

\noindent\textbf{Cylinders and strips}. There are two kinds of subsets that a translation surface might have and whose mere existence guarantees that the translation flow in certain directions, restricted to these subsets, is simple. These subsets are cylinders and (infinite) strips.

\begin{tcbdefinition}{Cylinders and strips}{CylindersAndStrips}
A \emphdef{horizontal cylinder} $C_{c,I}$ is a translation surface of the form $([0,c]\times I) / \sim$, where $I\subset\R$ is open (but not necessarily bounded), connected, and where $(0,s)$ is identified with $(c,s)$ for all $s\in I$. The numbers $c$ and $h=|I|$ are called the \emphdef[circumference (cylinder)]{circumference} and \emphdef[height (cylinder)]{height} of the cylinder respectively. A \emphdef{horizontal strip} $C_{\infty,I}$ is a translation surface of the form $\R\times I$, where $I$ is a bounded open interval. Analogously, the height of the horizontal strip is $h=|I|$.

An open subset of a translation surface $M$ is called a \emphdef{cylinder} (respectively a \emphdef{strip}) in direction $\theta$ if it is isomorphic to $e^{-i\theta}C_{c,I}$ (respect. to $e^{-i\theta}C_{\infty,I}$).
\end{tcbdefinition}

\begin{figure}[!ht]
\begin{minipage}{.45\textwidth}
\begin{center}\includegraphics[scale=0.8]{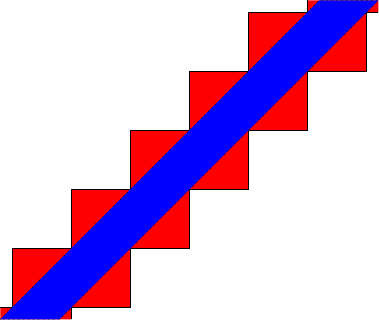}\end{center}
\subcaption{Two infinite strips in the infinite staircase from Figure~\ref{fig:StaircaseFirst}.}
\end{minipage}
\hspace{.1\textwidth}
\begin{minipage}{.45\textwidth}
\begin{center}\includegraphics[scale=0.8]{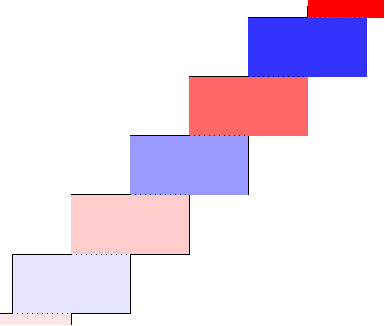}\end{center}
\subcaption{Cylinders in the horizontal direction in the infinite staircase.}
\end{minipage}
\bigskip \\
\begin{center}
\begin{minipage}{.45\textwidth}
\begin{center}\includegraphics{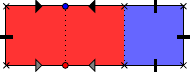}\end{center}
\subcaption{Two vertical cylinders in the genus 1 half-translation surface from Figure~\ref{fig:HalfTransSurf}.}
\end{minipage}
\end{center}
\caption{Cylinders and strips in translation and half-translation surfaces.}
\label{fig:cylindersinflatsurfaces}
\end{figure}

One can think of strips as cylinders of infinite circumference and finite height. For half-translation surfaces $M=(X,q)$ the definition of cylinder still makes sense, though its direction is well defined only up to change of sign. Examples of cylinders and strips are illustrated in Figure~\ref{fig:cylindersinflatsurfaces}.
\begin{tcbexercise}{}{PeriodicTrajectoriesAndCylinders}
Let $M$ be a translation surface.
\begin{compactenum}
\item Show that a cylinder in direction $\theta$ of $M$ consists of periodic trajectories of $F^t_{M,\theta}$ whose lengths are the circumference of the cylinder.
\item Conversely, show that every periodic trajectory of $F^t_{M,\theta}$ is contained in a maximal cylinder. \textit{hint: if $x\in M$ is a point where the translation flow $F^t_{M,\theta}$ is defined for some time in the past and in the future, then $x$ is the center of a neighbourhood in $M$ formed by points with this property.}
\end{compactenum}
\end{tcbexercise}

\noindent \textbf{Strata of finite-type translation surfaces.}
In this paragraph we review without details the definition of strata of finite-type translation surfaces. The reader is invited to consult details in~\cite[Section 3.3]{Zorich06}, ~\cite[Section 2.3]{Moeller13} (analytic version) and~\cite[Section 6]{Yoccoz10} (constructive version).

All singularities of a finite-type translation surface have finite angles. In particular one can define appropriate moduli spaces for the set of equivalence classes of translation surfaces (or quadratic differentials). If $\kappa_1, \kappa_2, \ldots, \kappa_m$ are non-negative integers then $\cH(\kappa_1, \kappa_2, \ldots, \kappa_m)$ denotes the set of tuples $(M, p_1, \ldots, p_m)$ where $M$ is a translation surface (up to isomorphism) with conical singularities of angles  $2 (1 + \kappa_1) \pi$, $2 (1 + \kappa_2) \pi$, \ldots $2 (1 + \kappa_m) \pi$ at the (distinct and enumerated) points $p_1$, $p_2$, \ldots, $p_m$ respectively. This set is called a \emphdef[stratum (of translation surfaces)]{stratum} of translation surfaces. From Riemann-Hurwitz formula, one can deduce that all surfaces in a stratum have the same genus $g$ given by the equation $\sum_{i=1}^m \kappa_i = 2g-2$. For example, the translation surface depicted in Figure~\ref{fig:LShapedOrigami} is a point in $\cH(2)$. When $\kappa_i=0$ the point $p_i$  is  a regular point: $\cH(2,0)$ is the stratum of translation surfaces with one cone angle singularity of angle $6\pi$ and one marked point (hence genus 2). Strata have a nice complex orbifold structure given by period coordinates. These coordinates permit to define the Masur--Veech measure on strata, and it is a classical result of Masur and Veech that the suborbifold formed by area 1 translation surfaces always has finite volume.

Consider now a tuple of integers $\kappa_i\geq -1$ and denote by $\cQ(\kappa_1, \kappa_2, \ldots, \kappa_m)$ the space of meromorphic quadratic differentials $(X,q)$ (up to isomorphism) for which $X$ is homeomorphic to a genus $g$ compact surface and all poles of $q$ are simple. This set is called a \emph[stratum (of quadratic differentials)]{stratum} of quadratic differentials. In this case $(X,q)$ has a conical singularity of angle $(2 + \kappa_i) \pi$ at each $p_i$ (and are flat elsewhere). The genus of these surfaces is given by $\sum_{i=1}^m \kappa_i = 4g-4$. Period coordinates and a finite measure can also be defined for $\cQ(\kappa_1, \kappa_2, \ldots, \kappa_m)$  by considering orientation double coverings.

\begin{tcbremark}{}{SpecialNotation}
Following tradition, the stratum of translation surfaces with $n$ conical singularities of the same angle, say $2(1+\kappa)\pi$, is abbreviated as $\cH(\ldots,\kappa^n,\ldots)$. The same applies to strata of quadratic differentials. For example $\cH(1^2) = \cH(1,1)$ while $\cQ(2,-1^2) = \cQ(2,-1,-1)$.
\end{tcbremark}

\begin{tcbremark}{}{}
  For infinite-type translation surfaces, there are many types of wild singularities. They are discussed in detail in Section~\ref{sec:Singularities}. One could classify  singularities according to their local geometry and in this way extend the definition of strata for infinite-type translation surfaces. We want to stress that there are major technical problems down this road, for example there is no clear candidate for a probability measure on strata. \emph{i.e.} for an analog of the Masur-Veech measure.


\end{tcbremark}

\section{Examples}
\label{sec:ExamplesInfiniteTranslationSurfaces}

In the rest of this chapter we present several mathematical contexts where infinite-type translation surfaces arise. Most of these are revisited later to illustrate general aspects about translation surfaces.


\subsection{Polygonal billiards}
\label{ssec:PolygonalBilliards}

Polygonal billiards are a recurrent example in this text. Roughly speaking, a polygonal billiard is the dynamical system defined by the frictionless motion of a point-particle inside a Euclidean polygon $P$ where all collisions with the boundary are elastic.  This means that each time that a point hits a side of the polygon the angle of incidence of its trajectory will be equal to the angle of reflection. By convention, the motion of a point ends when reaching a corner. One of the main motivations to study polygonal billiards is the following longstanding conjecture:

\begin{tcbquestion}{}{PolygonalBilliardPeriodicOrbits}
\label{CONJ:PerBill}
Does every polygonal billiard have a closed trajectory?
\end{tcbquestion}

A trajectory is closed if the point that follows it returns (after finitely many bounces) to its starting position \emph{with the same direction} in which it started its motion. According to R. E. Schwartz ``it is fair to say that this 200-year-old problem is widely regarded as impenetrable" \cite{Schwartz09}. Question~\ref{qu:PolygonalBilliardPeriodicOrbits} has a positive answer in the following situations.

It is a great achievement of the theory of compact translation surfaces that Question~\ref{qu:PolygonalBilliardPeriodicOrbits} has a positive answer in the case where all internal angles of the polygon are rational multiples of $\pi$. More precisely, in that situation, the set of tangent vectors to the polygon whose associated orbit is periodic for the billiard flow is dense in the tangent space. We refer the reader to the original article of H.~Masur~\cite{Masur86} or the survey~\cite{MasurTabachnikov02}.
Beyond rational polygons, only triangles have been studied intensively and the following statement is a summary of our knowledge.
\begin{tcbtheorem}{}{TriangularBilliardPeriodicOrbits}
\begin{enumerate}
\item For right triangles, the set of tangent vectors to the billiard table whose associated orbit is periodic is dense \cite{Troubetzkoy2005}, \cite{Hooper2007}.
\item For acute triangles, the orthic triangle is a periodic billiard trajectory of period three (Fagnano 1745).
\item Every obtuse triangle whose big angle is at most $5/9\pi$ admits a periodic trajectory \cite{Schwartz09}.
\end{enumerate}
\end{tcbtheorem}
The proof of Item 1 can be found in the second volume~\cite{DHV2} in the section about recurrence. For Item 2 and 3 we refer the reader to the article of R.~Schwartz~\cite{Schwartz09}.

Despite the evidence provided by Theorem~\ref{thm:TriangularBilliardPeriodicOrbits}, one has to be prudent: for any $\epsilon>0$ there exists a triangle whose two small angles are within $\epsilon$ radians of $\frac{\pi}{6}$ and $\frac{\pi}{3}$ respectively that has no periodic billiard path of combinatorial length\footnote{The \emph{combinatorial length} of a periodic billiard path is defined as the number of times the point bounces on the boundary of the polygon before returning to its initial position \emph{with the same direction} it started its motion.} less than $\frac{1}{\epsilon}$ \cite{Schwartz06}. We refer the reader to \cite{Schwartz09} and the references therein for more details on this conjecture.

A common trick to study billiards in polygons is to produce from $P$ a translation surface $M(P)$ using an unfolding process. We follow the approach of J.C. Yoccoz to explain this process, see
\cite{Yoccoz10}, though this trick appears in different guises in the literature, as early as 1936 in a paper by Fox and Kershner~\cite{FoxKershner36} and in a 1975 paper by Katok and Zemlyakov~\cite{KatokZemliakov75}.


Let $\{e_1, e_2, \ldots, e_n\}$ be the sides of $P$, $r(e_j)$ be the linear part of the element in $\rm Isom(\R^{2})$ given by the reflection with respect to the line containing $e_j$ for $j\in\{1,2,\ldots,n\}$ and $R< \operatorname{O}(2,\R)$ be the subgroup generated by $\{r(e_1),\ldots,r(e_n)\}$. Note that this group is not necessarily finite, because the product of some reflections can produce a rotation about an irrational angle. We define a topological space $M(P)$ as the quotient of $P\times R$ by the following equivalence relation. Two points $(z,r)$ and $(z',r')$ are equivalent if and only if $z=z'$ and:
\begin{enumerate}
\item $r^{-1}r'=\bold{1}_R$ if $z$ belongs to the interior of $P$ (\ie we do not identify points in the interior of a fixed copy of $P$),
\item $r^{-1}r'\in\{\bold{1}_R,r(e_j)\}$ if $z$ belongs to the side $e_j$,
\item $r^{-1}r'\in R_z$ if $z$ is a vertex of $P$, where $R_z$ is the subgroup of $R$ generated by $r(e_i)$ and $r(e_j)$ if $z$ is the vertex given by the intersection of the sides $e_i$ and $e_j$.
\end{enumerate}
We denote by $\Sigma$ the image in $M(P)$ of the vertices of $P$. It is always possible to define a translation structure on $M(P)\setminus\Sigma$. Indeed, if we denote the interior of $P$ by $\mathring{P}$ then:
\begin{itemize}
\item for each $r\in R$ the map $
\mathring{P}\times\{r\}\to\R^2$ given by $(z,r)\to r(z)$ defines a chart and

\item if $z$ belongs to the side $e_j$ (but is not a vertex) of $P$ and $r\in R$, let $U=B(z,\epsilon)\cap P$ with $\epsilon>0$ small enough so that $U$ is an open half-disk plus its diameter. The  map:
$$
U\times\{r,r\cdot r(e_j)\}\to\R^2
$$
sending $(z,r)$ to $r(z)$ and $(z,r\cdot r(e_j))$ to $f_{r(e_j)}(r(z))$, where $f_{r(e_j)}$ is the reflection with respect to line containing the image of the side $e_j\times\{r\}\subset P\times\{r\} $ in the plane, defines a chart around $z$.
\end{itemize}
\begin{tcbexercise}{}{PolygonsBilliards}
Let $P$ be a Euclidean polygon.
\begin{enumerate}
\item Show that the charts defined above define a translation tructure on $M(P)\setminus\Sigma$.
\item Suppose that all interior angles of $P$ are commensurable with $\pi$ and let $N$ be the smallest integer such that any interior angle of $P$ can be written as $\pi\frac{m}{N}$ for some $m\in\N$. Show that in this case $R$ is isomorphic to the dihedral group of order $2N$:
$$
\langle a,b\hspace{1mm}|\hspace{1mm}a^2=1,b^N=1,(ab)^2=1\rangle.
$$
\item Let $z\in P$ be a vertex defined by the intersection of the sides $e_i$ and $e_j$ and let  as before $R_z$ be the subgroup of $R$ generated by $r(e_i)$ and $r(e_j)$. Prove that:
\begin{enumerate}
\item If the interior angle of $P$ at $z$ is of the form $\pi\frac{p}{q}$ with $p,q\in\N$ relatively prime, then $R_z$ is isomorphic to the dihedral group of order $2q$. Deduce that if we add the image of the vertex $z$ to $M(P)\setminus \Sigma$ we obtain a cone angle singularity of total angle $2p\pi$.
\item If the interior angle of $P$ at $z$ is of the form $\lambda\pi$, with $\lambda\in\R\setminus\Q$, then $R_z$ (and hence $R)$ is infinite. Deduce that if we add the image of the vertex $z$ to $M(P)\setminus \Sigma$ we obtain an infinite cone angle singularity.
\end{enumerate}
\end{enumerate}
\end{tcbexercise}

From the preceding exercise, we can deduce that Euclidean polygons can be classified into two types: those for which every interior angle is commensurable with $\pi$ and those for which this is not the case. The former are called \emphdef{rational polygons} and the later \emphdef{irrational polygons}.

\begin{figure}[!ht]
\begin{center}\includegraphics[scale=0.7]{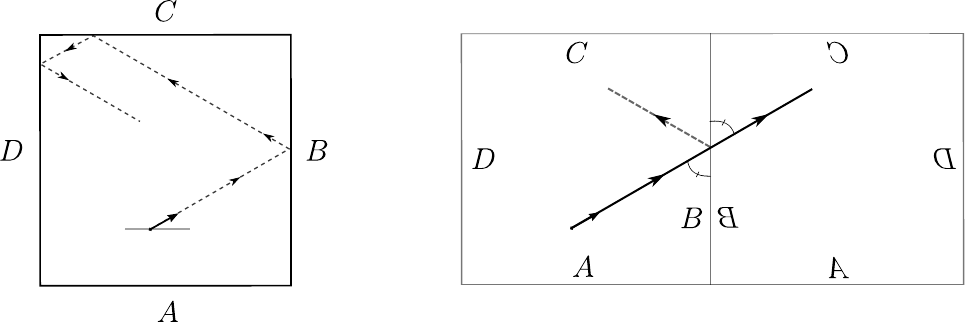}\end{center}
\caption{Billiard inside the unit square and the unfolding of a billiard trajectory.}
\label{fig:BilliardSquare}
\end{figure}

\textbf{Billiard trajectories and translation flows}. The main point of introducing the unfolding trick is to relate billiard trajectories on $P$ to trajectories of the translation flow on $M(P)\setminus \Sigma$. Let $\gamma(t)$, $t\in [0,T]$ be a billiard trajectory which bounces on the sides of $P$ at times $t_0=0<t_1<\ldots<t_N<t_{N+1}=T$. Denote by $e_{t_i}$ the side that $\gamma(t)$ encounters at $t=t_i$. Define $r_{t_0}:=\bold{1}_R$ and $r_{t_{i+1}}:=r_{t_{i}}\cdot r(e_{t_{i+1}})$. Then for every $r\in R$ the image on $M(P)\setminus \Sigma$ of the path defined by the formulas:
$$
\forall i \in \{0,1,\ldots,N\}, \quad \gamma_r(t) = (\gamma(t),r\cdot r_{t_i}) \quad \text{for $t_i \leq t\leq t_{i+1}$}
$$
defines a trajectory of a translation flow $F_{M(P),\theta}^t$ on $M(P)\setminus \Sigma$, for some $\theta\in\R/2\pi\Z$. Conversely, the image of every piece of trajectory of the translation flow on $M(P)\setminus \Sigma$ by the natural projection $\pi_P:M(P)\to P$ defines a billiard trajectory on $P$. Hence, we can reformulate many problems of billiards dynamics in terms of translation flows in $M(P)\setminus\Sigma$. In figure \ref{fig:BilliardSquare} we depict the unfolding process and the relation between billiard trajectories and the translation flow for the case of the billiard in the unit square.

\begin{figure}[!ht]
\begin{center}\includegraphics[scale=0.8]{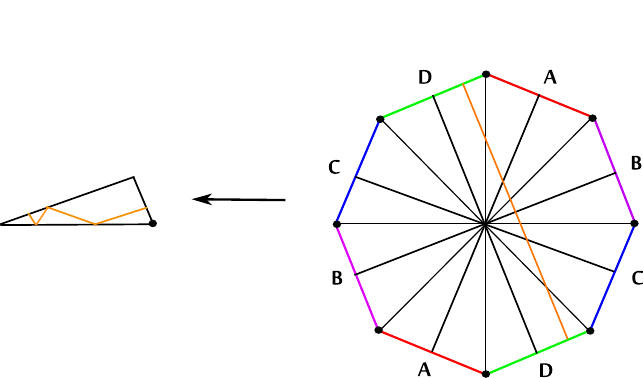}\end{center}
\caption{The billiard in the triangle $P$ with interior angles $(\frac{\pi}{2},\frac{\pi}{8},\frac{3\pi}{8})$ and the associated surface $M(P)$.}
\label{F:SPSurface}
\end{figure}

\textbf{The topology of $M(P)$}. If $P$ is a rational polygon, exercise \ref{exo:PolygonsBilliards} above implies that $M(P)$ is a compact translation surface tiled by finitely many copies of the polygon $P$. In figure \ref{F:SPSurface} we depict $M(P)$ for a triangle $P$ with interior angles $(\frac{\pi}{2},\frac{\pi}{8},\frac{3\pi}{8})$. Remark that $M(P)$ in this case is obtained by identifying opposite sides using translations on a regular octagon. By making a straightforward calculation on the number of copies of $P$ tilling $M(P)$ and applying the formula that relates the Euler characteristic of $M(P)$ to its genus, we deduce that this translation surface has genus 2 and only one conic singularity of total angle $6\pi$. In the same figure we illustrate how a billiard trajectory on $P$ unfolds in $M(P)$. In general, if $P$ has interior angles $\frac{p_i}{q_i}\pi$, a direct application of the Euler characteristic formula shows that

\begin{equation}
{\rm genus}(M(P))=1+\frac{N}{2} \left(n-2-\sum_{i=1}^{n}\frac{1}{q_{i}}\right)
\end{equation}
where  $N={\rm lcm}(q_1,\ldots,q_n)$. For a detailed proof of this formula see~\cite[Theorem 7.22]{Tabachnikov05}. In particular, the topology $M(P)$ when $P$ is a rational polygon \emph{depends only on the interior angles of $P$.}

If $P$ is an irrational polygon, Exercise~\ref{exo:PolygonsBilliards} above implies that $M(P)$ is not a manifold on points $p\in\Sigma$ coming from vertices where the interior angle is not commensurable with $\pi$. Indeed, since any such point $p$ is an infinite angle singularity, no neighbourhood of $p$ is compact and hence there cannot be a chart around this point into the plane. To avoid this nuisance we  \emph{remove first} from $P$ all vertices on which interior angles are not commensurable with $\pi$ and then we perform the unfolding construction. The result is a translation surface that we henceforth, abusing notation, denote by
 $M(P)$. Remark that in this case the group $R\cap \SO(2,\R)$ has index 2 in $R$ and is of the form $\Z^{\rho(P)}\times A$ where $A$ is finite Abelian group and the non-negative integer $\rho(P)$ is called the \emphdef[rank (of a polygon)]{rank} of $P$. If the interior angles of $P$ are $\lambda_j\pi$, $j=1,\ldots,N$ then:
\begin{equation}
	\label{eq:RankPolygon}
\rho(P)={\rm dim}_\Q(\sum_{j=1}^N\Q\lambda_j)-1
\end{equation}
The following theorem, whose proof will be discussed in Example~\ref{EXAMPLE:IrratBillAndProof} of Section~\ref{sec:CoveringSpaces} describes the topology of this surface.

\begin{tcbtheorem}{\cite{Valdez09_2}}{AnIrrationalBilliardIsLochNess}
Let $P$ be a Euclidean polygon with interior angles $\lambda_j\pi$, $j=1,\ldots,N$. If there exists $\lambda_j\in\R\setminus\Q$, then $M(P)$ is homemorphic to an infinite genus surface with one end.
\end{tcbtheorem}
In other words, \emph{all} translation surfaces stemming from irrational polygons are homeomorphic and of infinite type.\\

\begin{figure}[!ht]
\begin{minipage}{0.4\textwidth}
\begin{center}\includegraphics{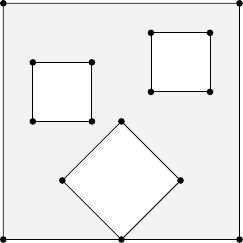}\end{center}
\subcaption{A rational generalized polygon.}
\end{minipage}
\hspace{0.1\textwidth}
\begin{minipage}{0.4\textwidth}
\begin{center}\includegraphics{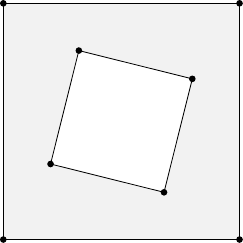}\end{center}
\subcaption{An irrational generalized polygon as considered in Exercise~\ref{exo:NonSimplePolygon}.}
\label{fig:InscribedSquares}
\end{minipage}
\caption{Two bounded generalized polygons.}
\label{fig:GeneralizedPolygons}
\end{figure}

It is fair to say that one can play billiards in domains of $\R^2$ which are more complicated than a Euclidean polygon. To fix ideas let us introduce a more general notion of polygon.
\begin{tcbdefinition}{}{GeneralizedPolygon}
A \emphdef{generalized polygon} is a closed subset $P$ of $\R^2$ whose boundary is a union of (Euclidean) segments such that for any compact set $K \subset \R^2$
the intersection $K \cap \partial P$ is a finite union of segments. The segments forming $\partial P$ are called the \emph{sides} of the generalized polygon.
\end{tcbdefinition}
All Euclidean polygons as well as the complement of their interior are generalized polygons. Two examples of generalized polygons are depicted in Figure~\ref{fig:GeneralizedPolygons}.

\begin{tcbdefinition}{}{RatIrratPolygon}
Let $P$ be a generalized polygon with sides $\{e_1,e_2,\ldots\}$ and $R<{\rm O}(2,\R)$ be the group generated by the linear parts of the reflections w.r.t. the lines containing $e_j$ for each $j\in\{1,2,\ldots\}$. Then $P$ is called a \emphdef[rational (generalized polygon)]{rational} generalized polygon if the group $R$ is finite and \emphdef[irrational (generalized polygon)]{irrational} otherwise.
\end{tcbdefinition}
\begin{tcbexercise}{}{NonSimplePolygon}
The purpose of this exercise is to illustrate that there are rather simple generalized polygons whose interior angles are all commensurable with $\pi$ but for which the unfolding process leads to an infinite-type translation surface. Consider two concentric and disjoint squares $P_{1}\subset P_{2}$ such that $P_{1}$ is obtained from $P_{2}$ by an homothety followed by a rotation by an angle $\lambda\pi$ with $\lambda\in\R\setminus\Q$ (see Figure~\ref{fig:InscribedSquares}). Let $P$ be the closure in the plane of $P_2\setminus P_1$. Prove that if we apply the unfolding process described in the preceding paragraph we obtain an infinite area infinite genus translation surface $M(P)$. Describe the set of conical singularities. 
\end{tcbexercise}

By definition, a compact generalized polygon $P$ is rational if and only if the translation surface $M(P)$ obtained by the unfolding process is compact (and hence of finite type). Exercise~\ref{exo:NonSimplePolygon} shows that there are \emph{irrational} compact generalized polygons whose interior angles are commensurable to $\pi$. As we will see later, the so-called wind-tree models are special classes of non-compact rational generalized polygons such that $M(P)$ is of infinite topological type and whose billiard dynamic is quite interesting.

\begin{tcbremark}{}{}
There are even more general spaces on which one can consider billiard dynamics and the unfolding trick. An example is the so-called \emph{parking garages} introduced by M.~Cohen and B.~Weiss, see~\cite{CohenWeiss2012}.
\end{tcbremark}

%
%
\subsection{Baker's surfaces and the  infinite Arnoux-Yoccoz surface}
\label{ssec:BakerBowmanArnouxYoccoz}

In this section we present two examples of finite area infinite-type translation surfaces
which appear when studying classical dynamical systems.

\emph{Baker's surfaces}. We construct for each $\alpha\in(0,1)$ a finite-area translation surface $B_\alpha$ of infinite genus. These examples are called \emph{baker's (or Chamanara) surfaces} because, as Chamanara, Gardiner and Lakic explain in~\cite{ChamanaraGardinerLakic}, they appear when considering the classical \emph{baker self-map} on the unit square. They are also examples of Hooper-Thurston-Veech surfaces, which we discuss in Sections~\ref{ssec:ThurstonVeechConstructionsIntro} and~\ref{sec:HooperThurstonVeechConstruction}.

\begin{figure}[!ht]
\begin{center}\includegraphics[scale=1.4]{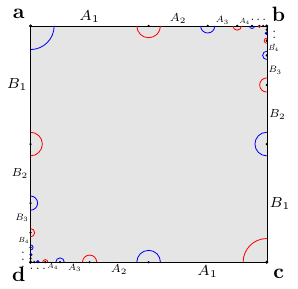}\end{center}
\caption{Baker's surface $B_{\alpha}$ for $\alpha=\frac{1}{2}$ and its wild singularity.}
\label{F:cham}
\label{fig:BakersSurface}
\end{figure}

Understanding the construction of baker's surfaces is easier when considering Figure~\ref{F:cham}.
Consider a square $abcd$ with $ab$ and $bc$ as its upper and right sides, respectively, each having a length of $\frac{\alpha}{1-\alpha}$. For each $i\in\N \setminus \{0\}$, define $a_i\in ab$ and $b_i\in bc$ as the point such that $|a_i b|=|b_i b|=\frac{\alpha^{i+1}}{1-\alpha}$. Analogously, let $c_i\in cd$ and $d_i\in da$ be such that $|c_i d|=|d_i d|=\frac{\alpha^{i+1}}{1-\alpha}$. Let us also denote $a_0 = d_0 = a$ and $b_0 = c_0 = c$. These points define a partition of the sides and each segment $\left[a_i,\, a_{i+1}\right]$, $\left[b_i,\, b_{i+1}\right]$, $\left[c_i,\, c_{i+1}\right]$ and $\left[d_i,\, d_{i+1}\right]$ has length $\alpha^{i+1}$ (see Figure~\ref{F:cham}). Now, identify by translation the horizontal sides $\left[a_i,\, a_{i+1}\right]$ with $\left[c_i,\, c_{i+1}\right]$ (labeled $A_i$ on Figure~\ref{F:cham}) and the vertical sides $\left[b_i,\, b_{i+1}\right]$ with $\left[d_i,\, d_{i+1}\right]$ (labeled $B_i$ on Figure~\ref{F:cham}). We denote by $\widehat{B_\alpha}$ the topological space that results from these identifications.

\begin{tcbexercise}{}{BakerSSurface}
Show that if we remove from
$\widehat{B_\alpha}$ the image of the points $\{a_i,b_i,c_i,d_i\}_{i=0}^\infty$, the result is a geometric translation surface, that is, it has a natural atlas whose change of coordinates are translations. We denote this surface by $B_\alpha$. Show that for each $\alpha\in(0,1)$ the metric completion of $B_\alpha$ w.r.t. the flat metric is obtained by adding just one point $x_\infty$. Moreover, show that this point is a wild singularity and that the metric completion of $B_\alpha$ is precisely the quotient $\widehat{B_\alpha}$ defined above.
\end{tcbexercise}

It is a result of Chamanara~\cite{Chamanara} that all $B_\alpha$ are homeomorphic: they all have infinite genus and only one end (see Section~\ref{Sec:Topological-classification-surfaces} for a detailed discussion on the topological classification of surfaces). Theorem~\ref{thm:AnIrrationalBilliardIsLochNess} implies that every $B_\alpha$
is homeomorphic to the topological surface obtained by unfolding an irrational polygon.

Each $B_\alpha$ is a translation surface formed only by flat points, \ie it contains no conical singularities. However, as the preceding exercise shows, a wild singularity appears when considering its metric completion, see Definition~\ref{def:Singularities}. We discuss wild singularities with more detail in the next chapter.

We revisit baker's surfaces in Section~\ref{ssec:VeechGroupBakersSurface}. There, we show that for rational parameters of the form
$\alpha=\frac{1}{q}$, $q\geq 2$, baker's surfaces have a lot of structure. More precisely, we show (see Proposition~\ref{prop:exa:TVBaker}) that $B_\frac{1}{q}$ is a Hooper--Thurston--Veech surface. These are, as discussed in Section~\ref{ssec:ThurstonVeechConstructionsIntro}, infinite-type analogs of translation surfaces obtained by a famous construction due to Thurston and Veech. As a consequence, following the work of Chamanara~\cite{Chamanara}, we prove in Theorem~\ref{thm:BakerVeechGroup} that the Veech group of $B_{\frac{1}{q}}$, $q\geq 2$, is (modulo conjugation) a group generated by two parabolic matrices (which together generate a free group) and the involution defined by rotating by $\pi$ the square in Figure~\ref{F:cham}.

\begin{tcbremark}{}{ChamanaraBakerIsTheSame}
\label{RMK:ChamanaraBakerIsTheSame}
Baker's surfaces were originally introduced by
R. Chamanara in \cite{Chamanara}, though the recipe for its construction is attributed to A. de Carvalho. By this reason some authors call these surfaces \emph{Chamanara's surfaces}.
It is important to remark that R. Chamanara constructed his surfaces in \cite{Chamanara} from a unit square, \emph{ergo}
if we apply the homothety of factor $\frac{1-\alpha}{\alpha}$
to the square $abcd$ introduced above and proceed with the construction of $B_\alpha$ we retrieve Chamanara's original construction.
\end{tcbremark}

\begin{tcbquestion}{}{BakerCylinders}
For which directions $\theta$ does the surface $B_\alpha$ have a cylinder, and for which directions does the translation flow $F_\theta^t$ have a dense trajectory?
\end{tcbquestion}

One of the consequences of $B_{\frac{1}{q}}$, $q\geq 2$, being a Hooper-Thurston-Veech surface is that the set of directions
$\theta$ for which $F_\theta^t$ decomposes the surface into infinitely many cylinders is infinitely countable. Moreover, using the work of Hooper~\cite{Hooper-infinite_Thurston_Veech}, one can show (though we do not discuss in this book how) that there are uncountably many directions $\theta$ for which $F_\theta^t$ has dense trajectories. Sadly, this describes only a measure zero set of directions.  I Chapter~\ref{chap:TranslationFlowsAndIET}, we prove that for any $\alpha$ and any $\theta$ the translation flow in direction $\theta$ in $B_\alpha$ has zero entropy (we also discuss the notion of entropy in that Chapter in detail). However, very little is known about the dynamics of the translation flow in $B_\alpha$ is a generic direction.

\begin{figure}[!ht]
\begin{center}\includegraphics[scale=.8]{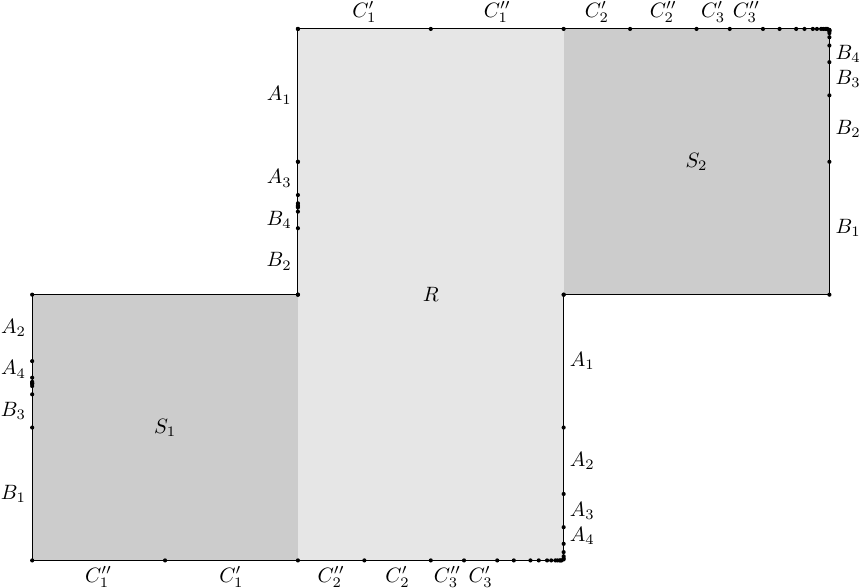}\end{center}
\caption{The infinite Arnoux-Yoccoz surface $(X_\infty,\omega_\infty)$.}
\label{infinite-arnoux-yoccoz}
\end{figure}

\emph{The infinite Arnoux-Yoccoz surface}. In~\cite{Bowman-Arnoux-Yoccoz}, Bowman provides an explicit description of the family of Arnoux-Yoccoz (translation) surfaces $(X_g,\omega_g)_{g\geq 0}$. Each element in this family is a genus $g$ translation surface having an affine involution $\rho_g$ and an affine pseudo-Anosov homeomorphism $\psi_g$ with dilatation $\alpha^{-1}>1$, where $\alpha$ is the unique positive root of $x^g+x^{g-1}+\cdots+x-1$. Moreover, he provides a triangulation for each $(X_g,\omega_g)_{g\geq 0}$ and, after some calculations, shows that as $g$ tends to infinity each triangle has a limiting position and it makes sense to consider the  geometric limit. This limit is an infinite-type translation surface $(X_\infty,\omega_\infty)$, which is illustrated in Figure~\ref{infinite-arnoux-yoccoz}. This figure is composed by two unit squares $S_1$ and $S_2$ plus a rectangle $R$ of width 1. Segments in the boundary with the same labels are identified by translations. Each segment whose label has $i\geq 1$ as subindex has length $\frac{1}{2^i}$.

As with Baker's surfaces, it can be shown that $(X_\infty,\omega_\infty)$ has only one singularity in its metric completion. Moreover, Bowman shows that $(X_\infty,\omega_\infty)$ is homeomorphic to an infinite genus surface with only one end and calculates explicitly two elements in its affine group: an orientation reversing involution $\rho_\infty$ and a `pseudo-Anosov' map $\psi_\infty$. The involution $\rho_\infty$ is easy to see from Figure~\ref{infinite-arnoux-yoccoz}: it sends the interior of the upper horizontal rectangle onto the interior of the inferior horizontal rectangle, each segment $A_i$ to a segment $B_i$, and segment
$C_i^{'}$ to $C_i^{''}$. On the other hand, the linear part of the hyperbolic element $\psi_\infty$ is $\begin{psmallmatrix}2 & 0\\0 & \frac{1}{2}\end{psmallmatrix}$, $\psi_{\infty}(R)$ is the rectangle formed by $S_1$ and the lower half of $R$, and $\psi(S_1\cup S_2)$ is the rectangle formed by $S_2$ and the upper part of $R$. As explained in [\emph{ibid.}], $\psi_\infty$ is a variant of the baker map.

Despite the resemblances, baker's surfaces and the infinite Arnoux-Yoccoz surfaces can be very different. More precisely, following Bowman, in Theorem~\ref{thm:InfiniteAYVeechGroup} we show that $\Aff(X_\infty,\omega_\infty)$ is generated by $\rho_\infty$ and $\psi_\infty$.
In particular, this group is virtually cyclic. In opposition, the Veech group of baker's surface $B_{\frac{1}{q}}$, for $q\geq 2$, always contains a free group of rank 2. We refer the reader to Theorem~\ref{thm:BakerVeechGroup} for a precise statement.

%
%
\subsection{The infinite staircase}
\label{ssec:InfiniteStaircase}

Infinite-type translation surfaces appear when considering infinite coverings of finite-type translation surfaces. This is a framework that we explore with more detail in Section~\ref{sec:CoveringSpaces} and we exploit recurrently along this text. The infinite staircase, introduced in Section \ref{SECC:WhatIsATranslationSurface}, was the first example of an infinite covering to be studied systematically, see
\cite{HooperHubertWeiss}. To be more precise, let $\T^2$ be the flat torus $\C/\Z^2$. The infinite staircase is a covering of $\T^2 \setminus 0$. Translation surfaces which are coverings of $\T^2$ (branched at most over $0$) are called \emphdef[square-tiled surface]{square-tiled surfaces} or \emphdef[origami]{origamis}. The infinite staircase is thus an infinite-type square-tiled surface and, as it can be seen from Figure \ref{fig:StaircaseFirst}, it has four infinite angle singularities.

\begin{figure}[!ht]
\begin{minipage}{.45\textwidth}
\begin{center}\includegraphics[scale=0.9]{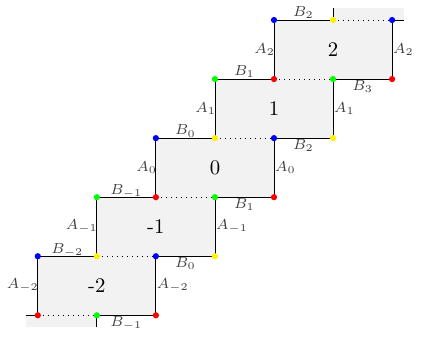}\end{center}
\subcaption{The infinite staircase. A (ramified) covering of a torus with
4 singularities of infinite conical type.}
\end{minipage}
\hspace{.1\textwidth}
\begin{minipage}{.45\textwidth}
\begin{center}\includegraphics[scale=0.8]{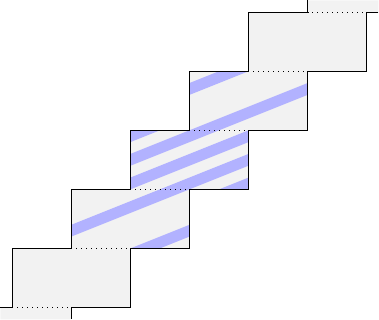}\end{center}
\subcaption{In direction $(p,q)$ with $p+q \equiv 1 \mod 2$, the staircase
decomposes into a countable union of cylinders. Here is the example $(p,q) = (5,2)$.}
\end{minipage}
\bigskip \\
\begin{minipage}{.45\textwidth}
\begin{center}\includegraphics[scale=0.9]{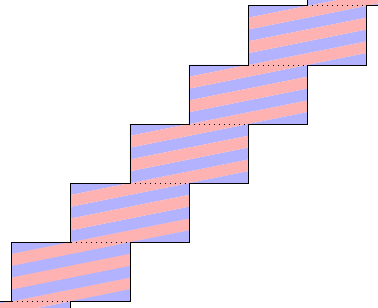}\end{center}
\subcaption{In direction $(p,q)$ with $p+q \equiv 0 \mod 2$ and relatively prime, the staircase
decomposes into a union of two strips. Here is the example of $(p,q) = (5,1)$.}
\end{minipage}
\hspace{.1\textwidth}
\begin{minipage}{.45\textwidth}
\begin{center}\includegraphics[scale=0.8]{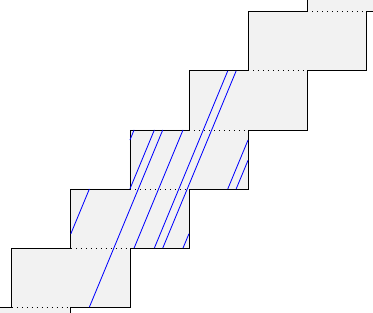}\end{center}
\subcaption{A piece of a trajectory in the direction stabilized by the affine map $\psi_h \psi_v$.}
\end{minipage}
\caption{The infinite staircase and some trajectories of its translation flows.}
\label{fig:StaircaseFlows}
\end{figure}

The infinite staircase has a lot of symmetries. Indeed, its affine group contains:
\begin{enumerate}
\item a group of translations isomorphic to $\Z$,
\item two elements\footnote{As we see later, these elements correspond to horizontal and vertical multitwists.} $\psi_h$ and $\psi_v$ with derivatives $D \psi_h = \begin{pmatrix}1&2\\0&1\end{pmatrix}$ and $D\psi_v = \begin{pmatrix}1&0\\2&1\end{pmatrix}$ and
\item a symmetry with derivative $-\bold{Id}$ that stabilizes the staircase at level $0$, see figure \ref{fig:StaircaseFlows} (a).
\end{enumerate}
\begin{tcbtheorem}{}{StaircaseProperties}
The infinite staircase has the following geometric and dynamical properties:
\begin{compactitem}
\item its affine group is generated by the elements listed above. In particular its Veech group
 is a lattice in $\SL(2,\R)$ generated by the matrices $\begin{pmatrix}1&2\\0&1\end{pmatrix}$ and $\begin{pmatrix}0&1\\-1&0\end{pmatrix}$,
\item if $\theta\in\R/2\pi\Z$ is a direction of rational slope, then the translation flow $F_\theta^t$ decomposes the infinite staircase into infinitely many cylinders or into two infinite strips, and
\item  if $\theta\in\R/2\pi\Z$ is a direction of irrational slope, the translation flow $F_\theta^t$
is ergodic with respect to Lebesgue measure.
\end{compactitem}
\end{tcbtheorem}

We refer the reader to~\cite{Dajani-Kalle2021} for background material on ergodic theory. In the result above, ergodicity in irrational directions is due to J.-P. Conze~\cite{Conze76}, \cite{ConzeKeane76} and K. Schmidt~\cite{Schmidt78}. It was generalized later by I.~Oren~\cite{Oren1983} and strengthened into a measure classification by J.~Aaronson, H.~Nakada, O.~Sarig and R.~Solomyak~\cite{AaronsonNakadaSarigSolomyak02}. The Veech group computation has been done by P.~Hooper, P.~Hubert and B.~Weiss~\cite{HooperHubertWeiss} and the presence of these affine symmetries provides an alternative way to prove ergodicity of the translation flow. All this will be discussed in
a forthcoming volume~\cite{DHV2}
dedicated mostly to dynamical aspects of infinite-type translation surfaces. Other results about this model include \cite{Kesten-uniform_distribution}, \cite{Ralston14-GenericDiscrepancy}, \cite{AvilaDolgopyatDuryevSarig-2015}.

\subsection{Wind-tree models}
\label{ssec:WindTreeIntro}
\label{ssec:Intro:WindTree}

\begin{figure}[!ht]
\begin{minipage}{.5\textwidth}
\begin{center} \includegraphics[scale=0.4]{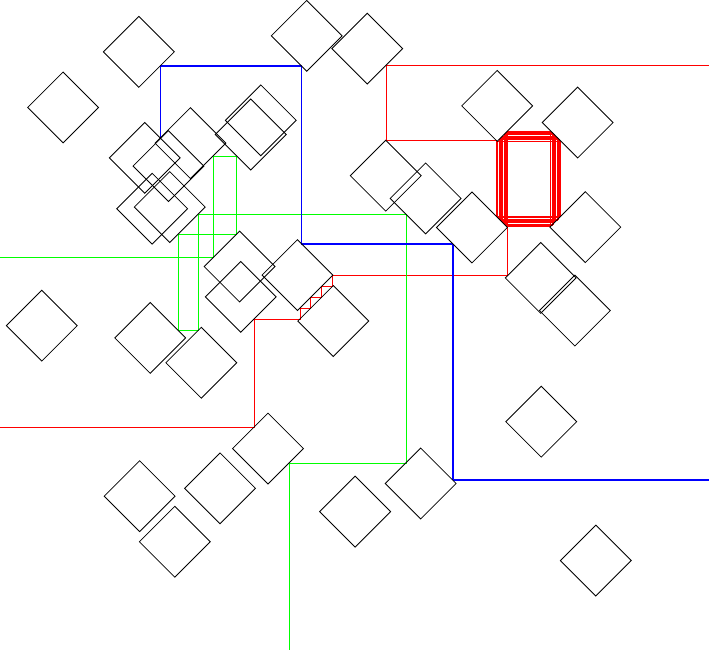} \end{center}
\subcaption{The Ehrenfests' wind-tree model}
\label{fig:EhrenfestOriginalWindtree}
\end{minipage}
\begin{minipage}{.5\textwidth}
\begin{center} \includegraphics[scale=0.3]{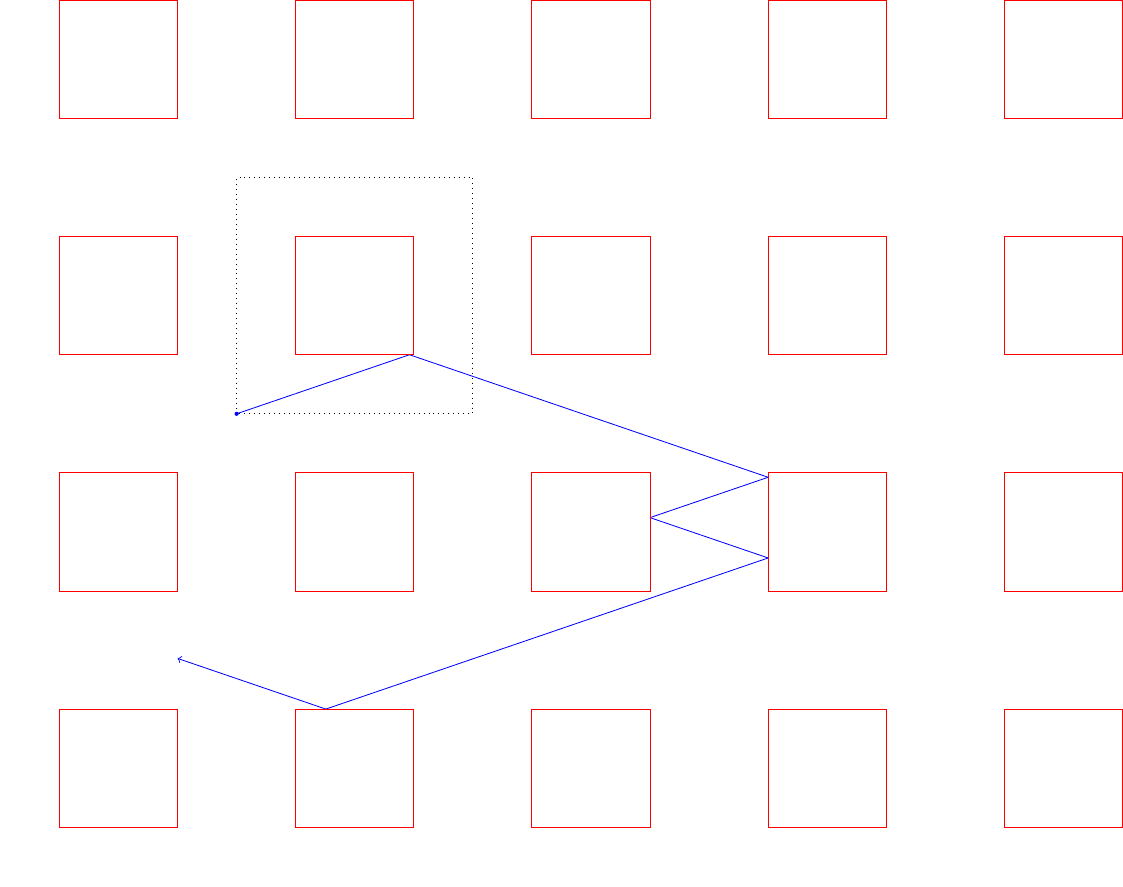} \end{center}
\subcaption{The Hardy-Weber periodic wind-tree model. In this picture the parameters are $a=b=1/2$.}
\label{fig:PeriodicWindtree}
\end{minipage}
\caption{The original Ehrenfests' wind-tree model and the periodic version from Hardy-Weber.}
\label{fig:WindtreeModels}
\end{figure}

In this section we present a special kind of billiards on generalized polygons: wind-tree models. Historically, it was Paul and Tatiana Ehrenfest \cite{EhrenfestEhrenfest12} who defined first the (random) wind-tree model as a simplified model of Lorenz gas. They were interested in the robustness
of Boltzmann's ergodic hypothesis. In this first version, square scatterers with random sizes are randomly placed in the plane, and the billiard ball could only move in a very specific direction (see Figure~\ref{fig:EhrenfestOriginalWindtree}).
Much later on, two other physicists, J.~Hardy and J.~Weber, studied a periodic version of the Ehrenfests' wind-tree model~\cite{HardyWeber80} (see Figure~\ref{fig:PeriodicWindtree}). However, they also considered the billiard flow only in some very specific directions. It wasn't until the work of P.~Hubert, S.~Leli\`evre and S.~Troubetzkoy~\cite{HubertLelievreTroubetzkoy11} that the study of the billiard flow in a generic direction was initiated. Since then a lot of work has been devoted to the study of periodic wind-tree models and their generalizations. See, for example,~\cite{AvilaHubert-recurrence}, \cite{DelecroixHubertLelievre14}, \cite{Delecroix-divergent}, \cite{FraczekUlcigrai-non_ergodicity}, \cite{FraczekUlcigrai-ergodic_examples}, \cite{DelecroixZorich-cries_and_whispers} or \cite{Pardo-windtree}. Despite a huge amount of work in the periodic case, very little is known about
the original Ehrenfests' wind-tree model. We discuss more on these random wind-tree models in Section~\ref{ssec:IntroRandomModels}. For the rest of this section we focus on Hardy-Weber's periodic wind-tree models.

\begin{tcbquestion}{}{DynamicsRandomEhrenfestModel}
Describe the dynamics of the billiard ball on the original (random) Ehrenfests' wind-tree model for a generic direction $\theta\in\R/2\pi\Z$.
\end{tcbquestion}

\emph{Periodic wind-tree models}. Let $\Z^2$ be the standard lattice in $\R^2$ and $(a,b)\in(0,1)^2$. For each $(m,n)\in\Z^2$ let us consider the rectangle $R_{m,n} = [m - a/2, m + a/2] \times [n - b/2, n+b/2]$ of size $a \times b$ centered at the lattice point $(m,n)$. We denote by $T_{a,b}$ the billiard table which is obtained by removing from the plane $\R^2$ the union of all rectangles $R_{m,n}$. We call the billiard table $T_{a,b}$ the \emph{periodic wind-tree model} with parameters $(a,b)$.  As discussed at the end of section \ref{ssec:PolygonalBilliards}, $T_{a,b}$ is a non-compact generalized polygon, and hence we can perform the unfolding trick on it.

\begin{figure}[!ht]
\begin{center}
\hspace{.02\textwidth}
\begin{minipage}{.65\textwidth}%
\begin{center}%
\includegraphics[scale=0.7]{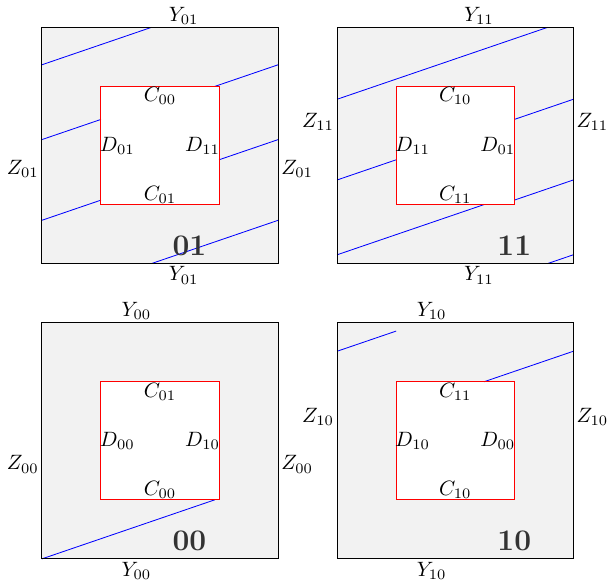}
\end{center}
\subcaption{The quotient translation surface $X_{a,b}$. This surface has genus $5$ and belongs to the stratum $\cH(2^4)$.}
\label{fig:WindtreeQuotientSurface}
\end{minipage}
\\
\medskip
\begin{minipage}{.65\textwidth}%
\begin{center}%
\includegraphics[scale=0.7]{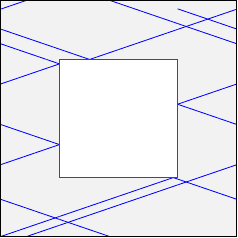}
\end{center}
\subcaption{The quotient billiard obtained by considering the fundamental domain from Figure~\ref{fig:PeriodicWindtree}.}
\label{fig:WindtreeQuotientBilliard}
\end{minipage}
\end{center}
\caption{The billiard and compact translation surface obtained as the quotient of the periodic wind-tree model with parameters $a=b=1/2$.
The quotient of the orbit shown in Figure~\ref{fig:PeriodicWindtree} is shown in blue in both figures.}
\label{fig:WindtreeUnfoldings}
\end{figure}

Given that the billiard flow on the table $T_{a,b}$ is invariant under the action by translations of $\Z^2$ and the interior angles of $T_{a,b}$ are all equal to $\frac{3\pi}{2}$, when we perform the unfolding construction on $T_{a,b}$ the result is an infinite-type surface $W_{a,b}$ which is a $\Z^2$-covering of the genus 5 translation surface $X_{a,b}$ illustrated in Figure~\ref{fig:WindtreeQuotientSurface}. It is easy to see that $W_{a,b}$ is a complete translation surface (for the flat metric) presenting only finite cone angle singularities. Note that the surface $X_{a,b}$ corresponds to the unfolding procedure applied to a billiard in a torus with a rectangular obstacle (see Figure~\ref{fig:WindtreeQuotientBilliard}). Using this, one can show that $X_{a,b}$ is a covering of a genus two surface $L_{a,b}$ in the stratum $\mathcal{H}(2)$. This seemingly trivial observation plays a fundamental role in the study of periodic wind-tree models. As we will see in Chapter~\ref{sec:CoveringSpaces}, it is possible to use standard arguments in covering space theory to prove that the infinite translation surface $W_\emph{a,b}$ always has infinite genus and one end.

The following statement summarizes the principal results that we prove in this book regarding the Veech groups of periodic wind-tree models. In the statement, $\Gamma(W_{a,b} \to L_{a,b})$ denotes the subgroup of the Veech group $\Gamma(W_{a,b})$ formed by the linear parts of elements in $\Aff(L_{a,b})$ which have a lift to $W_{a,b}$. As we see later, this is called the \emph{relative} Veech group of the covering $W_{a,b} \to L_{a,b}$. For more details see Definition~\ref{def:RelativeVeechGroup}.

\begin{tcbtheorem}{}{WindTreeVeechGroup0}
Let $(a,b)$ be parameters for the wind-tree model so that the Veech group of $L_{a,b}$ is a
lattice in $\SL(2,\R)$. Then
\begin{enumerate}
\item $\Gamma(W_{a,b} \to L_{a,b})$ has infinite index in $\Gamma(L_{a,b})$,
\item $\Gamma(W_{a,b} \to L_{a,b})$ contains a non-trivial normal subgroup,
in particular its limit set\footnote{See Appendix~\ref{app:FuchsianGroups} for more details about Fuchsian groups and their limit sets} is the whole circle, and
\item if furthermore $(a,b) \in \cE$, where $\cE$ is defined
in~\eqref{eq:DefinitionEParametersWindtree}, then $\Gamma(W_{a,b} \to L_{a,b})$
contains a parabolic element.
\end{enumerate}
\end{tcbtheorem}

In particular, for parameters $(a,b)$ such that that the Veech group of $L_{a,b}$ is a
lattice $\SL(2,\R)$, $\Gamma(W_{a,b} \to L_{a,b})$ is a Fuchsian group of the first kind which is not a lattice, and thus it is infinitely generated.

We now summarise below what is known about the dynamics of the billiard flow on $T_{a,b}$ (or equivalently the translation flow on $W_{a,b}$). We refer the reader to~\cite{KrerleyViana16} and~\cite{Aaronson} for a review of standard notions in dynamical systems such as recurrence or ergodicity in infinite-measure spaces.
\begin{tcbtheorem}{Generic and exceptional wind-tree behaviors}{WindtreeResults}
For every pair $(a,b)\in(0,1)^2$ the following is true for the translation flow $F^t_\theta$ on the translation surface $W_{a,b}$, obtained from the wind-tree model by unfolding, in a generic direction $\theta$:
\begin{compactitem}
\item the translation flow $F^t_\theta$ is conservative\footnote{Conservative in this context means that almost every orbit comes back arbitrarily close to its starting point.} and not ergodic,
\item for all point $x \in W_{a,b}$ with infinite forward orbit we have
\[
\limsup_{t \to \infty} \frac{\log \dist(x, F^t_\theta)}{\log t} = \frac{2}{3}.
\]
where the distance is the distance in $W_{a,b}$,

\end{compactitem}
However, for the parameters $a=b=1/2$
\begin{compactitem}
\item there are uncountably many directions $\theta$ for which the flow $F^t_\theta$ is completely dissipative
\item there are uncountably many directions $\theta$ for which the flow $F^t_\theta$ is ergodic.
\end{compactitem}
\end{tcbtheorem}

Recall that a dynamical system is \emph{dissipative} if the space is a union of wandering sets up to a zero measure. This is the opposite behavior of a conservative dynamical system.

In the result above, recurrence is attributed to A.~Avila and P.~Hubert~\cite{AvilaHubert-recurrence}. 
The rate of diffusion is a manifestation of the \emph{Kontsevich-Zorich phenomenon} and in this context, the result is due to V.~Delecroix, P.~Hubert and S.~Leli\`evre~\cite{DelecroixHubertLelievre14}.
The non-ergodicity is a phenomenon found by K.~Fraczek and C.~Ulcigrai~\cite{FraczekUlcigrai-non_ergodicity} and
the existence of divergent flow is due to V.~Delecroix~\cite{Delecroix-divergent}. Finally, the existence of ergodic directions is a consequence of the work of P.~Hooper~\cite{Hooper-infinite_Thurston_Veech} or a method relative to essential values. All these properties and most of these results are discussed in a (second) forthcoming volume~\cite{DHV2} dedicated mostly to dynamical aspects of infinite-type translation surfaces.


\subsection{Panov planes}
\label{ssec:PanovPlanes}


In \cite{Panov}, D. Panov provided one of the first studies of the dynamics of foliations $\mathcal{F}_\theta$ in an infinite area quadratic differential (see Section~\ref{ssec:Structures} for details on this foliations). His construction
was then revived recently in~\cite{JohnsonSchmoll14}, \cite{FraczekSchmoll-reflector}, \cite{Artigiani-Eaton_lenses} and~\cite{FraczekShiUlcigrai}.
\begin{tcbdefinition}{Panove plane}{PanovePlane}
A \emph{Panov plane} is a not identically zero meromorphic quadratic differential $\tilde{q}$ on
the plane $\C$ with at most simple poles and that is invariant under translation by a
lattice $\Lambda$. That is for all $v$ in $\Lambda$, if $\tilde{q}=fdz^2$, we have that $f(z + v) = f(v)$.

Equivalently, a Panov plane is the pull-back to the universal cover of a not identically zero meromorphic
quadratic differential $q$ on a torus.
\end{tcbdefinition}
As we discussed in Remark~\ref{rk:NoPoles}, if we want the Panov plane to be a holomorphic quadratic differential
one has to remove from the plane finitely many orbits of the lattice $\Lambda$ corresponding to the simple poles of the differential
$\widetilde{q}$.

\begin{tcbremark}{}{}
Any doubly periodic meromorphic function $f: \C \to \C$ can be written in terms of the Weierstrass $\cP_\Lambda(z)$
function and its derivative. 
\end{tcbremark}

\begin{figure}[!ht]
\begin{center}
\begin{minipage}{.75\textwidth}%
\begin{center}%
\includegraphics{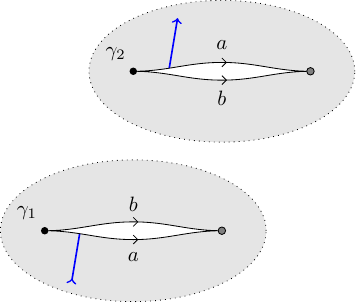}
\end{center}
\subcaption{Slit construction along $\gamma_1$ and $\gamma_2$. A geodesic segment
in the slit surface is drawn in blue.}
\label{fig:SlitConstruction}
\end{minipage}
\\
\medskip
\begin{minipage}{.75\textwidth}%
\begin{center}%
\includegraphics{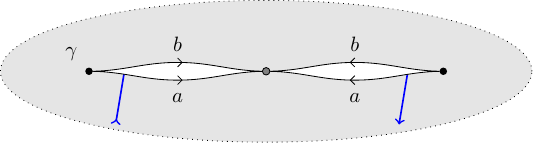}
\subcaption{The special case where $\gamma_1$ and
$\gamma_2$ are two oppositely oriented halves of a given segment $\gamma$.}
\end{center}
\end{minipage}
\caption{Two slit constructions.}
\end{center}
\end{figure}

In this introductory section we discuss in detail the case of differentials in the stratum $\cQ(2,-1^2)$
(that is, quadratic differentials on a genus 1 surface   with one
singularity of angle $4\pi$ and two singularities of angle $\pi$)
which includes Panov original example. Let us first start with a simple surgery
operation (see Figure~\ref{fig:SlitConstruction}).
\begin{tcbdefinition}{}{GluingSlits}
Let $M$ be half-translation surface (not necessarily connected) and $\gamma_1$ and $\gamma_2$ be two
open parallel oriented geodesic segments of the same length in $M$ with disjoint interiors and which are disjoint
from the singularities of $M$. The \emph{slit construction} performed on $M$
along $(\gamma_1,\gamma_2)$ is the half-translation surface $M'$ obtained by
cutting $M$ along $\gamma_1$ and $\gamma_2$ and gluing the left side of $\gamma_1$ to the right of $\gamma_2$ and vice-versa.
\end{tcbdefinition}

In this section we only consider slit constructions where $\gamma_1$ and $\gamma_2$
are the two oppositely-oriented halves of a given segment $\gamma$ as in the right
part of Figure~\ref{fig:SlitConstruction}.
In this situation the two endpoints of $\gamma$ become a simple zero of the quadratic differential defining $M'$ and the midpoint
(corresponding to the junction of $\gamma_1$ and $\gamma_2$) becomes two simple poles\footnote{Formally speaking these poles do not belong to $M'$ but to the metric completion $\widehat{M'}$, see remark \ref{rk:NoPoles}}.

\begin{figure}[!ht]
\begin{center}
\begin{minipage}{0.45\textwidth}
\begin{center}\includegraphics[scale=1]{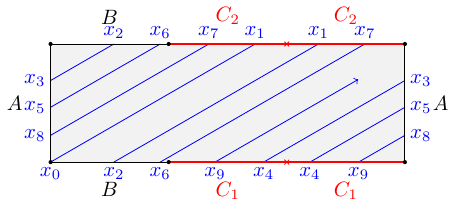}\end{center}
\subcaption{A flat picture of $(X,q)$ in $\cQ(2,-1^2)$.}
\end{minipage}
\\
\medskip
\begin{minipage}{0.45\textwidth}
\begin{center}\includegraphics{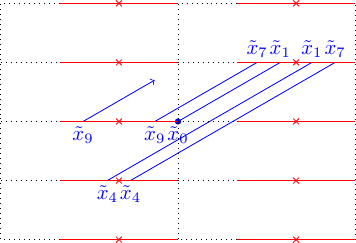}\end{center}
\subcaption{A flat picture of $(\widetilde{X},\widetilde{q})$.}
\end{minipage}
\caption{The original example of D. Panov. The torus (on the left) is built
from three unit squares with identifications given by the upper case letters.
Its universal cover (on the right) can be seen as the flat plane with slits
displaced periodically. A separatrix starting from $x_0$ and its lift starting
from $\tilde{x}_0$ are show in blue. Each segment $[x_i,x_{i+1}]$ corresponds
to a piece of the orbit. The lift of $x_i$ is denoted by $\tilde{x}_i$.}
\label{fig:OriginalPanovPlane}
\end{center}
\end{figure}

Now we arrive at the original Panov construction. Consider the quadratic
differential $(\C, dz^2)$, let
$\Lambda_0 = \Z (3,0) \oplus \Z (0,1)$ and $v_0 = (2,0)$. We perform a slit
construction along each segment $[t, t+v_0]$ for $t \in \Lambda$ and obtain
the original Panov plane as drawn in Figure~\ref{fig:OriginalPanovPlane} (b).
We denote by $(\C, \tilde{q})$ this quadratic differential.
Since this construction is invariant under the translation by the lattice
we obtain a differential $(X, q)$ in genus one that belongs to
the stratum $\cQ(2,-1^2)$ (see Section~\ref{sec:CoveringSpaces} for more about coverings).

Now let us consider an affine pseudo-Anosov element $\phi \in \Aff(X, q)$ of the torus covered by the original Panov plane. It has an induced action $\phi_*$ on $H_1(X; \Z) \simeq \Z^2$. This action is tightly linked to the behavior
of the lift of the invariant foliations of $\phi$ to the plane.

\begin{figure}[!ht]
\begin{center}
\begin{minipage}{0.8\textwidth}
\begin{center}\includegraphics[scale=0.8]{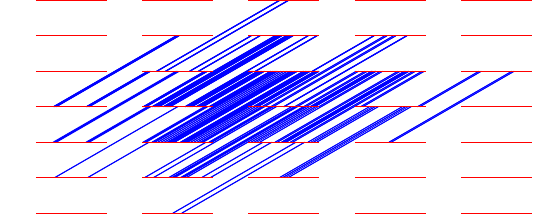}\end{center}
\subcaption{Elliptic case (original Panov example, $\phi_h \phi_v^{-1} \phi_h$) in direction $(\sqrt{3},1)$.\\}
\label{fig:KZPanovElliptic}
\end{minipage}
\\
\begin{minipage}{0.6\textwidth}
\begin{center}\includegraphics[scale=0.5]{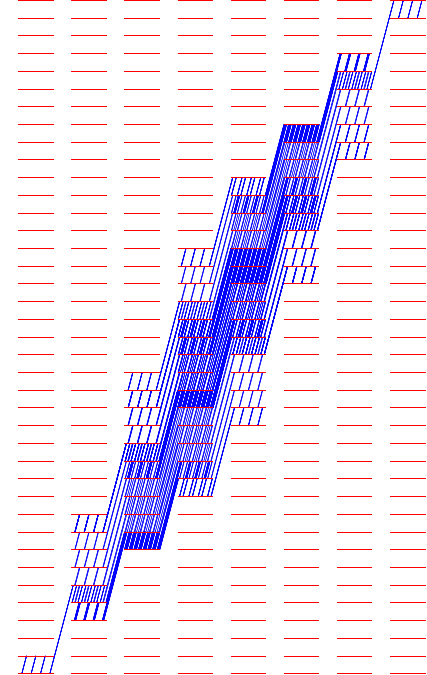}\end{center}
\subcaption{Hyperbolic case ($\phi_h \phi_v$) in direction $(6, -3 + \sqrt{21})$.}
\label{fig:KZPanovHyperbolic}
\end{minipage}
\caption{Kontsevich-Zorich phenomenon for lifting of invariant foliations of pseudo-Anosov homeomorphisms in Panov planes.}
\label{fig:KontsevichZorichPanovPlanes}
\end{center}
\end{figure}

The following exercise provides some computations of $\phi_*$ (see
Figure~\ref{fig:KontsevichZorichPanovPlanes} for an application).
\begin{tcbexercise}{}{PanovePlaneAffineGroup}
 \label{exercise:PanovPlaneAffineGroup}
Let $(\C,\widetilde{q}) \to (X,q)$ be the Panov plane of Figure~\ref{fig:OriginalPanovPlane}.
\begin{compactenum}
\item Let $a$ (respectively $b$) be a closed curve in $(X,q)$ obtained by joining
the left side denoted $A$ to the right one (resp. the bottom side denoted $B$ to the top one).
Show that $H_1(X; \Z) = \Z a \oplus \Z b$.
\item Show that there exists two elements $\phi_h$ and $\phi_v$ in $\Aff(X,q)$ so that
 $\phi_h$ acts as a horizontal twist, $\phi_v$ as a vertical twist and their derivatives
 are $D \phi_h = \displaystyle \begin{pmatrix}1&3\\0&1\end{pmatrix}$
and $D \phi_v = \displaystyle \begin{pmatrix}1&0\\1&1\end{pmatrix}$
(\textit{note: this is a special case of the Thurston-Veech construction of
Section~\ref{ssec:ThurstonVeechConstructionsIntro}}).
\item Let $\rho: \Aff(X,q) \to \SL(2,\Z)$ be defined by $\rho(\phi) = \phi_*$ after having
identified $H_1(X; \Z)$ to $\Z^2$ using the basis $\{a,b\}$ from the first question. Show
that
\[
\rho(\phi_h) = \begin{pmatrix}1&1\\0&1\end{pmatrix}
\quad \text{and} \quad
\rho(\phi_v) = \begin{pmatrix}1&0\\1&1\end{pmatrix}.
\]
\item Show that $\rho(\phi_h \phi_v^{-1} \phi_h)$ is elliptic of order 4.
\item Show that $\rho(\phi_h \phi_v)$ is hyperbolic.
\item Show that $\rho(\phi_h \phi_v^{-2} \phi_h)$ is parabolic.
\end{compactenum}
\end{tcbexercise}

\begin{tcbtheorem}{}{KZPAForPanovPlanes}
Let $(\C, \widetilde{q}) \to (X,q)$ be a Panov plane. Let $\phi:X \to X$
be an affine pseudo-Anosov homeomorphism of $(X,q)$ and let
$\phi_*: H_1(X; \Z) \to H_1(X; \Z)$ be its induced action on homology.
Then
\begin{compactenum}
\item if $\phi_*$ is the identity or elliptic, then the eigenfoliations of $\phi$ in $(X,q)$
lift to ergodic foliations of $\C$.
\item If $\phi_*$ is hyperbolic, there exist directions $v^+$ in $\R^2 \setminus \{0\}$ so that each
lift $\widetilde{\gamma}$ of a leaf of the unstable  eigenfoliation
of $\phi$ is at bounded distance away from $\R v^+$.
\end{compactenum}
\end{tcbtheorem}
Concrete examples of each case of the above result are provided in Figure~\ref{fig:KontsevichZorichPanovPlanes}. The
relevant computations were proposed in Exercise~\ref{exercise:PanovPlaneAffineGroup}. The second item in the statement above is a manifestation of the so called Kontsevich-Zorich phenomenon which will be discussed in~\cite{DHV2}. 

\begin{figure}
\begin{center}
\includegraphics{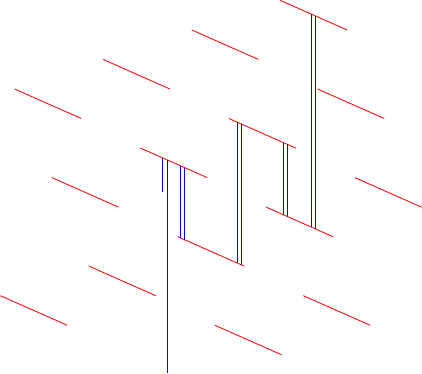}
\end{center}
\caption{A Panov plane coming from a quadratic differential in $\cQ(2,-1^2)$.}
\label{fig:GeneralizedPanov}
\end{figure}

We now study a natural generalization of Panov's example.
Let $\Lambda \subset \C$ be a lattice and $v \in \C^*$ so that
the open segment $\{tv: t \in (0,1)\}$ is disjoint from $\Lambda$.
We define the quadratic differential $\tilde{q}_{\Lambda, v}$ in the
plane by adding a slit parallel to $v$ based at each point of $\Lambda$. Because
this construction is invariant under the action of $\Lambda$ by translation we can
also make this slit construction in the torus $\C / \Lambda$ and obtain a
quadratic differential $q_{\Lambda,v}$. This is illustrated in Figure~\ref{fig:GeneralizedPanov}.

\begin{tcbexercise}{}{IsomQ2ppH00}
In this exercise we consider finite quotients and coverings of finite-type quadratic differentials. Recall that in a stratum
$\cQ(\kappa)$ of quadratic differentials the zeros and punctures are labeled (see Section~\ref{ssec:Structures}).
\begin{compactenum}
\item Show that any differential in $\cQ(2,-1^2)$ is hyperelliptic\footnote{Recall that a quadratic differential $(X,q)$ is hyperelliptic if it is a double covering of a quadratic differential on the sphere $(\C \P^1, \overline{q})$. That is to say, there exists an involution of the surface $M$ that preserves the quadratic differential and such that the quotient is (topologically) a sphere.}. Show that the quotient quadratic differential on the sphere has 4 poles and no zero (up to a choice of labeling of the poles, it is a differential in $\cQ(0,-1^4)$).
\item \label{quest:isom1} Show that the quotient studied in the previous question provides an $\SL(2,\R)$-equivariant map between a subspace of $\cQ(2, -1^2, 0^2)$ and $\cQ(0,-1^4)$. Show that most differentials in $\cQ(0,-1^4)$ admit 4 preimages under this map.
\item \label{quest:isom2} By considering a different half-translation covering of $\cQ(0,-1^4)$, show that there exists an isomorphism between $\cQ(0,-1^4)$ and a subspace of $\cH(0^6)$.
\item Express the maps of questions~\ref{quest:isom1}
and~\ref{quest:isom2} in terms of the lattices $\Lambda$ and non-zero vectors $v \in \C^*$
that we used to build the differential $q_{\Lambda,v}$ in $\cQ(2,-1^2)$.
\end{compactenum}
\end{tcbexercise}

\begin{figure}[!ht]
    \begin{center}
\begin{minipage}{.3\textwidth}
\begin{center}\includegraphics{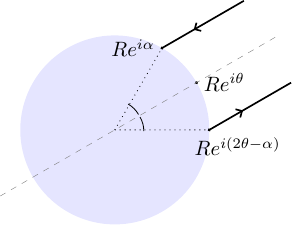}\end{center}
\subcaption{An Eaton lense is a perfect retroreflector: a trajectory entering
with direction $\theta$ exits with direction $-\theta$.}
\label{fig:EatonLensTrajectory}
\end{minipage}
\hspace{.2\textwidth}
\begin{minipage}{.3\textwidth}
\begin{center}\includegraphics{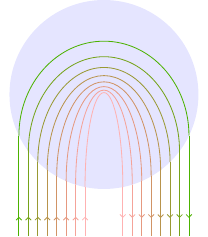}\end{center}
\subcaption{A Eaton lens with a beam of vertical light rays.}
\label{fig:EatonLens}
\end{minipage}
\\
\begin{minipage}{.3\textwidth}
\begin{center}\includegraphics{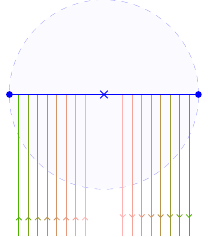}\end{center}
\subcaption{A quadratic slit (with two poles and a double zero) with
the same beam of rays.}
\label{fig:QuadraticSlit}
\end{minipage}
\caption{From a Eaton lens to its flat model.}
\label{fig:EatonLensOriginalAndFlat}
\end{center}
\end{figure}

The space $\cQ(2,-1^2)$ (or equivalently $\cH(0,0)$, see Exercise~\ref{exo:IsomQ2ppH00}) plays a distinguished role. First of all because of its relation to 2-dimensional affine lattices
\footnote{The moduli space of 2-dimensional lattices $\SL(2,\R) / \SL(2,\Z)$
is naturally identified to the moduli space of volume one flat tori $\cH_1(0)$. The
moduli space of affine lattices $\SL(2,\R) \rtimes \R^2 / (\SL(2,\Z) \rtimes \Z^2)$
is the space of pairs $(\Lambda, \Lambda + v)$ where the translate $v + \Lambda$
can be identified to a point on the torus $\C / \Lambda$. It can be identified to the
union $\cH_1(0,0) \sqcup \cH_1(0)$ where
$\cH_1(0,0)$ corresponds to affine lattices $v + \Lambda$ so that $v \not \in \Lambda$. This
identification is $\SL(2,\R)$-equivariant and permits to see $\cH_1(0,0)$ as a dense
$\SL(2,\R)$-invariant set in the space of affine lattices.}. Another
reason that motivates its study is a model that shares some analogy with the periodic
wind-tree model: a \emph{circular Eaton lens}. This is a physical gadget that
has the following property: every light ray that enters the lens is reflected in the
same direction but with opposite orientation (see Figure~\ref{fig:EatonLens}).
More precisely the Eaton lens of radius $R$ is a ball of radius $R$ such that a
trajectory entering at $R e^{i \alpha}$ in direction $e^{i \theta}$ will exit from
$R e^{i (2 \theta - \alpha)}$ with direction $e^{- i \theta}$ (see Figure~\ref{fig:EatonLensTrajectory}). The precise behavior
of the trajectory inside the ball is not important for our purposes.
Periodic configurations of Eaton lenses are discussed in the following exercise.

\begin{tcbexercise}{}{EatonLens}
Let $\Lambda \subset \C$ be a lattice and $R > 0$ so that $\displaystyle R < \min_{u,v \in \Lambda, u\not=v} \|u - v\|$.
Let $E_{\Lambda,R}$ be the plane $\C$ on which an Eaton lens of radius $R$ is placed at
each point of $\Lambda$ (the condition on $R$ ensures that the lenses do not intersect). We have
a well defined non-orientable foliation in $E_{\Lambda,R}$ for each direction
$\theta$.
\begin{compactenum}
\item Show that the vertical flow in $E_{\Lambda,R}$ can be identified with vertical flow
in a Panov plane $\tilde{q}_{\Lambda,v_0}$ for a unique $v_0 \in \C$
outside of the lenses (in other words each Eaton lens can be replaced by a slit as depicted in
Figure~\ref{fig:EatonLensOriginalAndFlat}).
\item Show that by applying a rotation one can naturally identify the
foliation in direction $\theta$ of $E_{\Lambda,R}$ to the vertical foliation of
$e^{-i\theta} E_{\Lambda,R}$.
\item Deduce that the one-parameter family of differentials in $\cQ(2,-1^2)$ obtained
by rotating the direction in $E_{\Lambda,R}$ is different from the family
$\{e^{i \theta} q_{\Lambda,v_0}\}$.
\end{compactenum}
\end{tcbexercise}

\begin{tcbtheorem}{generic and non-generic behavior for $\cQ(2,-1^2)$}{GenericBehaviorQ2pp}
Let $q_{\Lambda,v}$ be a quadratic differential in $\cQ(2,-1^2)$.  Then for
almost every direction $\theta$ each leaf of the foliation $\widetilde{\mathcal{F}_\theta}$ in the universal cover
stays at bounded distance from a line.

On the other hand, let $(\C,\widetilde{q}) \to (X,q)$ be the Panov plane of
Example~\ref{fig:OriginalPanovPlane}. Then there are uncountably
many directions in which the foliation $\widetilde{\mathcal{F}_\theta}$ is ergodic.
This set of directions contains all eigendirections of
pseudo-Anosov homeomorphisms in $\Aff(X,q)$ so
that their action on homology $\phi_*$ is elliptic.
\end{tcbtheorem}

This result is, in spirit, very similar to Theorem~\ref{thm:WindtreeResults}: one can describe the generic behavior of the dynamics, but for some exceptional cases the dynamics exhibit a dramatic contrast. We do not provide a proof in this text. The first item in Theorem~\ref{thm:GenericBehaviorQ2pp} is another manifestation of the Kontsevich-Zorich phenomenon that, in this form, is due to K.~Fracezk and M.~Schmoll~\cite{FraczekSchmoll-reflector}. A stronger form of
this theorem has been recently proved by K.~Fraczek, R.~Shi and C.~Ulcigrai~\cite{FraczekShiUlcigrai}.
The phenomenon described by the second item, is related to specific directions for which the associated
Kontsevich-Zorich cocycle is the identity. It can be obtained from two very different
sources. Either from general ergodic argument about $G$-coverings
using essential values
or from the Hooper-Thurston-Veech construction for infinite-type surfaces discussed
 in Section~\ref{ssec:ThurstonVeechConstructionsIntro}.

\smallskip


Ergodicity implies that almost every leaf is dense, but in general it is hard to provide a concrete example of such a leaf.
This is precisely one of the main contributions of D.~Panov~\cite{Panov}: the construction of a dense leaf. We will prove in a (second) forthcoming volume~\cite{DHV2},
dedicated mostly to dynamical aspects of infinite-type translation surfaces, a more general version due to C. Johnson and M.~Schmoll~\cite{JohnsonSchmoll14}.
\begin{tcbtheorem}{\cite{Panov}, \cite{JohnsonSchmoll14}}{PanovMinimality}
Let $(\C, \widetilde{q}) \to (X,q)$ be a Panov plane and assume
that $(X,q)$ admits a pseudo-Anosov affine element $\phi$ such that:
\begin{itemize}
\item $\phi$ has a fixed point $p$,
\item the action on homology $\phi_*$ is elliptic but not $\pm 1$.
\end{itemize}
Then the lifts of each separatrix issued from $p$ in the stable or unstable eigendirection of $\phi$ are dense in $(\C,\widetilde{q})$.
\end{tcbtheorem}
This result applies in particular to the separatrix drawn in Figure~\ref{fig:KZPanovElliptic}. We will not provide a proof of it in this text.

\smallskip

\begin{figure}[!ht]
\begin{minipage}{.45\textwidth}
\begin{center}
\includegraphics{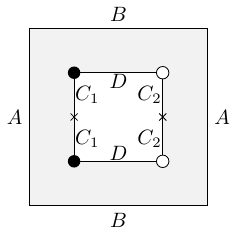}
\end{center}
\subcaption{The quotient $(X_h,q_h)$.}
\label{fig:WindtreeHypQuotient1}
\end{minipage}
\hspace{.05\textwidth}
\begin{minipage}{.45\textwidth}
\begin{center}
\includegraphics{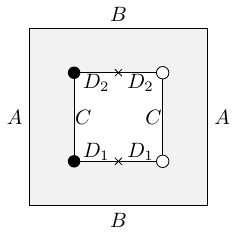}
\end{center}
\subcaption{The quotient $(X_v, q_v)$.}
\label{fig:WindtreeHypQuotient2}
\end{minipage}
\caption{Two quotients of the surface $X_{a,b}$ from Figure~\ref{fig:WindtreeQuotientSurface} obtained
as the quotient of the wind-tree model. Both are hyperelliptic quadratic differentials that belong
to $\cQ(1^2, -1^2)$.}
\label{fig:WindtreeHypQuotients}
\end{figure}

We finish this Section on Panov planes by explaining a relationship between those steaming from elements in $\cQ(1^2,-1^2)$ and wind-tree models. Recall that the wind-tree model can be thought as a $\Z^2$-periodic
translation surface $W_{a,b}$ whose quotient $X_{a,b}$ has genus $5$.
The surface $X_{a,b}$ has many symmetries and in particular two quotients $(X_h, q_h)$
and $(X_v, q_v)$ that are quadratic differentials on tori. These surfaces are shown
in Figure~\ref{fig:WindtreeHypQuotients}.

The quotients $(X_h, q_h)$ and $(X_v, q_v)$ will be of great use when computing the constant $2/3$ appearing in Theorem~\ref{thm:WindtreeResults}. This point of view on the wind-tree model is the starting point of the article~\cite{DelecroixZorich-cries_and_whispers} and~\cite{Pardo-windtree}. In these works, the authors study generalizations of the wind-tree model with obstacles more complicated than rectangles but for which a quotient still exists in $\cQ(1^{2n}, -1^{2n})$. Thanks to this fact, the diffusion rate of the translation flow can be computed in these more complicated surfaces.

\subsection{Hooper-Thurston-Veech construction}
\label{ssec:ThurstonVeechConstructionsIntro}
On a widely circulated preprint Thurston explained how to construct affine pseudo-Anosov homeomorphisms on any finite-type surface. Given two non-homotopic curves, $\alpha$ and $\beta$, that fill a topological surface $S$, there exists a unique half-translation structure $M$ on $S$ which admits two affine automorphisms $\phi_h$ and $\phi_v$ which act as a Dehn twist\footnote{For a precise definition of Dehn twist we refer the reader to Chapter 3 in \cite{FarbMargalit}.} along $\alpha$ and $\beta$ respectively. These have derivatives $D \phi_h = h_\lambda := \begin{psmallmatrix}1 & \lambda\\0 & 1\end{psmallmatrix}$ and $D \phi_v := v_{-\lambda} = \begin{psmallmatrix}1 & 0\\-\lambda & 1\end{psmallmatrix}$, where $\lambda=n(\alpha,\beta) > 0$ is some \emph{integer} number. In particular, $M$ is a square-tiled surface.

Because $\phi_h$ and $\phi_v$ are affine homeomorphisms of the same half-translation surface, they generate a subgroup of $\Aff(M)$ whose  elements are mostly pseudo-Anosov homeomorphisms. Thurston's construction was later extended by himself and Veech to the case of filling multicurves~\cite{Thurston88}, \cite{Veech89}, case in which the parameter $\lambda > 0$ is not necessarily an integer anymore. In particular this extended construction produces half-translation surfaces with affine pseudo-Anosov homeomorphisms which are not square-tiled surfaces. In this section we first recall in detail Thurston-Veech construction. Then we illustrate with a family of examples arising from the infinite staircase how this construction can be extended for infinite-type surfaces. The details of this extension, that we call the Hooper-Thurston-Veech construction for it was developed by P. Hooper \cite{Hooper-infinite_Thurston_Veech}, are explained in full detail in Chapter \ref{ch:Symmetries}.

\textbf{Historical note}. The Thurston-Veech construction seems to be forgotten by the community during the 1990's until a Thursday evening in July 2003 when, after the traditional \emph{Bouillabaisse}\footnote{Bouillabaisse is a traditional Proven\c{c}al fish stew originating from the port city of Marseille.} dinner at CIRM (France), John H. Hubbard revived the construction by explaining it to a large tipsy audience. For this reason, sometimes people refer to surfaces obtained by the Thurston-Veech construction as \emph{Bouillabaisse surfaces}.

\textbf{Thurston-Veech construction}. Throughout this section we use the following notation for parabolic matrices in $\SL(2,\R)$
\begin{equation} \label{eq:HLambdaVLambda}
h_\lambda := \begin{pmatrix}1&\lambda\\0&1\end{pmatrix}
\qquad \text{and} \qquad
v_\lambda := \begin{pmatrix}1&0\\\lambda&1\end{pmatrix}.
\end{equation}
Recall that the \emph{modulus} of a cylinder
$C$ is the real number: $$\mu(C) = \frac{\height(C)}{\circumference(C)}.$$
The basis of the Thurston-Veech construction is the following classical statement. For a proof see~\cite{HubertSchmidt}.

\begin{tcblemma}{}{ParabolicCylinderDecomposition}
Let $M$ be a finite-type translation or half-translation surface.

If the horizontal foliation of $M$ decomposes into $r$ cylinders of commensurable moduli $\mu_i =  \frac{n_i}{\lambda}$, then there exists a unique affine automorphism $\phi$ of $M$ that stabilizes pointwise the boundaries of the cylinders and has derivative $D\phi = h_{ \lambda}$.

Conversly, if $\phi$ is a parabolic affine homeomorphism of $M$ with derivative $h_\lambda$ then $M$ decomposes into a finite family of maximal horizontal cylinders $\{H_1,\ldots,H_r\}$ of commensurable moduli $\mu(H_i) = \frac{p_i}{q_i}\frac{1}{\lambda}$. Moreover, a finite positive power of $\phi$ stabilizes pointwise the boundary of each cylinder.
\end{tcblemma}
A parabolic automorphism that stabilizes the boundary of cylinders such as in Lemma~\ref{lem:ParabolicCylinderDecomposition} acts as a power of a Dehn twist in each cylinder. Such element is called a \emph{multitwist}. We introduce below the formal topological definition of multitwist. It is a bit technical since we cover both finite-type and infinite-type surfaces.
\begin{tcbdefinition}{}{Multitwist}
Let $S$ be a (topological) surface. A
simple closed curve in $S$ is
\emphdef[essential (curve)]{essential} if its
complement does not contain a component homeomorphic to
a disk or a punctured disk.
A disjoint union $\alpha = \cup \alpha_i$ of simple closed curves is
\emphdef{locally finite} if for any connected compact set
$K$ in $S$ the number of connected components of
$K \cap \alpha$ is finite.

A \emphdef{multicurve} in $S$ is a locally-finite union
$\displaystyle \alpha = \bigcup_{i \in I} \alpha_i$ of pairwise disjoint
and pairwise non-homotopic essential simple closed curves $\alpha_i$.

The \emphdef{multitwist} $T_\alpha$ associated to the multicurve
$\alpha$ is the mapping class element\footnote{A \emphdef{mapping class element} is
an homeomorphism of $S$ up to isotopy. In the context of multitwists $T_\alpha$
the Dehn twist $T_{\alpha_i}$ along a single component of $\alpha$ is a
well defined mapping class. One can choose the $T_{\alpha_i}$ to have disjoint
support in which case it is easy to see that the infinite product
$T_\alpha$ is well-defined.} defined by the (finite or infinite) product
$\displaystyle \prod_{i \in I} T_{\alpha_i}$
of simple Dehn twists $T_{\alpha_i}$ along the curve $\alpha_i$.
More generally, if $\{m_i\}_{i \in I}$ is a set of positive integers,
we define a multitwist along $\alpha$ with multiplicities $\{m_i\}_{i\in I}$ by
$\prod_{i \in I} T_{\alpha_i}^{m_i}$.
\end{tcbdefinition}

\begin{tcbdefinition}{}{AffineMultitwist}
Let $M$ be a translation or half-translation. A parabolic multitwist $f\in\Aff(M)$ which preserves the orbits of the translation flow \footnote{or leaves of the corresponding foliation in the case of a half-translation surface} parallel to the eigendirection of $Df$ is called an \emphdef{affine multitwist}.
\end{tcbdefinition}

The proof of the first part of Lemma~\ref{lem:ParabolicCylinderDecomposition}
is left to the reader (the argument is contained in the first part of Exercise~\ref{exo:FiniteThurstonVeech}). The second part is a bit more delicate, a proof can be found in~\cite[Lemma 4, Section 1.4]{HubertSchmidt}\footnote{Be careful that in the reference~\cite{HubertSchmidt} the authors used an opposite definition for modulus of a cylinder, namely $\operatorname{circumference}(C)/\operatorname{height}(C)$}.

Note that by applying the matrix $\begin{psmallmatrix} \cos(\theta) & -\sin(\theta) \\ \sin(\theta) & \cos(\theta) \end{psmallmatrix}$ to $M$, a similar statement holds in any direction $\theta$. In particular, if the vertical direction admits a cylinder decomposition with commensurable moduli, then the multitwist in the vertical direction has derivative
\[
\begin{pmatrix} 0 & -1\\1 & 0\end{pmatrix}
h_\lambda
\begin{pmatrix} 0 & -1\\1 & 0\end{pmatrix}^{-1}
=
v_{-\lambda}.
\]
Note that the derivative of a vertical multitwist is $v_{-\lambda}$ and not $v_{\lambda}$!

In the following exercise, we invite the reader to do the converse of the Thurston-Veech construction. Namely start from a translation surfaces with transverse multitwists and analyze the coordinates of the translation surface in terms of the topology of the surface.
\begin{tcbexercise}{}{CylinderDecomposition}
Let $M$ be a finite-type half-translation surface with conical singularities $\Sigma \subset \widehat{M}$ (recall
that $M$ is necessarily punctured at the simple poles, if any). Assume that $M$ admits vertical and horizontal
cylinder decompositions $H = \{H_1, \ldots, H_r\}$ and $V = \{V_1,\ldots,V_s\}$. Let
$\alpha = \{\alpha_1,\ldots,\alpha_r\}$ and $\beta = \{\beta_1,\ldots,\beta_s\}$ be the core curves of the
cylinders in the families $H$ and $V$ respectively. Prove that:
\begin{compactenum}
\item $\alpha$ and $\beta$ are multicurves in $M$ (that is, their components are simple closed curves, pairwise non-intersecting and pairwise non-homotopic).
\item The multicurves fill the surface, that is $M \setminus (\alpha \cup \beta)$ is a finite union of topological disks, possibly punctured.
\setcounter{EnumerateCounterContinuation}{\theenumi}
\end{compactenum}
Now assume furthermore that all the cylinders in $H$ and $V$ have modulus $1/\lambda$ and let
$E = (E_{ij})_{1\leq i \leq r, 1 \leq j \leq s}$ be the matrix whose coefficient
$E_{ij}$ is the geometric intersection number between $\alpha_i$ and $\beta_j$ \ie the number of time they intersect).
Let also $\textbf{h}_h \in \R^r_+$ and $\textbf{h}_v \in\R^s_+$ be the vectors of heights of cylinders in
$H$ and $V$.
\begin{compactenum}
\setcounter{enumi}{\theEnumerateCounterContinuation}
\item Show that the circumferences of horizontal cylinders (respectively, vertical cylinders)
are given by the coordinates of the vector $E \textbf{h}_v$ (respectively $E^t \textbf{h}_h$).
\item Deduce that $E \textbf{h}_v = \lambda \textbf{h}_h$ and $E^t \textbf{h}_h = \lambda \textbf{h}_v$.
\item Applying the Perron-Frobenius theorem to $E E^t$ and $E^t E$ show that $\lambda$ is uniquely determined by the matrix $E$
and that the pair $(\textbf{h}_h, \textbf{h}_v)$ is uniquely determined up to scaling.
\item Show that the vector $(\textbf{h}_h, \textbf{h}_v) \in \R^{r+s}$ is an eigenvector with eigenvalue $\lambda$ of the $(r+s) \times (r+s)$ matrix $\begin{psmallmatrix} 0 & E^t \\ E & 0\end{psmallmatrix}$.
\item Show that there exists a unique affine multitwist $\phi_h$ (respectively $\phi_v$) with $D \phi_h = h_\lambda$ (resp. $D \phi_v = v_{-\lambda}$)
that preserves the cylinder decomposition $H$ (resp. $V$) and acts in each cylinder as a Dehn twist.
\end{compactenum}
\end{tcbexercise}
The following exercise is precisely the Thurston-Veech construction.
\begin{tcbexercise}{}{FiniteThurstonVeech}
Let $\alpha = \{\alpha_1,\ldots,\alpha_r\}$ and $\beta = \{\beta_1,\ldots,\beta_s\}$ be two transverse multicurves that fill a finite-type surface $S$
(possibly with punctures). Let $E := (\iota(\alpha_i, \beta_j))_{1 \leq i \leq r, 1 \leq j \leq s}$ the matrix of geometric intersection numbers.
\begin{compactenum}
\item Show that there is a unique half-translation structure $M$ on $S$ (up to isotopy) and a unique positive real number $\lambda$ so that $M$ admits horizontal and vertical affine multitwists $\phi_h$ and $\phi_v$ with derivatives $D \phi_h = h_\lambda$ and $D \phi_v = v_{-\lambda}$ that act as the product of the Dehn twists along the curves $\alpha_i$ and $\beta_j$ respectively. \emph{Hint}: use exercise \ref{exo:CylinderDecomposition} to find the desired half-translation structure. The idea is to find the correct coordinates for a horizontal and a vertical cylinder decomposition on $S$.
\setcounter{EnumerateCounterContinuation}{\theenumi}
\end{compactenum}
Now, let us consider the general Thurston-Veech construction by fixing multiplicities \linebreak$\textbf{m}=(m_1,\ldots,m_r)\in \N^r$, $\textbf{n}=(n_1,\ldots,n_s)\in \N^s$. Define the two $r \times s$ matrices $F_h$ and $F_v$ by:
$$
(F_h)_{ij} = m_i E_{ij} \qquad (F_v)_{ij} = E_{ij} n_j.
$$
\begin{compactenum}
\setcounter{enumi}{\theEnumerateCounterContinuation}
\item Show that there exists $\lambda, \textbf{h}_h, \textbf{h}_v$ so that $F_h \textbf{h}_h = \lambda \textbf{h}_v$ and $F_v^t \textbf{h}_v = \lambda \textbf{h}_h$.
\item What can be said about uniqueness of $\lambda$, $\textbf{h}_h$ and $\textbf{h}_v$?
\item Deduce that there is a unique half-translation structure $M$ on $S$ (up to isotopy) that admits two affine multitwists $\phi_h$ and $\phi_v$ with derivative $h_\lambda$ and $v_{-\lambda}$ and act as the product of the $m_i$-th power of the Dehn twist along the curves $\alpha_i$ and the product of the $n_j$-th power of the Dehn twist along the curves $\beta_j$ respectively.
\end{compactenum}
\end{tcbexercise}
An avid reader can consult the details in the original articles~\cite{Thurston88}, \cite{Veech89} or the more recent~\cite{HubertLanneau06} or~\cite{Leininger04}.

\begin{tcbexercise}{}{ThurstonVeechNotTranslation}
Apply the Thurston-Veech construction of Exercise~\ref{exo:FiniteThurstonVeech}
to the two multicurves $C_h$ and $C_v$ and vectors $\textbf{m} = (1,1)$ and
$\textbf{n} = (1,1)$ in the following genus 2 surface
\begin{center}
\includegraphics{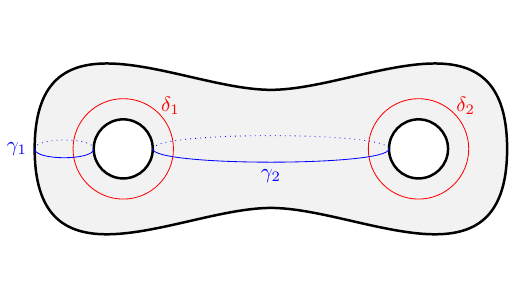}
\end{center}
Make a (flat) picture of the surface and determine its stratum.
\end{tcbexercise}

\begin{tcbexercise}{}{SingularitiesGeometryThurstonVeech}
Let $\alpha = \{\alpha_i\}_{i=1,\ldots,r}$ and $\beta = \{\beta_j\}_{j=1,\ldots,s}$ be two transverse multicurves filling a topological surface $S$.
Let $M$ be the half-translation surface with underlying topological surface $S$ obtained from the Thurston-Veech construction (see Exercise~\ref{exo:FiniteThurstonVeech}).
\begin{compactenum}
\item Show that the angles of the conical singularities of $M$ can be read from the multicurves: there is exactly one conical singularity in each component $U$ of the complement of $S \setminus (\alpha \cup \beta)$ and its angle is the number
of components of $\alpha$ (or equivalently $\beta$) that bound $U$ multiplied by $\pi$.
\item Deduce that the poles come from bigons and hence need to be in components of $S \setminus (\alpha \cup \beta)$ containing at least one puncture.
\item Show that $M$ is a translation surface if and only if there exists a coherent way of orienting the multicurves so that their geometric intersection coincides with the algebraic intersection.
\item Deduce that the stratum in which the half-translation surface belongs only depends on $\alpha$ and $\beta$ and is not affected by the possible multiplicity vectors $\mathbf{m}$ and $\mathbf{n}$.
\item Using Thurston-Veech's construction construct an example of a half-	translation surface in the principal stratum of genus 2 surfaces $\cQ(1,1,1,1)$.
\end{compactenum}
\end{tcbexercise}

The horizontal and vertical directions in the infinite staircase $S$ (see Section~\ref{ssec:InfiniteStaircase}) define two transverse cylinder decompositions. The moduli of all these cylinders is equal to $\frac{1}{2}$ and there exist "generalized multitwists" $\phi_h$ and $\phi_v$ with derivatives $h_2$ and $v_{-2}$ that act as a single Dehn twist in each cylinder. This is a particular case of an infinite Thurston-Veech construction that is due to Hooper~\cite{Hooper-infinite_Thurston_Veech} and that we call the Hooper-Thurston-Veech construction. All details regarding this construction are given in Section~\ref{sec:HooperThurstonVeechConstruction}.

\begin{figure}[H]
\begin{center}%
\includegraphics{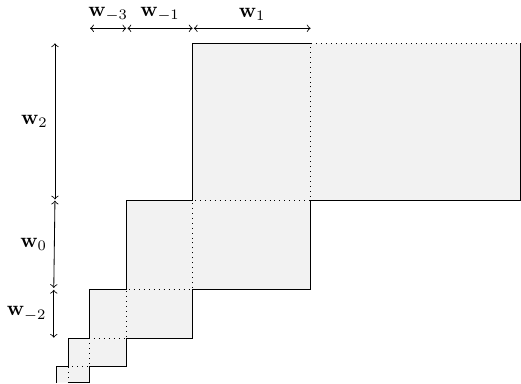}%
\hspace{1cm}%
\includegraphics{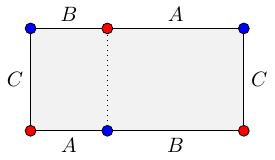}%
\end{center}
\caption{On top: the $\lambda$-staircase with $\lambda=27/13$ (opposite sides are identified by translation). At the bottom, the dilation surface obtained by quotienting the left picture by the affine dilation.}
\label{fig:LambdaStaircase}
\end{figure}

Let us recall that for finite-type surfaces the Thurston-Veech construction produces a translation (or half-translation) structure which is unique up to scaling (this is a consequence of Perron-Frobenius theorem as we saw in Exercise~\ref{exo:FiniteThurstonVeech}). This is not the case anymore for the Hooper-Thurston-Veech construction! Let us consider $\lambda > 2$. The equation $x^2 - \lambda x + 1$ has two distinct positive real roots $r_+ := \frac{\lambda + \sqrt{\lambda^2 - 4}}{2}$ and $r_- := \frac{1}{r_+}$. Let $\textbf{h}_n := r_+^n$ for $n \in \Z$. This sequence satisfies the equation
\begin{equation} \label{eq:WidthsLambdaStaircase}
\textbf{h}_{n-1} + \textbf{h}_{n+1} = \lambda \textbf{h}_n.
\end{equation}
Now we construct an infinite-type surface called the $\lambda$-staircase by gluing horizontal cylinders of height $\textbf{h}_{2n}$ with vertical cylinders of height $\textbf{h}_{2n+1}$ in the same pattern as the original infinite staircase, see Figure~\ref{fig:LambdaStaircase}. Because of equations~\eqref{eq:WidthsLambdaStaircase}, all horizontal and vertical cylinders in the translation surface $S_\lambda$ have modulus $1/\lambda$. For this reason they admit generalized multitwists with derivatives $h_\lambda$ and $v_{-\lambda}$. We denote $G_\lambda$ the group generated by $h_\lambda$ and $v_\lambda$. The following result describes the dynamics of the translation flow on  $\lambda$-staircases (we refer the reader to Appendix~\ref{appendix:sec:RosenGLambda} for a detailed study of $G_\lambda$).
\begin{tcbtheorem}{}{LambdaStaircase}
Let $\lambda > 2$ and $S_\lambda$ be the $\lambda$-staircase built from the sequence $\textbf{h}_n = \left(\frac{\lambda + \sqrt{\lambda^2 - 4}}{2}\right)^n$. Then the Veech group of $S_\lambda$ is the group generated by $G_\lambda$ and $\begin{psmallmatrix}r_+&0\\0&r_+\end{psmallmatrix}$. The Veech group is identified with the affine group via the derivative map. Moreover:
\begin{enumerate}
\item if $\theta$ is a parabolic direction of $G_\lambda$ (\ie fixed by a conjugate of $\phi_h$ or $\phi_v$)
then the translation flow in direction $\theta$ decomposes  $S_\lambda$ into an infinite union of cylinders
of modulus $1/\lambda$,
\item if the direction $\theta$ is in the limit set\footnote{See Appendix~\ref{app:FuchsianGroups} for more details about Fuchsian groups and their limit sets. To any non-vertical direction $\theta\in\R/2\pi\Z$ one can associate the \emph{co-slope} $\frac{\cos(\theta)}{\sin(\theta)}\in\partial\mathbb{H}^2=\R$.}  of $G_\lambda$ but neither parabolic nor fixed by
a conjugate of $\phi_h \phi_v$ then $F^t_\theta$ is ergodic. Moreover there
exists a direction $\theta'$ so that the translation flow in $S_\lambda$ in direction $\theta$ is conjugate
to the translation flow in $S$ in direction $\theta'$ up to a time change,
\item if the direction $\theta$ does not belong to the limit set of $G_\lambda$ or is fixed by
a conjugate of $\phi_h \phi_v$ then the translation flow $F^t_\theta$ in direction $\theta$ in $S_\lambda$ is
completely dissipative.
\end{enumerate}
\end{tcbtheorem}
The last item of Theorem~\ref{thm:LambdaStaircase} shares some similarity with the fact that the flow in the standard infinite staircase ($\lambda=2$) in the direction $\pi/4$ is completely dissipative (the surface is made of two infinite strips), see Theorem \ref{thm:StaircaseProperties}. In the situation of the $\lambda$-staircase the situation is somehow more drastic: there is an open set of directions for which the surface decomposes as a union of two "affine infinite strips". In order to study this phenomena it will be
convenient to consider the quotient of $S_\lambda$ by the unique affine
homeomorphism with derivative $\begin{psmallmatrix}r_+ &0\\0&r_+\end{psmallmatrix}$, see Figure~\ref{fig:LambdaStaircase}. The quotient is no longer a translation surface but a \emph{dilation surface}, that is, the transition maps of the atlas are dilations $z \mapsto a z + b$ with $a \in \R$ and $b \in \C$.

Dilation surfaces are discussed in Section~\ref{sec:CoveringSpaces}. In Section~\ref{sec:HooperThurstonVeechConstruction} we introduce the Hooper-Thurston-Veech construction and discuss the geometry of the surfaces obtained with it.

\subsection{Random models}
\label{ssec:IntroRandomModels}
In this section we discuss some random models of infinite-type translation surfaces. They are random in the sense that they depend on some parameters and we would like to understand the dynamics of the translation flow for a given random parameter. They might be thought of as ``translation flows in random environment".

There are two distinct ways of considering a random parameter. Either in the topological sense: we want to understand the dynamics for a dense $G_\delta$ set of parameters (\emph{generic parameter}), or in the probabilistic sense: for a full-measure set of parameters (\emph{typical parameter}). The topological version, and the measurable version, lead to different kinds of questions. For instance, in the context of circle rotations, the set of Liouville numbers\footnote{A real number in $\alpha \in [0,1]$ is called \emph{Liouville} if its continued fraction expansion $[0; a_1, a_2, \ldots]$ is so that $\lim a_n = +\infty$.} although forming a dense $G_\delta$ set in $[0,1]$, has zero Lebesgue measure. In more general situation, it is often the case that the parameters that are very well approximated by rational numbers form dense $G_\delta$ sets. For example, in~\cite{KerckhoffMasurSmillie86} S.~Kerckhoff, H.~Masur and J.~Smillie prove the existence of a dense $G_\delta$ set of polygons whose billiard flow is ergodic. This result is (easily) derived from the (hard) fact that ergodicity is prevalent in rational polygons.

\emph{M\'alaga map}. We now consider the maps $T_{\underline{\alpha}}$ introduced by A.~M\'alaga in her PhD thesis~\cite{Malaga14}. As we see later they can be thought as the first-return map of the translation flow on particular class of infinite-type translation surfaces
called \emph{generalized staircases}.

Let $I = [0,1) \times \Z$, and denote $I_n = [0,1) \times \{n\}$. Let $\underline\alpha = (\alpha_n)_{n\in\Z}$ be a sequence of numbers in $[0,1)$ and $f:[0,1) \to \{1,-1\}$ the piecewise constant map with value $1$ on $[0,1/2)$ and $-1$ on $[1/2,1)$. We define a map $T_{\underline\alpha}:I \to I$ by
\begin{equation}
T_{\underline{\alpha}} (x,n) = (x + \alpha_n \mod 1,\ n + f(x + \alpha_n\mod 1)).
\end{equation}
In other words, if we are on an interval $I_n$ at level $n$ we rotate first by an angle $\alpha_n$ and, depending on which half of $I_n$ it lands on, the point $x+\alpha_n$ either climbs to the interval $I_{n+1}$ or descends to $I_{n-1}$. When the parameter $\underline\alpha$ is a constant sequence then this map encodes the dynamics of the translation flow on the infinite staircase (see Section~\ref{ssec:InfiniteStaircase}) in a given direction.

It is easy to check that the Lebesgue measure $\lambda$ on $I$ is invariant by the transformation $T_{\underline{\alpha}}$. Therefore there is room for many standard questions from the infinite ergodic theory point of view: \eg conservativity and ergodicity. To go beyond these cases, it is often simpler to state results about the dynamics of the map $T_{\underline\alpha}$ that are valid only for a subset of parameters $\underline\alpha$. In that sense $T_{\underline{\alpha}}$ is a random model since we are interested in picking the sequence $\underline\alpha$ at random and then study the dynamics.

\begin{figure}[!ht]
\begin{center}
\includegraphics{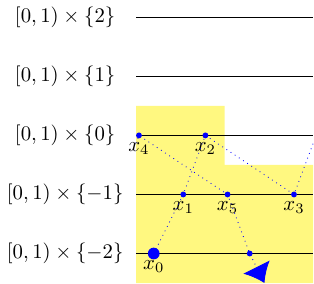}
\end{center}
\caption{A M\'alaga map with $\alpha_0=1/2$. Any trajectory starting at a negative (respectively positive) level $n$ will remain in non-positive levels (resp. non-negative levels) forever. An orbit labeled $x_0$, $x_1=T_{\underline{\alpha}}(x_0)$, $x_2=T^2_{\underline{\alpha}}(x_0)$, \ldots is shown in blue (we also draw artificial dotted segments between $x_i$ and $x_{i+1}$ to produce a better visualization of the orbit). The yellow region delimits a subset of $I = [0,1) \times \Z$ invariant under the map $T_{\underline{\alpha}}$.}
\label{fig:MalagaMapBarrier}
\end{figure}

The only known strategy to prove results on the maps $T_{\underline\alpha}$ is based on the fact that the value $\alpha_n=1/2$ creates a barrier (see Figure~\ref{fig:MalagaMapBarrier}).
\begin{tcbtheorem}{\cite{Malaga14}}{MalagaRecurrence}
Let $\underline\alpha$ be a sequence so that $1/2$ is both an accumulation point of $(\alpha_n)_{n \geq 0}$ and $(\alpha_n)_{n \leq 0}$. Then the map $T_{\underline \alpha}$ is conservative, that is, almost every orbit comes back arbitrarily close to its starting point.
\end{tcbtheorem}
This theorem is a consequence of the so-called \emph{boxes lemma} that will be introduced
and discussed in a (second) forthcoming volume~\cite{DHV2}
dedicated mostly to dynamical aspects of infinite-type translation surfaces.
As a result, one has
\begin{tcbcorollary}{}{MalagaMapRandomRecurrent}
Let $\mu$ be the Lebesgue measure on $[0,1)$. Then for $\mu^{\otimes \Z}$-almost every $\underline{\alpha} \in [0,1)^\Z$, the map $T_{\underline{\alpha}}$ is recurrent.
\end{tcbcorollary}
\begin{proof}[of Corollary~\ref{cor:MalagaMapRandomRecurrent}]
We claim that if $(X_n)_{n \geq 0}$ is a sequence of independent and identically distributed uniform random variables in $\left[0,1\right)$ then
for almost every $\omega$ and for any $s \in \left[0,1\right)$ there exists a divergent sequence of positive
integers $(n_k)_k$ so that
\[
\lim_{k \to \infty} X_{n_k}(\omega) = s.
\]
Assuming the claim with the value $s=1/2$ one can apply Theorem~\ref{thm:MalagaRecurrence}.

The proof of the claim is a standard probability argument. Any non-empty open set $U$
in $\left[0, 1 \right)$ has positive Lebesgue-measure. Hence almost surely,
$X_n \in U$ for infinitely many $n$. One then concludes by considering the
countable collection of open balls $B_{q,k}$ centered at rational points $q$ with radius $1/k$
and the fact that a countable intersection of almost sure events is almost sure.
\end{proof}
Contrarily to recurrence properties, ergodicity needs quantitative estimates. The only result in this direction for the maps $T_{\underline\alpha}$ is also due to M\'alaga.
\begin{tcbtheorem}{\cite{Malaga14}}{MalagaMinimalErgodic}
The set of parameters $\underline{\alpha}$ in $\left[0,1\right)^\Z$ for which $T_{\underline{\alpha}}$ is minimal and ergodic contains a dense $G_\delta$ set.
\end{tcbtheorem}

\smallskip

\begin{figure}[!ht]
\begin{center}%
\includegraphics{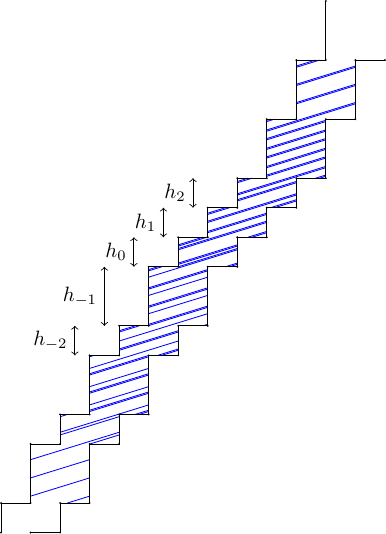}%
\end{center}
\caption{An orbit in some $\underline{h}$-staircase where each height $h_n$ is either $1$ or $2$. This example is an infinite square-tiled surface. The first return map of the translation flow in direction $\theta$ can be identified with a map $T_{\underline\alpha}$ for a certain sequence $\underline\alpha$.}
\label{fig:MalagaSurface}
\end{figure}

Recall the infinite staircase that was introduced in Section~\ref{ssec:InfiniteStaircase}. We now consider staircases where each step has width 2 but can have any height. Given a bi-infinite sequence of positive real numbers $\underline{h} = (h_n)_{n\in\Z}$ we consider the rectangles $R_n := [0,2] \times [0,h_n]$. From them we build the translation surface $M_{\underline{h}}$ by performing for each $n \in \Z$ the following identifications by translations:
\begin{itemize}
\item the left side of $R_n$ is glued to its right side,
\item the left half of the top side of $R_n$, $[0,1] \times \{h_n\}$ is glued to the right part of the bottom side of $R_{n-1}$, $[1,2] \times \{0\}$,
\item the right part of the top side of $R_n$, $[1,2] \times \{h_n\}$ is glued to the left part of the bottom side of $R_{n-1}$, $[0,1] \times \{0\}$.
\end{itemize}
We call the surface $M_{\underline{h}}$ the $\underline{h}$-infinite staircase or simply a generalized staircase. An example can be seen in Figure~\ref{fig:MalagaSurface}.

The following result makes the link between the generalized staircases and
M\'alaga maps. The definition of first return map can be found in
Chapter~\ref{chap:TranslationFlowsAndIET}.
\begin{tcblemma}{}{MalagaReturn}
Let $M_{\underline{h}}$ be a generalized staircase. Then for each direction $\theta$ that is
not horizontal, the first return map of the translation flow $F^t_{S_{\underline{h}},\theta}$ to the union of the horizontal sides of the rectangles $R_n$
 is canonically identified with the M\'alaga map
$T_{\underline{\alpha}}$ with angles $\alpha_n = h_n \tan(\theta) + 1/2 \mod 1$.
\end{tcblemma}
It is interesting to notice that a fixed generalized staircase $M_{\underline{h}}$ gives rise to a one-parameter
family of M\'alaga maps. A motivated reader can jump to Section~\ref{ssec:ReturnMapMalaga}, where we give several hints on how to proof this Lemma.

In the exercise below we just invite the reader to prove elementary properties of the
M\'alaga map $T_{\underline{\alpha}}$ and the generalized staircase $M_{\underline{h}}$.
\begin{tcbexercise}{}{MalagaAndStaircases}
\begin{enumerate}
\item Under which condition does the surface $M_{\underline{h}}$ have finite area?
\item Show that if $\mu$ is a Borel measure on $(0,+\infty)$ then for $\mu^{\otimes \Z}$-almost every sequence $\underline{h}$ the surface $M_{\underline{h}}$ has infinite area.
\item Show that if $\underline{\alpha}$ is an infinite sequence whose values are $1/3$ and $2/3$ (in any order) then in each interval $I_n$ at least $2/3$ of the points has a divergent orbit. (\textit{hint: the partition in intervals of the form $(k/6, (k+1)/6) \times \{n\}$ is preserved by $T_{\underline\alpha}$. One just needs to study the associated map on the discrete set $\{0,1,2,3,4,5\} \times \Z$}).
\end{enumerate}
\end{tcbexercise}

Theorem~\ref{thm:MalagaRecurrence} implies the following result for $M_{\underline{h}}$.
\begin{tcbcorollary}{}{MalagaRecurrence}
Let $\mu$ be a probability measure on $(0,+\infty)$ whose support contains a half-line. Then for $\mu^{\otimes \Z}$-almost every sequence $\underline{h}$ the translation flow in $M_{\underline{h}}$ is conservative in all directions.
\end{tcbcorollary}
Note that contrarily to the infinite staircase from Section~\ref{ssec:InfiniteStaircase}, the randomized version given in the above corollary is recurrent in \emph{all} directions!
However, the ergodic result of Theorem~\ref{thm:MalagaMinimalErgodic} does not apply to a surface $M_{\underline{h}}$ in almost every direction. The following two questions are
wide open.
\begin{tcbquestion}{}{MalagaMapErgodicity1}
Let $\mu$ be the Lebesgue measure on $[0,1)$. Is $T_{\underline\alpha}$ ergodic for $\mu^\Z$-almost every $\underline\alpha$?
\end{tcbquestion}

\begin{tcbquestion}{}{MalagaMapErgodicity2}
Let $\mu$ be the measure supported on $\{1,2\}$ with $\mu(\{1\}) = \mu(\{2\}) = 1/2$. Is the straightline flow of $M_{\underline{h}}$ ergodic in almost every directions for $\mu^\Z$-almost every $\underline{h}$?
\end{tcbquestion}

\emph{Random wind-tree models}. We finish this section by introducing random variations of the periodic wind-tree model introduced in Section~\ref{ssec:Intro:WindTree}.

Given $\underline\omega \in \{0,1\}^{\Z^2}$ we define a billiard on a generalized polygon $T(\underline\omega)$ as follows. For each position $(m,n)\in\Z^2$ so that $\omega_{m,n}=1$ we place a square obstacle with side length $1/2$ and consider the billiard in the complement of the obstacles. For the constant sequence $\omega_{m,n} = 1$ the generalized polygon $T(\underline\omega)$ is the usual wind-tree model table $T_{1/2,1/2}$ with dimensions $a=b=\frac{1}{2}$ (see Section~\ref{ssec:Intro:WindTree} for the geometric meaning of the parameters $a$ and $b$). While for $\omega_{m,n}=0$ we recover the plane $(\C,dz)$. As for the windtree, one can consider the 4-fold cover giving rise to a translation surface $W(\underline\omega)$.

The only sensible known result about this model is the existence of configurations in which the translation flow is conservative.
\begin{tcbtheorem}{\cite{Troubetzkoy-typical_recurrence_ehrenfest}}{}
The set of parameters $\underline\omega$ for which the translation flow on the translation surface $W(\underline\omega)$ is conservative in almost every direction contains a dense $G_\delta$ set.
\end{tcbtheorem}
The tools used to prove this theorem are the same as the ones used to prove Theorem~\ref{thm:MalagaRecurrence} and Corollary~\ref{cor:MalagaRecurrence}. However the conclusion is much weaker because in order to create barriers in this 2-dimensional coverings one needs a 1-dimensional crown of obstacles.

\begin{figure}[!ht]
\begin{center}\includegraphics{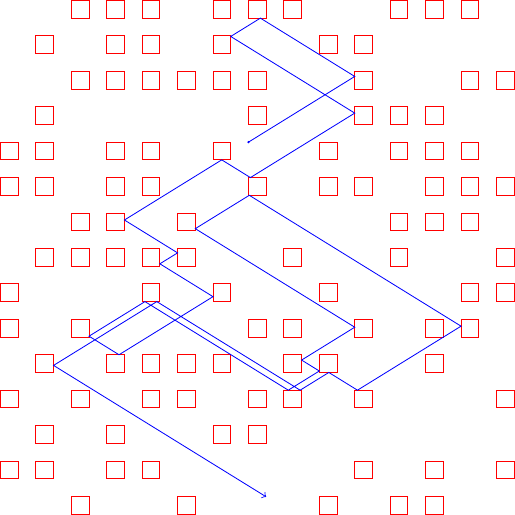}\end{center}
\caption{A trajectory in a Bernoulli wind-tree model with $p=1/2$.}
\label{fig:BernoulliWindtree}
\end{figure}

We now describe the Bernoulli wind-tree model (see Figure~\ref{fig:BernoulliWindtree}). For any given parameter $0 < p < 1$ we consider a sequence $\underline\omega \in \{0,1\}^{\Z^2}$ where each $\omega_{m,n}$ is chosen according to an independent Bernoulli random variable with parameter $p$. In other words, for each point $(m,n)$ in the lattice $\Z^2$ toss a coin which lands on your favourite side with probability $p$. If for $(m,n)$ the coin lands on your favourite side place an obstacle at position $(m,n)$. We denote this product of Bernoulli measures on $\{0,1\}^{\Z^2}$ by $\mu_p$.

\begin{tcbquestion}{}{RandomWindtreeDynamics}
For a given parameter $p$, what are the dynamical properties of $W(\underline\omega)$ where $\underline\omega$ is taken accordingly to the product of Bernoulli measures $\mu_p$?
\end{tcbquestion}

Recall that the original Ehrenfests' wind-tree model concerns obstacles that are not displaced along
a lattice. The natural way to choose at random a countable set of points in the plane is by
a Poisson point process\footnote{The \emph{Poisson point process with density $\delta$} is a way
of choosing randomly a countable set of points in the plane. It is defined by a probability measure
$\mu_\delta$ on countable subsets of $\R^2$ that can be defined by the following property.
For any open set $U$ in $\R^2$ the expected number of points in $U$ is $\delta \area(U)$, in other
words
\[
\int \# (\Gamma \cap U) d\mu_\delta(\Gamma) = \delta \area(U).
\]
The interested reader might want to consult~\cite{CoxIsham} or~\cite{DaleyVereJones}.
}. It is well known that the Poisson point process with density $\delta$ can be obtained
as the limit of the Bernoulli point processes with probability $\delta \epsilon^2$ on the
lattice $\epsilon \Z^2$. Any result relating the dynamics of the Bernoulli wind-tree model
to the Ehrenfests' wind-tree one would be really interesting.

\subsection{Holomorphic foliations}\label{SS:HF}

Every integral curve of a complex vector field can be endowed with a natural translation structure. We illustrate this principle in the following paragraphs by giving explicit examples of complex differential equations whose leaves are, as translation surfaces, precisely those we obtained by unfolding a triangular billiard table in Section~\ref{ssec:PolygonalBilliards}. As a consequence, one can reformulate the problem of finding periodic orbits on a triangular billiard in terms of a real quadratic vector field in $\R^4$.

Consider the complex differential equations in $\C^2$:
\begin{equation}
	\label{eq:ComplexVectorField}
\begin{array}{ccc}
\frac{\partial z_1}{\partial T} & = & \lambda_2 z_1^2+(\lambda_3-\lambda_2)z_1z_2 = P_\lambda(z_1,z_2)\\
\\
\frac{\partial z_2}{\partial T} & = & -\lambda_1z_2^2+(\lambda_3+\lambda_1)z_1z_2 = Q_\lambda(z_1,z_2)
\end{array}
\end{equation}
where $\lambda=(\lambda_1\pi,\lambda_2\pi,\lambda_3\pi)$ are the interior angles of a Euclidean triangle, in particular $\sum\lambda_i=1$. This kind of differential equations arise naturally when studying homogeneous holomorphic foliations on $\C^2$, for more details see \cite{Valdez09_1} and references within. The integral curves of the vector field $V_\lambda=P_\lambda\frac{\partial}{\partial z_{1}}+Q_\lambda\frac{\partial}{\partial z_{2}}$ defined by equation (\ref{eq:ComplexVectorField}) define a (singular) holomorphic foliation $\cF_\lambda$ on $\C^2$.  This foliation has $F(z_1,z_2)=z_1^{\lambda_1}z_2^{\lambda_2}(z_2-z_1)^{\lambda_3}$ as first integral, \ie, the irreducible components of $F^{-1}(w)$, $w\in\C$, in $\C^2\setminus\{(0,0)\}$ are precisely the leaves of $\cF_\lambda$ and the origin is the only singular leaf. Moreover, $\cF_\lambda$ is invariant by the action of the homothety group $\C^*$ and hence every two leaves in $\C^2\setminus\{(0,0)\cup F^{-1}(0)\}$ are diffeomorphic. We refer to any leaf of $\cF_\lambda$ in $\C^2\setminus\{(0,0)\cup F^{-1}(0)\}$ a \emph{generic leaf}.

Let us now endow a generic leaf $\cL\in\cF_\lambda$ with a translation structure. Let $V_{\lambda |\cL}$ be the restriction of the complex vector field $V_\lambda$ to its integral curve $\cL$. Consider the holomorphic 1--form $\eta$ on $\cL$ satisfying the equation $\eta(V_{\lambda |\cL})=1$ (this form is unique because the leaves have complex dimension 1). The pair $(\cL,\eta)$ is a translation surface. The following result relates the foliation $\cF_\lambda$ to triangular billiards introduced in Example~\ref{ssec:PolygonalBilliards}.
\begin{tcbtheorem}{\cite{Valdez09_1}}{FoliationBilliadsLeaves}
Let $(\lambda_1\pi,\lambda_2\pi,\lambda_3\pi)$ be the interior angles of a Euclidean triangle $P$ and $\cL\in\cF_\lambda$ be a generic leaf endowed with its natural translation structure. Then $(\cL,\eta)$ is isomorphic, as a translation surface, to the translation surface $M(P)$ obtained by performing the unfolding construction on P.
\end{tcbtheorem}
We will sketch a proof of this theorem in Section~\ref{sec:CoveringSpaces}.

This correspondence between integral curves of $V_\lambda$ and translation surfaces defined by triangular billiards can be extended to the level of billiard dynamics. Indeed, let $\cF_{\lambda,0}$ be the real foliation on $\C^2$ defined by the integral curves of the real analytic vector field $\Re(V_\lambda)$. This vector field is quadratic and homogeneous and the restriction of $\cF_{\lambda,0|\cL}$ is a real foliation of dimension 1 on $\cL$. For every generic leaf $\cL\in\cF_\lambda$ let $re^{i\theta}\cL$ denote the image of the restriction to $\cL$ of the homothety $T_{re^{i\theta}}(z_1,z_2):=re^{i\theta}(z_1,z_2)$.
\begin{tcbcorollary}{\cite{Valdez09_1}}{}
Let $(\lambda_1\pi,\lambda_2\pi,\lambda_3\pi)$ be the interior angles of a Euclidean triangle $P$. For every generic leaf $\cL\in\cL_\lambda$ there exists a direction $\theta\in\SS^1$ such that the real foliation $\cF_{\lambda,0|\cL}$ is analytically conjugated\footnote{That is, there exists a biholomorphism $\mathcal{L}\to M(P)$ which sends leaves of $\cF_{\lambda,0|\cL}$ to the trajectories of the flow $F_\theta^t$.} to the foliation on $M(P)$ defined by the translation flow $F^t_\theta$. Moreover, for every $\theta'\in\SS^1$,
$\cF_{\lambda,0|re^{i\theta'}\cL}$ is analytically conjugated to the foliation on $M(P)$ defined by the translation flow $F^t_{\theta+\theta'}$.
\end{tcbcorollary}
This corollary permits us to reformulate the open question about the existence of periodic trajectories in a polygonal billiard in terms of a homogeneous quadratic vector field (see Question~\ref{CONJ:PerBill}).
\begin{tcbquestion}{}{}
Does the homogeneous quadratic vector field in $\R^4$:
\begin{equation}
\label{eq:VectFieldBill}
\begin{array}{ccccc}
2\Re(V_\lambda) & = & [\lambda_2(x_1^2-y_1^2)-(\lambda_2+\lambda_3)(x_1x_2-y_1y_2)]\partial/\partial x_1 &+&\\
\\
& & [2\lambda_2x_1y_1-(\lambda_2+\lambda_3)(x_1y_2+x_2y_1)]\partial/\partial y_1 &+&\\
\\
& & [\lambda_1(x_2^2-y_2^2)-(\lambda_1+\lambda_3)(x_1x_2-y_1y_2)]\partial/\partial x_2 &+&\\
\\
& & [2\lambda_1x_2y_2-(\lambda_1+\lambda_3)(x_1y_2+x_2y_1)]\partial/\partial y_2\\
\\
\end{array}
\end{equation}
have a periodic orbit for any choice of parameters $(\lambda_1\pi,\lambda_2\pi,\lambda_3\pi)$ representing the interior angles of a Euclidean triangle $P$?
\end{tcbquestion}

\begin{tcbremark}{}{}
\label{Remark:Projectivize}
As a careful reader might have noticed, the vector field (\ref{eq:VectFieldBill}) is defined on $\R^4$ whereas the phase space of a triangular billiard has dimension 3. Since this vector field is defined by homogeneous polynomials of degree 2, it is possible to consider its projectivization to either $\mathbb{S}^3$ or $\mathbb{RP}(3)$.  For more details on this topic, we refer to \cite{Valdez09_1}.
\end{tcbremark}



\chapter{From topology to geometry}
\label{ch:TopologyGeometry}

An orientable surface with finitely-generated fundamental group is determined, up to homeomorphism, by two non-negative integers: its genus and the number of punctures it has. In Section~\ref{Sec:Topological-classification-surfaces}, we explain the topological classification of orientable surfaces with infinitely generated fundamental group. In a nutshell, the topological type of the surface is determined by its genus (which can be infinite), its space of ends and those ends that can be accumulated by genus.

Once the topological description is settled, we discuss geometric structures. One way to produce a translation structure on an infinite-type surface is to pull-back a translation structure under infinite covering of a finite-type surface. This is explained in Section~\ref{sec:CoveringSpaces}. Many of the surfaces encountered in the introduction are indeed of this form : unfolding of irrational polygonal billiards from Section~\ref{ssec:PolygonalBilliards}, the infinite staircase from Section~\ref{ssec:InfiniteStaircase}, wind-tree models from Section~\ref{ssec:WindTreeIntro}, Panov planes from Section~\ref{ssec:PanovPlanes} as well as some of the Hooper-Thurston-Veech constructions of Section~\ref{ssec:ThurstonVeechConstructionsIntro}. The general question of the existence of a translation structure on a given topologial surface is addressed in Section~\ref{Sec:Existence-translation-structure}. We even provide constructions for any given conformal structure.

In the introductory chapter we saw that infinite-type translation surfaces might present singularities that have no counterpart in the realm of finite type translation surfaces. Their analysis is the subject of Section~\ref{sec:Singularities}. The central object is the \emph{space of linear approaches} which is an extension of the tangent bundle to the metric completion of the surface. To illustrate this section, we introduce the icicled surface constructed by A.~Randecker.

\section{Topological classification of surfaces}
	\label{Sec:Topological-classification-surfaces}

The topological classification theorem for orientable surfaces with\linebreak finitely-generated fundamental group and empty boundary states that any such surface is determined up to homeomorphism by a pair of non-negative integers $(g,n)$ corresponding to the genus and the number of punctures, respectively \cite{Fulton95}. The topological classification of orientable surfaces with infinitely-generated fundamental group and empty boundary was done by I. Richards~\cite{Richards63} and is based on the work of Raymond and Freudenthal's theory of \emph{ends} of topological spaces. Roughly speaking, to a sufficiently nice non-compact topological space $X$ one can associate a topological invariant called the space of ends which encodes the different forms on which a point in $X$ can escape to infinity. In the particular case when $X$ is a surface, there are two kinds of ends: those that are accumulated by genus
and those which have genus zero (a.k.a. \emph{planar ends}). Ends together with genus define the topological invariants needed to classify orientable surfaces with infinitely-generated fundamental group and empty boundary.

\subsection{Ends of topological spaces}
\label{ssec:EndsTopologicalSpaces}

For the purposes of this text we use two definitions of the space of ends of a topological space, which in our context are equivalent. One uses nested sequences of open sets and the other proper rays. The space of ends was originally introduced by H. Freudenthal \cite{Freudenthal31}. For a modern approach (in English) and a more detailed discussion about these two definitios we refer to \cite{Diestel03}.
\begin{tcbdefinition}{\cite{Freudenthal31}}{EndsFreudenthal}
Let $X$ be a locally-compact, locally-connected, connected-Hausdorff space and $U_{1}\supseteq U_{2}\supseteq\ldots$ be an infinite sequence of non-empty connected open subsets of $X$ such that for each $i\in\N$ the boundary $\partial U_{i}$ of $U_{i}$ is compact and $\bigcap\limits_{i\in\N}\overline{U_{i}}=\emptyset$. Two such sequences $U_{1}\supseteq U_{2}\supseteq\ldots$ and $U'_{1}\supseteq U'_{2}\supseteq\ldots$ are equivalent if for every $i\in\N$ there exist $j$ such that $U_{i}\supseteq U'_{j}$ and conversely, that is, for every $i\in\N$ there exist $j$ such that $U'_{i}\supseteq U_{j}$. The corresponding equivalence classes are called \emphdef[end (topological space)]{ends} of $X$ and we denote by $\Ends(X)$ the set of all ends of $X$.
\end{tcbdefinition}
For every non-empty open subset $U$ of $X$ with compact boundary let
\begin{equation}
\label{E:BaseTopEnds}
U^{*}:=\{[U_{1}\supseteq U_{2}\supseteq\ldots]\in\Ends(X)\hspace{1mm}|\hspace{1mm}U_{j}\subset U\hspace{1mm}\text{for some j}\}.
\end{equation}
We endow $\Ends(X)$ with the topology generated by all sets $U^*$. The collection formed by all open sets of $X$ and all sets of the form $U\cup U^{*}$ is a base for a topology of $X':=X\cup\Ends(X)$.

\begin{tcbtheorem}{\cite{Raymond60}}{TopologicalSpaceEnds}
	\label{TH:TOPOLENDS}
Let $X'=X\cup\Ends(X)$ be the topological space defined above. Then,
\begin{compactenum}
\item $X'$ is Hausdorff, connected and locally-connected.
\item $\Ends(X)$ is closed and has no interior points in $X'$.
\item $\Ends(X)$ is totally-disconnected.
\item $X'$ is compact.
\item If $V$ is any open connected set in $X'$, then the set $V\setminus \Ends(X)$ remains connected.
\end{compactenum}
\end{tcbtheorem}
We stress that in his work Raymond defines the space of ends as an inverse limit obtained by considering complements of nested sequences of relatively-compact subsets of $X$. His approach is equivalent to the one presented here. In summary, the space $X'$ is a compactification of $X$ obtained by adding a point for each way one can escape to infinity in $X$. Even though Theorem
\ref{TH:TOPOLENDS} works for a large class of spaces, we will be using it in two simple contexts: surfaces and infinite (but locally-finite) graphs.
\begin{tcbdefinition}{}{EndsProperRays}
Let $X$ be a path-connected, locally-compact space. A proper continuous map $r:[0,\infty)\to X$ is called \emphdef[proper ray (topological space)]{proper ray} in $X$. Two proper rays $r_1,r_2$ are equivalent if for every compact $K\subset X$ there exists an $n\in\N$ such that $r_1([n,\infty))$ and $r_2([n,\infty))$ are contained in the same path-connected component of $X\setminus K$. The equivalence class of a proper ray $r$ is denoted by $\endray(r)$ and the set of equivalence classes $\{\endray(r):\ \text{$r$ is a proper ray in $X$}\}$
is also refered to as the space of ends of $X$.
\end{tcbdefinition}
A sequence of ends $\endray(r_n)$ converges to $\endray(r)$ if for every compact $K\subset X$ there exists a sequence of integers $M_n$ such that $r_n[M_n,\infty)$ and $r[M_n,\infty)$ lie in the same path component of $X\setminus K$ whenever $n$ is sufficiently large.


\begin{tcbexercise}{}{}
Prove that definitions \ref{def:EndsFreudenthal} and \ref{def:EndsProperRays} lead to homeomorphic spaces of ends when $X$ is a locally-compact path-connected Hausdorff space.
\end{tcbexercise}


\begin{tcbexercise}{}{}
Let $\beta_n=\{e_i\}_{i=1}^n$ be the standard basis of $\Z^n$. For example, $\beta_1=\{1\}$ and $\beta_2=\{(1,0), (0,1)\}$. Prove that the Cayley graph of $\Z^n$ w.r.t. $\beta_n$ has one end if $n\geq 2$ and two ends if $n=1$. Denote by $F_n$ the free group generated by $n\geq 2$ elements $\{a_1,\ldots,a_n\}$ and $\Gamma_n$ the corresponding Cayley graph. Prove that all spaces $\Ends(F_n)$, $n\geq 2$ are homeomorphic to the Cantor set.
\end{tcbexercise}


\subsection{Topological classification of orientable surfaces}
\label{ssec:TopologicalClassificationOfOrientableSurfaces}
By \emphdef{surface} we mean a connected topological manifold $S$ of real dimension two. We stress that in this text all surfaces are required to be second-countable topological spaces, \ie, to have a countable basis for their topology. In particular surfaces  are metrizable. Henceforth all the surfaces considered in this text are orientable and, unless explicitly stated, will have empty boundary. A simple closed \emph{curve} in $S$ is a continuous injective map $\alpha:\mathbb{S}^1\to S$, where $\mathbb{S}^1\subset\C$ is the unit circle. A curve is said to be \emph{essential} if it is not isotopic to the boundary curve of a neighbourhood of a puncture of $S$ nor to a point. We often abuse language and use the term curve to refer to the map $\alpha$, its image and its isotopy class in $S$. A simple curve $\alpha$  is \emph{separating} if $S\setminus \alpha$ has two connected components.

\begin{tcbdefinition}{}{GenusZeroInfiniteGenus}
	\label{DEF:GenusZeroInfiniteGenus}
A set of disjoint isotopy classes of curves $\{a_i\}_{i\in I}\subset S$ is called multicurve. The genus of $S$ is the maximal cardinality of
a multicurve $\{a_i\}_{i\in I}$ for which there exist representatives $\{\alpha_i\}_{i \in I}$ such that $S\setminus\cup_{i\in I}\alpha_i$ is connected. If such cardinality is infinite we say that $S$ has infinite genus.
 An end $[U_1\supseteq U_2\supset\ldots]\in\Ends(S) $ is called planar if there exists an $i\in\N$ such that $U_i$ has genus zero. We define $\Ends_\infty(S)\subset\Ends(S)$ as the subspace of all ends which are not planar, and we refer to them as \emph{ends accumulated by genus} or \emph{ends of infinite genus}.
\end{tcbdefinition}




Note that a surface $S$ has genus zero if every simple essential curve in $S$ is a separating curve. We have chosen the notation $\Ends_\infty(S)$, for if $[U_1\supseteq U_2\supset\ldots]\in\Ends_\infty(S)$
then every $U_i$ is an infinite genus surface. In other words, non-planar ends always have infinite genus.
It follows from the definitions that  $\Ends_\infty(S)$ forms a closed subset of $\Ends(S)$.

\begin{tcbtheorem}{Topological classification of orientable surfaces}{}
Two orientable surfaces $S$ and $S'$ with empty boundary are homeomorphic if and only if they have the same genus $g\in\{\Z_{\geq 0}\cup\infty\}$, and both $\Ends_\infty(S)\subset\Ends(S)$ and $\Ends_\infty(S')\subset\Ends(S')$ are homeomorphic as nested topological spaces, that is, there exists a homeomorphism $h:\Ends(S)\to\Ends(S)'$ whose restriction $h_|:\Ends_\infty(S)\to \Ends_\infty(S')$ is a homeomorphism as well.
\end{tcbtheorem}
This result was first announced by Ker\'ekj\'art\'o, but a complete proof was given by Richards \cite{Richards63}. It tells us that the complete topological invariant of an arbitrary orientable topological surfaces is the genus plus a couple of nested topological spaces characterising the ways a point in $S$ has to escape to infinity and, among these ways, which of them carry infinite genus. Note that if $S$ is a surface of genus $g$ having $n$ punctures, then $\Ends_\infty(S)=\emptyset$ and $\Ends(S)=\{1,\ldots,n\}$ endowed with the discrete topology.

\begin{figure}[!ht]
\begin{center}
\begin{minipage}{.35\textwidth}%
\begin{center}%
\includegraphics[scale=0.3]{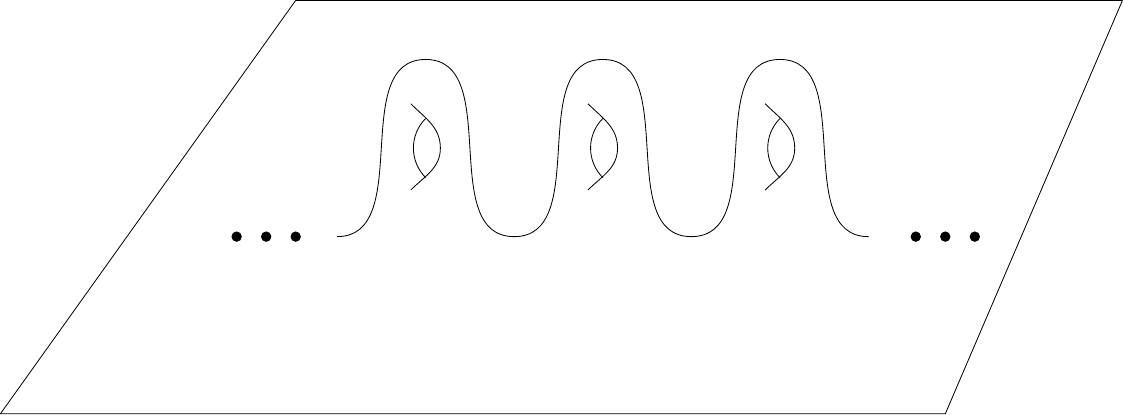}
\end{center}
\subcaption{The Loch Ness Monster}
\label{fig:LochNessMonster}
\end{minipage}
\hspace{.05\textwidth}
\begin{minipage}{.45\textwidth}%
\begin{center}%
\includegraphics[scale=0.35]{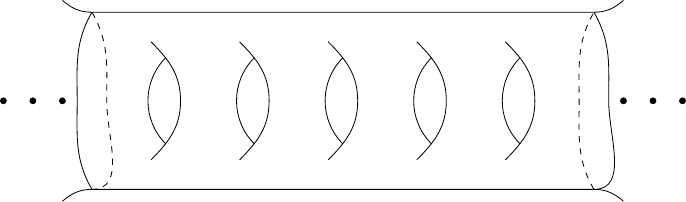}
\end{center}
\subcaption{Jacob's Ladder}
\label{fig:JacobsLadder}
\end{minipage}
\hspace{.05\textwidth}
\begin{minipage}{.45\textwidth}%
\begin{center}%
\includegraphics[scale=0.6]{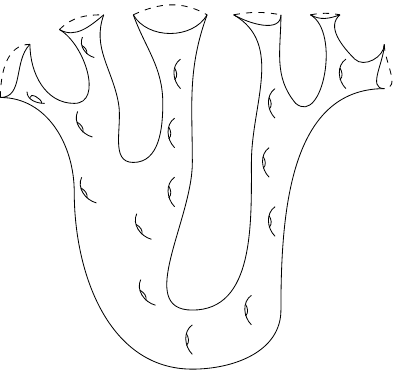}
\end{center}
\subcaption{Blooming Cantor Tree.}
\label{fig:BloomingCantorTree}
\end{minipage}
\hspace{.05\textwidth}
\begin{minipage}{.45\textwidth}%
\begin{center}%
\includegraphics[scale=0.6]{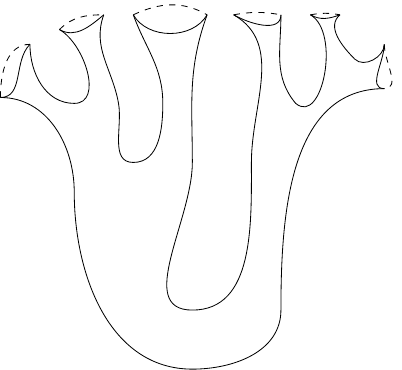}
\end{center}
\subcaption{Cantor Tree.}
\label{fig:CantorTree}
\end{minipage}\bigskip \\
\end{center}
\caption{Famous infinite-type surfaces}
\label{fig:FamousInfiniteTypeSurfaces}
\end{figure}

From Theorem \ref{TH:TOPOLENDS} we deduce that both $\Ends_\infty(S)$ and $\Ends(S)$ are closed, totally-disconnected and Hausdorff. Hence the couple $\Ends_\infty(S)\subset \Ends(S)$ is homeomorphic to a couple of nested closed subsets of $C'\subset C$ of the standard Cantor set obtained by removing middle thirds from an interval. On the other hand, the following result tells us that every such pair defines a homeomorphism class of orientable surfaces.
\begin{tcbtheorem}{\cite{Richards63}}{}
Let $C'\subset C$ be a nested pair of closed subsets of the Cantor set. Then there exist a surface $S$ such that $\Ends_\infty(S)\subset\Ends(S)$ is homeomorphic to $C'\subset C$ as nested pair of topological spaces.
\end{tcbtheorem}
\begin{proof}[Sketch of proof]
Consider $C$ as a subspace of points in the Sphere $\mathbb{S}^{2}$. Then $\Ends(S=\mathbb{S}^{2}\setminus C)$ is homeomorphic to $C$. An easy way to see this is to imagine that points in $C$ are pulled away to infinity. If $U_{1}\supseteq U_{2}\supseteq\ldots$ represents a point $c'\in C'$ then, by gluing a suficiently small 2-dimensional torus into a sufficiently small neighborhood of a point $p\in U_{j}\setminus U_{j+1}$ for every $j>>1$, we produce the desired surface $S$.
\end{proof}
From the preceding theorem we deduce that there are uncountably many different topological types of surfaces whose fundamental group is not finitely-generated. Some of these surfaces seem to appear more naturally than others, and hence have names of their own. The following nomenclature can be attributed to A. Phillips and D. Sullivan~\cite{PhillipsSullivan81} or E. Ghys~\cite{Ghys95}.
\begin{tcbdefinition}{}{MonsterLadder}
An infinite-genus surface $S$ with only one end is called a \emph{Loch Ness monster}. An infinite-genus surface with two ends, each accumulated by genus, is called a \emph{Jacob's Ladder}. If a surface $S$ has
only planar ends and $\Ends(S)$ is homeomorphic to the Cantor set, then $S$ is called a \emph{Cantor tree}. On the other hand, if $S$ has no planar ends and $\Ends(S)$ is homeomorphic to the Cantor set, then $S$ is called a \emph{blooming Cantor tree}. These are illustrated in figure \ref{fig:FamousInfiniteTypeSurfaces}.
\end{tcbdefinition}

\begin{figure}[!ht]
\begin{center}
\includegraphics[scale=1.5]{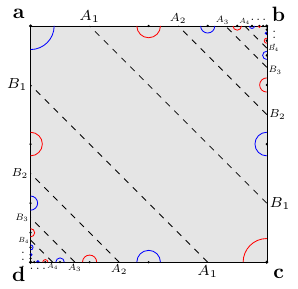}
\end{center}
\caption{$B_{\frac{1}{2}}\setminus\{c_1,\ldots,c_n\}$ is connected for any $n\in\N$, hence it has infinite genus.}
\label{fig:chamanarainfinitegenus}
\end{figure}

\begin{tcbexample}{}{BakerIsAMonster}
Recall that in Example~\ref{ssec:BakerBowmanArnouxYoccoz} we defined for every $\alpha\in(0,1)$ a surface $B_\alpha$ called baker's surface. We claim that for every parameter, $B_\alpha$ is a Loch Ness monster. It is sufficient to show this is true for $\alpha=\frac{1}{2}$, case in which $B_{\frac{1}{2}}$ has area 1. Consider figure \ref{fig:chamanarainfinitegenus} and for each $n\in\N$ the simple closed curve $c_n$ formed by the union of the segments $I_n$ and $J_n$ of slope -1 given by: $I_n$ joins the midpoint of $A_n$ in the lower side of the square to the midpoint of $B_n$ in the right side, and $J_n$ joint the midpoint of $B_n$ in the left side to the midpoint of $A_n$ in the top. All curves in $\{c_n\}_{n\in\N}$ are non-isotopic. Moreover, from the figure one can see that $B_{\frac{1}{2}}\setminus\{c_1,\ldots,c_n\}$ is connected for any $n\in\N$, therefore $B_{\frac{1}{2}}$ has infinite genus. To show that $B_{\frac{1}{2}}$ is one-ended we use the following lemma, whose proof is left as an exercise to the reader:
\begin{tcblemma}{}{CriterionInfiniteGenus}
Let $S$ be any surface. The space of ends $\Ends(S)$ has only one element if and only if for every compact $K\subset S$ there exists a compact subset $K'$ of $S$ containing $K$ and such that $S\setminus K'$ is connected.
\end{tcblemma}
Every compact subset $K\in B_{\frac{1}{2}}$ is contained in the complement of a small neighbourhood $U$ of the points in the square to which the translation surface structure cannot be extended. These points are the extrema of the segments labeled by $A_n$, $B_n$, $n\in\N$ in the figure plus the corners $a,b,c$ and $d$. If we denote by $K'$ this complement then by definition it is a compact set. Any point in $B_{\frac{1}{2}}\setminus K'$ can be joined by an arc to the extremity of one segment in $\{A_n,B_n\}_{n\in\N}$. Given that the extremities of these segments accumulate to the corners $b$ and $d$ we conclude that one can connected any two points in the complement of $K'$ through an arc.
\end{tcbexample}

Very similar arguments are used in the next section to prove that the infinite staircase and any infinite-type surface coming from the billiard in generic triangle are homeomorphic to the Loch-Ness monster.

\begin{tcbexercise}{}{}
Let $\eps>0$ be small (say less than $\frac{1}{10}$) and consider the orbit in $\R^2$ of the square $[0,\eps]^2$ by the translation group $(x,y)\to(x+n,y+m)$,   $(n,m)\in\Z\times\Z$. Consider the infinite polygon given by the complement of this orbit, and the translation surface $M$ obtained from identifying pairs of opposite sides in each square. Prove that $M$ is a Loch-Ness monster. Use the translation $(x,y)\to(x,y+1)$ to produce a quotient of $M$ homeomorphic to Jacob's ladder.
\end{tcbexercise}

\begin{tcbexercise}{}{}
By gluing polygons (see Definition~\ref{def:TranslationSurfaceConstructive} in the preceeding chapter), construct a translation surface homeomorphic to (A) the Cantor tree and (B) the Cantor blooming tree.
\end{tcbexercise}

\section{Covering spaces}
\label{sec:CoveringSpaces}



Many infinite-type translation surfaces in this book are infinite coverings of
finite-type translation surfaces. Let us mention the infinite staircase
(see Section~\ref{ssec:InfiniteStaircase}) and the wind-tree model
(see Section~\ref{ssec:WindTreeIntro}). As we see in the next
Chapters, many geometrical and dynamical properties of such infinite
coverings can be deduced from properties of the finite-type
base surface.


\begin{tcbdefinition}{}{TranslationCovering}
Let $M$ be a translation surface,
$\Sigma\subset M$ a (possibly empty) discrete subset of $M$ and
define $M^0:=M\setminus\Sigma$. Consider a non-ramified covering map
$p^0:\widetilde{M}^0\to M^0$ defined by a subgroup $\Gamma < \pi_1(M^0)$.
We endow the covering surface $\widetilde{M}^0$ with the translation structure
defined by the pullback of the translation structure of $M^0$ via $p^0$. Let
$\widehat{p}: \widehat{\widetilde{M}} \to \widehat{M}$
be the continuation of $p^0$ to the metric completions and $p:\widetilde{M} \to M$ be the
restriction of $\widehat{p}$ to $\widetilde{M} \subset \widehat{\widetilde{M}}$ that is
obtained by adding to $\widetilde{M}^0$ the conical singularities and regular points
in $\widehat{\widetilde{M}}$ that belong to preimages $\widehat{p}^{-1}(s)$ of points
$s$ in $\Sigma$. We call $p:\widetilde{M}\to M$ the \emphdef{translation covering}
of $M$ defined by the subgroup $\Gamma$.
\end{tcbdefinition}
It is important to notice that for a translation covering $p:\widetilde{M}\to
M$ the map $p$ is not necessarily surjective. However, it
extends continuously in a unique way to the surjective map
$\widehat{p}: \widehat{\widetilde{M}}\to\widehat{M}$.
The set of points in $M$ over which $\widehat{p}$ fails to be a covering
are called \emphdef[branching points (translation covering)]{branching points} and their union the \emphdef[branching locus (translation covering)]{branching locus}.
Their preimages in $\widehat{M}$ are called
\emphdef[ramification points (translation covering)]{ramification points}.

For the infinite staircase, illustrated in Figure~\ref{fig:staircase1}, the
base surface $M$ is a torus tiled by two squares and $\Sigma$ is formed by the
image of the vertices of these two squares. In this example
$\widetilde{M}^0=\widetilde{M}$ because no point in $\widehat{\widetilde{M}} \setminus \widehat{M}$
is either regular or a conical singularity. This is an example where
the translation covering $p:\widetilde{M}\to M$ is not surjective.

By a slight abuse of language we apply the standard terminology of
the regular covering $p^0: \widetilde{M^0} \to M^0$ to the translation
covering $p: \widetilde{M} \to M$. We call \emphdef[degree (translation covering)]{degree}
of $p$ the degree of $p^0$, we say that $p$ is \emphdef[normal (translation covering)]{normal}
or \emphdef[regular (translation covering)]{regular} if $p^0$ is normal, etc.

By definition all points in $\widetilde{M}^0$ are regular and
the deck transformation group $\deck(p)$ acts on $\widetilde{M}^0$ by
translations (see Definition \ref{def:TypologyAffineMaps}). We say that the
translation covering is \emph{finite} if the fibers of the covering map are
finite.
\begin{tcbexercise}{}{}
Prove that if $p:\widetilde{M}\to M$ is a finite translation covering over a finite-type translation surface $M$, then $\widetilde{M}$ is always a translation surface.
\end{tcbexercise}

\emph{Outline of this section}. We dedicate the rest of this section to discuss some relevant examples of translation coverings and related covering spaces. First we put the emphasis on normal coverings defined by infinite abelian groups. The second series of examples focus the attention on the geometry of the base surface, more precisely when $M$ is a torus and $\Sigma$ a singleton. These are the so called squared-tiled surfaces or origamis. Then we deal with coverings arising from irrational polygonal billiards. These fail to be translation coverings in the sense of Definition~\ref{def:TranslationCovering} because the base is not a translation surface (though it has a Euclidean structure) and the deck transformation group does not act by translations. Finally, we discuss dilation surfaces. This class of surfaces, as we see later, naturally appears when stuying the dynamics of the translation flow on infinite-type translation surfaces.


\textbf{G-coverings}. Let $p:\widetilde{M}\to M$ be a translation covering. Recall that the fibers of the covering map $p:\widetilde{M}^0 \to M^0$ are naturally identified with the cosets of the subgroup $p_*\pi_1(\widetilde{M}^0)<\pi_1(M^0)$ and that any covering of $M^0$ is determined (up to covering isomorphism) by  (the conjugacy class of) a subgroup of $\pi_1(M^0)$.  Given that the induced map $p_*:\pi_1(\widetilde{M}^0)\to\pi_1(M^0)$ is always injective, we will also write $\pi_1(\widetilde{M}^0)<\pi_1(M^0)$ when there is no ambiguity.

Most of the translation coverings we study in this text are normal coverings, that is to say $\pi_1(\widetilde{M}^0)$ is a normal subgroup of $\pi_1(M^0)$, hence the deck transformation group $\deck(p)$ of the translation covering can be identified with $G = \pi_1(M^0) / \pi_1(\widetilde{M}^0)$ and acts transitively on the fibers of the covering. \emph{In this situation $p:\widetilde{M}^0\to M^0$ is called a $G$-covering}, because the group $G$ also acts freely and properly discontinously on $\widetilde{M}^0$, $M^0=\widetilde{M}^0/G$ and $p$ is just the projection $\widetilde{M}^0\to\widetilde{M}^0/G$. For more details on $G$-coverings we refer the reader to the book of W.~Fulton~\cite{Fulton95}.


\begin{tcbexample}{$\Z$-coverings and the infinite staircase}{Zcoverings}
Let $M$ and $\Sigma$ as before. A $\Z$-covering of $M^0=M\setminus\Sigma$ is given by the kernel of a surjective morphism  $f: \pi_1(M^0) \to \Z$.
Since $\Z$ is Abelian, $f$ factorizes through a map $\phi: H_1(M^0; \Z) \to \Z$, \ie, $\phi \in H^1(M^0;
\Z)$. On the other hand, the algebraic intersection form
\[
\langle\cdot,\cdot\rangle: H_1(M^0; \Z) \times H_1(M , \Sigma; \Z) \to \Z
\]
is non-degenerate and hence induces a bijection $H^1(M^0; \Z) \simeq H_1(M, \Sigma; \Z)$.
In particular, the map $\phi$ is represented by a unique element $c$ of
$H_1(M , \Sigma; \Z)$, that is for every, $v\in H_1(M^0; \Z) $ we have $\phi(v) = \langle v, c\rangle$. Hence, if we define the submodule
$$
c^\bot:=\{ c'\in H_1(M^0; \Z) \hspace{1mm}|\hspace{1mm} \langle c',c\rangle=0 \}
$$
and $\textbf{ab}:\pi_1(M^0)\to H_1(M^0; \Z)$ denotes the abelianization homomorphism, then the $\Z$-covering in question is given by $\ker(f)=\textbf{ab}^{-1}(c^\bot)$.


\begin{tcbexercise}{}{}
Show that an element $\phi$ of $ H^1(M^0; \Z)$ is surjective if and only if it is represented by a primitive class $c$ in $H_1(M,\Sigma;\Z)$
(\ie, there is no $n \in \Z \setminus \{\pm 1\}$ such that $c = n c'$ with $c' \in H_1(M,\Sigma;\Z)$).
\end{tcbexercise}

The preceding exercise implies that every primitive element $\pm c\in H_1(M,\Sigma;\Z)\simeq H^1(M^0;\Z)$ determines a $\Z$-covering of $M$. In the following lines we explain in detail how the infinite staircase introduced in Section~\ref{ssec:InfiniteStaircase} can be seen as a $\Z$-covering. Consider the rectangle $P$ on the right hand side of Figure~\ref{fig:staircase1}
formed by two unit squares. If we identify using translations parallel sides labeled with the same
letter the result is a flat torus $M$ with two marked
points $\Sigma\subset M$ coming from the corners of the
aforementioned squares. Denote by $\{A,B,C\}$ the oriented basis of $H_1(M,\Sigma;\Z)$ defined
by the identification of the sides $\{A^\pm,B^\pm,C^\pm\}$ respectively (in the figure the orientation of these is illustrated). The
$\Z$-covering $p:\widetilde{M}\to M$ defined by the cycle $c=B-A$ can be concretely constructed from
the polygon $P$ and $c$ as follows. Let $\{\phi_A,\phi_B,\phi_C\}$ and $\{\gamma_A,\gamma_B,\gamma_C\}$
be the basis of $H^1(M^0;\Z)={\rm Hom}(H_1(M^0;\Z),\Z)$ defined by $\{A,B,C\}$ and its dual in $H_1(M^0;\Z)$, respectively. If we denote by $\phi_c$ the cohomology class determined by $c$, then $\widetilde{M}^0$ is given by performing the following identifications (using translations) in the infinite family $P\times\Z$: identify
the side $(A^+,n)$ with $(A^-,n+\phi_c(\gamma_A)) = (A^-,n-1)$, $(B^+,n)$ with $(B^-,n+\phi_c(\gamma_B)) = (B^-,n+1)$
and $(C^+,n)$ with $(C^-, n+\phi_c(\gamma_C)) = (C^-, n)$. The covering space $p:\widetilde{M}^0\to M^0$ is depicted on the left hand
side of Figure~\ref{fig:staircase1} and is exactly the \emph{infinite
staircase} introduced in Section~\ref{ssec:InfiniteStaircase}. It is important to remark that the metric completion of $\widetilde{M}^0$ is not a surface, because in order to obtain it we have to add four infinite angle singularities, none of which has a compact neighbourhood. Therefore in this case the extension of $p:\widetilde{M}^0\to M^0$ to the metric completion is not a map between translation surfaces.

\begin{figure}[H]
\begin{center}
\includegraphics[scale=0.775]{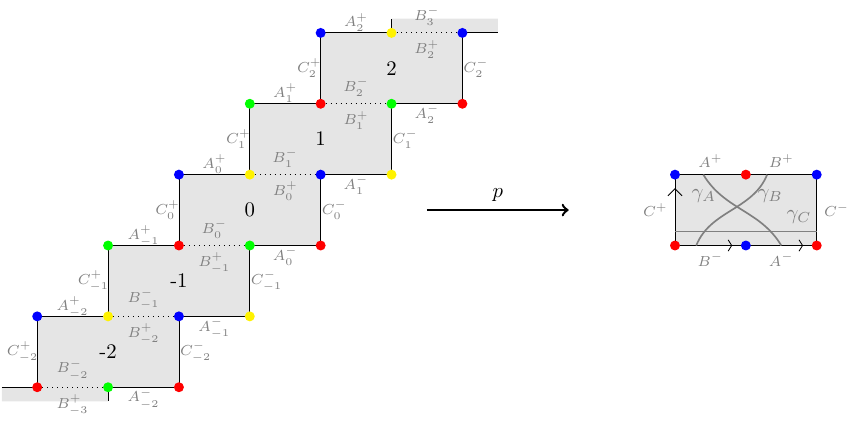}
\end{center}
\caption{The staircase as a $\Z$-covering of a torus.}
\label{fig:staircase1}
\end{figure}
\end{tcbexample}

\begin{tcbexercise}{The infinite staircase is a Loch Ness monster}{InfiniteStaircaseIsLochNess}
Let $p:\widetilde{M}\to M$ be the $\Z$-translation covering given by the infinite staircase as depicted in figure~\ref{fig:staircase1}.
\begin{enumerate}
\item For each $k\in\Z$, let $V_k\subset \widetilde{M}$ be the maximal vertical cylinder containing the horizontal saddle connection defined by the side $A_k^+$ and $\gamma_k$ its the core curve. Show that $S\setminus \bigcup_{k\in\Z}\gamma_k$ is connected and any two curves in $\{\gamma_k\}_{k \in \Z}$ are non-isotopic.
\item If $K\subset \widetilde{M}$ is compact, show that there exist a neighbourhood $U$ of $\Sing(\widetilde{M})$ whose complement contains $K$. Use that neighbourhood to construct a compact set $K'$ containg $K$ such that $\widetilde{M}\setminus K'$ is connected. \emph{Hint}: the set $K$ must be contained in a finite union of steps of the staircase; then proceed as in the case of baker's surface explaine in Example~\ref{exa:BakerIsAMonster}.
\end{enumerate}
Deduce from the preceding points that the infinite staircase is homeomorphic to the Loch Ness monster.
\end{tcbexercise}

\begin{tcbremark}{}{}
It is somehow surprising that the infinite staircase is a Loch Ness monster,
because it actually looks like it has two ends and the deck transformation
group of the covering $p:\widetilde{M}\to M$ it defines is $\Z$, which is a
group with two ends. Hence, by the \v{S}varc--Milnor\footnote{The \v{S}varc--Milnor lemma
says the following: given a discrete group $G$ acting properly-discontinuously and
cocompactly by isometries on a (quasi-geodesic) metric space $X$ then $X$ is quasi-isometric to
$G$. For a more precise statement, see~\cite{Loh17}. Informally, \emph{quasi-isometric} means that at large scale $X$
and $G$ look the same. Formally, a quasi-isometry between two metric spaces
$X$ and $Y$ is a (not-necessarily continuous) function $f: X \to Y$ such
that there exists real constants $A \geq 1$, $B \geq 0$ and $C \geq 0$ such
that the two following properties hold
\[
\begin{array}{l}
\forall x, x' \in X, \quad
\frac{1}{A} d_X(x,x') - B \leq d_Y(f(x), f(x')) \leq A d_X(x, x') + B, \\
\forall y \in Y, \quad
\exists x \in X, d_Y(f(x), y) \leq C.
\end{array}
\]
Two metric spaces $X$ and $Y$ are quasi-isometric if there is a quasi-isometry
from $X$ to $Y$ (this definition is symmetric!).

Given a finitely generate group $G$ with generating set $S$ (such as $\Z$ with
$\{-1, +1\}$), one can associate a metric $d_S$ on $G$ coming from the word
distance with respect to $S$.  It can be viewed as the distance on the Cayley
graph of $G$ with respect to the generating set $S$ (where the edge length is
1). Changing the generating set $S$ does not change the quasi-isometric type of
$G$ and this is why one often say "quasi-isometric to the group $G$". More
formally one should say "quasi-isometric to the group $G$ endowed with any of
its word length metric".

Two quasi-isometric geodesic spaces have the same set of ends. In particular,
if $p:\widetilde{S}\to S$ is a non-ramified normal covering of a compact
surface $S$, then $\Ends(\widetilde{S})$ is homeomorphic to the space of
ends of the Cayley graph of the deck group $\Deck(p)$.
This implies that the deck transformation group imposes restrictions on
the topology of a non-ramified covering. For example, given that a
finitely-generated group has either 0,1,2 or infinitely many ends, there exist
no non-ramified normal covering $\widetilde{S}\to S$  (with $S$ a compact
surface) for which $\Ends(\widetilde{S})$ has $3$ elements. For a more detailed discussion on the coarse geometry of infinite translation surfaces we refer the reader to C. Karg~\cite{Karg2020}.
}
lemma one would expect to have two ends. The catch here is that $\Deck(p)$ does not act cocompactly (because of the existence of infinite
cone angle singularities) and hence the infinite staircase is not quasi-isometric to the Cayley graph of $\Z$.
\end{tcbremark}


\begin{tcbexample}{$\Z^d$-coverings}{ZDCoverings}
	\label{Ex:Zd-coverings}
We now focus our attention on higher rank (torsion-free) translation coverings. Consider the the kernel of a surjective morphism $f:\pi_1(M^0)\to\Z^d$ where
$M^0 = M \setminus \Sigma$. Such morphism factorizes through a map
$\phi:H_1(M^0;\Z)\to\Z^d$, in other words an element of $H^1(M^0; \Z^d)$.
By duality, we can identify $\phi$ as an element of the relative
homology $H_1(M,\Sigma; \Z^d)$ that is a $d$-tuple
$c = (c_1, \ldots, c_d)$ of relative cycles in $H_1(M, \Sigma; \Z)$.
To such $d$-tuple, one can associate the vector subspace
$V = \Q c_1 + \Q c_2 + \ldots + \Q c_d$ of $H_1(M, \Sigma; \Q)$.
The surjectivity of the morphism $f$ implies that $\dim_\Q V = d$.

Conversely, given $V$ a vector subspace of $H_1(M, \Sigma; \Q)$ we reconstruct
a morphism $f: \pi_1(M^0) \to \Z^d$. Namely, we define a normal subgroup
of $\pi_1(M^0)$ as follows
\begin{equation}
\label{eq:FundamentalGroupSubgroupForAbelianCoverings}
\Gamma_V := \{\gamma \in \pi_1(M^0):\ \forall v \in V,\ \langle \gamma,v\rangle = 0\}.
\end{equation}
We claim that the subgroup $\Gamma_V$ of $\pi_1(M^0)$ defines a
$\Z^d$-covering of $M^0$. To see this, let us choose a
$\Z$-basis $(c_1, c_2, \ldots, c_d)$ of $V_\Z := V \cap H_1(M, \Sigma; \Z)$.
We define a morphism $f: \pi_1(M^0) \to \Z^d$ by setting
\[
f(\gamma) = (\langle c_1, \gamma\rangle, \ldots, \langle c_d, \gamma\rangle).
\]
It is clear that $\ker(f) = \Gamma_V$.

We summarize the discussion in the following lemma. Compare with~\cite[Proposition 3.2]{HooperWeiss09} for the case $d=1$.
\begin{tcblemma}{}{AbelianCoveringCorrespondence}
Let $M$ be a translation surface and $\Sigma \subset M$ a discrete subset.
The map that associates to a subspace $V$ of $H_1(M,\Sigma;\Q)$ of dimension $d$
the normal subgroup $\Gamma_V$ defined
in~\eqref{eq:FundamentalGroupSubgroupForAbelianCoverings} is a bijection
between subspaces of dimension $d$ of $H_1(M, \Sigma; \Q)$ and $\Z^d$-coverings
of $M$ ramified at most over $\Sigma$.
\end{tcblemma}

Given the translation surface $M$ and the set $\Sigma \subset M$ there are
several canonical subspaces $V \subset H_1(M, \Sigma; \Q)$ that gives rise to
interesting covers:
\begin{enumerate}
\item $V = H_1(M, \Sigma; \Q)$ defines the \emph{maximal (relative) Abelian cover},
\item If $M=(X,\omega)$, then $V = \ker( \gamma \in H_1(M, \Sigma; \Q) \mapsto \int_\gamma \omega )$ defines the \emph{maximal (absolute) Abelian cover without drift}. This example will reappear in Definition~\ref{def:Kerhol}, Chapter~\ref{ch:Symmetries}.
\end{enumerate}

\begin{figure}[H]
\begin{center}\includegraphics{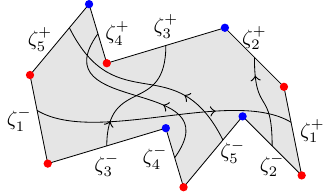}\end{center}
\caption{A polygon $P$ with a side pairing $\zeta_i^- \leftrightarrow \zeta_i^+$
gives rise to a translation surface $M$ with a canonical basis $\{\zeta_i\}_{i=1,\ldots,N}$
of $H_1(M, \Sigma; \Z)$ and a dual basis $(\gamma_1, \ldots, \gamma_N)$ of
$H_1(M \setminus \Sigma; \Z)$.}
\label{fig:PolygonAndHomology}
\end{figure}
We now generalize the construction of the infinite staircase presented above
to arbitrary $\Z^d$-coverings. Let $P$ be a (finite) polygon with
$2N$ sides labeled $\zeta_i^-, \zeta_i^+$ for $i=1,\ldots,N$ and where
$\zeta_i^- = - \zeta_i^+$ as vectors in $\C$ (the orientation of each side
comes from the standard orientation of the plane). Let $M$ be the translation surface
obtained by considering the pairing $f(\zeta_i^+)=\zeta_i^-$ on the sides of $P$, that is, contructed by
identifying $\zeta_i^+$ with $\zeta_i^-$ using a translation (see Definition \ref{def:TranslationSurfaceConstructive}, Chapter \ref{ch:Introduction}). Denote
by $\Sigma \subset M$ the images of the vertices of the polygon $P$.
Let $\zeta_1,\ldots,\zeta_N$ the image of $\zeta_1^+,\ldots,\zeta^+_N$ in $M$.
The set $\{\zeta_1,\ldots,\zeta_N\}$ is a basis of the relative homology $H_1(M,\Sigma;\Z)$.
Let $\{\gamma_1, \ldots, \gamma_N\}$ be the dual basis in $H_1(M \setminus \Sigma; \Z)$.
See also Figure~\ref{fig:PolygonAndHomology}.

Now let $V \subset H_1(M,\Sigma;\Q)$ be a $d$-dimensional subspace as in
Lemma~\ref{lem:AbelianCoveringCorrespondence} and let us denote
$V_\Z = V \cap H_1(M,\Sigma;\Z)$. Let us choose a
$\Z$-basis $\{c_1,\ldots, c_d\}$ of $V$ (that is $c_i \in H_1(M,\Sigma;\Z)$
and $\Z c_1 + \Z c_2 + \cdots + \Z c_d = V_\Z$). By duality, each $c_i$ can be
seen as an element $\phi_{c_i}$ in $H^1(M^0; \Z)$.
We construct the $\Z^d$-covering $p:\widetilde{M}^0\to M^0$ associated to $V$ as follows.
Let $P^0$ be the polygon $P$ with the vertices removed. In $P^0 \times \Z^d$
for each $i=1,\ldots, d$ and $(n_1,\ldots,n_d) \in \Z^d$ the side
$(\zeta_i^+,(n_1,\ldots,n_d))$ in $P^0 \times \{(n_1,\ldots,n_d)\}$ with the side
$(\zeta_i^-, (n_1 + \phi_1(\gamma_1), \ldots, n_d + \phi_d(\gamma_i)))$
in $P^0\times\{(n_1 + \phi_1(\gamma_i),\ldots, n_d + \phi_d(\gamma_i))\}$.


\begin{tcbexercise}{}{}
Let $M$ be a finite-type translation surface obtained from a polygon $P$ as in the construction of $\Z$-covering from Example~\ref{exa:Zcoverings}.
Let $G$ be a group together with a surjective group morphism $f:\pi_1(M\setminus\Sigma)\to G$. Describe in concrete geometrical terms how to construct the $G$-covering $\widetilde{M}\to M$ defined by the kernel of $f$.
\end{tcbexercise}

\end{tcbexample}

\begin{tcbexample}{Square-tiled surfaces}{SquareTiledSurfaces}

The infinite staircase introduced in Section~\ref{ssec:InfiniteStaircase}
belongs to the class of translation surfaces built from
unit squares that we discuss now in more detail.
\begin{tcbdefinition}{}{SquareTiledSurface}
Let $\T^2$ be the flat torus $\C/\Z^2$. A \emphdef{square-tiled surface} is a translation covering $p: M\to \T^2$ branched at most over $0$.
\end{tcbdefinition}

\begin{tcbremark}{}{}
Square-tiled surfaces are also sometimes called \emphdef{origami} in the literature.
\end{tcbremark}

\begin{tcbexercise}{}{}
Show that the following are equivalent for an origami $p: M \to \T^2$
\begin{compactenum}
\item $M$ is compact,
\item $M$ is of finite-type,
\item the degree $\deg(p^0)$ is finite.
\end{compactenum}
\end{tcbexercise}

Since $\pi_1(\T^2 \setminus \{0\}) \simeq F_2$, the free group on two
generators, square-tiled surfaces are in bijection with conjugacy classes of
subgroups of $F_2$ (see Definition~\ref{def:TranslationCovering}).
For example, the commutator subgroup $\left[F_2,F_2\right]$ is in
correspondence with the covering $p:\C\setminus \Z^2\to \T^2$ with deck
transformation group isomorphic to $\Z^2 = F_2 / \left[F_2,F_2\right]$.

\begin{tcbexercise}{}{}
Describe a subgroup $H< F_2$ that determines the infinite staircase as a square-tiled surface. Is H in this case a normal subgroup?
\end{tcbexercise}

From the constructive point of view (Definition~\ref{def:TranslationSurfaceConstructive}),
a square-tiled surface corresponds exactly to translation surfaces
built from copies of the unit square.
More precisely, let $C = [0,1] \times [0,1]$ be the unit square (thought
as a Euclidean polygon) and let $I$ be an at most countable set.
Let $r,u:I\to I$ (the "right" and "up" gluing rules). Then we construct a surface
from $\bigcup_{i \in I} C \times \{i\}$ by making for each $i \in I$ the following
identifications of edges by translation
\begin{itemize}
\item the right edge of $C \times \{i\}$ with the left edge of $C \times \{r(i)\}$ and
\item the top edge of $C \times \{i\}$ with the bottom edge of $C \times \{u(i)\}$.
\end{itemize}
If the permutations $r$ and $u$ act transitively on $I$ (that is for any pair
$i,j \in I$ there is a sequence of elements $s_1$, \ldots, $s_n$ in
$\{r,u,r^{-1},u^{-1}\}$ such that $s_n \cdots s_2 s_1 (i) = j$) then the
resulting topological space is connected. Let us recall than when
we construct a translation surface from polygons one needs to
remove the vertices adjacent to infinitely many squares
(see Definition~\ref{def:TranslationSurfaceConstructive}).
This is the case for the four infinite angle singularities appearing in
the infinite staircase.

\begin{tcbexercise}{}{}
Describe a pair of bijections $r,u:\Z\to \Z$ for which the above construction
produces the infinite staircase.
\end{tcbexercise}{}{}

\end{tcbexample}



\begin{tcbexample}{Irrational billiards}{IrrationalBilliardAndProof}
\label{EXAMPLE:IrratBillAndProof}
In Section~\ref{ssec:PolygonalBilliards} we explained how the dynamics of a billiard ball on a Euclidean polygonal table $P$ can be interpreted as the dynamics of the translation flow on a translation surface $M(P)$ obtained from $P$ by unfolding. In the following lines we sketch the proof of Theorem~\ref{thm:AnIrrationalBilliardIsLochNess} for totally-irrational triangles. More precisely, we show that for every totally-irrational triangle the surface $M(P)$ is homeomorphic to the Loch Ness monster. The general case is detailed in \cite{Valdez09_2}.

Let $\{\xi_1,\xi_2,\xi_3\}$ be the sides of a totally irrational triangle\footnote{A Euclidean triangle $P$ is said to be \emph{totally-irrational} if it has the following property: if $\lambda_i\pi, \lambda_j\pi$ are interior angles of $P$ and $n_i\lambda_i+n_j\lambda_j\in\Z$ for some $n_i,n_j\in\Z$, then $n_i=n_j=0$. Note that totally-irrational triangles are generic in the space of triangles.} $P$ and $\mathbb{S}^2(P)$ be the Euclidean surface\footnote{By Euclidean surface we mean a surface endowed with an atlas where transition functions are isometries of the Euclidean plane.}  obtained by identifying two copies of $P\setminus\text{Vertices}(P)$ along sides with the same labels using orientation preserving isometries. The Euclidean surface
$\mathbb{S}^2(P)$ is homeomorphic to a three-punctured sphere. With respect to its natural Euclidean metric (inherited by from $P\subset\R^2$) the sectional curvature of $\mathbb{S}^2(P)$ at every point is zero but $\mathbb{S}^2(P)$ is not a translation surface.
\begin{tcbexercise}{}{}
Let $P$ be a totally-irrational triangle.
\begin{compactenum}
\item Show that there is a natural normal covering map $p: M(P) \to \mathbb{S}^2(P)$ coming from the unfolding construction whose deck transformation group is isomorphic to $\Z^2$.
\item Show that $\deck(p)$ acts by isometries on $M(P)$.
\end{compactenum}
\end{tcbexercise}

Hence $p: M(P) \to \mathbb{S}^2(P)$ is a $\Z^2$-covering but its base \emph{is not} a translation surface. We include this example to illustrate how the discussion about $G$-coverings we had before can be extended to covering maps where the base is an Euclidean surface and the covering space a translation surface.

The surface $M(P)$ has infinite genus because, by choosing the appropiate subgroup of $F_2$, one can find an intermediate covering
$$
M(P)\to M_n \to \mathbb{S}^2(P)
$$
where $M_n$ is a finite-type positive genus surface. To prove that $M(P)$ has only one end we use lemma \ref{lem:CriterionInfiniteGenus} and argue in a similar way as we did for baker's surface or the infinite staircase: if $K\subset M(P)$ is compact, then $p(K)\subset \mathbb{S}^2(P)$ is contained in $K_\epsilon$, the complement of the union of a sufficiently small $\epsilon$-neighbourhood of $\text{Vertices}(P)$ in $\mathbb{S}^2(P)$. Since $K$ is compact, it is possible to find finitely many copies of $K_\epsilon$ in $p^{-1}(K_\epsilon)$ covering $K$. The union of these finitely many copies of $K_\epsilon$ defines
a compact set $K'\subset M(P)$. Think of $M(P)$ as $D\times\Z^2$, where $D$ is a fundamental domain for the action of the deck transformation group (isomorphic to $\Z^2$) of $p: M(P) \to \mathbb{S}^2(P)$. Since $K'$ is compact, it is contained in a ball of the form:
$$
B_N:=\{D\times\{(n,m)\}\hspace{1mm}|\hspace{1mm} n^2+m^2<N\}
$$
for some large $N\in\N$. Let us remark that $M(P)\setminus B_N$ is connected, because $\Z^2$ is a group with only one end. Every point $z\in M(P)\setminus K'$ can be connected by a path to a point $z'$ which projects to $p(z')\in\mathbb{S}^2(P)\setminus K_\epsilon$. Given that $P$ is a totally irrational triangle, it is possible to connect $z'$ through an arc to a copy of $D$ outside $B_N$, that is a fundamental domain of the form $D\times\{(n_0,m_0)\}$, for some $n_0,m_0>N$. Given that the choice of $z$ was arbitrary and $M(P)\setminus B_N$ is connected, one can connect any two points in $M(P)\setminus K'$ through an arc.
\end{tcbexample}

\begin{tcbquestion}{}{}
Let $P$ be a generalized irrational polygon (see Definitions~\ref{def:GeneralizedPolygon} and~\ref{def:RatIrratPolygon}) which is not simply connected. Describe the topology of $S(P)$. Is it always homeomorphic to a Loch Ness monster?
\end{tcbquestion}




\textbf{Dilation surfaces}.
The infinite staircase forms part of a larger class
of infinite-type translation surfaces called $\lambda$-staircases which we
defined in Section~\ref{ssec:ThurstonVeechConstructionsIntro}. Every
$\lambda$-staircase admits an affine automorphism with a derivative of the form
$r_+\rm{Id}$, where $r_+>1$ is the largest root of $x^2-\lambda x +1$ whose
action is discrete. Its quotient inherits a geometric structure that is called
a dilation surface that we discuss now. This situation is very similar to
irrational billiards $p:M(P)\to\mathbb{S}^2(P)$ that we described above.

A \emphdef{dilation} is a linear isomorphism $\phi: \C \simeq \R^2$ of the form $\phi(z)
= az + b$ with $a \in \R \setminus \{0\}$. The dilation is \emph{positive} if $a >
0$. We denote by $\Dil(\R^2)$ and $\Dil^+(\R^2)$ the group of respectively
dilations and positive dilations of $\R^2$.
\begin{tcbdefinition}{Constructive}{}
Let $\cP$ be a finite set of Euclidean polygons and $f: E(\cP) \to E(\cP)$ a pairing of
the edges of the polygons in $\cP$ such that edges that are paired together
have opposite outward normal vectors. For each pair $(e, f(e))$ there is a unique
orientable dilation $\phi_e$ such that $\phi_e(e) = f(e)$. We define $M$ as
$\bigsqcup_{P \in \cP} P / \sim$ where two points $x$ and $y$ are identified if they
belong to a pair of edges $(e, f(e))$ and $\phi_e(x) = y$. If $M$ is conncted, we
call it \emph{the (compact) dilation surface obtained from the family of polygons $\cP$
and pairing $f$}.
\end{tcbdefinition}
The example of the quotient of the $\lambda$-staircase can be visualized on
Figure~\ref{fig:LambdaStaircase}.
Note that the only difference with our constructive definition of translation surfaces
(Definition~\ref{def:TranslationSurfaceConstructive}) is that we do not require that
$e$ and $f(e)$ have the same length.

We now briefly discuss how to make a geometric definition (similar to
Definition~\ref{def:TranslationSurfaceGeometric}). A \emphdef{dilation atlas}
on a topological surface $S$ is an atlas $\mathcal{T}=\{\phi_i:U_i\to\C\}$ for
which the transition maps $\phi_j\circ\phi_i^{-1}:\phi_i(U_i\cap
U_j)\to\phi_j(U_i\cap U_j)$ are orientable dilations. A \emphdef{dilation
structure} on a surface $S$ is a maximal dilation atlas. It is clear that
given a (constructive) dilation surface $M$ it admits a dilation structure on
$S \setminus \Sigma$ where $S$ is the underlying topological surface and $\Sigma$ the
image of the vertices of $\cP$ in $S$. However, the converse is not
true. A dilation atlas defined on the complement of a finite set of points
of a compact surface $S \setminus \Sigma$ does not necessarily extend to a
dilation surface.
Similarly to the case of compact translation surfaces where only conical
points are allowed, the geometric definition needs to allow only certain
local models of singularities that correspond to finite covers of the
punctured disc. More precisely, each point $p \in \Sigma$ must have a
neighborhood isomorphic to a finite cover of the dilation structure
in $\C \setminus \{0\}$ whose holonomy around $0$ is generated by
$z \mapsto r z$ for some $r > 0$.
Note that dilation structure is a particular case of $(G,X)$-structures
discussed in Appendix~\ref{Appendix:GXStructures}. Dilation surfaces also admit an analytic definition, see \emph{twisted laminations} in~\cite{McMullen00} and
\emph{twisted quadratic differentials} in~\cite{Wang21}.

For more about the geometry and dynamics of dilation surfaces, we recommand the lecture
of~\cite{DuryevFougeronGhazouani19, BoulangerFougeronGhazouani20}.

\emph{Translation coverings}. Given a compact dilation surface $M$ with underlying
surface $S$ and singularities $\Sigma$, it has an associated \emph{holonomy map} $h: \pi_1(x_0, S \setminus \Sigma) \to \Dil^+(\R^2)$.
It is constructed as follows. For a closed rectifiable curve $\gamma: [0,1] \to S$
based at $x_0$, consider a set of charts $\phi_i: U_i \to \C$ whose union of domains contain
$\gamma$. There exists real numbers $t_0 = 0 < t_1 < \ldots < t_n = 1$ such that
for each $j=0, \ldots, n-1$ the segment $\gamma([t_j, t_{j+1}])$ is contained
in the domain $U_{i_j}$ of the chart $phi_{i_j}$. The holonomy $h(\gamma)$ is the composition of
the transition between these various charts, that is
\[
h(\gamma) := (\phi_{i_n} \circ \phi_{i_{n-1}}^{-1}) \cdots (\phi_{i_2} \circ \phi_{i_1}^{-1}) \cdot (\phi_{i_1} \circ \phi_{i_0}^{-1})
\]
where each transition map $(\phi_{i_n} \circ \phi_{i_{n-1}}^{-1})$ should be extended to a dilation of the whole plane $\R^2$ before multiplying them. It is not hard to see that this composition
does not depend on the choice of the charts (see Appendix~\ref{Appendix:GXStructures} for
a proof in the more general context of $(G, X)$-structures).

Let's consider the morphism $c: \pi_1(S \setminus \Sigma, x_0) \to \R^*$ obtained from composing the holonomy
map $h: \pi_1(x_0, S \setminus \Sigma) \to \Dil^+(\R^2)$ with the derivative $\Dil^+(\R^2) \to \R^*$. This
partial holonomy map only records the "dilation" part and not the "translation" part of
the structure. The kernel of $c$ defines a normal covering $\pi:\widetilde{S}\to S \setminus \Sigma$.
The covering space $\widetilde{S}$ inherits an atlas from $S$ where all transition
functions are translations and we hence
call it the \emph{translation covering associated to the dilation surface $S$}.
Given that $\R^*$ is an Abelian torsion free group, every non-trivial
translation covering over a dilation surface $S$ with finitely generated
fundamental group has deck transformation group isomorphic to $\Z^d$, for some
$d\geq 1$. The archetypical example of this kind of coverings are
$\lambda$-staircases, as illustrated in Figure~\ref{fig:LambdaStaircase} where
$d=1$.

\section{Existence of translation structures}
	\label{Sec:Existence-translation-structure}
In this chapter we discuss the existence of translation structures on
a given topological surface. We mostly consider the analytical point of
view where translation surfaces are pairs $(X, \omega)$ made of a Riemann
surface $X$ together with a non-zero holomorphic 1--form $\omega$.

Let us first begin with the existence of conformal structures on topological
surfaces.
\begin{tcbtheorem}{}{ExistenceComplexStructureOnTopologicalSurface}
Any topological surface $S$ admits a conformal structure.
\end{tcbtheorem}

\begin{tcbremark}{}{}
The object of \emph{Teichm\"uller theory} is to answer to the question: how
many conformal structures there are on a given topological surface $S$?
When $S$ is of finite type, the conformal structures on $S$ up to isotopy are parametrized by a
complex manifold homeomorphic to a ball, the so-called Teichm\"uller space. For an infinite-type surface
$S$ there are many different Teichm\"uller spaces associated to the various
quasi-conformal structures of $S$. Each of these Teichm\"uller spaces
is an infinite dimensional (Banach) manifold. For more on this subject and the
geometry of quasi-conformal mappings, we refer the reader
to~\cite{Ahlfors-quasiconformal_mappings}, \cite{Abikoff80}, \cite{Nag88},
\cite{ImayoshiTaniguchi}, \cite{GardinerLakic00} and~\cite{Hubbard-book1}.
\end{tcbremark}

\begin{proof}
By definition the topology of any surface has a countable basis and hence $S$
can be triangulated, see Section 46A in~\cite{AhlforsSario-book}. Consider a triangulation on $S$ and then declare that
each triangle is a Euclidean equilateral triangle with edge length one. This
defines a conformal structure on the complement of the vertices. At the
vertices, since there are finitely many edges adjacent to it, there is unique
way to extend this conformal structure to the whole $S$.
\end{proof}

We consider three types of constructions on a given Riemann
surfaces $X$ that produce translation surface structures: with a prescribed set of conical singularities in
Section~\ref{ssec:Strata}, with finite area in Section~\ref{ssec:Hodge} and
finally with no periods in Section~\ref{ssec:Exact}.

\begin{tcbremark}{}{}
We quickly discuss in this Remark whether given an infinite-type Riemann
surface $X$ one can construct a holomorphic 1--form with multiple constraints
among: its divisor (as in Theorem~\ref{thm:ExistenceHolomorphicOneForm1}), its
area (as in Theorem~\ref{thm:ExistenceHolomorphicOneForm2}) and its periods (as
in Theorem~\ref{thm:ExistenceHolomorphicOneForm3}).

For compact Riemann surfaces the constraints are very rigid. First of all,
holomorphic 1--forms are square-integrable. Moreover,
each holomorphic 1--form $\omega$ is determined up to scaling by its divisor (because the quotient of two differentials with the same divisor is a holomorphic function, which
is necessarily constant).
The form $\omega$ is also determined by its periods (Torelli theorem). Since both
the space of divisors and the space of periods have higher dimensions
(respectively $2g-2$ and $g(g+1)/2$) that makes them very rigid
objects. To some extent, this rigidity extends to the class of
parabolic Riemann surfaces. For more about the analytic study of parabolic surfaces we refer
to~\cite{FeldmanKnoerrerTrubowitz}.

Beyond parabolic surfaces, understanding the constraints associated to
divisors and periods of square-integrable 1--form seem a delicate problem.
\end{tcbremark}

\subsection{Translation structures with prescribed conical angles}
\label{ssec:Strata}
A \emphdef{divisor} on a Riemann surface $X$ is a map $D:X\to\Z$ with discrete
support, \ie, for any compact subset $K\subset X$ the set $\{x\in K:D(x)\neq 0\}$
is finite. If $\omega$ is a non-zero meromorphic differential on $X$ its order of
vanishing $z \to {\rm ord}_z(\omega)$ is a divisor called
the divisor of $\omega$. Similarly, one can associate a divisor to meromorphic
functions and more generally to any meromorphic section of a holomorphic line
bundle.
\begin{tcbtheorem}{}{ExistenceHolomorphicOneForm1}
Any Riemann surface $X$ different from the sphere admits a non-zero holomorphic
1--form. Moreover if $X$ is non-compact, given a divisor $D:X \to \Z_{\geq 0}$
there exists a 1--form $\omega$ on $X$ whose divisor is $D$.
\end{tcbtheorem}

Combining Theorem~\ref{thm:ExistenceComplexStructureOnTopologicalSurface}
and Theorem~\ref{thm:ExistenceHolomorphicOneForm1}, we see that translation
structures exist on any topological surface different from the sphere.
More precisely, recall from Exercise~\ref{exo:ConicalSingularities} that the
order of vanishing $d = {\rm ord}_z(\omega)$ of the differential $\omega$
determines the angle $\alpha$ of the conical singularity at $z$ of the
translation structure $(X, \omega)$ by the formula $\alpha = 2\pi(d + 1)$. In
particular, given a non-compact topological surface $S$, and any finite or
infinite discrete multiset $A$ of positive multiples of $2\pi$ there exists a translation
structure on $S$ whose multiset of conical singularities is equal
to $A$.

\begin{tcbremark}{}{rk:OneFormWithoutZero}
 Let us mention that if $\Sigma$ is chosen to be empty in
Theorem~\ref{thm:ExistenceHolomorphicOneForm1}, then the obtained translation
surface does not have conical singularities.
\end{tcbremark}

\begin{proof}

If $X$ is a compact Riemann
surface of genus $g>0$, then $X$ always admits a translation surface structure
because the space of holomorphic 1--forms on $X$ has complex
dimension\footnote{This can be deduced from the Riemann-Roch theorem, see for example Theorem 16.9 in~\cite{Forster77}.} $g$.

We now assume that $X$ is non-compact. Hence $D:X \to
\Z_{>0}$ on $X$ is the divisor of a holomorphic function $f$, see~\cite[Chapter 26]{Forster77}. Any holomorphic vector bundle $E$ on $X$ is (holomorphically)
trivial, see~\cite[Theorem~30.4]{Forster77}. In particular, the
canonical bundle\footnote{For every pair of charts $(U_i,z_i)$ and $(U_j,z_j)$
of $X$ the function $g_{ij}:=\frac{dz_j}{dz_i}$ is holomorphic in $U_i\cap U_j$
and does not vanish. In the language of vector bundles, the family $(g_{ij})$
is a cocycle. It defines a line bundle $K_X$ over $X$ called the
\emphdef{canonical line bundle} of $X$. The holomorphic sections of $K_X$ are the
holomorphic 1--forms on $X$.} $K_X$ is trivial. Let then $\omega$ denote a
non-zero constant section of $K_X$ in any trivialization. This holomorphic 1--form
is nowhere vanishing.

To obtain the desired translation surface structure one simply takes the product $f\omega$
as its divisor is $(f\omega) = (f) = D$.
\end{proof}

\subsection{Hodge decomposition and integrable holomorphic 1--forms}
\label{ssec:Hodge}
We now turn to the construction of translation structure of finite area via
Hodge theory. Two relevant sources for this section are~\cite[Chapter V]{AhlforsSario-book}
and~\cite[Chapter 1]{FeldmanKnoerrerTrubowitz}.

Given a Riemann surface $X$, we say that a non-zero holomorphic 1--form $\omega$ on $X$
is \emphdef[square-integrable 1--form]{square-integrable} if the area of the translation
surface $(X, \omega)$ is finite. This terminology is justified in Exercise~\ref{exo:SquareIntegrable}.

\begin{tcbtheorem}{}{ExistenceHolomorphicOneForm2}
Let $X$ be a Riemann surface which is not of genus $0$. Then $X$ admits a
non-zero square-integrable holomorphic 1--form. In other words, a translation
structure of finite area.
\end{tcbtheorem}

Let $X$ be a Riemann surface. The \emphdef{Hodge star operator} on $X$ is the
linear operator that maps a (real or complex) one form $\omega = f(x,y) dx +
g(x,y) dy$ in local coordinates\footnote{It is
a simple computation that this definition is independent on the choice of local coordinates.} $z = x + iy$ to $*\omega = -g(x,y) dx + f(x,y) dy$.
We denote by $L^2(X, T^*X)$ the space of measurable square-integrable real one
forms on $X$ with respect to the scalar product
\[
(\alpha, \beta) = \int_X \alpha \wedge *\beta.
\]
The space of \emphdef{square-integrable harmonic one forms} on $X$ is
\[
\cH(X) = \{\alpha \in L^2(X, T^*X): d \alpha = d *\alpha = 0\}.
\]
\begin{tcbexercise}{}{SquareIntegrable}
The Hodge star operator extends to the complex valued 1--form by linearity.
\begin{enumerate}
\item Show that in a local chart $z = x + iy$ we have $*dz = -i d\overline{z}$.
\item Show more generally that for a holomorphic 1--form $\omega$, we have $*\omega = -i \overline{\omega}$.
\item Deduce that the area form $\frac{i}{2} \int_X \omega \wedge
\overline{\omega}$ of a translation surface $(X, \omega)$, see
Exercise~\ref{exo:ConicalSingularityAndZeroOfDifferential}, can alternatively
be written as $\frac{1}{2} \int_X \omega \wedge *\omega$.
\item Conclude that a translation surface $(X, \omega)$ has finite area if and only if
$\omega$ belongs to $L^2(X, T^*X)$.
\end{enumerate}
\end{tcbexercise}
The Hodge decomposition theorem
(see~\cite[Section V.10 "Orthogonal decompositions"]{AhlforsSario-book}
or~\cite[Theorem 1.3]{FeldmanKnoerrerTrubowitz}) states that we have
an orthogonal decomposition
\[
L^2(X, T^* X) = \overline{d C_c^\infty(X)} \oplus \overline{* d C_c^\infty(X)} \oplus \cH(X)
\]
where $C_c^\infty(X)$ are the smooth functions with compact support on $X$.

\begin{proof}
We present now the proof of Theorem~\ref{thm:ExistenceHolomorphicOneForm2}.
A complex 1--form $\omega$ is holomorphic if and only if its real
part is harmonic (that is satisfies $d \alpha = d *\alpha = 0$) and
$\Im(\omega) = * \Re(\omega)$ (Cauchy-Riemann equations).
In particular, $(X, \omega)$ has finite area as a translation surface
if and only if its real part $\alpha = \Re(\omega)$ belongs to $\cH(X)$.
It is hence sufficient to show that $\cH(X)$ contains nonzero elements.

Now, to construct a non-trivial harmonic 1--form we rely on
orthogonal projections in $L^2(X,T^*X)$. Namely, given a closed
curve $\gamma$ in $X$ one can then construct a smooth 1--form with
support in a tubular neighborhood of $\gamma$ such that
$\int_\gamma \eta'_\gamma = 1$. We define $\eta_\gamma$ as the projection
of $\eta'_\gamma$ on $\cH(X)$ (the harmonic 1--forms). By Hodge theorem,
for any pairs of closed curves $\gamma_1$ and $\gamma_2$ we have
\begin{equation}
\label{eq:HodgeAndAlgebraicIntersections}
(\eta_{\gamma_1}, * \eta_{\gamma_2})
=
\langle \gamma_1, \gamma_2 \rangle
\end{equation}
where the inner product on the left is the inner product on $L^2(X, T^*X)$ while
the one on the right is the algebraic intersection form. Now the surface $X$ has
genus if there are two cycles $\gamma_1$ and $\gamma_2$ such that
$\langle \gamma_1, \gamma_2 \rangle \not= 0$. For such pairs the associated harmonic
1--forms $\eta_{\gamma_1}$ and $\eta_{\gamma_2}$ are non-trivial
by~\eqref{eq:HodgeAndAlgebraicIntersections}. This concludes the proof of the
first part.
\end{proof}

\subsection{Exact 1--forms and holomorphic functions}
\label{ssec:Exact}
Let $X$ be a Riemann surface. If there exists a non-constant holomorphic map
$\phi: X \to \C$ then the pull-back $\omega := \phi^* dz = d \phi$ is a
1--form on $X$. By construction, such form is exact. Conversly, if $\omega$
is an exact holomorphic 1--form on $X$, then choosing a point $z_0 \in X$ the
period map
\[
\begin{array}{lll}
X & \to & \C \\
z & \mapsto & \int_{z_0}^z \omega
\end{array}
\]
is well-defined and provides a non-constant holomorphic map to $\C$. In other
words we proved the following result.
\begin{tcbproposition}{}{}
Let $X$ be a Riemann surface. Then the following are equivalent.
\begin{enumerate}
\item $X$ admits an exact non-zero holomorphic 1--form,
\item there exists a non-constant holomorphic map $X \to \C$.
\end{enumerate}
\end{tcbproposition}
By the maximum principle, for compact Riemann surface $X$ any holomorphic map
$X \to \C$ is constant. The following result shows that in the non-compact
setting one can even impose to the map $\phi$ to have a nowhere vanishing differential,
or equivalently for $\omega := \phi^* dz$ to have no zeroes.
\begin{tcbtheorem}{}{ExistenceHolomorphicOneForm3}
Let $X$ be a non-compact Riemann surface. Then $X$ admits an exact holomorphic
1--form $\omega$ that vanishes nowhere.
\end{tcbtheorem}

Many examples of translation surfaces $(X, \omega)$ as in
Theorem~\ref{thm:ExistenceHolomorphicOneForm3} are obtained by considering
an open subset $U$ of the complex plane $\C$ endowed with the canonical
differential $dz$. Note also that the set of translation surfaces
$(X, \omega)$ where $\omega$ is an exact holomorphic 1--form
that vanishes nowhere is stable under taking subsurfaces and taking
non-ramified coverings.

We will only give a sketch of the proof of
Theorem~\ref{thm:ExistenceHolomorphicOneForm3}. We refer to the original
article~\cite{GunningNarasimhan67} for more details (see
also~\cite{KusunokiSainouchi71} for an extension).

The main tools in the construction are based on the fact that any open Riemann
surface admits a nowhere vanishing holomorphic 1--form (see proof of
Theorem~\ref{thm:ExistenceHolomorphicOneForm1}) and the following result:
\begin{tcbtheorem}{Mergelyan-Bishop theorem}{MergelyanBishop}
Let $X$ be a non-compact Riemann surface and $K \subset X$ a compact subset
whose complement has no relatively compact connected component. Then any
function on $K$ which is holomorphic on the interior of $K$ is the
uniform limit of restrictions of holomorphic functions on $X$.
\end{tcbtheorem}

\begin{proof}
We address now the proof of Theorem~\ref{thm:ExistenceHolomorphicOneForm3}. We provide a sketch of the original proof as in~\cite{GunningNarasimhan67}.

Recall that we are given an open Riemann surface $X$ and want to build
a differential $\omega$ which does not vanish anywhere and whose periods
are all zero. The proof proceeds by induction. We pick an exhaustion $U_0
\subset U_1 \subset \ldots $ of $X$ by relatively compact subsets with smooth
boundaries and with no compact connected component in the complement. Our
aim is to build a sequence of nowhere vanishing holomorphic differential $\omega_n$
such that
\begin{enumerate}
\item the sequence $\omega_n$ converges to a nowhere vanishing 1--form $\omega_\infty$,
\item the integral of $\omega_n$ along curves contained in $U_n$ are zero.
\end{enumerate}

Without loss of generality, we can assume $U_0$ to be a disk. We can then take
$\omega_0$ be any nowhere vanishing differential on $X$ which exists
by Theorem~\ref{thm:ExistenceHolomorphicOneForm1} 
This provides the initial step of the induction.
Now assume that we constructed $\omega_n$. We sketch how to build $\omega_{n+1}$.
Using Theorem~\ref{thm:MergelyanBishop}, one can construct a holomorphic function $f_n: X \to \C$
such that $|f_n| < e^{-n}$ inside $U_n$ and the integral of $\omega_{n+1}
:= e^{f_n} \omega_n$ vanishes along curves inside $U_{n+1}$.

The limiting differential is obtained by considering the pointwise limit
$\omega_\infty := e^{f_1 + f_2 + \ldots} \omega_0$ which exists by our requirements
on the smallness of $f_n$ inside $U_n$.
\end{proof}

\section{Singularities and the space of linear approaches}
\label{sec:Singularities}

In this section we explain how the natural extension of the tangent bundle of a translation surface $M$ to its metric completion provides the necessary invariants to describe the singularities of $M$. Let us begin by recalling two indispensable notions.\\

\textbf{Saddle connections and holonomy vectors}. As we saw in the preceding chapter, the term \emph{flow} is an abuse of language when referring to the translation flow $F_\theta^t$ on $M$ for it might not be defined for all points in $M$ for all $t\in\R$. We distinguish those orbits of $F_\theta^t$ which are not defined for all times in the future or in the past.

\begin{tcbdefinition}{}{SaddleConnection}
\label{def:saddle_connection} \label{def:holonomy_vector}
Let $M$ be a translation surface, $z\in M$ and $I\subset\R$ the \emph{maximal} domain of definition of $F_\theta^t(z)$.
\begin{itemize}
\item If $I=\R$, the orbit of $z$ under $F_\theta^t$ is called a \emph{regular geodesic}\footnote{This nomenclature is justified by remark~\ref{rk:GeodesicFlowVSTranslationFlow}.},
\item if $I$ is unbounded, but $I\neq \R$, this orbit is called a \emph{separatrix}. Furthermore, if $I$ is bounded below (respectively above), the orbit is called a forward separatrix (respectively backward), and
\item if $I$ is bounded, the orbit of $z$ under $F_\theta^t$ called a \emph{saddle connection} in direction $\theta$.
\end{itemize}
If $\gamma$ is the trace of a saddle connection in direction $\theta$, we denote by $v_\gamma$ the vector of length $|\gamma|$ and direction $\theta$.  This vector is called the  \emph{holonomy vector associated} to the saddle connection $\gamma$. The set of all holonomy vectors is denoted by $V_{hol}(M)$.
\end{tcbdefinition}

It is clear from the definition that saddle connections are geodesics whose extremities are points in $\Sing(M)$ or marked points\footnote{A marked point or a puncture is a regular point in $\widehat{M}$ that does not appear in $M$.}, and that $v\in V_{hol}(M)$  implies that $-v\in V_{hol}(M)$. On the other hand, affine automorphisms send saddle connections to saddle connections and hence $\Gamma(M)$, the Veech group of $M$, acts on $V_{hol}(M)$. As we see in the next chapter, this action can be used to deduce properties of $\Gamma(M)$, for example it is a well known fact (see \cite{MasurTabachnikov02}) that for finite-type translation surfaces $V_{hol}(M)$ is an infinite discrete subset of $\R^2$ and hence $\Gamma(M)$ must be discrete in this context.

\begin{tcbexercise}{}{HolonomyVectors}
Let $P$ be a polygon and $M(P)$ the translation surface obtained by unfolding the billiard  $P$ as described in Section~\ref{ssec:PolygonalBilliards}.
\begin{enumerate}
\item Show that the natural projection $\pi_P:M(P)\to P$ sends saddle connections to \emph{generalized diagonals}, that is, billiard trajectories whose extremities are vertices of $P$. Is $\pi_P$ restricted to a saddle connection a bijection? Is every generalized diagonal in $P$ the image of a saddle connection?
\item For which triangles $P$ does $M(P)$ has no conical singularities?
\item Show that if $P$ is a triangle whose interior angles are all irrational multiples of $\pi$, then $V_{hol}(M(P))$ is an infinite non-discrete subset of the plane that is contained in the complement of a disc $\{(x,y)\in\R^2\hspace{1mm}|\hspace{1mm} |z|< r \}$ for some $r>0$ that depends only on $P$.
\end{enumerate}
\end{tcbexercise}

\begin{tcbquestion}{}{}
Let $P$ be an irrational polygon. Is $V_{hol}(M(P))\subset\R^2$ bounded above? Is the interior of the closure of $V_{hol}(M(P))$ empty?
\end{tcbquestion}

In the preceding chapter we introduced conical, infinite angle and wild singularities (see Exercise~\ref{exo:ConicalSingularities} and Definition~\ref{def:Singularities}) and the notion of tame translation surface $M$: these are surfaces for which $\Sing(M)$ does not contain wild singularities. If $M$ is a tame translation surface we can describe a small neighbourhood of any point in its metric completion: it is either isometric to a neighbourhood of $0\in\C$ or to a cyclic (maybe infinite) covering of a neighbourhood of $0\in\C$ ramified over $0$. The class of tame translation surfaces includes all compact translation surfaces, those arising from billiards on irrational polygons, wind-tree models, and, as the following exercise shows, tameness is preserved when taking coverings.

\begin{tcbexercise}{}{}
Prove that if  $p:\widetilde{M}\to M$ is a $G$-covering whose base is a tame translation surface (not necessarily of finite topological type), then $\widetilde{M}$ is also a tame translation surface.
\end{tcbexercise}

\begin{tcbexample}{Step surface}{InfiniteStepSurface}
We now introduce step surfaces studied
in~\cite{DegliEspostiDelMagnoLenci-1998,DegliEspostiDelMagnoLenci-2000}.
These examples show that, contrary to intuition, infinite-type tame translation
surfaces of finite area exist. Consider the generalized polygon
$P:=\cup_{n\in\N}[n-1,n]\times [0,h_n]$, where $(h_n)_{n\in\N}$ is a
monotonically decreasing sequence converging to zero. A polygonal billiard
having $P$ as table is called an \emphdef[step billiard]{(infinite) step
billiard}\footnote{The term \emph{italian billiard} is also
used to refer to this family of tables.}. A billiard ball in $P$ can move in
at most 4 directions, therefore the translation surface $M(P)$, called a
\emphdef[step surface]{(infinite) step surface}, obtained by
unfolding (see Section~\ref{ssec:PolygonalBilliards}) is tiled by four copies
of $P$. As illustrated in Figure~\ref{fig:StepSurface}, these copies are $\{P,\sigma_xP,\sigma_y\sigma_xP,\sigma_y P\}$, where $\sigma_x$ and $\sigma_y$ are the reflections of the plane w.r.t. the horizontal and vertical lines passing through the origin respectively. The area of $M(P)$ is equal to
$4\sum_{n=1}^\infty h_n$, hence $M(P)$ has finite area if $(h_n)_{n\in\N}$ is
summable. It is not difficult to check that the diameter of $M$ is infinite,
and hence the only points defining singularities come from the vertices of $P$
where the interior angle is equal to $\frac{3\pi}{2}$. More precisely,
$\Sing(M)\subset M$ is discrete and formed exclusively by conic singularites of
total angle $6\pi$. In particular, $M(P)$ is always tame.
\begin{figure}[H]
\begin{minipage}{0.4\textwidth}
	\includegraphics[scale=0.7]{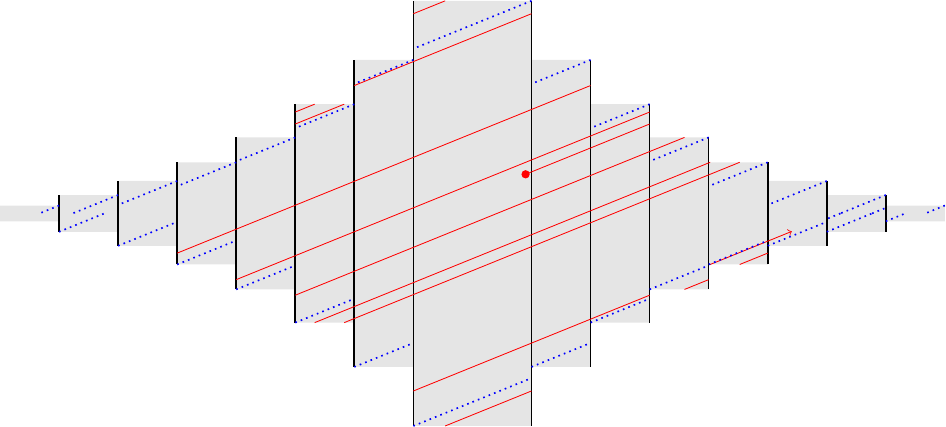}
\end{minipage}
\caption{A step surface with steps $h_n = \left(\frac{3}{4}\right)^n$.}
\label{fig:StepSurface}
\end{figure}

\end{tcbexample}

\subsection{The space of linear approaches}
\label{ssec:LinearApproaches}


Every translation surface $M$ can be endowed with a distance coming from the intrinsic flat metric which allows us to consider, as we defined in the preceding chapter, its metric completion $\widehat{M}$ and $\Sing(M)\subset\widehat{M}$ the set of singularities. Given that $M^0:=\widehat{M}\setminus \Sing(M)$ has no singularities, its tangent bundle $TM^0$ is isomorphic to the product $M^0\times\R^2$. In the following paragraphs we introduce a continuous extension of this bundle, with its zero section removed, to $\widehat{M}$. At the level of the \emph{unit} tangent (sub)bundle, this extension produces a topological space called the \emph{space of linear approaches} of $M$, which we use to describe the geometry of a translation surface near any singularity. The germ of the idea behind the space of linear approaches can be found in R. Chamanara's work \cite{Chamanara}. As seen in example \ref{ssec:BakerBowmanArnouxYoccoz}, the metric completion of baker's surface $B_\alpha$ with respect to its natural metric is achieved by adding one wild singularity $z_\infty$. R. Chamanara points out that ``\emph{[g]eometrically, the surface spirals infinitely many times around this point}''. The space of linear approaches formalizes this intuition. Moreover, one can extract enough information from this space to calculate the Veech group of baker's surface, as done in Section \ref{ssec:VeechGroupBakersSurface}. The approach and notations we use to present the space of linear approaches slightly differ from those used in \cite{BowmanValdez13} and \cite{Randecker16}, which are the principal references on the subject.

\begin{tcbdefinition}{}{ExtensionTangentBundle}
\label{Def:ExtensionTangentBundle}
Given $\eps>0$, let
\[
\mathcal{L}^{\eps}(M)
:=\{\text{geodesic trajectories}\ \gamma:(0,\eps)\to M^0\},
\]
where geodesics\footnote{Instead of speaking in terms of geodesics one can phrase the definitions in this section in terms of trajectories of the translation flows $F_\theta^t$ on $M$. We have decided not to do this to be more coherent with the exposition found in \cite{BowmanValdez13} and to keep the notation simple.} are taken with respect to the flat metric on $M$. Two geodesics $\gamma_1\in\mathcal{L}^{\eps_1}(M)$ and $\gamma_2\in\mathcal{L}^{\eps_2}(M)$ are said to be
equivalent if and only if $\gamma_1(t)=\gamma_2(t)$ for all
$t\in(0,\min\{\eps_1,\eps_2\})$. We denote this equivalence relation by
$\sim$ and define:
\begin{equation}
T\widehat{M}:=\bigsqcup_{\eps>0}\leX{\huge /}\sim
\end{equation}
The equivalence class of $\gamma$ will be denoted by $[\gamma]$.
\end{tcbdefinition}

To understand better what elements of $T\widehat{M}$ are, endow $TM^0=M^0\times\R^2$ with its natural product topology and suppose that there exists $\gamma:(0,\eps)\to M^0$ in $\mathcal{L}^\eps(M)$ for which $\lim_{t\to 0}\gamma(t)=z_0$ exists in $M^0$. Then, given that geodesics are parametrized with constant non-zero speed, the limit $\lim_{t\to 0}(\gamma(t),\gamma')=(z_0,v_0)$ exists in $TM^0$ (with its zero section removed) and is independent of the representative chosen within $[\gamma]$. Hence we can associate to $[\gamma]$ the tangent vector $\gamma'$ based at $z_0$ and think that this class encodes the fact that one can \emph{approach} the point $z_0$ using the a translation flow with speed $v_0$.

In general $\lim_{t\to 0}\gamma(t)$ only exists in the metric completion $\widehat{M}$ and is independent from the representative chosen in $[\gamma]$. Hence we can think of $(\lim_{t\to 0}\gamma(0),\gamma')$ as the tangent vector $\gamma'$ based at the point $z_0=\lim_{t\to 0}\gamma(t)$, keeping in mind that strictly speaking there is no tangent space at $z_0$ for in general one cannot extend the differentiable structure of $M^0$ to all points in the metric completion.

Let us now define a topology on $T\widehat{M}$. The uniform convergence of functions induces the uniform topology on $\leX$. For each $\eps,\eps'\in\R^+$ we say that $\eps\ll\eps'$ if $-\eps\leq-\eps'$, where $\leq$ is the standard order in $\R$. This is just a formality so that for each $\eps\ll\eps'$ the restriction of element in
$\leX$ to the interval $(0,\eps')$ defines a continuous injection:
\[
\rho_{\eps}^{\eps'}:\leX\to\leXp
\]
and hence $\langle\leX,\rho_{\eps}^{\eps'}\rangle$ is a direct system of
topological spaces over $(\R^+,\ll)$. Since for every
$\epsilon\ll\epsilon'$ the projection map $\gamma\mapsto [\gamma]$ from
$\leX$ to $T\widehat{M}$ is injective and commutes with $\rho_{\eps}^{\eps'}$,
we have the equality of sets
\begin{equation}\label{lX}
T\widehat{M}=\varinjlim \leX
\end{equation}
It is then natural to endow $T\widehat{M}$ with the limit topology. Let us describe more precisely this topological space. For every $x\in M^0$, $r,t>0$ define:
 \begin{equation}
	\label{E:BasisTopLA}
B(x,t,r):=\{[\gamma] \hspace{1mm} | \hspace{1mm} d(\gamma(t),x)<r\},
\end{equation}
where $d(\dot,\dot)$ denotes the distance in $M^0$ defined by the flat metric. Note that if $[\gamma]\in B(x,t,r)$, then there exist a representative $\gamma:(0,\eps)\to M^0$ for which $\eps>t$. We claim that the family of sets $\mathcal{B}:=\{B(x,t,r)\}_{x\in M^0;r,t>0}$ is a basis for the limit topology. Indeed, it is sufficient to note that for every fixed $\eps>0$, the family of sets
$$
\mathcal{B}_\eps:=\{B_\eps(x,t,r):=\{\gamma\in\leX \hspace{1mm} | \hspace{1mm} d(\gamma(t),x)<r\}\}_{x\in M^0;r,t>0}
$$
is a basis for the topology on $\leX$.


\begin{tcbexercise}{}{dirbp}
Let $M$ be a translation surface. To every $[\gamma]\in T\widehat{M}$ we can associate its base point ${\rm bp}[\gamma]:=\lim_{t\to 0}\gamma(t)\in\widehat{M}$ and its direction\footnote{Here the convention is to send the horizontal direction to zero.} ${\rm dir}[\gamma]=\gamma'\in\R^2$. This association defines two maps ${\rm bp}:T\widehat{M}\to\widehat{M}$ and ${\rm dir}:T\widehat{M}\to\R^2$ called the \emph{basepoint} and \emph{direction} map respectively. Prove that:
\begin{enumerate}
\item The space $T\widehat{M}$ is Hausdorff and second countable.
\item Let $(TM^0)^*$ denote the tangent bundle $TM^0$ with its zero section removed. Show that the map
$$
i:(TM^0)^*\to T\widehat{M}
$$
that associates to each $(z,v)$ the class  $[\gamma]$ defined by $({\rm bp}[\gamma],{\rm dir}[\gamma])=(z,v)$ is a topological embedding.
\end{enumerate}
\emph{Hint}: use the fact that $M$ is a Riemann surface and hence second countable; prove that the maps ${\rm bp}$ and ${\rm dir}$ are continuous.
\end{tcbexercise}


\emph{Functoriality}. Let $f:M_1\to M_2$ be an affine map whose derivative lies in $\GL(2,\R)$ and $\widehat{f}:\widehat{M_1}\to\widehat{M}_2$ its (unique) continous extension to the metric completion. We define the map $\widehat{Tf}:T\widehat{M}_1\to T\widehat{M}_2$ by setting $\widehat{Tf}[\gamma]:=[f\circ\gamma]$. It is clear from the defintion that $\widehat{TId}=Id_{T\widehat{M}}$ and that $\widehat{T(f\circ g)}= \widehat{Tf}\circ \widehat{Tg}$. Moreover, if $\lim_{t\to 0}\gamma(t)=z_0$ for any given $\gamma\in\mathcal{L}^\eps(M)$, then ${\rm bp}(\widehat{Tf}[\gamma])=\widehat{f}(z_0)$ and ${\rm dir}(\widehat{Tf}[\gamma])=Df\cdot\gamma'$. It is in this last sense that $\widehat{Tf}$ is an extension of the tangent map $Tf:TM_1^0\to TM_2^0$.


\begin{tcbexercise}{}{functoriality}
Prove that the tangent map $\widehat{Tf}$ defined above is continuous. \emph{Disclaimer}: it will be shown in Example~\ref{exa::NotRegularSpace} below the space $T\widehat{M}$ is not in general a regular ($T_3$) space and hence this exercise does not follow immediately from classical extension theorems.
\end{tcbexercise}

We now introduce the space of linear approaches of a translation surface. The key point here is that in order to understand the geometry near a singularity it is sufficient to consider approximations at unit speed.

\begin{tcbdefinition}{}{SpaceOfLinearApproaches}
	\label{Def:space of linear approaches}
Let $M$ be a translation surface. The space of linear approaches of $M$ is the subspace of $T\widehat{M}$ defined by:
\begin{equation}
T^1\widehat{M}:=\{[\gamma]\in T\widehat{M}\hspace{1mm}|\hspace{1mm} \text{$\gamma'$ is a unit vector}\}
\end{equation}
Every $[\gamma]\in T^1\widehat{M}$ is called a \emph{linear approach} to the point $\lim_{t\to 0}\gamma(t)\in\widehat{M}$. For every $z\in\widehat{M}$, the subspace $T^1_z\widehat{M}:={\rm bp}^{-1}(z)\cap T^1\widehat{M}$ is called the space of linear approaches to the point $z\in\widehat{M}$.
\end{tcbdefinition}
Recall that $\mathbb{S}^1\subset\C$ denotes the unit circle given by the image of the segment $[0,2\pi]$ under the exponential map. If $z$ is a regular point of $M$ the restriction ${\rm dir}_|:T^1_z\widehat{M}\to \mathbb{S}^1$ defines a homeomorphism and hence the space of linear approaches to a regular point is naturally parametrized by $\R/2\pi\Z$. If $z$ is a conical singularity of angle $2k\pi$ the map ${\rm dir}_|$ can be lifted to a homeomorphism between $T^1_z\widehat{M}$ and the natural $k:1$ covering of $\mathbb{S}^1$, hence in this case $T^1_z\widehat{M}$ is parametrized by $\R/2k\pi\Z$. In the same line of thought it can be proven that if $z$ is an infinite angle singularity, $T^1_z\widehat{M}$ is naturally parametrized by $\R$. A simple way to understand these facts is to pick a linear approach $[\gamma]$ in $T^1_z\widehat{M}$  and think of it as a geodesic segment anchored at the point $z={\rm bp}[\gamma]$. The aforementioned parametrizations are achieved by ``turning the geodesic segment around $x$" clockwise and  counter-clockwise until we come back to the initial position or we browse the whole space $T^1_z\widehat{M}$.

\begin{tcbexercise}{}{}
Is there a translation surface $M$ for which there exists a point $z$ in the metric completion $\widehat{M}$ such that $T^1_z\widehat{M}=\emptyset$?
\end{tcbexercise}

We summarize the preceding discussion in the following:

\begin{tcbtheorem}{}{SpaceOfLinearApproaches}
	\label{THM:SpaceOfLinearApproaches}
Let $M$ be a translation surface, $\widehat{M}$ its metric completion w.r.t. the flat metric and $M^0=M\setminus\Sing(M)$. There exists continuous extensions of $(TM^0)^*$ (the tangent bundle $TM^0$ with its zero section removed) and the unit tangent $T^1M^0$ to $\widehat{M}$, that we denote by $T\widehat{M}$ and $T^1\widehat{M}$ respectively. These spaces are Hausdorff, second countable but in general not regular.
Both $(TM^0)^*$ and $T^1M^0$ are dense in $T\widehat{M}$ and $T^1\widehat{M}$. Moreover, these extensions are functorial in the category of translation surfaces: the derivative of any every affine map $f:M_1\to M_2$ whose derivative lies in $\GL(2,\R)$ can be extended to a continuous map $\widehat{Tf}:T\widehat{M}_1\to T\widehat{M}_2$. This extension is such that $\widehat{TId}=Id_{T\widehat{M}}$ and that $\widehat{T(f\circ g)}= \widehat{Tf}\circ \widehat{Tg}$.
\end{tcbtheorem}

\begin{tcbremark}{}{}
The notions of linear approach and the space of linear approaches can be also defined using the formalism of germs and stalks coming from algebraic geometry.
\end{tcbremark}

\subsection{Rotational components}
	\label{SSEC:RotationalComponents}
The idea  of ``turning a geodesic segment around a point'' is formalized with the notion of  \emph{rotational component},
which we define in the following paragraphs. First we need the notion of an \emph{angular sector}. We follow the approach of R. Schwartz in~\cite{Schwartz11}.


Consider two infinite rays $r_1$ and $r_2$ emanating from the origin in $\R^2$. A \emph{sector} in $\R^2$ is one of the connected components of $\R^2\backslash\{r_1,r_2\}$ and its angle is defined as the angle between the corresponding rays measured from the interior of the sector. A Euclidean angular sector is a space obtained by gluing together \emph{along rays in the boundary} at most countably many angular sectors in a consecutive pattern. This consecutive pattern might be cyclic, case in which we say that the Euclidean angular sector is \emph{closed}. For example, for every $n\in\N\cup\{\infty\}$,  the $n$-to-$1$ translation covering of $\R^2$ ramified over the origin is considered to be a Euclidean angular sector. The cone point of a Euclidean angular sector is the equivalence class of the origin(s) under gluing and it is the only point in the Euclidean angular sector that does not have a neighbourhood that is locally isometric to the plane. The \emph{angle} of a Euclidean angular sector is the sum of the angles of the sectors forming it (which might be divergent).
\begin{tcbexercise}{}{}
Prove that a \emph{closed} Euclidean angular sector defines a translation surface if and only if its total angle is an integer multiple of $2\pi$. Is this translation surface of finite type?
\end{tcbexercise}

\begin{tcbdefinition}{Angular sector}{AngularSector}
\label{DEF:GeomAngularSectorSurf}
An \emphdef{angular sector} $A$ in a translation surface $M$ is a subset that is isometric to a \emph{punctured} neighbourhood of the cone point of a Euclidean angular sector.
\end{tcbdefinition}

\begin{tcbremark}{}{}
It is important to stress that our definition requires the angular sector to be a subset of a translation surface. If $A\subset M$ is an angular sector of total angle $\theta\in(0,\infty]$, the point $z\in\widehat{M}$ in the metric completion
corresponding to the cone point of the Euclidean angular sector in question is well-defined. For this reason we denote sometimes $A=A(\theta,z)$ and call $z$ the cone point of the angular sector. It is important to remark that the cone point of an angular sector \emph{is not necessarily a cone angle singularity}, for it might be a wild singularity.
\end{tcbremark}




\begin{tcbdefinition}{Rotational component}{RotationalComponent}
\label{Definition:RotationalComponent}
Let $[\gamma_1]$ and $[\gamma_2]$ be two linear approaches in
$T^1_z\widehat{M}$. We say that $[\gamma_1]$ and $[\gamma_2]$ belong to the same rotational component
if and only if there exist representatives $\gamma_i:(0,\eps_i)\to M^0$,
$i=1,2$ and an angular sector $A(\theta,z)$ containing the images of $\gamma_1$ and $\gamma_2$ to $M^0$.
We denote by $Rot[\gamma]$ the set of all linear approaches belonging to the same rotational component as
$[\gamma]$ and we call it the \emphdef{rotational
component} of $T^1_z\widehat{M}$ containing $[\gamma]$. The \emph{total angle} of $Rot[\gamma]$ is the supremum of all $\theta$ for which there exist $[\gamma_1],[\gamma_2]\in Rot[\gamma]$ with representatives $\gamma_1,\gamma_2$ whose traces are contained in an angular sector $A(\theta,z)$. By convention, the total angle of a singleton is zero. A rotational component of total angle $2\pi$ is called regular.  Let us define $IrrRot(M)$ to be the set:
\begin{equation}
	\label{EQ:NonRegRot}
\{Rot[\gamma]\hspace{1mm}|\hspace{1mm} Rot(\gamma)\hspace{1mm}\text{is a non-regular rotational component of $M$}\}
\end{equation}
and $IrrRot(M)=: IrrRot(M)_+\sqcup IrrRot(M)_0$, where $ IrrRot(M)_+$ is the set of non-regular rotational components whose total angle is positive and $IrrRot(M)_0$ is its complement in $IrrRot(M)$.
\end{tcbdefinition}

Loosely speaking, if  two linear approaches  $[\gamma_1]$ and $[\gamma_2]$ in $T^1_z\widehat{M}$ are in the same rotational component and correspond to rays $r_1$ and $r_2$ in some Euclidean angular sector isometric to $A(\theta,z)$ then one can ``rotate $[\gamma_1]$ into $[\gamma_2]$",
by considering an angular displacement that takes $r_1$ into $r_2$. This intuition can be formalized using the function $\rm dir$ defined in exercise \ref{exo:dirbp}. Indeed, if we fix $[\gamma_0]$ in a  rotational component $Rot[\gamma]$ and ${\rm dir[\gamma_0]}=\theta_0\in\R/2\pi\Z$, then the choice of a point $\widetilde{\theta_0}\in\R$ in the fiber over $\theta_0$ of the universal covering $\R\to\R/2\pi\Z$ defines an injective lift
\begin{equation}
	\label{eq:InjectiveLift}
\widetilde{\rm dir}:Rot[\gamma]\to\R
\end{equation}
In this context if $\widetilde{\theta_1}$ and $\widetilde{\theta_2}$ are the images of $[\gamma_1]$ and $[\gamma_2]$ by the map $\widetilde{\rm dir}$, then the angular displacement from one into the other is achieved by a translation in $\R$ sending $\widetilde{\theta_1}$ to $\widetilde{\theta_2}$. In particular, $Rot[\gamma]$ is regular if and only if ${\rm bp}[\gamma]$ is a regular point.

If $M$ is a tame translation surface, then $IrrRot(M)_0=\emptyset$ and the cardinality of $IrrRot(M)$ is equal to the cardinality of the set of singularities $\Sing(M)$. If $M$ is wild then the cardinality of $\Sing(M)$ is not necessarily related to the cardinality of $IrrRot(M)$: we will see in Example~\ref{exa:IcicledSurface} that there exist a translation surface with only two (wild) singularities but which has an uncountable number of rotational components. In more general terms, given  that translation surfaces are second-countable (\emph{i.e.} the topology has a countable base) if $\Sing(M)$ is countable then $IrrRot(M)_+$ is at most countable. In particular, if  $\Sing(M)$ is countable and $IrrRot(M)$ is uncountable then the set of rotational components with total angle zero $IrrRot(M)_0$ is uncountable.


Given that the map $\widetilde{\rm dir}$ defined above is injective, one can endow any rotational component with the standard topology of the real line. Moreover, using this map it is not difficult to prove that
every rotational component admits a $(G,X)$-structure, where $X=\R$ and $G=2\pi\Z$ acts by translations,
possibly with non-empty boundary, see \cite{BowmanValdez13} and Appendix \ref{Appendix:GXStructures}. However this $(G,X)$-structure is somehow imposed artificially as every rotational component carries a natural subspace topology coming from the limit topology of $T\widehat{M}$. As we see in the following example, these topologies in general do not agree. As a matter of fact even in simple cases rotational components with their subspace topology are not metrizable spaces.


\begin{tcbexample}{}{NotRegularSpace}
\label{exa::NotRegularSpace}
The following description relies on Figure~\ref{fig:StackBoxes} and was first described by Bowman~\cite{Bowman12} and Randecker~\cite{Randecker16}. Let $W=\{w_n\}_{n\geq 0}$ be a strictly monotonically decreasing sequence of positive real numbers that converges to 0. For each $n\in \N$ consider the rectangle $R_n$ of height 1 and width $w_n$. For each $n\geq 1$ glue the lower edge of $R_{n+1}$ to the left part of the upper side of $R_n$ and the vertical sides of $R_n$ together. The result looks like the pile of rectangles shown in the figure. For each $n$ there is still a segment on the top right of each box $R_n$ to be identified. To do this we subdivide the lower edge of the bottom box $R_0$ into segments of size $w_{n+1}-w_n$ starting from $w_0 - w_1$ on the right side. After doing this operation, we have vertically aligned horizontal segments that we identify by translations (they carry the labels $A_n$ on Figure~\ref{fig:StackBoxes}). The result of these gluings is a wild translation surface $M$ having only one singularity $z$. Moreover, $T_z^1\widehat{M}$ is formed by just one rotational component that can be (set theoretically) identified using the map (\ref{eq:InjectiveLift}) with $(0,\infty)$. The convention in this identification is that the linear approach in the vertical direction to the point $z$ (in the lower left corner of Figure~\ref{fig:StackBoxes}) corresponds to $\pi$. Using this correspondence we can write
\begin{equation}
	\label{EQ:closed_set_linear_approaches}
T_z^1\widehat{M}=\{[\gamma_s]\hspace{1mm}|\hspace{1mm}s\in(0,\infty)\}.
\end{equation}

The claim is that $T_z^1\widehat{M}$ with the topology induced from $T^1\widehat{M}$ is not a regular space. That is, there exist a closed subset $F$ and a point $[\gamma]$ in $T_z^1\widehat{M}$ which cannot be separated by neighbourhoods. We define first the closed subset. For every fixed $\eps>0$ let $\rho_\eps:\mathcal{L}^\eps(M)\to T^1\widehat{M}$ be the natural projection. By definition $\rho_\eps(\mathcal{L}^\eps)$ is an open subset. For every fixed $\eps'>0$ let
$$
F_{\eps'}:=T_z^1\widehat{M} \setminus \bigcup_{\eps>\eps'}\rho_\eps(\mathcal{L}^\eps)
$$
The set $F_{\eps'}$ is closed and is formed by all linear approaches $[\gamma]$ to $z$ for which the length of any representative is bounded above by $\eps'$. The rest of the proof of our claim is left to the reader in the following:

\begin{tcbexercise}{}{NotRegularSpace}
Let $[\gamma_\pi]$ be the linear approach to $z$ in the horizontal direction defined according to (\ref{EQ:closed_set_linear_approaches}). Show that if $\eps'>0$ is small enough, every neighbourhood of $[\gamma_\pi]$ and $F_{\eps'}$ intersect. Conclude that $T_z^1\widehat{M}$ is not a regular space and hence not metrizable. \emph{Hint}: show first that for small $\eps'>0$ the set $F_{\eps'}$ is formed by the union of horizontal saddle connections of the form $\sqcup_{n\geq k} A_n$ for a sufficiently large $k$ (see Figure~\ref{fig:StackBoxes})).
\end{tcbexercise}
\end{tcbexample}


\begin{figure}[!ht]
\begin{center}
\includegraphics[scale=1]{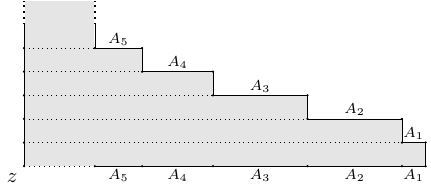}
\end{center}
\caption{The stack of boxes surface.}
\label{fig:StackBoxes}
\end{figure}

In what follows we illustrate with a series of examples how simple or complicated the spaces $T_z^1\widehat{M}$ (and their rotational components) can be when $z$ is in the singular locus $\Sing(M)$. Except for the double parabola and baker's surface, all examples are due to A. Randecker. See \cite{Randecker16} for details.


\begin{tcbexample}{Double parabola}{DoubleParabola}
\label{Ex:doubleparabola}
As we see in the following lines, rotational components can be singletons. Let $\pm I_n$ be a family of segments in the $xy$-plane whose endpoints are $(\pm 2^n,2^{2n})$ and $(\pm 2^{n+1},2^{2(n+1)})$,
$n\in\Z$. Let $\pm J_n$ be the family of segments whose endpoints are
$(\pm 2^n,-2^{2n})$ and $(\pm 2^{n+1},-2^{2(n+1)})$, $n\in\Z$. Let $P_-$ be closure of the connected component of
\begin{equation}
	\label{DoubleParabola}
\R^2\setminus\{\{\pm I_n\}_{n\in\Z}\cup(0,0)\cup\{\pm J_n\}_{n\in\Z}\}
\end{equation}
containing the negative $x$-axis. Analogously, let $P_+$ be the closure of the connected component of
(\ref{DoubleParabola}) containing the positive $x$-axis. By construction $\partial P_-=\{-I_n\}_{n\in\Z}\cup(0,0)\cup \{-J_n\}_{n\in\Z}$ and
$\partial P_+=\{I_n\}_{n\in\Z}\cup(0,0)\cup \{J_n\}_{n\in\Z}$. Remove all vertices
(and the origin) from $P_-$ and $P_+$ and identify this two disjoint domains along parallel sides of the same length using translations. This produces a translation surface $M$ that we call the \emph{double parabola}.

\begin{figure}[H]
\begin{center}
\includegraphics[scale=1]{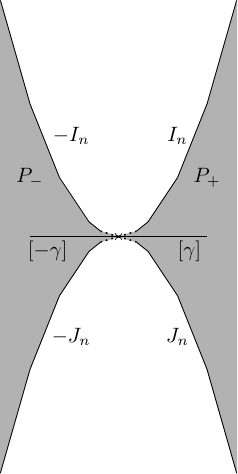}
\end{center}
\caption{The double parabola. The linear approaches $[\gamma]$ and $[-\gamma]$ are the only elements in the rotational components they define.}
\label{fig:StackBoxes}
\end{figure}
\begin{tcbexercise}{}{}
Let $M$ be the double parabola defined above.
\begin{enumerate}
\item Show that $\Sing(\widehat{M})$ is just one wild singularity $z_\infty$. \emph{Hint}: $\lim_{n\to-\infty}|\pm I_n|=\lim_{n\to-\infty}|\pm J_n|=0$ and these segments accumulate to the origin.
\item Determine all rotational components in $T^1_{z_\infty}\widehat{M}$. In particular, show that each rotational component defined by
$\pm\gamma(t)=(\pm t,0),$ $t\in(0,1)$ consists of only one point.
\end{enumerate}
\end{tcbexercise}
\end{tcbexample}

\begin{figure}[H]
\begin{center}
\includegraphics{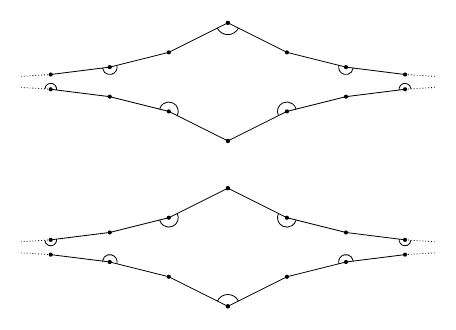}
\end{center}
\caption{The space of linear approaches to the wild singularity $x_\infty$ of the exponential surface.}
\label{fig:ExpoSurface2}
\end{figure}

\begin{tcbexample}{Exponential surface}{ExponentialSurface}
There exist finite-area translation surfaces with only one wild singularity $z_\infty$ and $T^1_{z_\infty}\widehat{M}$ is formed by two rotational components parametrized by $\R$. Indeed,  consider the infinite polygon $P$ depicted in figure \ref{fig:ExpoSurface1}: its vertices are $\{(n,\pm2^{-|n|})\hspace{1mm}|\hspace{1mm}n\in\Z\}$ and its interior contains the real axis. The sides of $P$ can be grouped into pairs of parallel sides; these are labeled in the figure with the same labels. The translation surface $M$ obtained by gluing the parallel sides of $P$ by translation is called the \emph{exponential surface}. By construction the metric completion $\widehat{M}$ is obtained by adding to $M$ all vertices of $P$, and by the way we glued the sides of $P$ to obtain $M$, these vertices merge \emph{a priori} to two points $z_1,z_2\in\widehat{S}$. However, since every saddle connection in the vertical direction joins $z_1$ to $z_2$ and the length of these accumulates to zero we conclude that $z_1=z_2$ is the only singularity in the metric completion $\widehat{M}$. Let us denote this singularity by $z_\infty$. We claim that $z_\infty$ is a wild singularity. Indeed, given that the total area of $S$ is finite, $z_\infty$ cannot be an infinite angle singularity and it cannot be a finite angle singularity since the number of saddle connections emerging from it is infinite. The space of linear approaches $T^1_{z_\infty}\widehat{M}$ is formed by two rotational components which are isometric to $\R$, as depicted in Figure~\ref{fig:ExpoSurface2}.

\begin{figure}[H]
\begin{center}
\includegraphics{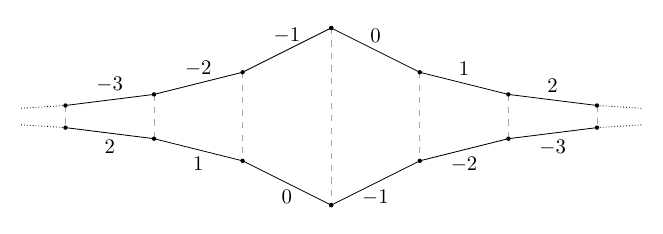}
\end{center}
\caption{The exponential surface.}
\label{fig:ExpoSurface1}
\end{figure}
\end{tcbexample}
\begin{tcbexample}{Baker Surface}{RotationalComponentsBakersSurface}
As we saw in~\ref{ssec:BakerBowmanArnouxYoccoz}, baker's surface $B_\alpha$ has only one wild singularity $z_\infty\in\widehat{B_\alpha}$. Chamanara's remark ``[g]eometrically, the surface spirals infinitely
many times around this point'' can now be rephrased in the terminology we just introduced as \emph{the space of linear approaches $T^1_{z_\infty}\widehat{M}$ contains an unbounded rotational component}. As a matter of fact we now have the tools to be completely precise. For the sake of clarity we set $\alpha=\frac{1}{2}$ and we refer henceforth to figure \ref{fig:chamanarainfinitegenus} in the preceding section. Suppose that the intersection of the diagonals in the unit square defining baker's surface is
the origin and the corners {\Large\textbf{a}} and {\Large\textbf{c}}
have coordinates $(-\frac{1}{2},\frac{1}{2})$ and $(\frac{1}{2},-\frac{1}{2})$ respectively. In this coordinates the geodesic segments $\gamma_1(t):=(1-t)( -\frac{1}{2},\frac{1}{2})$ and $\gamma_2(t):=-\gamma_1(t)$, $t\in(0,\varepsilon)$,  are linear approaches to $z_\infty$ that define two bi-infinite rotational components $\overline{[\gamma_1]}$ and $\overline{[\gamma_2]}$, each isometric to $\R$. These rotational components are drawn in red and blue in figure \ref{F:cham}  in section \ref{ssec:BakerBowmanArnouxYoccoz}. On the other hand,
$\eta_1(t)=(1-t)(\frac{1}{2},\frac{1}{2})$ and $\eta_2(t)=-\eta_1(t)$, define two bounded rotational
components whose total angle is $\pi/4$. The boundary of these rotational components is formed by the
horizontal and vertical saddle connections labeled $A_i$ and $B_i$, $i\in\N$ in the same figure. In the next section we see how the action of the Veech group on $T^1_{z_\infty}\widehat{M}$ can be used to describe all rotational components, see Lemma~\ref{lem:RotCompChamanara} in Section \ref{ssec:VeechGroupBakersSurface}.
\end{tcbexample}


\begin{tcbexample}{Icicled surface}{IcicledSurface}
One of the main issues when studying the family of translation flows $\{F_\theta^t\}_\theta$ on $M$ is to determine whether they are well defined. Namely, for each angle $\theta$ we would like a flow defined on a set of full measure in $M$. The only obstruction for these flows to be defined for all times are points in $\Sing(M)$. Intuition tells us that if $\Sing(M)$ is a singleton, then the translation flow should be defined for all directions on sets of full measure. In what follows we present an example that shows that this intuition is wrong: there exist a finite area translation surface $M$ of infinite type with only two singularities such that the \emph{only} direction on which the translation flow $F^t_\theta$ is defined on a set of full measure is the horizontal direction. In other words, for infinite-type translation surfaces of finite area, having a finite set of singularities does not imply that the translation flow is defined in a set of full measure for almost every direction.

The example we present is known as the \emph{icicled surface}. The construction of this example starts with the closed rectangle $R=[0,1]\times[0,2]$. Remove from $R$ the top and bottom sides $[0,1]\times\{2\}$, $[0,1]\times\{0\}$ and identify vertical sides using a translation. For every $n\geq 1$ and odd $i\in\{1,\ldots,2^n-1\}$ consider the vertical segments $I_{i,n}^{top}$, $I_{i,n}^{bot}$  inside $R$ of length $2^{-n}$ with one extremity at $(\frac{i}{2^n},2)$ and $(\frac{i}{2^n},0)$ respectively. These are the \emph{icicles} and are depicted in Figure~\ref{fig:IcicledSurface1}.

\begin{figure}[H]
\begin{center}\includegraphics[scale=0.8]{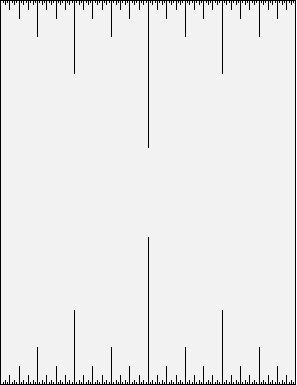}\end{center}
\caption{The icicles in the icicled surface.\\}
\label{fig:IcicledSurface1}
\end{figure}

We now define gluings on icicles with one extremity in $[0,\frac{1}{2}]\times\{2\}$:
\begin{enumerate}
\item First subdivide each side of the icicle $I^{top}_{1,1}$ at $(\frac{1}{2},2)$ following the geometric progression $\frac{1}{2^k}$, $k>1$ from bottom to top. This creates an infinite subdivision of the icicle $I^{top}_{1,1}$ into segments of length $\frac{1}{2^k}$.
\item Glue the left side of the lower half of the icicle $I^{top}_{1,1}$ to the right side of the icicle $I^{top}_{1,2}$ (which has one extremity at $(\frac{1}{4},1)$).
\item For every $n>2$ and ever odd $i\in\{3,\ldots,2^{n-1}-1\}$ such that $I^{top}_{i,n}$ has one extremity in $[0,\frac{1}{2}]\times\{1\}$, glue the left side of $I^{top}_{i,n}$ to the right side of $I^{top}_{i-2,n}$.
\item For every $n>2$, the left side of the icicle at $(\frac{1}{2^n-1},2)$ is cut into two segments of the same length. The lower segment obtained by this subdivision is glued to the right side of the icicle at $(\frac{2^{n-1}-1}{2^n},2)$ and the upper segment to the right side of the only segment at $(\frac{1}{2},2)$ that has the same length.
\end{enumerate}
The gluings described above are sketched in Figure~\ref{fig:IcicledSurface2}. The gluings on the icicles with one extremity on $[\frac{1}{2},1]\times\{2\}$ are defined in an analogous way, the only difference is that we need to change \emph{left} for \emph{right} in all instructions above. The gluings for icicles with one extremity on $[0,1]\times\{0\}$ are
just a mirror image (w.r.t. the real axis) of the gluings defined above. The result of these gluings is the \emph{icicled surface}.


\begin{figure}[H]
\begin{center}\includegraphics[scale=1.2]{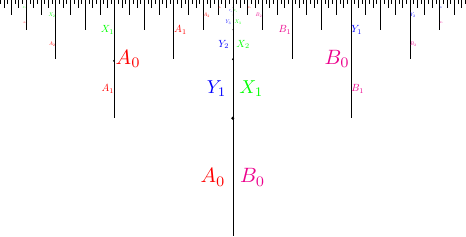}\end{center}
\caption{Gluings of the icicles to produce the icicled surface.}
\label{fig:IcicledSurface2}
\end{figure}

\begin{tcbexercise}{}{}
Prove that the icicled surface $M$ has only two wild singularities $\{\sigma_{top},\sigma_{bot}\}$. \emph{Hint}: by the way the gluings were defined, all the tips of the icicles $I^{top}_{i,n}$ are identified into a point that we denote by $\sigma_{top}$. Show that every non-dyadic point in $[0,1]\times\{2\}$, \ie, which is not an extremity of an icicle, is the limit of a sequence of tips of icicles. In other words, every non-dyadic point in $[0,1]\times\{2\}$ is identified with $\sigma_{top}$ in the metric completion $\widehat{M}$. The same argument works for dyadic points in $[0,1]\times\{0\}$

\end{tcbexercise}

\begin{tcbexercise}{}{}
Prove that the space of linear approaches $T^1_{\sigma_{top}}\widehat{M}$ contains \emph{uncountably} many rotational components which are singletons, and at least two rotational components of total angle $\frac{\pi}{4}$. Does it contain unbounded rotational components (\ie, of infinite total angle)?
\end{tcbexercise}

Note that for every flat point $z$ in the icicled surface the geodesic flow $F_\theta^t(x)$ is not defined in the future and in the past whenever $\theta\neq 0$. Indeed, by the way the gluings defining the icicled surface were performed, every such trajectory reaches the top or the bottom of the rectangle $R$ in finite time. On the other hand, for $\theta=0$ the geodesic flow $F_\theta^t$ decomposes the icicled surface into an infinite family of cylinders and saddle connections.
\end{tcbexample}

\subsubsection{Proof of Lemma~\ref{lem:triangulation}}
	\label{SSEC:Debt}
We finish this Chapter paying a debt we had from Chapter~\ref{CH:INTRODUCTION}, namely, let us prove using the language we developed
Lemma~\ref{lem:triangulation}. Recall that this lemma states that a finite area simply connected translation
surface with finitely many segments as boundary can be triangulated.
\begin{proof}[Proof of Lemma~\ref{lem:triangulation}]
\label{proof:lem:Triangulation}
Let $M$ be a simply connected translation surface of finite area whose boundary consists
of finitely many Euclidean segments.

We consider the set $L$ of straight line segments whose ends are either singularity or vertex on the
boundary. Now consider the subsets of segments in $L$ whose elements have pairwise disjoint interior.
Such subsets are ordered by inclusion and we claim that the maximal such subsets correspond exactly
to edges of a triangulation. As we want to proceed by induction, we just need to show that if the
surface $M$ is not already a triangle the set $L$ is not empty.

It is well known that an convex polygon can be triangulated. This corresponds when $M$ has no
conical singularity and each angular sector of points on the boundary has angle $\leq \pi$. To
conclude it is enough to prove the following claim.

\begin{tcblemma}{}{FiniteAreaAndSaddleConnections}
Let $M$ be a translation surface of finite area whose boundary consists of
finitely many Euclidean segments. Let $A$ be an angular sector in $T \widehat{M}$
based at a point $p$ whose angle is greater than or equal to $\pi$. We claim that among the germs of
straight-line segment in $A$ there is a saddle connection that joins to a conical singularity
in $M$ or a point on the boundary whose visual sector is not $\pi$.
\end{tcblemma}

\begin{proof}[Proof of Lemma~\ref{lem:FiniteAreaAndSaddleConnections}]
Let $A$ be an angular sector with angle $> \pi$. There exists a germ of a map from the half-plane
into $A$. Either we discover a singularity, and we are done. Otherwise the boundary of this piece of
half-plane is made of a union of segments. Since the surface $M$ is simply connected, the map is an embedding and
because the area is finite, it needs to have corners.
\end{proof}

\end{proof}

\chapter{Symmetries}
\label{ch:Symmetries}

We introduced in Section~\ref{ssec:Structures} the group of affine homeomorphisms of a translation surface as well as the subgroups of isometries and translations. These groups are the center of attention of this Chapter. We start in Section~\ref{sec:AnalyticAutomorphisms} with the group of isometries. Because they are conformal automorphisms, we can rely on powerful tools from hyperbolic geometry. This allows us to provide a complete classification of translation surfaces having a group of translations that does not act properly discontinuously.

Recall from Definition~\ref{def:VeechGroup} that the Veech group of a translation surface $M$ is the subgroup of $\GL(2,\R)$ consisting of the linear parts of the affine homeomorphisms. The Veech group of a finite-type translation surface is a discrete subgroup of $\SL(2,\R)$. Moreover, a classical result by W.~Veech shows that in the (extremely rare) case that the Veech group is a lattice, the translation flows behave similarly to the translation flows on a torus. For a precise statement see Theorem~\ref{thm:VeechDichotomy}. In Section~\ref{sec:EveryGroupIsVeech}, we prove that any countable subgroup of $\GL^+(2,\R)$ is the Veech group of some infinite-area translation surface with trivial dynamics : any non-singular orbit in any direction escapes to infinity. In particular Veech groups for infinite-area surfaces do not immediately provide consequences for the dynamics. However, Veech groups turn out to be useful to study translation flows of Hooper-Thurson-Veech surfaces and infinite coverings of finite-type surfaces as we will see in~\cite{DHV2}.

As in Chapter~\ref{ch:TopologyGeometry}, infinite coverings of finite-type surfaces play a special role. In Section~\ref{sec:AffineGroupsCoverings}, we explore the relationship between the affine group of a covering and the affine group of its base using the so-called lifting criterion from algebraic topology. It becomes particularly concrete in the situation of Abelian coverings and we apply it to study Veech groups of Abelian coverings of genus 2 translation surfaces in Section~\ref{ssec:HooperWeissTheorem} and the wind-tree models in Section~\ref{ssec:SymmetriesWindTree}.

Then, we present in Section~\ref{sec:HooperThurstonVeechConstruction} an extension of the Thurston-Veech construction from Section~\ref{ssec:ThurstonVeechConstructionsIntro} to the realm of infinite-type surfaces following P.~Hooper's work~\cite{Hooper-infinite_Thurston_Veech}. The main novelty of the infinite-type construction is that the Perron-Frobenius theorem does not apply and there is a lot of freedom in the choice of the eigenvalue and positive eigenvector to build the translation structure.

We end the chapter with the computation of the Veech groups of the baker and Bowman-Arnoux-Yoccoz surfaces from Section~\ref{ssec:BakerBowmanArnouxYoccoz}. These are some of the few non-trivial examples on which Veech groups of finite-area infinite-type surfaces are known.

\section{Analytic affine automorphisms: isometries and translations}
\label{sec:AnalyticAutomorphisms}
Through this section, we adopt the analytic point of view for translation
surfaces. That is, a translation surface $M$ will be thought as a couple $(X,
\omega)$ where $X$ is a Riemann surface and $\omega$ a non-zero holomorphic
one-form (see Definition~\ref{def:TranslationSurfaceAnalytic} and
Theorem~\ref{thm:EquivalencesInDefinitions}).
We denote the group of analytic automorphisms of the
Riemann surface $X$ by $\Aut(X)$ . The following elementary lemma, whose proof is left to the
reader, characterizes the intersection between $\Aff(X,\omega)$ and $\Aut(X)$.
\begin{tcblemma}{}{}
An affine homeomorphism $f \in \Aff(X, \omega)$ is analytic if and only if $D(\phi)$ is a similitude, i.e. a matrix of the form $\begin{psmallmatrix} a&-b\\ b&a \end{psmallmatrix}$ that corresponds to multiplication by $a+ \sqrt{-1} b$ in $\C$ under the identification $\R^2 = \C$.
\end{tcblemma}
In particular orientable isometries, and hence translations, are automorphisms of the underlying Riemann surface: $\Tr(X,\omega) \subset \Isom^+(X,\omega) \subset \Aut(X)$.

A classical result of A. Hurwitz gives upper bounds on $\Aut(X)$ when $X$ is a compact Riemann surface.
\begin{tcbtheorem}{Hurwitz's automorphism theorem~\cite{Hurwitz1893}}{HurwitzBoundAutomorphisms}
Let $X$ be a compact Riemann surface of genus $g \geq 2$. Then the group $\Aut(X)$ of analytic automorphism of $X$ is finite and $|\Aut(X)| \leq 84(g-1)$.
\end{tcbtheorem}
In genus $0$ or $1$ the automorphism groups are infinite.
The following result provides an analogue of Hurwitz's theorem for the subgroup of translations $\Tr(X,\omega)$.
\begin{tcbtheorem}{\cite{SchlagePuchtaWeitzeSchmithuesen-translations}}{BoundTranslationAutomorphisms}
Let $(X,\omega)$ be a compact translation surface of genus $g \geq 2$. Then $|\Tr(X,\omega)| \leq 4g-4$.

Moreover, the order is exactly $4g - 4$ if and only if $(X,\omega)$ is a normal square-tiled surface\footnote{That is, a square-tiled surface $p: M \to \T^2$ defined by a normal subgroup of $\pi_1(\T^2 \backslash \{0\})$, see Example~\ref{exa:SquareTiledSurfaces}.}  in the principal stratum $\cH(1,\ldots,1)$, that is,  when $\omega$ has $2g-2$ simple zeroes. And this could happen only if $g$ is odd or $g-1$ is divisible by 3.
\end{tcbtheorem}

As a consequence of Hurwitz theorem, for every compact Riemann surface $X$ of genus $g \geq 2$, the group $\Aut(X)$ acts properly discontinuously on $X$ (see Definition~\ref{def:ProperlyDiscontinuousAction}). The following result provides the list of surfaces for which automorphisms do not act properly discontinuously (or equivalently, for which the group of automorphism is not countable).
\begin{tcbtheorem}{}{thm:BigAnalyticAut}
Let $X$ be a Riemann surface such that its group $\Aut(X)$ of analytic
automorphism does not act properly discontinuously on $X$.
Then $X$ is conformally equivalent to one of the following surfaces:
\begin{compactitem}
\item the Riemann sphere $\hat{\C}$ for which $\Aut(\hat{\C}) = \SL(2,\C)$,
\item the plane $\C$ for which $\Aut(\C) = \C^* \ltimes \C$,
\item the half-plane $\H^2$ for which $\Aut(\H^2) = \PSL(2,\R)$,
\item an annulus $\{z \in \C; r < |z| < R\}$ with $0 < r < R \leq \infty$ for which $\Aut(X)$ contains a subgroup isomorphic to $\SS^1$,
\item the punctured plane $\C^*$ for which $\Aut(\C^*)=\C^*$ or
\item a torus $\R^2 / \Lambda$, in which case $\Aut(\R^2/\Lambda)$ contains the group of translations $\R^2/\Lambda$.
\end{compactitem}
\label{thm:BigAnalyticAut}
\end{tcbtheorem}
The main ingredient in the proof is the Poincar\'e-Koebe's uniformization theorem.
\begin{proof}
Let $X$ be a Riemann surface, $\widetilde{X}$ its universal cover and $G = \Aut(\widetilde{X})$.
Let $H \subset G$ be the subgroup acting freely and properly discontinuously on $\widetilde{X}$ so that $X = \widetilde{X} / H$. In particular the action of $H$ on $\widetilde{X}$ has no fixed points.

Let $N_G(H) = \{g \in G\ | \ g H = H g\}$ be the normalizer of $H$ in $G$. We claim that $\Aut(X)$ is in bijection with $N_G(H) / H$. Indeed, every element $\varphi\in \Aut(X)$ can be lifted to an element $\widetilde{\varphi}$ which lives in $N_G(H)$ for $\widetilde{\varphi}$ must send $H$-orbits in the universal cover to $H$-orbits. On the other hand, every element $\widetilde{\psi}\in N_G(H)$ descends to an analytic automorphism $\psi\in\Aut(X)$. Moreover $\psi=Id_X$ if and only if $\widetilde{\psi}$ leaves each $H$-orbits in the universal cover invariant, that is, if $\widetilde{\psi}\in H$.

We want to list the possible group $H$ for which $\Aut(X) \simeq N_G(H) / H$ is non-discrete.
Using the uniformization theorem, we proceed by cases considering the three different universal covers that $X$ can have:
\begin{itemize}
\item $\widetilde{X} = \SS^2$. Given that every non-trivial analytic automorphism of the Riemann sphere has a fixed point and $H$ cannot have fixed points we have that
 $X = \widetilde{X}$.

\item $\widetilde{X} = \C$. The only analytic automorphisms of $\C$ without fixed points are the translations.
Hence $H = \{1\}$, $H = \Z u$ or $H = \Z u \oplus \Z v$, where $u$ and $v$ are complex numbers different from zero such that $u/v\notin\R$.
The resulting quotients are, respectively, the plane, a cylinder conformally equivalent to $\C^*$ and a torus.

\item $\widetilde{X} = \H^2$.
Let us first remark that the group $N_G(H)$ is closed in $\PSL(2,\R)$. Hence it is either discrete or contains a one-parameter subgroup.
The one-parameter subgroups of $\PSL(2,\R)$ are well known and are of one of the following type:
\begin{compactitem}
\item parabolic, i.e. conjugate to $P = \{h_t\}$ where $h_t = \begin{pmatrix} 1 & t \\ 0 & 1\end{pmatrix}$, with infinitesimal generator
$\displaystyle \mathfrak{h} := \frac{d}{dt}|_{t=0} h_t = \begin{pmatrix}0&1\\0&0\end{pmatrix}$,
\item hyperbolic, i.e. conjugate to $A = \{g_t\}$ where $g_t = \begin{pmatrix}e^t & 0 \\0 & e^{-t}\end{pmatrix}$, with infinitesimal generator
$\displaystyle \mathfrak{a} := \frac{d}{dt}_{t=0} g_t = \begin{pmatrix}1&0\\0&-1\end{pmatrix}$,
\item rotational, i.e. conjugate to $R = \{r_\theta\}$ where $r_\theta = \begin{pmatrix} \cos(\theta) & -\sin(\theta) \\ \sin(\theta) & \cos(\theta) \end{pmatrix}$ with infinitesimal generator $\displaystyle \mathfrak{r} := \frac{d}{d\theta}|_{\theta=0} r_\theta = \begin{pmatrix}0&1\\-1&0\end{pmatrix}$.
\end{compactitem}

From now on we assume that $\widetilde{X} = \H^2$, that $H \not= \{1\}$
(that corresponds to $X = \H^2$) and that the normalizer $N_G(H)$ contains a
one-parameter subgroup $\{s_t\}_{t \in \R}$. For all $m \in H$ and all $t \in \R$,
the conjugate $s_t m s_{-t}$ belongs to $H$. As $H$ is discrete, for $t$ small enough
$s_t m s_{-t} = m$ or equivalently $s_t m = m s_t$. Taking the derivative at $t=0$ we obtain that
\begin{equation}
\label{eq:InfinitesimalRelation}
\mathfrak{s} m = m \mathfrak{s}
\end{equation}
where $\mathfrak{s} = \frac{d}{dt}|_{t=0} s_t$ is the infinitesimal generator of the one-parameter subgroup.
We compute the equations obtained in the three cases of subgroups.  Let $m = \begin{pmatrix}a & b \\c & d\end{pmatrix}$
be a matrix in $H$. Up to conjugation we can assume that $\{s_t\}$ is one of $P$, $A$ or $R$ defined above.

In the parabolic case, equation~\eqref{eq:InfinitesimalRelation} writes
$\mathfrak{h} m = m \mathfrak{h}$, that is
\[
\begin{pmatrix}c&d\\0&0\end{pmatrix}
=
\begin{pmatrix}0&a\\0&c\end{pmatrix}.
\]
Hence, $c=0$ and $a=d$. This implies that either $m$ acts like the identity on
$\mathbb{H}^2$ or  $m \in P$. In the former case we conclude that $X=\H^2$. In
the later that $H = \langle h_{t_0} \rangle$ for some $t_0$ (because $H$ is a
discrete subgroup of $\PSL(2,\R)$). In this case, $X$ is conformally equivalent
to the punctured unit disc (which is conformally equivalent to an infinite
annulus $\{z\in\C; r<|z|\leq\infty\}$). Moreover, in this case $P/H$ is isomorphic to $\SS^1$.

Now, assume that $\{s_t\} = A$. Then equation~\eqref{eq:InfinitesimalRelation}
gives $\mathfrak{a} m = m \mathfrak{a}$, that is
\[
\begin{pmatrix}a&b\\-c&-d\end{pmatrix}
=
\begin{pmatrix}a&-b\\c&-d\end{pmatrix}
\]
and hence $b = c = 0$. In other words $m \in A$. Hence $H = \langle g_{t_0}
\rangle$ for some $t_0$ and $X$ is conformally equivalent to a finite annulus
$\{z \in \C; r < z < R\}$ with $0 < r < R < \infty$. Moreover, in this case
$P/H$ is isomorphic to $\SS^1$.

Finally, let us consider the case $\{s_t\} = R$.
Equation~\eqref{eq:InfinitesimalRelation} gives $\mathfrak{r} m = m \mathfrak{r}$
which writes
\[
\begin{pmatrix}c&d\\-a&-b\end{pmatrix}
=
\begin{pmatrix}-b&a\\-d&c\end{pmatrix}.
\]
Hence, $a = d$ and $b = -c$. In other words $m \in R$. Given that every
non-trivial element of $R$ has a fixed point in $\H^2$ and the action of $H$ on
$\widetilde{X}$ has no fixed points we conclude that $H = \{1\}$, which was excluded.
\end{itemize}
\end{proof}

We will now use Theorem~\ref{thm:BigAnalyticAut} to classify translation
surfaces with uncountable translation group.
\begin{tcbtheorem}{}{TranslationSurfacesBigAut}
Let $(X,\omega)$ be a
translation surface whose translation group does not act properly
discontinuously. Then there exists $\lambda \in \C^*$ such that $(X, \lambda \omega)$
is isomorphic as a translation surface to one of the following:
\begin{compactenum}
\item the plane $(\C,dz)$ or the half plane $(\H^2,dz)$, \label{case:plane}
\item a torus $(\C / (\Z u \oplus \Z v),\pi_*dz)$, \label{case:flat_torus} where $\pi:\C\to\C / (\Z u \oplus \Z v)$ is the natural projection.
\item a strip $(\{z \in \C\ | \ -\frac{\pi}{2} < \Im(z) < \frac{\pi}{2}\},dz)$, \label{case:strip} isomorphic to  $(\H^2,\frac{dz}{z})$,\label{case:half_plane}
\item a finite cylinder $(\{z \in \C\ | \ -\frac{\pi}{2} < \Im(z) < \frac{\pi}{2}\} / (z \sim z+t),dz)$, with $0<t<\infty$.\label{case:anulus}
\item an infinite cylinder $(\{z \in \C\ | \ -\frac{\pi}{2} \leq \Im(z) \leq \frac{\pi}{2}, \ R<\Re(z)\} / (z \sim z+i\pi),dz)$, where $-\infty\leq R$.
\label{case:InfiniteCylinder}
\end{compactenum}
\end{tcbtheorem}

\begin{proof}
The strategy of the proof is the following: we determine first by direct calculation all simply-connected translation surfaces $(X,\omega)$ for which $\Trans(X,\omega)<\Aut(X)$ does not act properly discontinuously. Then we use the list of  Theorem \ref{thm:BigAnalyticAut} to get all possible (non-simply connected) quotients.

The sphere is excluded since it supports no Abelian differential. One-parameter subgroups of $\Aut(\C)$ are of the form $z\to z+tw$ or $z\to e^{t \alpha }z$ for some fixed $\alpha, w\in\C^*$ and $t\in \R$. Let $(\C,\omega=f(z)dz)$, with $f$ an entire function, be an Abelian differential. This form has to be invariant under the one-parameter subgroups of $\Aut(\C)$ mentioned before. The first possibility yields:

\begin{equation}
\label{eq:TranslationInvariance}
f(z) dz = f(z+tw) dz.
\end{equation}
As this holds for arbitrarily small $t$, necessarily $f$ is constant so that $\omega = \alpha dz$ for $\alpha\in\C^*$. For the second possibility we get
\begin{equation}
\label{eq:LoxodromicInvariance}
f(z) dz = f(e^{t \alpha} z) e^{t\alpha} dz.
\end{equation}
Taking derivatives with respect to $t$ and setting $t=0$ we obtain
\[
0 = z f'(z) + f(z).
\]
The only non-identically zero solutions to this equation are meromorphic functions  $f(z) = c/ z$ with $c \in \C^*$. This case yields forms that are not Abelian differentials, hence it is excluded.

We now determine all translation surfaces $(\H^2, f(z) dz)$ which are invariant under one-parameter subgroups of $\Aut(\H^2)$. As in the proof of Theorem~\ref{thm:BigAnalyticAut}, there are three one-parameter groups to consider: $P$ (parabolic), $A$ (hyperbolic) and $R$ (rotational). Forms invariant under the parabolic subgroup must be translation invariant, hence must satisfy equation~\eqref{eq:TranslationInvariance} and therefore are of the form $\alpha dz$, for some $\alpha\in\C^*$. Forms invariant under $R$ are also of the form $\alpha dz$, for some $\alpha\in\C^*$. This is best seen when performing a direct calculation of the invariance equation in the Poincar\'e disc. Forms invariant under the hyperbolic one-parameter group $A$ must satisfy equation~\eqref{eq:LoxodromicInvariance}, for $\alpha=2$. Hence in this case $\omega=\frac{\alpha dz}{z}$.

In summary: simply connected translation surfaces for which $\Trans(X,\omega)<\Aut(X)$ does not act properly discontinuously are all projectively isomorphic to either: $(\C, dz)$, $(\H^2, dz)$ or $(\H^2, \frac{ dz}{z})$.

We now use Theorem~\ref{thm:BigAnalyticAut}. We deduce that the only possible non-trivial (translation) quotients of $(\C,dz)$ are a bi-infinite cylinder  of the form $\{z \in \C\ | \ -\frac{\pi}{2} \leq \Im(z) \leq \frac{\pi}{2}, \ -\infty<\Re(z)\} / (z \sim z+i\pi)$ or a torus. On the other hand, the only possible non-trivial (translation) quotients of $(\H^2,dz)$ are a infinite half-cylinder  $\{z \in \C\ | \ -\frac{\pi}{2} \leq \Im(z) \leq \frac{\pi}{2}, \ R<\Re(z)\} / (z \sim z+i\pi)$, where $-\infty< R$. Finally, the non-trivial quotients of $(\H^2,\frac{dz}{z})$ are finite cylinders of the form
$\{z \in \C\ | \ -\frac{\pi}{2} < \Im(z) < \frac{\pi}{2}\} / (z \sim z+t)$, with $0<t<\infty$, see Exercise~\ref{exo:plane_infinite_band} below.
\end{proof}

\begin{tcbexercise}{}{plane_infinite_band}
Show that $(\H^2,\frac{dz}{z})$ is isometric to an infinite strip in $(\C,dz)$ and describe such strip. For each case in Theorem~\ref{thm:TranslationSurfacesBigAut} compute $\Trans(X,\omega)$ and $\Isom(X,\omega)$.
\end{tcbexercise}

\begin{tcbexercise}{}{}
Give examples that are not in the list of Corollary~\ref{thm:TranslationSurfacesBigAut} of surfaces for which $\Trans(X,\omega)$ acts properly discontinuously on $(X,\omega)$ but $\Isom(X,\omega)$ does not.
\end{tcbexercise}

Only two cases of Theorem~\ref{thm:TranslationSurfacesBigAut} have finite area. Hence we obtain the following:
\begin{tcbcorollary}{}{}
Except the case of tori and finite cylinders, any finite area translation surface has a finite group of translations.
\end{tcbcorollary}


\section{Every group is a Veech group}
\label{sec:EveryGroupIsVeech}

Recall that the Veech group $\Gamma(M)$ of a translation surface $M$ is the subgroup of $\GL(2,\R)$ formed by the derivatives of elements in $\Aff(M)$ and $\Gamma^+(M)=D(\Aff^+(M))$ are the derivatives of orientation-preserving elements. For a finite-type surface $M$ the Veech group $\Gamma(M)$ is Fuchsian, \ie a discrete subgroup of $\SL(2,\R)$, and is never co-compact. There are still many open questions regarding Veech groups of finite-type surfaces. For example, it is still unknown whether among such Fuchsian groups one can find one which is generated by a hyperbolic matrix, or if there exist one of the second kind (\emph{i.e.} whose limit set in $\partial\H^2$ is a Cantor set).
As we will see later in this Chapter, for finite-type translation surfaces, the Veech group plays an important role in the dynamics of the translation flows.

In this section we discuss two simple and na\"ive questions in the context of infinite-type translation surfaces: are there any restrictions for a Fuchsian group to be a Veech group? Do Veech groups play a role in the dynamics of the translation flows? As we explain in what follows, the answer to both of these questions is no.

\subsection{Veech groups of Loch Ness monsters}
\begin{tcbtheorem}{\cite{PrzytyckiSchmithuesenValdez11}}{AnyGroupIsVeech}
Let $A$ be an at most countable group and $F: A \to \GL^+(2,\R)$ a group morphism.
Then there exists a translation surface surface $M$ homeomorphic to the
Loch Ness monster and a group isomorphism
\[
\begin{array}{ccc}
A & \to & \Aff^+(M) \\
a & \mapsto & f_a
\end{array}
\]
such that $F(a) = D(f_a)$ where $D: \Aff^+(M) \to \GL^+(2,\R)$ is the derivative map.
\end{tcbtheorem}
Recall from Chapter~\ref{ch:TopologyGeometry} that the Loch Ness monster is the (homeomorphism class of) topological surface of infinite genus with only one end (see Definition~\ref{def:MonsterLadder}). Loosely speaking, this is the simplest infinite genus surface. It is the homeomorphism type of the translation surface $M(P)$ associated to an irrational polygon $P$ via unfolding (Theorem~\ref{thm:AnIrrationalBilliardIsLochNess} proven in Example~\ref{exa:IrrationalBilliardAndProof}) as well as the infinite staircase (Exercise~\ref{exo:InfiniteStaircaseIsLochNess}).

\begin{proof}
The proof we present is constructive and is inspired by ~\cite[Construction 4.6]{PrzytyckiSchmithuesenValdez11}.
We make the construction under the hypothesis that $A$ is not the trivial group.
In the case $A$ is the trivial group, one has to construct a translation
surface homeomorphic to the Loch Ness monster with trivial affine group.

The idea of the construction is to consider a normal ramified covering
$\pi:M \to \R^2$ with Deck group $A$. The construction is done so that
the deck group action is done via affine elements of $M$ with derivatives
prescribed by $F:A \to \GL^+(2,\R)$. Note that unless the image of $F$ is
identically $Id$, the map $\pi: M \to \R^2$ will not be a translation covering.
This construction will ensure that $A$ is contained in the affine group. In order
to prevent any additional affine automorphism we perform in a second step some
additional surgeries.

We will make use of the slit operation that was introduced in
Definition~\ref{def:GluingSlits} Section~\ref{ssec:PanovPlanes} to describe the Panov
planes and that we recall now. Let $M_1$ and $M_2$ be two translation surfaces and
$\gamma_1\subset M_1$ and $\gamma_2\subset M_2$ two closed parallel straight-line
segments of equal length  which do not encounter any singularities. Then cut the
translation surfaces $M_1$ and $M_2$ along $\gamma_1$ and $\gamma_2$ respectively
and glue these cuts back so that $M_1$ becomes attached to $M_2$ (see the left hand
side of Figure~\ref{fig:SlitConstruction}). The translation surface $M$ obtained this
way presents two singularities of angle $4\pi$ at the extremities of the cuts.

\begin{figure}[!ht]
\begin{center}
    \includegraphics[scale=0.7]{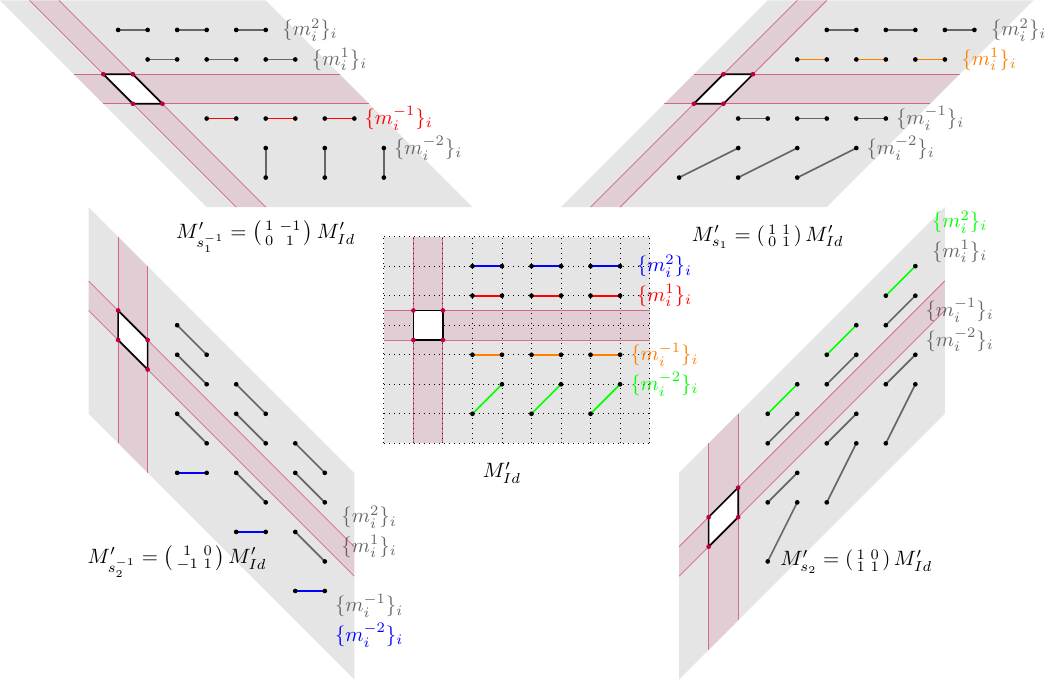}
\end{center}
\caption{An example of the ramified covering $\pi: M \to M_{Id}$ as in
the proof of Theorem~\ref{thm:AnyGroupIsVeech}. The group is $A = \SL(2,\Z)$ and
$F: \SL(2,\Z) \to \SL(2,\R)$ is the inclusion.}
\label{fig:EveryGroup}
\end{figure}

\emph{Step 1}. We consider $\R^2$ as a translation surface and construct a normal cover with
deck group $A$. Let us fix a generating set for $A$ with a given ordering
$U = \{a_1, a_2, \ldots\}$. For each $a_j$ in this generating set let
$(u_j, v_j) := F(a_j)^{-1}(1,0)$
and define
\[
x_j := |u_j| + 1 \qquad y_j := - (j + |v_1| + |v_2| + \ldots + |v_j|).
\]
For $i \in \N$, we consider the following segments in $\R^2$
\[
m^{j}_{i} := [(2i, j),\, (2i+1, j)]
\qquad
m^{-j}_{i} := [(i x_j, y_j),\, (i x_j + u_j, y_j + v_j)].
\]
See Figure~\ref{fig:EveryGroup} for an illustration.
Let us remark that these segments are pairwise disjoint.

For each $a \in A$, consider a copy $M'_a \simeq \R^2$ of the Euclidean plane.
For each $a \in A$, $a_j \in U$ and $i \in \N$
glue the segment $F(a) \cdot m_i^j$ in $M'_a$ to the segment
$F(a a_j) \cdot m^{-j}_i$ in $M'_{a\,a_j}$. Note that by definition of $m^j_i$
and $m^{-j}_i$ these segments are parallel and of equal length.
We denote by $M'$ the resulting translation surface. It is a covering
$\pi: M' \to \R^2$ ramified over the endpoints of the slits $m^j_i$ and
$m^{-j}_i$ and the copies $M'_a$ of the Euclidean plane constitute sheets
of this covering.

By construction the group $A$ acts by deck group on $M'_a$. This action is
affine and we denote by $f_a \in \Aff^+(M')$ the corresponding element. We
have $D(f_a) = F(a)$.

We now prove that $M'$ is homeomorphic to the Loch Ness monster.
The surface $M'$ is first of all connected. Next, the surfaces
$M'_a$ and $M'_{a\, a_j}$ are glued along two infinite
families of parallel slits. And these slits do not accumulate in
the metric completion of $M'_a$ and $M'_{a\,a_j}$.
For each slit $m^j_i$, the boundary of a sufficiently small regular
neighbourhood of $F(a) \cdot m^j_i$ in $M'_a$ defines a non-separating
essential curve in $M'_a$. Because the slit are disjoint, one can choose
these essential curves to be disjoint. Therefore $M'$ has infinite genus.
To see that $\Ends(M')$ has only one point we use Lemma~\ref{lem:CriterionInfiniteGenus}
from Section~\ref{sec:CoveringSpaces}. If $K\subset M'_a$ is compact then there exists
a ball $B$ in $\R^2 \simeq M'_{Id}$ such that $K$ is contained in
$\displaystyle \bigcup_{a \in A} s_a(B)$ where $s_a: \R^2 \to M'_a$ are
the sections to the sheets of $M'$. But since we add infinitely many slits,
the complement of this union is still connected. This concludes the fact
that $M'$ is a Loch Ness monster.

\emph{Step 2}.
The problem with the surface $M'$ is that it might admit more symmetries
than its deck group $A$. We perform a surgery on $M'$ to avoid this. By
doing so, we are careful to preserve the fact that we keep a Loch Ness
monster.

Let $[-2,-1] \times [-1/2, 1/2]$ be the square in $M'_Id \simeq \R^2$.
We remove $Z$ from the surface and glue back the opposite sides left open in the surface
$M_{Id}$. This creates a new surface $M_{Id}$ with a conical singularity $p_{Id}$ of
angle $6\pi$ with an
horizontal and a vertical saddle connection joining $p_{Id}$ to itself. Moreover,
if we draw the eight half-lines starting from the singularity in the direction
of these two saddle connections, we separate $M'_{Id}$ into 4 quadrants
and 2 infinite bands. Among these quadrants, two of them contain no
singularities. See Figure~\ref{fig:EveryGroup}.

We perform the same surgery in $M'_a$ for $a \in A$ by removing the rectangle
$F(a) \cdot Z$. We denote this new surface $M_a$ and by $p_a \in M_a$ the singularity
of angle $6\pi$ coming from the surgery in the surface $M'_a$.
We denote by $M$ the surface obtained by gluing the slits $m_i^j$ as in step 1.
$M$ is homeomorphic to the Loch Ness monster (the same argument we
gave in step 1 applies). The deck group $A$ acts by affine transformation
and we keep the notation $f_a$ to denote the action of $a$.
Note that $M$ is a ramified covering of the plane $\R^2$ from which we
removed the rectangle $[-2,-1] \times [-1/2,1/2]$.

Now let $f \in\Aff^+(M)$ be an affine automorphism. Because the only conical
singularities of angle $6\pi$ in $M$ are the $p_a$ we need to have that
$f(p_{Id}) = p_{a_0}$ for some $a_0 \in A$. We claim that $f = f_{a_0}$.
Indeed, the horizontal and vertical saddle connections in $M_{Id}$ joining
$p_{Id}$ to itself must be mapped to saddle connections joining $p_{a_0}$
to itself. Moreover, the 2 infinite bands parallel to these saddle connections and the
four quadrants are mapped to infinite bands and quadrants. Since two of
these quadrants do not contain any singularity, the only possibility
is that $f_{a_0} = f$.
\end{proof}

\textbf{Veech groups and the translation flow}. The following classical result of Veech, relates the translation flow on a finite-type surface with its Veech group:

\begin{tcbtheorem}{\cite{Veech89}}{VeechDichotomy}
Let $M$ a finite-type translation surface whose Veech group $\Gamma(M)$ is a lattice\footnote{That is, such that $\H/\Gamma(M)$ has finite hyperbolic area.}.
Then for any direction $\theta$ the translation flow $F_{M,\theta}^t$ in direction
$\theta$ in $M$ is either completely periodic\footnote{A direction is called completely periodic if $F_{M,\theta}^t$ decomposes $M$ into a (finite) family of maximal cylinders. In the context of this theorem the moduli of two cylinders appearing in the decomposition are always commensurable.} and the direction $\theta$ is stabilized
by a multitwist in $\Aff(M)$, or it is uniquely ergodic.
\end{tcbtheorem}

\begin{tcbdefinition}{}{definition-of-veech-surface}
  A finite-type translation surface $M$ whose Veech group is a lattice is  called a \emph{Veech surface} or a \emph{lattice surface}.
\end{tcbdefinition}

It is important to remark that finite-type Veech translation surfaces are quite rare. More precisely:
\begin{tcbtheorem}{\cite{Moeller09}}{TH:VeechGroupsAreTrivialGenerically}
The affine group $\Aff^+(M)$ of genus $g(M)\geq 2$ generic translation surface is $\Z/2\Z$ or trivial, depending whether $M$ belongs to the hyperelliptic component or not.
\end{tcbtheorem}
Here generic means that $M=(X,\omega)$ lies in the complement of a countable union of real codimension one submanifolds of its stratum.

Theorem~\ref{thm:VeechDichotomy} generalizes a classical result that tells us that the translation flow $F_\theta^t$ on the flat torus $\R^2/\Z^2$ is either periodic or uniquely ergodic depending on whether $\theta$ is a rational multiple of $\pi$ or not. We invite the reader to show periodicity for rational slopes.
Unique ergodicity for irrational slope is an early XXth century result by H. Weyl, W. Sierpiński and P. Bohl. Their proofs do not involve any action of the Veech group, contrarily to the one of W.~Veech in~\cite{Veech89} for the general case. Note however that finer dynamical properties of the translation flows might differ between tori, origamis and general Veech surfaces such as weak mixing~\cite{AvilaDelecroix} or rigidity~\cite{Ferenczi-Hubert}. For a more detailed discussion on Veech surfaces of finite type, we refer the reader to~\cite{HubertSchmidt}.



The purpose of the exercise below is to show that for infinite-type surface knowing the Veech group is not enough to determine the dynamics of the translation flow.

\begin{tcbexercise}{}{VeechSaysNothingAboutDyn}
In this exercise we invite the reader to make a modification of the construction
in the proof of Theorem~\ref{thm:AnyGroupIsVeech} so that we can control the
translation flows of the surface $M$.

Let $A$ be a finitely generated group with generating set $U = \{a_1, \ldots, a_k\}$.
Let $F: A \to \GL_+(2,\R)$ be a morphism. Show that there exist a translation surface
$M$ as in the conclusion of Theorem~\ref{thm:AnyGroupIsVeech} that has the following
additional property: for every direction
$\theta$ and for any $x\in M$ then either the orbit $F_{\theta}^t(x)$ hits a
singularity in finite time or it escapes every compact set $K \subset M$.

\emph{Hint}: in step 1 of the proof of Theorem~\ref{thm:AnyGroupIsVeech} instead of
considering slits in $\R^2$ you could consider a cyclic covering of degree $2k$
of $\R^2$ ramified over the origin and place each family $\{m_i^j\}_i$ of slits
in the different sheets of this covering.
\end{tcbexercise}
Exercise~\ref{exo:VeechSaysNothingAboutDyn} should warn the reader that the Veech
group does not capture interesting properties of translation flows of
infinite-type translation surfaces. However, recall that in the proof
of Theorem~\ref{thm:AnyGroupIsVeech} the surface $M$ is a $A$-covering of the plane
with a square removed. Hence the affine group $\Aff^+(M) \simeq A$ (which is the
same thing as the deck group of the covering) acts properly discontinuously
(see Definition~\ref{def:ProperlyDiscontinuousAction}).
In many examples in this book such as the infinite staircase, wind-tree models and
the Hooper-Thurston-Veech surfaces (that we introduce at the end of this Chapter)
the affine group does not act properly discontinuously and there is a rich interaction
between the affine group and the dynamics of the translation flows.

\subsection{Veech groups of tame translation surfaces}
Recall from Definition~\ref{def:Singularities} that a translation surface $M$ is tame if $\Sing(M) \subset \widehat{M}$ contains only conical and infinite angle singularities. We now discuss how this geometric condition reflects and imposes constraints on the Veech group.

\begin{tcbdefinition}{}{ContractiveMatrix}
A matrix $g \in \GL(2,\R)$ is \emphdef[contracting (element of $\GL(2,\R)$]{contracting} if $\displaystyle \lim_{n \to +\infty} g^n$ exists and is the zero matrix.
\end{tcbdefinition}
Equivalently, a matrix $g \in \GL(2,\R)$ is contracting if its eigenvalues are smaller than one in absolute value.

\begin{tcbtheorem}{}{TameVeechGroup}
Let $M$ be a tame translation surface with at least two singularities.
Then its Veech group $\Gamma(M)$ does not contain contracting elements.
\end{tcbtheorem}

\begin{proof}
Let $\|.\|$ denote the Euclidean norm on $\R^2$.
We proceed by contradiction. Let $\phi \in \Aff(M)$ such that $g=D(\phi)$ is contracting.
Then, there exists $n_0 > 0$ such that for all $v \in \R^2$ we have $\|g^{n_0} v\| < \|v\|$.
Because of this property, the map $\psi := \phi^{n_0}$ is a contraction for the metric
on $\widehat{M}$. Hence it has a unique fixed point $x_0 \in \widehat{M}$
and for all $x \in \widehat{M}$ we have that
$\psi^n(x)$ converges to $x_0$ as $n \to \infty$ (Banach-fixed point theorem).

By assumption there exists a singularity $x_1 \neq x_0$. Hence $\psi^n(x_1)$ is a
sequence of singularities that converges to $x_0$. Therefore $x_0$ must be a wild
singularity. This is the desired contradiction.
\end{proof}

A natural question from Theorem~\ref{thm:TameVeechGroup} is whether we can
extend the construction of Theorem~\ref{thm:AnyGroupIsVeech} to guarantee
that the result $M$ is tame. As shown in the following result, the answer
is positive.

\begin{tcbtheorem}{\cite{PrzytyckiSchmithuesenValdez11}}{VeechGroupLochNessTame}
Let $G$ be a countable subgroup of $\GL^+(2,\R)$ without contracting elements. Then there
exists a tame translation surface $M$ homeomorphic to the Loch Ness monster such that
$\Gamma^+(M)=G$.
\end{tcbtheorem}
Recall that in the proof of Theorem~\ref{thm:AnyGroupIsVeech} the surface $M$ is built out of planes with slits $M'_a$ and each of this piece is tame.
However, this is not enough to guarantee that the construction produces a tame translation surface (think about the case where $F(A)$ is a group
of diagonal matrices). In order to prevent accumulation of singularities the strategy of~\cite{PrzytyckiSchmithuesenValdez11} consists
in using what they call ``buffer surfaces" to separate $M_a$ from the other pieces $M_{a s_j}$ it is glued to. We refer the reader to the
original article for details about the construction.

From Theorems~\ref{thm:AnyGroupIsVeech} and~\ref{thm:VeechGroupLochNessTame} we conclude that there are no restrictions (other than the obvious ones) to realize countable subgroups of matrices as Veech groups of Loch Ness monsters, even if we impose the translation surface structure to be tame. The following result tells us that are no restrictions either for the ``most complicated'' infinite genus surface which is, informally speaking, the Blooming Cantor tree (the surface of infinite genus without planar ends and whose space of ends is homeomorphic to the Cantor set).
\begin{tcbtheorem}{\cite{Ramirez-Valdez16}}{VeechGroupBloomingCantorTreeEtAl}
Let $G<\GL_+(2,\R)$ be a countable subgroup without contracting elements. Then there exists a tame translation surface $M$ homeomorphic to the Blooming Cantor tree for which $\Gamma^+(M)=G$.
\end{tcbtheorem}

Again, the proof of this result is similar in spirit to the proof of Theorem~\ref{thm:VeechGroupLochNessTame}: given a set of generators $U$ of $G$ one constructs a family of surfaces $\{M_g\}$, each homeomorphic to the Blooming Cantor tree, and glues them following the structure of the Caley graph of $G$ relative to $U$. Although the proof is simple, it is far from elementary, for the gluings have to be done with care. Results similar to the preceding Theorem hold for an infinite family of topological types of infinite-type translation surfaces. For details we refer the reader to~\cite{Ramirez-Valdez16}.

Recall from Chapter~\ref{ch:TopologyGeometry} that to any $g\in\{\Z_{\geq 0}\cup \infty\}$ and
a nested couple of closed subsets $C'\subset C$ of the
Cantor set corresponds a unique topological type of surface: the one with genus
$g$ and such that $C'\subset C$ is homeomorphic to
$\Ends_{\infty}(S)\subset\Ends(S)$. The following general situation is till an open problem
\begin{tcbquestion}{}{VeechGroupAnySurface}
\label{Q:VeechGroupAnySurface}
Let $C'\subseteq C$ be a nested couple of closed subsets of the Cantor set and $g\in\{\Z_{\geq 0}\cup \infty\}$. Is it possible to realise any countable subgroup of $\GL_+(2,\R)$ (with or without contracting elements) as the Veech group of a (tame) translation surface $M$ such that $C'\subset C$ is homeomorphic to $\Ends_{\infty}(S)\subset\Ends(S)$.
\end{tcbquestion}
This question remains open even for simple cases such as Jacob's Ladder, that is when $C'=C$ has only two elements. It is also unknown whether there exist a translation surface homeomorphic to Jacob's Ladder whose Veech group is $\SL(2,\Z)$.

Recall that in Section~\ref{sec:AnalyticAutomorphisms} we classified translation surfaces $M$ for which
the subgroup $\Tr(M)$ of translations in $\Aff(M)$ is uncountable (see Theorem~\ref{thm:TranslationSurfacesBigAut}).
Here we classify tame translation surfaces for which the Veech group is uncountable.
\begin{tcblemma}{}{}
Let $M$ be a tame translation surface such that its Veech group $\Gamma^+(M)$ is uncountable.
Then $\Gamma^+(M)$ is conjugate to either $\GL_+(2,\R)$,
$$
P:=\{\left(\begin{smallmatrix} 1 & t \\ 0 & s \end{smallmatrix}\right)\hspace{1mm}|\hspace{1mm}\text{$t\in\R$ and $s\in\R^+$}\}
$$
or to $P'<\GL_+(2,\R)$, the group generated by $P$ and $-Id$.
\end{tcblemma}
\begin{proof}
Suppose first that $M$ has no saddle connections. There are two subcases in this regime: $M$ has singularities or it only one singularity. If $M$ has no singularities then, by the tameness assumption, $M$ is either $\R^2$, or a flat infinite cylinder, \eg $\R^2/\Z$ or a flat torus. The Veech group of a flat torus is conjugated within $\GL(2,\R)$ to $\SL(2,\Z)$, so it is countable. On the other hand if $M=\R^2$ then $\Gamma^+(M)=\GL_+(2,\R)$ and if $M$ is an infinite cylinder then its Veech group of conjugated to $P'$. If $M$ has only one singularity, tameness implies that $M$ must be a cyclic covering of the plane ramified over a point, case in which $\Gamma^+(M)=\GL_+(2,\R)$.

If $M$ carries saddle connections and the group $\Gamma^+(M)$ is uncountable, then $G$ is conjugate to either $P$ or $P'$. First remark that tameness implies that all holonomy vectors (and hence saddle connections) must be parallel. Indeed, tameness implies that the set $V_{hol}(M)$ of holonomy vectors of $M$ is at most countable. If there were two linearly independent holonomy vectors $v_1$, $v_2$ in $V_{hol}(M)$ then the map $\eta:\Gamma^+(M)\to V_{hol}(M)\times V_{hol}(M)$ defined by the linear action of matrices of vectors  by $\eta(g)=(g(v_1),g(v_2))$ would be an embedding, which is impossible since $V_{hol}(M)$ is countable and $\Gamma(M)$ uncountable. Hence $P$ is a subgroup of $\Gamma^+(M)$. To conclude we claim that  $P<\Gamma^+(M)<P'$ and since $P$ is of index 2 in $P'$ we have $G=P$ or $G=P'$. Indeed, the claim follows from the fact that $Spine(M)$, which is the union of all singularities of $\widehat{M}$, and all horizontal separatrices and saddle connections, is a complete metric space w.r.t. its intrinsic path metric; then if $g\in\Gamma^+(M)$ and $e_1$ is a unit horizontal vector, then from the preceding discussion we have that $g(e_1)=\pm e_1$, for in any other case $g$ restricted to $Spine(M)$ would be a contraction having a unique fixed point to which singularities must accumulate, contradicting the surface's tameness.
\end{proof}
\begin{tcbexercise}{}{GL2RTopology}
Prove that if $M$ is a translation surface with Veech group $\Gamma^+(M)=\GL_+(2,\R)$, then $M$ is either the plane or a covering of the plane ramified over one point.
\end{tcbexercise}

In the light of the preceding exercise, it is natural to ask if there are topological restrictions for a tame translation surface to have a Veech group conjugated to either $P$ or $P'$. As the following theorem illustrates, this is not the case for surfaces without planar ends.
\begin{tcbtheorem}{\cite{Ramirez-Valdez16}}{}
Let $C$ be any closed subset of the Cantor set. Then there exists tame translation surfaces $M$ and $M'$ for which:
\begin{itemize}
\item $\Ends(M)=\Ends_\infty(M)$, $\Ends(M')=\Ends_\infty(M')$ and both of these spaces are homeomorphic to $C$, and
\item the Veech groups of $M$ and $M'$ are conjugated to $P$ and $P'$ respectively.
\end{itemize}
\end{tcbtheorem}

\begin{proof}
We only provide a sketch of the proof to illustrate the main ideas.
First, let's consider the case where $C$ is just one point (\ie when the surfaces are homeomorphic to a Loch Ness monster) and discuss the key points of the general case.
The construction starts from the Euclidean plane $\R^2$. We consider the family $L^+$ of horizontal segments
$\ell^+_i := [((4i-1), 0), (4i, 0)]$ for $i \in \N$. Let $M_P$ be the tame translation surface obtained by gluing the segment $l^+_{2i-1}$ to the segment $l^+_{2i}$ for each $i\in\N$. Given that $L^+$ is a discrete family of parallel segments contained in the real axis, the surface $M_P$ is a tame Loch Ness monster. Moreover, its Veech group is precisely $P$ because the set of all holonomy vectors of $M_P$ is $\{\pm e_1\}$ and there are no saddle connections on the half plane $x<0$. To produce a Loch Ness monster with Veech group $P'$ consider on $\R^2$ the additional families of parallel segments $L^-$ made of $\ell_i^- := [((1-4i), 0), (-4i, 0)]$ for $i \in \N$.
and define $M_{P'}$ as the surface
obtained by gluing $l^+_{2i-1}$ to $l^+_{2i}$ and $l^-_{2i-1}$ to $l^-_{2i}$, for each $i\in\N$. The Veech group of $M_{P'}$ is $P'$ for reasons that are analogous to the ones presented for the surface $M_P$.

The general case follows from~\cite[construction 1.1]{Ramirez-Valdez16}. This construction takes as input data a closed subset $C$ of the Cantor set and a tame Loch Ness monster $M$ and produces a tame translation surface $M_{elem}$ (called an \emph{elementary piece}) satisfying $\Ends_\infty(M_{elem})=\Ends(M)=C$. When
$M=M_P$ (respectively $M_{P'}$) it is possible to assure that $P$ (respectively $P'$) persists through the steps of the construction so that $\Gamma(M_{elem})=P$ (respectively $\Gamma(M_{elem})=P'$). The surface $M_{elem}$ is obtained by a cut-and-paste construction similar to the one presented in the proof of Theorem \ref{thm:AnyGroupIsVeech}, the main difference being that instead of following the Cayley graph of a group as guide for the gluings one emulates the structure of a subtree $T_C$ of the infinite regular tree of degree 3 which satisfies $Ends(T_C)=C$.
\end{proof}

\begin{tcbquestion}{}{}
Let $C'\subseteq C$ be a nested couple of closed subspaces of the Cantor set. Is it possible to realise the group $P$ (respectively $P'$) as Veech group of a tame translation surface $M$ (respectively $M'$) such that $C'\subseteq C$ is homeomorphic to $\Ends_\infty(M)\subset\Ends(M)$?
\end{tcbquestion}

\section{Affine groups and Veech groups of coverings}
\label{sec:AffineGroupsCoverings}
We now study affine groups of translation coverings $p: \widetilde{M} \to M$ where $M$ is a finite-type translation surface (see Section~\ref{sec:CoveringSpaces}, Definition~\ref{def:TranslationCovering}). More precisely, we will study which affine elements of $M$ can be \emph{lifted} to $\widetilde{M}$; that is whether there exists $\widetilde{f} \in \Aff(\widetilde{M})$ so that $p \circ \widetilde{f} = \widetilde{f} \circ p$.

In Section~\ref{ssec:LiftingCriterion} we present the general lifting criterion which follows from standard algebraic topology. Next, in Section~\ref{ssec:AffineGroupsOfAbelianCoverings} we study the case of Abelian coverings and discuss the importance of the affine group action on homology of the base surface. In Section~\ref{ssec:LiftingParabolicElements} we focus on parabolic elements of $\Aff(M, \Sigma)$ that, as we already have seen in the Thurston-Veech construction from Section~\ref{ssec:ThurstonVeechConstructionsIntro}, are intimately related to cylinder decompositions of $M$. Then, in Section~\ref{ssec:HooperWeissTheorem} we prove a theorem due to Hooper and Weiss that provides a powerful method to show that the Veech group of a covering is non-elementary\footnote{Recall that a subgroup $G$ of $\SL(2,\R)$ is called \emph{elementary} if its limit set $\Lambda(G)\subset\partial\mathbb{H}$ has at most two elements, see also Appendix~\ref{app:FuchsianGroups}} Finally, in Section~\ref{ssec:SymmetriesWindTree}
we apply these results to the wind-tree model.

Most statements concern translation surfaces though they
can be applied to half-translation surfaces
(see Section~\ref{ssec:HalfTranslationSurfaces}) after an
appropriate double-covering construction.

\subsection{The lifting criterion}
\label{ssec:LiftingCriterion}
Let us start with the general lifting criterion. Let $p: \widetilde{M} \to M$ be a translation covering branched over $\Sigma \subset M$. In order to lift affine elements from $M$ to $\widetilde{M}$, the branching locus $\Sigma$ must be preserved. Hence we only deal with elements of the group $\Aff(M, \Sigma)$, which is the subgroup of $\Aff(M)$ preserving the set $\Sigma$. The obstruction for a map $f$ in $\Aff(M,\Sigma)$ to be lifted to $\widetilde{M}$ can be measured by its action on the fundamental group of $M \setminus \Sigma$. Recall that if $x_0$ is any point in $M$ then $f$ induces a group morphism $f_*: \pi_1(M, x_0) \to \pi_1(M, f(x_0))$. Given that $f$ can change basepoints, $f$ does not necessarily induces an automorphism of $\pi_1(M, x_0)$. However, the action of $f$ on conjugacy classes is well defined. To lighter notations, we will often write $\pi_1(M)$ without reference to the base point when the properties we require are conjugacy invariant.

\begin{tcbtheorem}{Lifting criterion}{LiftingCriterion}
Let $p:\widetilde{M}\to M$ be a translation covering with branching locus $\Sigma\subset M$ and define $\widetilde{\Sigma} := p^{-1}(\Sigma)$. Let $f\in\Aff(M,\Sigma)$ (\ie $f \in \Aff(M)$ and $f(\Sigma)=\Sigma$). Then there exists a lift $\widetilde{f}:\widetilde{M}\to\widetilde{M}$, (\ie $f \circ p = p \circ \widetilde{f}$ on $\widetilde{M}$) if and only if the induced outer\footnote{Recall that outer automorphism group of a group, G, is the quotient, $\Aut(G) / \Inn(G)$, where $\Aut(G)$ is the automorphism group of G and $\Inn(G)$ is the subgroup consisting of inner automorphisms. } morphism $f_* \in \Out(\pi_1(M \setminus \Sigma))$ preserves the conjugacy class of $p_*(\pi_1(\widetilde{M} \setminus \widetilde{\Sigma}))$ in $\pi_1(M \setminus \Sigma)$.
\label{LiftingCriterion}
\end{tcbtheorem}

\begin{proof}
From Proposition 1.33, p. 61 in~\cite{Hatcher} we have the following topological criterion:
\begin{tcblemma}{Topological lifting criterion}{}
Let $S$ be a topological surface, $p:\widetilde{S}\to S$ a covering map and $f:S\to S$ a continuous function. Let $x_0\in S$, $x_1=f(x_0)$ and $\tilde{x}_0\in p^{-1}(x_0)$, $\tilde{x}_1\in p^{-1}(x_1)$. Then there exists a lift $\widetilde{f}:\widetilde{S}\to\widetilde{S}$ satisfying $\widetilde{f}(\tilde{x}_0)=\tilde{x}_1$ if and only if the induced morphisms $f_*: \pi_1(S, x_0) \to \pi_1(S, x_1)$, $p_{0*}: \pi_1(\widetilde{S}, \tilde{x}_0) \to \pi_1(S, x_0)$ and $p_{1*}: \pi_1(\widetilde{S}, \tilde{x}_1) \to \pi_1(S, x_1)$ on the corresponding fundamental groups satisfy
$$
f_* p_{0*}\pi_{1}(\widetilde{S}, \widetilde{x_0}) \subset p_{1*}\pi_1(\widetilde{S}, \widetilde{x_1})
$$
\end{tcblemma}
In order to apply the above Lemma, we fix a basepoint $x_0$ in $M^0=M \setminus \Sigma$ and let $x_1$ be its image under $f$. Since $f$ is an element of $\Aff(M,\Sigma)$ the image $x_1$ is not in $\Sigma$. Now, a pair of points $\tilde{x}_0$ and $\tilde{x}_1$ in $\widetilde{M}^0$ exist if and only if $f$ preserves the conjugacy class of $(p_0)_*(\pi_1(\widetilde{M}^0,\tilde{x}_0))$. If such a lifting exists, it can be uniquely extended to $\widetilde{M}$ in a continuous way.
\end{proof}
\begin{tcbexercise}{}{}
Let $p:\widetilde{M}\to M$ be a translation covering. We define the subgroups $G_1(p) \subset \Aff(M)$ and $G_2(p) \subset \Aff(\widetilde{M})$ as follows
\begin{align}
\label{eq:DefinitionG1}
G_1(p)&:=\{f\in\Aff(M):\text{there exists $\widetilde{f}\in\Aff(\widetilde{M})$ such that $p\circ\widetilde{f}=f\circ p$}\}\\
\label{eq:DefinitionG2}
G_2(p)&:=\{\widetilde{f}\in\Aff(\widetilde{M}): \text{there exists $f\in\Aff(M)$ such that $p\circ\widetilde{f}=f\circ p$}\}
\end{align}
\begin{enumerate}
\item Prove that the composition with $p$ provides a surjective map
\[
\begin{array}{ccc}
G_2(p)    & \to & G_1(p) \\
\tilde{f} & \mapsto & p \circ \tilde{f}
\end{array}.
\]
\item Show that the kernel of the map $G_2(p) \to G_1(p)$ considered above is $\Deck(p)$.
\item Assume that $p: \widetilde{M} \to M$ is a normal cover, that is $\Deck(p)$ acts transitively on the fibers of $p$. Then prove that $G_2(p)$ is the normalizer of $\Deck(p)$ in $\Aff(\widetilde{M})$ (\ie $\tilde{f}$ in $\Aff(\widetilde{M})$ belongs to $G_2(p)$ if and only if $\tilde{f} \Deck(p) \tilde{f}^{-1}=\Deck(p)$).
\item Prove that if $\tilde{f} \in G_2(p)$ maps to $f \in G_1(p)$ then the derivatives $D(\tilde{f})$ and $D(f)$ are equal.
\end{enumerate}
\end{tcbexercise}
\begin{tcbdefinition}{}{RelativeVeechGroup}
Let $p: \widetilde{M} \to M$ be a translation covering. The \emphdef{relative Veech group} $\Gamma(p)$ or $\Gamma(\widetilde{M} \to M)$
of the covering $p$ is the group $D(G_1(p)) = D(G_2(p))$, that is, the image by the derivative homomorphism of the groups $G_1(p)$ and $G_2(p)$ defined in~\eqref{eq:DefinitionG1} and \eqref{eq:DefinitionG2}.
\end{tcbdefinition}

The following direct consequence of Theorem~\ref{LiftingCriterion} is already a rich source of examples.
\begin{tcbcorollary}{}{LiftingCharacteristicCoverings}
Let $M$ be a translation surface, $\Sigma\subset M$ a discrete subset and $M^0=M\setminus\Sigma$. Let $p:\widetilde{M}\to M$ be a translation covering defined by a characteristic subgroup\footnote{A subgroup $H$ of a group $G$ is \emphdef[characteristic (subgroup)]{characteristic} if for any automorphism $\phi:G \to G$ we have $\phi(H) = H$.} of $\pi_1(M^0)$. Then every element of $\Aff(M,\Sigma)$ lifts to $\Aff(\widetilde{M})$.
In particular, $\Gamma(M, \Sigma) \subset \Gamma(\widetilde{M})$.
\end{tcbcorollary}
\begin{tcbexample}{}{CharacteristicCoverings}
Let $G$ be a group. We consider the following families of subgroups of $G$:
\begin{enumerate}
\item The lower central series:
\begin{equation}
\label{eq:CentralSeries}
\gamma_0(G)\unrhd \gamma_1(G)\unrhd\ldots\unrhd \gamma_n(G)\unrhd\ldots
\end{equation}
where $\gamma_0(G) := G$ and $\gamma_{n+1}(G):=[\gamma_n(G),G]$ for each $n>0$.
\item The derived series:
\begin{equation}
\label{eq:Derivedseries}
\delta_0(G) \unrhd \delta_1(G) \unrhd\ldots\unrhd \delta_n(G) \unrhd\ldots
\end{equation}
where $\delta_0(G) = G$ and $\delta_{n+1}(G):=[\delta_n(G), \delta_n(G)]$ for $n>0$.
\item The subgroup of $n$-th powers:
\begin{equation}
\label{eq:nthpowers}
G^n:=\langle g^n:g\in G\rangle
\end{equation}
\end{enumerate}
It is not difficult to check that all subgroups of $G$ listed above are characteristic. In particular Corollary~\ref{cor:LiftingCharacteristicCoverings} applies. Let us consider the case of square-tiled surfaces where everything can be made more explicit.
\begin{tcbcorollary}{}{LiftingCharacteristicCoveringsConcrete}
Let $p: \widetilde{M} \to \T^2$ be a square-tiled surface determined by a subgroup of $F_2 \simeq \pi_1(\T^2 \setminus \{0\})$ in one of the families~\eqref{eq:CentralSeries}, \eqref{eq:Derivedseries} or~\eqref{eq:nthpowers} defined above. Then every element of $\Aff(\T^2,\{0\}) \simeq \SL(2,\Z)$ lifts to $\widetilde{M}$, in particular $\Gamma(\widetilde{M})$ contains $\SL(2,\Z)$. Moreover, if $p$ is ramified then $\Gamma(\widetilde{M}) = \SL(2,\Z)$.
\end{tcbcorollary}
\begin{proof}
The first part of the statement directly follows from Corollary~\ref{cor:LiftingCharacteristicCoverings}. The Veech group of the torus with a marked point is $\SL(2,\Z)$. In particular the Veech group of $\widetilde{M}$ contains $\SL(2,\Z)$. However, it might not be equal to it: the Euclidean plane admits affine automorphisms that permute the fibers of the covering map $p: \C \to \T^2$. When the covering is ramified, the fact that the covering subgroup is characteristic ensures that the corners of each square in $\widetilde{M}$ are conical singularities (when the ramification order is finite) or punctures that gives rise to infinite angle singularity in the metric completion (when the order is infinite). All other points of $\widetilde{M}$ are regular. Hence, any element of $\Aff(\widetilde{M})$ preserves $\widehat{p}^{-1}(\{0\})$, where $\widehat{p}$ is the unique extension of $p$ to the metric completion (see Definition~\ref{def:TranslationCovering} and the paragraph that follows). As a consequence every element of $\Aff(\widetilde{M})$ descends to $M$.
\end{proof}

To illustrate let us consider the subgroup $\gamma_1(F_2) = \delta_1(F_2) = [F_2,F_2]$ of the fundamental group $F_2=\pi_1(\T^2\setminus\{0\})$. The associated (non-ramified) covering space $p:\C\setminus\Z^2\to\T^2\setminus\{0\}$ is a $\Z^2$-covering, because $\Z^2=F_2/\delta_1(F_2)$. Now, recall from Definition~\ref{def:TranslationCovering} in Chapter~\ref{ch:TopologyGeometry} that the corresponding translation covering is obtained by adding the conical singularities and regular points in the metric completion of $\C\setminus\Z^2$ to which $p$ can be extended continuously. In the case of $\C \setminus \Z^2$, the metric completion is obtained by adding $\Z^2$ which are all regular points and $p$ is the (unramified) universal cover $\widehat{p}:\C\to\T^2$.
An other explicit example is given by the second group in the central series: $F_2/\gamma_2(F_2)$ is isomorphic to the \emph{Heisenberg group}
$$
\left\{\left(\begin{smallmatrix} 1 & a & c \\ 0 & 1 & b \\ 0 & 0 & 1 \end{smallmatrix}\right) \hspace{1mm}|\hspace{1mm} a,b,c\in\Z \right\},
$$
More generally, subgroups of the central series gives rise to a nilpotent coverings.

The second group in the derived series $\delta_2(F_2)$ is however less explicit. The quotient $F_2 / \delta_2(F_2)$ is called the \emph{free metabelian group of rank 2}.
\begin{tcbexercise}{}{}
\begin{enumerate}
\item Show that for $n \geq 0$, $[\gamma_{n+1}(F_2):\gamma_n(F_2)] = +\infty$ and $[\delta_{n+1}(F_2):\gamma_n(F_2)] = +\infty$.
\item Show that the subgroups $\gamma_n(F_2)$ and $\delta_n(F_2)$ are pairwise distinct except: $\gamma_0(F_2) = \delta_0(F_2) = F_2$ and $\gamma_1(F_2) = \delta_1(F_2) = [F_2,F_2]$.
\end{enumerate}
(\emph{hint}: use the fact that subgroups of $F_2$ are free groups)
\end{tcbexercise}

The study of quotients by the $n$-th power subgroup is much more delicate and has a long history. The quotient $B(m,n) := F_m/(F_m)^n$ is called the \emph{free $m$-generator Burnside group of exponent $n$} (each element in $B(m,n)$ has order at most $n$). It has been known for a while that for $n=2,3,4,6$ these groups are finite. However,
it is a deep result that there exists $n_0$ such that $B(2,n)$ is infinite for every $n\geq n_0$ (Adian and Novikov (1968), Ol'shanski\u{i} (1982), Ivanov (1994), Lysenok (1996), Delzant and Gromov (2008), and Coulon (2018)). Let us also mention that it is still unknown whether $B(2,5)$ is infinite. More recently, a versatile and easy-to-apply tool for constructing examples of finitely generated infinite $n$-periodic groups (such as the Burnside groups) has been developed by Coulon and Gruber, see~\cite{Coulon-Gruber} for details.
\end{tcbexample}

We now present two examples of coverings that are not translation coverings but for which the deck group acts by affine transformations. The first example are irrational billiards (for which the deck group acts by isometries) and the second one dilation surfaces (for which the Veech group acts by dilations).

\begin{tcbexample}{}{}
In Section \ref{ssec:PolygonalBilliards} we discussed how to translate the dynamics of a billiard ball in a polygon $P$ into the dynamics of the translation flow on a translation surface $M(P)$. Moreover, Theorem~\ref{thm:AnIrrationalBilliardIsLochNess} proven in Example~\ref{exa:IrrationalBilliardAndProof} shows that if $P$ is an irrational polygon (\ie none of its interior angles is commensurable with $\pi$) then $M(P)$ is homeomorphic to the Loch Ness monster.
\begin{tcbtheorem}{\cite{Valdez12}}{NonDiscreteVeechGroupsExist}
\label{THM:NonDiscreteVeechGroupsExist}
Let $P$ be a simply connected Euclidean polygon with interior angles $\{\lambda_j\pi\}_{j=1}^N$. If there exists $\lambda_j\in\R\setminus\Q$ then the Veech group $\Gamma^+(M(P))$ is a subgroup of $\SO(2,\R)$. Moreover, the subgroup of $\SO(2,\R)$ generated by the rotations $z\to e^{2\lambda_j\pi}z$, $j=1,\ldots,N$ is a subgroup of $\Gamma^+(M(P))$ of finite index.
\end{tcbtheorem}

\begin{proof}
We discuss the case of triangles $P$, for the case of a general polygon is analogous. From Section~\ref{ssec:PolygonalBilliards} and the proof of Theorem~\ref{thm:AnIrrationalBilliardIsLochNess} presented in Example~\ref{EXAMPLE:IrratBillAndProof} we can deduce that there exists a non-ramified covering $p:M(P)\to \mathbb{S}^2(P)$ whose base is a three-punctured sphere made of two copies of $P$ (with the vertices removed) and whose deck transformation group is isometric to $\Z\times A$ or $\Z^2\times A$ (where $A$ is a finite abelian group) depending on whether the rank of $P$ is equal to 1 or 2, see (\ref{eq:RankPolygon}) in Section~\ref{ssec:PolygonalBilliards}. Let us emphasize that here $\mathbb{S}^2(P)$ is not a translation surfaces and $p$ is not a translation covering.
The surface $\mathbb{S}^2(P)$ has a $\Isom^+(\R^2)$-structure (see Appendix \ref{Appendix:GXStructures}). Let $\hol:\pi_1(M(P))\to \Isom^+(\R^2)$ be its holonomy. Then the covering $p$ is defined by the normal subgroup $\hol^{-1}(\Trans(\R^2))$.

By construction, the image of $\hol$ is contained in the group generated by the reflections with respect to the sides of $P$. Let $R(P)$ be the group generated by the linear part of the reflections with respect to the sides of $P$ and $R^+(P) = R(P) \cap \SO(2)$. The deck group of $p$ acts on the translation surface $M(P)$ by affine homeomorphisms whose linear parts generate the group $R^+(P)$. We claim that actually $\Gamma^+(M(P))$ is formed only by rotations. Indeed, in $M(P)$ there are saddle connections of minimal length, say $L>0$, corresponding to geodesics of minimal length in the Euclidean surface $\mathbb{S}^2(P)$. Since one of the $\lambda_j$ is irrational, $R^+(P)$ contains an irrational rotation. In particular, the holonomy vectors of length $L$ are dense in the circle $\{z\in\C:|z|=L\}$. Now, if we had an element $A\in\Gamma^+(M(P))$ different from a rotation, then the image of $\{z\in\C:|z|=L\}$ under $A$ would be an ellipse containing infinitely many points corresponding to holonomy vectors of norm less than $L$, which is a contradiction to the fact that $L$ is the minimal length of a saddle connection. To finish the proof remark that the set of holonomy vectors of minimal length of $M(P)$ is of the form $\sqcup_{i=1}^k R^+(P)v_i$, where $v_i$ are minimal length holonomy vectors and the group $\Gamma^+(M(P))/R^+(P)$ acts freely on the set of $R^+(P)$-orbits of holonomy vectors of minimal length.
\end{proof}
\end{tcbexample}

\emph{Veech groups of translation coverings of dilation surfaces}. At the end of Section~\ref{sec:CoveringSpaces} we explained that each dilation surface $S$ has a natural normal covering $\pi:\widetilde{S}\to S$ defined by the inverse image of $\Trans(\R^2)$ under the holonomy $\hol: \pi_1(S)\to \Dil^+(\R^2)$ of the $\Dil^+(\R^2)$-structure of $S$. The pullback of the dilation surface structure of $S$ defines thus a translation surface structure on $\widetilde{S}$ and the deck transformation group of $\pi$ acts on this translation surface by dilations. In particular, if $\Delta$ denotes the image of $\pi_1(S)\to \Dil^+(\R^2) \to \R^*$, where the second arrow is just the derivative map, then the Veech group of $\widetilde{S}$ contains an element of the form $\begin{psmallmatrix} a & 0 \\ 0 & a \end{psmallmatrix}$, for each $a\in\Delta$.

\begin{tcbexercise}{}{}
Let $S$ be a finite type dilation surface and $\pi:\widetilde{S}\to S$ its canonical translation cover. Let $\Aff(S)$ be the set of homeomorphisms of $S$ which are affine in local coordinates. Show that:
\begin{enumerate}
  \item $\Aff(S)$ is a group.
  \item Every element in $\Aff(S)$ lifts to an element of $\Aff(\widetilde{S})$.
\end{enumerate}
\end{tcbexercise}

Veech groups can also be defined for dilation surfaces and have been studied in~\cite{DuryevFougeronGhazouani19}.
 As with translation surfaces, the characterization of Fuchsian groups that are Veech groups of dilation surfaces is still unknown.

\subsection{Abelian coverings and $\Aff(M,\Sigma)$-action on $H_1(M,\Sigma;\Z)$}
\label{ssec:AffineGroupsOfAbelianCoverings}
The lifting criterion, as stated in Theorem~\ref{LiftingCriterion}, is defined
in terms of the action of $\Aff(M)$ on the fundamental group. In the case of
Abelian coverings this can naturally be reformulated in terms of the action
of $\Aff(M,\Sigma)$ on the relative homology space $H_1(M, \Sigma; \Q)$.
The central theme of this section and the next one is to study this action.

First, let's mention that the properties of the action of the affine group on
$H_1(M; \Z)$ are tightly linked to dynamical properties of translation flows or
foliations in $\widetilde{M}$ (see for example Theorem~\ref{thm:KZPAForPanovPlanes}
in Chapter~\ref{CH:INTRODUCTION}). Dynamical properties of Abelian coverings will be investigated
in a (second) forthcoming volume~\cite{DHV2}
dedicated mostly to dynamical aspects of infinite-type translation surfaces.

Recall from Lemma~\ref{lem:AbelianCoveringCorrespondence} that an Abelian covering of a translation surface $M$ at most branched over $\Sigma$ is determined by a subspace $V \subset H_1(M, \Sigma; \Q)$.
On the other hand, the group $\Aff(M,\Sigma)$ has an induced action
on $H_1(M, \Sigma; \Q)$. As a consequence of the lifting criterion
(Theorem~\ref{LiftingCriterion}) we have the following.
\begin{tcbtheorem}{Abelian lifting criterion}{AbelianLiftingCriterion}
Let $M$ be a compact translation surface and $\Sigma \subset M$ a finite set.
Let $V$ be a vector subspace of $H_1(M,\Sigma; \Q)$ and $p: \widetilde{M} \to M$ the
Abelian translation covering it determines. Then an element $f \in \Aff(M,\Sigma)$ lifts to $\widetilde{M}$ if and only
if $f_* V = V$.
\end{tcbtheorem}

In Example~\ref{exa:CharacteristicCoverings}, we already encountered an Abelian covering. Namely, the covering determined by the vector subspace $V=H_1(M,\Sigma;\Q)$
which corresponds to the derived subgroup $\delta_1(\pi_1(M \setminus \Sigma))$
of the fundamental group $\pi_1(M \setminus \Sigma)$. In this situation, all
elements of $\Aff(M,\Sigma)$ lift to $\widetilde{M}$ (see
Corollary~\ref{cor:LiftingCharacteristicCoverings}).

In view of Theorem~\ref{thm:AbelianLiftingCriterion}, studying the affine groups of Abelian translation
covers of a compact surface $M$ amounts to studying the stabilizers of vector subspaces of $H_1(M, \Sigma; \Q)$ under the action of the affine group $\Aff(M,\Sigma)$.
An important subspace that is stabilized by $\Aff(M, \Sigma)$ is the kernel of the holonomy
map, which we define as follows:
\begin{tcbdefinition}{}{Kerhol}
Let $M=(X,\omega)$ be a finite-type translation surface and $\Sigma$ be a finite subset in $M$ (possibly empty). We define $\kerhol(M, \Sigma)$ as the kernel of the holonomy map
\[
\begin{array}{ccc}
{\rm hol}: H_1(M, \Sigma; \Q) & \to & \C \\
c & \mapsto & \int_c \omega
\end{array}.
\]
When $\Sigma = \emptyset$, we will simply use the simpler notation $\kerhol(M)$ for
$\kerhol(M, \emptyset)$.
\end{tcbdefinition}
Let us emphasize that in the definition of $\kerhol$ we used the homology with
rational coefficients $H_1(M, \Sigma; \Q)$ and not $H_1(M, \Sigma; \R)$. By the rank-nullity theorem, the kernel of $\hol: H_1(M, \Sigma; \R)\to\C$ has always codimension 2. On the other hand, as we see
in Theorem~\ref{thm:CohomologyNonElementaryVeechGroup}, the codimension of
$\kerhol(M)$ is intimately related to the arithmetic properties of the
holomy vectors of $M$.

\begin{tcbexercise}{}{InvarianceKerhol}
Let $M$ be a finite-type translation surface and $\Sigma \subset M$ a finite subset.
\begin{enumerate}
\item Prove that $\kerhol(M,\Sigma)$ is invariant under the action of $\Aff(M,\Sigma)$.
\item Prove that all elements of $\Aff(M, \Sigma)$ lift to the Abelian covering
defined by $V = \kerhol(M, \Sigma)$ (\emph{hint:} use Theorem~\ref{thm:AbelianLiftingCriterion}).
\end{enumerate}
\end{tcbexercise}

Most Abelian coverings that we study in this book are defined by subspaces of $\kerhol(M,\Sigma)$. This is the case of the
infinite staircase: as it can be seen in Figure \ref{fig:staircase1} in Example \ref{exa:Zcoverings} this $\Z$-covering
is defined by the cycle $B-A \in H_1(M,\Sigma;\Z)$ and $\hol(A) = \hol(B) = 1$. It is also the case, for
the wind-tree models $W_{a,b}\to X_{a,b}$ for any $(a,b)\in(0,1)\times(0,1)$. Let us consider Figure~\ref{fig:WindtreeQuotientSurface}
in which the cycles $Z_{ij}$ and $Y_{ij}$ appear. Let us define the two absolute cycles
\begin{equation}
\label{eq:WindTreeCycleCovering}
Y = Y_{00} + Y_{10} - Y_{01} - Y_{11} \quad \text{and} \quad Z = Z_{00} + Z_{01} - Z_{10} - Z_{11}
\end{equation}
in $H_1(X_{a,b}; \Z)$. The subspace of $H_1(X_{a,b}; \Q)$ generated by $Y$ and $Z$ defines the unramified covering $W_{a,b} \to X_{a,b}$.
To see that it belongs to $\kerhol(X_{a,b})$ it is enough to notice that
\[
\forall i,j \in \{0,1\}, \quad \hol(Z_{i,j}) = 1 \quad \text{and} \quad \hol(Y_{i,j}) = \sqrt{-1}.
\]

As shown in the following exercise, half-translation surfaces (or quadratic differentials)
provides examples of $\Z^d$-coverings defined by subspaces in $\kerhol$.
\begin{tcbexercise}{}{HalfTranslationKerhol}
Let $M$ be the original example of D. Panov as illustrated in the left-hand side of Figure~\ref{fig:OriginalPanovPlane} and $p:\widetilde{M}\to M$ the $\Z^2$-covering that results from pulling the half-translation\footnote{Recall that a half-translation surface is defined as a pair $(X,q)$ where $X$ is a Riemann surface and $q$ a quadratic differential which is not identically zero. See also Definitions~\ref{def:HalfTranslationSurfaceConstructive} and~\ref{def:FlatSurfaceGeometric}} structure on $M$ to its universal cover.
\begin{enumerate}
\item Show that we have a commutative square
\begin{equation}
    \label{Eq:QDPanov}
\begin{CD}
\widetilde{M}' @>>> \widetilde{M} \\
@Vp'VV @VVpV \\
M' @>>> M
\end{CD}
\end{equation}
where $M'$ and $\widetilde{M}'$ are respectively the double orientation coverings
of $M$ and $\widetilde{M}$ and $p'$ is a $\Z^2$-covering. Show that the $\Z^2$-covering $p'$ is defined by a subspace of $\kerhol(M', \Sigma')$
where $\Sigma'$ is the union of the preimages of $\Sigma$, the preimages of the
zeros of odd order in $M$, and the preimages of poles in $M$.
\item (\emph{General case}).  Let now $M$ be a half-translation surface and $p: \widetilde{M} \to M$ a
$\Z^d$-half-translation covering branched over $\Sigma$ such that
$\widetilde{M}$ is not a translation surface (\eg Panov planes from~\ref{ssec:PanovPlanes}).
Show that we have a commutative diagram as in (\ref{Eq:QDPanov}), where $M'$ and $\widetilde{M}'$ are respectively the double orientation coverings
of $M$ and $\widetilde{M}$ and $p'$ is a $\Z^d$-covering. Show that the $\Z^d$-covering $p'$ is defined by a subspace of $\kerhol(M', \Sigma')$
where $\Sigma'$ is the union of the preimages of $\Sigma$, the preimages of the
zeros of odd order in $M$, and the preimages of poles in $M$.
\end{enumerate}
\end{tcbexercise}

In the following exercise we show that for a square-tiled surface $M$, the
space $\kerhol(M, \Sigma)$ admits a natural complement in $H_1(M,\Sigma; \Q)$.
This result is generalized in
Theorem~\ref{thm:CohomologyNonElementaryVeechGroup} and
Lemma~\ref{lem:KerholAndCovers} to any finite-type translation surface with
non-elementary Veech group.
\begin{tcbexercise}{}{TorusKerhol}
Let $\T^2 = \R^2 / \Z^2$ be the square torus and $\Sigma \subset \T^2$ a finite
set that contains the origin $0$.
\begin{enumerate}
\item Show that $\dim \kerhol(\T^2, \Sigma) \leq |\Sigma| - 1$ with
equality if and only if $\Sigma$ are rational points, that is $\Sigma \subset \Q^2 / \Z^2$.
\item Prove that if $\Sigma \subset \Q^2 / \Z^2$ then there is a decomposition
\[
H_1(\T^2, \Sigma; \Q) = H_1(\T^2; \Q) \oplus \kerhol(\T^2, \Sigma).
\]
\item Let $p: M \to \T^2$ be a finite-type square-tiled surface.
Let $U_M \subset H_1(M; \Q)$ be the
subspace generated by the preimages under $p$ of closed curves in $\T^2$.
Show that $p_*: U_M \to H_1(\T^2; \Q)$ is an isomorphism and
that we have a decomposition
\begin{equation}
\label{eq:TangentSpaceOrigami}
H_1(M; \Q) = U_M \oplus \kerhol(M)
\end{equation}
that is orthogonal with respect to the algebraic intersection form on $H_1(M; \Q)$.
\item Show that $U_M$ defined in the previous question is $\Aff(M)$-invariant.
\item Let $M$ be a finite-type square-tiled surface as before. Let $\Sigma' \subset M$ be
a finite set contained in the preimage of the rational
points of the torus which does not necessarily contains the conical singularities of $M$. Show that
\begin{equation}
\label{eq:OrigamiHomologyDecomposition}
H_1(M, \Sigma'; \Q) = U_M \oplus \kerhol(M, \Sigma').
\end{equation}
\end{enumerate}
\end{tcbexercise}

We now study in full generality a finite-type surface $M$ with non-elementary Veech group.
It will clarify the special role that $\kerhol(M,\Sigma)$ plays when considering the
Veech groups of a $\Z^d$-covering of $M$. The discussion requires
some technical concepts from compact translation surfaces. A reader that is not familiar
with these notions should not be afraid of not understanding everything as most of
the results developed here will not be used intensively in the following sections.

Recall that a subgroup of $\SL(2,\R)$ is called non-elementary if
it contains two hyperbolic elements with disjoint axes (see Appendix~\ref{app:FuchsianGroups}). By extension, we say that
the affine group $\Aff(M,\Sigma)$ of a finite-type translation
surface is non-elementary if its image $\Gamma(M,\Sigma)$ by the derivative homomorphism
$D: \Aff(M,\Sigma) \to \SL(2,\R)$ is non-elementary (since we assume that
$M$ is finite-type, the affine elements necessarily preserve the area and hence
image is necessarily in $\SL(2,\R)$).

An important source of examples of finite-type translation surface $M$
with non-elementary Veech groups are the ones obtained from the
Thurston-Veech construction (see Section~\ref{ssec:ThurstonVeechConstructionsIntro}).
In this case, the Veech group $\Gamma(M)$ contains two parabolic elements
with distinct fixed points in $\partial \H^2$. Among these Thurston-Veech
surfaces, there are some very exceptional ones for which the Veech group $\Gamma(M)$
is of finite covolume: these are the so-called Veech surfaces. Finite-type
square-tiled surfaces (see Example~\ref{exa:SquareTiledSurfaces}) are Veech
surfaces and their Veech groups are commensurable with $\SL(2,\Z)$.
An other example of surface with non-elementary Veech group but which does not
come from a Thurston-Veech construction is the
Arnoux-Yoccoz surface~\cite{HubertLanneauMoeller09}. We terminate our list with mentioning a general construction from~\cite{HubertSchmidt04}, for which we need some
definitions.

\begin{tcbdefinition}{}{DEF:holonomyfield}
Let $M$ be a finite-type translation surface and $v_1,v_2$ two non-parallel vectors in $\hol(H_1(M;\Z))$. The smallest subfield $k$ of $\R$ for which any element in $\hol(H_1(M;\Z))$ can be written as $av_1+bv_2$, $a,b\in k$ is called the \emph{holonomy field} of $M$.
\end{tcbdefinition}

Equivalently, let $T^{id}_M := \R \Re(\omega) + \R \Im(\omega) \subset H^1(M; \R)$ be the
\emphdef[tautological plane (translation surface)]{tautological plane}. Then $k$ as defined above is also the field of definition of $T^{id}_M$, that is,
the smallest field such that one can write $T^{id}_M$ as a kernel of a matrix with
coefficients in $k$ (in a $\Q$-basis of $H^1(M; \Q)$). A direct consequFence of these definitions is that the Veech group of a finite type translation surface $M$ with holonomy field $k$ is a subgroup of $\SL(2,k)$.

\begin{tcbtheorem}{\cite{KenyonSmillie00}}{holonomyfieldtracefield}
Let $M$ be a finite-type translation surface with holonomy field $k$ and suppose that there exists a hyperbolic element $f\in\Aff(M)$, then $k = \Q[\tr(D(f))]$.
\end{tcbtheorem}

For a more detailed discussion on holonomy fields and its applications we refer the reader also to~\cite{GutkinJudge00}, \cite{HubertLanneau06}
and~\cite{McMullen03}. For a discussion about how the concept of holonomy field can be extended to infinite-type surfaces and some new phenomena that appear in this context see Ñ \cite{SchmithuesenValdez14}. Remark that if $M$ is a square-tiled surface then its holonomy field is $k=\Q$ and $U_M$, the subspace of $H_1(M;\Q)$ defined in Exercise~\ref{exo:TorusKerhol}, is the Poincar\'e dual\footnote{If $D:H^1(M)\to H_1(M)$ is the isomorphism given by Poincar\'e-duality, then we say that a subspace $W\subset H_1(M)$ is the Poincar\'e dual of $V\subset H^1(M)$ if $D^{-1}(W)=V$.} of the tautological plane $T^{id}_M $.

\begin{tcbdefinition}{}{RationalPoint}
Let $M$ be a compact translation surface with holonomy field $k$ and $\Sigma$ be its
set of conical singularities. Let
$\Lambda = \{\int_\gamma \omega: \gamma \in H_1(M, \Sigma; \Q)\}$.
A point in $M$ is called \emphdef[$k$-rational point]{$k$-rational}
if the holonomy of any path from $x$ to a point in $\Sigma$ belongs to $\Lambda$
\end{tcbdefinition}

\begin{tcbdefinition}{}{PeriodicPoint}
Let $M$ be a surface with non-elementary Veech group.
A point $x$ in $M$ is called
\emphdef[periodic point (translation surface)]{periodic} if its stabilizer
in $\Aff(M)$ has finite index (or equivalently, if its
orbit is finite).
\end{tcbdefinition}
For square-tiled surfaces, the rational points are the elements of $\Q^2/\Z^2$
which are all periodic. This justifies our terminology. For a surface with
non-elementary Veech group the periodic points are rational points, however
the converse is not true: in a Veech surface the set of periodic
points is finite~\cite{GutkinHubertSchmidt03,Moeller06-periodic}. The fact
that there exist rational points which are not periodic can be used to construct
finite-type translation surfaces whose Veech group do not admit a finite set
of generators, see~\cite{HubertSchmidt04}.

We now discuss how the homology of a translation surface with non-elementary Veech group decomposes, but first let us discuss and recall some useful facts.

For each embedding $\sigma: k \to \R$ or a pair of complex conjugate embeddings
$\sigma, \overline{\sigma}: k \to \C$ we can define a Galois conjugate
$T_M^\sigma \subset H^1(M; \R)$ of the tautological plane by applying $\sigma$
on the defining equations of $T_M^{id}$. For real embeddings $T_M^\sigma$
has (real) dimension 2 while for a pair of complex conjugate embeddings it has (real) dimension 4.

Recall also that for any compact translation surface $M$ of genus $g$ the algebraic intersection number:
$$
\langle\cdot,\cdot\rangle : H_1(M,\Z)\times H_1(M,\Z)\to\Z
$$
extends uniquely to a symplectic form on $H_1(M,\R)$. Moreover, if $\eta_a,\eta_b$ are the Poincar\'e duals of $a,b\in H_1(M,\Z)$ then $\langle a,b\rangle = \int_M \eta_a\wedge\eta_b$. On the other hand, the natural action of $\MCG(M)$, the mapping class group of $M$ on $H_1(M,\Z)$ defines a (surjective) \emph{symplectic} representation:
$$
\MCG(M)\to \Sp(H_1(M,\Z))
$$
We refer the reader to Chapter 6 in~\cite{FarbMargalit} for a more detailed dicussion.

\begin{tcbtheorem}{}{CohomologyNonElementaryVeechGroup}
Let $M$ be a compact translation surface whose affine group is
non-elementary and $k$ be its holonomy field. Then
each subspace $T_M^\sigma$ is symplectic and $\Aff(M)$-invariant.
There is direct sum decomposition orthogonal with respect to
the intersection form
\begin{equation}
\label{eq:DefinitionTangent}
T_M := \bigoplus_{\sigma \in \Emb(k)} T_M^\sigma
\end{equation}
where $\Emb(k)$ is the set of real embeddings and pairs of complex conjugate
embeddings of $k$.

The Poincar\'e dual of $T_M$ is a subspace of $H_1(M; \R)$ defined over $\Q$ and the
anihilator of $T_M$ is exactly $\kerhol(M)$.
As a consequence $\kerhol(M)$ is
symplectic and, if we henceforth identify $T_M$ to its Poincar\'e dual in $H_1(M; \Q)$ (and keep the notation $T_M$) we have a direct sum decomposition orthogonal with respect
to the algebraic intersection form
\begin{equation}
\label{eq:SumTangentKerhol}
H_1(M; \R) = T_M \oplus (\kerhol(M) \otimes \R)
\end{equation}
Any subspace $V \subset H_1(M; \R)$ that has a non-elementary
stabilizer in $\Aff(M)$ decomposes as
\begin{equation}
\label{eq:NonElementaryInvariantDecomposition}
V = \bigoplus_{\sigma \in I} T_M^\sigma \oplus (\kerhol(M) \otimes \R \cap V)
\end{equation}
where $I \subset \Emb(k)$.
\end{tcbtheorem}

One of the main ingredients in the proof of the above result is the following Theorem of Thurston:

\begin{tcbtheorem}{}{TranslationPseudoAnosovHomologyAction}
Let $M$ be a finite-type translation surface and $f \in \Aff(M)$
a hyperbolic element. Then the dominant eigenvalue of  $D(f)$
is also the dominant eigenvalue of the action of $f$ in
$H_1(M; \Z)$. Moreover, this eigenvalue is simple.
\end{tcbtheorem}

\begin{proof}
We provide a proof only for the case when $k$ is totally real. The
general case introduce some technical details
that would obfuscate the proof.

First of all, the
action of $\Aff(M)$ on the tautological plane $T_M^{id} = \R \Re(\omega) \oplus
\R \Im(\omega)$ identifies with the derivative on $\R^2$. Indeed, by definition of
the affine action
\[
\begin{pmatrix} f^* \Re(\omega) \\ f^* \Im(\omega) \end{pmatrix}
=
D(f) \begin{pmatrix} \Re(\omega) \\ \Im(\omega) \end{pmatrix}.
\]
The restriction of the symplectic form $\langle \eta_1, \eta_2 \rangle = \int_M \eta_1 \wedge \eta_2$
on the tautological plane is non-degenerate since $\int_M \Re(\omega) \wedge \Im(\omega)$ is the
area of the surface

More generally, for each embedding $\sigma: k \to \R$, the induced action of $\Aff(M)$ on
$T_M^\sigma$ identifies with the action of $\sigma \circ D$ on $\R^2$. This shows
that each $T_M^\sigma$ is symplectic and $\Aff(M)$-invariant.

We now prove that the sum in~\eqref{eq:SumTangentKerhol} is direct and symplectic orthogonal. Let us consider the fact that the holonomy field is the field
generated by the traces of the matrices $D(f)$ for $f \in \Aff(M)$, see Theorem~\ref{thm:holonomyfieldtracefield}. Because
$\Aff(M)$ is non-elementary, it admits a hyperbolic element.
Its trace $t = \tr(D(f))$ is a generator of $k$ and hence the
embeddings $\sigma(t)$ are pairwise distinct. The eigenvalues
of $D(f)$ is a pair $\lambda$, $1/\lambda$ such that
$t = \lambda + 1/\lambda$. Moreover $\lambda$ is simple as
an eigenvalue of the induced action $f_*: H_1(M;\R) \to H_1(M; \R)$, see Theorem~\ref{thm:TranslationPseudoAnosovHomologyAction}. Hence
$T_M^{Id} = E_\lambda(f_*) \oplus E_{1/\lambda}(f_*)$ where $E_\lambda(f_*)$
and $E_{1/\lambda}(f_*)$ are the eigenspaces of $f_*$ with eigenvalues
$\lambda$ and $1/\lambda$.
This is also true for each Galois conjugate:
$T_M^{\sigma} = E_{\sigma(\lambda)}(f_*) \oplus E_{1/\sigma(\lambda)}(f_*)$.
Because the eigenspaces are in direct sums, the sum in~\eqref{eq:DefinitionTangent}
is direct. Since the matrix $f_*$ is symplectic, this sum is symplectic orthogonal\footnote{If $M$ is a symplectic matrix (preserving a symplectiv form $\omega$) and $\lambda,\mu$ are eigenvalues such that $\lambda\mu\neq 1$ then for every $\lambda$-eigenvector $v$ and $\mu$-eigenvector $w$ we have that $\omega(v,w)=\omega(Mv,Mw)=\lambda\mu\omega(v,w)$ and hence $\omega(v,w)=0$.}.

Since $\kerhol(M)$ is a $\Q$-subspace of $H_1(M; \R)$ it
is also anihilated by the $T_M^\sigma$. The direct sum
in equation~\eqref{eq:SumTangentKerhol} follows.

Now we want to prove the last part of the Theorem. Let $V$ be a subspace
of $H^1(M; \R)$ with a non-elementary stabilizer $G$ under the action of
$\Aff(M)$. Because $G$ is non-elementary, it contains two hyperbolic elements $f_1$
and $f_2$ such that their derivatives $D(f_1)$ and $D(f_2)$ have no common
fixed points on $\partial \H^2$ (see Theorem~\ref{thm:ElementaryFuchsianGroup}).
Equivalently, the derivatives $D(f_1)$ and $D(f_2)$ have no common eigenspace.
Hence, the only invariant subspaces of $\R^2$ preserved by both $D(\Aff(M))$
are only the trivial subspaces $0$ and $\R^2$. The same property holds for
any $T_M^\sigma$.

Now the invariant subspace $V$ must decompose along the eigenspaces of any
element in its stabilizer. Since the decompositions in eigenspaces of $f_1$
and $f_2$ lead to distinct decompositions of $T_M^\sigma$, we obtain
the decomposition~\eqref{eq:NonElementaryInvariantDecomposition}.
\end{proof}

\begin{tcbcorollary}{}{}
Let $M$ be finite type translation surface and $\Sigma \subset M$ a finite set
such that $\Aff(M, \Sigma)$ is non-elementary. Let $T_M = \bigoplus_{\sigma \in \Emb(k)} T_M^\sigma$
as in Theorem~\ref{thm:CohomologyNonElementaryVeechGroup}

Let $V$ be a vector subspace of $H_1(M, \Sigma; \Q)$. If the
stabilizer of $V$ in $\Aff(M,\Sigma)$ is non-elementary
then either $V$ contains\footnote{Here we are abusing notation by writing $T_M\subset H_1(M,\Q)$ for the Poincar\'e dual of $T_M$ as defined in Theorem~\ref{thm:CohomologyNonElementaryVeechGroup}. } $T_M$ or $V$ is contained in $\kerhol(M, \Sigma)$.

In the particular case of $\Z$-coverings, that is $\dim_\Q V = 1$, only the second
possibility occurs.
\end{tcbcorollary}

\begin{proof}
The decomposition~\eqref{eq:SumTangentKerhol}
from Theorem~\ref{thm:CohomologyNonElementaryVeechGroup}
is also valid for relative homology
\begin{equation}
\label{eq:RelativeHomologyDecomposition}
H_1(M,\Sigma; \Q) = \bigoplus_{\sigma \in \Emb(k)} T_M^\sigma \oplus \kerhol(M,\Sigma).
\end{equation}
The reason is that the class of $\Re(\omega)$ and $\Im(\omega)$ are well
defined in $H^1(M \setminus \Sigma; \Q)$ which by Poincar\'e duality
can be identified to subspaces of $H_1(M,\Sigma;\Q)$.

Now, the second part of Theorem~\ref{thm:CohomologyNonElementaryVeechGroup} shows that an
invariant subspace must respect the decomposition~\eqref{eq:RelativeHomologyDecomposition}.
Since $V$ is assumed to be defined over $\Q$ it must either contain the whole sum
$\bigoplus_{\sigma \in \Emb(k)} T_M^\sigma$ (that is the smallest subspace defined
over $\Q$ that contains $T_M^{id}$) or be contained in $\kerhol(M,\Sigma)$

With the additional hypothesis that $\dim_Q V = 1$, it can not contain
$T_M$ which has non-zero even dimension.
\end{proof}

Let us mention two important elementary results about $\kerhol(M,\Sigma)$.
The first one concerns the behavior under finite coverings.

\begin{tcblemma}{}{KerholAndCovers}
Let $M$ be a compact translation surface with non-elementary Veech
group, $k$ its holonomy field and
$\Sigma \subset M$ a finite set of $k$-rational points. Let $p: \widetilde{M} \to M$ a finite
covering branched over $\Sigma$ and $\widetilde{\Sigma} := p^{-1}(\Sigma)$. Then
$p_*: H_1(\widetilde{M}, \widetilde{\Sigma}; \Q) \to H_1(M, \Sigma; \Q)$ induces
an isomorphism between $T_{\widetilde{M}}$ and $T_M$ as defined in Theorem~\ref{thm:CohomologyNonElementaryVeechGroup}; moreover
$\codim \kerhol(M, \Sigma) = \codim \kerhol(\widetilde{M}, \widetilde{\Sigma})$
(where the codimensions are respectively measured in $H_1(M, \Sigma; \Q)$
and $H_1(\widetilde{M}, \widetilde{\Sigma}; \Q)$).
\end{tcblemma}
The second lemma states that the surfaces for which $\kerhol(M)$ is maximal
are exactly the square-tiled surfaces.
\begin{tcblemma}{}{KerholAndTori}
Let $M$ be a compact translation surface and $\Sigma \subset M$ a (possibly
empty) finite set of points.  Then $\codim \kerhol(M, \Sigma) = 2$ (where the
codimension is taken inside $H_1(M, \Sigma; \Q)$) if and only if $p:M \to \T^2$ is a
square-tiled surface and $\Sigma \subset p^{-1}(0)$.
\end{tcblemma}

\begin{proof}[Proof of Lemma~\ref{lem:KerholAndCovers}]
Let $\omega$ and $\widetilde{\omega}$ be respectively the one forms on $M$ and $\widetilde{M}$. By definition of translation covering, $\widetilde{\omega} = p^* \omega$. It follows that $T_{\widetilde{M}}$ and $T_M$ are isomorphic.

Let us define
\[
L := \left\{ \int_\gamma \omega:\ \gamma \in H_1(M, \Sigma; \Z) \right\}
\qquad \text{and} \qquad
\widetilde{L} := \left\{ \int_{\gamma} \widetilde{\omega}:\ \gamma \in H_1(\widetilde{M}, \widetilde{\Sigma}; \Z) \right\}.
\]
By paths projection and paths lifting we have that
\[
\widetilde{L} \subset L \subset \frac{1}{\deg(p)} \widetilde{L}.
\]
In particular $L \otimes \Q = \widetilde{L} \otimes \Q$.
Now, by definition of $\kerhol(M,\Sigma)$ and $\kerhol(\widetilde{M},\widetilde{\Sigma})$, we have that
$\codim \kerhol(M,\Sigma) = \rk_\Z L$ and
$\codim \kerhol(\widetilde{M},\widetilde{\Sigma}) = \rk_\Z \widetilde{L}$ from which the lemma follows.
\end{proof}

\begin{proof}[Proof of Lemma~\ref{lem:KerholAndTori}]
As in the proof Lemma~\ref{lem:KerholAndCovers}, let us consider
\[
L := \left\{\int_\gamma \omega:\ \gamma \in H_1(M, \Sigma; \Z) \right\} \subset \C.
\]
Fix a reference point $x_0 \in \Sigma$ if $\Sigma$ is non-empty and any point in $M$ otherwise.
Then consider the period mapping
\[
P: \begin{array}{ccc}
M & \to & \C / L \\
x & \mapsto & \int_{x_0}^x \omega.
\end{array}.
\]
The map is well defined since two paths from $x_0$ to $x$ differ by an
element in absolute homology.

Now $\codim \kerhol(M,\Sigma) = 2$ if and only if $L$ is a lattice in $\C$.
Hence if $\codim \kerhol(M,\Sigma) = 2$, the map $P$ provides a translation covering with
a torus target. And by construction, all points of $\Sigma$ are mapped to $0$ in the
quotient.

The proof of the converse is a particular case of Lemma~\ref{lem:KerholAndCovers}.
\end{proof}

A ``generic" surface $M$ in a given stratum $\cH(\kappa)$ of Abelian differentials
satisfies $\dim\kerhol(M) = 0$. However, Lemmas~\ref{lem:KerholAndCovers}
and~\ref{lem:KerholAndTori} show that when $M$ comes from a covering
then $\kerhol(M)$ is non-trivial.

\begin{tcbexercise}{}{LAffineAction}
Let $L$ be the L-shaped surface made from 3 squares as in Figure~\ref{fig:LShapedExo} (after identifying opposite sides using translations)
and let $\alpha = \gamma_2 - \gamma_1$ and $\beta = \gamma_3 - \gamma_4$ (as elements
of the absolute homology $H_1(L; \Z)$).
\begin{figure}[H]
\begin{center}
\includegraphics[scale=0.8]{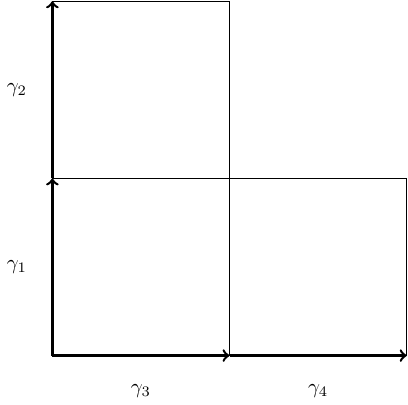}
\end{center}
\caption{}
\label{fig:LShapedExo}
\end{figure}

\begin{enumerate}
\item Show that $\kerhol(L) = \Q \alpha \oplus \Q \beta$,
\item Show that $L$ admits no translation automorphism (so that we can identify
the Veech group and the affine group).
\item Show that the following three matrices belong to the Veech group of $L$
\[
s = \begin{pmatrix}0&1\\-1&0\end{pmatrix},
\quad
h_2 = \begin{pmatrix}1&2\\0&1\end{pmatrix},
\quad
v_2 = \begin{pmatrix}1&0\\2&1\end{pmatrix}.
\]
(actually these 3 matrices generate the Veech group)
\item Let $\rho: \Aff(M) \to \GL(\kerhol(M))$. Show that in the basis $(\alpha, \beta)$
of $\kerhol(M)$ one has
\[
\rho(s) = \begin{pmatrix}0&1\\-1&0\end{pmatrix},
\quad
\rho(h_2) = \begin{pmatrix}1&1\\0&1\end{pmatrix},
\quad
\rho(v_2) = \begin{pmatrix}1&0\\1&1\end{pmatrix}.
\]
\item Conclude that the subgroup of $\Aff(L)$ that consists of elements that act
trivially on $\kerhol(L)$ is a non-trivial normal subgroup of $\Aff(L)$ of infinite
index. Does it contain a parabolic element?
\end{enumerate}
\end{tcbexercise}

\subsubsection{Affine action on absolute and relative homologies}
\label{sssec:AbsoluteVSRelative}
In this section we study the differences between the action of affine groups on
absolute and relative homology groups.

Let $M$ be a finite-type translation surface and $\Sigma \subset M$ be a finite
set so that $\Aff(M,\Sigma)$ is non-elementary. Recall that:
\begin{itemize}
  \item both $\kerhol(M)$
 and $\kerhol(M,\Sigma)$ are $\Aff(M,\Sigma)$-invariant subgroups of
 $H_1(M,\Sigma;\Q)$ (see Exercise~\ref{exo:InvarianceKerhol}).
  \item To any discrete (\ie Fuchsian) subgroup $\Gamma$ of $\SL(2,\R)$ one can associate its limit set $\Lambda(\Gamma)\subset\partial\H$: this is the set of limit points of $\Gamma z$ for any $z\in\H$. For more details see Appendix~\ref{app:FuchsianGroups}.
\end{itemize}

\begin{tcbtheorem}{\cite{HooperWeiss09}}{KerPsiVSkerPsi0}
Let $M$ be a finite-type translation surface and let
$$\psi: \Aff(M,\Sigma) \to GL(\kerhol(M,\Sigma))$$
and $$\psi_0: \Aff(M,\Sigma) \to GL(\kerhol(M))$$ be the induced actions. Then, $\ker(\psi)$ is a normal
subgroup of $\ker(\psi_0)$ and the quotient is virtually Abelian of rank at most $r_0(r-r_0)$
where $r_0 = \dim(\kerhol(M))$ and $r = \dim(\kerhol(M,\Sigma))$.

In particular, if $\Gamma(M)$ is non-elementary and $\ker(\psi)$ is
non-trivial then
$$
\Lambda (D\ker \psi) = \Lambda (D\ker \psi_0) = \Lambda \Gamma(M),
$$
where $D:\Aff(M,\Sigma)\to\SL(2,\R)$ denotes the derivative.
\end{tcbtheorem}

\begin{proof}
Because $\kerhol(M) \subset \kerhol(M,\Sigma)$, in an appropriate basis of
$\kerhol(M,\Sigma)$, $\psi$ has the form
\begin{equation}
\label{eq:ActionHomologyVSRelativeHomology}
\psi(f) =
\begin{pmatrix}
\psi_0(f) & \alpha(f) \\
0         & \sigma(f)
\end{pmatrix}
\end{equation}
where $\sigma: \Aff(M,\Sigma) \to \GL(V)$ is the action induced by permutations
on $V = \{v \in \Q^\Sigma: \sum v_i = 0\}$ and $\alpha(f)$ is some $r_0
\times (r - r_0)$ matrix with integer coefficients.

Hence, $\ker(\sigma)$ is of finite index (the image of $\sigma$ is a permutation
group) and the image of $\alpha$ is an Abelian group whose rank is given by
the dimensions of the matrix: $r_0 (r - r_0)$.

The last part of the proof concerning the limit sets is a direct application of
Lemma~\ref{lem:NormalSubgroupLimitSet} from the Appendix~\ref{Appendix:FuchsianGroups}.
\end{proof}

\begin{tcbexercise}{}{TorusAffineAction}
This exercise is a continuation of Exercise~\ref{exo:TorusKerhol}.
Recall that $\Aff(\T^2, \{0\}) \simeq \SL(2,\Z)$. Instead of considering $\Aff(\T^2, \Sigma)$ (that can admit translations) we consider subgroups of $\Aff(\T^2, \{0\})$ that we identify to subgroups of $\SL(2,\Z)$. Let $n$ be a positive integer and denote $\Sigma_n$ the set of points $x$ in $\T^2 = \R^2/\Z^2$ such that such that $nx = 0$. A point in $\Sigma_n$ is called a \emph{$n$-torsion point}.
\begin{enumerate}
\item Show that $\Sigma_n$ is a subgroup of $(\T^2, +)$ isomorphic to $\Z/n\Z \times \Z/n\Z$.
\item Show that $\SL(2,\Z)$ preserves $\Sigma_n$.
\setcounter{EnumerateCounterContinuation}{\theenumi}
\end{enumerate}
The \emph{principal congruence subgroup} of degree $n$ is
\[
\Gamma(n) :=\{
A \in \SL(2,\Z):
A \equiv Id \mod n\}
\]
\begin{enumerate}
\setcounter{enumi}{\theEnumerateCounterContinuation}
\item Show that $\Gamma(n)$ is the pointwise stabilizer of $\Sigma_n$, \ie $A s = s$ for each $s$ in $\Sigma_n$.
\setcounter{EnumerateCounterContinuation}{\theenumi}
\end{enumerate}
Recall that we have a map $\delta: H_1(\T^2, \Sigma_n; \Z) \to \Z^{\Sigma_n}$ which to a relative cycle associates its signed boundary. On the other hand $\Sigma_n$ is a group and we have a map $\Z^{\Sigma_n} \to \Sigma_n$ which to a formal linear combination of points associates its sum in $\Sigma_n$. We denote $\delta': H_1(\T^2, \Sigma_n; \Z) \to \Sigma_n$ the composition of these two maps.
\begin{enumerate}
\setcounter{enumi}{\theEnumerateCounterContinuation}
\item Show that $\delta$ is equivariant with respect to the action of $\SL(2,\Z)$.
\item Show that $\delta'$ coincide with the composition of $\hol: H_1(\T^2, \Sigma_n; \Z) \to \R^2$ with the quotient map $\R^2 \to \T^2$.
\item Deduce from the two previous items that $\Gamma(n)$ preserves $\kerhol(\T^2, \Sigma_n)$ and that for each element in $\Gamma(n)$ its induced action on $\kerhol(\T^2, \Sigma_n)$ is trivial.
\item Deduce that if $V \subset \kerhol(\T^2, \Sigma_n)$ is a vector subspace and $M$ the translation Abelian covering associated to $V$ then the affine group of $M$ contains $\Gamma(n)$.
\end{enumerate}
\end{tcbexercise}
In the above exercise, we have seen a family of examples where the affine group
$\Aff(M, \Sigma)$  acts via a finite group on $\kerhol(M, \Sigma)$. However, as shown
in the following example this is not the case in general.
\begin{tcbexample}{}{OrigamiH11-1}
We now present an example of a finite-type square-tiled surface $Z_{3,1}$ made of 4 squares in the stratum $\cH(1,1)$ of
Abelian differentials. We compute explicitly the action of $\Aff(Z_{3,1}) = \Aff(Z_{3,1},\Sigma)$ on
$H_1(Z_{3,1}; \Z)$ and $H_1(Z_{3,1}, \Sigma; \Z)$ where $\Sigma$ is the set made of the two conical
singularities of $O$. We use the notation introduced in the proof of Theorem~\ref{thm:KerPsiVSkerPsi0}.
\newsavebox{\smlmat}
\savebox{\smlmat}{$\left(\begin{smallmatrix}1&0\\1&1\end{smallmatrix}\right)$}

\begin{figure}[H]
\begin{center}
\includegraphics[scale=0.7]{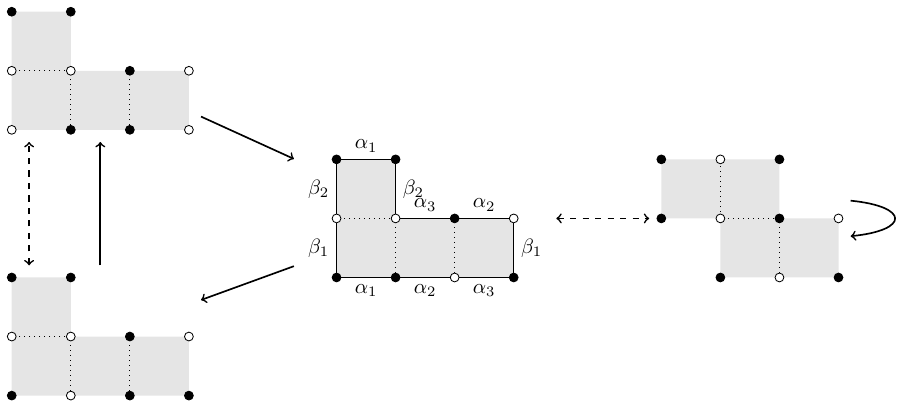}
\end{center}
\caption{The $\SL(2,\Z)$-orbit of the square-tiled surface $Z_{3,1}$ (shown here in the middle). The plain arrows correspond to the
action of the horizontal parabolic element~\usebox{\smlmat}, while the dotted ones to the rotation by $\pi/2$.}
\label{fig:SL2ZOrbit-H11}
\end{figure}
The Veech group (which in that case identifies to the affine group) is
an index 4 subgroup of $\SL(2,\Z)$ generated by the matrices
\[
g_1 = \begin{pmatrix}
1&3 \\
0&1
\end{pmatrix}
\quad
g_2=\begin{pmatrix}
1&0\\
1&1
\end{pmatrix}
\quad
g_3=\begin{pmatrix}
1&-3\\
1&-2
\end{pmatrix}.
\]
Then in the basis\footnote{Oriented in the Figure from left to right and bottom to top.} $\{2\alpha_1 - \alpha_2 - \alpha_3, 2\alpha_2 - 2\alpha_1 + \beta_2 - \beta_1, \alpha_2 - \alpha_1\}$ of $\kerhol(Z_{3,1},\Sigma)$ depicted
in Figure~\ref{fig:SL2ZOrbit-H11} the map $\psi_0$ is given by
\[
\psi_0(g_1) =
\begin{pmatrix}
1&1&0 \\
0&1&0 \\
0&0&1
\end{pmatrix}
\quad
\psi_0(g_2)=
\begin{pmatrix}
-2&1&1\\
-1&0&0\\
0&0&1
\end{pmatrix}
\quad
\psi_0(g_3)=
\begin{pmatrix}
-2&3&1\\
-1&1&0\\
0&0&1
\end{pmatrix}.
\]
One can show that for
$g = \begin{pmatrix}
37 & 60 \\
8 & 13
\end{pmatrix} = g_1 g_3 (g_2)^2 g_3 g_2$
then $\alpha(g) = (-1,0)$. Which shows that $\ker(\psi)$ and $\ker(\psi_0)$
are different. Moreover, taking its conjugate
$g' = \begin{pmatrix}
13 & 24 \\
20 & 37
\end{pmatrix} = g_2 g_3 g_2 g_1 g_3 g_2$ we have $\alpha(g') = (-1,-1)$. Hence, the Abelian group
of Theorem~\ref{thm:KerPsiVSkerPsi0} has rank $2$ in this example (the maximal possible).

Let us also mention that this example shows that the $\Aff(M,\Sigma)$-invariant subspace
$H_1(M;\R)$ inside $H_1(M,\Sigma;\R)$ does not necessarily admits a complement that is
also invariant.
\end{tcbexample}


\subsection{Lifting parabolic elements in $G$-coverings}
\label{ssec:LiftingParabolicElements}
In this section we study parabolic elements in more detail. Because a parabolic element
in the affine group of a finite-type translation surface is associated to a cylinder decomposition
their properties are often much easier to understand than for hyperbolic elements.
As we will see in the next chapters, cylinders are extremely useful to study recurrence
and ergodicity of coverings.

Cylinders and strips in translation surfaces were defined in Definition~\ref{def:CylindersAndStrips}. First, let's consider the general question of preimages of cylinders.
Let $p: \widetilde{M} \to M$ be a translation covering. Then, the preimages in $\widetilde{M}$ of cylinders in $M$ are either cylinders or strips. The following statement
gives a precise description.
\begin{tcblemma}{}{LiftingCylindersInCoverings}
Let $p: \widetilde{M} \to M$ be a $G$-covering branched over $\Sigma$ and determined
by the kernel of the morphism $\rho: \pi_1(M \setminus \Sigma) \to G$. Let $C$ be
a maximal cylinder in $M \setminus \Sigma$ with core curve $\gamma$ seen as a conjugacy
class in $\pi_1(M \setminus \Sigma)$.
Let $k \in \{1, 2, \ldots, \} \cup \{\infty\}$ be the order of $\rho(\gamma)$ in $G$. Then,
if $k$ is finite, the preimage of $C$ in $\widetilde{M}$ consists of a disjoint collection of copies of cylinders
with same height and whose circumference is $k$-times longer than the circumference of $C$. If $k$
is infinite, the preimage of $C$ consists of a disjoint union of strips of the same height as $C$.
\end{tcblemma}

\begin{proof}
The argumentation is based on the path-lifting property from the theory
of covering spaces. Namely, let us fix a basepoint $x$ on the curve $\gamma$.
Then this curve can be
lifted from any point $\widetilde{x}$ in the preimage of $x$. The endpoint
of the lifting is then $\widetilde{x} \cdot \gamma$ where $\cdot$ represents
the right action of the fundamental group $\pi_1(M \setminus \Sigma, x)$
on the fiber of $x$. Repeating this lifting from $\widetilde{x} \cdot \gamma$
we end up in $\widetilde{x} \cdot (\gamma)^2$, etc. This procedure comes back
for the first time to $\widetilde{x}$ in $k$ steps if and only if the order of
$\rho(\gamma)$ is $k$. This proves the part of the Lemma concerning the circumferences.
The preimage of a cylinder $C$ in $\widetilde{M}$ has the same height because
the lift of a saddle connection in the boundary of $C$ is a saddle connections (of the same length).
\end{proof}

\begin{tcbdefinition}{}{CompletelyPeriodicDirection}
Let $M$ be a translation or half-translation surface. We say that $\theta\in\R/2\pi\Z$ is:
\begin{enumerate}
\item \emphdef[completely periodic (direction)]{completely periodic} if there exists a collection of maximal cylinders $\{C_i\}_{i\in I}$ parallel to $\theta$ such that $\cup_{i\in I} C_i$ is dense in $M$. Such a collection of cylinders is called a \emphdef{cylinder decomposition} in direction $\theta$.
\item an \emphdef[affine multitwist direction (translation surface)]{affine multitwist direction} if it is completely
periodic and there exists an affine multitwist\footnote{See Definition~\ref{def:AffineMultitwist}.} $f$ in $\Aff(M)$ that acts as a product
of Dehn twists in each cylinder $C_i$, $i\in I$.
\item A \emphdef[parabolic direction (translation surface)]{parabolic direction} if there exists $f\in\Aff(M)$ parabolic such that $D(f)$ fixes the direction $\theta$.
\end{enumerate}
\end{tcbdefinition}

\begin{tcbremark}{}{}
For finite-type translation surface $M$, Lemma~\ref{lem:ParabolicCylinderDecomposition} shows parabolic directions and affine multitwist directions are the same thing. Moreover, the existence of an affine multitwist for a given direction is equivalent to the commensurability of the moduli of the corresponding cylinders.

As shown in Example~\ref{exa:BandsInfiniteStaircase} below, a parabolic direction in a translation surface is not necessarily an affine multitwist direction. However, as we see in Section~\ref{sec:HooperThurstonVeechConstruction}, the condition for a completely periodic direction to be stabilized by an affine multitwist can still be formulated in terms of moduli of the cylinders.
\end{tcbremark}

As a consequence of Lemma~\ref{lem:LiftingCylindersInCoverings} we have the
following general statement for lifting parabolic affine automorphisms.
\begin{tcbcorollary}{}{LiftingParabolicsViaCylinders}
Let $M$ be a finite-type translation surface and let $p: \widetilde{M} \to M$
be a $G$-translation covering branched over $\Sigma$ defined by
a morphism $\rho: \pi_1(M \setminus \Sigma) \to G$.
Let $f \in \Aff(M,\Sigma)$ be a parabolic
element for which $C_1$, \ldots, $C_k$ are the maximal cylinders in the corresponding decomposition of $M$. Let $\gamma_i$ be the
core curve of $C_i$ considered as a conjugacy class in
$\pi_1(M \setminus \Sigma)$. If all the elements $\rho(\gamma_1)$, \ldots, $\rho(\gamma_k)$
have finite order $m_1$, \ldots, $m_k$ in $G$, then $f^{\lcm(m_1, \ldots, m_k)}$ lifts to $\widetilde{M}$.
\end{tcbcorollary}

\begin{proof}[of Corollary~\ref{cor:LiftingParabolicsViaCylinders}]
Let us assume that $\rho(\gamma_1)$, \ldots, $\rho(\gamma_k)$ have finite order $m_1$, \ldots, $m_k$.
Then, by Lemma~\ref{lem:LiftingCylindersInCoverings}, the preimage of the cylinder $C_i$ in
$\widetilde{M}$ is a disjoint union of cylinders with same height and whose circumference is
$m_i$ times the one of $C_i$. In particular, the modulus of any of these cylinders in $\widetilde{M}$
is $\frac{\mu(C_i)}{m_i}$. Let $f \in \Aff(M)$ be the multitwist whose derivative stabilizes the parabolic
direction. Let $m = \lcm(m_1, \ldots, m_k)$. Then $f^m$ lifts as a multitwist in $\widetilde{M}$.
\end{proof}

We now turn to the case of Abelian coverings. As we already saw in
Theorem~\ref{thm:AbelianLiftingCriterion} from
Section~\ref{ssec:AffineGroupsOfAbelianCoverings} the lifting
of elements of $\Aff(M, \Sigma)$ to the covering $\widetilde{M}$
can be expressed in terms of the action of this element on
the relative homology group $H_1(M,\Sigma; \Q)$.
\begin{tcbdefinition}{}{}
Let $M$ be a finite-type translation surface and $\theta\in\R/2\pi\Z$ be a completely
periodic direction. Let $\gamma_1$, \ldots, $\gamma_k$ be the core curves
of the maximal cylinders in the corresponding decomposition of $M$. The \emph{homological dimension} of the completely
periodic direction $\theta$ in $M$ is the dimension of the subspace
spanned by $\gamma_1$, \ldots, $\gamma_k$ in absolute homology $H_1(M; \Q)$.
\end{tcbdefinition}
Note that the homological dimension is always smaller than or equal to the
number of cylinders.
\begin{tcbtheorem}{one-cylinder trick}{OneCylinderTrick}
Let $M$ be a finite-type translation surface and $p:\widetilde{M} \to M$
be a $\Z^d$-translation covering defined by a subspace of $\kerhol(M,\Sigma)$ and branched over $\Sigma$. Let $f\in\Aff(M)$ be a multitwist
whose derivative stabilizes a completely periodic direction $\theta\in\R/2\pi\Z$.
\begin{enumerate}
\item If there is \emph{only one} cylinder $C$ in the cylinder decomposition of $M$
in direction $\theta$, then the preimages of
this cylinder in $\widetilde{M}$ are disjoint isometric copies of $C$.
\item The direction $\theta$ has homological dimension 1 in $M$ if and only if $f_*$ acts as the
identity on $\kerhol(M,\Sigma)$.
\end{enumerate}
In both situations above $f$ lifts to an element $\widetilde{f}\in\Aff(\widetilde{M})$.
\end{tcbtheorem}
\begin{tcbexample}{}{BandsInfiniteStaircase}
\label{exa:BandsInTheInfiniteStaircase}
First, let's consider the situation of the infinite staircase from
Section~\ref{ssec:InfiniteStaircase} to which Theorem~\ref{thm:OneCylinderTrick} applies.
Recall that the infinite staircase is a $\Z$-translation covering $p:\widetilde{M} \to M$,
where
the base $M$ is a torus and the covering map $p$ is branched
over two points (see Example~\ref{exa:Zcoverings}).
In the horizontal
and vertical directions, the cylinder decompositions consist of a single cylinder
and hence the first item of Theorem~\ref{thm:OneCylinderTrick} applies. In the
diagonal direction $\theta = \pi/4$ there are two cylinders but their core
curves are homologous and the second item of Theorem~\ref{thm:OneCylinderTrick}
applies. Note that in this diagonal direction, the staircase decomposes
into two strips. The parabolic element $\widetilde{f}$ that stabilizes
this strip decomposition is not a multitwist in the sense of Definition~\ref{def:Multitwist}.
\end{tcbexample}

We will actually prove results stronger than Theorem~\ref{thm:OneCylinderTrick}. In particular, we will describe explicitly the action of multitwists on homology.
The algebraic intersection form
\[
\langle\cdot,\cdot\rangle:H_1(M \setminus \Sigma;\R)\times H_1(M,\Sigma;\R)\to\R
\]
will play an important role.
\begin{tcblemma}{}{MultitwistLift}
Let $M$ be a finite-type translation surface and $V$ a subspace of $\kerhol(M, \Sigma)$ defining a $\Z^d$-covering
$p:\widetilde{M} \to M$ branched over $\Sigma \subset M$. Let $\theta\in\R/2\pi\Z$ be an affine multitwist direction in $M \setminus \Sigma$ with corresponding multitwist $f \in \Aff(M, \Sigma)$. Let $\gamma_1$, \ldots, $\gamma_k$ in $H_1(M \setminus \Sigma; \Q)$ be the core curves of the cylinders $C_1,\ldots,C_k$ in the and $\mu_1$, \ldots, $\mu_k$ their moduli.
Let us define
\[
\phi:
\begin{array}{ccc}
H_1(M,\Sigma; \Q) & \to & H_1(M,\Sigma; \Q) \\
v & \mapsto & \sum_{i=1}^k \mu_i\ \langle \gamma_i,v \rangle\ \overline{\gamma_i}
\end{array}
\]
where $\overline{\gamma_i}$ denotes the image of $\gamma_i$ in $H_1(M, \Sigma; \Q)$.
Then
\begin{enumerate}
\item $\langle  \gamma_i, V \rangle = 0$ for each $i=1,\ldots,k$ if and
only if the preimage of each $C_i$ in $\widetilde{M}$ is a disjoint union
of cylinders isometric to $C_i$.
\item $\phi(V) = 0$ if and only if $f_*$ acts trivially on $V$,
\item $\phi(V) \subset V$ if and only if $f$ lifts to the cover defined by $V$.
\end{enumerate}
\end{tcblemma}
The main ingredient for the proof of Lemma~\ref{lem:MultitwistLift} is the description of the action of a multitwist on homology given by the following exercise.

\begin{tcbexercise}{}{MultitwistHomologyAction}
Let $\theta\in\R/2\pi\Z$ be an affine multitwist direction in $M \setminus \Sigma$ with corresponding multitwist $f \in \Aff(M, \Sigma)$. Let $\gamma_1$, \ldots, $\gamma_k$ in $H_1(M \setminus \Sigma; \Q)$ be the core curves of the corresponding cylinders $C_1,\ldots,C_k$. The multitwist $f$ acts as a power of a Dehn twist in each
cylinder. Let $m_i$ be the multiplicity of the twist
in the cylinder $C_i$. Show that the action of $f$ on
$H_1(M, \Sigma; \Z)$ is given by:
\begin{equation}
  \label{EQ:FormulaMultitwistHomologyAction}
f_*(v) = v + \sum_{i=1}^k m_i \langle \gamma_i, v \rangle\ \overline{\gamma_i}
\end{equation}
\end{tcbexercise}

The last ingredient in the proof of Theorem~\ref{thm:OneCylinderTrick} is the explicit expression of $\Im(\omega)$ for
horizontal completely periodic directions.
\begin{tcblemma}{}{BaseHomologyAndIntersection}
Let $M$ be a finite-type translation surface whose horizontal direction is
completely periodic. Let $\Sigma \subset M$ be a finite set of points
and let $\gamma_1$, \ldots, $\gamma_k$ in $H_1(M \setminus \Sigma; \R)$
be the core curves of maximal cylinders in $M \setminus \Sigma$ and
$h_1$, $h_2$, \ldots, $h_k$ their heights. Then the cycle
$\gamma := h_1 \gamma_1 + h_2 \gamma_2 + \ldots + h_k \gamma_k$ is the
Poincar\'e dual of $\Im(\omega)$. That is to say
\[
\forall \eta \in H_1(M, \Sigma; \R), \quad \langle \gamma, \eta \rangle = \int_\eta \Im(\omega).
\]
\end{tcblemma}

We first proove Theorem~\ref{thm:OneCylinderTrick} as a consequence
of Lemma~\ref{lem:MultitwistLift} and Lemma~\ref{lem:BaseHomologyAndIntersection}.
\begin{proof}[Proof of Theorem~\ref{thm:OneCylinderTrick}]
Let us assume that the surface $M$ decomposes as a single cylinder $C$ with core curve $\gamma$ in direction $\theta$. By Lemma~\ref{lem:BaseHomologyAndIntersection} the core curve
is the Poincar\'e dual of an element $\eta$ in the tautological plane
$\R \Re(\omega) \oplus \R \Im(\omega)$, that is $\eta=a\Re(\omega)+b\Im(\omega)$, for some $a,b\in\R$. Hence for any $\delta\in\kerhol(M,\Sigma)$ we have that
$\langle \gamma,\delta \rangle = a\int_\delta \Re(\omega) + b\int_\delta \Im(\omega) = 0$. By the first item in Lemma~\ref{lem:MultitwistLift} we conclude that the preimages of $C$ in $\widetilde{M}$ are disjoint isometric copies of $C$.

Now let us prove the second item. We assume that the horizontal direction
has homological dimension one. That is to say for each $i=1,\ldots,k$
there exists $\alpha_i \in \Q_{>0}$ so that
$\overline{\gamma_i} = \alpha_i \overline{\gamma_1}$.
Since the core curves are primitive elements of $H_1(M; \Z)$ we have
$\alpha_i = 1$ for all $i=1,\ldots,k$.  Now recall that the modulus of the cylinder $C_i$ is $\mu_i = h_i / w_i$ where
$h_i$ and $w_i$ are the heights and circumferences of $C_i$ respectively. But
since the $\gamma_i$ are homologous, the $w_i$ are equal and $\mu_i = h_i / w_1$.
In particular, we have for any $v \in H_1(M, \Sigma; \Q)$
\[
\phi(v)=\sum_{i=1}^k \mu_i \langle \gamma_i, v \rangle \overline{\gamma_i} =
\frac{1}{w_1} \left\langle \sum_{i=1}^k  h_i \gamma_i, v \right\rangle \overline{\gamma_1}.
\]
Now, similarly to the one-cylinder case we just apply
Lemma~\ref{lem:BaseHomologyAndIntersection} to conclude that
$\left\langle \sum_{i=1}^k  h_i \gamma_i, v \right\rangle$ vanishes
for vectors $v \in \kerhol(M, \Sigma)$. The conclusion follows from the second item in Lemma~\ref{lem:MultitwistLift} and the Abelian lifting criterion in Theorem~\ref{thm:AbelianLiftingCriterion}.
\end{proof}
\begin{proof}[Proof of Lemma~\ref{lem:BaseHomologyAndIntersection}]
Let us consider the following generating system of $H_1(M, \Sigma; \R)$. Let $\eta_1, \eta_2, \ldots, \eta_n$ be the saddle connections on the boundary
of the cylinders appearing in the horizontal direction. We orient them in such way that $\int_{\eta_i}\omega>0$, for all $i=1,\ldots,n$.
Next, for each cylinder choose a saddle connection that goes from the top
to the botton. We denote $\xi_i$ the transverse saddle connection in the $i$-th cylinder.
The set $\{\eta_1,\ldots,\eta_n,\xi_1,\ldots,\xi_k\}$ generates $H_1(M, \Sigma; \R)$ and
\begin{itemize}
\item for all $i=1,\ldots,n$ and $j=1,\ldots,k$ we have $\langle \eta_i, \gamma_j \rangle = 0$
\item for all $i=1,\ldots,k$ and $j=1,\ldots,k$ we have $\langle \xi_i, \gamma_j \rangle = \delta_{ij}$
where $\delta_{ij}$ is the Kronecker symbol.
\end{itemize}
As $\int_{\eta_i} \Im(\omega)=0$ and $\int_{\xi_i} \Im(\omega)=h_i$ we obtain the result.
\end{proof}
\begin{proof}[Proof of Lemma~\ref{lem:MultitwistLift}]
Suppose that the $\Z^d$-covering $p:\widetilde{M}\to M$ is defined by the kernel of $\rho:\pi_1(M\setminus\Sigma)\to G\simeq\Z^d$. Then $\langle \gamma_i,V\rangle=0$ if and only if $\rho(\gamma_i)=0$ for every $i=1,\ldots,k$. From Lemma~\ref{lem:LiftingCylindersInCoverings} we conclude that this last condition happens if and only if for every $i=1,\ldots,k$ the preimage of each $C_i$ in $\widetilde{M}$ is a disjoint union of cylinders isometric to $C_i$.

The multitwist $f$ acts as a power of a Dehn twist in each
cylinder.

Let $m_i$ be the multiplicity of the twist
in the cylinder $C_i$.
Now, recall that this multiplicity $m_i$ is proportional to the moduli $\mu_i$
of the cylinders. Namely if the parabolic stabilizes the horizontal direction and is hence written as
$\displaystyle
\begin{psmallmatrix}
1 & \lambda \\
0 & 1
\end{psmallmatrix}$
then $\lambda \mu_i = m_i$. By formula (\ref{EQ:FormulaMultitwistHomologyAction}) from Exercise~\ref{exo:MultitwistHomologyAction} we have that $f_*(v)= v + \frac{1}{\lambda}\phi(v)$ where
$\phi$ was defined in the statement of Lemma~\ref{lem:MultitwistLift}.
Hence $f_*$ acts trivially on $V$ if and only
if $\phi(V) = 0$. From here the third item follows trivially and the last item is deduced
from the Abelian lifting criterion (Theorem~\ref{thm:AbelianLiftingCriterion}).
\end{proof}
Recall from Definition~\ref{def:PeriodicPoint} and the discussion after that
a periodic point in a square tiled surface $p: M \to \T^2$ is a point that
projects to $\Q^2/\Z^2$ under $p$.
\begin{tcbcorollary}{one-cylinder trick for square-tiled surfaces}{OneCylinderTrickSquareTiled}
Let $M$ be a square-tiled surface and let $\Sigma \subset M$ be a finite set of periodic points. Assume that $M$ admits a completely periodic direction
of homological dimension one. Let $\widetilde{M} \to M$ be a $\Z^d$-covering
determined by a subspace $V$ of $\kerhol(M,\Sigma)$. Then the relative Veech group\footnote{See Definition~\ref{def:RelativeVeechGroup}}
$\Gamma(\widetilde{M}\to M)$ contains a non-trivial normal subgroup of $\Gamma(M)$
generated by parabolic elements.
\end{tcbcorollary}

A non-trivial normal subgroup of the lattice $\Gamma(M)$ is
 "big" in the following sense: its limit set is the whole circle.
This is a direct consequence of Lemma~\ref{lem:NormalSubgroupLimitSet}, see also Appendix~\ref{app:FuchsianGroups} for details. A particular case of
Corollary~\ref{cor:OneCylinderTrickSquareTiled} appeared in~\cite{HubertSchmithuesen10}.

\begin{proof}
Let $G$ be the subgroup of $\Aff(M)$ generated by the matrices $Df$ where $f$ is a multitwist that stabilizes a completely periodic direction with homological dimension one in $M$. Then by assumption $G$ is non-trivial. It also follows from the one-cylinder trick (Theorem~\ref{thm:OneCylinderTrick}) that any element in $G$ lifts to $\widetilde{M}$. Let us prove that $G$ is a normal subgroup. If $f$ is a multitwist in $G$ stabilizing a direction $\theta$ and $h \in \Aff(M)$, then $h f h^{-1}$ is the multitwist stabilizing the direction $h \theta$. Moreover, the direction $\theta$ and $h \theta$ are affinely equivalent. In particular the completely periodic direction $\theta$ and $h \theta$ in $M$ have the same homological dimensions and $h f h^{-1}$ lifts to $\widetilde{M}$.
\end{proof}

\begin{tcbexample}{}{EierlegendeWollmilchsau}
The eierlegende Wollmilchsau\footnote{The literal translation for this German word could be \emph{egg-laying woolly milk sow}. This of course does not capture its meaning, and for this reason we added the following explanation, which is taken from die Deutsche Welle (a German public state-owned international broadcaster funded by the German federal tax budget, for more details see \href{https://www.dw.com/en/eierlegende-wollmilchsau/a-6616972}{here}): Germans have, apparently, never subscribed to the premise that you can't make all the people happy all the time. That's the idea behind the word \emph{Eierlegende Wollmilchsau}. Who wouldn't be happy to have a pig that lays eggs, gives milk and produces wool before it's turned into a tasty Sunday roast? These days, the term has come to mean any person or device that's able to be used in a wide variety of ways -- like an alarm clock that makes coffee and toast or a cell phone that reminds you to pickup the dry cleaning.} $p:W\to \T^2$ is the square-tiled surface depicted in figure~\ref{fig:Eierlegendewollmilchsau}. It was discovered independently by Herrlich and Weitze-Schmith\"{u}sen~\cite{HerrlichSchmithusen08} and Forni~\cite{Forni06}. $W$ is a genus 3 translation surface with $4$ conical singularities of total angle $4\pi$ (\ie it belongs to the stratum $\cH(1,1,1,1)$).
\begin{figure}[H]
\begin{center}\includegraphics{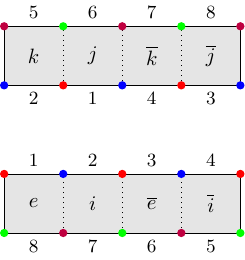}\end{center}
\caption{The eierlegende Wollmichsau square-tiled surface $W$. The integers on the outside of the polygons
correspond to the edge pairing. The bold labels in the squares show that $W$ has Deck group
transformation the quaternion group.}
\label{fig:Eierlegendewollmilchsau}
\end{figure}
\begin{tcbexercise}{}{Wollmilchsau}
Let $W$ be the eierlegende Wollmilchsau square-tiled surface and $\Sigma$ be its set of four conical singularities.
\begin{enumerate}
\item
Show that the Deck group of $p:W\to \T^2$ can be identified with the quaternion group
\[
Q := \{e, i, j, k, \overline{e}, \overline{i}, \overline{j}, \overline{k}\}
\]
with neutral element $e$ and multiplication rules
\begin{itemize}
\item $i^2 = j^2 = k^2 = \overline{e}$
\item $ij = k$, $ki = j$ and $jk = i$
\item $\forall g \in Q$, $g^{-1} = \overline{g}$ where $\overline{\overline{g}} = g$.
\end{itemize}
\item Show that the subgroup of $\pi_1(\T^2 \setminus \{0\})$ that defines $W$
is characteristic. Conclude that the Veech group of $W$ is $\SL(2,\Z)$
(see Corollary~\ref{cor:LiftingCharacteristicCoverings}).
\item Use the one-cylinder trick (Theorem~\ref{thm:OneCylinderTrick}) to prove that the horizontal multitwist
in $W$ acts trivially on $\kerhol(W, \Sigma)$.
\item
Using the fact that the Veech group of $W$ is $\SL(2,\Z)$, prove that any multitwist in
$\Aff(W, \Sigma)$ acts trivially on $\kerhol(W, \Sigma)$. \emph{Hint}: there is only one conjugacy class of parabolic elements in $\Gamma(W) = \SL(2,\Z)$.
\item
Show that the derivative map induces an isomorphism between the subgroup of $\Aff(W, \Sigma)$ generated
by the multitwists of the previous question and $\Gamma(4)$
(\emph{hint:} for surjectivity, use the fact that $\Gamma(4)$ is generated by its parabolic elements)
\item
Conclude that the Veech group of any $\Z^d$-covering $\widetilde{W}$ of $W$ at most
branched over $\Sigma$ contains the congruence group $\Gamma(4)$
(\emph{hint:} use Theorem~\ref{thm:AbelianLiftingCriterion}).
\end{enumerate}
\end{tcbexercise}
In the above exercise we explicited the action of a finite index subgroup of $\Aff(M,\Sigma)$ on $\kerhol(W, \Sigma)$.
The complete description of the action of $\Aff(W,\Sigma)$ on $\kerhol(W, \Sigma)$ is given
in~\cite{MatheusYoccoz10}. Other properties of the eierlegende Wollmilchsau are
discussed in~\cite{Forni06} and \cite{HerrlichSchmithusen08}.

Let us also mention that the eierlegende Wollmilchsau $W$ is not the only known example of
finite-type square-tiled surface such that a finite index subgroup of its affine group acts trivially
in homology. There is a cousin of the eierlegende Wollmilchsau: a genus 4 square-tiled surface called
the \emph{ornithorynque} (the french for platypus). Similarly to what is proposed
in Exercise~\ref{exo:Wollmilchsau}, one can prove that a subgroup of its affine
group isomorphic to $\Gamma(3)$ acts trivially on homology. We refer to~\cite{MatheusYoccoz10}
and~\cite{ForniMatheusZorich10}, \cite[Chapter 7.2]{ForniMatheus14} for more details on this example.
\end{tcbexample}

\subsection{The Hooper-Weiss theorem}
\label{ssec:HooperWeissTheorem}
Let us recall from Definition~\ref{def:PeriodicPoint} and the remark that follows
that a periodic point in a square-tiled surface $p:M \to \T^2$ is a point that projects
under $p$ to a point in $\Q^2/\Z^2$. The aim of this section is to prove the following result.
\begin{tcbtheorem}{\cite{HooperWeiss09}}{HooperWeissFirstKind}
Let $M$ be a square-tiled surface of genus $2$ and let $\Sigma \subset M$ be a finite set of periodic points. Let $\widetilde{M}\to M$ be a $\Z^d$-covering determined by a subspace $V$ of $\kerhol(M, \Sigma)$. Then the relative Veech group\footnote{See Definition~\ref{def:RelativeVeechGroup}} $\Gamma(\widetilde{M} \to M)$ contains a non-trivial normal subgroup of $\Gamma(M, \Sigma)$\footnote{See Definition~\ref{def:VeechGroup}}.
\end{tcbtheorem}

The reader should notice the similarity of Theorem~\cite{HooperWeiss09} with Corollary~\ref{cor:OneCylinderTrickSquareTiled} of the one-cylinder trick. In both situations, the morphism $\Aff(M) \to \GL(\kerhol(M,\Sigma))$ has an infinite kernel and all elements in this kernel lift to $\Aff(\widetilde{M})$. In particular, Theorem~\cite{HooperWeiss09} implies that $\Gamma(\widetilde{M} \to M)$ is a Fuchsian group of the first kind (see Appendix~\ref{app:FuchsianGroups} for more details on Fuchsian groups).

\begin{proof}[Proof of Theorem~\ref{thm:HooperWeissFirstKind}]
Since $\Sigma$ are periodic points, the relative Veech group $\Gamma(M,\Sigma)$
is a lattice of finite index in $\Gamma(M)$. Let
$\psi: \Aff(M,\Sigma) \to \GL(\kerhol(M,\Sigma))$ and
$\psi_0: \Aff(M,\Sigma) \to GL(\kerhol(M))$ be the morphisms induced
by the action of the affine group on respectively $H_1(M, \Sigma; \Z)$ and $H_1(M; \Z)$.
By Theorem~\ref{thm:KerPsiVSkerPsi0}, it suffices
to show that the image of $\ker(\psi_0)$ in $\Gamma(M)$ is non-trivial. We
suppose that $D\ker(\psi_0)$ is trivial and derive a contradiction.

Because $M$ has genus 2, $\kerhol(M)$ has dimension 2. We know from
Theorem~\ref{thm:CohomologyNonElementaryVeechGroup} that $\kerhol(M)$ is
symplectic. Up to passing to a finite index subgroup $A$ of $\Aff(M,\Sigma)$ we
can further suppose that:
\begin{enumerate}
\item the image $\psi_0(A)$ is torsion-free
(\eg consider $\psi_0^{-1}(\Gamma(2))$ where $\Gamma(2)$ denote the principal
congruence subgroup)
\item that both $D$ and $\psi_0$ are faithful representations of
$A$ into $\SL(2,\Z)$.
(because in this case $\ker(\Aff(M) \to \Gamma(M))$ is finite since it is a subgroup of the group of conformal automorphisms of a compact Riemann surface).
\end{enumerate}
From now on we fix such finite index subgroup $A \subset \Aff(M,\Sigma)$.

Now, let $f \in A$. Then
\begin{enumerate}
\item if $D(f)$ is parabolic then $\psi_0(f)$ is parabolic. Indeed, by
Exercise~\ref{exo:MultitwistHomologyAction}, the action of $f$ on homology has
only eigenvalues which are roots of unity. Furthemore, $\psi_0(f)$ can not be
identity or elliptic since by assumption $\psi|_{A}$ is injective.
\item If $D(f)$ is hyperbolic, then $2 \leq |\tr(\psi_0(f))| < |\tr(D(f))|$
(the dominant eigenvalue of $D(f)$ is the stretch-factor of the pseudo-Anosov
and dominates any other eigenvalue of the action of $f$ on $H_1(M;\R)$,
see Theorem~\ref{thm:TranslationPseudoAnosovHomologyAction}).
\end{enumerate}
Now consider the hyperbolic surfaces
$S=\H^2/D(A)$ and $S_0 = \H^2/\psi_0(A)$. Because both $\psi_0$
and $D$ are faithful, there exists an homotopy equivalence $\phi: S \to S_0$
induced by the isomorphism at the level of fundamental groups. Because
parabolic elements in $A$ are mapped to parabolic elements by
$\psi_0$ we know that the number of cusps $n$ of $S$ is smaller than or equal
to the one $n_0$ of $S_0$. We claim that there are actually equal.
Indeed, let us denote by $\gamma_1$, \ldots, $\gamma_n$ in $H_1(S; \Z)$
the generators of peripheral curves in $S$. If we had $n_0 > n$ then these linearly
dependent elements would be mapped to linearly independent elements
of $H_1(S_0; \Z)$ contradicting the homotopy equivalence. Hence, $S$ and $S_0$
have the same number of cusps and, because of the homotopy equivalence, the same
genus. Thus $\psi_0(A)$ is a lattice.

Now for each non-peripherial $\gamma\in\pi_1(S)$ let $\ell(\gamma)$
denote the length of the geodesic representative on $S$, and let $\ell_0(\gamma)$
denote the length of the geodesic representative of $\phi_*(\gamma)$. It is a deep
result by Thurston (see~\cite[Theorem 3.1]{Thurston86}) that
\begin{equation}
\label{eq:ThurstonExpansionFact}
\sup_{\gamma\in\pi_1(S)}\frac{\ell_0(\phi_*(\gamma))}{\ell(\gamma)}\geq 1,
\end{equation}
with equality only if $\phi$ is homotopic to an isometry of the hyperbolic
surfaces $S$ and $S_0$.

Now a closed geodesic $\gamma$ in $S$ corresponds to the
conjugacy class $D(f)$ of a hyperbolic affine element $f$ in $\Aff(M)$.
Its image $\gamma_0 = \phi(\gamma)$ corresponds to the conjugacy class
of $\psi_0(f)$. The hyperbolic lengths of the loops $\gamma$ and
$\gamma_0$ are respectively the logarithms of the dominant
eigenvalues of $D(f)$ and $\psi_0(f)$.
Thus, equality in~(\ref{eq:ThurstonExpansionFact}) implies that $\phi$ is isotopic to an isometry, thus the logarithms of the dominant
eigenvalues of $D(f)$ and $\psi_0(f)$ are equal and this contradicts the
fact that $|\tr(\psi_0(f))| < |\tr(D(f))|$.
\end{proof}

\begin{tcbremark}{}{GeneralHooperWeissFirstKind}
The proof of Theorem~\ref{thm:HooperWeissFirstKind} is actually more general. Along the same line one can prove the following generalization.
\begin{tcbtheorem}{}{HooperWeissFirstKindII}
Let $M$ be a Veech surface. Assume that there exists a rank 2 subspace
$E \subset \kerhol(M)$ defined over $\Q$ and preserved by a subgroup $A$ of
finite index of $\Aff(M)$. Then the kernel of the induced action $A \to \GL(E)$
is infinite or, equivalently, the image of this kernel in $\Gamma(M)$ is non-trivial.
\end{tcbtheorem}
We will even see a more general version in Lemma~\ref{lem:WindTreeAffineGroupKernel} and
Remark~\ref{rk:GeneralHooperWeissFirstKind2}.
\end{tcbremark}

\begin{tcbexercise}{}{}
Let $M = (X,q)$ be a meromorphic quadratic differential with at most simple poles on a torus
such that its orientation cover is a
Veech surface. Let $\widetilde{M} \to (X,q)$ be the $\Z^2$-half-translation covering
obtained by taking the universal cover of the torus (\ie a Panov plane as in
Section~\ref{ssec:PanovPlanes}).Show that Theorem~\ref{thm:HooperWeissFirstKindII} can be applied to
$\Aff(M) \to \GL(H_1(M, \Z))$. \textit{Hint}: use Exercise~\ref{exo:HalfTranslationKerhol}.
\end{tcbexercise}

\begin{tcbexample}{}{}
We now come back to the square-tiled surface $Z_{3,1}$ of Example~\ref{exa:OrigamiH11-1}.
This example also appears as~\cite[Example 7.2]{HooperWeiss09}.
As can be seen on Figure~\ref{fig:SL2ZOrbit-H11} any cylinder decomposition of
the square-tiled surface $Z_{3,1}$ has homological dimension 2. As a consequence, the
kernel of the action $\Aff(Z_{3,1}) \to \GL(\kerhol(Z_{3,1})$
\emph{does not contain any parabolic element} (see Theorem~\ref{thm:OneCylinderTrick}).
However, the Hooper-Weiss theorem does apply: this kernel is infinite.

The square-tiled surface $Z_{3,1}$ belongs to the stratum $\cH(1,1)$. Let us mention that
in the other stratum of genus $2$, $\cH(2)$ any square-tiled surface admits a completely
periodic direction made of a single maximal cylinder~\cite{HubertLelievre06}.
\end{tcbexample}

The following is a delicate open question.
\begin{tcbquestion}{}{}
Does there exist a square-tiled surface $p:M \to \T^2$ of genus 3 such that the kernel of
the induced action $\Aff(M) \to \GL(\kerhol(M))$ is finite?
\end{tcbquestion}
Note that Theorem~\ref{thm:OneCylinderTrick} discards any square-tiled surface that admits a completely periodic of homological dimension one. The reader interested by this question should consult~\cite{HubertMatheus}.


\subsection{Affine symmetries of wind-tree models}
\label{ssec:SymmetriesWindTree}
\label{ssec:CylinderDecompositionWindtree_1/2_1/2}

\begin{figure}[ht!]
\begin{center}%
\includegraphics{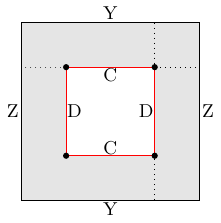}
\hspace{1cm}
\includegraphics{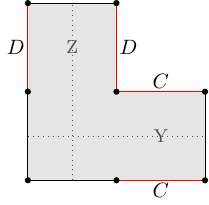}
\end{center}
\caption{Two views of $L_{1/2,1/2}$. The left picture is inspired from
Figure~\ref{fig:WindtreeQuotientSurface} for $X_{a,b}$. The right picture
is obtained from the left one by cutting along the dotted lines and gluing
along the sides $Y$ and $Z$.}
\label{fig:L1/2}
\end{figure}

In this section we study the affine symmetries of the wind-tree models $W_{a,b}$, which were defined in Section~\ref{ssec:Intro:WindTree}. Recall that $W_{a,b}$
is a $\Z^2$-covering of a finite-type translation surface $X_{a,b}$ that belongs to the
stratum $\cH(2^4)$. Moreover, $X_{a,b}$ is a $K$-covering of $L_{a,b}$ in $\cH(2)$ where
$K = \Z/2\Z \times \Z/2\Z$ is the Klein group (see Figure~\ref{fig:L1/2}).

By important work of K. Calta~\cite{Calta04} and
C. McMullen~\cite{McMullen03}, we know the list of $(a,b)$ for which $L_{a,b}$
are Veech surfaces.
\begin{tcbtheorem}{\cite{Calta04}, \cite{McMullen03}}{Calta:McMullen:Veech:Surfaces:H2}
Let $a,b \in (0,1)$, then the following are equivalent:
\begin{enumerate}
\item \label{item:thm:Calta:McMullen:Veech}
The surface $L_{a,b}$ is a Veech surface.
\item \label{item:thm:Calta:McMullen:ThurstonVeech}
The surface $L_{a,b}$ comes from a Thurston-Veech construction. In particular, it admits
horizontal and vertical affine multitwists.
\item \label{item:thm:Calta:McMullen:PseudoAnosov}
The surface $L_{a,b}$ admits an affine pseudo-Anosov homeomorphism.
\item \label{item:thm:Calta:McMullen:Explicit:Parameters}
The parameters $a$ and $b$ are rational, or there exists a square-free positive integer $D$
and rational numbers $x,y$ such that $1/(1-a) = x + y \sqrt{D}$ and $1/(1-b) = (1-x) + y \sqrt{D}$.
\end{enumerate}
\end{tcbtheorem}

Some of the implications in Theorem~\ref{thm:Calta:McMullen:Veech:Surfaces:H2} are straightforward:
\ref{item:thm:Calta:McMullen:Veech} implies \ref{item:thm:Calta:McMullen:PseudoAnosov},
\ref{item:thm:Calta:McMullen:Explicit:Parameters} is equivalent to \ref{item:thm:Calta:McMullen:ThurstonVeech}
and \ref{item:thm:Calta:McMullen:ThurstonVeech} implies \ref{item:thm:Calta:McMullen:PseudoAnosov}. The
reader non-familiar with this topic might want to prove them. All the other implications are non trivial
and we refer to the original articles.

The following subset of parameters plays an important role
\begin{equation}
\label{eq:DefinitionEParametersWindtree}
\cE := \left\{
\left(\frac{p}{q},\frac{r}{s}\right) \in \Q^2:\ 0<p<q,\ 0<r<s,\ \text{$p$, $q$ odd and $r$, $s$ even}
\right\}.
\end{equation}
\begin{tcbtheorem}{}{WindTreeVeechGroup}
For the wind-tree model parameters $(a,b)$, such that $L_{a,b}$ is a Veech surface, the following holds:
\begin{enumerate}
\item $\Gamma(W_{a,b} \to L_{a,b})$ has\footnote{See Definition~\ref{def:RelativeVeechGroup}.} infinite index in $\Gamma(L_{a,b})$,
\item $\Gamma(W_{a,b} \to L_{a,b})$ contains a non-trivial normal subgroup,
in particular its limit set is the whole circle, and
\item if furthermore $(a,b) \in \cE$, where $\cE$ is defined
in~\eqref{eq:DefinitionEParametersWindtree}, then $\Gamma(W_{a,b} \to L_{a,b})$
contains a parabolic element.
\end{enumerate}
\end{tcbtheorem}
The rest of the section is dedicated to the proof of Theorem~\ref{thm:WindTreeVeechGroup}.
We will prove the first item about the index of
$\Gamma(W_{a,b} \to L_{a,b})$ in $\Gamma(L_{a,b})$. Then, we consider
parameters $(a,b)$ in $\cE$ and use the one-cylinder trick
(Theorem~\ref{thm:OneCylinderTrick}) to prove the third item.
Finally, we apply a variant of the Hooper-Weiss theorem
(Theorem~\ref{thm:HooperWeissFirstKind}) for all
parameters $(a,b)$ which provides a proof of the second item.

\subsubsection{$\Gamma(W_{a,b} \to L_{a,b})$ has infinite index in $\Gamma(L_{a,b})$}
We begin the proof of Theorem~\ref{thm:WindTreeVeechGroup} by proving that
$\Gamma(W_{a,b} \to L_{a,b})$ is not a lattice. We use an argument
that appears in~\cite{Cabrol}.

Let $a,b$ be parameters so that $L_{a,b}$ is a Veech surface. Let
$f \in \Aff(X_{a,b})$ be the multitwist stabilizing the horizontal direction.
Since $X_{a,b}\to L_{a,b}$ is a finite covering (the deck transformation group is the Klein group), it is sufficient to show that no power of $f$ lifts to $W_{a,b}$.
Let $V$ be the subspace of $\kerhol(X_{a,b})$ generated by the elements
$Y$ and $Z$ as in~\eqref{eq:WindTreeCycleCovering}. We want to show that no power of
$f_*$ preserves $V$. By a direct computation, one can see that $f_*$ fixes
$Y$ but for any $n \geq 1$, we have $(f_*)^n Z \not\in V$. As a consequence
no power of $f$ lifts to $W_{a,b}$ and the index of $\Gamma(W_{a,b} \to L_{a,b})$
in $\Gamma(L_{a,b})$ is infinite.

\subsubsection{Parabolic elements for parameters in $\cE$}
Let us recall that in~\eqref{eq:DefinitionEParametersWindtree} we defined
a subset of rational parameters $\cE$ for the windtree-model. For
$(a,b)$ in $\cE$ we consider parabolic elements in $L_{a,b}$ and show that some of
them can be lifted to $W_{a,b}$.

\begin{tcblemma}{}{WindTreeOneCylinderTrick}
Let $(a,b) \in \cE$. Then $L_{a,b}$ admits a completely periodic direction made of
a single cylinder. For any such direction, if the surface $L_{a,b}$ decomposes
into a cylinder of modulus $\mu$ then the surface $X_{a,b}$ decomposes into two
cylinders of modulus $\mu/2$ whose core curves are homologous.
\end{tcblemma}
As a consequence of the above Lemma and Theorem~\ref{thm:OneCylinderTrick}
we have.
\begin{tcbcorollary}{}{WindTreeOneCylinderTrick}
Let $(a,b) \in \cE$. The square of any multitwist stabilizing a one
cylinder direction in $L_{a,b}$ lifts to $X_{a,b}$ and $W_{a,b}$.
In particular, $\Gamma(W_{a,b} \to X_{a,b})$ contains a parabolic
element.
\end{tcbcorollary}

\begin{proof}[Proof of Lemma~\ref{lem:WindTreeOneCylinderTrick}]
We consider the quotient $\overline{L_{a,b}}$ of $L_{a,b}$ by the hyperelliptic
involution (the unique affine element in $\Aff(L_{a,b})$ whose derivative is minus identity). The surface
$\overline{L_{a,b}}$ is depicted in Figure~\ref{fig:LQuotient}. The advantage of
this surface is that it has genus 0, and describing a covering of a genus 0 surface
is very convenient (the fundamental group is generated by loops around the punctures).
\begin{figure}[!ht]
\begin{minipage}{0.4\textwidth}
\begin{center}%
\includegraphics{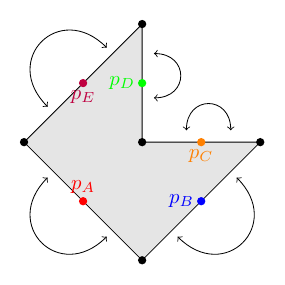}%
\end{center}
\subcaption{The quotient $\overline{L_{a,b}}$ of $L_{a,b}$ is a quadratic differential
on the sphere (stratum $\cQ(1,-1^5)$). Each side of the pentagon is glued to itself
via a 180° rotation and contains a pole in its center. The five poles are labeled $p_A$,
$p_B$, $p_C$, $p_D$, $p_E$ in counterclockwise order.}
\end{minipage}
\hspace{0.1\textwidth}%
\begin{minipage}{0.4\textwidth}
\begin{center}%
\includegraphics[scale=0.85]{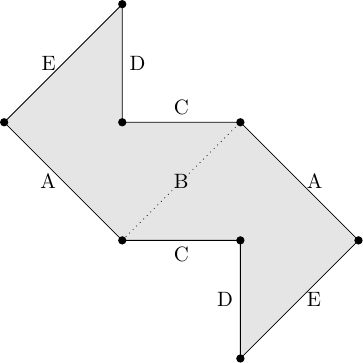}%
\end{center}
\subcaption{The surface $L_{a,b}$ as two copies of a pentagon. The labelling of the
sides is compatible with Figure~\ref{fig:L1/2}: the preimage of the pole $p_A$ belongs
to side $A$, etc.}
\end{minipage}
\caption{$\overline{L_{a,b}}$ and $L_{a,b}$.}
\label{fig:LQuotient}
\end{figure}
The first step, is to describe the coverings $L_{a,b} \to \overline{L_{a,b}}$
and $X_{a,b} \to \overline{L_{a,b}}$. Both involve homology with
coefficients in $\Z/2\Z$ relative to the set of singularities
$\overline{\Sigma}$ of $\overline{L_{a,b}}$ that consists of the five poles $p_A$, \ldots, $p_E$ and
a conical singularity of angle $3\pi$.
We denote by
$\overline{A}$, $\overline{B}$, $\overline{C}$, $\overline{D}$ and
$\overline{E}$ the elements in
$H_1(\overline{L_{a,b}}, \overline{\Sigma}; \Z/2\Z)$ that
corresponds to the quotients of the sides $A$, $B$, $C$, $D$ and $E$
of the figure $L_{a,b}$ (see Figure~\ref{fig:LQuotient}). Since
we consider $\Z / 2\Z$ coefficients, the orientation of these
curves does not matter. These 5 elements form a basis of
$H_1(\overline{L_{a,b}}, \overline{\Sigma}; \Z/2\Z)$.

In this basis the cycle $c_L \in H_1(\overline{L_{a,b}}, \overline{\Sigma}; \Z/2\Z)$ that corresponds to the degree 2
covering $L_{a,b} \to \overline{L_{a,b}}$ and $c_X \in H_1(\overline{L_{a,b}}, \overline{\Sigma}; (\Z/2\Z)^3)$ that
corresponds to the degree 8 covering $X_{a,b} \to \overline{L_{a,b}}$ are
\begin{equation}
\label{eq:XLcocycle}
c_L = \overline{A} + \overline{B} + \overline{C} + \overline{D} + \overline{E}
\qquad
c_X = (1,1,0) \overline{D} + (1,0,1) \overline{C} + (1,0,0) (\overline{A} + \overline{B} + \overline{E}).
\end{equation}
The formula for $c_L$ and $c_X$ can be deduced by analyzing the action of the Deck group of
$X_{a,b} \to L_{a,b}$ on $H_1(X_{a,b}; \Z)$.

\begin{figure}[!ht]
\begin{minipage}{0.3\textwidth}
\begin{center}\includegraphics{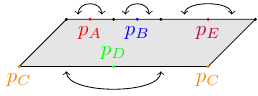}\end{center}
\subcaption{A one cylinder direction in $\overline{L_{a,b}}$ when $(a,b) \in \cE$
has one of its sides with the poles $\{p_C,p_D\}$ while the other three
are on the other side.}
\label{fig:Lbaronecyl}
\end{minipage}
\hspace{0.1\textwidth}
\begin{minipage}{0.6\textwidth}
\begin{center}\includegraphics{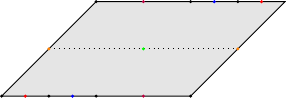}\end{center}
\subcaption{The preimage of the one cylinder direction of Figure~\ref{fig:Lbaronecyl}
in $L_{a,b}$.}
\label{fig:Lonecyl}
\end{minipage}
\begin{minipage}{\textwidth}
\begin{center}\includegraphics{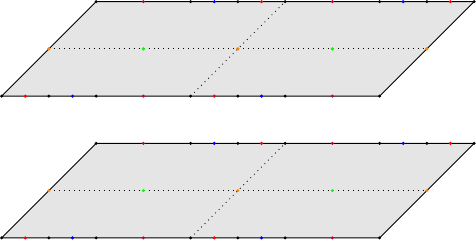}\end{center}
\subcaption{The preimage of the one cylinder direction of Figure~\ref{fig:Lbaronecyl}
in $X_{a,b}$.}
\label{fig:Xonecyl}
\end{minipage}
\caption{A one cylinder direction in $\overline{L_{a,b}}$ and its preimages in
$L_{a,b}$ and $X_{a,b}$.}
\label{fig:LXonecyl}
\end{figure}
A direction is completely periodic in $L_{a,b}$ if and only if it is completely periodic
in $\overline{L_{a,b}}$. Furthermore, if the direction in $L_{a,b}$ is made of a single
cylinder then it is also the case in $\overline{L_{a,b}}$ which is a quotient.

A single cylinder direction in $\overline{L_{a,b}}$ has two poles on one boundary
and the three other together and the $3\pi$ singularity on the other.
By~\cite{HubertLelievreTroubetzkoy11} for any completely periodic direction in
$\overline{L_{a,b}}$ made of a single cylinder, one side of this cylinder
contains the two poles $p_C$ and $p_D$ while the three other poles and the
conical singularity of angle $3\pi$ are on the other side
This fact holds exactly
for the parameters $(a,b)$ in $\cE$ and this is where this hypothesis is crucial.

From now on, we analyze a single cylinder direction in $\overline{L_{a,b}}$
where one boundary is made of a saddle connection between $p_C$ and $p_D$.
Such one cylinder direction has the form shown in
Figure~\ref{fig:Lbaronecyl} (up to permutation of the poles $p_A$, $p_B$, $p_E$
located on the top of the cylinder). The core curve
$\gamma$ seen as an element of
$H_1(\overline{L_{a,b}} \setminus \{p_A, p_B, p_C, p_D, p_E, q\}; \Z/2\Z)$ intersects the cycles $\overline{C}$ and $\overline{D}$ but not the other three.
From formulas~\eqref{eq:XLcocycle} we compute
\[
\langle c_X, \gamma \rangle = 1 + 1 = 0
\quad \text{and} \quad
\langle c_L, \gamma \rangle = (1,1,0) + (1,0,1) = (0,1,1).
\]
Hence, the preimage in $L_{a,b}$ and $X_{a,b}$ of this cylinder
consists respectively of a single cylinder of modulus $2 \overline{\mu}$
(see Figure~\ref{fig:Lonecyl}) and two cylinders of modulus $\overline{\mu}$
with homologous core curves (see Figure~\ref{fig:Xonecyl}). Hence the
conclusion of Lemma~\ref{lem:WindTreeOneCylinderTrick} holds for any
one-cylinder direction in $\overline{L_{a,b}}$ one of its boundary
consisting of a saddle connection between $p_C$ and $p_D$.

As already mentioned, for $(a,b) \in \cE$ all one-cylinder directions
are of this form. This concludes the proof of the lemma.
\end{proof}

\subsubsection{Homology action of $\Aff(X_{a,b})$}
We now prove the second item in Theorem~\ref{thm:WindTreeVeechGroup} by using a
variant of the Hooper-Weiss Theorem~\ref{thm:HooperWeissFirstKind}.

First of all the surface $X_{a,b}$ admits a decomposition of $H_1(X_{a,b}; \Q)$
which is finer than the discussed for a general surface in
Theorem~\ref{thm:CohomologyNonElementaryVeechGroup}. This is a very general phenomenon
for normal covers: the deck group of the cover acts on homology, and the affine group
preserves the isotypical components of this action~\footnote{The isotypical components
of a finite group
are defined via the representation theory of finite groups, see~\cite{Serre77}. For a definition in
the context of translation surfaces see~\cite{MatheusYoccozZmiaikou14}.}

Let $\tau_h$ be the translation of $X_{a,b}$ that exchanges
the polygon $00$ with $10$ and $01$ with $11$ of
Figure~\ref{fig:WindtreeQuotientSurface}. In other words,
$\tau_h$ exchanges the two copies drawn on the same horizontal.
Similarly, let $\tau_v$
be the one exchanging the copies $00$ with $01$ and $10$ with $11$.
\begin{tcblemma}{\cite{DelecroixHubertLelievre14}}{WindTreeKZSplitting:Veech}
Let $(a,b)$ be any parameter for the surface $L_{a,b}$. Let
$\tau_h$ and $\tau_v$ be the generators of the Deck transformations of
the covering $X_{a,b} \mapsto L_{a,b}$. Let us define for each
$s,t \in \{+1,-1\}$ the following subspaces:
\[
E^{st} := \{v \in H_1(X_{a,b}; \R):\ \tau_h v = s v \quad \tau_v v = t v\}.
\]
Then
\begin{enumerate}
\item $\dim E^{++} = 4$ and $\dim E^{+-} = \dim E^{-+} = \dim E^{--} = 2$,
\item there is a symplectic orthogonal splitting $H_1(X_{a,b}; \R) = E^{++} \oplus E^{+-} \oplus E^{-+} \oplus E^{--}$
and $E^{+-} \oplus E^{-+} \oplus E^{--} \subset \kerhol(X_{a,b})$,
\item each element of $\Aff(X_{a,b})$ preserves $E^{++}$ and $E^{--}$,
\item each element of $\Aff(X_{a,b})$ either preserves or exchanges $E^{+-}$ and $E^{-+}$.
\end{enumerate}
\end{tcblemma}

\begin{proof}
We refer to~\cite{DelecroixHubertLelievre14} for details. All these results
follow from the computation of the topology of the various quotients. First of
all $L_{a,b} = X_{a,b} / \langle \tau_h, \tau_v\rangle$ hence $E^{++} \simeq
H_1(L_{a,b}; \Q)$ which has rank 4. Next $X_{a,b} / \langle \tau_h \rangle$,
$X_{a,b} / \langle \tau_v \rangle$ and $X_{a,b} / \langle \tau_h, \tau_v
\rangle$ have all genus 3.
\end{proof}

\begin{tcblemma}{\cite{Pardo-commutator}}{WindTreeAffineGroupKernel}
Let $(a,b)$ be parameters so that $L_{a,b}$ is a Veech surface. Let
$H_1(X_{a,b}; \R) = E^{++} \oplus E^{+-} \oplus E^{-+} \oplus E^{--}$
be as in Lemma~\ref{lem:WindTreeKZSplitting:Veech}. Then the kernel of
the group action $\Aff(X_{a,b}) \to \GL(E^{+-} \oplus E^{-+} \oplus E^{--})$
is non-trivial.
\end{tcblemma}

\begin{proof}
Up to passing to a subgroup $A$ of index 2 of $\Aff(X_{a,b})$, we can assume that
each individual component $E^{st}$ for $(s,t) \in \{+-, -+, --\}$
is preserved. Since each component $E^{st}$ has rank 2 and is contained
in $\kerhol(X_{a,b})$ the general Hooper-Weiss theorem
(Theorem~\ref{thm:HooperWeissFirstKindII})
implies that the kernels of $A \to \GL(E^{st})$ are infinite.
Their images in the lattice $\Gamma(X_{a,b})$ are non-trivial normal subgroups and
hence non-elementary of the first kind (see Lemma~\ref{lem:NormalSubgroupLimitSet}). In particular
there exist $\phi_{+-}$, $\phi_{-+}$ and $\phi_{--}$ in $A$ so that $\phi_{st}$ acts
trivially on $E^{st}$ for each $st \in \{+-,-+,--\}$.

Let us first show that one can choose $\phi_{+-}$ and $\phi_{-+}$ so
that $g_1 := D(\phi_{+-})$ and $g_2 := D(\phi_{-+})$ do not have
common fixed points on $\partial \H^2$. Assume they have a common
fixed points. Because the group is non-elementary, there exists an
element $\psi$ that does not have any common fixed point with
$\phi_{+-}$. The elements $\psi^n \phi_{+-} \psi^{-n}$ for $n \geq 1$
are all pseudo-Anosovs which act trivially on $E_{+-}$. Moreover,
the fixed points of their derivatives $D(\psi^n \phi_{+-} \psi^{-n})$
are pairwise disjoint. Hence, there is one, say for $n = n_0$, whose
derivative has disjoint fixed points from $\phi_{-+}$ and we can
replace $\phi_{+-}$ by its conjugate $\psi^{n_0} \phi_{+-} \psi^{-n_0}$.

From now on, we assume that $g_{+-} := D(\phi_{+-})$ and
$g_{-+} := D(\phi_{-+})$ have disjoint fixed points. By the
ping-pong lemma (see Theorem~\ref{thm:PingPongLemma}
from Appendix~\ref{app:FuchsianGroups}) there exists an integer $m$
so that $(g_{+-})^m$ and $g_{-+}$ generate a free-group whose
non-trivial elements are hyperbolic. In particular the commutator
$[(\phi_{+-})^{m}, \phi_{-+}]$ whose derivative
is $D([(\phi_{+-})^{m}, \phi_{-+}]) = [(g_{+-})^n, g_{-+}]$
is a hyperbolic affine element. Moreover it acts trivially
on both $E^{+-}$ and $E^{-+}$.

In order to construct a hyperbolic element that acts trivially
on all of $E^{+-}$, $E^{-+}$ and $E^{--}$ one can repeat the same
argument with the two affine elements $[(\phi_{+-})^{m}, \phi_{-+}]$
and $\phi_{--}$.

Since $L_{a,b}$ is a Veech surface, so is $X_{a,b}$ and hence the limit of any non-trivial normal subgroup of
$\Gamma(X_{a,b})$ is the whole circle, see Lemma~\ref{lem:NormalSubgroupLimitSet}. To finish the proof recall that
the subspace $V<H_1(X_{a,b},\Q)$ defining the non-ramified covering $W_{a,b}\to X_{a,b}$ is generated by the cycles $Y$
and $Z$ described in equation (\ref{eq:WindTreeCycleCovering}). Both $X$ and $Y$ lie in $E^{+-} \oplus E^{-+} \oplus E^{--}$, hence by the Abelian lifting criterion every element in the kernel of $\Aff(X_{a,b}) \to \GL(E^{+-} \oplus E^{-+} \oplus E^{--})$
lifts to $\Aff(W_{a,b})$.

\end{proof}

\begin{tcbremark}{}{GeneralHooperWeissFirstKind2}
The proof of Lemma~\ref{lem:WindTreeAffineGroupKernel} does actually provide
an extension of the general Hooper-Weiss theorem. If $M$ is a Veech surface
and there exist subspaces $E_1$, \ldots, $E_k$ of kerhol in $H_1(M; \Z)$ such that
\begin{enumerate}
\item each $E_i$ has rank 2 and they are in direct sum
\item a subgroup of finite index $A$ of $\Aff(M)$ preserves each $E_i$
\end{enumerate}
Then the kernel of the induced action $A \to \GL(E_1) \times \ldots \times \GL(E_k)$
is infinite and its image in $\Gamma(A)$ is a non-trivial normal subgroup.
\end{tcbremark}

\section{Hooper-Thurston-Veech construction}
\label{sec:HooperThurstonVeechConstruction}

In Section~\ref{ssec:ThurstonVeechConstructionsIntro} we introduced the
Thurston-Veech construction, which given two multicurves filling a finite-type
topological surface $S$ produces a half-translation surface $M$ on $S$ admiting
a horizontal and a vertical cylinder decomposition where all cylinders involved
have the same modulus. As a consequence the surface $M$ has a lot of affine
pseudo-Anosov homeomorphisms. In the following paragraphs we detail an
extension of the Thurston-Veech construction based on the ideas of P.
Hooper~\cite{Hooper-infinite_Thurston_Veech} and illustrate it with several
examples. Infinite-type quadratic differentials constructed this way have
non-elementary Veech groups and their translation flows exhibit rich dynamics,
which we study in detail in a (second) forthcoming volume~\cite{DHV2}
dedicated mostly to dynamical aspects of infinite-type translation surfaces.

\begin{tcbremark}{}{halftranslationequalqdiff}
In this Chapter we construct surfaces where changes of coordinates are translations or half-translations. For the sake of a simple nomenclature, we will henceforth use the equivalent term quadratic differential (see Chapter~\ref{chap:Introduction}).
\end{tcbremark}

The whole construction is based in the following statement, which is the extension of Lemma~\ref{lem:ParabolicCylinderDecomposition} to infinite-type surfaces.
\begin{tcblemma}{}{CylinderDecompImpliesParabolic}
Let $M$ be a half-translation surface for which the horizontal direction is completely periodic and suppose that every maximal cylinder in the corresponding cylinder decomposition has modulus equal to $\frac{1}{\lambda}$ for some $\lambda\in\R^*$. Then there exists a unique affine automorphism $\varphi_h$ which fixes the boundaries of the cylinders and has derivative $D\varphi_h=h_\lambda=\begin{psmallmatrix}1 & \lambda\\ 0 & 1\end{psmallmatrix}$. Moreover, the automorphism $\varphi_h$ acts as a Dehn twist along the core curve of each cylinder.
\end{tcblemma}

\begin{tcbremark}{}{}
If $M$ is a half-translation surface and $\theta\in\R/2\pi\Z$ is a completely periodic direction for which every cylinder has modulus equal to $\frac{1}{\lambda}$ for some $\lambda\in\R^*$, then one can apply to $M$ the rotation $R_\theta\in\SO(2,\R)$ that takes $\theta$ to the horizontal direction and apply the preceding lemma. For example, if $\theta=\frac{\pi}{2}$, then there exists a unique affine automorphism $\varphi_v$ which fixes the boundaries of the vertical cylinders, acting on each cylinder as a Dehn twists and with derivative $D\varphi_v=v_{-\lambda}=\begin{psmallmatrix}1 & 0\\ -\lambda & 1\end{psmallmatrix}$.
\end{tcbremark}

\begin{tcbexercise}{}{}
Let $M$ be a half-translation surface and suppose that $\theta\in\R/2\pi\Z$ is a direction in which the translation flow decomposes $M$ into maximal cylinders $\{C_i\}_{i\in I}$ of modulus $\{\mu_i\}$. Moreover, suppose that for every $i\in I$ there exists $s\in\R^*$ such that $m_i:=s\mu_i \in \Z$. Show that there exist a unique $f\in\Aff(M)$ fixing the boundary of each cylinder $C_i$ whose derivative is, up to conjugation in $\SO(2,\R)$, equal to $\begin{psmallmatrix}1 & s\\ 0 & 1\end{psmallmatrix}$ .
\end{tcbexercise}

The proof of Lemma~\ref{lem:CylinderDecompImpliesParabolic} can be deduced from the preceding exercise and we leave it to the reader.

Let us remark that Lemma~\ref{lem:ParabolicCylinderDecomposition} gives an equivalence
between parabolic directions and affine multitwist directions for finite-type
surfaces. Though for infinite-type surfaces, only the implication stated in
Lemma~\ref{lem:CylinderDecompImpliesParabolic} is valid. Indeed, we already
encountered many counterexamples to the converse. First of all, in
Example~\ref{exa:BandsInTheInfiniteStaircase} we showed that the Veech group of
the infinite staircase has a parabolic element fixing the direction
$\theta=\frac{\pi}{4}$ and in this direction the staircase decomposes into two strips.
Secondly, if one considers the plane or any cyclic cover ramified over the origin,
its Veech group is $\SL(2,\R)$. In this case, any direction $\theta$ is parabolic but
none of them contains a cylinder. Finally, more exotic examples can be constructed from
Theorem~\ref{thm:AnyGroupIsVeech} and Exercise~\ref{exo:VeechSaysNothingAboutDyn}.
Note that all these examples have infinite area.

\begin{tcbquestion}{}{ParabolicMultitwitsCylDecom}
\label{Q:ParabolicMultitwitsCylDecom}
Let $M$ be an infinite-type translation surface of finite area.
Are parabolic directions in $M$ affine multitwist directions?
\end{tcbquestion}

\begin{tcbdefinition}{}{HooperThurstonVeechSurface}
A quadratic differential $M$ is called a \emphdef[Hooper-Thurston-Veech surface]{Hooper-Thurston-Veech} surface if it admits horizontal and vertical cylinder decompositions where all cylinders have the same modulus. The common modulus of the cylinders is called the modulus of the Hooper-Thurston-Veech surface.
\end{tcbdefinition}
\begin{tcbremark}{}{}
If $M$ is a Hooper-Thurston-Veech surface of modulus $\frac{1}{\lambda}$ then $\Aff(M)$ has two affine multitwists $\phi_h$ and $\phi_v$ with $D\phi_h=\begin{psmallmatrix}1 & \lambda\\0 & 1\end{psmallmatrix}$ and $D\phi_v=\begin{psmallmatrix}1 & 0\\ -\lambda & 1\end{psmallmatrix}$. Moreover, if $\lambda\geq 2$ then
\begin{equation}
\label{eq:GroupGLambda}
G_{\lambda}=\langle\begin{psmallmatrix}1 & \lambda\\0 & 1\end{psmallmatrix},\begin{psmallmatrix}1 & 0\\ -\lambda & 1\end{psmallmatrix}\rangle
\end{equation}
is a free subgroup of the Veech group $\Gamma(M)$ (see Theorem~\ref{thm:GLambdaGeometry} from Appendix~\ref{Appendix:FuchsianGroups}).
\end{tcbremark}

Note that finite-type Hooper-Thurston-Veech surfaces are those given by the
classical Thurston-Veech construction discussed in
Section~\ref{ssec:ThurstonVeechConstructionsIntro}. Translation coverings
provide an elementary way of obtaining Hooper-Thuston-Veech surfaces. More
precisely the lifting criterion for cylinders stated in
Lemma~\ref{lem:LiftingCylindersInCoverings} implies the following:
\begin{tcblemma}{}{LiftingHTV}
Let $M$ be a finite-type translation surface obtained by the Thurston-Veech construction
for the multicurves $\alpha=\cup_{i\in I}\alpha_i$ and $\beta=\cup_{j\in J}\beta_j$ and $G$ any group.
Let $p:\widetilde{M}\to M$ be the $G$-covering branched at most over the singularities $\Sigma$
of $M$ defined by a morphism $\rho:\pi_1(M\setminus\Sigma)\to G$ such that for all $i \in I$
we have $\rho(\alpha_i) = 1$ and for all $j \in J$ we have $\rho(\beta_j) = 1$.
Then $\widetilde{M}$ is a Hooper-Thurston-Veech surface.
\end{tcblemma}
The rest of the section is dedicated to the general construction of
Hooper-Thurston-Veech surfaces starting from a pair of transverse
multicurves in an infinite-type surface.

\subsection{From cylinder decompositions to bipartite graphs}
Let $S$ be a
topological surface. Recall from Definition~\ref{def:Multitwist} that a
multicurve in $S$ is a locally finite family of essential simple
closed curves in $S$, pairwise non-intersecting and pairwise non-homotopic.
A pair of simple closed curves $c$ and $c'$ are in \emphdef{minimal position} if
$|c\cap c'|=\iota(c,c')$, where $\iota(\cdot,\cdot)$ denotes the
geometric intersection between isotopy classes of curves.
\begin{tcbdefinition}{}{}
We say that a pair of multicurves $\alpha = \cup \alpha_i$, $\beta = \cup \beta_j$
\begin{itemize}
\item is in \emph{minimal position} if every pair $\alpha_i, \beta_j$ of
components of $\alpha$ and $\beta$ are in minimal position,
\item \emph{fills} $S$ if every connected component of
the complement of the union $\alpha \cup \beta$ is a topological disc with at most
one puncture.
\end{itemize}
\end{tcbdefinition}

\begin{tcbexercise}{}{}
Let $\alpha$ and $\beta$ two multicurves in a topological surface $S$. Show
that each component $\alpha_i$ of $\alpha$ intersects finitely
many components of $\beta$.

\emph{hint:} use the fact that a multicurve is by definition locally finite.
\end{tcbexercise}

\begin{tcbdefinition}{}{ConfigurationGraph}
Let $\alpha=\cup_{i\in I}\alpha_i$ and $\beta=\cup_{j\in J}\beta_j$ be two multicurves in a topological surface
$S$ in minimal position (not necessarily filling).
The \emphdef{configuration graph} of the pair $(\alpha,\beta)$ is the bipartite graph\footnote{We allow multiple edges since $\alpha_i$ and $\beta_j$ can intersect more than once. Strictly speaking $\cG(\alpha\cup\beta)$ is a multigraph, but we keep the term graph to facilitate the reading of this section.} $\cG(\alpha\cup\beta)$
whose vertex set is $I\sqcup J$ and where there is an edge between two vertices $i\in I$ and $j\in J$ for every
intersection between the corresponding curves $\alpha_i$ and $\beta_j$.
\end{tcbdefinition}

Before elaborating the construction for infinite-type surfaces, we explain how
the configuration graph already appeared in disguise form in the
Thurston-Veech contruction from
Section~\ref{ssec:ThurstonVeechConstructionsIntro} (more precisely
Exercises~\ref{exo:CylinderDecomposition} and~\ref{exo:FiniteThurstonVeech}). Let $S$ be a
finite-type surface and $\alpha=\cup \alpha_i$, $\beta=\cup \beta_j$ a pair of
filling multicurves in minimal position. As we have seen there is a
half-translation structure on $S$ which admits an horizontal
and a vertical affine multitwist along the multicurves $\alpha$ and $\beta$
respectively. To construct such half-translation structure, we considered the
intersection matrix $E = (\iota(\alpha_i, \beta_j))_{i\in I,j\in J}$ and shown that the heights of
the horizontal cylinders $\textbf{h}_h$ and the
heights of the vertical cylinders $\textbf{h}_v$ must satisfy $E \textbf{h}_v =
\lambda \textbf{h}_h$ and $E^t \textbf{h}_h = \lambda \textbf{h}_v$. Let us
consider the following matrix
$$
A=
\begin{pmatrix}
0 & E \\
E^t & 0
\end{pmatrix}.
$$
Then $A$ is the adjacency matrix of the configuration graph $\cG(\alpha \cup \beta)$
and the vector $\textbf{h} = (\textbf{h}_h, \textbf{h}_v)$ is an eigenvector of $A$.

Now, let us consider a Hooper-Thurston-Veech surface with
modulus $1/\lambda$. Let $H=\{H_i\}_{i\in I}$ and $V=\{V_j\}_{j\in J}$
be respectively the horizontal and vertical cylinders of $M$.
For every $i \in I$ let $\alpha_i$ be the
core curve of $H_i$ and for every $j \in J$ let $\beta_j$ be
the core curve of $V_j$.
Then $\alpha=\cup \alpha_i$ and $\beta=\cup \beta_j$ are multicurves
in $M$. Let $\textbf{h}: I \cup J \to \R_{>0}$ be the function which to
an index associates the width of the corresponding cylinder.
Then it satisfies the equation $A\textbf{h} = \lambda \textbf{h}$ where $A$ is the
adjacency operator of the graph $\cG(\alpha \cup \beta)$ defined by
\begin{equation}
\label{eq:AdjacencyOperator}
(A\textbf{h})(v) := \sum_{w \sim v} \textbf{h}(w)
\end{equation}
where the sum above is taken over the
vertices $w$ adjacent to $v$ in $\cG(\alpha \cup \beta)$
(with the multiplicity given by the number of edges
between $v$ and $w$).
\begin{tcbdefinition}{}{LambdaHarmonicFunction}
Let $\cG = (V,E)$ be a graph with vertices of finite degree.
Let $A: V^\R \to \R$ its adjacency operator.
A function $\textbf{h}: V \to \R$ that satisfies $A\textbf{h}=\lambda \textbf{h}$
is called a \emphdef{$\lambda$-harmonic function}.
\end{tcbdefinition}

\begin{tcbexercise}{}{}
Let $H=\{H_i\}$ and $V=\{V_j\}$ be horizontal and vertical cylinder decompositions
of a quadratic differential $M$, and let $\cG(\alpha\cup\beta)$ be the configuration
graph of the corresponding core curves.
Suppose that each horizontal cylinder intersects finitely many vertical cylinders and \emph{vice versa}. Show that the union of the core curves of the cylinders in $H\cup V$ fills $M$, \ie their complement in $M$ is a union of discs with at most one puncture. On the other hand, show that the core curves of the horizontal and vertical cylinder decompositions of the infinite-type translation surface $M$ in Example~\ref{exa:InfiniteStepSurface}, Chapter~\ref{ch:TopologyGeometry}, do not fill.
\end{tcbexercise}

We now present the main construction that extends to infinite-type surfaces
the Thurston-Veech construction of Exercise~\ref{exo:FiniteThurstonVeech}.
\begin{tcbtheorem}{}{HooperThurstonVeechConstruction}
Let $S$ be topological surface, $\alpha=\cup_{i\in I}\alpha_i$ and $\beta=\cup_{j\in J}\beta_j$ two multicurves in minimal position whose union fills $S$.
Suppose that:
\begin{itemize}
\item \textbf{C1}. The configuration graph $\cG(\alpha\cup\beta)$ has finite degree (\ie there is an upper bound on the degree of the vertices in $\cG(\alpha \cup \beta)$),
\item \textbf{C2}. For every component $D$ of the complement of $\alpha\cup\beta$ in $S$, its boundary $\partial D$ in $S$ is connected.
\item \textbf{C3}. For every component $D$ of the complement of $\alpha\cup\beta$ for which $\partial D$ intersects infinitely many curves in $\alpha\cup\beta$, then $D$ is a disk \emph{without} punctures.
\end{itemize}
Let $\textbf{h}: \cG(\alpha \cup \beta) \to \R_{>0}$ be a positive $\lambda$-harmonic function, where $\lambda > 1$.

Then, there exist a half-translation structure $M=M(\alpha,\beta,\textbf{h})$
on $S$ which is a Hooper-Thurson-Veech surface of modulus $1/\lambda$ whose
horizontal cylinders (respectively vertical cylinders) have core curves
$\{\alpha_i\}_{i \in I}$ (respectively $\{\beta_j\}_{j \in J}$).
\end{tcbtheorem}

\begin{tcbremark}{}{ExistenceEigenfunction}
Any connected infinite graph of finite degree has a positive $\lambda$-harmonic function if $\lambda$ is large enough. Moreover, we always have that $\lambda\geq 2$. A more detailed discussion of this fact can be found in Appendix~\ref{Apendix:EigenvaluesInfiniteGraphs}.
For the configuration graph $\cG(\alpha\cup\beta)$ connectedness and finite degree are guaranteed by the fact
that $\alpha\cup\beta$ fills the surface $S$ and contidion
\textbf{C1} above respectively.
\end{tcbremark}

\begin{tcbremark}{}{}
Note that if $S$ is of finite-type then the conditions \textbf{C1}, \textbf{C2} and \textbf{C3} above are satisfied automatically and the result follows from the Thurston-Veech construction as explained in Section~\ref{ssec:ThurstonVeechConstructionsIntro}.
\end{tcbremark}

\begin{proof}
The union $\alpha\cup\beta$ of the multicurves $\alpha$ and $\beta$ define a graph embedded in $S$: the vertices are points in $\bigcup_{(i,j)\in I\times J}\alpha_i\cap\beta_j$ and edges are the connected components of $\alpha\cup\beta$ minus its vertices. 
It is important not to confuse the (geometric) graph $\alpha\cup\beta$ with the (abstract) configuration graph $\cG(\alpha\cup\beta)$.

To define the quadratic differential we consider a dual embedded graph $(\alpha\cup\beta)^*$ in $S$. When all complementary components of $S\setminus\alpha\cup\beta$ have a compact boundary $(\alpha\cup\beta)^*$ is just the dual graph of $\alpha\cup\beta$. In the case when there are complementary components of $S\setminus\alpha\cup\beta$ whose boundary is formed by infinitely many vertices and edges we need to define $(\alpha\cup\beta)^*$ more carefully. The vertices of $(\alpha\cup\beta)^*$ are defined as follows:
\begin{itemize}
\item \textbf{A1}. For every connected component $D$ of $S\setminus \alpha\cup\beta$ such that $\partial D$ is formed by finitely many edges of $\alpha\cup\beta$ choose a unique point $v_D$ inside the connected component. If $D$ is a punctured disc, then choose $v_D$ to be the puncture.
\item \textbf{A2}. If $D$ is a connected component of $S\setminus \alpha\cup\beta$ such that $\partial D$ is formed by infinitely many edges of $\alpha\cup\beta$, then by hypothesis $\textbf{C3}$ above $D$ is homeomorphic to a disc. Moreover, $D\cup\partial D$ has only one end to which we associate an ideal point $v_D$. We think of $v_D$ as a vertex ``at infinity''  of $(\alpha\cup\beta)^*$, see Figure~\ref{fig:DualGraph}.
\end{itemize}

\begin{figure}[!ht]
\begin{center}
\begin{minipage}{.45\textwidth}
\begin{center}
\includegraphics[scale=0.7]{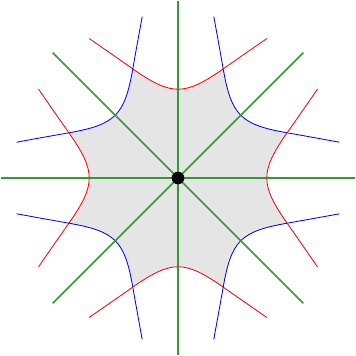}
\end{center}
\subcaption{$\partial D$ has finitely many sides}
\end{minipage}
\begin{minipage}{.45\textwidth}
\begin{center}
\includegraphics[scale=0.8]{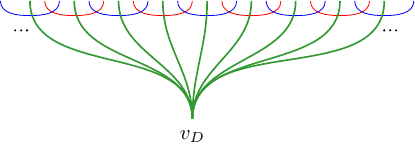}
\end{center}
\subcaption{$\partial D$ has infinitely many sides}
\end{minipage}
\end{center}
\caption{Edges of the graph $(\alpha\cup\beta)^*$. }
\label{fig:DualGraph}
\end{figure}

The points $v_D$, as chosen above, are the vertices of $(\alpha\cup\beta)^*$. Two vertices in this graph are joined by an edge in $S$ if the closures of the corresponding connected components of $S\setminus \alpha\cup\beta$ share an edge of $\alpha\cup\beta$. Edges are chosen to be pairwise disjoint. For edges with vertices at infinity, edges are chosen to be proper rays defining the end of $D\cup\partial D$ associated to $v_D$, see Figure~\ref{fig:DualGraph}.

Given that $\alpha\cup\beta$ fills and the boundary $\partial D$ of each complementary component $D$ of $\alpha\cup\beta$ is connected (hypothesis \textbf{C2}), every connected component $S\setminus (\alpha\cup\beta)^*$ is a topological quadrilateral which contains a unique vertex of $\alpha\cup\beta$. Hence there is a well defined bijection between edges in the abstract graph $\mathcal{G}(\alpha\cup\beta)$ and the set of these quadrilaterals. This way, for every edge $e\in E(\mathcal{G}(\alpha\cup\beta))$ we denote by $R_e$ the closure in $S$ of the corresponding topological quadrilateral with the convention to add to $R_e$ vertices $v_D$ corresponding to punctures
in $S$ and points at infinity (as chosen above) so that $R_e$ has four vertices.

Note that there are only two sides of $R_e$ intersecting the multicurve $\alpha$, which henceforth are called \emph{vertical sides}. The other two sides are in consequence called \emph{horizontal}, see Figure~\ref{fig:MakingRectanglesEuclidean}.

We now build a half-translation atlas on $S$ by mapping the topological
quadrilaterals $R_e$ of the dual graph $(\alpha \cup \beta)^*$ into Euclidean
rectangles. As in the finite-type case,
any positive $\lambda$-harmonic function $\textbf{h}: \cG(\alpha \cup \beta) \to \R_{>0}$ gives
us compatible heights of horizontal and vertical cylinders. As noted in Remark~\ref{rk:ExistenceEigenfunction}, given that the union of the multicurves $\alpha$ and $\beta$ fill the surface and hypothesis \textbf{C1}, we can always count on the existence of such a function. Let us define the maps
\begin{equation}
    \label{eq:HTVProjections}
p_\alpha:E(\cG(\alpha\cup\beta)) \to V(\cG(\alpha\cup\beta))
\qquad \text{and} \qquad
p_\beta: E(\cG(\alpha\cup\beta)) \to V(\cG(\alpha\cup\beta))
\end{equation}
which associates to an edge $e$ of the configuration graph $\cG(\alpha\cup\beta)$
its endpoints $p_\alpha(e)$ in $I$ and $p_\beta(e)$ in $J$. Our aim is to
construct a homeomorphism
\begin{equation}
\label{eq:HTVpsie}
R_e\longrightarrow [0,\textbf{h}\circ p_\beta(e)] \times [0,\textbf{h}\circ p_\alpha(e)]
\end{equation}
that will determine a unique half-translation atlas for $S$ with the required
properties (see Figure~\ref{fig:MakingRectanglesEuclidean}).

\begin{figure}[!ht]
\begin{center}
\includegraphics{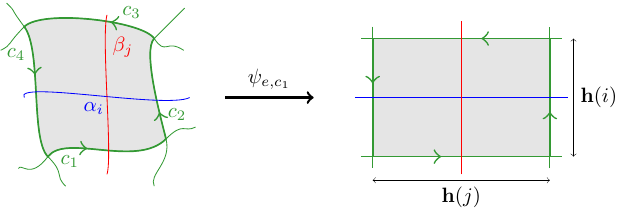}
\end{center}
\caption{Transforming topological Rectangles into Euclidean rectangles.}
\label{fig:MakingRectanglesEuclidean}
\label{fig:MakingRectanglesEuclidean}
\end{figure}

In order to ensure compatibility of the atlas we first consider mapping
of the edges. For each oriented edge $c$ of $(\alpha \cup \beta)^*$ crossing
a component $\alpha_i$ of $\alpha$ we choose an homeomorphism
$\theta_c: c \to [0, \textbf{h}(i)]$ such that
\begin{itemize}
\item the start point of $c$ is mapped to $0$ and its endpoint to $\textbf{h}(i)$,
\item $c^{-1}(\textbf{h}(i)/2)$ is the intersection of $c$ with $\alpha_i$,
\item if $c'$ is the same edge with opposite orientation $\theta_{c'}(x) = \textbf{h}(i) - \theta_c(x)$.
\end{itemize}

Then for each rectangle $R_e$ (corresponding to an
edge $e$ of $\cG(\alpha \cup \beta)$) and one of its horizontal
sides $c_1$ we construct a map $\psi_{e,c_1}$ as follows.
Choose a cyclic ordering of the sides $c_1, c_2, c_3, c_4$ following
the orientation of the surface. Then define on the boundary $\partial R_e$
\[
\psi_{e,c_1}: x \mapsto
\begin{cases}
(\theta_{c_1}(x), 0) & \text{if $x \in c_1$} \\
(\textbf{h}(j), \theta_{c_2}(x)) & \text{if $x \in c_2$} \\
(\textbf{h}(j) - \theta_{c_3}(x), \textbf{h}(i)) & \text{if $x \in c_3$} \\
(0, \textbf{h}(i) - \theta_{c_4}(x)) & \text{if $x \in c_4$}
\end{cases}
\]
Note that by our choice of compatibilities of the maps $\theta_c$ we have that
$\psi_{e,c_1}|_{\partial R_e} = - \psi_{e, c'_3} |_{\partial R_e}$ where
$c'_3$ is the edge $c_3$ with reversed orientation.

Now, we extend each $\psi_{e,c_1}$ into a homeomorphism such that
\begin{itemize}
\item $\alpha_i \cap R_e$ is mapped to $[0,\textbf{h}(j)] \times \{\textbf{h}(i)/2\}$,
\item $\beta_j \cap R_e$ is mapped to $\{\textbf{h}(j)/2\} \times [0,\textbf{h}(i)]$,
\item $\psi_{e,c_1} = - \psi_{e,c'_3}$.
\end{itemize}

We claim that the family $\{\psi_{e,c}\}$ we constructed above extends to a unique
half-translation atlas on $S$. First of all, the coordinate change between
$\psi_{e,c_1}$ and $\psi_{e,c'_3}$ is a point reflection. Notice also that each
point in the interior of a rectangle is covered by two charts. Now, for each
edge $c$ in the intersection of two neighbouring rectangles $R_e$ and $R_{e'}$, the
compatibility of the maps $\psi_e$ on the boundaries via the maps $\theta_c$ can be
used to provide charts compatible with the atlas.

Note that the argument is similar to the one we used in the introduction to show that a
constructive translation surface (obtained from gluing polygons) is also a geometric
translation surface (endowed with a translation atlas). See
item~\ref{item:ConstructiveGeometric} of Theorem~\ref{thm:EquivalencesInDefinitions}.

We denote the resulting half-translation surface $M(\alpha, \beta, \textbf{h})$.

Now, for every $i\in I$, the curve $\alpha_i$ is the core curve of the horizontal
cylinder $H_i := \cup_{e\in p_{\alpha}^{-1}(i)} R_e$. Because $\textbf{h}$ is harmonic we have
$$
\sum_{e\in p_\alpha^{-1}(i)}\textbf{h}(p_{\beta}(e))=\sum_{j\sim i}\textbf{h}(j)=\lambda \textbf{h}(i).
$$
In other words, the circumference $\sum_{e\in p_\alpha^{-1}(i)}\textbf{h}(p_{\beta}(e))$ of
$H_i$ is $\lambda$ times its height $\textbf{h}(i)$. That is, the modulus $\mu(H_i)$ of $H_i$
is equal to $\frac{1}{\lambda}$. The same computation with $\beta_j$ shows that the vertical cylinders
$V_j := \cup_{e\in p_{\beta}^{-1}(j)} R_e$ have core curve $\beta_j$
and modulus $\frac{1}{\lambda}$.

\end{proof}

\begin{tcbremark}{}{ParabolicAutomorphismsSingularities}
\label{Lemma:ParabolicAutomorphismsSingularities}
Let $M=M(\alpha,\beta,\textbf{h})$ be a Hooper-Thurston-Veech surface as in Theorem~\ref{thm:HooperThurstonVeechConstruction}. For every connected component $D$ of $S\setminus\alpha\cup\beta$ we have that:
\begin{enumerate}
\item If $\partial D$ is formed by $2k$ edges of $\alpha\cup\beta$ then the vertex $v_D$ of $(\alpha\cup\beta)^*$ is a conical singularity of total angle $k\pi$ in\footnote{For $k=2$, the vertex $v_D$ is a regular point of $M$.} $\Sing(M)$ and all other points of $D$ are regular.
\item If $\partial D$ is formed by infinitely many edges of $\alpha\cup\beta$ then the vertex $v_D$ of $(\alpha\cup\beta)^*$ is either an infinite angle singularity or a wild singularity.
\end{enumerate}
This follows from the fact that in the proof of Theorem~\ref{thm:HooperThurstonVeechConstruction} we are glueing Euclidean rectangles with a common vertex at $v_D$. When $v_D$ is an infinite angle or wild singularity, it has a bi-infinite rotational component\footnote{\ie the map $\widetilde{\rm dir}$ restricted to each rotational component defines a bijection with $\R$, see (\ref{eq:InjectiveLift}), in Section~\ref{SSEC:RotationalComponents}.}. Moreover, the existence of such bi-infinite rotational component among the rotational components of $v_D$ does not depend on the positive $\lambda$-harmonic function \textbf{h}. However, not all rotational components of linear approaches to $v_D$ are bi-infinite. In the next section we see that baker's surfaces are examples of Hooper-Thurston-Veech surfaces. All surfaces in this family have a single wild singularity which has bi-infinite rotational components but also bounded rotational components if $\lambda>2$.
\end{tcbremark}

In the rest of this Chapter we illustrate with several examples how to determine explicitly the quadratic differential $M(\alpha,\beta,\textbf{h})$.

\begin{tcbexample}{}{HTVStaircase}
In the following paragraphs we show that the infinite staircase and all $\lambda$-staircases are Hooper-Thurston-Veech surfaces homeomorphic to the Loch Ness monster. Figure~\ref{Fig:MulticurvesInfiniteStaircase} shows a Loch Ness monster with two multicurves $\alpha$ and $\beta$ in minimal position satisfying contidions \textbf{C1}--\textbf{C4} of Theorem~\ref{thm:HooperThurstonVeechConstruction} and the corresponding configuration graph $\mathcal{G}(\alpha\cup\beta)$, to which we refer hereafter as the $\Z$-graph because its vertices are indexed by $\Z$.
\begin{figure}[H]
\begin{center}
\begin{minipage}{.45\textwidth}
\begin{center}
\includegraphics[scale=.5]{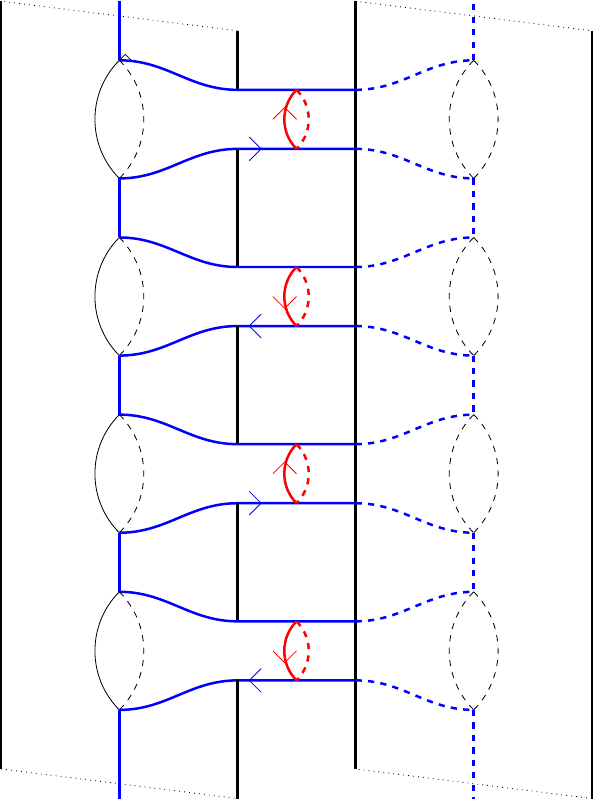}
\end{center}
\subcaption{Two oriented multicurves $\alpha$ (in blue) and $\beta$ (in red) in the Loch Ness monster for which the Hooper-Thurston-Veech construction produces the infinite staircase.}
\end{minipage}
\begin{minipage}{.45\textwidth}
\begin{center}
\includegraphics[scale=1]{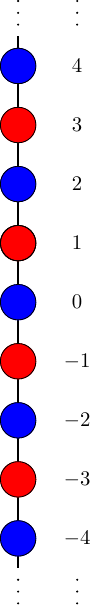}
\end{center}
\subcaption{The $\Z$-graph}
\end{minipage}
\begin{minipage}{.45\textwidth}
\begin{center}
\includegraphics[scale=1.2]{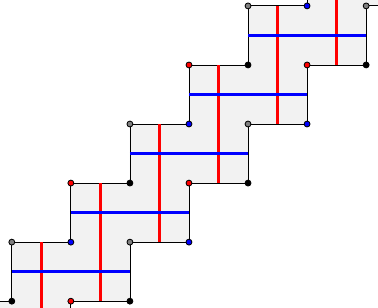}
\end{center}
\end{minipage}
\end{center}
\caption{The infinite staircase as a Hooper-Thurson-Veech surface.}
\label{Fig:MulticurvesInfiniteStaircase}
\end{figure}

A $\lambda$-harmonic function on the $\Z$-graph satisfies
\begin{equation}
	\label{eq:RecurrenceRelationStaircase}
\textbf{h}(n+1)=\lambda \textbf{h}(n)-\textbf{h}(n-1)
\end{equation}
for every $n\in\Z$. This is a linear recurrence relation whose characteristic polynomial is $x^2-\lambda x+1$. By Remark~\ref{rk:ExistenceEigenfunction} (see also Appendix~\ref{Apendix:EigenvaluesInfiniteGraphs}) we know that $\lambda\geq 2$. Hence let $r_+:=\frac{\lambda+\sqrt{\lambda^2-4}}{2}$ and $r_-:=\frac{1}{r_+}$ be the corresponding roots. There are two cases to consider:
\begin{itemize}
\item $\lambda=2$. Since $r_+=r_-=1$ the general solution of (\ref{eq:RecurrenceRelationStaircase}) is of the form $\textbf{h}(n)=C+Dn$. Positive solutions only exist if $D=0$ and hence the only
positive $\lambda$-harmonic function in this case (up to scaling) is $\textbf{h}(n)=1$ for all $n\in\Z$.
\item $\lambda>2$. For each fixed value of $\lambda$, the general solution of (\ref{eq:RecurrenceRelationStaircase}) is a linear combination of the positive eigenfunctions $\textbf{h}_1(n)=r_+^n$ and $\textbf{h}_2(n)=r_+^{-n}$.
\end{itemize}

We now explain a systematic way to picture $M(\alpha,\beta,\textbf{h})$. In general it can be difficult to visualize the graph $\alpha\cup\beta$ and how the rectangles lie in $S$. The following steps are intented to overcome this difficulty.
\begin{enumerate}
\item \emph{Orientations}. Fix an orientation on each of the curves in $\alpha\cup\beta$ and the standard orientation in $\R^2$. Let $\psi_{e,c_1}$ be an orientation preserving homeomorphism as defined in the proof of Theorem~\ref{thm:HooperThurstonVeechConstruction}.
The orientation of each component of $\alpha$ and $\beta$,
  induces an orientation of the arcs $\psi_{e,c1}(\alpha \cap R_e)$ and $\psi_{e,c1}(\beta \cap R_e)$.
We say that $R_e$ is \emph{positive} (respectively \emph{negative}) if the algebraic intersection (w.r.t. the orientations just described) of $\psi_{e,c_1}(\alpha\cap R_e)$ with $\psi_{e,c_1}(\beta\cap R_e)$ is positive (respect. negative). In Figures \ref{fig:OrientationRectangles1} and \ref{fig:OrientationRectangulos} we illustrate positive and negative rectangles.

\item \emph{Gluing rectangles}. The orientation of elements in $\alpha\cup\beta$ defines two permutations
\begin{equation}
	\label{eq:RibbonStructureFromOrientation}
r:E(\mathcal{G}(\alpha\cup\beta))\to E(\mathcal{G}(\alpha\cup\beta))\hspace{2mm}\text{and}\hspace{2mm} u:E(\mathcal{G}(\alpha\cup\beta))\to E(\mathcal{G}(\alpha\cup\beta))
\end{equation}
as follows. Let $e\in E(\mathcal{G}(\alpha\cup\beta))$ be determined by a point $p_{ij}\in\alpha_i\cap\beta_j$ and let $p_{ij'}\in\alpha_i\cap\beta_j$ be the point of $\alpha_i\cap\beta$ that follows from $p_{ij}$ when considering the orientation of $\alpha_i$. If $e'\in E(\mathcal{G}(\alpha\cup\beta))$ denotes the edge determined by $p_{ij'}$ then $r(e)=e'$. The permutation $u$ is defined in an analogous way by considering points in $\alpha\cap\beta_j$. Note that for every vertex $i\in I$ (respectively $j\in J$) of $\cG(\alpha\cup\beta)$, $r$ defines a cyclic permutation of edges adjacent to $i$ (respect. $u$ for edges adjacent to $j$).

As discussed in the proof of Theorem~\ref{thm:HooperThurstonVeechConstruction}, the rectangles $[0,\textbf{h}\circ p_\beta(e)]\times[0,\textbf{h}\circ p_\alpha(e)]$ are glued along parallel edges following the configuration of the multicirves $\alpha\cup\beta$. More precisely: for each $e\in E$ is glued to the right side of the Euclidean rectangle $[0,\textbf{h}\circ p_\beta(e)]\times[0,\textbf{h}\circ p_\alpha(e)]$ to the left side of $R_{r(e)}$ using a translation if both rectangles $R_e$ and $R_{r(e)}$ are positive or using a half-translation if both rectangles are negative. On the other the upper side of $[0,\textbf{h}\circ p_\beta(e)]\times[0,\textbf{h}\circ p_\alpha(e)]$ is glued to the lower side of $R_{u(e)}$ always using a translation.
\item \emph{Singularities}. Remove all corners of rectangles glued as above which correspond to marked points in $S$ or which produce singularities which are not conical.
\end{enumerate}

\begin{figure}[H]
\begin{center}
\includegraphics[scale=1]{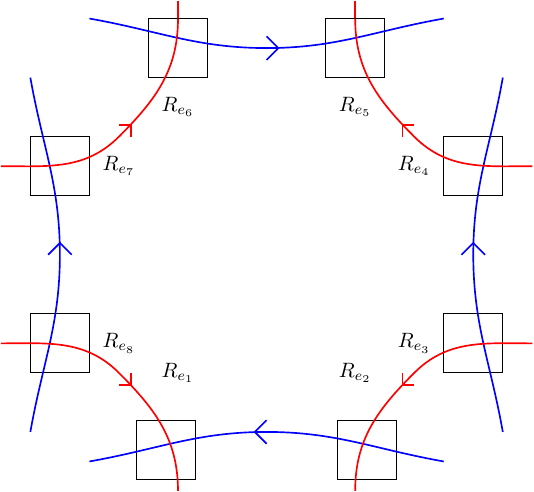}
\end{center}
\caption{Orientation of the multicurves $\alpha\cup\beta$ inducing an orientation on the rectangles $R_e$ in $S$.\\}
\label{fig:OrientationRectangles1}
\end{figure}

\begin{figure}[H]
\begin{center}
\includegraphics[scale=1]{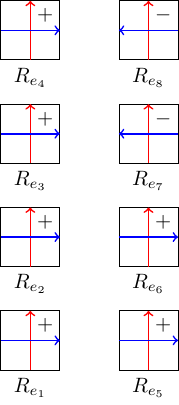}
\end{center}
\caption{Orientation of the Euclidean rectangles $R_e$ induced by the homeomorphism $\psi_e$.}
\label{fig:OrientationRectangulos}
\end{figure}

To visualize concretely the steps explained above consider the oriented multicurves depicted in Figure~\ref{Fig:MulticurvesInfiniteStaircase} and the $\lambda$-harmonic function $\textbf{h}(n)=1$ for all $n\in\Z$ of the adjacency operator of the $\Z$-graph (\ie case $\lambda=2$). For each edge $e$ in the $\Z$-graph the rectangle $R_e$ is positive and identified with a unit square, therefore all gluings of these are done using translations. Therefore $M(\alpha,\beta,\textbf{h})$ is a translation surface. Moreover, it is one of our favorites examples in this book: the infinite staircase. See Figure~\ref{Fig:MulticurvesInfiniteStaircase}.

\end{tcbexample}

\begin{tcbexercise}{}{LambdaStairCasesHTV}
For $\lambda>2$, consider any positive linear combination $\textbf{h}$ of the $\lambda$-harmonic functions of the $\Z$-graph $\textbf{h}_1(n)=r_+^n$ and $\textbf{h}_2(n)=r_+^{-n}$. Show that $M(\alpha,\beta,\textbf{h})$ is a translation surface. Draw the rectangles corresponding rectangles and identifications between them, and compare your results with Figure~\ref{fig:LambdaStaircase} in Chapter~\ref{chap:Introduction}.
\end{tcbexercise}

\begin{tcbexercise}{}{}
Consider the oriented multicurves $\alpha$ and $\beta$ on the Loch Ness monster shown in the Figure below. Show that the configuration graph $\mathcal{G}(\alpha\cup\beta)$ is isomorphic to the $\Z$-graph. Construct the quadratic differential $M(\alpha,\beta,\textbf{h})$ following the steps explained in Example \ref{exa:HTVStaircase} for the $2$-harmonic function $\textbf{h}(n)=1$, for all $n\in\Z$. Show that $M(\alpha,\beta,\textbf{h})$ \emph{is not} isomorphic to the infinite-staircase and explain what went wrong. \emph{Hint:} determine the number of infinite angle singularities that $M(\alpha,\beta,\textbf{h})$ has in this case and compare with the case of the infinite-staircase. Is $M(\alpha,\beta,\textbf{h})$ homeomorphic to the Loch Ness monster?
\begin{center}
\includegraphics[scale=.8]{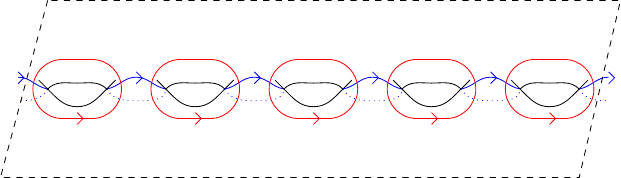}
\end{center}

\end{tcbexercise}

\begin{tcbexercise}{}{}
Consider the oriented multicurves $\alpha=\{\alpha_i\}_{i\geq 0}$ and $\beta=\{\beta_j\}_{j\geq 0}$ on the surface $S$ depicted in the figure below. Let $\textbf{h}$ be any positive $\lambda$-harmonic function of $\mathcal{G}(\alpha\cup\beta)$. Show that
the quadratic differential $M(\alpha,\beta,\textbf{h})$ is not the square of an Abelian differential, \ie there are changes of coordinates in the corresponding atlas that are half-translations.
\begin{center}
\includegraphics[scale=1]{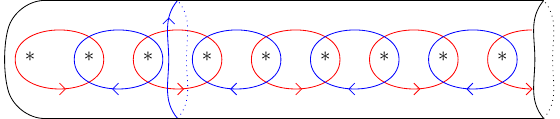}
\end{center}
\end{tcbexercise}


\subsubsection{Constructing Hooper-Thurston-Veech surfaces from bipartite graphs}
In the preceding section, we have shown that a Hooper-Thurston-Veech surface is encoded by a pair of multicurves $\alpha$ and $\beta$ (the core curves of the horizontal and vertical cylinders) and a positive $\lambda$-harmonic function $\textbf{h}: \cG(\alpha\cup\beta) \to \R_{>0}$ of the configuration graph (which encodes the heights of cylinders). As explained in Example~\ref{exa:HTVStaircase}, an orientation of the components of $\alpha$ and $\beta$ defines a pair of permutations $r,u$ on the edges of $\cG(\alpha \cup \beta)$ (see~\eqref{eq:RibbonStructureFromOrientation}) which encode how cylinders intersect. In this section we explain how to construct a Hooper-Thurston-Veech from an abstract bipartite graph of finite degree, a pair of permutations and a positive $\lambda$-harmonic function. This is the way the construction is presented in~\cite{Hooper-infinite_Thurston_Veech}.

\begin{tcbdefinition}{}{RibbonStructure}
Let $\cG$ be a bipartite graph with vertices $I \sqcup J$. A \emphdef{ribbon structure} on $\cG$ is a choice
for each vertex $v$ of a cyclic ordering to the edges adjacent to $v$. In other words, it is a pair of
permutations $r: E(\cG) \to E(\cG)$ and $u: E(\cG) \to E(\cG)$ such that for each edge adjacent to
$i \in I$ (respectively $j \in J$) the orbit of $i$ by $r$ (respectively of $j$ by $u$) are the edges adjacent
to $i$ (respectively $j$).
\end{tcbdefinition}

\begin{tcbtheorem}{\cite{Hooper-infinite_Thurston_Veech}}{}
Let $\cG$ be a connected bipartite graph of finite degree endowed with a ribbon
structure. Then for every positive $\lambda$-harmonic function $\textbf{h}$
of $\cG$ there is an Abelian differential $M=M(\cG,h)$ which is a
Hooper-Thurston-Veech surface of modulus $\frac{1}{\lambda}$. Let $H=\{H_i\}_{i\in I}$ and
$V=\{V_j\}_{j\in J}$ be the corresponding horizontal and vertical cylinder
decompositions; and $\alpha$ and $\beta$ be the multicurves given by their core curves. Then $\cG$ is the configuration graph $\cG(\alpha\cup\beta)$.
\end{tcbtheorem}
\begin{proof}
We sketch the main ingredients of the proof. Let $V(\cG)=I \sqcup J$ be the set of vertices and $E(\cG)$ the set of edges of the bipartite graph $\cG$. We have natural projections:
$$
p_{I}:E(\cG)\to I\hspace{1cm}\text{and}\hspace{1cm}p_{J}:E(\cG)\to J
$$
which assign to an edge $e=ij$ the vertices $p_I(e)=i$ and $p_j(e)=j$. Consider a positive $\lambda$-harmonic function $\textbf{h}$ of $\cG$ and for each $e\in E(\cG)$ define $R_e=[0,\textbf{h}\circ p_J(e)]\times[0,\textbf{h}\circ p_I(e)]$. The ribbon structure of $\cG$ defines a pair of permutations
$$
r:E(\cG)\to E(\cG)\hspace{1cm}\text{and}\hspace{1cm}u:E(\cG)\to E(\cG)
$$
as follows: $r(e)=\sigma_{p_I(e)}(e)$ and $u(e)=\sigma_{p_J(e)}(e)$, where $\sigma_{p_I(e)}$ and $\sigma_{p_J(e)}$ are the permutations given by the ribbon structure on edges adjacent to $p_I(e)$ and $p_J(e)$ respectively. Thus, the $r$-orbit (respectively $u$-orbit) of $e$ is formed by the all edges adjacent to $I(e)$ (respect. $J(e)$). To obtained the desired translation surface glue using a translation the right side of the Euclidean rectangle $R_e$ to the left side of $R_{r(e)}$, and the upper side of $R_e$ to the bottom side of $R_{u(e)}$. The rest of the proof follows from arguments analogous to the ones used in the proof of Theorem~\ref{thm:HooperThurstonVeechConstruction}.
\end{proof}
\begin{tcbremark}{}{}
Note that the steps of the construction of $M(\cG,\textbf{h})$ are very similar to those explained in Example~\ref{exa:HTVStaircase}, except that we do not have to care about orientations for there is no surface involved in the input data. On the other hand, if we consider $\alpha$ and $\beta$ the multicurves formed by the core curves of the cylinder decompositions $H$ and $V$ of $M(\cG,\textbf{h})$ respectively, then $M(\cG,\textbf{h})$ is equal to $M(\alpha,\beta,\textbf{h})$, the Hooper-Thurston-Veech surface given by Theorem~\ref{thm:HooperThurstonVeechConstruction}.
\end{tcbremark}

\begin{tcbexercise}{}{FiniteAreaHTV}
Let $\textbf{h}$ be a positive $\lambda$-harmonic function used to construct a Hooper-Thurston-Veech surface $M$ of the form $M(\alpha,\beta,\textbf{h})$ or $M(\cG, \textbf{h})$. Show that the area of $M$ is given by
\[
\area(M) = \frac{\lambda}{2} \sum_{v \in V(\cG)} \textbf{h}_v^2.
\]
Deduce that $\textbf{h}$ is a sequence in $\ell^2(V)$ if and only the surface $M(\cG, \textbf{h})$ has finite area.
\end{tcbexercise}

\begin{tcbexample}{The $\N$-graph}{NCayleyGraph}
Let $\cG$ be the graph depicted in Figure~\ref{fig:NNIGraph}. Henceforth we call $\cG$ the $\N$-graph and as indicated in the Figure vertices are indexed by the non-negative integers.
Given that every vertex except one has degree 2, there is a unique ribbon structure on this graph. As in Example~\ref{exa:HTVStaircase}, a positive $\lambda$-harmonic function $\textbf{h}$ on $\cG$ satisfies the (linear recurrence) relation $\textbf{h}(n)=\lambda \textbf{h}(n-1)-\textbf{h}(n-2)$, for $n\geq 2$ with initial condition $\lambda \textbf{h}(0)=\textbf{h}(1)$. The corresponding characteristic polynomial is $x^2-\lambda x+1$. For $\lambda=2$ there is only one root and the general solution is of the form  $\textbf{h}(n)=A+Bn$, $A,B\in\R$.

\begin{figure}[H]
\begin{center}
\includegraphics[scale=1]{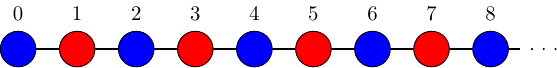}
\end{center}
\caption{The $\N$-graph.}
\label{fig:NNIGraph}
\end{figure}
Let's now consider the case of simple roots $r_+>r_-$ corresponding to $\lambda>2$. The general solution in this case is of the form $\textbf{h}(n)=Ar_+^n+Br_-^n$.

\begin{tcbexercise}{}{}
Draw the infinite area surface $M(\cG,\textbf{h})$ for $\cG$ equal to $\N$-graph and $\textbf{h}$ is the $2$-harmonic function $\textbf{h}(n)=n+1$. For $\lambda>2$ and $\textbf{h}(0)=r_+-r_-$ show that $\textbf{h}(n)=(r_+^{n+1}-r_-^{n+1})$ is a $\lambda$-harmonic function and draw the corresponding Hooper-Thurston-Veech surface $M(\cG,\textbf{h})$. You should get in both cases something that resembles a staircase going upwards. One can turn the $\N$-graph into a metric graph by declaring that each edge has length 1. Consider then the subset of each surface corresponding to the ball of radius $n$ in the $\N$-graph centered at the origin. Denote this subset $Y_n(\lambda)$.  Compare the growth of the area of $Y_n(\lambda)$ for parameters $\lambda=2$ and $\lambda>2$.
\end{tcbexercise}
\end{tcbexample}

\begin{tcbexample}{}{ModifiedNGraph}
\emph{Modified $\N$-graph}. Let $\cG_k$ be the graph obtained by adding $k \geq 1$ vertices of valence 1 to the vertex $0$ in the $\N$-graph, see Figure~\ref{fig:MNNIGraph}.

\begin{figure}[H]
\begin{center}
\includegraphics[scale=.8]{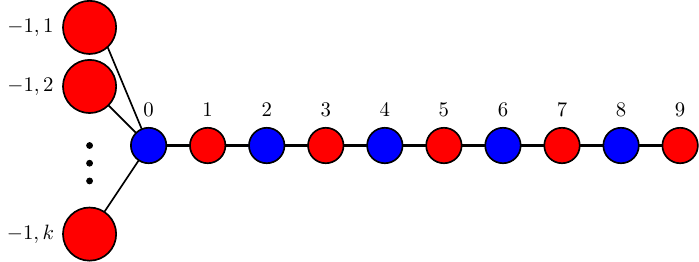}
\end{center}
\caption{Modified $\N-$graph $\cG_k$.}
\label{fig:MNNIGraph}
\end{figure}

Henceforth we call $\cG_k$ the modified $\N$-graph. Label by $\{(-1,1),\ldots,(-1,k)\}$ the $k$ vertices of valence 1 adjacent to $0$ the only vertex of valence $n+1$. Any positive $\lambda$-harmonic function $\textbf{h}$ satisfies $\lambda \textbf{h}((-1,j)) = \textbf{h}(0)$ for all $j=1,\ldots,k$ and $\textbf{h}(1)=\textbf{h}(0)(\lambda-\frac{k}{\lambda})$. For every $n\geq 1$ we have the same linear recurrence relations as in Example~\ref{exa:NCayleyGraph}, namely $\textbf{h}(n)=\lambda \textbf{h}(n-1) - \textbf{h}(n-2)$. Its characteristic polynomial is $x^2-\lambda x+1$. There are hence two cases to consider:
\begin{enumerate}
\item $\lambda=2$. Here the general solution of the linear recurrence relation that a $2$-harmonic function must satisfy is of the form $\textbf{h}(n)=A+Bn$. By evaluating in $n=0$ and $n=1$, and substituing the value for $\textbf{h}(1)$ calculated above we obtain that $\textbf{h}(n)=\textbf{h}(0)+(\lambda \textbf{h}(0)-\frac{k}{\lambda})n$. The only possible values that yield a positive $2$-harmonic function are $k=1$ (which is just the case of the $\N$-graph treated before) and $k=2$.
\begin{tcbexercise}{}{}
Show that for $k=2$ the corresponding Hooper-Thurston-Veech surface is tiled by infinitely many unit squares but is not isometric to the infinite staircase.
\end{tcbexercise}
\item $\lambda>2$. In this case the general solution of the linear recurrence relation
that $\textbf{h}$ satifies must satisfy of the form
\begin{equation}
	\label{eq:GenSolRecurrenceRel}
\textbf{h}(n)=Ar_+^n+Br_-^n
\end{equation}
(where $r_+>r_-$ are the roots of the characteristic polynomial). Therefore:
$$
\begin{pmatrix}
A\\
B
\end{pmatrix}
=
\frac{\textbf{h}(0)}{r_+-r_-
}\begin{pmatrix}
-r_- & 1\\
r_+ & -1
\end{pmatrix}
\begin{pmatrix}
1\\
\lambda-\frac{k}{\lambda}
\end{pmatrix}
$$
Hence:
\begin{equation}
	\label{eq:GenSolRecurrenceRel1}
\textbf{h}(n)=\frac{\textbf{h}(0)}{r_+-r_-}\left((\lambda-\frac{k}{\lambda}-r_-)r_+^n-(r_+-\lambda+\frac{k}{\lambda})r_-^n\right)
\end{equation}
\begin{tcbexercise}{}{}
Show that for $k\geq 2$ fixed, the condition $(\lambda-\frac{k}{\lambda})r_+>1$ is sufficient to guarantee that we have a positive eigenvector for the adjacency operator. \textit{hint: the sequence $\mu_n:=\frac{r_+^n-r_-^n}{r_+^{n-1}-r_-^{n-1}}$, $n\geq 2$ is decreasing and $\lim_{n\to\infty}\mu_n=r_+$.}

For fixed $k$ find $\lambda$ such that $A=0$. Use~\eqref{eq:GenSolRecurrenceRel} to obtain $\textbf{h}(n)$ and use Exercise~\ref{exo:FiniteAreaHTV} to prove that in this case the corresponding Hooper-Thurston-Veech surface has finite area.
\end{tcbexercise}
\end{enumerate}

\end{tcbexample}

\section{Some Veech groups of finite area surfaces}
  \label{sec:VeechGroupsChamanaraArnouxYoccoz}

This Chapter ends with the computation of the Veech groups for baker surfaces $B_\frac{1}{q}$, $q\geq 2$, and the Bowman-Arnoux-Yoccoz surface, both introduced in Section~\ref{ssec:BakerBowmanArnouxYoccoz}.
These are some of the few non-trivial examples on which the full Veech groups of finite-area infinite-type surfaces are known.

\subsection{The Veech group of baker's surfaces $B_{\frac{1}{q}}$.}
\label{ssec:VeechGroupBakersSurface}

In this section we revisit baker's surfaces $B_\alpha$, with $\alpha>0$. These finite area infinite-type surfaces were introduced in Example~\ref{ssec:BakerBowmanArnouxYoccoz} from Chapter~\ref{chap:Introduction}. For parameters $\alpha=\frac{1}{q}$, $q\geq 2$ we show that $B_{\alpha}$ is a Hooper-Thurston-Veech surface and we describe its Veech group explicitly.

\begin{tcbproposition}{}{exa:TVBaker}
\label{exa:TVBaker}
For $\alpha=\frac{1}{q}$, $q\geq 2$, the baker's surface $B_{\frac{1}{q}}$ is affinely equivalent to a Hooper-Thurston-Veech surface of modulus $\frac{\sqrt{q}}{q+1}$. That is, there exists a matrix $N\in\SL(2,\R)$ such that\footnote{Here $N\cdot$ denote the natural action of $\GL(2,\R)$ on translation surfaces.} $N\cdot B_{\frac{1}{q}}$ is a Hooper-Thurston-Veech surface.
\end{tcbproposition}

\begin{figure}[!ht]
\begin{center}
\begin{minipage}{.45\textwidth}
\begin{center}\includegraphics[scale=0.8]{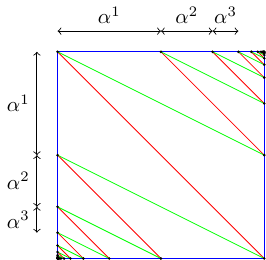}\end{center}
\subcaption{The baker's surface $B_\alpha$}
\end{minipage}%
\begin{minipage}{.45\textwidth}
\begin{center}\includegraphics[scale=0.8]{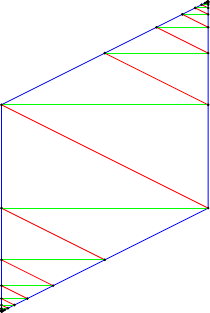}\end{center}
\subcaption{$A_1 B_\alpha$}
\end{minipage} \bigskip \\
\begin{minipage}{.45\textwidth}
\begin{center}\includegraphics[scale=0.4]{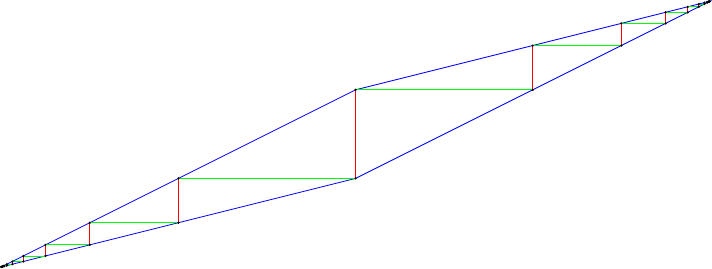}\end{center}
\subcaption{$A_2 A_1 B_\alpha$}
\end{minipage}
\begin{minipage}{.45\textwidth}
\begin{center}\includegraphics[scale=0.8]{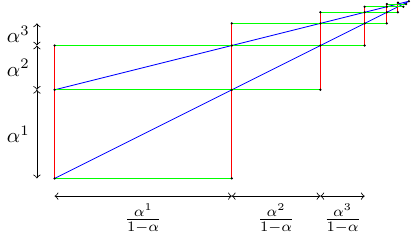}\end{center}
\subcaption{An infinite staircase obtained after cut and paste operations on $A_2 A_1 B_\alpha$}
\end{minipage}
\end{center}
\caption{The baker's surface as an infinite staircase: an example of the Hooper-Thurston-Veech construction.}
\label{fig:BakerAsTV}
\label{fig:BakerAsTVInfStair}
\end{figure}

\begin{proof}
In the top left part of Figure~\ref{fig:BakerAsTV}, we depict two families of parallel saddle connections that bound two transverse cylinder decompositions of baker's surface $B_\alpha$. The first step to see that $B_\alpha$ can be obtained by a Hooper-Thurston-Veech construction is to render cylinders in these two decompositions vertical and horizontal respectively. This is achieved by first applying an affine map with derivative $A_1=\begin{psmallmatrix}1 & 0\\ \alpha & 1\end{psmallmatrix}$ and then another with derivative $A_2=\begin{psmallmatrix}1 & \frac{1}{1-\alpha}\\ 0 & 1\end{psmallmatrix}$. The result is the image depicted at the bottom left part of figure \ref{fig:BakerAsTV}. By cutting and pasting we can present the obtained surface as an infinite (half) staircase going upwards, as shown in the bottom right part of figure \ref{fig:BakerAsTVInfStair}. Let us denote this (half) infinite staircase $B'_\alpha$. From the figure it is clear that the horizontal and vertical translation flow decompose $B'_\alpha$ into two infinite families of cylinders that we denote $H=\{H_0,H_1,\ldots\}$ and $V=\{V_0,V_1,\ldots\}$ respectively. Here, the enumeration of horizontal cylinders goes from bottom to top, and for vertical cylinders, it goes from left to right. Let us denote by $\mu(\_)=\frac{\text{height}(\_)}{\text{circumference}(\_)}$ the function associating to a cylinder its modulus. A direct calculation shows that
$$
\mu(H_0)=1-\alpha, \quad \mu(H_j)=\frac{\alpha(1-\alpha)}{1+\alpha} \ \text{for $j\geq1$} \quad \text{and}\quad \mu(V_k)=\frac{1}{(1+\alpha)(1-\alpha)},
$$
for $k\geq0$.
In other words, $A_2A_1B_{\alpha}$ is not a Thurston-Veech surface because the moduli of the horizontal and vertical cylinders are not equal. To fix this issue one can apply an affine map with hyperbolic derivative of the form $A_3(X)=\begin{psmallmatrix} X & 0 \\ 0 & X^{-1}
\end{psmallmatrix}$. This hyperbolic affine map changes moduli of horizontal and vertical cylinders by a factor of $X^{-2}$ and $X^2$ respectively. A direct calculation shows that
 if we pick $A_3(\sqrt[4]{\alpha(1-\alpha)^2})$, then $\mu(A_3H_j)=\mu (A_3V_k)=\frac{\sqrt{\alpha}}{1+\alpha}$ for all $j\geq 1$, $k\geq 0$ but $\mu(A_3H_0)=\frac{1}{\sqrt{\alpha}}$. To fix this last issue we can artificially subdivide the horizontal cylinder $H_0$ into $q$ horizontal cylinders $H_{0,i}$, $i=1,\ldots, n$ of equal height and modulus $\mu(H_{0,i})=\frac{1}{q\sqrt{\alpha}}$. Hence, $A_3A_2A_1B_{\alpha}$ after this subdivision of $A_3H_0$ is a Hooper-Thurston-Veech surface if and only if $\alpha=\frac{1}{q}$, $q\geq 2$.
\end{proof}

\begin{tcbexercise}{}{}
From Figure~\ref{fig:BakerAsTVInfStair} deduce that baker's surface
$B_\alpha$ for  $\alpha=\frac{1}{q}$ is a Hooper-Thurston-Veech surface
obtained by a construction on the modified $\N$-graph $\cG_{q+1}$ as explained
in Example \ref{exa:ModifiedNGraph}. Moreover, show that in this case the eigenvalue of the
adjacency operator is $\lambda=\frac{\sqrt{q}}{q+1}$.
\end{tcbexercise}

\begin{tcbexercise}{}{ratpargiveTV}
Prove that the surface $B_\alpha$ for $\alpha\in\mathbb{Q}$ is given by a Hooper-Thurston-Veech construction, up to application of an affine transformation and a subdivision of its horizontal and vertical cylinders.
\end{tcbexercise}

We now determine the Veech groups of baker's surfaces for parameters $\alpha=\frac{1}{q}$, $q\geq 2$.  Let us recall that from (\ref{eq:GroupGLambda}) that for every $\lambda\in\R$, $G_\lambda$ denotes the subgroup of $\SL(2,\R)$ generated by:
\[
h_\lambda = \begin{pmatrix}1 & \lambda \\ 0 & 1 \end{pmatrix}
\qquad
v_\lambda = \begin{pmatrix}1 & 0 \\ \lambda & 1 \end{pmatrix}.
\]

\begin{tcbtheorem}{}{BakerVeechGroup}
Let $q \geq 2$. The Veech group of baker's surface $B_{\frac{1}{q}}$ is conjugated in $\SL(2,\R)$ to the group generated by $-\rm Id$ and $G_\lambda$, where $\lambda=\frac{q+1}{\sqrt{q}}$.
\end{tcbtheorem}

\begin{tcbremark}{}{}
Theorem~\ref{thm:BakerVeechGroup} was originally stated
by Chamanara in~\cite{Chamanara}. However there is an issue with the proof Chamanara presents: he claimed that for any translation surface of finite area the translation flow on the invariant direction of a parabolic affine automorphism always defines a cylinder decomposition of the surface. The proof of this (Proposition 2 in~\cite{Chamanara}) assumes that a generalized version of Keane's theorem holds (Proposition 1 in
~\cite{Chamanara}). We will provide a counterexample to this assumption in
Section~\ref{sec:CounterexampleKeane}. Herrlich and Randecker have also computed the Veech group of $B_{\frac{1}{2}}$, see \cite{HerrlichRandecker}.
\end{tcbremark}
\begin{proof}[Proof of Theorem~\ref{thm:BakerVeechGroup}]
From Theorem~\ref{exa:TVBaker} we know that Baker's surface $B_{\frac{1}{q}}$ is, up to composition with an affine map, a Hooper-Thurston-Veech surface. Indeed, if we apply to $B_{\frac{1}{q}}$ the matrix
\begin{equation}
    \label{EQ:A-change-coordinates}
 A = \begin{pmatrix} \sqrt[4]{\alpha(1-\alpha)^2} & 0\\ 0 & \frac{1}{\sqrt[4]{\alpha(1-\alpha)^2}}\end{pmatrix}\begin{pmatrix}1 & \frac{1}{1-\alpha}\\ 0 & 1\end{pmatrix}
\begin{pmatrix}1 & 0\\ \alpha & 1\end{pmatrix}
\end{equation}
where $\alpha=\frac{1}{q}$ using the $\SL(2,\R)$-action on translation surfaces, the result is
a Hooper-Thurston-Veech surface $M(\mathcal{G}_{q+1},h)$, where $\mathcal{G}_{q+1}$ is a modified $\N$-graph and $\lambda=\frac{q+1}{\sqrt{q}}$. This surface can be depicted as a (half) staircase and we denote it $M_\frac{1}{q}$, see figure \ref{fig:BakerAsTV}. Given that $B_{\frac{1}{q}}$ and $M_{\frac{1}{q}}$ are affinely equivalent, their Veech groups are conjugated. More precisely $\Gamma(M_{\frac{1}{q}})=A\Gamma(B_{\frac{1}{q}})A^{-1}$.

Rotational components were introduced in Section~\ref{SSEC:RotationalComponents} in Chapter~\ref{ch:TopologyGeometry} in order to describe wild singularities of translation surfaces. Moreover, we started a discussion on rotational components of baker's surfaces in Example~\ref{exa:RotationalComponentsBakersSurface} in Chapter~\ref{ch:TopologyGeometry}.
In the following lemma we continue this discussion by describing the set of rotational components of $B_\frac{1}{q}$.

From the Hooper-Thurson-Veech construction we know that $G_\lambda$ is a subgroup of $\Gamma(M_{\frac{1}{q}})$ and $h'=\begin{psmallmatrix}q & 0\\ 0 & \frac{1}{q}\end{psmallmatrix}\in\Gamma(M_{\frac{1}{q}})$. Hence $G'_\lambda:=A^{-1}G_\lambda A<\Gamma(B_\frac{1}{q})$. Define $H$ as the quotient of the group generated by $\{G_\lambda',-Id\}$ modulo the subgroup generated by the hyperbolic element $h'$.
\begin{tcblemma}{}{RotCompChamanara}
The rotational components of $B_{\frac{1}{q}}$ are:
\begin{itemize}
\item two bi-infinite rotational\footnote{\ie the map $\widetilde{\rm dir}$ restricted to each rotational component defines a bijection with $\R$, see (\ref{eq:InjectiveLift}).}  components isometric to $\R$.
\item a countable collection of open rotational components\footnote{\ie the map $\widetilde{\rm dir}$ restricted to each rotational component $R_i$ defines a bijection with an open subinterval of $\R$, see (\ref{eq:InjectiveLift}).} of finite total angle given by the $H$-orbit of a rotational component of length $\frac{\pi}{2}$.
\end{itemize}
Moreover, for every direction $\theta\in\R/2\pi\Z$ there is at most one rotational component of finite total angle containing a linear approach parallel to $\theta$.
\end{tcblemma}
\begin{proof}

From Figure~\ref{F:cham} in Section~\ref{ssec:BakerBowmanArnouxYoccoz} we have that:
\begin{enumerate}
    \item There exist only two bi-infinite rotational components (in red and blue in the Figure). These are the only rotational components of infinite total angle.
    \item In $B_{\frac{1}{q}}$ there exist two open rotational components $R_0'$ and $-R_0'$ of total angle $\frac{\pi}{2}$: they are
    defined by segments approaching the vertices $\textbf{b}$ and $\textbf{d}$ in the Figure, respectively. Moreover, these are the only bounded rotational components of total angle $\frac{\pi}{2}$ whose co-slopes\footnote{Recall that the co-slope of a vector $(x,y)$ is defined as $x/y$.} parametrize the ray $(0,\infty)$. This means, roughly speaking, that when turning around the singularity in any of these two rotational components we go from the horizontal to the vertical direction (swaping a total angle of $\frac{\pi}{2}$).
\end{enumerate}

The rotational components $R_0'$ and $-R_0'$ are fixed by a (hyperbolic) affine automorphism with derivative $h'=\begin{psmallmatrix}q & 0\\0 & \frac{1}{q}\end{psmallmatrix}$.
Let us denote by $R_0$ the finite angle rotational component in $M_\frac{1}{q}$ corresponding to $R_0'$ via the affine map whose derivative is the matrix A with $\alpha=\frac{1}{q}$, see~(\ref{EQ:A-change-coordinates}).

If set of co-slopes of unit vectors defining linear approaches in  $R_0$ is the interval $I_0=(r_2,r_1)$, where $r_2=\frac{1}{\sqrt{q}}$ and $r_1=\sqrt{q}$. Indeed, this follows from the fact that $A\begin{psmallmatrix}1 \\ 0  \end{psmallmatrix}=\begin{psmallmatrix}\beta\frac{q}{q-1} \\ \beta^{-1}\frac{1}{q}  \end{psmallmatrix}$ and $A\begin{psmallmatrix}0 \\ 1  \end{psmallmatrix}=\begin{psmallmatrix}\beta\frac{q}{q-1} \\ \beta^{-1}  \end{psmallmatrix}$, where $\beta = \sqrt[4]{\frac{1}{q}(1-\frac{1}{q})^2}$.
If we consider the action of matrices by linear transformations on co-slopes, a direct computation shows that the open interval $I_0$ is the $h=Ah'A^{-1}$ orbit of the interval $]\frac{2}{\lambda},\frac{2}{\lambda}+\frac{\lambda}{2}]$, where $\lambda=\frac{q+1}{\sqrt{q}}$. Now, Theorem~\ref{thm:GLambdaLimitSet} from Appendix~\ref{app:FuchsianGroups} we know that the $G_\lambda$-orbit of the interval $I_0$ is equal to $\R\backslash\Lambda(G_\lambda)$, where $\Lambda(G_\lambda)$ denotes the limit set of $G_\lambda$. In particular, the complement in $\R$ of the $G_\lambda$-orbit of $I_0$ has empty interior.

Now we use again the flat geometry of baker's surface illustrated in Figure~\ref{F:cham}, Section~\ref{ssec:BakerBowmanArnouxYoccoz}. There are no rotational components of total angle equal to zero. Indeed, if $\gamma$ is a linear approach in $B_{\frac{1}{q}}$ defining a rotational component of angle zero then the co-slope $\theta$ of vector defined by its derivative lies in $\Lambda(G)$. Then there exists an element of $H$ sending $\theta$ to the vertical or horizontal direction. From the Figure it is clear that any linear approach in the vertical or horizontal direction lies in a bi-infinite rotational component, which is a contradiction since rotational components of total angle zero cannot be sent by the affine group into rotational component of positive of infinite angle. Therefore the set of rotational components with finite total angle of $B_{\frac{1}{q}}$ is given by the $H$-orbit of $R_0'$ and $-R_0'$.
Indeed, if $R$ is such a rotational component then the interval of co-slopes of its image via $A$ in $M_{\frac{1}{q}}$ is an open interval. Hence it intersects $\R\backslash\Lambda(G_\lambda)$. Therefore there is an element of $H$ sending $R$ to a rotational component intersecting either  $R_0'$ or $-R_0'$ an the claim follows. In particular, we have that for every $\theta\in\R/2\pi\Z$  there is at most one rotational component of finite total angle containing a linear approach parallel to $\theta$.
\end{proof}

We now finish the proof of Theorem \ref{thm:BakerVeechGroup}. Let $g$ be an element of the Veech group of $M
_\frac{1}{q}$ different from $-Id$. By the preceding lemma, $G_\lambda$ acts transitively on the set of  rotational components of finite total angle and hence there exists $g'\in G_\lambda$ such that $g'g$ fixes the interval $I_0$, which is the set of directions defining the rotational component of finite total angle fixed by the hyperbolic map $A\begin{psmallmatrix} q & 0 \\ 0 & \frac{1}{q}
\end{psmallmatrix}A^{-1}
$. Therefore $A^{-1}g'gA$ must be a hyperbolic matrix that fixes the horizontal and vertical direction. Since in $B_{\frac{1}{q}}$ all saddle connections in the vertical and horizontal direction have lengths in the set $\{q^n\}_{n\in\Z}$ we conclude that $A^{-1}g'gA$ must be a power of the hyperbolic matrix $\begin{psmallmatrix}q & 0\\ 0& \frac{1}{q}\end{psmallmatrix}$ and hence $g\in G_\lambda$.
\end{proof}


\subsection{The Veech group of the infinite Arnoux-Yoccoz surface}

The infinite Arnoux-Yoccoz surface $(X_\infty,\omega_\infty)$ was introduced in Section~\ref{ssec:BakerBowmanArnouxYoccoz}. Using Figure~\ref{infinite-arnoux-yoccoz} we explained there that $\Aff(X_\infty,\omega_\infty)$ contains an orientation reversing involution $\rho_\infty$ and a map $\psi_\infty$ with hyperbolic derivative $\begin{psmallmatrix}2 & 0\\0 & \frac{1}{2}\end{psmallmatrix}$. The rest of this section is dedicated to sketch the proof of the following:

\begin{tcbtheorem}{\cite{Bowman-Arnoux-Yoccoz}}{InfiniteAYVeechGroup}
The affine group of the infinite Arnoux-Yoccoz surface is generated by the involution $\rho_\infty$ and the hyperbolic map $\psi_\infty$.
\end{tcbtheorem}

\begin{proof}
We present a small variation of the original proof. However, we do use the following two facts about the infinite Arnoux-Yoccoz surface (see Lemmas 4.7 and 4.10 in~\cite{Bowman-Arnoux-Yoccoz}):
\begin{enumerate}
\item Every singular vertical trajectory of the translation flow in $\Aff(X_\infty,\omega_\infty)$ is a saddle connection and saddle connections are dense in this surface.
\item Every element of $\Aff(X_\infty,\omega_\infty)$ leaves the vertical direction invariant.
\end{enumerate}
As detailed in Section 4.1 in~\cite{Bowman-Arnoux-Yoccoz} the first return map of the vertical translation flow to the horizontal boundary of the polygon depicted in Figure~\ref{infinite-arnoux-yoccoz} defines an exchange map $f_\infty$ on an infinite number of intervals; (1) above is obtained by Bowman by a detailed analysis of the dynamics of $f_\infty$. Infinite interval exchange transformations are discussed in detail in the next Chapter. On the other hand, (2) above follows from the fact that any vertical saddle connection in $(X_\infty,\omega_\infty)$ lies inside a bi-infinite rotational component, whereas for any $\theta\neq\pm\frac{\pi}{2}$ there are singular trajectories contained in bounded rotational components.

Let now $\phi\in\Aff(X_\infty,\omega_\infty)$. By (2) above $\phi$ fixes the horizontal direction. Given that all vertical saddle connections have lengths in $(2^n)_{n\in\Z}$ we can suppose, modulo
composition by a power of $\rho_\infty$ and $\psi_\infty$, that $D\phi$ is the identity in the vertical direction. We claim that $D\phi$ cannot be a translation automorphism or a parabolic element, thus $\phi$ lies in the group generated by $\rho_\infty$ and $\psi_\infty$. Indeed, if $D\phi$ was parabolic, then up to passing to an appropiate power, $\phi$ should fix the saddle connection in Figure~\ref{infinite-arnoux-yoccoz} labeled by $A_1$. Then $\phi$ must fix all vertical saddle connections, which by (2) above are dense, implying that $\phi$ is the identity. On the other hand, from Figure~\ref{infinite-arnoux-yoccoz} we can see that vertical saddle connections of a given length occur in pairs and there is no translation automorphism taking one to the other.
\end{proof}

In general there is no recipe to calculate the Veech group a given translation surface, not even for Hooper-Thurston-Veech surfaces. So far, no lattice example has been found among Veech groups of infinite-type finite area translation surfaces. We conclude this Chapter with the following open problem.

\begin{tcbquestion}{}{}
Does there exist a finite area infinite-type translation surface whose Veech group is a lattice in $\SL(2,\R)$?
\end{tcbquestion}

\chapter{Translation flows}
\label{chap:TranslationFlowsAndIET}

With this last chapter begins the study of the dynamics of the translation flows. It is a prelude to our second volume~\cite{DHV2} dedicated to dynamics.

The rich idea to study flows via first return maps goes back to Poincar\'e. We see in Section~\ref{sec:IETAndTranslationFlows} how this construction can be performed for infinite-type translation surfaces and give rise to infinite interval exchange transformations, or infinite IET for short. We discuss these first return maps for baker's surfaces $B_\alpha$. We next discuss more general first return maps adapted to the study of infinite coverings, Hooper-Thurston-Veech surfaces or infinite step billiards.

Section~\ref{sec:Jungle} is dedicated to Arnoux-Ornstein-Weiss construction (Theorem~\ref{thm:VershikArnouxOrnsteinWeiss}): every aperiodic measure-preserving transformation is measurably isomorphic to an interval exchange transformation on $(0,1)$. In plain words: anything can happen in the world of infinite IETs. The main ingredient of this construction is the so-called cutting and stacking construction.

For a finite type surface, a classical result of Keane~\cite{Keane75} shows that the absence of vertical saddle connection implies minimality of the (vertical) translation flow. We show in Section~\ref{sec:CounterexampleKeane} that this does not hold for infinite-type surfaces. More precisely, we construct a non-minimal infinite interval exchange transformation without connections.

We finish this chapter with Section~\ref{sec:entropy} discussing the (metric) entropy of interval exchange transformations. The main result is Theorem~\ref{thm:EntropyInfiniteIET} following the ideas of F.~Blume~\cite{Blume2012} which provides an upper bound for the entropy of an infinite IET in terms of the lengths of the subintervals it permutes. We apply the latter to show that some of our favorite examples have zero entropy : baker and Bowman-Arnoux-Yoccoz surfaces from Section~\ref{ssec:BakerBowmanArnouxYoccoz} and the infinite step billiards from Example~\ref{exa:InfiniteStepSurface}. However, as we prove in Example~\ref{exa:PositiveEntropyIET} following the work of P.~Shields~\cite{Shields73}, there exist infinite IETs with positive entropy.

We assume that the reader has some basic knowledge of ergodic theory
and topological dynamics. If not, most notions discussed in this chapter can be found in~\cite{KrerleyViana16, Dajani-Kalle2021, EinsiedlerWard11} or references therein.

\section{First return maps of translation flows}
\label{sec:IETAndTranslationFlows}

We assume that the reader is familiar with the basic aspects of interval
exchange transformations on a finite number of intervals. We refer to
the surveys of~\cite{Yoccoz10} and~\cite{Viana06}. The aim of this
section is to study first-return maps in the setting of infinite translation
surfaces.

\subsection{Interval exchange transformations}
\label{ssec:FirstReturnMaps}
Let us recall that given a translation surface $M$ there is a well-defined
(partial) translation flow $F^t: M \dashrightarrow M$ (see Section~\ref{ssec:Structures}).
Given a regular point $x \in M$ let
$$
c^+_x=\sup\{t\in \left(0,+\infty\right]\hspace{1mm}:\hspace{1mm} \text{$F^t(x)$ is defined}\}.
$$
The parametrized curve $\alpha: t \in [0,c^+_x) \mapsto F^t(x)$ is
called the \emph{forward orbit} of $x$ or the \emph{forward ray}. Similarly, one
can define a backward orbit by looking at
$$
c^-_x=\inf\{t\in \left[-\infty,0\right)\hspace{1mm}:\hspace{1mm} \text{$F^t(x)$ is defined}\}.
$$
The map $\alpha: t \in (c^-_x, c^+_x) \mapsto F^t(x)$ is called a
\emph{two sided orbit} or the \emph{leaf through $x$}. Often, the maps $\alpha$
are identified with their embedded image in the surface. Recall from
Definition~\ref{def:saddle_connection} that a two sided orbit is called:
\begin{itemize}
\item a saddle connection if both $c^-_x$ and $c^+_x$ are finite,
\item a (forward) separatrix if $c^-_x$ is finite and $c^+_x$ is infinite,
\item a (backward) separatrix if $c^-_x$ is infinite and $c^+_x$ is finite,
\item a regular leaf (or geodesic) if both $c^-_x$ and $c^+_x$ are infinite.
\end{itemize}

In order to define a first-return map, we need a choice of a curve transverse
to the translation flow.

\begin{tcbdefinition}{}{TransverseCurve}
Let $M$ be a translation surface and $\gamma: I \to M$ a piecewise smooth embedded curve
where $I \subset \R$ is a bounded open interval.
We say that $\gamma$ is \emphdef[transverse (curve)]{transverse} to the vertical translation
flow if $\gamma$ has a continuous extension $\overline{\gamma}:\overline{I} \to \widehat{M}$
and its derivative $\gamma'$, in translation charts, is never vertical.
\end{tcbdefinition}

Remark that if $M=(X,\omega)$ and $\gamma$ is transverse to the vertical translation flow, then the one-form $\Re(\omega)$ restricted to $\gamma(I)$ is non-degenerate.

Henceforth we sometimes abuse notation and write $\gamma$ to denote the image of $\gamma: I \to M$. We will mostly encounter curves $\gamma$ that are union of straight-line segments in
$M$. A union of segments is transverse if and only if none of its components is vertical.

We now define the first-return map on a transverse curve $\gamma$. The
construction is simple: start with a point $x \in \gamma$, follow
the vertical translation flow in $M$ until it reaches again $\gamma$ at $y$ and
then define the first-return map as $x \mapsto y$. Let us introduce a formal definition.

\begin{tcbdefinition}{}{}
Let $M$ be a translation surface and $\gamma: I \to M$ be a smooth
curve transverse to the flow.
The \emphdef[first-return time]{(forward) first-return time} of the translation flow $F^t: M \to M$
to $\gamma$ is the map $r: M \rightarrow (0,+\infty]$ defined by
\[
r(x) = r^+(x) := \inf \{t \in (0,c^+_x):\ F^t(x) \in \gamma\}.
\]
Here $r(x)=+\infty$ if the forward orbit of $x$ never meets $\gamma$ again
(even though $c^+_x$ might be finite). The \emph{backward first-return time}
$r^-(x)$ is the return time of the backward flow
\[
r^-(x)  := \sup \{t \in (c^-_x,0):\ F^t(x) \in \gamma\}.
\]
Analogously, if the backward orbit of $x$ never meets $\gamma$ again we set $r^-(x)=-\infty$ (even though $c^-_x$ might be finite).

The \emphdef{first-return map} (or \emphdef{Poincar\'e map}) of the
vertical translation flow on $\gamma$ is the (partial) map
$f: \gamma \dashrightarrow \gamma$ defined as
$f(x) = F^{r(x)}(x)$ if $r(x) < \infty$ and undefined otherwise.
We will denote by $D \subset \gamma$ the set of points where $r$ is finite and $R =
f(D)$ (that is $D$ and $R$ are respectively the domain and range of $f$).
\end{tcbdefinition}
\pagebreak
\begin{tcbexample}{}{TorusWithTwoMarkedPoints}
In what follows we describe a classical example: the first return map to a transverse interval in a compact translation surface.
\begin{figure}[H]
\begin{minipage}{0.4\textwidth}\begin{center}\includegraphics{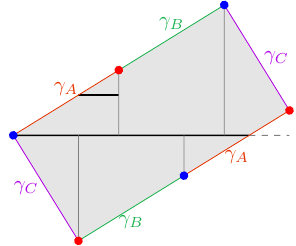}\end{center}
\subcaption{A torus with two marked points drawn with a horizontal transversal segment.}
\label{fig:TorusWithSegment}
\end{minipage}
\hspace{0.1\textwidth}
\begin{minipage}{0.4\textwidth}\begin{center}\includegraphics{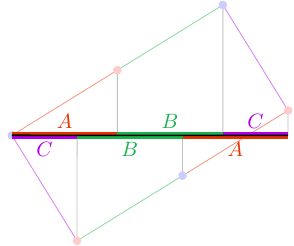}\end{center}
\subcaption{The first-return map of the vertical translation flow from Figure~\ref{fig:TorusWithSegment}.}
\label{fig:TorusFirstReturnMap}
\end{minipage}
\caption{First return maps in a finite type surface.}
\label{fig:TwoTorusIET}
\end{figure}
Let $M$ be the torus on the left hand side
of Figure~\ref{fig:staircase1}. Let us fix an angle $\theta \in (0, \pi/2)$. A first return map
for $M$ is shown in Figure~\ref{fig:TwoTorusIET} and consists of an interval
exchange on three subintervals labeled $A$, $B$ and $C$. This first return map $f_\theta: D \to R$ is represented on
Figure~\ref{fig:TorusFirstReturnMap} following a standard convention. The domain is represented on top and each component
of $D$ is labelled with a letter. Each subinterval $f(D_i)$ in the codomain is labeled with the label of $D_i$.
The map $f_\theta$ can be encoded by a permutation the $\left(\begin{array}{l}A\ B\ C\\C\ B\
A\end{array}\right)$ and its length data
\[
\lambda_A = \sin(\theta),\
\lambda_B = \sin(\theta) \text{ and }
\lambda_C = \cos(\theta).
\]
\end{tcbexample}

We now introduce interval exchange transformations that are piecewise translation of the interval\footnote{Interval exchange transformations can, in general, be defined using a large class of (piecewise) homeomorphisms of the interval.}. We
then prove that the first return map of translation flows as considered above are naturally conjugated to
an interval exchange transformation. The definitions are meant to capture first return map to intervals of the translation flow on an infinite-type translation surface.
\begin{tcbdefinition}{}{iet}
Let $I$ be a bounded open interval.
A \emphdef{partial interval exchange transformation on $I$} is an homeomorphism $f: D \rightarrow R$ between open subsets
$D \subset I$ and $R \subset I$ for which $f$ is a translation on each connected component of $D$. If furthermore
$D$ has full measure in $I$ then we say that $f$ is an \emphdef{interval exchange transformation}.

The \emphdef[singularities (of iet)]{forward and backward singularities} of an
interval exchange transformation $f: D \rightarrow R$ are respectively
$\overline{I \setminus D}$ and $\overline{I \setminus R}$ where the closure
is taken in $\overline{I}$.
\end{tcbdefinition}
Let us notice that a partial interval exchange transformation $f$ preserves the Lebesgue measure: given a measurable set
$A \subset R$ we have $\Leb(f^{-1} (A)) = \Leb(A)$. Since $f$ is bimeasurable we also have for $A \subset D$
$\Leb(f(A)) = \Leb(A)$).

A technical tool will be the following parametrization of curves.
\begin{tcbdefinition}{}{CurveNaturalParametrization}
Let $M=(X,\omega)$ be a translation surface and $\gamma: I \to M$ a transverse curve to the vertical translation
flow. We say that $\gamma$ has a
\emphdef[natural parametrization (transverse curve)]{natural parametrization} if
$\gamma^* \Re(\omega) = du$, where $u$ is the standard coordinate on $I \subset \R$. Equivalently,
the parameter is natural if in translation charts we have $\gamma'(u) = (1,y(u))$ for some function $y:I \to \R$.
\end{tcbdefinition}
\begin{tcbremark}{}{}
The natural parametrization of $\gamma: I \to M$ is a priori different from
the unit-speed parametrization, \emph{i.e.} when $\|\gamma'(u)\| = 1$.
\end{tcbremark}

\begin{tcbexercise}{}{}
Let $M$ be a translation surface and $\gamma: I \to M$ be a transverse curve.
Show that there is a unique (possibly orientation reversing) diffeomorphism $\phi: I \to I$
such that $\gamma \circ \phi$ has natural parametrization.
\end{tcbexercise}

We now present the most general version on how to construct (partial) interval exchange
transformations from translation flows. In many situations however, it is more convenient to use Corollary~\ref{cor:TranslationFlowFirstReturn}.
\begin{tcbdefinition}{}{}
Let $\gamma$ be a curve transverse to the vertical flow in a translation surface $M$.
Its left (or right) endpoint $x$ is called an
\emphdef[admissible endpoint (of
transverse curve)]{admissible endpoint}
if $x$ is either a singularity of $M$ or its forward-return time $r^+(x)$ is infinite or
its backward-return time $r^-(x)$ is infinite.
\end{tcbdefinition}
\begin{tcbtheorem}{}{TranslationFlowFirstReturn}
Let $M$ be a translation surface, $\gamma: I \to M$ be a curve transverse to the vertical flow with its natural parametrization, $r^+$ and $r^-$ the forward and backward first-return times
to $\gamma(I)$ and $T: \gamma(I) \dashrightarrow \gamma(I)$ the corresponding first-return
map. Let $D$ and $R$ be the subsets of $I$ such that $\gamma(D)$ (resp. $\gamma(R)$)
is the domain (resp. range) of $T$. Let $P$ be the subset of $M$ consisting of the
endpoints of $\gamma$ that are not admissible (it consists of either zero, one or two points).
Then
\begin{enumerate}
\item the map $f$ obtained as the restriction of $\gamma^{-1} \circ T \circ \gamma: I \dashrightarrow I$
on the domain $D \setminus \{F^{r^-(x)}(x): x \in P\}$ and range $R \setminus \{F^{r^+(x)}(x): x \in P\}$
is a partial interval exchange transformation.
\item For any connected component $C$ of the domain of $f$ and any pair of points $x,y \in C$ we have
\begin{equation}
\label{eq:ReturnTimeFlowBox}
r^+(x) + \int_{f(x)}^{f(y)} \Im(\omega) = r^+(y) + \int_x^y \Im(\omega)
\end{equation}
where the paths of integration between $x$ and $y$ and $f(x)$ and $f(y)$ is taken on the transverse
curve $\gamma(I)$.
\end{enumerate}
\end{tcbtheorem}

\begin{tcbcorollary}{}{TranslationFlowFirstReturn}
Let $M$ be a translation surface and $\gamma: I \to M$ be a transverse straight-line segment in
$M$ whose both ends are admissible. Let $T: \gamma(I) \dashrightarrow \gamma(I)$ be the first return map of
the translation flow on $M$ to $\gamma$ and $r^+$ the first-return time. Then
\begin{enumerate}
\item the map $f = \gamma^{-1} \circ T \circ \gamma: I \dashrightarrow I$ is a partial interval exchange map,
\item on any connected component of the domain of $f$, the function $r^+ \circ \gamma$ is constant.
\end{enumerate}
\end{tcbcorollary}

\begin{proof}[Proof of Theorem~\ref{thm:TranslationFlowFirstReturn}]
Let $u_0 \in D$, that is, $x_0 := \gamma(u_0)$ is such that $r^+(x_0) < \infty$. We will show that there is
a neighborhood of $x_0$ such that $f$ is a translation around $x_0$ and that~\eqref{eq:ReturnTimeFlowBox}
holds in this neighborhood.
For a flat point $x \in M$, we denote by $b(x)$ the radius of the largest disk centered at $x$ that contains
only flat points and no points of  $\partial \gamma$. The map $x \to b(x)$ is positive and continuous
on $M \backslash (\Sing(M) \cup \partial \gamma)$. Indeed, one can see $b(x)$
as the distance to $\Sing(M) \cup \partial \gamma$ in the metric completion $\widehat{M}$. Now by
compactness, $\epsilon := \inf_{0 \leq t \leq r(x_0)} b(F^t(x_0)) > 0$. Let
$J$ be a neighborhood of $u_0$ in $I$ so that $\gamma(J)$ is contained in the ball of radius
$\epsilon$ around $x_0$.

We claim that for any point $y$ in $\gamma(J)$, the image by the first-return
map $T(y)$ is well-defined. By assumption $T(x_0)$ is well-defined. The map $r \circ \gamma$
has a discontinuity (or cease to be defined) at points that hit a singularity or whose orbit runs
into the endpoints of $\gamma$ before returning to $\gamma$. Because the latter points were removed
from the domain, this proves that $r \circ \gamma$ is continuous on $J$.

Now consider the (topological) rectangle $R := \{(u,t) \in J \times \R: 0 \leq t \leq r(\gamma(u)))\}$.
We have
an immersion
\[
\psi:
\begin{array}{lll}
\widetilde{U} & \to & M \\
(u,t) & \mapsto & F^t(\gamma(u))
\end{array}.
\]
that satisfies $\psi^* \omega = du + i dt$ and that is embedded in restriction
to $t < r(\gamma(u))$. See Figure~\ref{fig:FlowBox}.
\begin{figure}[!ht]
\begin{center}\includegraphics{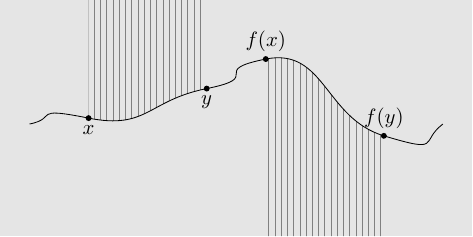}\end{center}
\caption{A flow box in a translation surface (drawn in a neighborhood of $\gamma$).}
\label{fig:FlowBox}
\end{figure}

Now notice that inside the flow box $r(x) = \int_x^{f(x)} \Im(\omega)$. Hence~\eqref{eq:ReturnTimeFlowBox} corresponds
to the fact that the integral of $\Im(\omega)$ around the boundary of the flow box vanishes because of Stoke's theorem.
\end{proof}

\begin{proof}[of Corollary~\ref{cor:TranslationFlowFirstReturn}]
By the first item of Theorem~\ref{thm:TranslationFlowFirstReturn}, the first return map is a partial
interval exchange transformation.

To prove the second part, one uses the second item of Theorem~\ref{thm:TranslationFlowFirstReturn} together
with the fact that since $\gamma$ is a straight-line segment, the integral
$\int_x^y \Im(\omega)$ is proportional to $\int_x^y \Re(\omega)$.
\end{proof}
\begin{tcbexample}{}{BakerIET}

We now study some first-return maps of the translation flows in the baker's surfaces $B_\alpha$, where $\alpha$ is a parameter in $(0,1)$. These surfaces were introduced in Section~\ref{ssec:BakerBowmanArnouxYoccoz} (see also Example~\ref{exa:BakerIsAMonster} and Section~\ref{ssec:VeechGroupBakersSurface}). We recall that the surface $B_\alpha$ is built from a square with side length $\frac{\alpha}{1-\alpha}$.


\begin{figure}[H]
\begin{minipage}{0.4\textwidth}
\begin{center}\includegraphics{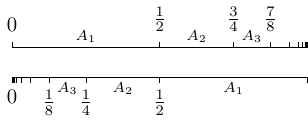}\end{center}
\subcaption{The top interval represents the domain of $f$ and the bottom interval its range.}
\end{minipage}
\hspace{0.1\textwidth}
\begin{minipage}{0.4\textwidth}
\begin{center}\includegraphics{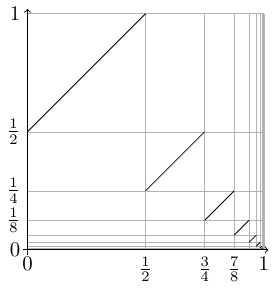}\end{center}
\subcaption{The graph of $f$.}
\end{minipage}
\caption{Two views of the interval exchange transformation $f: (0,1) \dashrightarrow (0,1)$  induced by the the first-return map of the vertical translation flow of the baker's surface $B_{1/2}$ on the transversal formed by the horizontal saddle connections labeled by $\{A_i\}_{i=1}^\infty$ in Figure~\ref{F:cham}.}
\label{fig:VerticalBakerIET}
\end{figure}

Now let us consider an angle $\theta \in (0,\pi/2)$. We will construct the first-return map of the translation flow in $B_\alpha$ in direction $\theta$ on the curve $\gamma$ that corresponds to the union of the vertical and the horizontal side\footnote{This curve is not a piecewise smooth embedded curve because it is formed by countably many disjoint open intervals. Therefore it does not satisfy the hypothesis of Definition~\ref{def:TransverseCurve}. However, for any $\theta\in(0,\pi/2)$ the curve $\gamma$ is homotopic along the leaves of the translation flow $F_\theta^t$ to $\delta\setminus (x_n)$, where $\delta$ is the diagonal in Figure~\ref{fig:BakersSurface} joining the points $a$ and $c$, and $(x_n)$ is a closed infinite sequence of points in $\delta$. The diagonal $\delta$ satisfies the hypothesis of Definition~\ref{def:TransverseCurve} and Corollary~\ref{cor:TranslationFlowFirstReturn}, hence it can be taken to compute a first return map. Nonetheless it is simpler for calculations to use the curve $\gamma$. This amounts to the same first return map because first return maps are invariant if we perform a homotopy along leaves of the translation flow.}. Let $\omega = dx + \sqrt{-1} dy$ be the one form associated to the translation structure of $B_\alpha$. We have
\[
\Re(e^{\sqrt{-1}(\pi/2-\theta)} \omega) = \sin(\theta) dx + \cos(\theta) dy.
\]
Hence the natural parametrization (see Definition~\ref{def:CurveNaturalParametrization}) of the curve $\gamma$ for the translation flow in direction
$\theta$ "multiplies" distances by $\sin(\theta)$ on the horizontal side and $\cos(\theta)$ on the vertical side. In this example, the components of the domain and range of $f_{\alpha,\theta}$ corresponds exactly to the side pairing of Figure~\ref{fig:BakersSurface}. If we label the components according to the numbering on that figure they appear as
$A_1\ A_2\ A_3\ \cdots\ B_3\ B_2\ B_1$ in the domain
and $B_1\ B_2\ B_3\ \cdots\ A_3\ A_2\ A_1$ in the range.
Now, the lengths of each of these components are given by the formulas
\[
\lambda_{A_i} = \sin(\theta)\ \alpha^i
\quad \text{and} \quad
\lambda_{B_i} = \cos(\theta)\ \alpha^i.
\]
We illustrate this in Figure~\ref{fig:GeneralBakerIET}.
Note that as $\theta$ tends to $\pi/2$ the interval exchange transformation $f_{\alpha,\theta}$ degenerates to the example of Figure~\ref{fig:VerticalBakerIET}.

\begin{figure}[H]
\begin{center}\includegraphics[scale=0.94]{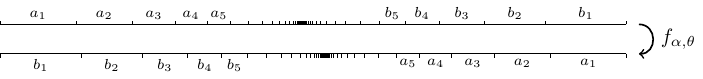}\end{center}
\caption{An example of interval exchange transformation $f_{\alpha,\theta}$ obtained as first return map of the translation flow in baker's surfaces $B_\alpha$.}
\label{fig:GeneralBakerIET}
\end{figure}

\begin{tcbquestion}{}{}
Let $\alpha \in (0,1)$. Does the interval exchange $f_{\alpha,\theta}$ admit a dense orbit for almost every $\theta$?
\end{tcbquestion}
There is no single $\alpha$ for which the answer to this question is known! However, as we will discuss in~\cite{DHV2}, it follows from the work of P. Hooper~\cite{Hooper-infinite_Thurston_Veech} that for $\alpha \in \Q$ one can find parameters $\theta$ for which most orbits are dense (more precisely, the corresponding interval exchange transformation is \emph{ergodic}). The corresponding subset of parameters $\theta$ has positive Hausdorff measure but zero Lebesgue measure.

\end{tcbexample}

We refer the reader to~\cite{HooperRafiRandecker2020} and ~\cite{Bruin-Lukina-2101} and reference within for more examples of infinite interval exchange transformations whose dynamics have been studied recently.

\subsubsection{Suspension}
For an interval exchange transformation $f$ with finitely many components, the
Veech zippered rectangle~\cite{Veech82} and Masur's polygon~\cite{Masur82}
provide concrete constructions of a translation surface so that the first return map of the translation flow on some direction is precisely $T$ (see also~\cite[Section 4 \emph{Suspension of i.e.m.: the zippered
rectangle construction}]{Yoccoz10}). The following exercise shows that one can
always construct a surface from an interval exchange transformation.
\begin{tcbexercise}{}{}
Let $f: D \to R$ be an interval exchange transformation on $I=(0,1)$. Let
\[
S := \Big( \big((0,1) \times [0,1]\big)\ \backslash\ \big((I \backslash D) \times \{0\} \cup (I \backslash R \times \{1\})\big)\Big)\ /\ (x,1) \sim (f(x),0).
\]
Prove that $S$ is a translation surface and that the first-return map of the vertical flow on $(0,1) \times \{1/2\}$
is equal to $f$.
\end{tcbexercise}

\subsubsection{Encoding}
Finite interval exchange transformations are completely described by:
\begin{compactitem}
\item a finite \emph{alphabet} $\cA$ used to label the (finitely many) connected components of the domain $D$ and the range $R$,
\item a \emph{combinatorial data} which is a pair of bijections $\pi_t, \pi_b$ from $\cA$ to $\{1,2,\ldots,d\}$
that describes how the components of $D$ are permuted,
\item a \emph{length data} $\lambda \in \R_+^\cA$ which are the lengths of the components of $D$.
\end{compactitem}
For more details on the finite case the reader is invited to
read~\cite[Section 3: \emph{Interval exchange maps : basic definitions}]{Yoccoz10}. We now give an extension of this
description to handle infinite interval exchange transformations.
\begin{tcbdefinition}{}{ConstructiveIET}
Let $\cA$ be an at most countable set and $(<_{top},<_{bot})$ two total orders on $\cA$. Let $\lambda=(\lambda_i)_{i\in \cA}$ a summable family of positive real numbers and let us define for $i \in \cA$
$$
\alpha_i^{top}:=\sum_{j:\ j<_{top} i}\lambda_j
\quad \text{and} \quad
\alpha_i^{bot}:=\sum_{j:\ j<_{bot} i}\lambda_j
$$
$$
I_i^{top}:=(\alpha_i^{top},\alpha_i^{top}+\lambda_i)
\quad \text{and} \quad
I_i^{bot}:=(\alpha_i^{bot},\alpha_i^{bot}+\lambda_i).
$$
We define the interval
exchange transformation $f=f_{(A,<_{top},<_{bot},\lambda)}$ as the map that sends
$I_i^{top}$ to $I_i^{bot}$ by translation. In other words, for $x\in
I_i^{top}$, $f(x):=x-\alpha_i^{top}+\alpha_i^{bot}$.
\end{tcbdefinition}
\begin{tcbexercise}{}{FiniteAndInfiniteIET}
\begin{itemize}
\item Show that definitions \ref{def:iet} and \ref{def:ConstructiveIET} are equivalent. \emph{Hint}: in definition \ref{def:iet} the family of open intervals $(D_i)_i$ is at most countable.
\item Show that in Definition~\ref{def:ConstructiveIET}, if $\cA$ is finite of cardinality $d$ then a total
order $<_{top}$ on $\cA$ is equivalent to a bijection $\pi_{top}: \cA \to \{1,2,\ldots,d\}$.
\item Describe a set $\cA$, two orders $<_{top}$, $<_{bot}$ and a length data $\lambda \in \cR^\cA$ for the interval exchange transformations $f_{\alpha,\theta}$ of Example~\ref{exa:BakerIET}.
\end{itemize}
\end{tcbexercise}

\subsubsection{Iteration of iet and periodic points}

In the following paragraphs we discuss how to define the domains of the iterates of a partial interval exchange transformation and we recall some basic properties of their periodic orbits.

\begin{tcbdefinition}{}{}
Given a partial interval exchange transformation $f: D \to R$ on an interval $I$ and $n\in\Z$ let $D_n := D_n(T)$ be the set
of points where $f^n$ is defined. Formally $D_0 := I$, $D_1 := D$, $D_{-1} := R$ and for
$n \geq 1$
\[
D_{n+1} := f^{-1} (D_n \cap D_{-1})
\qquad \text{and} \qquad
D_{-n-1} := f (D_{-n} \cap D_1).
\]
We also set
$\displaystyle D_{-\infty} := \bigcap_{n \geq 0} D_{-n}$ and $D_{+\infty} := \bigcap_{n \geq 0} D_n$
\end{tcbdefinition}
Note that by definition $f^n (D_n) = D_{-n}$.

\begin{tcbdefinition}{}{IETConnection}
Let $f: D \to R$ be a partial interval exchange transformation on an interval $I$.
A \emphdef[connection (of an iet)]{connection} for $f$ is a triple $(m, x, y)$ where $m$ is a
non-negative integer, $x \in D_m \backslash D_{m+1}$, $y \in D_{-m} \backslash D_{-(m+1)}$ and $f^m x = y$.
\end{tcbdefinition}

\begin{tcbexercise}{}{}
Let $M$ be a translation surface and $\gamma:I \to M$ a transverse curve. Assume that
the first return map $T: \gamma(I) \dashrightarrow \gamma(I)$ is defined almost everywhere. Let $f: I \dashrightarrow I$
be map conjugated to $T$ by $\gamma$ as in Theorem~\ref{thm:TranslationFlowFirstReturn}. Let $x$ be a point in $D(f)$ such that $\gamma(x)$ belongs to a vertical saddle connection in $M$. Show that the orbit of $x$ under $f$ is finite in both the future and the past.
\end{tcbexercise}

As in the case of finite
interval exchange transformations, periodic orbits always come in families. In
other words, as the following Lemma states, being periodic this is an "open property".
\begin{tcblemma}{}{PeriodicOrbitsIET}
Let $f:D \to R$ be a partial interval exchange transformation on an interval
$I$. If $x \in D$ is a periodic orbit of period $n$ then there exists an open
interval $U$ containing $x$ so that $f^n$ is the identity on $U$.
The maximal such interval $U$ is delimited by points that belong to a connection,
that is both extremities belong to $I \setminus (D_{+\infty} \cup D_{-\infty})$.
\end{tcblemma}
The proof is left to the reader (\emph{hint}: just use the fact that $f^n = f \circ f \circ \cdots \circ f$ is an interval exchange
transformation).
Let us also remark that this fact was already proven for translation surfaces, see Exercise~\ref{exo:PeriodicTrajectoriesAndCylinders}.

We finish this section with a characterization of partial interval exchange transformations that do not have periodic orbits. The result is a direct consequence of Lemma~\ref{lem:PeriodicOrbitsIET}
and is used in Section~\ref{sec:entropy} when we study entropy.
\begin{tcbcorollary}{}{PeriodicOrbitsIETObstruingShrinking}
Let $f:D \to R$ be a partial interval exchange transformation on an interval $I$ and let $f^n: D_n \to D_{-n}$ be its $n$-th iterate.
Let $\cD^{(n)}$ be the set of connected components of $D_n$. Then the following are equivalent
\begin{compactitem}
\item $f$ does not have periodic orbits,
\item $\displaystyle \lim_{n \to \infty} \max_{J \in \cD^{(n)}} |J| = 0$.
\end{compactitem}
\end{tcbcorollary}

\subsection{More general first return maps}
\label{ssec:Coding}

Infinite interval exchange transformations are useful for studying the translation flow but they have an important limitation, namely, the domain of definition has to be a finite length interval. For this reason it is necessary to consider first return maps with infinite-lenght domains. In what follows we discuss several contexts on which such first return maps arise. The dynamics of the translation flow on these examples will be discussed thoroughly in the second Volume of this text.

\subsubsection{Skew-products of iet as return maps in $G$-coverings}
\label{sec:SkewProductsAndGCoverings}
In this section we discuss a construction of first return maps in
$G$-coverings of finite-type translation surfaces. The maps we obtain
are particular case of skew-products. Dynamical properties of general
skew-products in the context of infinite measure dynamical systems
have received a lot of attention since the 1970s due to the work of J.-P.~Conze and
K.~Schmidt.

Let us recall that $G$-coverings were introduced in Section~\ref{sec:CoveringSpaces}.
We start with the definition of skew-product restricted to the context of
finite interval exchange transformations
\begin{tcbdefinition}{}{def:SkewProductIET}
Let $f: D \to R$ be a finite interval exchange transformation on an interval $I$, that is
$D$ has finitely many connected components. Let $(G, \cdot)$ be an at most countable
group and let $\phi: D \rightarrow G$ be a locally constant function. The map
$$
f_\phi :
\begin{matrix} D \times G &\to & R \times G \\
(x,g) & \to & (f(x), g \cdot \phi(x))
\end{matrix}
$$
is called a \emphdef{skew-product} with \emphdef[fiber (skew-product)]{fiber} $G$ and \emphdef[base (skew-product)]{base} $I$.

The associated $f$-\emphdef[cocycle (skew-product)]{cocycle} is the sequence of maps $\{\phi_n:D_n \to G\}_{n \in \Z}$ defined by
\[
\phi_n(x) =
\left\{ \begin{array}{ll}
\phi(x) \cdot \phi(f(x)) \cdots \phi(f^{n-1}(x)) & \text{if $n \geq 0$} \\
(\phi(f^{n}(x)) \cdot \phi(f^{n+1}(x)) \cdots \phi(f^{-1}(x)))^{-1} & \text{if $n < 0$}
\end{array} \right.
\]
and satisfies
\[
f_\phi^n(x, g) = (f^n(x), \phi_n(x)).
\]
\end{tcbdefinition}
Note that a skew-product $f_\phi$ preserves the product of the Lebesgue measure on $I$ and the
counting measure on the group $G$. If the group $G$ is infinite, this measure has infinite
total mass.

Our aim is now to show a correspondence between translation flows
of $G$-coverings of finite-type translation surfaces and skew-products of finite
interval exchange transformations in a similar way as in Theorem~\ref{thm:TranslationFlowFirstReturn}.

\begin{tcbexample}{}{InfiniteStaircaseAndSkewProducts}
We start with the construction of a skew-product for the infinite staircase
that was introduced in Section~\ref{ssec:InfiniteStaircase} and further
studied in Example~\ref{exa:Zcoverings}.
\begin{figure}[H]
\begin{center}\includegraphics[scale=0.9]{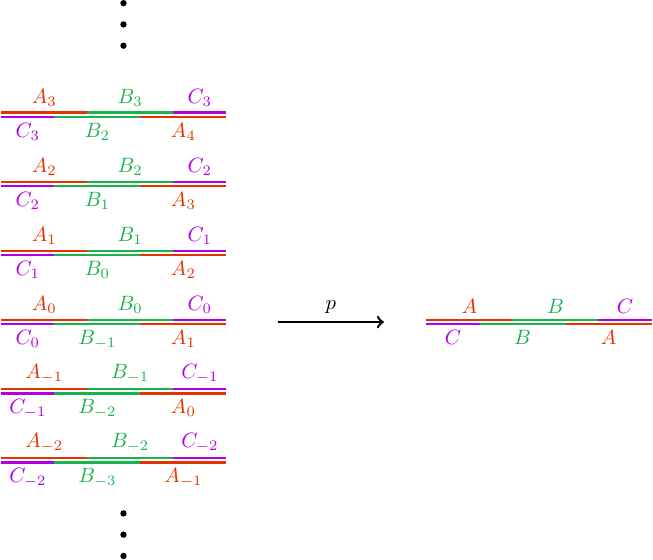}\end{center}
\caption{The skew-product $(f_\theta)_\phi$ for the infinite staircase. The map
$p: I \times \Z \to I$ represented on the picture is the forgetful map. It commutes
with the dynamics: $p \circ (f_\theta)_\phi = f_\theta \circ p$.}
\label{fig:FirstReturnInfiniteStaircase}
\end{figure}

Let $\widetilde{M}$ be the infinite staircase. Recall that $\widetilde{M}$ is the
$\Z$-covering of a torus $M$ with two marked points defined by a cocycle $c \in
H_1(M, \Sigma; \Z)$. Using the notations from Figure~\ref{fig:staircase1} the
cocycle $c$ is $B-A$.

In Example~\ref{exa:TorusWithTwoMarkedPoints} we considered the torus $M$
defining the base of this cover. We built a first return $f_\theta: I \to I$
of the translation flow in direction $\theta$ in $M$. See in particular
Figure~\ref{fig:TorusWithSegment}.

Let $\widetilde{I}$ be the preimage of the segment $I$ by the covering map
$p: \widetilde{M} \to M$. The preimage $\widetilde{I}$ is the union of countably
many segments. The first return map on $\widetilde{I}$ of the translation flow in $\widetilde{M}$
in direction $\theta$ is a skew-product with base $f_\theta: I \to I$, fiber $\Z$ and
cocycle
\[
\phi: \begin{array}{lll}
I & \to & \Z \\
x & \mapsto & \left\{
\begin{array}{lll}
1 & \text{if $x \in I_B$} \\
-1 & \text{if $x \in I_A$} \\
0 & \text{if $x \in I_C$}
\end{array} \right.
\end{array}
\]
The values of the function $f_\theta$ are the algebraic intersection of the cycles
$\gamma_A$, $\gamma_B$ and $\gamma_C$ dual to the sides of the fundamental
polygon, see Figure~\ref{fig:staircase1}, with the cycle $c$ defining the
covering $p: \widetilde{M} \to M$.

Here the group $\Z$ is additive and the $f_\theta$-cocycles form a
\emphdef{Birkhoff sum}
\[
\phi_n(x) = \phi(x) + \phi(f_\theta(x)) + \ldots + \phi(f_\theta^{n-1}(x)).
\]

Figure~\ref{fig:FirstReturnInfiniteStaircase} represents the skew-product $(f_\theta)_\phi$.
\begin{tcbexercise}{}{}
Given $f_\theta$ as considered above, let $g_\theta$ be the map induced
on the union $J$ of the subintervals $I_A$ and $I_B$. That is $g_\theta(x) = f_\theta^{r(x)}(x)$
where $r(x) = \min \{m \geq 1: f_\theta^m(x) \in I_A \cup I_B\}$ is the first return
time to $I_A \cup I_B$.
\begin{enumerate}
\item Show that $g_\theta$ is a rotation,
\item Let $\psi: J \to \Z$ defined by $\psi(x) = \sum_{i=0}^{r(x)-1} \phi(x) = \phi_{r(x)}(x)$.
Show that the values of $\psi$ are constant on $I_A$ and $I_B$ and does not
depend on $\theta$.
\item Show that the first return map on the preimage $\widetilde{J}$ of $J$ by the covering
map $p: \widetilde{M} \to M$ identifies with the skew product with base $g_\theta$ and
cocycle $\psi$.
\end{enumerate}
\end{tcbexercise}
The alternative first return map provided in the above exercise is a particular case of
skew-products over rotations as considered in the work of J.~Aaronson, H.~Nakada,
O.~Sarig and R.~Solomyak~\cite{AaronsonNakadaSarigSolomyak02}.
\end{tcbexample}

We now turn to the general construction of first return maps in $G$-coverings.
Let $M$ be a compact translation surface and $I$ a horizontal segment in $M$.
Let $f: I \to I$ be the first return map to $I$ of the translation flow in $M$.
By Theorem~\ref{thm:TranslationFlowFirstReturn}, $f$ is (conjugated to) a finite interval
exchange transformation.

Let $x_0 \in I$ a reference point. For each point $x \in I$ we associate an element of
$\pi_1(M \setminus \Sigma; x_0)$ as follows. Consider the closed curve $\gamma(x)$ made
of the concatenation of
\begin{itemize}
\item the horizontal segment from $x_0$ to $x$,
\item the vertical segment $\{F^t(x)\}_{x \in [0, r(x)]}$ from $x$ to $f(x) = F^{r(x)}(x)$,
\item the horizontal segment from $f(x)$ to $x_0$.
\end{itemize}
The class of $\gamma(x)$ in $\pi_1(M \setminus \Sigma; x_0)$ only depends on the subinterval
of the interval exchange transformation $f$ the point $x$ belongs to.

Let $p: \widetilde{M} \to M$ be a $G$-covering of the compact translation surface $M$
given by the surjective cocycle $c: \pi_1(M \backslash \Sigma, x_0) \to G$. We construct
a map $\phi: I \to G$ by setting $\phi(x) := c(\gamma(x))$.

\begin{tcbtheorem}{}{}
With the above notations, the first return map of the translation flow
$\widetilde{F}^t$ in $\widetilde{M}$ on the preimage $\widetilde{I}$ of
the segment $I$ by the covering map $p: \widetilde{M} \to M$ identifies
with the skew-product with basis $f: I \to I$ and cocycle $\phi: I \to G$.
\end{tcbtheorem}

\begin{proof}
We make an arbitrary choice of a base point in the preimage $p^{-1}(x_0)$ so that we can identify $p^{-1}(x_0)$ with $G$.
It induces an identification of $\widetilde{I}$ with $I \times G$. Given a point $(x, g) \in I \times G$, its image under the translation flow $\widetilde{F}^t$ on $\widetilde{M}$ is of the form $(f(x), g')$.
By construction $(f(x), g') = (f(x), g) \cdot \gamma(x)$ where $\cdot$ represents the right action of $\pi_1(M \setminus \Sigma; x_0)$ on $I \times G$. Formally the action only exists in the fiber $p^{-1}(x_0)$ but since $I$ is simply connected we extend it to a locally constant action on $\widetilde{I}$. By definition of the $G$-covering associated to $c$, this action is given by $(f(x), g) = (f(x), g c(\gamma(x)))$.
\end{proof}

\subsubsection{M\'alaga maps}
\label{ssec:ReturnMapMalaga}
The M\'alaga maps $T_{\underline{\alpha}}: [0,1) \times \Z \to [0,1) \times \Z$ were introduced in Section~\ref{ssec:IntroRandomModels}. The aim of this section is to provide a proof of Lemma~\ref{lem:MalagaReturn} that identifies $T_{\underline{\alpha}}$ with a first return map in a generalized staircase. This construction presents some similarity with the skew-product constructions of Section~\ref{sec:SkewProductsAndGCoverings}.

Recall that in Section~\ref{ssec:IntroRandomModels} we introduced the
generalized staircase $S_{\underline{h}}$ where $\underline{h} \in (0,\infty)^\Z$.
The surface $S_{\underline{h}}$ was constructed from a particular gluing of the
rectangles $R_n = [0,2] \times [0,h_n]$.

We consider the disjoint union $\widetilde{J}$ of bottom sides of the
rectangles and we identify the bottom of the $n$-th rectangle $R_n$ to $J_n :=
[0,2) \times \{n\}$ (seen as a circle). We map each $J_n$ to $[0,1) \times
\{n\}$ by $x \mapsto \frac{x - 1}{2} \mod 1$. We obtain a continuous map $p:
\widetilde{J} \to \left [0,1\right)^\Z$ by setting
\[
p(x, n) = \left( \frac{x-1}{2} \mod 1, n \right).
\]

\begin{tcbexercise}{}{MalagaReturn}
Prove Lemma~\ref{lem:MalagaReturn} using the above identification of $\widetilde{J}$ with $[0,1)^\Z$.
\end{tcbexercise}

\subsubsection{Hooper-Thurston-Veech surfaces}
Let us recall from Section~\ref{sec:HooperThurstonVeechConstruction}
that a Hooper-Thurston-Veech surface $M$ is a translation surface
that admits horizontal and vertical cylinder
decompositions $H=\{H_i\}_{i\in I}$ and $V=\{V_j\}_{j\in J}$, where all
cylinders involved have the same moduli, see
Definition~\ref{def:HooperThurstonVeechSurface}. Moreover, the intersection of
these cylinders decomposes $M$ into a family of rectangles $\{R_e\}_{e \in E}$.

\begin{figure}[!ht]
\begin{center}%
\begin{minipage}{0.25\textwidth}\includegraphics[scale=0.8]{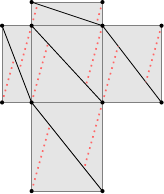}
\subcaption{A local picture of the first return map in a Hooper-Thurston-Veech surface.}
\label{fig:HTVFirstReturnDiag}\end{minipage}%
\hspace{0.09\textwidth}%
\begin{minipage}{0.6\textwidth}\includegraphics[scale=1.05]{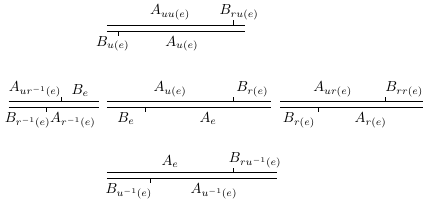}
\subcaption{The first return map seen on a union of the intervals appearing in Figure~\ref{fig:HTVFirstReturnDiag}.}
\end{minipage}%
\end{center}
\caption{The construction of the first return map in a Hooper-Thurston-Veech construction.}
\label{fig:HTVFirstReturn}
\end{figure}

We consider in each rectangle $R_e$ of the surface $M$ the diagonal $\gamma_e$
joining the top left and bottom right corners. We denote $\Gamma$ the union of
these diagonals. For an angle $\theta \in (0,\pi/2)$, the first return map of
the translation flow in direction $\theta$ in $M$ on $\Gamma$ is represented in
Figure~\ref{fig:HTVFirstReturn}. In this picture  $r: E \to E$ and
$u: E \to E$ denote the two bijections that indicate which rectangle
is glued to the right and up side of $R_e$ respectively.

\begin{tcbexercise}{}{}
Let $M$ be a Hooper-Thurston-Veech surface defined by a pair of multicurves $\alpha=\{\alpha_i\}_{i\in I}
$, $\beta=\{\beta_j\}_{j\in J}$ and $\textbf{h}$ a positive $\lambda$-harmonic function as in Theorem~\ref{thm:HooperThurstonVeechConstruction}. Let $\Gamma=\cup_{e\in E(\cG(\alpha\cup\beta))}\gamma_e$ be the set of diagonals defined above. Show that:
\begin{enumerate}
\item For every edge $e\in E(\cG(\alpha\cup\beta))$ between vertices $i\in I$ and $j\in J$ the
length $|\gamma_e|$ is bounded above by $(1+\lambda^2)^\frac{1}{2}\textbf{h}(\alpha_i)$.
\item For every $i\in I$ and $j\in J$ we have that $\textbf{h}(\alpha_i)\leq \min_{p_\alpha(e)=\alpha_i}|\gamma_e|$ and $\textbf{h}(\beta_j)\leq \min_{p_\beta(e)=\beta_j}|\gamma_e|$. Here $p_\alpha$ and $p_\beta$ are the natural projections sending an edge in $\cG(\alpha\cup\beta)$ to the vertices that define it, see~(\ref{eq:HTVProjections}).
\end{enumerate}
Conclude using both items above that $\Gamma$ is of finite length if and only if $\textbf{h}$ belongs to $\ell^1(V(\cG(\alpha\cup\beta)))$. Recall from Exercise~\ref{exo:FiniteAreaHTV} that $\textbf{h}$ belongs to $\ell^2(V(\cG(\alpha\cup\beta)))$ if and only if $M$ has finite area.
\end{tcbexercise}

\begin{tcbquestion}{}{L2NotL1}
Does there exist an infinite graph $\cG$ of finite degree having a positive $\lambda$-harmonic function belonging to $\ell^2(V(\cG))\setminus\ell^1(V(\cG))$? If so, describe an explicit example.
\end{tcbquestion}

\subsubsection{Step billiards}
In this section we study the translation flows in step surfaces from
Example~\ref{exa:InfiniteStepSurface}. The maps we obtain are very similar to
the M\'alaga maps defined above.

\begin{figure}[!ht]
\begin{center}\includegraphics{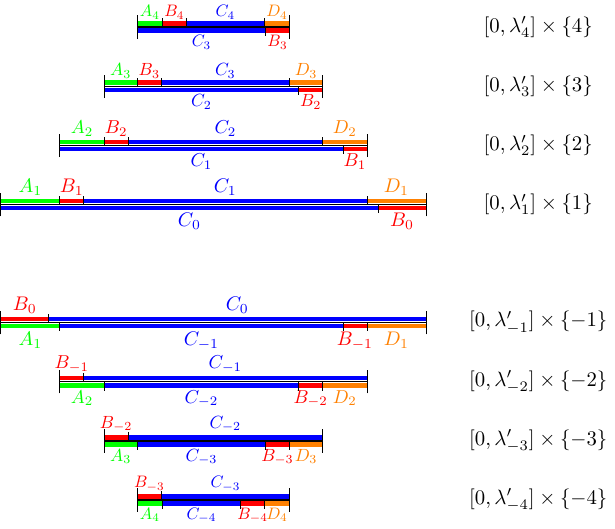}\end{center}
\caption{The first return map of the step surface of Figure~\ref{fig:StepSurface}. For each level $[0, \lambda'_i] \times \{i\}$ the domain and codomain are respectively represented by colored rectangle on top and bottom.}
\label{fig:StepSurfaceReturnMap}
\end{figure}

We first introduce a generalized interval exchange transformation $f_{\underline{h}, \underline{\alpha}}$
that depends on a decreasing sequence of real numbers $\underline{h} = (h_n)_{n \geq 1}$ and
real numbers $\underline{\alpha} = (\alpha_n)_{n \geq 0}$ that satisfy $\alpha_n \in [0, h_{n+1}/2]$.
The subdomains will be identified with the letters $B_i$, $C_i$ for $i\in \Z$ and $A_i$, $D_i$ for $i \geq 1$.
Let
\begin{align}
\label{eq:IETGeneralizedMalaga}
\lambda_{A_i} &= \lambda_{D_i} = h_i - h_{i+1}, \quad n \geq 1 \\
\quad
\lambda_{B_i} &= \lambda_{B_{-i}} = \alpha_n, \quad n \geq 0 \\
\quad
\lambda_{C_i} &= \lambda_{C_{-i}} = 2 h_{i+1} - \alpha_n, \quad n \geq 0
\end{align}
For $i \in \Z \setminus \{0\}$ we set
\begin{equation}
\lambda'_i = \left\{\begin{array}{ll}
   \lambda_{B_{i+1}} + \lambda_{C_{i+1}} & \text{if $i < 0$} \\
   \lambda_{B_{i-1}} + \lambda_{C_{i-1}} & \text{if $i > 0$}
\end{array} \right. .
\end{equation}
They also satisfy
\[
\lambda'_i = \left\{\begin{array}{ll}
   \lambda_{A_{-i}} + \lambda_{B_i} + \lambda_{C_i} + \lambda_{D_{-i}} & \text{if $i < 0$} \\
   \lambda_{A_i} + \lambda_{B_i} + \lambda_{C_i} + \lambda_{D_i} & \text{if $i > 0$} \\
\end{array} \right. .
\]
The domain of $f_{\underline{h}, \underline{\alpha}}$ is the infinite union
\[
\bigcup_{i \in \Z \setminus \{0\}} [0, \lambda'_i] \times \{i\}.
\]
The way $f_{\underline{h}, \underline{\alpha}}$ is defined can be seen on Figure~\ref{fig:StepSurfaceReturnMap}. More
precisely, each level $[0,\lambda'_i] \times \{i\}$ is represented with a partition of its domain on top and a
partition of its range at the bottom.

\begin{tcblemma}{}{}
Let $M$ be the step surface with step heights $\underline{h} = (h_n)_{n \geq 1}$.
Then for any direction $\theta \in (0, \pi/2)$, the first return
map on the transversal drawn on the left hand side of
Figure~\ref{fig:StepSurfaceReturnMap} is canonically identified with
the generalized inteval exchange $f_{\underline{h}, \underline{\alpha}}$
where $\underline{\alpha} = (\alpha_n)_{n \geq 0}$ is defined as
$\alpha_n = \cot(\theta) \mod 2h_{n+1}$.
\end{tcblemma}

\section{The jungle of infinite interval exchange transformations}
\label{sec:Jungle}
The main point of this section is to warn the reader of the following fact: the dynamical properties of infinite interval exchange transformations (and thus translation flows) are too rich and diverse to be studied as a whole. More precisely, any aperiodic\footnote{Through this text all dynamical systems are considered on
standard Borel spaces $(X,\mathcal{B},\mu)$, that is, we suppose that $X$ is
Polish space and $\mathcal{B}$ is the corresponding Borel $\sigma$-algebra. A measure-preserving transformation $T$ on a measure
space
$(X,\mathcal{B},\mu)$ is called \emphdef{aperiodic} if $\mu(\{x\in X :
T^k(x)=x\})=0$ for every integer $k\geq 0$.} prescribed dynamical behaviour can be realized by an interval exchange transformation. This is formally stated in following result of P.~Arnoux, D.~Ornstein and B.~Weiss.


\begin{tcbtheorem}{Arnoux-Ornstein-Weiss~\cite{ArnouxOrnsteinWeiss}}{VershikArnouxOrnsteinWeiss}
Every invertible aperiodic measure-preserving transformation $T:X \rightarrow X$ is (measurably) isomorphic to an interval exchange transformation $f:D \to R$ on $I=[0,1]$. Moreover, $f$ can be chosen so that the only point in $D$ and in $R$ to which the extremities of intervals accumulate is $1$.
\end{tcbtheorem}
 The main ideas for the proof of this result can already be found in the work of Vershik~\cite{Vershik-MarkovCompacta}. One of the main ingredients used in its proof is the cutting and stacking construction. Most of this section is dedicated to explain this construction, its relation to Bratteli-Vershik diagrams and how this can be used to construct an interval exchange transformations. These, together with the
Kakutani-Rokhlin Lemma (see Lemma~\ref{lem:Rokhlin}), give to a motivated reader enough elements to understand the proof as presented in~\cite{ArnouxOrnsteinWeiss}.


\begin{tcbremark}{}{}
There are analogs of Theorem~\ref{thm:VershikArnouxOrnsteinWeiss} in the category of topological dynamical systems. That is, instead of $T: X \to X$ being an invertible measurable dynamical systems, we can consider an homeomorphism of the Cantor set. In this settings, using  the construction of R.~H.~Herman, I.~.F.~Putman and C.~F.~Skau~\cite{HermanPutnamSkau1992} one can build an interval exchange transformation $f: D \to R$ whose points with biinfinite orbits are dense. Moreover, one can construct a continuous surjective map $\pi: X \to [0,1]$ such that $\pi \circ T = f \circ \pi$ and $\pi$ is $1$-to-$1$ on points of $D$ with biinfinite orbits under $f$. Note that here $\pi$ can not be an homeomorphism as the Cantor set $X$ and $[0,1]$ are not homeomorphic and $f$ is not continuous. This is a main difference with Theorem~\ref{thm:VershikArnouxOrnsteinWeiss}.
\end{tcbremark}

\textbf{Towers and stacks}. A tower is a simple concept in the theory of dynamical systems.

\begin{tcbdefinition}{}{Tower}
Let $(X, \cB, \mu)$ be a probability space and $T: X \to X$ a measurable measure-preserving invertible system. A \emphdef{tower} for $(X,T)$ is a finite orbit $(Y, TY, T^2 Y, \ldots T^{n-1}Y)$ of a subset $Y \in \cB$ such that $Y$, $T Y$, $T^2 Y$, \ldots, $T^{n-1} Y$ are disjoint in measure. The numbers $\mu(Y)$ and $n$ are called the \emph{width} and the \emph{height} of the tower respectively. The subsets $T^i Y$ are called the \emph{levels}. The $0$-th level $Y$ and $(n-1)$-th level $T^{n-1}Y$ are called the \emph{base} and the \emph{top} respectively .
\end{tcbdefinition}

Despite its simplicity, towers are a very important tool in dynamical systems. Dynamical properties are often easier to understand in terms of towers and they allow to produce examples of dynamical systems with exotic properties using the cutting and stacking construction. Furthermore, as in the case of finite interval exchange transformations, towers become a powerful tool to prove generic dynamical properties.

\begin{tcbexample}{}{TowerIET}
Let $M$ be a translation surface of finite area for which the vertical direction has been fixed. Let $I\subset M$ be an open and relatively compact horizontal interval. We consider the first-return map $f: I \to I$ of the vertical translation flow and assume that its domain has full measure in $I$.
Let $R \subset S$ be an embedded rectangle with horizontal and vertical sides so that it does not contain any extremity of $I$ in its interior. Then $R \cap I$ is a finite union of intervals $J_0$, $J_1$, \ldots, $J_{n-1}$ that form a tower for the transformation $f$. See Figure~\ref{fig:RectangleAsTower}.
\begin{figure}[H]
\begin{center}\includegraphics{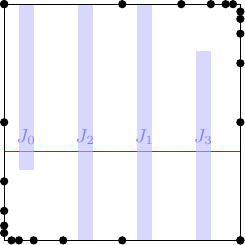}\end{center}
\caption{On this figure one can see the baker's surface from Section~\ref{ssec:BakerBowmanArnouxYoccoz} with a vertical blue rectangle that intersects four times the interval $I$ in red. It gives rise to a tower of height 4 for the first return map $f: I \to I$.}
\label{fig:RectangleAsTower}
\end{figure}
\end{tcbexample}

\begin{tcbexample}{}{TowerOdometer}
In the next paragraphs we give the elements to describe towers for the dyadic odometer. This is a dynamical system defined on $\Z_2$, the ring of 2-adic integers. For the purpose of this example we only use the additive structure of $\Z_2$, hence multiplication is not discussed. There are several ways to define this ring. For our purposes we consider the discrete topology on $\{0,1\}$ and the corresponding product topology on $\Z_2:=\{0,1\}^\N$. As a consequence $\Z_2$ is a totally disconnected and perfect, hence homeomorphic to the Cantor set. The family of subsets (a.k.a.  cylinders):
$$
C_m(k_0,\ldots,k_l):=\{(x_n)\in \Z_2\hspace{1mm}|\hspace{1mm} x_m=k_0, x_{m+1}=k_1,\ldots, x_{m+l}=k_l\}
$$
where $m\in\Z_{\geq 0}$ and $k_i\in\{0,1\}$ for $i=0,\ldots,k$ generate the Borel $\sigma$-algebra of $\Z_2$. By Caratheodory's extension theorem, $\mu_2(C_m(k_0,\ldots,k_l))=\frac{1}{2^{l+1}}$ admits a unique extension to a probability measure $\mu_2$ on $\Z_2$.
Addition on $\Z_2$ is defined as $(x_n)+(y_n)=(x_n+y_n+k_n \mod 2)$, where $k_0=0$ and $k_n$ is defined inductively by:
$$
k_n=
\begin{cases}
1 & \text{if $x_{n-1} = y_{n-1} = 1$} \\
0 & \text{otherwise.}
\end{cases}
$$
\begin{tcbremark}{}{}
The ring $\Z_2$ can also be defined as the inverse limit $\varprojlim_n \Z / 2^n \Z$, for $n\geq 0$ where each space in the corresponding inverse system has been given the discrete topology. The measure $\mu_2$ in this context is just the Haar measure on $\Z_2$, which is the projective limit of the (normalized) counting measure on $\Z / 2^n \Z$
\end{tcbremark}

Let $\underline{1}:=(1,0,0,\ldots)$. The dyadic odometer $T:\Z_2\to \Z_2$ is defined as $(x_n)=(x_n)+\underline{1}$. In particular, $f(1,\ldots,1,0,x_k,\ldots)=(0,\ldots,0,1,x_k,\ldots)$.
Let $\underline{1}:=(1,0,0,\ldots)$. The dyadic odometer $T:\Z_2\to \Z_2$ is defined as $T(x_n)=(x_n)+\underline{1}$. In particular, $T(1,\ldots,1,0,x_k,\ldots)=(0,\ldots,0,1,x_k,\ldots)$.

\begin{tcbexercise}{}{TowersForOdometer}
In this exercise we construct towers for the dyadic odometer.
\begin{itemize}
\item Show that the dyadic odometer is a homeomorphism that preserves the measure $\mu_2$.
\item For any non-negative integer $l\geq 0$ consider the cylinder $Y_l:=C_0(0,\ldots,0)$. Prove that $Y_l, T(Y_l),\ldots,T^{2^{l+1}-1}(Y_{l})$ are disjoint clopen subsets whose union is $\Z_2$. In particular, for any $l\geq 0$ the dyadic integers form a tower of $2^{l+1}$ levels of width $\frac{1}{2^{l+1}}$.
\item Prove that the dyadic odometer does not have any periodic orbit and is \emph{minimal}, that is, any orbit is dense in $\Z_2$.
\end{itemize}
\end{tcbexercise}
\end{tcbexample}

\textbf{The cutting and stacking construction}. We explain an inductive method which produces an interval exchange transformation on the interval $(0,1)$.

Up to isomorphism in measure, there is a unique tower of given width $w$ and height $h$.
Let $S_{w,h} := (0,w) \times \{0, 1, \ldots, h-1\}$ and consider the partial map
\[
T_{w,h}:
\begin{array}{lll}
S_{w,h} & \to & S_{w,h} \\
(x, i) & \mapsto & \left\{ \begin{array}{ll}
(x, i+1) & \text{if $i < h-1$} \\
\text{undefined} & \text{if $i = h-1$}
\end{array} \right.
\end{array}
\]
Here $S_{w,h}$ is endowed with the Lebesgue measure on each $(0,w)$ so that it has
total mass $w \times h$. A tower for the system $(X, T, \mu)$ can be alternatively
defined as a measure preserving injection $\iota: S_{w,h} \to X$ such that on
the domain of $T_{w,h}$ we have $\iota \circ T_{w,h} = T \circ \iota$. The
subsets $\iota( (0,w) \times \{i\})$ are the levels of the tower.
We call $(S_{w,h}, T_{w,h})$ the \emphdef{stack of with $w$ and height $h$}.

\begin{figure}
\begin{center}%
\includegraphics{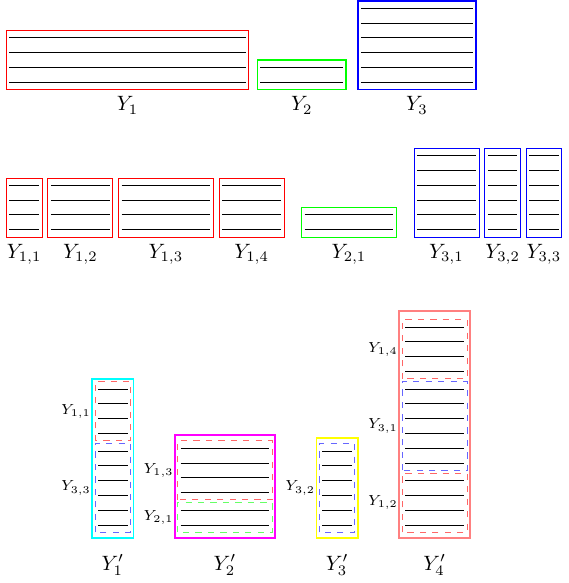}%
\end{center}
\caption{One cutting and stacking step. On the top part we start with 3 stacks that
form a certain step of the construction. In the middle, they are cut
in respectively $a_1=4$, $a_2=1$ and $a_3=3$ stacks. At the bottom, these $4+1+3$
stacks are stacked according to the words $W_1 = (3,3)\, (1,1)$, $W_2 = (2,1)\, (1,3)$, $W_3=(3,2)$ and $W_4=(1,2)\,(3,1)\,(4,1)$.}
\label{fig:CuttingAndStacking}
\end{figure}

We now consider two basic operations on a set of stacks: the cutting and the
stacking operations\footnote{Some authors, including~\cite{ArnouxOrnsteinWeiss}, consider an
extra operation in the cutting and stacking construction which consists in
adding intervals above and below stacks obtained after cutting. This operation
is often referred to as \emph{adding spacers}. The usage of these spacers is
convenient in certain construction such as the Chacon's example~\cite{Chacon69}.
In the present section, we use spacers in a rather hidden way in a single
place, namely Exercise~\ref{exo:WorkAroundSpacers}.}. It will be useful for the
reader to keep an eye on the left-hand side of
Figure~\ref{fig:CuttingAndStacking} as a reference.

Let $S_{w,h}$ be a stack and $w_1, w_2, \ldots, w_a$ positive real numbers so that
$w_1 + w_2 + \ldots + w_a = w$. We define the \emphdef{cutting operation} as the map
\[
\begin{array}{ccc}
S_{w_1, h} \sqcup S_{w_2, h} \sqcup \cdots \sqcup S_{w_a, h} & \to & S_{w, h} \\
 (x, i) \in S_{w_j, h} & \mapsto & (w_1 + \ldots + w_{j-1}  + x, i)
\end{array}
\]
In more visual terms the stack $S_{w, h}$ is cut vertically into pieces of
widths $w_1$, \ldots, $w_a$ (in this order) and the same height $n$.

Let $S_{w,h_1}$, $S_{w,h_2}$, \ldots, $S_{w,h_a}$ be a sequence of stacks
with the same width $w$ and let $h=\sum_{i=1}^ah_i$. We define the \emphdef{stacking operation} as the map
\[
\begin{array}{ccc}
S_{w, h} & \to & S_{w, h_1} \sqcup S_{w, h_2} \sqcup \cdots \sqcup S_{w, h_a} \\
(x, h_1 + \ldots + h_{j-1} + i) & \mapsto & (x, i) \in S_{w, h_j}
\end{array}
\]
In more visual terms, the stacks $S_{w, h_1}$, $S_{w, h_2}$, \ldots, $S_{w, h_a}$
are placed on top of each other to form a unique stack of height $h_1 + \ldots + h_a$.

The map associated to either a cutting or a stacking operation, is a bijection in measure.

We now define the cutting and stacking construction which is a sequence of stacks
built from the previous two operations. More precisely we build inductively a sequence
of finite collections of stacks $C^{(k)} = (S_{w_1^{(k)}, h_1^{(k)}}, \ldots, S_{w_{n_k}^{(k)}, h_{n_k}^{(k)}})$ together with an injective map
\[
\phi^{(k)}: S^{(k)}_{w_1^{(k)}, h_1^{(k)}} \sqcup S^{(k)}_{w_2^{(k)}, h_2^{(k)}} \sqcup \cdots \sqcup S^{(k)}_{w_{n_k}^{(k)}, h_{n_k}^{(k)}} \to (0,1)
\]
which is a local translation and a bijection in measure. More concretely, each level of each stack $S^{(k)}_{w_i^{(k)}, h_i^{(k)}}$ is identified to a subinterval of $(0,1)$.
The collection $C^{(k)}$ and the map $\phi^{(k)}$ induce a partial map $f^{(k)}$ of
$(0,1)$ by making the following diagram commute
\[
\begin{CD}
C^{(k)} @>T^{(k)}>> C^{(k)} \\
@V\phi^{(k)}VV @ VV\phi^{(k)}V \\
(0,1) @>>f^{(k)}> (0,1)
\end{CD}
\]
and where $T^{(k)} = T_{w_1^{(k)}, h_1^{(k)}} \sqcup \cdots \sqcup T_{w_{n_k}^{(k)}, h_{n_k}^{(k)}}$. The domain of $f^{(k)}$ is the image of $C^{(k)}$ under $\phi^{(k)}$ where we
remove the top of each stack. The main compatibility condition of this sequence
is: the domain of $f^{(k+1)}$ contains the domain of $f^{(k)}$ up to
finitely many points. On the intersection of their domains the maps $f^{(k)}$
and $f^{(k+1)}$ coincide.

The initial step is $C^{(0)} := (S_{1,1})$, a stack of height and width 1, where the
injection $\phi^{(0)} : S_{1,1} \to (0,1)$ is the trivial one and
the associated partial map $f^{(0)}: (0,1) \to (0,1)$ is the empty map.

Suppose that at the $k$-th step we have a collection
$$
C^{(k)} = (S_{w_1, n_1},
S_{w_2, n_2}, \ldots, S_{w_k, n_k})
$$
    of stacks and injections $\phi^{(k)}_i:
S_{w_i, n_i} \to (0,1)$ whose images form a partition of $(0,1)$ in measure. To
define a cutting and stacking operation in order to build $C^{(k+1)}$ from
$C^{(k)}$ we need
\begin{itemize}
\item positive integers $a_1$, $a_2$, \ldots, $a_{n_k}$
\item for each $i \in \{1, \ldots, k\}$ a probability vector $(p_{ij})_{j=1}^{a_i}$,
that is $p_{ij}$ are positive real numbers and $\sum_{j=1}^{a_i} p_{ij}=1$,
\item a list of non-empty words $(W_1, W_2, \ldots, W_{n_{k+1}})$ on the alphabet
$\{(i,j): i \in \{1, \ldots, n_k\}, j \in \{1, \ldots, a_i\}\}$ so that each letter
appears exactly once in the list of words and two letters $(i,j)$ and $(i',j')$ in the same
word $W_\ell$ satsify $w_i p_{ij} = w_{i'} p_{i'j'}$.
\end{itemize}
From the above data we explain how build $C^{(k+1)}$. We first cut each stack
$S_{w_i, n_i}$ according to the weights $\{p_i w_{ij}\}_j$. Then we stack them
according to the words $W_i$ in order
to build $n_{k+1}$ new stacks of widths $(w'_1, w'_2, \ldots, w'_{n_{k+1}})$,
where $w'_\ell$ is the quantity $w_i p_{ij}$ where $(i,j)$ is any letter
appearing in $W'_\ell$, see Figure~\ref{fig:CuttingAndStacking} for an example.

Finally $\phi^{(k+1)}: C^{(k+1)} \to (0,1)$ is obtained as
$\phi^{(k+1)} := \phi^{(k)} \circ \phi_{cut} \circ \phi_{stack}$
by composing $\phi^{(k)}$ with the maps defining the cutting and stacking
operations above. The new partial map $f^{(k+1)}$ obtained from $C^{(k+1)}$
and $\phi^{(k+1)}$ clearly satisfies the compatibility condition. It differs
from $\phi^{(k)}$ by the fact that it is not defined anymore at the cut
points (these are finitely many points) and some part of the tops of the
towers in $C^{(k)}$ are part of the domain of $T^{(k+1)}$ if in the
stacking part of the construction they do not appear on top.

\begin{tcblemma}{}{WellDefinedCutAndStack}
Let $C^{(k)} = (S_{w_1^{(k)}, h_1^{(k)}}, \ldots, S_{w^{(k)}_{n_k}, h^{(k)}_{n_k}})$ be a sequence of lists of stacks obtained by the cutting and stacking construction
from $C^{(0)} = (S_{1,1})$. Let $f^{(k)}: (0,1) \to (0,1)$ be the associated sequence of partial maps.
Let $f$ be the map obtained as the union of all the $f^{(k)}$, that is
$f(x) = f^{(k)}(x)$ if $x$ belongs to the domain of $k$ (the compatibility condition ensures that $f^{(k)}(x) = f^{(k')}(x)$ if $x$
belongs to both domains of $T^{(k)}$ and $T^{(k')}$).
If
\begin{equation}
\label{eq:ConditionDefiningLimit}
\lim_{k \to \infty} \sum_{i=1,2,\ldots,n_k} w^{(k)}_i = 0
\end{equation}
then $f$ is an aperiodic interval exchange transformation.
\end{tcblemma}

\begin{proof}
By construction, $f$ is a local translation: for each $x$ in the domain
of $f$ there exists an open set $U$ containing $x$ such that $f(y) = y + f(x) - x$
for $y$ in $U$.

The subset of $(0,1)$ where $f^{(k)}$ is not defined are the images under
$\phi^{(k)}$ of the top of the towers in $C^{(k)}$ and the finitely points in
$(0,1)$ which are not in the image of $\phi^{(k)}$.
The total measure of this subset is $\sum_{i=1,2,\ldots,n_k} w^{(k)}_i$ . Hence, the measure
of points where $f$ is defined is at most
$\lim_{k \to \infty} \sum_{i=1,2,\ldots,n_k} w^{(k)}_i$. This shows the
first part of the statement.

We now show by contradiction that $f$ does not admit periodic points. Let us
assume that $x_0 = x$, $x_1 = f x$, \ldots, $x_n = f^p x = x$ is a periodic orbit
where $p > 0$. Then there exists an open neighborhood $U$ of $x$ which is periodic
with the same translations: $f^i U = U + x_i - x_0$. None of the points in $U$
can be contained in a tower of height larger than $p$ because the image of a point
in a tower just go one level up. This means that for all $y$ in $U$ there is
an iterate $f^m(y)$ with $m \in \{0, 1, \ldots, p-1\}$ which is undefined. That
contradicts the fact that $f$ had a domain with full measure.
\end{proof}

The above lemma shows how one can obtain an interval exchange transformation
from a cutting and stacking construction. Another ingredient for the proof of
Theorem~\ref{thm:VershikArnouxOrnsteinWeiss} is a lemma to control the
accumulation points of the intervals. We propose this step as an exercise.

\begin{tcbexercise}{}{WorkAroundSpacers}
Let $C^{(k)} = (S_{w_1^{(k)}, h_1^{(k)}}, \ldots, S_{w^{(k)}_{n_k}, h^{(k)}_{n_k}})$,
$\phi^{(k)}: C^{(k)} \to (0,1)$, $f^{(k)}: (0,1) \to (0,1)$, and $f: (0,1) \to (0,1)$ be as in Lemma~\ref{lem:WellDefinedCutAndStack}.

Furthermore assume the following three hypotheses: for each $k$ the $n_k$-th
tower has height one, that is $h^{(k)}_{n_k} = 1$,
$\phi^{(k)}([0,w_{n_k}^{(k)}] \times \{0\}) = (1 - w_{n_k}^{(k)}, 1)$ and in
the cutting and stacking construction, all words $x_\ell$ starts and ends with a
letter $(n_k, j)$ for some $j$.

Then, show that the only accumulation points of the singularities in the domain and
codomain of $f$ are at $1$.
\end{tcbexercise}

We described how one can associte an interval exchange transformation
to a cutting and stacking construction. The second half of the proof of
Theorem~\ref{thm:VershikArnouxOrnsteinWeiss} consists in showing that any
invertible aperiodic measure-preserving transformation
admits a cutting and stacking construction. The main ingredient of
this part is the Kakutani-Rokhlin lemma proved simultaneously by
S.~Kakutani~\cite{Kakutani43} and V.~A.~Rokhlin~\cite{Rokhlin48}.

\begin{tcblemma}{Kakutani-Rokhlin}{Rokhlin}
Let $(X,T,\mu)$ be a measurable aperiodic measure-preserving transformation, $n > 0$ be an integer and $\epsilon > 0$ a real number. Then there exists a measurable set $Y \subset X$ so that
\begin{itemize}
\item $(Y, TY, \ldots, T^{n-1} Y)$ is a tower,
\item $\mu (Y \cup TY \cup \ldots \cup T^{n-1} Y) > 1 - \epsilon$.
\end{itemize}
In particular, $\mu(Y)=\frac{1-\epsilon}{n}$.
\end{tcblemma}
The above Lemma says that on the complement of a set of measure $\epsilon>0$ the dynamical system $T$ behaves like a tower of height $n$ and width $\mu(Y)$. In other words, the map
$T$ restricted to $Y \cup TY \cup \ldots T^{n-1}Y$ is measurably conjugated
to $T_{w,n}: S_{w,n} \to S_{w,n}$ where $w = \mu(Y)$.

In their proof of Theorem~\ref{thm:VershikArnouxOrnsteinWeiss}, Arnoux, Ornstein and Weiss use the Kakutani-Rokhlin lemma
to build one step of the cutting and stacking construction. Given that their exposition is well written and clear, we refer the reader directly to~\cite{ArnouxOrnsteinWeiss} for the rest of the details.



\medskip
\emph{Bratteli-Vershik diagrams}. The cutting and stacking construction we
presented can also be formally described using fully ordered Bratteli-Vershik
diagrams which are infinite graphs. For a more detailed discussion about these diagrams we refer the
reader to~\cite{BezuglyiKarpel16} and~\cite{LindseyTrevino16}.
\begin{tcbdefinition}{}{BratteliDiagram}
A \emphdef{Bratelli-Vershik diagram} is a quadruple $B = (V, E, s: E \to V, t: E \to V)$ where
\begin{itemize}
\item $V$ is a countable set endowed with a partition $V=\sqcup_{i\geq 0} V_i$ into
finite subsets $V_i$ called the \emph{vertex set},
\item $E$ is a countable set endowed with a partition $E=\sqcup_{i\geq 0}E_i$
into finite subsets $E_i$ called the \emph{edge set},
\item the maps $s$ and $t$ called the \emph{source} and \emph{target}
\end{itemize}
which satisfy the following conditions
\begin{enumerate}
\item $V_0 = \{v_0\}$,
\item for every $i\geq 0$ we have $s(E_i)=V_{i}$ and $t(E_i)=V_{i+1}$.
\end{enumerate}
A Bratteli-Vershik diagram $B$ is \emph{fully ordered} if for every vertex
$v \in V$ there is a total order $<_s$ on its outoing edges $s^{-1}(v)$ and
a total order $<_t$ on its incoming edges $t^{-1}(v)$.
\end{tcbdefinition}

\begin{tcbremark}{}{NoVershikMap}
The standard definition of a Bratteli-Vershik diagram requires also the existence of a homeomorphism (the Vershik map) on the space of infinite paths rooted at the top vertex $v_0$ which satisfies a series of technical conditions. We refer the reader to Definition~2.7 in~\cite{BezuglyiKarpel16}. For our purposes will not require such a map as our aim is to build infinite interval exchange transformations.
\end{tcbremark}

We now describe the fully ordered Bratteli-Vershik diagram associated to a cutting and stacking construction. The vertex set $V_k$ is the set of $n_k$ stacks at the step $k$ of the construction. It is a finite set and $V_0$ is reduced to a single point since we started the initial step consists of the single trivial stack $S_{1,1}$. The edges and the total orders $<_s$ and $<_t$ between the level $k$ and $k+1$ describe a cutting and stacking step. Let $S_{w_1^{(k)}, h_1^{(k)}}, \ldots, S_{w_{n_k}^{(k)}, h_{n_k}^{(k+1)}}$ the stacks at the $k$-th step.
More precisely, $E_k$ are the set of stacks obtained after the cutting. The source of an edge is the stack of $V_k$ it comes from and its target is the stack of $V_{k+1}$ it is stacked into.
The total order $<_s$ corresponds to the "horizontal" ordering of the interval $(0,1)$ while the ordering $<_t$ corresponds to the "vertical" ordering. In Figure~\ref{fig:CuttingAndStacking} we illustrate part of the Bratteli-Vershik diagram encoding an explicit cutting and stacking construction.

\begin{figure}[H]
\begin{center}
\includegraphics[scale=0.8]{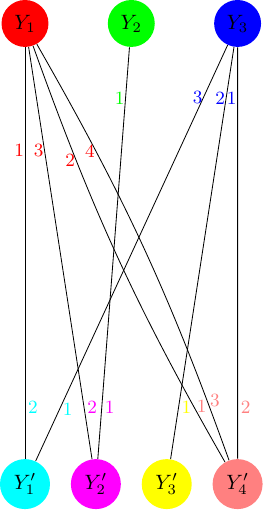}%
\end{center}
\caption{One level of the Bratteli-Vershik diagram encoding the combinatorial information of the cutting and stacking step of Figure~\ref{fig:CuttingAndStacking}.
\label{fig:BVDiagram}}
\end{figure}
In Figure~\ref{fig:BVDiagram} edges are oriented from top to bottom (\ie $Y_1$, $Y_2$ and $Y_3$ are sources and $Y'_1$, $Y'_2$, $Y'_3$ and $Y'_4$ are targets). The edges correspond to towers after the cutting step, namely $Y_{1,1}$, $Y_{1,2}$, $Y_{1,3}$, $Y_{1,4}$ (with source $Y_1$), $Y_{2,1}$ (with source $Y_2$) and $Y_{3,1}$, $Y_{3,2}$, $Y_{3,3}$ (with source $Y_3$)). The orders $<_s$ and $<_t$ appear on each edge near respectively the sources and targets.

\begin{tcbexercise}{}{CuttingStackingOdometer}
The idea of this exercise is first to show that the first return map of the vertical translation flow on baker's surface $B_{\frac{1}{2}}$ (see Figure~\ref{fig:VerticalBakerIET}) and dyadic odometer introduced in Example~\ref{exa:TowerOdometer} are essentially the same dynamical system. Second, to show that this dynamical system can be explicitly obtained by a cutting and stacking construction which can be encoded by a rather simple Bratteli-Vershik diagram.
\begin{enumerate}
\item Let $\Phi: \Z_2 \to [0,1]$ be the map that to a dyadic number $(x_n)\in\prod_{n\geq 0}\{0,1\}$ associates the real number $\Phi(x_n) = \sum \frac{x_n}{2^{-(i+1)}} $. Show that $\Phi$ is continuous and surjective.
\item Describe the points where $\Phi$ is not injective.
\item Let $T:\Z_2\to\Z_2$ be the diadic odometer. Show that $\Phi(T(x_n))= \Psi(\Phi(x_n))$, where $\Psi:(0,1)\to(0,1)$ is the first-return map to the transversal $I=(0,1)$ depicted in figure \ref{fig:RectangleAsTower} of the vertical translation flow on baker's surface $B_{\frac{1}{2}}$. In other words, show that the map $\Phi$ provides an explicit conjugation between the dyadic odometer and the interval exchange transformation from the baker map.
\item 4. Using Figure~\ref{fig:OdometerBV} as a reference, describe a cutting and stacking construction which produces $\Psi$. Does equation (\ref{eq:ConditionDefiningLimit}) hold in this case?
\end{enumerate}
\begin{figure}[H]
\begin{center}\includegraphics{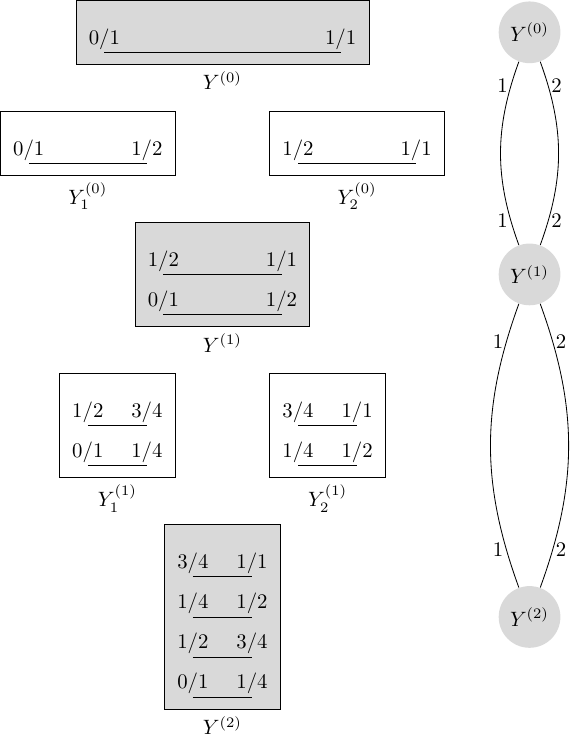}\end{center}
\caption{The cutting and stacking construction for the vertical flow in the Baker's surface $B_{1/2}$ (on the left) and its associated Bratteli-Vershik diagram (on the right).}
\label{fig:OdometerBV}
\end{figure}
\end{tcbexercise}

\section{Minimality and Keane criterion}
\label{sec:CounterexampleKeane}
Keane's theorem asserts that for \emph{finite} interval exchange transformations,
the absence of connection implies minimality (see~\cite{Keane75}). The aim of
this section is to define minimality and complete periodicity and show that
Keane's theorem does not hold for infinite interval exchange transformations.

\begin{tcbdefinition}{}{MinimalityPeriodicity}
Let $T: X \to X$ be a homeomorphism on a compact space $X$. We say that $T$ is
\emph{minimal} if for all $x \in X$ its orbit $\{T^n(x)\}_{n \geq 0}$ under $T$
is dense in $X$. We say that $T$ is \emph{completely periodic} if for all $x \in X$
its orbit is finite, in other words, there exists $n > 0$ such that $T^n(x) = x$.
\end{tcbdefinition}
If $x$ has finite orbit, the length of the orbit is called the \emph{period}.

\begin{tcbexercise}{}{CompletePeriodicity}
Show that in the definition of completely periodic we can invert quantifiers. Namely,
given $T: X \to X$ a homeomorphism of a compact space the following are equivalent
\begin{itemize}
\item for all $x \in X$ there exists $n > 0$ such that $T^n(x) = x$,
\item there exists $n > 0$ such that for all $x \in X$ we have $T^n(x) = x$.
\end{itemize}
\end{tcbexercise}

Note that interval exchange transformations are not homeomorphisms and Definition~\ref{def:MinimalityPeriodicity}
does not apply. By extension we use the following definition.
\begin{tcbdefinition}{}{}
Let $I$ be a bounded open interval and $f: D \to R$ be a finite or infinite interval exchange transformation
where $D, R \subset I$. Let $D_\infty$ the set of points with infinite forward orbit under $f$.
We say that $f$ is \emphdef[minimal (iet)]{minimal} if for all $x \in D_\infty$
its orbit $\{f^n(x)\}_{n \geq 0}$ is dense in $I$. We say that $f$ is \emph{completely periodic}
if for all $x \in D_\infty$ its orbit is finite, in other words, there exists $n > 0$ such that $f^n(x) = x$.
\end{tcbdefinition}

\begin{tcbexercise}{}{}
This is a continuation of Exercise~\ref{exo:CompletePeriodicity}. Build an
example of completely periodic interval exchange transformation $f: D \to R$ but
for which periods of points are unbounded. In other words
\begin{itemize}
\item for all $x \in D_\infty$ there exists $n > 0$ such that $f^n(x) = x$,
\item there does not exist $n > 0$ such that for all $x \in D_\infty$ we have $f^n(x) = x$.
\end{itemize}
\end{tcbexercise}

Recall from Definition~\ref{def:IETConnection} that a connection for
an interval exchange transformation is a triple $(m, x, y)$ such that
$m \geq 0$ is an integer, $x \in D_m \setminus (R \cup D_{m+1})$
and $f^m(x) = y$. In particular, $x$ is a backward singularity of $f$,
$y$ is a forward singularity of $f$.

Recall that in the case of finite interval exchange transformations if all
singularities are part of a connection then the interval exchange transformation
is completely periodic. This breaks dramatically in the infinite setting.
\begin{tcbtheorem}{}{InfinitePeriodicCounterexample}
There exists an interval exchange transformation $f: D \to R$ where
$D, R \subset (0,1)$ such that
\begin{compactitem}
\item the accumulation points in $[0,1]$ of the singularities of
$f$ and $f^{-1}$ are contained respectively in $\{1\}$ and $\{0\}$,
\item the forward orbits of backward singularities of $f$ are finite, \ie each of them is part of a connection,
\item the backward orbits of forward singularities of $f$ are finite,
\item for any $x \in D_\infty$, its forward orbit is dense.
\end{compactitem}
\end{tcbtheorem}

\begin{proof}
We already constructed an example of this kind, namely the interval exchange transformation
obtained as the first return map of the vertical flow in the baker's surface,
see Section~\ref{sec:Jungle} and more specifically Example~\ref{exa:BakerIET}.
The forward singularities of $f_{1/2,\pi/2}$ are the rational numbers
$\{0, 1/2, 3/4, 7/8, \ldots, 1\}$. Its backward singularites are
$\{1, 1/2, 1/4, 1/8, \ldots, 0\}$. For each $n \geq 1$ we have a
connection $(2^{-n}, 1 - 2^{-n}, 2^{n-1})$
\end{proof}

The objective for the remaining of this section is to prove the following theorem.
\begin{tcbtheorem}{}{InfiniteKeaneCounterexample}
Let $C_3 \subset [0,1]$ be the triadic Cantor set.
There exists an infinite interval exchange transformation $f: D \to R$ with $D, R \subset (0,1)$ that satisfies the following properties.
\begin{enumerate}
\item the accumulation points in $[0,1]$ of the singularities of
$f$ and $f^{-1}$ are contained respectively in $\{1\}$ and $\{0\}$,
\item $C_3 \subset D \cap R$ and $f(C_3) = C_3$,
\label{propitem:CantorPreserved}
\item $f$ has no connection.
\label{propitem:NoConnection}
\label{propitem:Accumulation}
\end{enumerate}
\end{tcbtheorem}

\begin{tcbexercise}{}{}
Show that neither constructions of Theorem~\ref{thm:InfinitePeriodicCounterexample} nor
Theorem~\ref{thm:InfiniteKeaneCounterexample} are possible with finite interval
exchange transformations.
\end{tcbexercise}

\begin{proof}
We first define an interval exchange transformation that satisfies~\ref{propitem:CantorPreserved} and~\ref{propitem:Accumulation}. Next, we will modify this construction to get rid of the saddle connections so that the resulting interval exchange transformation satisfies~\ref{propitem:NoConnection}.

We now describe the first construction. A picture of the resulting interval
exchange transformation can be seen in Figure~\ref{fig:OdometerCantor}.
Let $\cA$ be the infinite alphabet made of the letters $a_n$ and $b_n$,
where $n$ are positive integer indices.
For $n \geq 1$, let $I^{top}_{a_n}$ and $I^{top}_{b_n}$ be respectively the
intervals of real numbers in $[0,1]$ whose 3-adic expansions start with
$0.2^{n+1}0$ and $0.2^{n+1}1$. For example $I^{top}_{a_1} = [0,1/3]$,
$I^{top}_{b_1} = [1/3,2/3]$, $I^{top}_{a_2} = [2/3,8/9]$, $I^{top}_{b_2} = [8/9,10/9]$.
Let also $I^{bot}_{a_n} := 1 - I^{top}_{a_n}$ and $I^{bot}_{b_n} = 1 - I^{top}_{b_n}$.
The list of intervals $\{I^{top}_x: x \in \cA\}$ and $\{I^{bot}_x: x \in \cA\}$ define
two partitions of $(0,1)$ labeled by $\cA$. We denote $f_1: D \to R$ the associated
interval exchange transformation.
\begin{figure}[!ht]
\begin{center}
\includegraphics{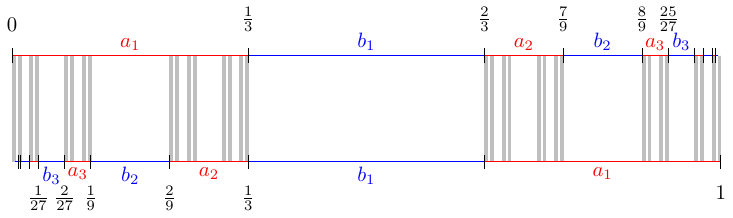}
\end{center}
\subcaption{The interval exchange $f_1$. Its
restriction to the triadic Cantor set $C_3$ is isomorphic to the diadic odometer. Each singularity is part of
a connection.}
\label{fig:OdometerCantor}
\vspace{.5cm}
\begin{center}
\includegraphics{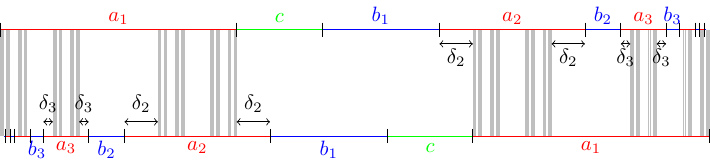}
\end{center}
\subcaption{The interval exchange $f_2$. Its restriction to $C_3$ coincide with $f_1$.}
\label{fig:OdometerCantorPerturbed}
\caption{The interval exchanges $f_1$ and $f_2$ from the proof of Theorem~\ref{thm:InfiniteKeaneCounterexample}.
We have ${f_1}|_{C_3}={f_2}|_{C_3}$.}
\end{figure}
Note that the order of the intervals is the exact same as in the odometer
$\pi_{top} = (a_1 b_1 a_2 b_2 \ldots)$ and $\pi_{bot} = (\ldots b_2 a_2 b_1
a_1)$.

It is easy to see that:
\begin{enumerate}
\item $f_1$ restricted to the triadic Cantor set is conjugate to the odometer of
Example~\ref{exa:TowerOdometer}, and
\item all other orbits are periodic and the period of their codings have one of the following forms:
$b_1$, $b_2 a_1$, $b_3 a_1 a_2 a_1$, $b_4 a_1 a_2 a_1 a_3 a_1 a_2 a_1$, \ldots they can be
inductively defined as $b_i u_i$ with $u_1 = \emptyset$ and $u_{i+1} = u_i a_i u_i$.
\end{enumerate}

For the interval exchange $f_1$, each singularity is part of a connection. We now perturb $f_1$ as follows.
Let $(\delta_i)_{i \geq 1}$ be a decreasing sequence of positive real numbers so that $\delta_{i} + \delta_{i+1} < \frac{1}{3^i}$.
Let us define the following operation on intervals $B([x,y], a, b) := [x-a,y+b]$.
We define a new interval exchange transformation $f_2$ by setting
\[
I'^{top}_{a_1} := I^{top}_{a_1} = \left[0,1/3\right]
\hspace{5mm}
I'^{top}_c     := \left[1/3,1/3+\delta_1\right]
\]

\[
I'^{top}_{a_i} := B(I^{top}_{a_i}, \delta_i, \delta_i)
\hspace{5mm}
I'^{top}_{b_i} := B(I^{top}_{b_i}, -\delta_i, -\delta_{i+1})
\]

\[
I'^{bot}_{a_1} := I^{bot}_{a_1} = \left[2/3,1\right]
\hspace{5mm}
I'^{bot}_c     := \left[2/3-\delta_1, 2/3\right]
\]

\[
I'^{bot}_{a_i} := B(I^{bot}_{a_i}, \delta_i, \delta_i)
\hspace{5mm}
I'^{bot}_{b_i} := B(I^{bot}_{b_i}, -\delta_{i+1}, -\delta_i).
\]
These new intervals define an interval exchange transformation $f_2$ on $[0,1]$ (see Figure~\ref{fig:OdometerCantorPerturbed}).
By construction, $f_1 = f_2$ on each of $I^{top}_{a_i}$. In particular, $f_2$ preserves the triadic Cantor set.

Now we claim that if the $\delta_i$
\begin{compactitem}
\item are rationally independent,
\item $\delta_{i+1} < \delta_i$, and
\item $\lambda_{b_i} = \frac{1}{3^i} - \delta_i - \delta_{i+1} < 2 (\delta_i - \delta_{i+1})$
\end{compactitem}
then $f_2$ has no saddle connection.

On each interval $I'^{top}_{x}$ we denote by $\tau_x$ the corresponding translation modulo $\Q$. We have
\[
\tau_{a_i} = 0
\qquad
\tau_c = -\delta_1
\qquad
\tau_{b_i} = \delta_{i+1} - \delta_i
\]

Let us denote by $I^{top}_{a_i,\ell}$ (respectively $I^{top}_{a_i,r}$) the left
and right part added to the interval $I^{top}_{a_i}$ to become $I'^{top}_{a_i}$
and let us consider the induced interval exchange transformation on the unions
of the $a_{i,\ell}$, $a_{i,r}$ for $i \geq 2$.
\begin{enumerate}
\item all points in $I'^{top}_{a_i,r}$ go into $I'^{top}_{a_{i-1},r}$ except for $i=2$ for which $T(I'^{top}_{a_{2,r}}) \subset I'^{top}_c$.
\item the points in $I'^{top}_{a_i,\ell}$ either go in one step in $I^{top}_{a_{i+1},\ell}$ or in many steps through $b_i$ and $a_i$ in $I^{top}_{a_i,r}$.
\end{enumerate}

\begin{figure}[!ht]
\begin{center}\includegraphics{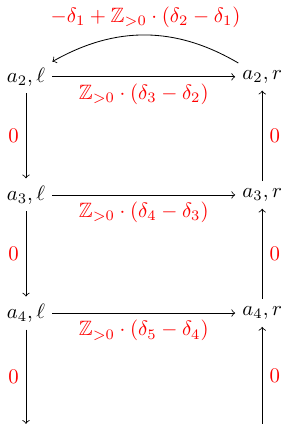}\end{center}
\caption{The translations modulo $\Q$ of the induced map on the added extremities of $I^{top}_{a_i}$. The value of $n_i$ is always a number larger than $2$.}
\label{fig:OrbitCantorPerturbed}
\end{figure}
We now prove that $f_2$ has no connection. Recall that a connection is a triple
$(x,y,m)$ where $x$ is a backward singularity, $y$ is a forward singularity and
$(f_2)^m(x) = y$. The backward or forward singularity modulo $\Q$ are
$-\delta_i$ and $+\delta_i$ for respectively the left and right side of the
interval labelled $a_i$.  We denote $\delta(a_{i,\ell}) = -\delta_i$ and
$\delta(a_{i,r}) = \delta_i$.

Recall that a connection is a finite orbit starting at a backward singularity $x$
and ending at a forward singularity $y$. This could be written $x + \tau = y$
where $\tau$ is the sum of the translations realized along the orbit of $x$. For our
purpose we consider this equation modulo $\Q$. As we already mentioned the values of
$x$ and $y$ are $\pm \delta_i$ depending whether they belong to $a_{i,\ell}$ or
$a_{i,r}$. The value of $\tau$ is obtained by following the orbit along the graph
of~\ref{fig:OrbitCantorPerturbed} and summing the weights on the edges. By the assumption on the linear independence of the $\delta_i$ no path exists that would lead to
an equation of the form $\delta(u) + \tau(\gamma) = \delta(v)$ where
$\gamma$ is an oriented path in the graph from $u$ to $v$ and $\tau(\gamma)$
is the sum of the weights along $\gamma$.
\end{proof}

\section{Entropy} \label{sec:entropy}
The entropy of a dynamical system is a non-negative real number that quantifies the degree of disorder of a transformation. We recall its definition in order to fix notations and refer the reader to ~\cite{CornfeldFominSinai}, ~\cite{Downarowicz-book}, ~\cite{KrerleyViana16} or~\cite{Dajani-Kalle2021} for more details.

Let $(X,\cB,\mu)$ be a standard Borel probability space. Given a (measurable) countable partition $\cP$ of $X$ we define its \emph{entropy} relative to $\mu$ as
\[
H(\mu, \cP) = \sum_{A \in \cP} - \mu(A) \log(\mu(A)).
\]
A finite partition of cardinality $n$ has entropy less than $\log(n)$, because the function $x\to-x\log(x)$ is concave on $(0,1)$. On the other hand, an infinite partition might have infinite entropy (\emph{e.g.} Example 9.1.4 in~\cite{KrerleyViana16}). For any countable family of partitions $\cP_n$ define:
$$
\bigvee_{n}\cP_n:=\{\cap_n P_n\hspace{1mm}|\hspace{1mm} P_n\in\cP_n\hspace{1mm}\text{for each $n$}\}.
$$
Given a measurable transformation $T:X \rightarrow X$ preserving the measure $\mu$ and a partition $\cP$ with finite entropy
, the (metric) \emph{entropy of $f$ relative to $\cP$} is given by
\[
h(T,\mu,\cP) = \lim_{n \to \infty} \frac{H\left(\mu, \bigvee_{k=0}^{n-1} T^{-k} \cP \right)}{n}.
\]
The above limit is equal to $\inf \frac{H\left(\mu, \bigvee_{k=0}^{n-1} T^{-k} \cP \right)}{n}$ as the sequence \linebreak  $n \mapsto H(\mu, \bigvee_{k=0}^{n-1} \cP^{(n)})$ is subadditive. The \emph{entropy} of the transformation $T$ is
\begin{equation}
	\label{eq:DefEntropy}
h(T,\mu):=\sup_\cP h(T,\mu,\cP)
\end{equation}
where $\cP$ ranges over all partitions with finite entropy\footnote{As a matter of fact this definition is not affected if the supremum is taken over only over finite partitions.}. When the measure is implicit in the discussion we write just $h(T)$ for the entropy of the transformation $T$, \emph{e.g.} for IETs Lebesgue is the implicit measure.
\begin{tcbexercise}{}{PeriodicSystemsHaveZeroEntropy}
\begin{enumerate}
\item Let $T$ be a measure preserving transformation on $(X,\mathcal{B},\mu)$ and suppose that $X = A \sqcup B$ with $A$ and $B$ measurable and $T$-invariant. Show that $h(T, \mu) = \mu(A) h(T|_A, \nu_A) + \mu(B) h(T|_B, \nu_B)$ where
$\nu_A$ and $\nu_B$ are the probability measures on $A$ and $B$ obtained by restricting $\mu$ and renormalizing, respectively.
\item Let $T$ be a measure preserving transformation on $(X,\mathcal{B},\mu)$ such that  for every $x\in X$ the orbit $\{T^n(x)\}_{n\in\N}$ is finite. Show that $h(T,\mu)$ is zero.
\item Let $T_i:(X_i,\mathcal{B}_i,\mu_i)$, $i=1,2$, be measurable-equivalent measure preserving transformations on two probability spaces\footnote{That is to say there exist measurable sets of full measure $X_i'\subset X_i$ and a  measurable bijection $\phi:X_1'\to X_2'$ with measurable inverse such that $\phi_*\mu=\nu$ and $\phi\circ T_1=T_2\circ\phi$ }. Prove that $h(T_1,\mu_1)=h(T_2,\mu_2)$.
\end{enumerate}
\end{tcbexercise}

As we show later (see Corollary~\ref{cor:EntropyFiniteIET}), finite interval exchange transformations have zero entropy. The main aim of this section is to give an upper bound estimate for the measured entropy of an infinite interval exchange. More precisely:

\begin{tcbtheorem}{}{EntropyInfiniteIET}
Let $f = f_{\pi,\lambda}$ be an interval exchange transformation on $[0,1]$.
Let $\lambda_1 \geq \lambda_2 \geq \ldots$ be the lengths of the subintervals (ordered by decreasing sizes) and $\Lambda_m := \lambda_{m+1} + \lambda_{m+2} + \ldots$. Then
\[
h(f) \leq \liminf_{m \to \infty} \log(m) \Lambda_m.
\]
\end{tcbtheorem}

We also recall the construction (see Example~\ref{exa:PositiveEntropyIET}) of an infinite IET having infinite entropy.

Theorem~\ref{thm:EntropyInfiniteIET} is a slightly modified version of a result that appear in~\cite{Blume2012}. For its proof we use some corollaries of Kolmogorov-Sinai's theorem (see~\cite{KrerleyViana16} for details) and an estimate on the combinatorial complexity of interval exchange transformations. We address these tools in what follows.

\begin{tcbtheorem}{Kolmogorov-Sinai}{KolmogorovSinai}
Let $T$ be a measurable, measure-preserving transformation of a probability space $(X,\mathcal{B},\mu)$ and
$\mathcal{P}_1\prec\ldots\prec\mathcal{P}_n\prec\ldots$ a non-decreasing sequence of partitions with finite entropy such that $\cup_{n=1}^\infty\mathcal{P}_n$ generates the $\sigma$-algebra $\mathcal{B}$, up to measure zero. Then
$$
h(T,\mu)=\lim_{n\to\infty} h(T,\mu,\mathcal{P}_n)
$$
\end{tcbtheorem}

Recall that given two partitions $\cP_1$ and $\cP_2$ we say that $\cP_2$ is finer than $\cP_2$ (and denote it by $\cP_1\prec\cP_2$) if every element of $\cP_2$ is contained in some element of $\cP_1$. On the other hand, a subset $\mathcal{B}'\subset\mathcal{B}$ generates the $\sigma$-algebra $\mathcal{B}$ up to measure zero if the $\sigma$-algebra generated by $\mathcal{B}'$ is equal to $\mathcal{B}$ up to measure zero\footnote{Let $\mathcal{B}_1$ and $\mathcal{B}_2$ be any two subsets of a $\sigma$-algebra $\mathcal{B}$ is a probability space $(X,\mathcal{B},\mu)$. Then $\mathcal{B}_1\subset\mathcal{B}_2$ up to measure zero if for every $B_1\in\mathcal{B}_1$ there exists $B_2\in\mathcal{B}_2$ such that $\mu(B_1\Delta B_2)=0$, where $B_1\Delta B_2=(B_1\cup B_2)\setminus(B_1\cap B_2)$}.

\begin{tcbcorollary}{}{EntropyDecreasingDiamaterPartition}
Suppose that $X$ is a metric space with a Borel probability measure $\mu$. Let $T:X\to X$ be a measurable, measure-preserving transformation and $\mathcal{P}_1\prec\ldots\prec\mathcal{P}_n\prec\ldots$ a non-decreasing sequence of partitions with finite entropy.  Let $\mathcal{P}_n(x)$ denote the element of $\mathcal{P}_n$ containing $x\in X$. Suppose that for $\mu$-almost every $x\in X$:
$$
\lim_{n\to\infty}{\rm diam}(\mathcal{P}_n(x))=0.
$$
Then $h(T,\mu)=\lim_{n\to\infty}h(T,\mu,\mathcal{P}_n)$
\end{tcbcorollary}



\emph{Combinatorial complexity}. Let $f:D\to R$ be a partial interval exchange transformation on $I$ and $\cP = \{P_a\}_{a \in A}$ a finite partition of $I$. For every $n\geq 1$ one can consider the map $D_n \to A^n$ which associates to an orbit of length $n$, $x$, $fx$, $f^2 x$, \ldots, $f^{n-1} x$ its coding $w_0 w_1 w_2 \ldots w_{n-1}$ with respect to the partition $\cP$. That is, $w_i$ is the unique element of $A$ so that $f^i x \in P_{w_i}$. The \emph{combinatorial complexity} of a partial interval exchange transformation $f: D \to R$ is the function $p: \N \to \N$ where $p(n)$ is the number of different codings of orbits of length $n$ of $f$. In other words, $p(n)$ the cardinality of the image of the  ``coding map" $D_n \to A^n$.
\begin{tcbproposition}{}{CodingFiniteIET}
Let $f: D \to R$ be a partial interval exchange transformation with finitely many intervals, let say $d$. Then its combinatorial complexity with respect to the set of connected components of $D$ satisfies $p(n) \leq (d-1)n + 1$. Moreover, if it is an interval exchange transformation without connections then $p(n) = (d-1)n + 1$.
\end{tcbproposition}
\begin{proof}
Let $\cI = \{I_1, I_2, \ldots, I_d\}$ be the set of connected components of $D$. We consider the coding of $f$ with respect to this natural partition. Let $f^n: D_n \to D_{-n}$ denote the $n$-th iterate of $f$. Now to each coding $w = w_0 w_1 \ldots w_{n-1}$ of length $n$ with respect to $\cI$, one can associate the interval $I_w\subset D$ formed by all points whose orbit of length $n$ has exactly the coding $w$. Namely
\[
I_w = I_{w_0} \cap f^{-1} (D_{-1} \cap I_{w_1}) \cap f^{-2} (D_{-2} \cap I_{w_2}) \cap \ldots f^{-n-1} (D_{-n-1} \cap I_{w_{n-1}}).
\]
$I_w$ is a maximal interval of the domain $D_{n-1}$ of $f^{n-1}$. Therefore $p(n)$ is equal to the number of connected components of $D_{n-1}$ (or equivalently of $D_{-n+1}$). Note that $D_{n} = f^{-n+1} (D_{-n+1} \cap D)$. Hence, as $f^{-n+1}$ is continuous on $D_{-n+1}$, the number of connected components of $D_{n}$ is exactly the number of connected components of $D_{-n+1} \cap D$. As $D$ is made of $d$ intervals, this number increases at most by $d-1$ as we iterate $f$. Therefore $p(n)$ is bounded above by $d+(n-1)(d-1)=(d-1)n+1$.
If $f$ is an interval exchange transformation without connections, each time we iterate $f$ exactly $d-1$ new separation points appear in the partition of $D_n$.
\end{proof}
\begin{tcbcorollary}{}{EntropyFiniteIET}
Finite interval exchange have zero entropy.
\end{tcbcorollary}
\begin{proof}
A classical result by Keane~\cite{Keane75} (that can be traced back to the work of Mayer~\cite{Mayer43}, see also~\cite{ZemljakovKatok-transitivity1} Section 2) implies that the domain of any IET on a finite number of intervals is decomposed into subsets where the IET is either periodic or minimal.  By Exercise~\ref{exo:PeriodicSystemsHaveZeroEntropy} it is sufficient to show that the restriction of $f$ to a minimal componet has entropy zero. Hence let $f:D\to R$ be a minimal interval exchange transformation and $\mathcal{I}=\{I_1,I_2,\ldots,I_d\}$ be the corresponding intervals. Since there are no periodic components, Corollary~\ref{cor:PeriodicOrbitsIETObstruingShrinking} implies that for every $\varepsilon>0$ there exists $n\in\N$ such that every interval $J\subset \bigvee_{k=0}^{n-1}f^{-k}(\mathcal{I})$ has length $|J|<\varepsilon$. Hence $\{\bigvee_{k=0}^{n-1}f^{-k}(\mathcal{I})\}_{n\geq 1}$ generates the Borel $\sigma$-algebra and therefore $h(f,\mu)=h(f,\mu,\mathcal{I})$. By Proposition~\ref{prop:CodingFiniteIET}, $\displaystyle \# \bigvee_{k=0}^{n-1} f^{-k}(\cI) \leq (d-1)n + 1$ and thus
\[
\frac{H(\mu, \bigvee_{k=0}^{n-1} \cI)}{n} \leq \frac{\log((d-1)n+1)}{n}\xrightarrow[n\to\infty]{} 0
\]
\end{proof}


\begin{tcbexample}{}{PositiveEntropyIET}
We present a cutting and stacking construction for a map of positive entropy that appears in~\cite{Shields73}. On the first step as usual we consider the stack $\mathbf{S}(1,0)=(0,1)$. In the second step we cut the unit interval in two pieces of equal length but we do not stack them. This produces two stacks $\mathbf{S}_1(1,\frac{1}{2})=(0,\frac{1}{2})$ and $\mathbf{S}_2(1,\frac{1}{2})=(\frac{1}{2},1)$. Now suppose that at the $k^{th}$ step, $k\geq 3$, one has stacks $\{\mathbf{S}_1,\ldots,\mathbf{S}_{q_k}\}$ having all the same width $w_k$. For each $1\leq i\leq q_k$ cut $\mathbf{S}_i$ into stacks $\mathbf{S}_i^j$, $1\leq j\leq 2q_k$, according to the parameters $p_{ij}=\frac{w_k}{2q_k}$. Then stack them into $q_k^2$ new stacks according to the rule:
\begin{equation}
	\label{eq:PositiveEntropyExample}
\mathbf{S}_{i,j}:=\mathbf{S}_j^i\mathbf{S}_{q_k+i}^j,\hspace{3mm}1\leq i\leq q_k,\hspace{3mm}1\leq j\leq q_k.
\end{equation}
Since each of the new $q_k^2$ stacks has width $\frac{w_k}{2q_k}$ the term defining the limit in~\ref{eq:ConditionDefiningLimit} satisfies
\[
\sum_{i=1,2,\ldots,q_k^2} w^{(k+1)}_i=\frac{1}{2}\sum_{i=1,2,\ldots,q_k} w^{(k)}_i
\]
and therefore the construction defines an infinite IET $f$ on $I=(0,1)$.
\begin{figure}[H]
\begin{center}
\includegraphics{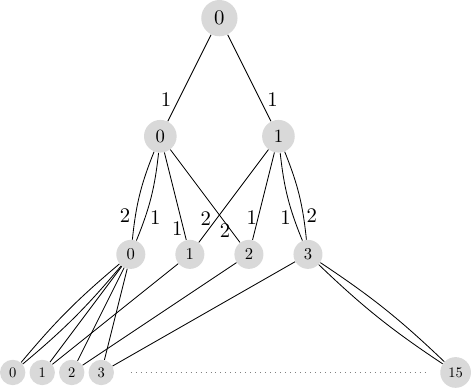}
\end{center}
\caption{The first four levels of a Bratteli-Vershik diagram of a map of positive entropy.}
\label{fig:BVPositiveEntropy}
\end{figure}
The corresponding Bratteli-Vershik diagram is depicted in Figure~\ref{fig:BVPositiveEntropy}. Here $V_0 = \{0\}$, for $k\geq 1$ we set $V_k = \{0, 1, \ldots, 2^{2^{k-1}}-1\}$ and arrows are defined by (\ref{eq:PositiveEntropyExample}). There are several ways to calculate que entropy of $f$. For example, from~\cite{Shields73} explains that $h(f)$ is given by $\lim_{k\to\infty}W_nq_k$, where $W_k$ is the sum of the widths of the $q_k$ stacks at the $k^{th}$ step of the construction. Therefore
$$
h(f)=\lim_{k\to\infty}(q_kw_k)\log(q_k)=\lim_{k\to\infty}\frac{1}{2^k}\log(2^{2^k})=\log(2)=0.69314...
$$
\end{tcbexample}

\begin{proof}[Proof of Theorem~\ref{thm:EntropyInfiniteIET}]
We assume that
\begin{equation} \label{eq:entropy_definition_C}
C := \liminf_{m \to \infty} \log(m) \Lambda_m < +\infty
\end{equation}
otherwise there is nothing to prove.

Let us denote by $I_i$ the subinterval of $[0,1]$ of size $\lambda_i$.
We consider the finite partitions induced by the domain of the IET of the form $\alpha_m = \{I_1, I_2, \ldots, I_m, J_m\}$ where
\[
J_m := \bigcup_{m' > m} I_{m'}.
\]
Note that $\mu(J_m) = \Lambda_m$.

Let $\alpha_m^{(n)} := \bigvee_{k=0}^{n-1} f^{-k} \alpha_m$. We claim that
 $$
\displaystyle h(f) = \lim_{m \to \infty} \lim_{n \to \infty} \frac{H(\mu, \alpha_m^{(n)})}{n} = \sup_{m \to \infty} \inf_{n \to \infty} \frac{H(\mu, \alpha_m^{(n)})}{n}.
$$
Indeed, since $H(\mu, \alpha_m^{(n)})$ is subadditive we have $\lim_{n \to \infty} \frac{H(\mu, \alpha_m^{(n)})}{n} = \inf_{n \to \infty} \frac{H(\mu, \alpha_m^{(n)})}{n}$. The sequence of partitions $\alpha_m$ is nested, that is $\alpha_{m+1}$ is a refinement of $\alpha_m$. Hence $h(\mu, \alpha_m) \leq h(\mu, \alpha_{m+1})$ and the limit in $m$ is a sup. This justifies the right hand side equality above.

We now justify the left hand side equality above. First remark that the entropy of $f$ is equal to the entropy of $f$ restricted to the complement of periodic components. Now consider for every $k\in\N$ the partition $\alpha_k^{(k)}$. Corollary~\ref{cor:PeriodicOrbitsIETObstruingShrinking} implies that for $\mu$-almost every $x$ not in a periodic component $\lim_{k\to\infty}{\rm diam} \alpha_k^{(k)}(x)=0$. The claim then follows from Corollary~\ref{cor:EntropyDecreasingDiamaterPartition}.

As a consequence, if $(m_k)$ is a diverging sequence and $(n_k)$ is any sequence, then
\[
h(f) \leq \liminf_{k \to \infty} \frac{H(f,\alpha_{m_k}^{n_k})}{n_k}.
\]
(simply because $H(f, \alpha_{m_k}^{n_k}) / n_k \geq h(f, \mu, \alpha_{m_k})$).

We code the orbits of the interval with respect to the partition $\alpha_m$ by using the letter $i$ for the interval $I_i$ and $*$ for the interval $J_m$. Let $\Delta_{m,n,k}$ be the set of finite codings of length $n$ that contains exactly $k$ times the letter $*$.
Let also $\Omega_{m,n,k}$ be the subset of $[0,1]$ that consists of points whose codings of length $n$ belongs to $\Delta_{m,n,k}$. For example, $\Omega_{m,n,0}$ is the set of points whose first $n$-th iterates never enter $J_m$.
Note that $\cA_n = \{\Omega_{m,n,k}\}_{k=0}^n$ is a partition of $[0,1]$ that is coarser than $\alpha_m^{(n)}$.

It follows from Proposition~\ref{prop:CodingFiniteIET} that $\# \Delta_{m,n,0} \leq (m-1)n+1 \leq mn$. As a consequence, we deduce a naive estimate for $\# \Delta_{m,n,k}$. By decomposing over all possible positions for the symbol $*$ we have
\[
\# \Delta_{m,n,k} \leq \sum_{\ell_1, \ell_2, \ldots, \ell_{k+1}} \prod_{i=1}^{k+1} (m \ell_i)
\]
where the sum runs over all non-negative integer vectors $(\ell_1, \ell_2, \ldots, \ell_{k+1})$ such that $\ell_1 + \ell_2 + \ldots + \ell_{k+1} = n-k$. The number of terms in the right hand side is $\binom{n}{k}$ for which an upper bound is given by $\binom{n}{k} = \prod_{i=0}^{k-1} \frac{n-i}{k-i} \leq n^k$. Using the trivial bound $\ell_i \leq n$, we obtain
\begin{equation} \label{eq:EntropyCounting1}
\# \Delta_{m,n,k} \leq n^k \left( m n \right)^{k+1} \leq (n^2m)^{k+1}.
\end{equation}

Now we claim that
\begin{equation} \label{eq:EntropyCounting2}
\sum_{k=0}^n k\, \mu(\Omega_{m,n,k}) = n \Lambda_m.
\end{equation}
Indeed, $\cA_n = \{\Omega_{m,n,k}\}_{k=0}^n$ is a partition of $[0,1]$ and the left
hand side in~\eqref{eq:EntropyCounting2} is the expectation of the number of
letter $*$ in a coding of length $n$. This is equal to $n$ times the
expectation of a letter $*$ which is $\Lambda_m = \mu(J_m)$. More formally, we have
\begin{align*}
\sum_{k=0}^n k \mu(\Omega_{m,n,k})
&= \int_0^1 \left( \mathbbm{1}_{J_m} + \mathbbm{1}_{J_m} \circ f + \ldots + \mathbbm{1}_{J_m} \circ f^{n-1} \right) d \mu \\
&= \sum_{k=0}^{n-1} \int_0^1 \mathbbm{1}_{J_m} \circ f^k d \mu \\
&= n\, \int_0^1 \mathbbm{1}_{J_m} d \mu = n\, \mu(J_m)\\
\end{align*}
where in the second line we used the linearity of the integral and in the third line the $f$-invariance of $\mu$. This concludes the proof of the claim.

Let us denote $g(x) = - x \log(x)$. In order to reduce entropy estimates to counting, we use the following inequality that follows from the concavity of the function $g$: if $A_i$, $i \in \cI$ are disjoint measurable sets whose union is $A$ then $\sum_{i \in \cI} g(\mu(A_i)) \leq \mu(A) \log(|\cI|) + g(\mu(A))$ \footnote{this can also be shown as a consequence of the following formula for conditional entropy: if $\cP = \{A_1,A_2,\ldots\}$ and $\cQ = \{B_1,B_2,\ldots\}$ are (measurable) partitions then $H(\cP|\cQ) = \sum_{B \in \cQ} \mu(B) H_B(\cP)$ where $H_B(\cP)$ is the entropy of $\cP$ for the measure $\mu_B$ induced on $B$.}.
Now, using~\eqref{eq:EntropyCounting1}, we obtain:
\begin{align}
H(\alpha_m^{(n)})
&= \sum_{k=0}^n \sum_{w \in \Delta_{m,n,k}} g(\mu(I_w)) \label{EQ:one} \\
&\leq \sum_{k=0}^n \left(\mu(\Omega_{m,n,k}) \log\# \Delta_{m,n,k} + g(\mu(\Omega_{m,n,k}))\right) \label{EQ:two} \\
&\leq \sum_{k=0}^n \left( \mu(\Omega_{m,n,k}) \log \# \Delta_{m,n,k}\right) + H(\{\Omega_{m,n,k}\}_{k=0}^n)  \label{EQ:three}\\
&\leq \sum_{k=0}^n (k+1) \mu(\Omega_{m,n,k}) \log(n^2 m) + \log(n+1)\label{EQ:four} \\
&\leq \log(n^2m) (n \Lambda_m + 1 ) + \log(n+1) \label{EQ:five}
\end{align}
Indeed, (\ref{EQ:one}) and (\ref{EQ:three}) follow by definition, (\ref{EQ:two}) uses the the counting bound with $A_i = I_w$ and $A = \Omega_{m,n,k}$, (\ref{EQ:four}) uses \eqref{eq:EntropyCounting1} and $H(\cP) \leq \log \# \cP$, and (\ref{EQ:five}) follows from~\eqref{eq:EntropyCounting2}.

Now, let $m_k$ so that $\log(m_k) \Lambda_{m_k}$ tends to $C$ defined in~\eqref{eq:entropy_definition_C}.
Pick $n_k$ so that
\[
\lim_{k \to \infty} \frac{\log(n_k)}{\log(m_k)} = \lim_{k \to \infty} \frac{\log(m_k)}{n_k} = 0.
\]
(For example, one can choose $n_k = \lfloor \log(m_k) \rfloor^2$). Then
\begin{align*}
h(f,\alpha)
&\leq \liminf_{k \to \infty} \frac{H \left(\alpha_{m_k}^{(n_k)}\right)}{n_k} \\
&\leq \liminf_{k \to \infty} \log(n_k^2 m_k) (\Lambda_{m_k} + 1/n_k) +\frac{\log(n_k+1)}{n_k} \\
&\leq \liminf_{k \to \infty} 2 \log(n_k) \Lambda_{m_k} + \log(m_k) \Lambda_{m_k} + 2\log(n_k) / n_k\\
& + \log(m_k) / n_k + \log(n_k+1)/n_k.
\end{align*}

By our choice of $m_k$ and $n_k$ we have $\log(m_k) \Lambda_{m_k} \to C$ and all other terms tend to 0:
\[
\lim_{k \to \infty} \log(n_k) \Lambda_{m_k} = \lim_{k \to \infty} \frac{\log(n_k)}{\log(m_k)} \log(m_k) \Lambda_{m_k}= \left(\lim_{k \to \infty}\frac{\log(n_k)}{\log(m_k)}\right) C = 0,
\quad
\]
and $\lim_{k \to \infty} \frac{\log(m_k)}{n_k} = 0$. Hence $h(f,\alpha) \leq C$.
\end{proof}



\begin{tcbexercise}{}{}
Find an example of an infinite interval exchange transformation such that the estimate of Theorem~\ref{thm:EntropyInfiniteIET} is positive but the entropy is null.
\end{tcbexercise}

\begin{tcbexample}{}{EntropyBakerSurface}
We apply the estimates given in Theorem~\ref{thm:EntropyInfiniteIET} to baker's surfaces. Let $\alpha\in(0,1)$ be fixed and $I\subset B_\alpha$ be a horizontal segment of length $\frac{\alpha}{1-\alpha}$ with no singularities in its interior. Consider $f:I\to I$ the infinite IET defined by the first return map to $I$ of the (vertical) translation flow on $B_\alpha$. The lengths of the corresponding subintervals are $\{\alpha^i\}_{i\geq 1}$, thus $\Lambda_m=\frac{\alpha^{m+1}}{1-\alpha}$. Therefore
$$
h(f)\leq\liminf \frac{\alpha^{m+1}}{1-\alpha}\log(m)=0.
$$
\end{tcbexample}

\begin{tcbexercise}{}{}
    Show that the infinite interval exchange transformation $f_{\alpha,\theta}$ arising from in Example~\ref{exa:BakerIET} by the translation flow in direction $\theta$ on baker's surface $B_\alpha$ has entropy zero.
\end{tcbexercise}

\begin{tcbexample}{}{exa:ItalianBilliard}
In Section~\ref{ssec:Coding} we discussed first return maps on a translation surface $M(P)$ associated to an infinite step billiard table $P$. Recall from Example~\ref{exa:InfiniteStepSurface} that $P=\cup_{n\in\N}[n-1,n]\times [0,h_n]$ where $(h_n)_{n\in\N}$ is a monotonically decreasing sequence of positive numbers. Hence, if $\sigma_x$ and $\sigma_y$ are the reflections of the plane w.r.t. the horizontal and vertical lines passing through the origin respectively, then $M(P)$ is tiled by $\{P,\sigma_xP,\sigma_y\sigma_xP,\sigma_y P\}$. Moreover, if the sequence $(h_n)_{n\in\N}$ is summable then $M(P)$ has finite area and $\gamma:=(\cup_{n\geq 1}\{n\}\times[-h_{n+1},h_{n+1}])\sqcup(\cup_{n\geq 1}\{-n\}\times[-h_{n+1},h_{n+1}])$ is a finite length disconnected smooth curve which is transversal to the translation flow on $M(P)$ in any direction different from the vertical direction. By the formulae following equation (\ref{eq:IETGeneralizedMalaga}), if $\theta\in(0,\frac{\pi}{2})$, the first return map $f_{\gamma,\theta}:\gamma\to\gamma$ of the translation flow $F_\theta^t$ to $\gamma$ defines a generalized infinite IET. We use the following general result and Theorem~\ref{thm:EntropyInfiniteIET} to show that many of these infinite IETs have zero entropy.

\begin{tcblemma}{}{}
Let $\lambda_1 \geq \lambda_2 \geq \lambda_3 \geq \ldots$ be summable positive numbers and
$\Lambda_m = \lambda_{m+1} + \lambda_{m+2} + \ldots$.
Let $\{J_1, J_2, J_3, \ldots\}$ be a partition of $\{1,2,3,\ldots\}$ and $K > 0$ a constant such that
$|J_i| \leq K$ for all $i$.
Denote
$\lambda'_i := \sum_{j \in J_i} \lambda_j$ and assume that $\lambda'_1 \geq \lambda'_2 \geq \ldots$.
Denote also $\Lambda'_m := \lambda'_{m+1} + \lambda'_{m+2} + \ldots$
and $d_i := \min \{j \in \N: K \lambda'_j < \lambda'_i\}$. Then
\begin{equation}
	\label{EQ:EstimateEntropy}
\liminf_{m \to \infty} \log(m) \Lambda_m
\leq
\liminf_{i \to \infty} \log(d_i) (K\ i\ \lambda'_i + \Lambda'_i).
\end{equation}
\end{tcblemma}

\begin{proof}
For each $i$ let us define $m_i$ as the minimum number $m$ such that
for all $j \leq i$, $J_j \cap \{1, \ldots, m\} \not= \emptyset$.

We claim that we have the two following inequalities
\begin{enumerate}
\item $m_i \leq K d_i$,
\item $\Lambda_{m_i} \leq K \ i \ \lambda'_i + \Lambda'_i$.
\end{enumerate}
The conclusion of the lemma follows directly from the above two inequalities.

Let us prove the first inequality.
Because each $J_i$ contains at most $K$ elements, $m_i$ is bounded by the
number of elements $j$ such that $\lambda_j \geq \lambda'_i / K$.
Now, among the $K d_i$ first numbers $j$ there exists
necessarily one that belongs to a partition $J_n$ with $n \geq d_i$.
Hence $\lambda_{K d_i} \leq \lambda'_{d_i} < \lambda'_i / K$ (where the last
inequality holds by definition of $d_i$). This shows the first inequality.

We now prove the second inequality.
We consider the intervals
$\lambda_m$ with $m > m_i$ and decompose them into two subsets:
the one for which $m$ belongs to $J_1 \cup J_2 \cup \ldots \cup J_i$ and
the rest. On the first hand, each subset $J_j$ contains at most $K$ elements.
Hence the sum of the elements of the first family is bounded by
$K \cdot i \cdot \lambda_i$. On the other hand, the other sum is bounded
by $\Lambda'_i$. This terminates the proof of the lemma.
\end{proof}

In the case of the step billiard we can take $\{\lambda_i\}$ as the sequence $\{\lambda_{A_i},\lambda_{D_i}\}\cup\{\lambda_{B_i},\lambda_{C_i}\}$ defined by the formulae following equation (\ref{eq:IETGeneralizedMalaga}) and $J_i$ is just the partition induced by considering all elements in a pair of transversal segments in $\gamma$ of the form $\{\pm n\}\cup[-h_n,h_n]$. Hence we have that $K = 8$. Moreover, if the size of the steps $h_n$ defining the infinite step billiard table $P$ decays exponentially fast then $d_i=i+1$. From $(\ref{EQ:EstimateEntropy})$ we obtain that $\liminf_{m \to \infty} \log(m) \Lambda_m=0$ and from Theorem~\ref{thm:EntropyInfiniteIET} that the first return map $f_{\gamma_\theta}$ has entropy zero. Remark that this fact depends only on the exponential decay of the sequence $(h_n)_{n\in\N}$ and not on the direction $\theta$.
\end{tcbexample}

\subsubsection{Abramov's formula} Let $M$ be a finite area translation surface and suppose that for a fixed direction $\theta$ the translation flow $F_\theta^t$ is defined almost everywhere on $M$. The entropy of this flow is by definition (see~\cite{CFS82}) the (metric) entropy of the time-1 map $F_\theta^1$ w.r.t. the Lebesgue measure on $M$. For example, if $M$ is obtained by suspending the IET from Example~\ref{exa:PositiveEntropyIET}, then the entropy of the (vertical) translation flow on $M$ is $log(2)$. On the other hand, the entropy of the translation flow on baker's surface is equal to zero.

Abramov's formula (see~\cite{Abramov661} and \cite{Abramov662}) which we describe below, relates the entropy of the translation flow to the entropy of the first return map to an interval (globally) transverse to it. More precisely, let $\gamma:I\to M$ be a smooth curve transverse to $F_\theta^t$. For simplicity let $\gamma$ also denote the image of this map and $dx$ the push-forward of Lebesgue measure of the interval $I$ to $\gamma$.  We suppose that $\gamma$ intersects every orbit of $F_\theta^t$. Since the area of $M$ is finite the first return time $r:\gamma\to(0,+\infty]$  of $F_\theta^t$ to $\gamma$ is defined almost everywhere on $\gamma$ as well as the first-return map $f:\gamma\to\gamma$. Abramov's formula states that:
\begin{equation}
	\label{eq:Abramov}
h(F_\theta^1,Leb_M)=\frac{h(f,dx)}{\int_\gamma rdx}
\end{equation}
Since $M$ is of finite area and the transversal curve $\gamma$ is global we have that $\int_\gamma rdx=\Area(M)$.

\appendix

\chapter{The hyperbolic plane and Fuchsian groups}
\label{app:FuchsianGroups}
\label{Appendix:FuchsianGroups}

Let us recall from Section~\ref{ssec:Structures} that translation surfaces admit a $\SL(2,\R)$-action. The quotient $\PSL(2,\R) = \SL(2,\R) / \{\pm 1\}$ identifies to the isometry group of the hyperbolic plane. This relation with geometry provides a very efficient tool to study $\SL(2,\R)$ and its discrete subgroups. In this appendix we recall classical definitions and results about hyperbolic geometry, discrete subgroups of $\PSL(2,\R)$ and Rosen continued fractions. For most of the proofs we refer to the literature such as~\cite{Katok92} or~\cite{Beardon}.

\section{Fuchsian groups}

Let $\hat{\C} = \C \cup \{\infty\}$ denote the Riemann sphere. It admits an action of $\PSL(2,\C)$ by \emphdef[homography]{homographies} also called \emphdef[M\"obius transformation]{M\"obius transformations}
\begin{equation}
    \label{eq:moebius-transformations}
\begin{pmatrix}
a & b \\
c & d
\end{pmatrix}
\cdot
z
= \frac{az + b}{cz + d}.
\end{equation}
It is a standard result that the M\"obius transformations are the only analytic automorphisms of the Riemann sphere $\hat{\C}$, see for example Theorem 1.8.2 in~\cite{Hubbard-book1}. The M\"obius transformations stabilizing respectively the affine plane $\C$ and the upper half space $\H^2 = \{z = x + iy \in \C : y > 0\}$ are respectively the subgroups $\Aff(\C) = \left\{\begin{pmatrix}a&b\\0&1\end{pmatrix}\right\}$ and $\PSL(2,\R)$. These two groups are all analytic automorphisms of $\C$ and $\H^2$ respectively, see also Theorem 1.8.2 in~\cite{Hubbard-book1}. The Poincar\'e-Koebe's uniformization theorem states that these are the only three simply connected Riemann surfaces up to analytic automorphisms.

In this section we focus on the upper half plane $\H^2$. It is somehow the richest of the three simply connected Riemann surfaces as it is the universal cover of any Riemann surface which is not homeomorphic to a sphere, a torus, a plane or a cylinder. It can be endowed with the \emphdef{Poincar\'e metric} or \emphdef{hyperbolic metric}:
\[
ds^2 = \frac{dx^2+dy^2}{y^2}.
\]
The geodesics in the metric space $(\H^2, ds)$ are the half circles and lines perpendicular to the boundary $\partial \H^2 := \R \cup \{\infty\}$, see Theorem 1.2.1 in~\cite{Katok92}. This metric space is uniquely geodesic: by every pair of points of $\H^2$ passes a unique geodesic. This property also extends to $\H^2 \cup \partial \H^2$ and a geodesic is entirely determined by its two endpoints in $\partial \H^2$. The group $\PSL(2,\R)$ acts by isometries with respect to this metric, see Theorem 1.1.2 in~\cite{Katok92}.

Elements of $\PSL(2,\R)$ are classified with respect to their trace as follows. We say that a matrix
$A \in \PSL(2,\R)$ different from $\pm 1$ is
\begin{itemize}
\item \emphdef[elliptic (matrix)]{elliptic} if $|\tr(A)|<2$.
\item \emphdef[parabolic (matrix)]{parabolic} if $|\tr(A)|=2$.
\item \emphdef[hyperbolic (matrix)]{hyperbolic} if $|\tr(A)|>2$.
\end{itemize}
As explained in Chapter 2 in~\cite{Katok92}, every elliptic element fixes exactly one point in $\H^2$, every parabolic fixes exactly one point in $\partial \H^2$ and every hyperbolic fixes exactly two points $\partial \H^2$. The \emphdef[axis (hyperbolic transformation)]{axis} of a hyperbolic element is the unique geodesic in $\H^2$ joining the corresponding fixed points.

A \emphdef{Fuchsian group} is a discrete subgroup of $\PSL(2,\R)$. Fuchsian groups can be characterized in more geometric terms using the following notion.
\begin{tcbdefinition}{}{PropDiscFreeAct}
\label{def:ProperlyDiscontinuousAction}
Let $X$ be a locally compact topological space and $\Gamma$ a discrete group acting on $X$.
The action is said to be \emphdef[properly discontinuous (action)]{properly discontinuous} if for every compact $K\subset X$ the set
\[
\{g \in \Gamma: g K \cap K \not= \emptyset\}
\]
is finite.
\end{tcbdefinition}

\begin{tcbtheorem}{\cite[Theorem~2.2.6 and Corollary~2.2.7]{Katok92}}{FuchsianGroupCharacterizations}
Let $\Gamma$ be a subgroup of $\PSL(2,\R)$. Then the following are equivalent
\begin{itemize}
\item $\Gamma$ is Fuchsian
\item $\Gamma$ acts properly discontinuously on $\H^2$,
\item for all $z\in\H^2$, the orbit $\Gamma \cdot z$ is discrete.
\end{itemize}
\end{tcbtheorem}

Let $\Gamma$ be a Fuchsian group. By the discretness of item 3 in Theorem~\ref{thm:FuchsianGroupCharacterizations}, for a point $z \in \H^2 \cup \partial \H^2$ the set of accumulation points of its orbit $\Gamma \cdot z$ in $\H^2 \cup \partial \H^2$ is contained in $\partial \H^2$. The latter is independent of the point $z \in \H^2 \cup \partial \H^2$ and is called \emph{the limit set of $\Gamma$} and denoted $\Lambda(\Gamma)$.

As explained in~\cite[Chapter~3.4]{Katok92}, there is the following trichotomy for the limit set of a Fuchsian group $\Gamma$
\begin{enumerate}
\item either $\Lambda(\Gamma)$ contains at most two points,
\item or it is homeomorphic to the Cantor set,
\item or it is equal to $\R\cup\{\infty\}$.
\end{enumerate}
In the first case we say that $\Gamma$ is \emphdef[elementary (Fuchsian group)]{elementary}. An elementary Fuchsian group either stabilizes a point in $\H^2$ (in which case $\Gamma$ consists of elliptic elements) or a single point in $\partial \H^2$ (in which case $\Gamma$ consists of parabolic elements) or a pair of points in $\partial \H^2$ (in which case $\Gamma$ consists of hyperbolic elements and possibly an elliptic involution). Fuchsian groups whose limit set is nowhere dense in $\R\cup\{\infty\}$, \ie cases (1) and (2), are said to be of the \emphdef[second kind (Fuchsian group)]{second kind} and those satifying (3) are of the \emphdef[first kind (Fuchsian group)]{first kind}.

A Fuchsian group for which the quotient $\Gamma \backslash \H^2$ has finite area is called a \emphdef[lattice (Fuchsian group)]{lattice} and, if in addition this su is compact the lattice is \emphdef[uniform (lattice)]{uniform}. A Fuchsian group is a lattice if and only if it is finitely generated and of the first kind. The following Theorem gives several characterizations for elementary Fuchsian groups.
\begin{tcbtheorem}{}{ElementaryFuchsianGroup}
Let $\Gamma$ be a Fuchsian group. Then the following are equivalent:
\begin{enumerate}
\item $\Lambda(\Gamma)$ is infinite,
\item $\Gamma$ contains a free group on two generators,
\item there is a pair of hyperbolic elements in $\Gamma$ whose sets of fixed points in  $\partial \H^2$ are disjoint,
\item there is a pair of hyperbolic elements in $\Gamma$ whose axes are disjoint.
\end{enumerate}
\end{tcbtheorem}
A proof of this Theorem can be deduced, for example, from the discussion on non-elementary Fuchsian groups presented in~\cite{Beardon}.

The following result is used in Sections~\ref{ssec:HooperWeissTheorem} and~\ref{ssec:SymmetriesWindTree} to show that Fuchsian groups of the first kind appear as Veech groups of wind-tree models.
\begin{tcblemma}{}{NormalSubgroupLimitSet}
Let $\Gamma$ be a non-elementary Fuchsian group and let $\Gamma'$ be a non-trivial normal subgroup of $\Gamma$. Then $\Gamma'$ is non-elementary and the limit sets of $\Gamma$ and $\Gamma'$ coincide.
\end{tcblemma}
\begin{proof}
Let us first show that $\Gamma'$ is non-elementary. Since it is non-trivial it
contains an element $h$ which is not the identity. Then, for an appropriate $g \in \Gamma$
the group $\langle h, g h g^{-1} \rangle$ is non-elementary. In particular $\Gamma'$
is non-elementary.

For the rest of the proof we use the following two facts, which can be found in~\cite{Katok92}: every non-elementary Fuchsian group contains infinitely many hyperbolic elements, no two of which have a common fixed point, and every limit set $\Lambda(\Gamma)$
containing more than one element is the closure of the set of fixed points of the hyperbolic elements in $\Gamma$.
Let then $g\in\Gamma$ be hyperbolic. It is sufficient to show that the attractive fixed
point of $g$ is accumulated by fixed points of elements in $\Gamma'$.
Let us consider any hyperbolic element $h \in \Gamma'$. If its attractive fixed point $x$
is fixed by $g$, then there is nothing to do. Otherwise, the attractive fixed
point of the hyperbolic element $g^n h g^{-n}\in\Gamma'$ is $g^n x$ which converges to the attractive fixed point of $g$.
\end{proof}



We conclude this section with a particular version of the ping-pong lemma (\emph{a.k.a.} the table-tennis lemma), adapted to the context of hyperbolic elements acting on $\H^2$. We refer the reader to II.B in~\cite{delaHarpe00} for more details about the proof.
\begin{tcbtheorem}{Ping-pong lemma}{PingPongLemma}
Let $G$ be a group acting on a topological space $X$.
Let $g_1$, $g_2$, \ldots, $g_k$ be elements of $G$ and
assume that there exist disjoint non-empty open subsets
$X^+_1, X^-_1, X^+_2, X^-_2, \ldots, X^+_k, X^-_k \subset X$ such that for all
$i$ in $\{1, \ldots, k\}$ we have
\[
h_i (X \setminus X_i^-) \subset X^+_i
\qquad \text{and} \qquad
h_i^{-1} (X \setminus X_i^+) \subset X^-_i
\]
Then
\begin{itemize}
\item all the $g_i$ have infinite order
\item the subgroup of $G$ generated by the $g_i$ is canonically isomorphic to the free product $\langle g_1
\rangle * \langle g_2 \rangle * \cdots * \langle g_k \rangle$.
\item the action of $\langle g_1, g_2, \ldots, g_k \rangle$ on $X$ is free and properly discontinuous.
\end{itemize}
\end{tcbtheorem}

As a consequence of the preceding discussion and the ping-pong lemma we have the following:
\begin{tcbtheorem}{}{PingPongCorollary}
Let $g,h$ be two hyperbolic elements in $\Isom(\H^2)$ such that $\Fix(g) \cap \Fix(h) = \emptyset$. Then, for large enough integer $m$ the group generated by $g^m$ and $h$ is free and all its elements different from the identity are hyperbolic.
\end{tcbtheorem}

\section{The group $G_\lambda$}
\label{APPENDIX:RosenContinuedFractions}
\label{appendix:sec:RosenGLambda}
In this section we investigate a family of Fuchsian group that we denote $G_\lambda$ where $\lambda > 2$ is a parameter. These groups play a central role in finite and infinite Thurston-Veech constructions, for example in the proofs of Section~\ref{ssec:VeechGroupBakersSurface}.

Our presentation mainly follows the original article of D.~Rosen~\cite{Rosen54}.
The latter is a continuation of earlier work of Hecke~\cite{Hecke1936}.

For $\lambda > 2$ let
\[
h_\lambda = \begin{pmatrix}1&\lambda\\0&1\end{pmatrix}
\qquad \text{and} \qquad
v_\lambda = \begin{pmatrix}1&0\\\lambda&1\end{pmatrix}.
\]
We define $G_\lambda := \left\langle h_\lambda, v_\lambda \right\rangle$. The main two results about the groups $G_\lambda$ are Theorem~\ref{thm:GLambdaGeometry} and Theorem~\ref{thm:GLambdaLimitSet} that we state now.
\begin{tcbtheorem}{}{GLambdaGeometry}
Let $\lambda > 2$ and $G_\lambda$ be as above.
Then $G_\lambda$ is a free group generated by $h_\lambda$ and $v_\lambda$ which is discrete and of infinite covolume in $\PSL(2,\R)$. The quotient $SO(2) \backslash \PSL(2,\R) / G_\lambda$ is a hyperbolic pair of pants with boundary lengths $0$, $0$ and $2 \operatorname{acosh}\left( \frac{\lambda^2-2}{2} \right)$. In other words, the quotient is a sphere with two cusps and a funnel.
\end{tcbtheorem}
Theorem~\ref{thm:GLambdaGeometry} is proven in Section~\ref{ssec:GLambdaFundamentalDomain} below. See also Figure A.1 for a picture of the fundamental domain of the action of $G_\lambda$ on $\H^2$.
\begin{tcbtheorem}{}{GLambdaLimitSet}
Let $\lambda > 2$ and $G_\lambda$ as above. Let $I = (r,-r+\lambda)$
where $r := \frac{\lambda + \sqrt{\lambda^2 - 4}}{2}$.
Then the intervals $\{g \cdot I\}_{g \in G_\lambda}$ in $\mathbb{P}^1 \R$ are disjoint and
the limit set $\Lambda(G_\lambda)$ of $G_\lambda$ is a Cantor set whose complement
is the union $\bigcup_{g \in G_\lambda} g \cdot I$.

For $\lambda = 2$, we have $\Lambda(G_2) = \R \cup \{\infty\}$.
\end{tcbtheorem}
Theorem~\ref{thm:GLambdaLimitSet} is proven in Section~\ref{ssec:Rosen}. See also Figure A.2 for a picture illustrating the gap $(r,-r+\lambda)$ in the limit set of $G_\lambda$.
\begin{tcbremark}{}{}
In the limit situation $\lambda=2$, $G_2 < \SL(2,\Z)$ is still a free group generated
by $h_2$ and $v_2$ which is discrete. However it has finite index in $\SL(2,\Z)$
and hence finite covolume. The quotient is a hyperbolic pair of pants with three cusps
(or rather "the" hyperbolic pair of pants with three cusps since it is unique up
to biholomorphism).

It was proven by Hecke that for $\lambda < 2$ the group $G_\lambda$ is discrete if and only if
$\lambda = 2 \cos(\pi/n)$ for some integer $n \geq 3$. The group $G_{2 \cos(\pi/n)}$
has finite covolume in $\SL(2,\R)$ and has for unique relation
\[
\left(h_{2 \cos(\pi/n)} (v_{2 \cos(pi/n)})^{-1}\right)^n = (-1)^n.
\]
The quotient is a pair of pants with two cusps and an orbifold point with angle $2\pi/n$.
\end{tcbremark}

\subsection{A fundamental domain for $G_\lambda$}
\label{ssec:GLambdaFundamentalDomain}
In this section we prove Theorem~\ref{thm:GLambdaGeometry} by finding an explicit fundamental domain for the action of $G_\lambda$. A useful ingredient in the proof is the ping-pong lemma (Theorem~\ref{thm:PingPongLemma}).
\begin{figure}%
\begin{center}%
\includegraphics{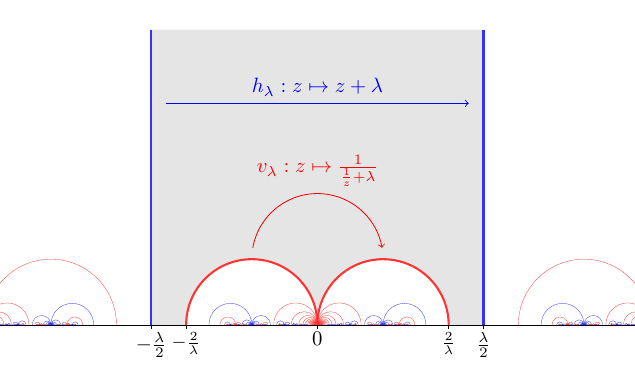}%
\caption{A fundamental domain for $G_\lambda$ (in gray with bold boundaries)
in the upper half plane model. To form the quotient
surface, the vertical blue geodesics and red geodesics on the boundary of the
fundamental domain are identified with respectively $h_\lambda$ and $v_\lambda$.
According to Theorem~\ref{thm:GLambdaGeometry} this quotient is a pair of pants
(a three punctured sphere).
The thinner lines are translates by $G_\lambda$ of the boundary of the fundamental domain.}
\end{center}
\label{fig:GLambdaFundamentalDomain}
\end{figure}

\begin{tcbtheorem}{}{GLambdaFundamentalDomain}
Let $\lambda > 2$. Then $G_\lambda$ is a free group generated by $h_\lambda$ and $v_\lambda$.
Let us define
\begin{align*}
X_1^+ &= \{\Re(z) > \lambda/2\} \\
X_1^- &= \{\Re(z) < -\lambda/2\} \\
X_2^+ &= \{|z - 1/\lambda| < 1/\lambda\} \\
X_2^- &= \{|z + 1/\lambda| < 1/\lambda\}
\end{align*}
Then $D_\lambda := \H^2 \setminus (X_1^+ \cup X_1^- \cup X_2^+ \cup X_2^-)$ is a fundamental region for the action of $G_\lambda$ on $\H^2$.
\end{tcbtheorem}

\begin{proof}
Recall that $h_\lambda \cdot z = z + \lambda$ and $v_\lambda \cdot z = \frac{1}{\lambda + \frac{1}{z}}$.
These are two parabolic elements that stabilize respectively $\infty$ and $0$.

The subsets $X_1^+, X_1^-, X_2^+, X_2^-$ do not exactly satisfy the assumption of Theorem~\ref{thm:PingPongLemma} for $g_1=h_\lambda$ and $g_2=v_\lambda$ but almost. The only problem comes from the geodesic boundaries $\{\Re(z) = \lambda/2\}$, $\{\Re(z)=-\lambda/2\}$, $\{|z-1/\lambda|=1/\lambda\}$ and $\{|z+1\/\lambda|=1/\lambda\}$. More precisely
\begin{align*}
h_\lambda \cdot \{\Re(z) = -\lambda/2\} &= \{\Re(z) = \lambda/2\} \\
v_\lambda \cdot \{|z+1/\lambda|=1/\lambda\} &= \{|z-1/\lambda|=1/\lambda\}.
\end{align*}
Indeed, $h_\lambda \cdot \frac{-\lambda}{2} = \frac{\lambda}{2}$ and $v_\lambda \cdot \frac{-2}{\lambda} = \frac{2}{\lambda}$. For that purpose, let us extend a bit the domains
\begin{align*}
Y_1^+ &= \{\Re(z) > \lambda/2 - \epsilon\} \\
Y_1^- &= \{\Re(z) < -\lambda/2 + \epsilon\} \\
Y_2^+ &= \{|z - (1/\lambda +\epsilon)| < 1/\lambda + \epsilon\} \\
Y_2^- &= \{|z + (1/\lambda + \epsilon) | < 1/\lambda + \epsilon\}
\end{align*}
where $\epsilon$ is a positive real such that $2 \epsilon < |\lambda/2 - 2/\lambda|$.
The sets $Y_1^+$, $Y_1^-$, $Y_2^+$ and $Y_2^-$ do satisfy the assumption of Theorem~\ref{thm:PingPongLemma} and this shows that $G_\lambda$ is a free group generated by $h_\lambda$ and $v_\lambda$.

We now show that $D_\lambda$ is a fundamental domain for $G_\lambda$.

Let us first show that no two elements of the interior $\mathring{D}_\lambda$ of $D_\lambda$ are congruent. This follows
from the ping-pong property : for any $g \in G_\lambda \setminus \{1\}$, we have $g \cdot \mathring{D}_\lambda \subset (X_1^+ \cup X_1^- \cup X_2^+ \cup X_2^-)$.

Now, let $z \in \H^2$. We construct an element $g \in G_\lambda$ such that $g z \in D_\lambda$.
For that purpose, we consider the function $f(z) = \dist(z, D_\lambda)$ (where the distance is with respect to the hyperbolic metric).
Our first aim is to show that for all $z \in \H^2$ there exists $g \in G_\lambda$ such that $f(g \cdot z) = 0$.
Let us first notice that the band $\H^2 \setminus (X_1 \cup X_1^-)$ is a fundamental domain for $\langle h_\lambda \rangle$.
In other words, if $z \in X_1^+ \cup X_1^-$ there exists $n \not= 0$ such that $h^n_\lambda \cdot z \not\in (X_1^+ \cup X_1^-)$.
For this particular $n$ we have $f(h^n_\lambda \cdot z) < f(z)$.
Similarly, if $z \in X_2^+ \cup X_2^-$, we can play the same game with $v_\lambda$.
This shows that whenever $z \not\in D_\lambda$ one can find an element in $\langle h_\lambda \rangle \cup \langle v_\lambda \rangle$
that makes the distance to $D_\lambda$ decrease. Since $G_\lambda$ acts freely and properly discontinuously this process
must stop after finitely many steps. But when this stops, the point must be in $G_\lambda$.
\end{proof}

\begin{proof}[Proof of Theorem~\ref{thm:GLambdaGeometry}]
By Theorem~\ref{thm:GLambdaFundamentalDomain}, the quotient $G_\lambda \backslash \H^2$ is biholomorphic
to the domain $D_\lambda$ quotiented by the two identifications of the boundaries $\{\Re(z)=-\lambda/2\}$
and $\{Re(z)=\lambda/2\}$ by $h_\lambda$ and the boundaries $\{|z+1/\lambda|=1/\lambda\}$ and
$\{|z-1/\lambda|=1/\lambda\}$ by $v_\lambda$.

The quotient is indeed a pair of pants with two cusps (corresponding to the points at infinity $0$ and $\infty$) and a funnel (corresponding to the
boundary segments $[-\lambda/2, -2/\lambda]$ and $[2/\lambda,\lambda/2]$. The boundary length of the funnel is the hyperbolic length of the simple geodesic
turning once around the funnel which corresponds to the conjugacy class of $g = h_\lambda (v_\lambda)^{-1}$. The hyperbolic
length $\ell$ of this curve is equal to half the leading eigenvalue. This can be deduced from the trace of the matrix $g$ via the relation
$\tr(g)/2 = \cosh(\ell/2)$. Indeed if $z$ is on the axis of a hyperbolic element $g \in \PSL(2,\R)$
then the distance between $z$ and $g \cdot z$ is $\ell$ such that $\tr(g)/2 = \cosh(\ell/2)$.
\end{proof}

\subsection{Rosen continued fraction}
\label{ssec:Rosen}
We now describe the Rosen continued fraction which is the main ingredient
in our proof of Theorem~\ref{thm:GLambdaLimitSet}. The Rosen map
is a map defined on an interval which corresponds to the action of certain elements
of $G_\lambda$ on the line at infinity $\partial \H^2$ in the upper half plane. Passing
from the Rosen continued
fraction to $G_\lambda$ will be done in the next section. The construction depends on
$\lambda > 2$ but we omit this dependency in our notations.

\begin{figure}%
\begin{center}%
\includegraphics{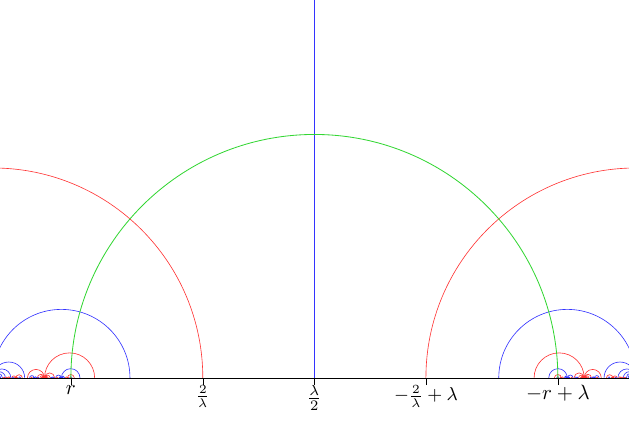}%
\caption{A zoom in the region $(r, -r+\lambda)$ that does not contain
any limit point of $G_\lambda$, see Theorem~\ref{thm:GLambdaLimitSet}.
The geodesic from $r$ to $-r+\lambda$ (in green) is a periodic geodesic on the
quotient which corresponds to the simple curve around the funnel.}
\end{center}
\label{fig:GLambdaLimitSet}
\end{figure}

Let $T: \left[-\frac{2}{\lambda}, \frac{2}{\lambda}\right] \to \left[-\frac{\lambda}{2}, \frac{\lambda}{2}\right]$ be the map defined
by
\[
T(x) =
\left\{ \begin{array}{ll}
\frac{1}{\epsilon(x) \cdot x} - \lambda
\left\lfloor \frac{1}{\lambda \epsilon(x) \cdot x} + \frac{1}{2} \right\rfloor
& \text{if $x \not= 0$} \\
0 & \text{if $x = 0$}.
\end{array} \right.
\]
where $\epsilon(x)$ is the sign of $x$.
For $x \not= 0$, we also note $\displaystyle a(x) = \left\lfloor \frac{1}{\lambda \epsilon(x) x} + \frac{1}{2} \right\rfloor \in \{1, 2, \ldots, +\infty\}$.

Note that $T$ has a codomain larger than its domain as $2/\lambda < 1 < \lambda/2$. We aim to study orbits under the map $T$ that is the list $(x, T(x), T^2(x), \ldots)$. Because the codomain is smaller than the domain, some orbits are finite and can not be extended further.
Let $\Lambda(T) \subset \left[-\frac{2}{\lambda}; \frac{2}{\lambda}\right]$ be the subsset of
points with infinite orbit, that is
\[
\Lambda(T) := \bigcap_{n \ge 0} \Lambda_n(T)
\quad \text{where} \quad
\Lambda_n(T) := T^{-n}
\left(
\left[-\frac{2}{\lambda}, \frac{2}{\lambda}\right]
 \right).
\]
Below, we give a description of $\Lambda(T)$ and explore the link with the limit set of $G_\lambda$.

For $n \geq 1$ and $x \in \Lambda_{n-1}(T)$, we define $a_n(x) = a(T^{n-1} x)$ and $\epsilon_n(x) = \epsilon(T^{n-1} x)$.
The \emph{Rosen continued fraction expansion} of $x \in \Lambda(T)$ is the infinite sequence of pairs $(\epsilon_n(x), a_n(x))_{n \ge 1}$.

\begin{figure}[!ht]
\begin{center}\includegraphics{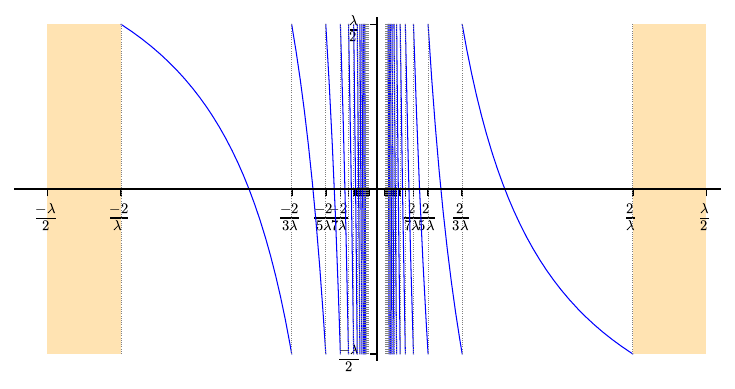}\end{center}
\caption{Graph of the Rosen map for $\lambda = 2.27$. The map is not defined
on $[\frac{-\lambda}{2}, \frac{-2}{\lambda}] \cup [\frac{2}{\lambda}, \frac{\lambda}{2}]$ (yellow region).}
\end{figure}

Recall that (invertible) $2 \times 2$ real matrices act by homography on the $\H^2$ and its
boundary $\R \cup \{\infty\}$ via
\[
\begin{pmatrix}a&b\\c&d\end{pmatrix} \cdot z = \frac{az + b}{cz + d}.
\]
Matrix multiplication corresponds to function composition. The only matrices that
act trivially are $\pm Id$ and this action of $\PSL(2,\R)$ on $\R \cup \{\infty\}$
corresponds to the action of $\PSL(2,\R)$ the boundary of the Poincar\'e half-plane. With
that notation in mind, we have for $x \in [-2/\lambda; 2/\lambda] \setminus \{0\}$ the expression
\[
T(x) =
\begin{pmatrix}
-\lambda r(x) & \epsilon \\
1 & 0
\end{pmatrix}
\cdot x
=
\begin{pmatrix}
0 & \epsilon(x) \\
1 & \lambda a(x)
\end{pmatrix}^{-1}
\cdot x.
\]
In the above formula, there is a slight abuse since the determinant of the matrix is $-\epsilon(x)$ and hence not an element of $\SL(2,\R)$ if $x > 0$. In order to make sense
we consider the action of the larger group $\SL^{\pm 1}(2,\R)$ on $\partial \H^2$.

We see that the inverse branches of $T$ are parametrized by the pairs $(\epsilon_1,a_1) \in \{\pm 1\} \times \Z_{> 0}$. Given such a pair, we define the corresponding inverse branch $H_{\epsilon_1:a_1}: \left[-\frac{\lambda}{2}; \frac{\lambda}{2} \right) \to \left[-\frac{\lambda}{2}; \frac{\lambda}{2} \right)$ by
\begin{equation} \label{eq:matrixHEpsilonA}
H_{\epsilon_1:a_1}(x) =
\begin{pmatrix}
0 & \epsilon_1 \\
1 & \lambda a_1
\end{pmatrix}
\cdot x
=
\frac{\epsilon_1}{\lambda a_1 + x}.
\end{equation}
Then a direct computation gives the following lemma.
\begin{tcblemma}{}{}
For any pair $(\epsilon_1,a_1) \in \{\pm 1\} \times \Z_{> 0}$, the map $H_{\epsilon_1:a_1}$ is a bijection on its image
\[
I_{\epsilon_1:a_1}
:=
H_{\epsilon_1:a_1}\left( \left( -\frac{\lambda}{2}, \frac{\lambda}{2} \right] \right)
=
\epsilon_1 \cdot \left( \frac{2}{(2a_1+1) \lambda}; \frac{2}{(2a_1-1) \lambda} \right]
\]
Furthermore, for $x \in \Lambda(T)$ we have $(\epsilon_1(x), a_1(x)) = (\epsilon_1, a_1)$ if and only if $x \in I_{\epsilon_1:a_1}$.
\end{tcblemma}

By iterating the inverse branches, we obtain a complete characterization of $\Lambda(T)$.
Before stating the main result we introduce further convenient notation.
We denote $\cS$ the set of finite sequences of elements in $\{\pm 1\} \times \Z_{> 0}$.
For $s = ((\epsilon_1, a_1), \ldots, (\epsilon_n, a_n)) \in \cS$ we denote
\begin{align*}
H_s(x)
&:=
H_{\epsilon_1:a_1} \circ H_{\epsilon_2:a_2} \circ \cdots \circ H_{\epsilon_n:a_n}(x) \\
&=
\begin{pmatrix}
0 & \epsilon_1 \\
1 & \lambda a_1
\end{pmatrix}
\begin{pmatrix}
0 & \epsilon_2 \\
1 & \lambda a_2
\end{pmatrix}
\cdots
\begin{pmatrix}
0 & \epsilon_n \\
1 & \lambda a_n
\end{pmatrix} \cdot x \\
&=
\cfrac{\epsilon_1}{\lambda a_1 + \cfrac{\epsilon_2}{\ldots + \cfrac{\epsilon_n}{\lambda a_n + x}}}.
\end{align*}
We also denote $I_s := H_s([-\lambda/2, \lambda/2))$ the image of $H_s$.
Let $\overline{\cS}$ the set of finite and infinite sequences of elements in $\{\pm 1\} \times \Z_{> 0}$.
\begin{tcblemma}{}{RosenCFBijection}
The continued fraction expansion
\[
\begin{array}{lll}
\Lambda(T) & \to & \overline{\cS} \\
x & \mapsto & ((\epsilon_n(x), a_n(x)))_{n \ge 1}
\end{array}
\]
is a bijection.
\end{tcblemma}

In order to prove the above lemma, we compute the diameters
of the intervals $I_s = H_s([-\lambda/2, \lambda/2))$. Given
$s = ((\epsilon_1, a_1), \ldots, (\epsilon_n, a_n)) \in \cS$ we
define
\[
\begin{pmatrix}
p_{n-1}(s) & p_n(s) \\
q_{n-1}(s) & q_n(s)
\end{pmatrix}
=
\begin{pmatrix}
0 & \epsilon_1 \\
1 & \lambda a_1
\end{pmatrix}
\begin{pmatrix}
0 & \epsilon_2 \\
1 & \lambda a_2
\end{pmatrix}
\cdots
\begin{pmatrix}
0 & \epsilon_n \\
1 & \lambda a_n
\end{pmatrix}
\]
or equivalently
\[
H_s(x) = \frac{p_{n-1}(s) x + p_n(s)}{q_{n-1}(s) x + q_n(s)}.
\]
The numbers $p_n(s)$ and $q_n(s)$ can be computed inductively by the following formulas
\[
\begin{array}{l@{;\ }l@{;\ }l}
p_{-1}(s) = 1  & p_0(s) = 0 & p_n(s) = \lambda a_n p_{n-1}(s) + \epsilon_n p_{n-2}(s) \\
q_{-1}(s) = 0  & q_0(s) = 1 & q_n(s) = \lambda a_n q_{n-1}(s) + \epsilon_n q_{n-2}(s)
\end{array} .
\]

\begin{tcblemma}{}{RosenGrowth}
Let $s = ((\epsilon_1, a_1), \ldots, (\epsilon_n, a_n)) \in \cS$.
For all $k \in \{1, \ldots, n\}$, we have
\begin{equation} \label{eq:RosenNaiveInequality}
|q_k(s)| > r \cdot |q_{k-1}(s)|
\end{equation}
where $r = \frac{\lambda + \sqrt{\lambda^2 - 4}}{2}$
and
\begin{equation} \label{eq:RosenSign}
\sign(q_k(s)) = \sign(a_k a_{k-1} \cdots a_1).
\end{equation}
If furthermore, $\sign(a_k a_{k-1}) = \epsilon_k$ or $|a_k| \geq 2$, then
\[
|q_k(s)| > \lambda \cdot |q_{k-1}(s)|.
\]
\end{tcblemma}

\begin{tcbcorollary}{}{RosenGrowthEstimate}
Let $((\epsilon_1, a_1), \ldots, (\epsilon_n, a_n)) \in \cS$.
Let $m$ be the number of indices
$k \in \{1, \ldots, n\}$ such that either $\sign(a_k a_{k-1}) = \epsilon_k$ or $|a_k| \geq 2$. Then we have
\[
q_n(s) > \lambda^{m} \cdot r^{n-m}
\]
where $r = \frac{\lambda + \sqrt{\lambda^2 - 4}}{2}$.
\end{tcbcorollary}

\begin{proof}
We prove the first inequality by induction on $k$.
Since $q_0 = 1$, $q_1 = \lambda a_1$ and $\lambda > r$, the
inequality holds for $k=1$. Assume that it holds for $k-1$.
Then
\begin{align*}
q_k &= \lambda a_k q_{k-1} + q_{k-2} \\
 & \geq \lambda |q_{k-1}| - |q_{k-2}| \\
 & > |q_{k-1}| (\lambda - 1/r)
\end{align*}
where in the last step we used the induction hypothesis.
Since $\lambda - 1/r = r$ this finishes the proof of the first inequality.

We now prove the formula for $\sign(q_k)$. By the previous inequality, we
got that $|\lambda a_k q_{k-1}| > |q_{k-2}|$. Hence
$\sign(q_k) = \sign(a_k q_{k-1})$ from which
$\sign(q_k) = \sign(a_k a_{k-1} \cdots a_1)$ follows.

Next, if $\sign(a_k a_{k-1}) = \epsilon_k$ then by~\eqref{eq:RosenSign} we have
$\sign(a_k q_{k-1}) = \sign(\epsilon_k q_{k-2})$ and hence
$|q_k| = \lambda |a_k| |q_{k-1}| + |q_{k-2}|$.
In particular, it follows that
$|q_k| \geq \lambda |q_{k-1}| + |q_{k-2}| \geq (\lambda + 1/r) |q_{k-1}|$
where we applied the inequality~\eqref{eq:RosenNaiveInequality} to $q_{k-2}$.
Finally, if $|a_k| \geq 2$ then
\[
|q_k| \geq 2 \lambda |q_{k-1}| - |q_{k-2}| \geq (2 \lambda - 1/r) |q_{k-1}|.
\]
Since $2 \lambda - 1/r > \lambda$ we also deduce that $|q_k| \geq \lambda |q_{k-1}|$
in that case.
\end{proof}

\begin{tcblemma}{}{RosenDiameterBound}
Let $s = ((\epsilon_1, a_1), \ldots, (\epsilon_n, a_n)) \in \cS$ and $y \in [-\lambda/2; \lambda/2)$.
Set $x = H_s(y)$. Then
\begin{equation}
x - \frac{p_n(s)}{q_n(s)}
=
\frac{(-1)^n \epsilon_1 \cdots \epsilon_n \cdot y}{q_n(s) (q_{n-1}(s) y + q_n)}.
\end{equation}
In particular
\begin{equation} \label{eq:RosenApproximationInequality}
\left| x - \frac{p_n(s)}{q_n(s)} \right|
\leq
\frac{\lambda}{\sqrt{\lambda^2 - 4}} \cdot \frac{1}{|q_n \cdot q_{n-1}|}.
\end{equation}
\end{tcblemma}

\begin{proof}
Let us use $p_n, p_{n-1}, q_n, q_{n-1}$ for $p_n(s), p_{n-1}(s), q_n(s), q_{n-1}(s)$ respectively.
From
\[
x = \begin{pmatrix}p_{n-1}&p_n\\q_{n-1}&q_n\end{pmatrix} \cdot y
\]
we can write
\[
x - \frac{p_n}{q_n}
=
\frac{p_{n-1} y + p_n}{q_{n-1} y + q_n} - \frac{p_n}{q_n}
=
\frac{ (q_n p_{n-1} - p_n q_{n-1}) y}{q_n (q_{n-1} y + q_n)}.
\]
Now, $q_n p_{n-1} - p_n q_{n-1}$ is the determinant of the matrix representation
of the homography $H_s$ which is a product of $n$ matrices with determinant respectively
$-\epsilon_1$, $-\epsilon_2$, \ldots, $-\epsilon_n$. Hence
\[
q_n p_{n-1} - p_n q_{n-1}
=
(-1)^n \epsilon_1 \cdots \epsilon_n
\]
from which we deduce the equality.

Let us now derive the inequality. By~\eqref{eq:RosenNaiveInequality} from Lemma~\ref{lem:RosenGrowth}, we obtain that
$|q_{n-1} y + q_n| > (r - \lambda/2) |q_{n-1}| = \frac{\sqrt{\lambda^2-4}}{2} |q_{n-1}|$.
The inequality~\eqref{eq:RosenApproximationInequality} follows.
\end{proof}

\begin{proof}[Proof of Lemma~\ref{lem:RosenCFBijection}]
For $x \in \Lambda \setminus \{0\}T$, we have $x = H_{\epsilon_1(x):a_1(x)}(T(x))$. Hence
\[
\Lambda(T) = \{0\} \cup \bigcup_{(\epsilon_1, a_1) \in \{\pm 1\} \times \Z_{> 0}} H_{\epsilon_1:a_1}(\Lambda(T)).
\]
By iterating this equality we see that, for $x \in \Lambda(T)$, either its orbit reaches zero after
finitely many steps and we have $x = H_{\epsilon_1(x):a_1(x)} \circ \cdots \circ H_{\epsilon_n(x):a_n(x)}(0)$
which corresponds to a finite continued fraction expansion. Otherwise, the orbit of $x$ never reaches
$0$ and we obtain an infinite continued fraction expansion $(\epsilon_n(x), a_n(x))_n$. From lemma~\ref{lem:RosenDiameterBound},
the diameter of $I_{\epsilon_1(x):a_1(x), \ldots, \epsilon_n(x):a_n(x)}$ tends to zero as $n \to \infty$.
Hence
\[
\{x\} = \bigcap_{n \ge 1} I_{\epsilon_1(x):a_1(x), \ldots, \epsilon_n(x):a_n(x)}.
\]
In other words $x \mapsto ((\epsilon_n(x), a_n(x)))_{n \ge 1}$ is injective.
\end{proof}

\begin{tcbcorollary}{}{}
Let $I = [-\lambda/2,-2/\lambda) \cup (2/\lambda, \lambda/2]$.
Then the sets $\{h_s(I)\}_{s \in \cS}$ are disjoint and their union is
the complement of $\Lambda(T)$ in $[-\lambda/2, \lambda/2]$.

In particular
$\min \Lambda(T) = -r$ and $\max \Lambda(T) = r$ where
$r := \frac{\lambda + \sqrt{\lambda^2-4}}{2}$.
\end{tcbcorollary}

\subsection{The limit set of $G_\lambda$}
We now prove Theorem~\ref{thm:GLambdaLimitSet}. To do so, we show that the limit set $\Lambda(G_\lambda)$ and the set $\Lambda(T)$ of the Rosen map are intimately related. The Rosen continued fraction expansion of points in $\Lambda(T)$ allows to determine the nature of points in the limit set.

\begin{tcbtheorem}{}{GLimitSetVersusTLimitSet}
Let $\lambda > 2$ and $T_\lambda: \left[-\frac{2}{\lambda}, \frac{2}{\lambda}\right] \to \left[-\frac{\lambda}{2}, \frac{\lambda}{2}\right]$ be the associated Rosen map as in Section~\ref{ssec:Rosen}.
We have $\Lambda(T_\lambda) = \Lambda(G_\lambda) \cap [-\lambda/2; \lambda/2]$.
Furthermore, a point $x \in \Lambda(T_\lambda)$ is
\begin{itemize}
\item stabilized by a parabolic element in $G_\lambda$ if and only if its continued fraction expansion is finite,
\item stabilized by an hyperbolic element in $G_\lambda$ if and only if its continued fraction expansion is periodic.
\end{itemize}
\end{tcbtheorem}

In order to tighten the link between the group $G_\lambda$ and Rosen continued fraction we introduce a slightly
larger group $\widehat{G_\lambda}$ in $\PSL^{\pm}(2,\R)$ where $\PSL^{\pm}(2,\R)$ denotes $2\times2$ matrices with
determinant $\pm 1$. The group $\PSL^{\pm}(2, \R)$ is identified with isometries of $\H^2$ that do not necessarily
preserves the orientation as follows. If $\begin{pmatrix}a & b \\c & d\end{pmatrix} \in \PSL^{\pm}(2, \R)$ with
determinant $-1$ we set
\[
\begin{pmatrix}a&b\\c&d\end{pmatrix}
\cdot z
:=
\frac{-a \overline{z} + b}{-c\overline{z} + d}.
\]
The group $\widehat{G_\lambda}$ is defined as the group generated by the three reflections
\[
s = \begin{pmatrix}
0 & 1 \\
1 & 0
\end{pmatrix}
\quad
t = \begin{pmatrix}
-1 & 0 \\
0 & 1
\end{pmatrix}
\quad \text{and} \quad
u = \begin{pmatrix}
-1 & \lambda \\
0 & 1
\end{pmatrix}.
\]
\begin{figure}[!ht]
\centering
\includegraphics{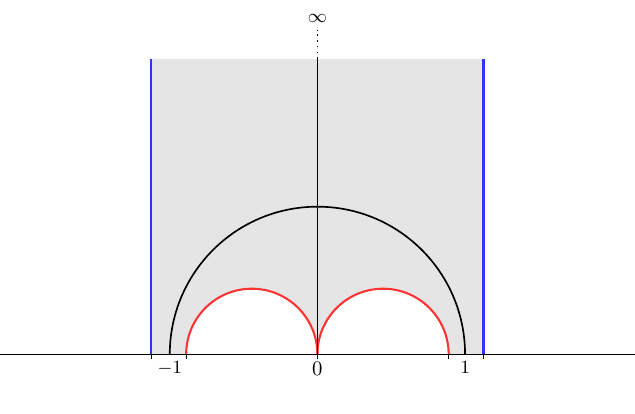}
\caption{The fundamental domain of $G_\lambda$ together with the reflection axis of $s$ and $t$. A fundamental domain for
$\widehat{G_\lambda}$ is obtained by considering a quarter of the one of $G_\lambda$.}
\label{fig:SAndTSymmetries}
\end{figure}

\begin{tcblemma}{}{}
The group $\widehat{G_\lambda}$ contains $G_\lambda$ and more precisely $h_\lambda = ut$ and $v_\lambda = s h_\lambda s$. It can be defined by generators and relations as $\widehat{G_\lambda} = \langle s, t, u | s^2, t^2, u^2, (st)^2 \rangle$\footnote{Such presentation is called a $(2, \infty, \infty)$-triangle group.}. Moreover $\widehat{G_\lambda}$ is a semi-direct product of $G_\lambda$ with the Klein four group $K := \langle s, t \rangle$ and we have
$s\,h_\lambda\,s = v_\lambda$ and $t\,h_\lambda\,t=(h_\lambda)^{-1}$
The fundamental domain of $D_\lambda$ of $G_\lambda$ from Theorem~\ref{thm:GLambdaFundamentalDomain} is preserved by the group $\langle s, t \rangle$.
\end{tcblemma}

\begin{proof}
The relations $h_\lambda = u$ and $v_\lambda = s h_\lambda s$ and $(h_\lambda)^{-1} = t h_\lambda t$ are straightforwrd matrix multiplications.
One can check that the elements $s$ and $t$ generate a group isomorphic to the Klein 4 group $\Z/2\Z \times \Z/2\Z$ that correspond to
a dihedral group of isometries of $\H^2$ fixing $\sqrt{-1}$. The group $K$ stabilizes the fundamental region $D_\lambda$ from Theorem~\ref{thm:GLambdaFundamentalDomain} and permutes the 4 sub-regions delimited by the axes of $s$ and $t$ as can be seen on Figure~\ref{fig:SAndTSymmetries}.
\end{proof}

The last ingredient needed for the proof of Theorem~\ref{thm:GLimitSetVersusTLimitSet} is the relation between $\widehat{G_\lambda}$ and Rosen continued fractions. The homographies that define the branches of the Rosen map are elements of $\widehat{G_\lambda}$. Indeed, the matrices $H_{\epsilon_1,a_1}$ defined in~\eqref{eq:matrixHEpsilonA} satisfy
\begin{equation} \label{eq:RosenAndG}
H_{+1:a_1} = s \, (h_\lambda)^{a_1}
\quad \text{and} \quad
H_{-1:a_1} = s \, t \, (h_\lambda)^{a_1}.
\end{equation}
The following lemma, whose proof is left to the reader, shows that the group $\widehat{G_\lambda}$ decomposes as finitely many translates of continued fraction matrices.
\begin{tcblemma}{}{RosenAndGroup}
For $x, y \in K$ let $\Delta_{x,y} := \{x (h_\lambda)^{a_1} x_2 (h_\lambda)^{a_2} x_3 \cdots x_{n} (h_\lambda)^{a_n} y : x_i \in \{s, st\}, a_i > 0 \}$. Then
\[
\widehat{G_\lambda}
=
K \sqcup \bigsqcup_{x, y \in K} \Delta_{x,y}.
\]
\end{tcblemma}

\begin{proof}[Proof of Theorem~\ref{thm:GLimitSetVersusTLimitSet}]
Let $x \in \Lambda(G_\lambda) \subset \R \cup \{\infty\}$, that is there exists a sequence $g_n \in G_\lambda$ such that for any $z \in \H^2$, $g_n z \to x$ as $n \to \infty$.

Since $K = \langle s, t \rangle$ stabilizes $\Lambda(G_\lambda)$ the limit set coincide : $\Lambda(\widehat{G_\lambda}) = \Lambda(G_\lambda)$. This shows in particular that $\Lambda(T) \subset \Lambda(G_\lambda)$. From Lemma~\ref{lem:RosenAndGroup}, one deduces that $\Lambda(G_\lambda) = \cup_{g \in K} g \Lambda(T)$.
\end{proof}

\begin{proof}[Proof of Theorem~\ref{thm:GLambdaLimitSet}]
The limit set $\widehat{G_\lambda}$ has a single cusp which is the orbit of $0$. A point $x \in \Lambda(G_\lambda)$ with a finite continued fraction expansion $s$ can be written as
$x = H_s(0)$. But $H_s$ is an homography in $\widehat{G_\lambda}$, so $x$ is stabilized by $H_s p H_s^{-1}$ where $p(z) = z/(z+1)$ is the parabolic element fixing $0$.

Conversely, if $x$ is stabilized by $g p g^{-1}$ then one can write $g = H_s$ for some $s \in \cS$.

Next, if $x$ has a periodic continued fraction expansion then $x = H_s(x)$ for some $H_s$. Since all $H_s$ are hyperbolic, we are done.

Conversely, if $g x = x$ then we can write down $g = H_s$.
\end{proof}

\chapter{Eigenvalues of infinite graphs. Geometric Structures}
\label{Apendix:EigenvaluesInfiniteGraphs}

In the first section of this Appendix we review some results about the spectrum of infinite graphs, more precisely about
positive $\lambda$-harmonic functions. These considerations are crucial for the
Hooper-Thurston-Veech construction considered in
Section~\ref{sec:HooperThurstonVeechConstruction}. Most of the material is well known and can be found in the
survey by B.~Mohar and W.~Woess~\cite{MoharWoess}. Our discussion follows Appendix~C in P.~Hooper~\cite{Hooper-infinite_Thurston_Veech}. In the second and last section, we briefly review the general notion of $(G,X)$-structures, developing map and holonomy.

\section{Spectral radius and Martin boundary}
Let $\cG = (V,E)$ be a connected graph (or multigraph) consisting of possibly infinitely many vertices
but such that each vertex $v$ has finitely many neighbors. We consider
the space of functions $\R^V$ from the vertices $V$ of the graph to the real numbers.
For each $p \in [1,\infty]$ we denote $\ell^p(V)$ the subspace of power $p$
integrable functions. Recall that we have the norm:
\[
\|f\|_{p} = \left( \sum_{v \in V} f(v)^p \right)^{1/p}.
\]
Because $V$ is discrete we have $\ell^p(V) \subset \ell^{p'}(V)$ for
$p \leq p'$. In particular $\ell^1(V) \subset \ell^2(V) \subset \ell^\infty(V)$.

We denote by $\{e_v\}_{v \in V}$ the functions
\[
e_v(u) =
\left\{\begin{array}{lll}
    1 & \text{if $u=v$} \\
    0 & \text{otherwise}
\end{array} \right. .
\]
This collection $\{e_v\}_{v \in V}$ is a Hilbert basis of $\ell^2(V)$ and for
any $f\in\R^V$ we have $f = \sum_{v \in V} f(v) e_v$.

The \emph{adjacency operator} $A = A_\cG$ of $\cG$ is the operator on $\R^V$ defined by
\[
(A f)(u) = \sum_{uv \in E} f(v).
\]
It is important to recall that in the case that $\cG$ is a multigraph, the set $E$ is a multiset, hence the same term can appear several times in the sum above.

\begin{tcbexercise}{}{}
Let $\cG$ be a graph and $A$ its adjacency operator.
\begin{enumerate}
\item For $p \in [1, \infty]$ and $f \in \ell_p(V)$ show that
\[
\|A f\|_p \leq \left(\max_{v \in V} \deg_\cG(v)\right) \cdot \|f\|_p.
\]
\item For each $v \in V$ show that
\[
\|A e_v\|_p \geq \deg_\cG(v)^{1/p}
\]
\item For each $v$ let $e'_v = \sum_{w \in N_\cG(v)} e_w$ where $N_\cG(v)$ are the neighbors
of $v$ in $\cG$. Show that
\[
\|A e'_v\|_p \geq \deg_\cG(v)^{1 - 1/p} \cdot \|e'_v\|_p.
\]
\item Deduce that the following are equivalent
\begin{enumerate}
\item $\cG$ has bounded degree,
\item $A$ preserves $\ell^p(V)$ for some $p \in [1,\infty]$,
\item $A$ preserves $\ell^p(V)$ for all $p \in [1,\infty]$.
\end{enumerate}
\end{enumerate}
\end{tcbexercise}

The \emph{coefficients} of $A$ are the numbers
$\{\langle A e_u, e_v \rangle\}_{u,v}$\footnote{Recall that $\langle (x_i), (y_i) \rangle=\sum x_iy_i$}. Note that the coefficient
$u,v$ is nothing else than the multiplicity of edges from $u$ to $v$.
The \emph{spectral radius} $\lambda_0(\cG)$ of the graph is the growth rate of coefficients of the powers of
\begin{equation}
\label{eq:DefinitionSpectralRadius}
\lambda_0(\cG) := \limsup_{n \to \infty} \left( \langle A^n e_u, e_v \rangle \right)^{1/n}
\end{equation}
where $u$ and $v$ are any two vertices of $\cG$. When the context is clear, we abbreviate $\lambda_0(\cG)$ by $\lambda_0$.

\begin{tcbexercise}{}{}
\begin{enumerate}
\item Show that $\langle A^n e_u, e_v \rangle$ is the number of oriented paths from $u$ to $v$ of length $n$.
\item Show that the quantity in~\eqref{eq:DefinitionSpectralRadius} is indeed independent of the choice of
the vertices $u$ and $v$ (hint: use the connectedness of the graph).
\end{enumerate}
\end{tcbexercise}

\begin{tcblemma}{}{EstimatesLambda0}
Let $\cG$ be a connected infinite graph of bounded degree, that is $\deg(\cG):=\max_{v \in V} \deg_\cG(v)<\infty$. Then $2 \sqrt{\deg(\cG) - 1} \leq \lambda_0(\cG) \leq \deg(\cG)$. In particular $\lambda_0(\cG) \geq 2$,

\end{tcblemma}

\begin{proof}
 The upper bound can be shown in various manners, we use a classical result of Schur (see Theorem 6.12-A~\cite{Taylor58}): if $(\alpha_{uv})$ is an infinite matrix of real scalars such that:
$$
\sum_v|\alpha_{uv}|\leq M_1\hspace{1cm}\text{for all $u$}
$$
and
$$
\sum_u|\alpha_{uv}|\leq M_2\hspace{1cm}\text{for all $v$}
$$
then $(\alpha_{uv})$ represents a bounded operator $A:\ell^2\to\ell^2$ such that $||A||\leq\sqrt{M_1M_2}$. For the lower bound let $T\to\cG$ be the universal cover and let  $T_{min}\hookrightarrow T$ a regular tree whose degree is $d=\min_{v \in V} \deg_{T}(v)$. Then in $T_{min}$
for any $n \geq 1$ and any $v \in V(T_{min}) $ we have:
\[
\langle A^{2n} e_v, e_v \rangle \geq C_n \cdot d \cdot (d-1)^{n-1}.
\]
where $C_n$ is the nth Catalan number. Since
\[
\lim_{n \to \infty} \left( C_n \cdot d \cdot (d-1)^{n-1} \right)^{1/(2n)}
=
2 \sqrt{d-1}
\]
we obtain the lower bound for $T_{min}$. Since we are working with a regular tree in the universal covering of $\cG$, this lower bound also works for $\cG$.
\end{proof}

Now we distinguish more carefully the different possible asymptotic behavior. For any connected multigraph $\cG$ let:
\[
R_\lambda := \frac{1}{\lambda} \sum_{n=0}^\infty \left( \frac{1}{\lambda} A \right)^n.
\]
It is known that all entries of $R_\lambda$ are real numbers if $\lambda>\lambda_0(\cG)$ and every entry of $R_\lambda$ is equal to $+\infty$ if $\lambda<\lambda_0(\cG)$.

\begin{tcbdefinition}{}{}
Let $\cG$ be a connected multigraph of bounded degree.
We say that $\cG$ is transient if all entries of $R_{\lambda_0}$ are in $\R$ and recurrent otherwise.
When $\cG$ is recurrent, we define
\[
\delta_{u,v} := \limsup_{n \to \infty} \frac{1}{\lambda_0^n} \left( \langle A^n e_u, e_v \rangle \right)
\]
We say that $\cG$ is \emph{null-recurrent} if $\delta_{u,v} = 0$ and \emph{positive recurrent} otherwise.
These properties are independent of the choice of the pair $u,v$ of vertices.
\end{tcbdefinition}

The matrix $R_\lambda$ satisfies $\lambda R_\lambda = A R_\lambda + I$. Hence, when $R_\lambda$ has entries in the real numbers, one can think of its columns $R_\lambda e_v$ as being 'almost' eigenfunctions of $A$.

For the reader familiar with random walk, the terminology in the above definitions is
explained in the exercise below.
\begin{tcbexercise}{}{}
Assume that $\cG$ is a $d$-regular graph, \ie each vertex has degree $d$.
Let $\{Z_n\}_{n \geq 0}$ be a random walk in $\cG$ starting from a root vertex $o$, \ie
$Z_0 = 0$, and where at each step the walk jumps randomly one of its neighbor with equal
probability. Let
\[
\tau := \inf \{n > 0: Z_n = o\}
\]
\begin{enumerate}
\item Show that $\P(\tau < +\infty) = 1$ if and only if $\cG$ is recurrent.
\item Show that $\E(\tau) < +\infty$ if and only if $\cG$ is positive recurrent.
\end{enumerate}
\end{tcbexercise}

Infinite graphs of bounded degree which are recurrent can always be used to produce Hooper-Thurston-Veech surfaces:

\begin{tcbtheorem}{}{SpectralRadiusRecurrentCase}
Let $\cG$ be a graph of bounded degree which is recurrent. Then
there exists a positive solution to the equation $A f = \lambda_0(\cG) f$ which
is unique up to scaling. Moreover, $\cG$ is positive-recurrent if and only if
this solution belongs to $\ell^2(V)$.
\end{tcbtheorem}

We refer to Theorem 6.2 in~\cite{MoharWoess} and references therein for a proof of the result above.
For transient graphs, one can parametrize the positive eigenfunctions by the so-called
Martin boundary. We briefly recall this notion in what follows.

For $o \in V$ a root vertex we define the \emph{$\lambda$-Martin kernel}:
\[
K_{\lambda,o}:
\begin{array}{lll}
V^2 & \to & \R \\
(u,v) & \mapsto & \frac{R_\lambda(e_v)(u)}{R_\lambda(e_v)(o)} .
\end{array}
\]
Now, to a vertex $v$ in $V$, we associate the function
$\phi_v: u \to K_{\lambda,o}(u,v)$. This defines an embedding
of $V\hookrightarrow\R^V$.
The Martin compactification $V_\lambda$ is the smallest compactification of $V$
to which all $\phi_u$ extend continuously, or, equivalently, for which the embedding $V\hookrightarrow\R^V$ has a continuous extension.
The \emph{$\lambda$-Martin boundary} is $M_\lambda := V_\lambda \setminus V$.
By construction, given $\zeta \in M_\lambda$ the quantity $\phi_v(\zeta)$ is
well defined. To obtain the correspondence that interest us we need to introduce the notion of minimal elements: these are points $\zeta\in M_\lambda$ such that, whenever $\eta_1$, $\eta_2\in M_\lambda$ satisfty $t k_{\eta_1}+(1-t)k_{\eta_2}=k_\zeta$ for some $0<t<1$, one has that $\zeta=\eta_1=\eta_2$.

\begin{tcbtheorem}{}{MartinBoundary}
Let $\cG$ be a transient graph of bounded degree and $\lambda$ such
either $\lambda > \lambda_0$ or $\lambda = \lambda_0$.
Then the function
\begin{equation}
    \label{eq:correspondence-Martin-boundary}
\begin{array}{lll}
M_\lambda & \to & \R^V \\
\zeta & \mapsto & (v \mapsto \phi_v(\zeta))
\end{array}
\end{equation}
provides a bijection between the minimal points of the Martin boundary and the extremal positive
$\lambda$-harmonic functions on $\cG$ which take value $1$ at the root vertex $o$.
\end{tcbtheorem}
\begin{tcbcorollary}{}{ExistenceEigenfunctions}
Let $\cG$ be a graph of bounded degree and $\lambda$ such
either $\lambda > \lambda_0$ or $\lambda = \lambda_0$ and $\cG$ is transient.
Then there exists a positive $\lambda$-harmonic function on $\cG$.
\end{tcbcorollary}

\subsection{Minimal $\lambda$-harmonic functions of regular trees}

In the following paragraphs we illustrate the correspondence between positive $\lambda$-harmonic eigenfunctions and points in the minimal Martin boundary when\footnote{The description of all positive $\lambda$-harmonic functions for the regular tree $\mathbb{T}_2$ was carried out in Example~\ref{exa:HTVStaircase}}  $\cG=\mathbb{T}_{q+1}$ is an regular tree of degree $d=q+1$, $q\geq 2$. The following general result (see Appendix F in~\cite{Hooper-infinite_Thurston_Veech} and references therein) stablishes the correspondence between extremal positive $\lambda$-harmonic functions on $\mathbb{T}_d$ and its Gromov boundary.


\begin{tcbtheorem}{}{}
If $\cG$ is a hyperbolic infinite graph and $\lambda > \lambda_0$, then every point in the Martin boundary
$M_\lambda$ is minimal and homeomorphic to the Gromov boundary of $\cG$. Moreover, if $\cG$ is a transient tree, then $M_{\lambda_0}$ is homeomorphic to the Gromov boundary of $\cG$.
\end{tcbtheorem}

 The graph $\mathbb{T}_{q+1}$ is transient, hyperbolic and its Gromov boundary is homeomorphic to its space of ends, which in turn is homeomorphic to a Cantor set. We refer the reader to~\cite{Loh17} and references therein for a discussion on Gromov hyperbolic spaces.

Let $o\in\mathbb{T}_{q+1}$ be a fixed vertex and $\omega\in\partial \mathbb{T}_{q+1}$ be a point in its Gromov boundary. We identify $\omega$ with a (equivalence class of) path $[x_0=o,x_{-1},x_{-2}\ldots]$ going to infinity. The end $\omega$ defines a partition $V(\mathbb{T}_{q+1})=\sqcup_{n\in\Z} H_n$ into horocycles, where $o\in H_0$. Every vertex in $H_n$ has a set of adjacent vertices formed by $q$ vertices in $H_{n+1}$ (which are 'further away' from $\omega$) and one vertex in $H_{n-1}$ (which is closer). From the work of Cartier (see Section 9 in~\cite{Cartier73}) we deduce that every minimal positive $\lambda$-harmonic function $\textbf{h}$ is constant along horocycles and grows exponentially along the ray $[x_0=0,x_1,\ldots]$. These facts, together with the equation $A\textbf{h}=\lambda \textbf{h}$, lead to the linear recurrence relation:
\[
\lambda \textbf{h}_n = q \textbf{h}_{n+1} + \textbf{h}_{n-1}
\]
Here $\textbf{h}_n$ is the value of $\textbf{h}$ along the horocyle $H_n$. The general solution of the equation above is of the form $\textbf{h}_n=Ar_+^n+Br_-^n$ where $r_+=\frac{1}{2q}(\lambda+\sqrt{\lambda^2-4q})$ and $r_-=\frac{1}{2q}(\lambda-\sqrt{\lambda^2-4q})$. Hence, we must have that $\lambda\geq 2\sqrt{q}$. For the extreme value
$\lambda=2\sqrt{q}$ the normalization  $\textbf{h}(o)=1$ leads to $\textbf{h}_n=r_+^n$. For $\lambda>2\sqrt{q}$ the same normalization leads to $\textbf{h}_n=r_+^n$ being the solution that corresponds to $\omega$, because $r_-^n$ decreases exponentially as we move along $[x_0=o,x_{-1},x_{-2}\ldots]$. We finish this Section with an exercise for a motivated reader.

\begin{tcbexercise}{}{}
In the preceding paragraph, the functions $\textbf{h}_+$ and $\textbf{h}_-$ depend on the choice of a basepoint $o\in\mathbb{T}_{q+1}$ and an end $\omega\in\partial\mathbb{T}_{q+1}$. Denote this dependency by $\textbf{h}_{+,o,\omega}$ and $\textbf{h}_{-,o,\omega}$, respectively. Show that , if we normalize $\textbf{h}_{+,o,\omega}(o)=\textbf{h}_{-,o,\omega}=1$,
then there exists a probability measure $\mu$ on $\partial\mathbb{T}_{q+1}$ such that:
$$
\textbf{h}_{-,o,\omega}=\int_{\partial\mathbb{T}_{q+1}} \textbf{h}_{+,o,\omega'} d\mu(\omega').
$$
\end{tcbexercise}





\section{Geometric structures}
  \label{Appendix:GXStructures}

We briefly review the general notion of $(G,X)$-structure, developing map and holonomy. For a more detailed discussion, we refer the reader to Chapter 43 in~\cite{Thurston97}.
Let $X$ be a manifold and $G$ a group acting on $X$ by homeomorphisms. A $(G,X)$-structure on a manifold $M$ is an atlas $\mathcal{A}=\{\varphi_i:U_i\to X\}$ whose transition functions are restrictions of elements in $G$. The ideas leading to this notion seem to be as old as Poincar\'e's studies of second order holomorphic ODEs on Riemann surfaces. For a detailed discussion on the history of $(G,X)$-structures, we refer the reader to \href{https://mathoverflow.net/questions/332012/history-of-the-notion-of-g-x-structure}{this post in Mathoverflow}\footnote{https://mathoverflow.net/questions/332012/history-of-the-notion-of-g-x-structure}.
The following are examples of $(G,X)$-strutures appering in this book.
\begin{enumerate}
\item A translation surface $M$ is a $\Trans(\R^2)$-structure on the complement of the set of conical singularities (see Definition~\ref{def:TranslationSurfaceGeometric}).
\item A half-translation surface (\emph{a.k.a.} quadratic differential) $M$ is a $\Trans(\R^2) \rtimes \{\pm 1\}$-structure on the complement of the set of conical singularities (see Definition~\ref{def:FlatSurfaceGeometric}).
\item A flat conical metric is a $\Isom^+(\R^2)= \Trans(\R^2) \rtimes \SO(2)$-structure (defined on the complement of the conical singularities). For example, consider the flat conical metric on the sphere given by the triangular pillowcase
$\mathbb{S}^2(P)$ used in Example~\ref{exa:IrrationalBilliardAndProof}.
\item A dilation surface is a $\Trans(\R^2) \rtimes \R_{>0}$-structure.
\item Every rotational component (see~\ref{SSEC:RotationalComponents}) admits a  $\Trans(\R)$-structure.
\end{enumerate}

\textbf{Developing maps and holonomy}. To construct the developing map of a $(G,X)$-structure on a manifold $M$ we need the following technical condition: if $g,h$ are two elements of $G$ which coincide on a non-emtpy open set, then $g=h$. This is true for all examples considered above. The developing map:
\begin{equation}
  \label{def:defmap}
  dev:\widetilde{M}\to X,
\end{equation}
where $\pi:\widetilde{M}\to M$ is the universal cover, plays the role of a global coordinate chart for the $(G,X)$-structure. To define it remark first that if $\varphi_i:U_i\to X$ are coordinate charts of $M$, $i=1,2$, such that $U_1\cap U_2\neq \emptyset$ then by the technical condition above, there exists a unique $g\in G$ such that $g\circ\varphi_2=\varphi_1$ on $U_1\cap U_2$. This allows us to extend $\varphi_1$ to a map $\varphi:U_1\cup U_2 \to X$ satisfying $\varphi = g\circ\varphi_2$ on $U_2$ and $\varphi=\varphi_1$ on $U_1$. Using this principle, it is easy to see that each lift $\widetilde{\varphi_1}$ of a chart of $M$ to the universal cover defines a map $dev:\widetilde{M}\to X$ extending $\widetilde{\varphi_1}$. This map is well defined up to the choice of the chart and its lift, and is called the \emphdef[developing map (G-structure)]{developing map} relative to the chart $\phi_1$. Different choices of charts and lifts define developing maps which differ by composition by an element of $G$. This can be expresed in the following equivariance condition:
\begin{equation}
  \label{holonomy-representation}
dev \circ \gamma = hol(\gamma) \circ dev
\end{equation}
for each deck transformation of the universal cover $\gamma\in\pi_1(M)$. Here
$$
hol:\pi_1(M)\to G
$$
is a homomorphism called the \emph{holonomy} (representation) of the $(G,X)$-structure. Remark that $h$ is also well defined, up to conjugation by an element of $G$.

\backmatter

\bibliographystyle{alpha} 
\bibliography{biblio}

\def\polhk#1{\setbox0=\hbox{#1}{\ooalign{\hidewidth
  \lower1.5ex\hbox{`}\hidewidth\crcr\unhbox0}}}
\begin{thebibliography}{DEDML00}

\bibitem[Aar97]{Aaronson}
Jon Aaronson.
\newblock {\em An introduction to infinite ergodic theory}, volume~50 of {\em
  Mathematical Surveys and Monographs}.
\newblock American Mathematical Society, Providence, RI, 1997.

\bibitem[{Abi}80]{Abikoff80}
William {Abikoff}.
\newblock {\em {The real analytic theory of Teichm\"uller space}}, volume 820.
\newblock Springer, Cham, 1980.

\bibitem[{Abr}66a]{Abramov661}
L.~M. {Abramov}.
\newblock {On the entropy of a flow.}
\newblock {\em {Transl., Ser. 2, Am. Math. Soc.}}, 49:167--170, 1966.

\bibitem[{Abr}66b]{Abramov662}
L.~M. {Abramov}.
\newblock {The entropy of a derived automorphism.}
\newblock {\em {Transl., Ser. 2, Am. Math. Soc.}}, 49:162--166, 1966.

\bibitem[AD16]{AvilaDelecroix}
Artur Avila and Vincent Delecroix.
\newblock Weak mixing directions in non-arithmetic {V}eech surfaces.
\newblock {\em J. Amer. Math. Soc.}, 29(4):1167--1208, 2016.

\bibitem[ADDS15]{AvilaDolgopyatDuryevSarig-2015}
A.~Avila, D.~Dolgopyat, E.~Duryev, and O.~Sarig.
\newblock The visits to zero of a random walk driven by an irrational rotation.
\newblock {\em Israel J. Math.}, 207(2):653--717, 2015.

\bibitem[AH20]{AvilaHubert-recurrence}
A.~Avila and P.~Hubert.
\newblock Recurrence for the wind-tree model.
\newblock {\em Ann. Inst. H. Poincar\'{e} C Anal. Non Lin\'{e}aire},
  37(1):1--11, 2020.

\bibitem[{Ahl}66]{Ahlfors-quasiconformal_mappings}
Lars~V. {Ahlfors}.
\newblock {\em {Lectures on quasiconformal mappings}}.
\newblock 1966.

\bibitem[ANSS02]{AaronsonNakadaSarigSolomyak02}
Jon Aaronson, Hitoshi Nakada, Omri Sarig, and Rita Solomyak.
\newblock Invariant measures and asymptotics for some skew products.
\newblock {\em Israel J. Math.}, 128:93--134, 2002.

\bibitem[AOW85]{ArnouxOrnsteinWeiss}
P.~Arnoux, D.~Ornstein, and B.~Weiss.
\newblock Cutting and stacking, interval exchanges and geometric models.
\newblock {\em Israel J. of Math.}, 50, 1985.

\bibitem[Arn06]{ArnoldODE}
Vladimir~I. Arnold.
\newblock {\em Ordinary differential equations}.
\newblock Universitext. Springer-Verlag, Berlin, 2006.
\newblock Translated from the Russian by Roger Cooke, Second printing of the
  1992 edition.

\bibitem[Art17]{Artigiani-Eaton_lenses}
Mauro Artigiani.
\newblock Exceptional ergodic directions in {E}aton lenses.
\newblock {\em Israel J. Math.}, 220(1):29--56, 2017.

\bibitem[AS60]{AhlforsSario-book}
Lars~V. Ahlfors and Leo Sario.
\newblock {\em Riemann surfaces}.
\newblock Princeton Mathematical Series, No. 26. Princeton University Press,
  Princeton, N.J., 1960.

\bibitem[Bea83]{Beardon}
Alan~F. Beardon.
\newblock {\em The geometry of discrete groups}, volume~91 of {\em Graduate
  Texts in Mathematics}.
\newblock Springer-Verlag, New York, 1983.

\bibitem[BFG20]{BoulangerFougeronGhazouani20}
Adrien {Boulanger}, Charles {Fougeron}, and Selim {Ghazouani}.
\newblock {Cascades in the dynamics of affine interval exchange
  transformations.}
\newblock {\em {Ergodic Theory Dyn. Syst.}}, 40(8):2073--2097, 2020.

\bibitem[BK16]{BezuglyiKarpel16}
S.~{Bezuglyi} and O.~{Karpel}.
\newblock {Bratteli diagrams: structure, measures, dynamics.}
\newblock In {\em {Dynamics and numbers. A special programm: June 1 -- July 31,
  2014. International conference: July 21--25, 2014, Max-Planck Institute for
  Mathematics, Bonn, Germany. Proceedings}}, pages 1--36. Providence, RI:
  American Mathematical Society (AMS), 2016.

\bibitem[BL23]{Bruin-Lukina-2101}
Henk Bruin and Olga Lukina.
\newblock Rotated odometers.
\newblock {\em J. Lond. Math. Soc., II. Ser.}, 107(6):1983--2024, 2023.

\bibitem[Blu12]{Blume2012}
Frank Blume.
\newblock An entropy estimate for infinite interval exchange transformations.
\newblock {\em Math. Z.}, 272(1-2):17--29, 2012.

\bibitem[{Bow}12]{Bowman12}
Joshua~P. {Bowman}.
\newblock {Finiteness conditions on translation surfaces.}
\newblock In {\em {Quasiconformal mappings, Riemann surfaces, and Teichm\"uller
  spaces. AMS special session in honor of Clifford J. Earle, Syracuse, NY, USA,
  October 2--3, 2010}}, pages 31--40. Providence, RI: American Mathematical
  Society (AMS), 2012.

\bibitem[Bow13]{Bowman-Arnoux-Yoccoz}
Joshua~P. Bowman.
\newblock The complete family of {A}rnoux-{Y}occoz surfaces.
\newblock {\em Geom. Dedicata}, 164:113--130, 2013.

\bibitem[BV13]{BowmanValdez13}
Joshua~P. Bowman and Ferr{\'a}n Valdez.
\newblock Wild singularities of flat surfaces.
\newblock {\em Israel J. Math.}, 197(1):69--97, 2013.

\bibitem[Cab12]{Cabrol}
J.~Cabrol.
\newblock {\em Origamis infinis: groupe de Veech et flot lin\'eaire}.
\newblock PhD thesis, Universit\'e d'Aix-Marseille, 2012.
\newblock PhD Thesis.

\bibitem[Cal04]{Calta04}
K.~Calta.
\newblock Veech surfaces and complete periodicity in genus two.
\newblock {\em J. of Amer. Math. Soc.}, 17(4):871--908, 2004.

\bibitem[Car73]{Cartier73}
Pierre Cartier.
\newblock G{\'e}om{\'e}trie et analyse sur les arbres.
\newblock Sem. {Bourbaki} 1971/72, {No}. 407, {Lect}. {Notes} {Math}. 317,
  123-140 (1973)., 1973.

\bibitem[CFS82a]{CornfeldFominSinai}
I.~P. Cornfeld, S.~V. Fomin, and Ya.~G. Sina{\u\i}.
\newblock {\em Ergodic theory}, volume 245 of {\em Grundlehren der
  Mathematischen Wissenschaften [Fundamental Principles of Mathematical
  Sciences]}.
\newblock Springer-Verlag, New York, 1982.
\newblock Translated from the Russian by A. B. Sosinski{\u\i}.

\bibitem[CFS82b]{CFS82}
I.~P. {Cornfeld}, S.~V. {Fomin}, and Ya.~G. {Sinai}.
\newblock {\em {Ergodic theory. Transl. from the Russian by A. B.
  Sossinskii.}}, volume 245.
\newblock Springer, Berlin, 1982.

\bibitem[CG19]{Coulon-Gruber}
R\'{e}mi Coulon and Dominik Gruber.
\newblock Small cancellation theory over {B}urnside groups.
\newblock {\em Adv. Math.}, 353:722--775, 2019.

\bibitem[CGL06]{ChamanaraGardinerLakic}
R.~Chamanara, F.~P. Gardiner, and N.~Lakic.
\newblock A hyperelliptic realization of the horseshoe and baker maps.
\newblock {\em Erg. Th. and Dyn. Syst.}, 26(6):1749--1768, 2006.

\bibitem[{Cha}69]{Chacon69}
R.~V. {Chacon}.
\newblock {Weakly mixing transformations which are not strongly mixing.}
\newblock {\em {Proc. Am. Math. Soc.}}, 22:559--562, 1969.

\bibitem[Cha04]{Chamanara}
R.~Chamanara.
\newblock {\em Affine automorphism groups of surfaces of infinite type}, volume
  355, pages 123--145.
\newblock American Math. Soc., 2004.

\bibitem[CI80]{CoxIsham}
David~Roxbee Cox and Valerie Isham.
\newblock {\em Point processes}.
\newblock Chapman \& Hall, London-New York, 1980.
\newblock Monographs on Applied Probability and Statistics.

\bibitem[CK76]{ConzeKeane76}
J.-P. Conze and M.~Keane.
\newblock Ergodicit\'e d'un flot cylindrique.
\newblock In {\em S\'eminaire de {P}robabilit\'es, {I} ({U}niv. {R}ennes,
  {R}ennes, 1976), {E}xp. {N}o. 5}, page~7. D\'ept. Math. Informat., Univ.
  Rennes, Rennes, 1976.

\bibitem[Con76]{Conze76}
J.-P. Conze.
\newblock Equir\'epartition et ergodicit\'e de transformations cylindriques.
\newblock In {\em S\'eminaire de {P}robabilit\'es, {I} ({U}niv. {R}ennes,
  {R}ennes, 1976), {E}xp. {N}o. 2}, page~21. D\'ept. Math. Informat., Univ.
  Rennes, Rennes, 1976.

\bibitem[CW12]{CohenWeiss2012}
Meital Cohen and Barak Weiss.
\newblock Parking garages with optimal dynamics.
\newblock {\em Geom. Dedicata}, 161:157--167, 2012.

\bibitem[{de }00]{delaHarpe00}
Pierre {de la Harpe}.
\newblock {\em {Topics in geometric group theory}}.
\newblock Chicago: The University of Chicago Press, 2000.

\bibitem[DEDML98]{DegliEspostiDelMagnoLenci-1998}
Mirko Degli~Esposti, Gianluigi Del~Magno, and Marco Lenci.
\newblock An infinite step billiard.
\newblock {\em Nonlinearity}, 11(4):991--1013, 1998.

\bibitem[DEDML00]{DegliEspostiDelMagnoLenci-2000}
Mirko Degli~Esposti, Gianluigi Del~Magno, and Marco Lenci.
\newblock Escape orbits and ergodicity in infinite step billiards.
\newblock {\em Nonlinearity}, 13(4):1275--1292, 2000.

\bibitem[Del13]{Delecroix-divergent}
V.~Delecroix.
\newblock Divergent directions in some periodic wind-tree models.
\newblock {\em J. of Modern Dyn.}, 7:1--29, 2013.

\bibitem[DFG19]{DuryevFougeronGhazouani19}
Eduard {Duryev}, Charles {Fougeron}, and Selim {Ghazouani}.
\newblock {Dilation surfaces and their Veech groups.}
\newblock {\em {J. Mod. Dyn.}}, 14:121--151, 2019.

\bibitem[DHL14]{DelecroixHubertLelievre14}
Vincent Delecroix, Pascal Hubert, and Samuel Leli{\`e}vre.
\newblock Diffusion for the periodic wind-tree model.
\newblock {\em Ann. Sci. \'Ec. Norm. Sup\'er. (4)}, 47(6):1085--1110, 2014.

\bibitem[DHV]{DHV2}
V.~Delecroix, P.~Hubert, and F.~Valdez.
\newblock Infinite translation surfaces in the wild. volume 2: dynamics.
\newblock in preparation.

\bibitem[DK03]{Diestel03}
Reinhard Diestel and Daniela K{\"u}hn.
\newblock Graph-theoretical versus topological ends of graphs.
\newblock {\em J. Combin. Theory Ser. B}, 87(1):197--206, 2003.
\newblock Dedicated to Crispin St. J. A. Nash-Williams.

\bibitem[DK21]{Dajani-Kalle2021}
Karma {Dajani} and Charlene {Kalle}.
\newblock {\em {A first course in ergodic theory}}.
\newblock Boca Raton, FL: CRC Press, 2021.

\bibitem[Don11]{Donaldson11}
Simon Donaldson.
\newblock {\em Riemann surfaces}, volume~22 of {\em Oxford Graduate Texts in
  Mathematics}.
\newblock Oxford University Press, Oxford, 2011.

\bibitem[Dow11]{Downarowicz-book}
Tomasz Downarowicz.
\newblock {\em Entropy in dynamical systems}, volume~18 of {\em New
  Mathematical Monographs}.
\newblock Cambridge University Press, Cambridge, 2011.

\bibitem[DVJ03]{DaleyVereJones}
D.~J. Daley and D.~Vere-Jones.
\newblock {\em An introduction to the theory of point processes. {V}ol. {I}}.
\newblock Probability and its Applications (New York). Springer-Verlag, New
  York, second edition, 2003.
\newblock Elementary theory and methods.

\bibitem[DZ20]{DelecroixZorich-cries_and_whispers}
Vincent Delecroix and Anton Zorich.
\newblock Cries and whispers in wind-tree forests.
\newblock In {\em What's next?---the mathematical legacy of {W}illiam {P}.
  {T}hurston}, volume 205 of {\em Ann. of Math. Stud.}, pages 83--115.
  Princeton Univ. Press, Princeton, NJ, 2020.

\bibitem[EE12]{EhrenfestEhrenfest12}
P.\ Ehrenfest and T.\ Ehrenfest.
\newblock Begriffliche grundlagen der statistischen auffassung in der mechanik.
\newblock {\em Encykl. d. Math. Wissensch.}, IV 2 II, Heft 6, 90 S, 1912.
\newblock (in German, translated in:) \textit{The conceptual foundations of the
  statistical approach in mechanics,} (trans. Moravicsik, M. J.), 10-13 Cornell
  University Press, Itacha NY, (1959).

\bibitem[EW11]{EinsiedlerWard11}
Manfred {Einsiedler} and Thomas {Ward}.
\newblock {\em {Ergodic theory. With a view towards number theory.}}, volume
  259.
\newblock London: Springer, 2011.

\bibitem[FaSU18]{FraczekShiUlcigrai}
Krzysztof Fr\polhk~aczek, Ronggang Shi, and Corinna Ulcigrai.
\newblock Genericity on curves and applications: pseudo-integrable billiards,
  {E}aton lenses and gap distributions.
\newblock {\em J. Mod. Dyn.}, 12:55--122, 2018.

\bibitem[FH19]{Ferenczi-Hubert}
S\'{e}bastien Ferenczi and Pascal Hubert.
\newblock Rigidity of square-tiled interval exchange transformations.
\newblock {\em J. Mod. Dyn.}, 14:153--177, 2019.

\bibitem[FK36]{FoxKershner36}
R.~Fox and R.~Kershner.
\newblock Concerning the transitive properties of geodesics in rational
  polyhedron.
\newblock {\em Duke}, 2(1):147--150, 1936.

\bibitem[FK92]{FarkasKra80}
H.~M. Farkas and I.~Kra.
\newblock {\em Riemann surfaces}, volume~71 of {\em Graduate Texts in
  Mathematics}.
\newblock Springer-Verlag, New York, second edition, 1992.

\bibitem[FKT03]{FeldmanKnoerrerTrubowitz}
Joel {Feldman}, Horst {Kn\"orrer}, and Eugene {Trubowitz}.
\newblock {\em {Riemann surfaces of infinite genus.}}, volume~20.
\newblock Providence, RI: American Mathematical Society (AMS), 2003.

\bibitem[FLP79]{FathiLaudenbachPoenaru}
Albert Fathi, Fran\c{c}ois Laudenbach, and Valentin Po\'enaru.
\newblock {\em Travaux de {T}hurston sur les surfaces}, volume~66 of {\em
  Ast\'erisque}.
\newblock Soci\'et\'e Math\'ematique de France, Paris, 1979.
\newblock S{\'e}minaire Orsay, With an English summary.

\bibitem[FM12]{FarbMargalit}
Benson Farb and Dan Margalit.
\newblock {\em A primer on mapping class groups}, volume~49 of {\em Princeton
  Mathematical Series}.
\newblock Princeton University Press, Princeton, NJ, 2012.

\bibitem[FM14]{ForniMatheus14}
Giovanni Forni and Carlos Matheus.
\newblock Introduction to {T}eichm\"uller theory and its applications to
  dynamics of interval exchange transformations, flows on surfaces and
  billiards.
\newblock {\em J. Mod. Dyn.}, 8(3-4):271--436, 2014.

\bibitem[FMZ11]{ForniMatheusZorich10}
G.~Forni, C.~Matheus, and A.~Zorich.
\newblock Square-tilde cyclic covers.
\newblock {\em J. of Modern Dyn.}, 5(2):285--318, 2011.

\bibitem[For91]{Forster77}
Otto Forster.
\newblock {\em Lectures on {R}iemann surfaces}, volume~81 of {\em Graduate
  Texts in Mathematics}.
\newblock Springer-Verlag, New York, 1991.
\newblock Translated from the 1977 German original by Bruce Gilligan, Reprint
  of the 1981 English translation.

\bibitem[For06]{Forni06}
G.~Forni.
\newblock {\em On the {L}yapunov exponents of the {K}ontsevich-{Z}orich
  cocycle}, pages 549--580.
\newblock Elsevier, 2006.

\bibitem[{Fre}31]{Freudenthal31}
Hans {Freudenthal}.
\newblock {\"Uber die Enden topologischer R\"aume und Gruppen.}
\newblock {\em {Math. Z.}}, 33:692--713, 1931.

\bibitem[FS14]{FraczekSchmoll-reflector}
Krzysztof Fr{\polhk{a}}czek and Martin Schmoll.
\newblock Directional localization of light rays in a periodic array of
  retro-reflector lenses.
\newblock {\em Nonlinearity}, 27(7):1689--1707, 2014.

\bibitem[FU14a]{FraczekUlcigrai-ergodic_examples}
Krzysztof Fr{\polhk{a}}czek and Corinna Ulcigrai.
\newblock Ergodic directions for billiards in a strip with periodically located
  obstacles.
\newblock {\em Comm. Math. Phys.}, 327(2):643--663, 2014.

\bibitem[FU14b]{FraczekUlcigrai-non_ergodicity}
Krzysztof Fr{\polhk{a}}czek and Corinna Ulcigrai.
\newblock Non-ergodic {$\Bbb{Z}$}-periodic billiards and infinite translation
  surfaces.
\newblock {\em Invent. Math.}, 197(2):241--298, 2014.

\bibitem[Ful95]{Fulton95}
William Fulton.
\newblock {\em Algebraic topology}, volume 153 of {\em Graduate Texts in
  Mathematics}.
\newblock Springer-Verlag, New York, 1995.
\newblock A first course.

\bibitem[GHS03]{GutkinHubertSchmidt03}
Eugene {Gutkin}, Pascal {Hubert}, and Thomas~A. {Schmidt}.
\newblock {Affine diffeomorphisms of translation surfaces: periodic points,
  Fuchsian groups, and arithmeticity.}
\newblock {\em {Ann. Sci. \'Ec. Norm. Sup\'er. (4)}}, 36(6):847--866, 2003.

\bibitem[Ghy95]{Ghys95}
{\'E}tienne Ghys.
\newblock Topologie des feuilles g\'en\'eriques.
\newblock {\em Ann. of Math. (2)}, 141(2):387--422, 1995.

\bibitem[GJ00]{GutkinJudge00}
Eugene {Gutkin} and Chris {Judge}.
\newblock {Affine mappings of translation surfaces: Geometry and arithmetic.}
\newblock {\em {Duke Math. J.}}, 103(2):191--213, 2000.

\bibitem[GL00]{GardinerLakic00}
Frederick~P. {Gardiner} and Nikola {Lakic}.
\newblock {\em {Quasiconformal Teichm\"uller theory}}, volume~76.
\newblock Providence, RI: American Mathematical Society, 2000.

\bibitem[GN67]{GunningNarasimhan67}
R.~C. {Gunning} and R.~{Narasimhan}.
\newblock {Immersion of open Riemann surfaces.}
\newblock {\em {Math. Ann.}}, 174:103--108, 1967.

\bibitem[Hat02]{Hatcher}
Allen Hatcher.
\newblock {\em Algebraic topology}.
\newblock Cambridge University Press, Cambridge, 2002.

\bibitem[Hec36]{Hecke1936}
E.~Hecke.
\newblock Über die bestimmung dirichletscher reihen durch ihre
  funktionalgleichung. (german) jfm 63.0264.03.
\newblock {\em Math. Ann.}, 112:664--699, 1936.

\bibitem[HHW13]{HooperHubertWeiss}
W.~Patrick Hooper, Pascal Hubert, and Barak Weiss.
\newblock Dynamics on the infinite staircase.
\newblock {\em Discrete Contin. Dyn. Syst.}, 33(9):4341--4347, 2013.

\bibitem[HL06a]{HubertLanneau06}
Pascal Hubert and Erwan Lanneau.
\newblock Veech groups without parabolic elements.
\newblock {\em Duke Math. J.}, 133(2):335--346, 2006.

\bibitem[HL06b]{HubertLelievre06}
Pascal {Hubert} and Samuel {Leli\`evre}.
\newblock {Prime arithmetic Teichm\"uller discs in ${\mathcal H}(2)$.}
\newblock {\em {Isr. J. Math.}}, 151:281--321, 2006.

\bibitem[HLM09]{HubertLanneauMoeller09}
Pascal {Hubert}, Erwan {Lanneau}, and Martin {M\"oller}.
\newblock {The Arnoux-Yoccoz Teichm\"uller disc.}
\newblock {\em {Geom. Funct. Anal.}}, 18(6):1988--2016, 2009.

\bibitem[HLT11]{HubertLelievreTroubetzkoy11}
Pascal Hubert, Samuel Leli{\`e}vre, and Serge Troubetzkoy.
\newblock The {E}hrenfest wind-tree model: periodic directions, recurrence,
  diffusion.
\newblock {\em J. Reine Angew. Math.}, 656:223--244, 2011.

\bibitem[HMS20]{HubertMatheus}
Pascal Hubert and Carlos Matheus~Santos.
\newblock An origami of genus 3 with arithmetic {K}ontsevich-{Z}orich
  monodromy.
\newblock {\em Math. Proc. Cambridge Philos. Soc.}, 169(1):19--30, 2020.

\bibitem[Hoo07]{Hooper2007}
W.~Patrick Hooper.
\newblock Periodic billiard paths in right triangles are unstable.
\newblock {\em Geom. Dedicata}, 125:39--46, 2007.

\bibitem[Hoo15]{Hooper-infinite_Thurston_Veech}
W.~Patrick Hooper.
\newblock The invariant measures of some infinite interval exchange maps.
\newblock {\em Geom. Topol.}, 19(4):1895--2038, 2015.

\bibitem[HPS92]{HermanPutnamSkau1992}
Richard~H. Herman, Ian~F. Putnam, and Christian~F. Skau.
\newblock Ordered {B}ratteli diagrams, dimension groups and topological
  dynamics.
\newblock {\em Internat. J. Math.}, 3(6):827--864, 1992.

\bibitem[HR]{HerrlichRandecker}
F.~Herrlich and A.~Randecker.
\newblock Notes on the {V}eech group of the {C}hamanara surface.
\newblock arXiv:1810.05257.

\bibitem[HRR20]{HooperRafiRandecker2020}
W.~Patrick Hooper, Kasra Rafi, and Anja Randecker.
\newblock Renormalizing an infinite rational {IET}.
\newblock {\em Discrete Contin. Dyn. Syst.}, 40(9):5105--5116, 2020.

\bibitem[HS04]{HubertSchmidt04}
Pascal {Hubert} and Thomas~A. {Schmidt}.
\newblock {Infinitely generated Veech groups.}
\newblock {\em {Duke Math. J.}}, 123(1):49--69, 2004.

\bibitem[HS06]{HubertSchmidt}
Pascal {Hubert} and Thomas~A. {Schmidt}.
\newblock {An introduction to Veech surfaces.}
\newblock In {\em {Handbook of dynamical systems. Volume 1B}}, pages 501--526.
  Amsterdam: Elsevier, 2006.

\bibitem[HS08]{HerrlichSchmithusen08}
Frank Herrlich and Gabriela Schmith{\"u}sen.
\newblock An extraordinary origami curve.
\newblock {\em Math. Nachr.}, 281(2):219--237, 2008.

\bibitem[HS10]{HubertSchmithuesen10}
P.~Hubert and G.~Schmithüsen.
\newblock Infinite translation surfaces with infinitely generated {V}eech
  groups.
\newblock {\em J. of Modern Dyn.}, 4(4):715--732, 2010.

\bibitem[Hub06]{Hubbard-book1}
John~Hamal Hubbard.
\newblock {\em Teichm\"uller theory and applications to geometry, topology, and
  dynamics. {V}ol. 1}.
\newblock Matrix Editions, Ithaca, NY, 2006.
\newblock Teichm{\"u}ller theory, With contributions by Adrien Douady, William
  Dunbar, Roland Roeder, Sylvain Bonnot, David Brown, Allen Hatcher, Chris
  Hruska and Sudeb Mitra, With forewords by William Thurston and Clifford
  Earle.

\bibitem[Hur92]{Hurwitz1893}
A.~Hurwitz.
\newblock Ueber algebraische {G}ebilde mit eindeutigen {T}ransformationen in
  sich.
\newblock {\em Math. Ann.}, 41(3):403--442, 1892.

\bibitem[HW80]{HardyWeber80}
J.~Hardy and J.~Weber.
\newblock Diffusion in a periodic wind-tree model.
\newblock {\em J. of Math. Phys.}, 21(7):1802--1808, 1980.

\bibitem[HW12]{HooperWeiss09}
W.~Patrick Hooper and Barak Weiss.
\newblock Generalized staircases: recurrence and symmetry.
\newblock {\em Ann. Inst. Fourier (Grenoble)}, 62(4):1581--1600, 2012.

\bibitem[IT92]{ImayoshiTaniguchi}
Y.~Imayoshi and M.~Taniguchi.
\newblock {\em An introduction to {T}eichm\"uller spaces}.
\newblock Springer-Verlag, Tokyo, 1992.
\newblock Translated and revised from the Japanese by the authors.

\bibitem[JS14]{JohnsonSchmoll14}
Chris Johnson and Martin Schmoll.
\newblock Pseudo-{A}nosov eigenfoliations on {P}anov planes.
\newblock {\em Electron. Res. Announc. Math. Sci.}, 21:89--108, 2014.

\bibitem[{Kak}43]{Kakutani43}
Shizuo {Kakutani}.
\newblock {Induced measure preserving transformations.}
\newblock {\em {Proc. Imp. Acad. Tokyo}}, 19:635--641, 1943.

\bibitem[{Kar}20]{Karg2020}
Christoph {Karg}.
\newblock {\em On the Coarse Geometry of Infinite Regular Translation
  Surfaces}.
\newblock PhD thesis, Karlsruher Institut f\"{u}r Technologie (KIT), 2020.

\bibitem[{Kat}92]{Katok92}
Svetlana {Katok}.
\newblock {\em {Fuchsian groups}}.
\newblock Chicago: The University of Chicago Press, 1992.

\bibitem[Kea75]{Keane75}
M.~Keane.
\newblock Interval exchange transformations.
\newblock {\em Math. Zeitschr.}, 141:77--102, 1975.

\bibitem[Kes60]{Kesten-uniform_distribution}
Harry Kesten.
\newblock Uniform distribution {${\rm mod}\,1$}.
\newblock {\em Ann. of Math. (2)}, 71:445--471, 1960.

\bibitem[KMS86]{KerckhoffMasurSmillie86}
S.~Kerckhoff, H.~Masur, and J.~Smillie.
\newblock Ergodicity of billiard flows and quadratic differentials.
\newblock {\em Annals of Math.}, 124(2):293--311, 1986.

\bibitem[KS71]{KusunokiSainouchi71}
Y.~{Kusunoki} and Y.~{Sainouchi}.
\newblock {Holomorphic differentials on open Riemann surfaces.}
\newblock {\em {J. Math. Kyoto Univ.}}, 11:181--194, 1971.

\bibitem[KS00]{KenyonSmillie00}
Richard {Kenyon} and John {Smillie}.
\newblock {Billiards on rational-angled triangles.}
\newblock {\em {Comment. Math. Helv.}}, 75(1):65--108, 2000.

\bibitem[KZ75]{KatokZemliakov75}
A.~Katok and Z.~Zemliakov.
\newblock Topological transitivity of billiards in polygons.
\newblock {\em Maths notes}, 18:760--764, 1975.

\bibitem[Lei04]{Leininger04}
Christopher~J. Leininger.
\newblock On groups generated by two positive multi-twists: {T}eichm\"uller
  curves and {L}ehmer's number.
\newblock {\em Geom. Topol.}, 8:1301--1359 (electronic), 2004.

\bibitem[L{\"o}h17]{Loh17}
Clara L{\"o}h.
\newblock {\em Geometric group theory. {An} introduction}.
\newblock Universitext. Cham: Springer, 2017.

\bibitem[LT16]{LindseyTrevino16}
Kathryn {Lindsey} and Rodrigo {Trevi\~no}.
\newblock {Infinite type flat surface models of ergodic systems.}
\newblock {\em {Discrete Contin. Dyn. Syst.}}, 36(10):5509--5553, 2016.

\bibitem[M\"06]{Moeller06-periodic}
Martin M\"{o}ller.
\newblock {Periodic points on Veech surfaces and the Mordell-Weil group over a
  Teichm\"uller curve.}
\newblock {\em {Invent. Math.}}, 165(3):633--649, 2006.

\bibitem[M\"09]{Moeller09}
Martin M\"{o}ller.
\newblock Affine groups of flat surfaces.
\newblock In {\em Handbook of {T}eichm\"uller theory. {V}ol. {II}}, volume~13
  of {\em IRMA Lect. Math. Theor. Phys.}, pages 369--387. Eur. Math. Soc.,
  Z\"urich, 2009.

\bibitem[M\"13]{Moeller13}
Martin M\"{o}ller.
\newblock Teichm\"{u}ller curves, mainly from the viewpoint of algebraic
  geometry.
\newblock In {\em Moduli spaces of {R}iemann surfaces}, volume~20 of {\em
  IAS/Park City Math. Ser.}, pages 267--318. Amer. Math. Soc., Providence, RI,
  2013.

\bibitem[M{\'a}l14]{Malaga14}
A.~M{\'a}laga.
\newblock {\em \'Etude d’une famille de transformations pr\'eservant la
  mesure de $\Z\times\mathbb{T}$}.
\newblock PhD thesis, Universit\'e Paris Sud - Paris XI, 2014.
\newblock Fran\c{c}ais. (NNT:2014PA112413). (tel-01127218).

\bibitem[Mas82]{Masur82}
H.~Masur.
\newblock Interval exchange transformations and measured foliations.
\newblock {\em Annals of Math.}, 115:169--200, 1982.

\bibitem[Mas86]{Masur86}
Howard Masur.
\newblock Closed trajectories for quadratic differentials with an application
  to billiards.
\newblock {\em Duke Math. J.}, 53(2):307--314, 1986.

\bibitem[May43]{Mayer43}
A.~Mayer.
\newblock Trajectories on the closed orientable surfaces.
\newblock {\em Rec. Math. [Mat. Sbornik] N.S.}, 12(54):71--84, 1943.

\bibitem[{McM}00]{McMullen00}
Curtis~T. {McMullen}.
\newblock {Polynomial invariants for fibered 3-manifolds and Teichm\"uller
  geodesics for foliations}.
\newblock {\em {Ann. Sci. \'Ec. Norm. Sup\'er. (4)}}, 33(4):519--560, 2000.

\bibitem[McM03]{McMullen03}
C.~McMullen.
\newblock Billiards and {T}eichm\"uller curves on {H}ilbert modular surfaces.
\newblock {\em J. of Amer. Math. Soc.}, 16(4):857--885, 2003.

\bibitem[MT02]{MasurTabachnikov02}
H.~Masur and S.~Tabachnikov.
\newblock {\em Handbook of dynamical systems 1A}, chapter Rational billiards
  and flat structures, pages 289--307.
\newblock Elsevier, 2002.

\bibitem[MW89]{MoharWoess}
B.~Mohar and W.~Woess.
\newblock A survey on spectra of infinite graphs.
\newblock {\em Bull. London Math. Soc.}, 21:209--234, 1989.

\bibitem[MY10]{MatheusYoccoz10}
Carlos Matheus and Jean-Christophe Yoccoz.
\newblock The action of the affine diffeomorphisms on the relative homology
  group of certain exceptionally symmetric origamis.
\newblock {\em J. Mod. Dyn.}, 4(3):453--486, 2010.

\bibitem[MYZ14]{MatheusYoccozZmiaikou14}
Carlos {Matheus}, Jean-Christophe {Yoccoz}, and David {Zmiaikou}.
\newblock {Homologie des origamis avec sym\'etries.}
\newblock {\em {Ann. Inst. Fourier}}, 64(3):1131--1176, 2014.

\bibitem[{Nag}88]{Nag88}
Subhashis {Nag}.
\newblock {\em {The complex analytic theory of Teichm\"uller spaces}}.
\newblock New York etc.: Wiley, 1988.

\bibitem[Ore83]{Oren1983}
Ishai Oren.
\newblock Ergodicity of cylinder flows arising from irregularities of
  distribution.
\newblock {\em Israel J. Math.}, 44(2):127--138, 1983.

\bibitem[Pan09]{Panov}
D.~Panov.
\newblock Foliations with unbounded deviation on $t^2$.
\newblock {\em J. of Modern Dyn.}, 3:589--594, 2009.

\bibitem[Par]{Pardo-commutator}
A.~Pardo.
\newblock A remark on $\mathbb{Z}^d$-covers of {V}eech surfaces.
\newblock arXiv:1810.05257.

\bibitem[Par18]{Pardo-windtree}
Angel Pardo.
\newblock Counting problem on wind-tree models.
\newblock {\em Geom. Topol.}, 22(3):1483--1536, 2018.

\bibitem[PS81]{PhillipsSullivan81}
Anthony Phillips and Dennis Sullivan.
\newblock Geometry of leaves.
\newblock {\em Topology}, 20(2):209--218, 1981.

\bibitem[PSV11]{PrzytyckiSchmithuesenValdez11}
Piotr Przytycki, Gabriela Schmith{\"u}sen, and Ferr{\'a}n Valdez.
\newblock Veech groups of {L}och {N}ess monsters.
\newblock {\em Ann. Inst. Fourier (Grenoble)}, 61(2):673--687, 2011.

\bibitem[Ral14]{Ralston14-GenericDiscrepancy}
David Ralston.
\newblock Generic {$\frac12$}-discrepancy of {$\{n\theta+x\}$}.
\newblock {\em New York J. Math.}, 20:195--208, 2014.

\bibitem[Ran16]{Randecker16}
Anja Randecker.
\newblock {\em Geometry and topology of wild translation surfaces}, volume 151
  S. of {\em 10}.
\newblock KIT Scientific publishing, 2016.

\bibitem[Ray60]{Raymond60}
Frank Raymond.
\newblock The end point compactification of manifolds.
\newblock {\em Pacific J. Math.}, 10:947--963, 1960.

\bibitem[Ric63]{Richards63}
Ian Richards.
\newblock On the classification of noncompact surfaces.
\newblock {\em Trans. Amer. Math. Soc.}, 106:259--269, 1963.

\bibitem[RMV17]{Ramirez-Valdez16}
Camilo Ram\'{\i}rez~Maluendas and Ferr\'{a}n Valdez.
\newblock Veech groups of infinite-genus surfaces.
\newblock {\em Algebr. Geom. Topol.}, 17(1):529--560, 2017.

\bibitem[{Rok}48]{Rokhlin48}
V.~A. {Rokhlin}.
\newblock {Die ''allgemeine'' Transformation mit invariantem Mass ist keine
  Durchmischung.}
\newblock {\em {Dokl. Akad. Nauk SSSR, n. Ser.}}, 60:349--351, 1948.

\bibitem[Ros54]{Rosen54}
David Rosen.
\newblock A class of continued fractions associated with certain properly
  discontinuous groups.
\newblock {\em Duke Math. J.}, 21:549--563, 1954.

\bibitem[Sch78]{Schmidt78}
K.~Schmidt.
\newblock A cylinder flow arising from irregularity of distribution.
\newblock {\em Compositio Math.}, 36(3):225--232, 1978.

\bibitem[Sch06]{Schwartz06}
Richard~Evan Schwartz.
\newblock Obtuse triangular billiards. {I}. {N}ear the {$(2,3,6)$} triangle.
\newblock {\em Experiment. Math.}, 15(2):161--182, 2006.

\bibitem[Sch09]{Schwartz09}
Richard~Evan Schwartz.
\newblock Obtuse triangular billiards. {II}. {O}ne hundred degrees worth of
  periodic trajectories.
\newblock {\em Experiment. Math.}, 18(2):137--171, 2009.

\bibitem[Sch11]{Schwartz11}
Richard~Evan Schwartz.
\newblock {\em Mostly surfaces}, volume~60 of {\em Student Mathematical
  Library}.
\newblock American Mathematical Society, Providence, RI, 2011.

\bibitem[Ser77]{Serre77}
Jean-Pierre Serre.
\newblock {\em Linear representations of finite groups}.
\newblock Springer-Verlag, New York-Heidelberg, 1977.
\newblock Translated from the second French edition by Leonard L. Scott,
  Graduate Texts in Mathematics, Vol. 42.

\bibitem[Shi73]{Shields73}
Paul Shields.
\newblock Cutting and independent stacking of intervals.
\newblock {\em Math. Systems Theory}, 7:1--4, 1973.

\bibitem[SPWS17]{SchlagePuchtaWeitzeSchmithuesen-translations}
Jan-Christoph Schlage-Puchta and Gabriela Weitze-Schmith\"{u}sen.
\newblock Finite translation surfaces with maximal number of translations.
\newblock {\em Israel J. Math.}, 217(1):1--15, 2017.

\bibitem[Tab05]{Tabachnikov05}
Serge Tabachnikov.
\newblock {\em Geometry and billiards}, volume~30 of {\em Student Mathematical
  Library}.
\newblock American Mathematical Society, Providence, RI; Mathematics Advanced
  Study Semesters, University Park, PA, 2005.

\bibitem[Tay58]{Taylor58}
Angus~E. Taylor.
\newblock Introduction to functional analysis.
\newblock New {York}: {John} {Wiley} \& {Sons}, {Inc}.; {London}: {Chapman} \&
  {Hall}, {Ltd}. {XVI}, 423 p. (1958)., 1958.

\bibitem[Thu]{Thurston86}
W.~Thurston.
\newblock Minimal stretch maps between hyperbolic surfaces.
\newblock 1998 eprint of 1986 preprint. arXiv:9801039.

\bibitem[Thu88]{Thurston88}
William~P. Thurston.
\newblock On the geometry and dynamics of diffeomorphisms of surfaces.
\newblock {\em Bull. Amer. Math. Soc. (N.S.)}, 19(2):417--431, 1988.

\bibitem[Thu97]{Thurston97}
William~P. Thurston.
\newblock {\em Three-dimensional geometry and topology. {V}ol. 1}, volume~35 of
  {\em Princeton Mathematical Series}.
\newblock Princeton University Press, Princeton, NJ, 1997.
\newblock Edited by Silvio Levy.

\bibitem[Tro86]{Troyanov86}
Marc Troyanov.
\newblock Les surfaces euclidiennes \`a singularit\'es coniques.
\newblock {\em Enseign. Math. (2)}, 32(1-2):79--94, 1986.

\bibitem[Tro05]{Troubetzkoy2005}
Serge Troubetzkoy.
\newblock Periodic billiard orbits in right triangles.
\newblock {\em Ann. Inst. Fourier}, 55(1):29--46, 2005.

\bibitem[Tro10]{Troubetzkoy-typical_recurrence_ehrenfest}
Serge Troubetzkoy.
\newblock Typical recurrence for the {E}hrenfest wind-tree model.
\newblock {\em J. Stat. Phys.}, 141(1):60--67, 2010.

\bibitem[Val09a]{Valdez09_1}
Ferr\'{a}n Valdez.
\newblock Billiards in polygons and homogeneous foliations on {$\bold C^2$}.
\newblock {\em Ergodic Theory Dynam. Systems}, 29(1):255--271, 2009.

\bibitem[Val09b]{Valdez09_2}
Ferr\'{a}n Valdez.
\newblock Infinite genus surfaces and irrational polygonal billiards.
\newblock {\em Geom. Dedicata}, 143:143--154, 2009.

\bibitem[Val12]{Valdez12}
Ferr{\'a}n Valdez.
\newblock Veech groups, irrational billiards and stable abelian differentials.
\newblock {\em Discrete Contin. Dyn. Syst.}, 32(3):1055--1063, 2012.

\bibitem[Vee82]{Veech82}
W.~A. Veech.
\newblock {G}auss measures for transformations on the space of interval
  exchange maps.
\newblock {\em Annals of Math.}, 115:201--242, 1982.

\bibitem[Vee89]{Veech89}
W.~Veech.
\newblock {T}eichm\"uller curves in the moduli space, {E}isenstein series and
  an application to triangular billiards.
\newblock {\em Invent. Math.}, 97:533--683, 1989.

\bibitem[Ver82]{Vershik-MarkovCompacta}
A.~M. Ver{\v{s}}ik.
\newblock A theorem on periodical {M}arkov approximation in ergodic theory.
\newblock In {\em Ergodic theory and related topics ({V}itte, 1981)}, volume~12
  of {\em Math. Res.}, pages 195--206. Akademie-Verlag, Berlin, 1982.

\bibitem[Via06]{Viana06}
Marcelo Viana.
\newblock Ergodic theory of interval exchange maps.
\newblock {\em Rev. Mat. Complut.}, 19(1):7--100, 2006.

\bibitem[VO16]{KrerleyViana16}
Marcelo {Viana} and Krerley {Oliveira}.
\newblock {\em {Foundations of ergodic theory.}}, volume 151.
\newblock Cambridge: Cambridge University Press, 2016.

\bibitem[VWS14]{SchmithuesenValdez14}
Ferr\'{a}n Valdez and Gabriela Weitze-Schmith\"{u}sen.
\newblock On the geometry and arithmetic of infinite translation surfaces.
\newblock {\em J. Singul.}, 9:226--244, 2014.

\bibitem[Wan22]{Wang21}
Jane Wang.
\newblock The realization problem for dilation surfaces.
\newblock {\em Int. Math. Res. Not.}, 2022(21):16672--16708, 2022.

\bibitem[Wri15]{Wright15}
Alex Wright.
\newblock Translation surfaces and their orbit closures: an introduction for a
  broad audience.
\newblock {\em EMS Surv. Math. Sci.}, 2(1):63--108, 2015.

\bibitem[Yoc10]{Yoccoz10}
Jean-Christophe Yoccoz.
\newblock Interval exchange maps and translation surfaces.
\newblock In {\em Homogeneous flows, moduli spaces and arithmetic}, volume~10
  of {\em Clay Math. Proc.}, pages 1--69. Amer. Math. Soc., Providence, RI,
  2010.

\bibitem[ZK75]{ZemljakovKatok-transitivity1}
A.~N. Zemljakov and A.~B. Katok.
\newblock Topological transitivity of billiards in polygons.
\newblock {\em Mat. Zametki}, 18(2):291--300, 1975.

\bibitem[Zor06]{Zorich06}
A.~Zorich.
\newblock {\em Flat surfaces}.
\newblock Springer, 2006.
\newblock in "Frontiers in Number Theory, Physics and Geometry. Volume 1: On
  random matrices , zeta functions and dynamical systems".

\end{thebibliography}

\end{document}